\documentclass[11pt,reqno]{amsart}

\usepackage{a4wide}

\usepackage{amssymb}
\usepackage[mathscr]{eucal}
\usepackage{enumerate}
\usepackage[all]{xy}
\ifx\pdfoutput\undefined \xyoption{dvips}\fi

\newtheorem{thm}{Theorem}
\newtheorem{lemma}[thm]{Lemma}
\newtheorem{cor}[thm]{Corollary}

\theoremstyle{definition}
\newtheorem{defn}[thm]{Definition}
\newtheorem{notation}[thm]{Notation}
\newtheorem{obs}[thm]{Observation}
\newtheorem{recall}[thm]{Recall}

\newenvironment{block}{{}{}}

\newcount\domsdepth
\newenvironment{domsitem}{%
    \advance\domsdepth by 1 %
    \begin{list}{%
        \ifnum \domsdepth = 1 $\bullet$ \else $-$ \fi } %
      {
        \setlength{\leftmargin}{1.5em}
        \setlength{\itemsep}{1ex}
        \setlength{\labelwidth}{1em}
        \setlength{\labelsep}{0.5em}
        \setlength{\itemindent}{0em}
     } %
  }{ %
  \end{list}}%

\def\parhead[#1]{\vspace{1ex}%
\noindent{\bf\boldmath #1.}}

\def\mlaux#1#2#3{\setbox0=\hbox{$\mathsurround=0pt #2{#3}$}%
  \dimen0=\dp0\advance\dimen0 by \ht0\lower#1\dimen0\box0}
\def\mlower#1#2{\mathpalette{\mlaux{#1}}{#2}}

\def\killdp#1{\setbox0=\hbox{#1} \dp0=0pt \box0}
\def\category#1{\underline{\mathchoice{\killdp{\rm\normalsize #1}}%
    {\killdp{\rm\normalsize #1}}{\killdp{\rm\scriptsize #1}}%
    {\killdp{\rm\tiny #1}}}}

\def\makellapm#1#2{\hbox to 0pt{\hss$\mathsurround=0pt #1{#2}$}}

\def\makerlapm#1#2{\hbox to 0pt{$\mathsurround=0pt #1{#2}$\hss}}
\def\rlapm{\relax\mathpalette\makerlapm}

\def\bicat#1{\underline{\mathcal{#1}}}
\def\twocat#1{\bicat{\text{\rm\em #1}}}
\let\tcat=\bicat

\def\OSimp{\tilde\Delta}

\def\coliminj{i}
\def\cprprod{c}
\def\cprprods{\bar{c}}

\def\orelse{\mathrel{\text{ or }}}
\def\also{\mathrel{\text{ and }}}

\def\ifundef#1{\expandafter\ifx\csname#1\endcsname\relax}

\ifundef{Top}
  \DeclareSymbolFont{dsymbolsC}{U}{txsyc}{m}{n}
  \SetSymbolFont{dsymbolsC}{bold}{U}{txsyc}{bx}{n}
  \DeclareFontSubstitution{U}{txsyc}{m}{n}

  \def\re@DeclareMathSymbol#1#2#3#4{%
      \let#1=\undefined
      \DeclareMathSymbol{#1}{#2}{#3}{#4}}

  \re@DeclareMathSymbol{\Top}{\mathord}{dsymbolsC}{120}
  \re@DeclareMathSymbol{\Bot}{\mathord}{dsymbolsC}{121}
\fi

\let\partinj=\Bot
\let\partproj=\Top

\def\mathonl#1{#1}

\def\im{\mathonl{\text{\rm im}}}
\def\coim{\mathonl{\text{\rm coim}}}

\def\op#1{#1^{\mathord{\text{\rm op}}}}
\def\co#1{#1^{\mathord{\text{\rm co}}}}
\def\coop#1{#1^{\mathord{\text{\rm coop}}}}
\def\vop#1{#1^{\mathord{\text{\rm vop}}}}
\def\hop#1{#1^{\mathord{\text{\rm hop}}}}
\def\hvop#1{#1^{\mathord{\text{\rm hvop}}}}
\def\refld#1{#1^{\mathord{\text{\rm r}}}}

\let\boundary=\partial

\def\cod{\mathonl{\text{\rm cod}}}
\def\t{\mathonl{\boldsymbol{t}}}

\def\dom{\mathonl{\text{\rm dom}}}
\def\s{\mathonl{\boldsymbol{s}}}

\def\comp{\mathbin{\boldsymbol{*}}}
\def\pb#1#2{\mathbin{{}_{#1}\mathord\times_{#2}}}

\def\idt{\mathonl{\boldsymbol{i}}}
\def\cotens{\mathonl{\pitchfork}}

\def\suptilde#1{{#1}^{\sim}}
\def\suphat#1{{#1}^{\wedge}}

\let\triangle=\triangleleft
\let\filledtri=\blacktriangleleft

\def\oriental{\mathonl{\mathcal{O}}}
\def\norient{\mathonl{\mathcal{N}}}

\def\functor#1{\mathonl{\mathrm{#1}}}
\let\function=\functor

\def\Cosk{{\functor{Cosk}}}
\def\Cosup{{\functor{Cos}}}

\def\Th{{\functor{Th}}}
\def\Dec{{\functor{Dec}}}

\def\Sup_#1(#2){\left|#2\right|_{#1}}
\def\Sk_#1(#2){\left\|#2\right\|_{#1}}

\def\Al#1{#1\text{-}\category{Alg}}
\def\CoAl#1{#1\text{-}\category{CoAlg}}

\def\rhv{^r}
\def\lhv{^l}

\def\ihoml#1{\homout{#1}\lhv}
\def\ihomr#1{\homout{#1}\rhv}
\def\homout#1{\functor{K}_{#1}}
\let\kan=\homout
\def\kanladj#1{\functor{F}_{#1}}

\def\Alg<#1;#2>{\Al{#1}(#2)}
\def\CoAlg<#1;#2>{\CoAl{#1}(#2)}

\def\LCoAlg<#1;#2>{\CoAl{#1}_L(#2)}

\let\decact=\odot
\def\Ndec{\nerv_{\Dec,\cpt_0}}

\def\ECat#1{{#1}\text{-}\Cat}

\def\CompCat{\ECat\Comp}

\def\dim{\function{dim}}

\def\lim{\functor{lim}}
\def\colim{\functor{colim}}

\def\obj{\functor{obj}}
\def\arr{\functor{arr}}
\def\sq{\functor{sq}}

\def\dis{\function{dis}}
\def\cpt{\mathonl{\Pi}}

\let\svec=\underline
\def\id{\mathord{\text{\rm id}}}

\def\inerve{\nerv_\inf}
\def\inladj{\functor{F}_\inf}
\def\nnerve{\nerv_n}
\def\nnladj{\functor{F}_n}
\def\groth{\mathbb{G}}
\def\interv{\functor{I}}
\def\interdc{\mathbb{I}}
\def\pathcat{\mathbb{P}}
\def\pathecat{\functor{P}}
\def\recons{\functor{R}}

\let\radj=\recons
\def\forget{\functor{U}}
\def\refl{\functor{L}}
\def\dgm{{D}}
\def\ladj{\functor{F}}
\def\gfunc{\functor{G}}
\let\ffunc=\ladj

\def\efunc{\functor{E}}
\def\kfunc{\functor{K}}
\def\lfunc{\functor{L}}
\def\incl{\functor{I}}
\def\nerv{\functor{N}}

\def\pathsq{\pathcat^{\scriptscriptstyle 2}}
\def\pathesq{\pathecat^{\scriptscriptstyle(2)}}

\def\xcube#1{\mathonl{1\text{\rm -Cube}}}

\def\squares{\mathonl{\text{\rm Sq}}}

\def\glob{\mathonl{\text{\rm Glob}}}

\def\Cpathcat{\Delta_{\pathcat}}

\let\iend=\partial

\let\diag=\nabla

\def\adjoint{\dashv}

\def\yoneda#1{\ulcorner{#1}\urcorner}

\def\Yoneda{\mathscr{Y}}

\def\defeq{\mathrel{\mlower{0.15}{\stackrel{\text{def}}=}}}

\def\p@ir#1<#2;#3>{#1{{\mskip 1.5mu}#2,\,#3{\mskip 1.5mu}}}
\def\br@ck#1{({#1})}

\def\pair{\p@ir\br@ck}
\def\satomp{\p@ir\satom}
\def\statomp{\p@ir\statom}
\def\atomp{\p@ir\atom}
\def\triple<#1;#2;#3>{({\mskip 1.5mu}#1,\,#2,\,#3{\mskip 1.5mu})}
\def\quadruple<#1;#2;#3;#4>{({\mskip 1.5mu}#1,\,#2,\,#3,\,#4{\mskip 1.5mu})}
\def\quintuple<#1;#2;#3;#4;#5>{({\mskip
    1.5mu}#1,\,#2,\,#3,\,#4,\,#5{\mskip 1.5mu})}
\def\sextuple<#1;#2;#3;#4;#5;#6>{({\mskip 1.5mu}#1,\,#2,\,#3,\,#4,\,#5,\,#6{\mskip 1.5mu})}
\def\septuple<#1;#2;#3;#4;#5;#6;#7>{({\mskip 1.5mu}#1,\,#2,\,#3,\,#4,\,#5,\,#6,\,#7{\mskip 1.5mu})}
\def\ntuple<#1;#2>{({\mskip 1.5mu}#1_1,\,#1_2,\,\dots,\,#1_{#2}{\mskip 1.5mu})}

\let\parvone=\varphi
\let\parvtwo=\xi

\let\sc=\mathbb
\def\scobj#1{\sc{#1}_o}
\def\scarr#1{\sc{#1}_a}

\let\dc=\mathbb
\def\dcobj#1{\dc{#1}_o}
\def\dcvarr#1{\dc{#1}_v}
\def\dcharr#1{\dc{#1}_h}
\def\dcsq#1{\dc{#1}_s}

\def\atom#1{\langle{#1}\rangle}
\def\satom#1{\lbrack\!\lbrack{#1}\rbrack\!\rbrack}
\let\statom=\satom

\def\cattuple#1{\sextuple<\scarr{#1};\scobj{#1};\comp;\idt;\s;\t>}
\def\dcattup<#1;#2;#3>{\sextuple<#1;#2;\comp_{#3};\idt_{#3};\s_{#3};\t_{#3}>}

\newdir{ >}{{}*!/-10pt/@{>}}
\newdir{u(}{{}*!/-5pt/@^{(}}
\newdir{|>}{%
  !/4.5pt/@{|}*:(1,-.2)@^{>}*:(1,+.2)@_{>}*!/-3pt/@{ }}

\def\arrow#1:#2->#3.{{#1\colon #2}\xy (0,0)*+{} %
  \ar (7,0)*+{}\endxy{#3}}

\def\epi#1:#2->#3.{{#1\colon #2}\xy (0,0)*+{} %
  \ar@{->>} (7,0)*+{}\endxy{#3}}

\def\overepi#1:#2->#3.{{#2}\xy (0,0)*+{} %
  \ar@{->>}^-{\smash{#1}} (7,0)*+{}\endxy{#3}}

\def\cover#1:#2->#3.{{#1\colon #2}\xy (0,0)*+{} %
  \ar@{-|>} (7,0)*+{}\endxy{#3}}

\def\overinc#1:#2->#3.{{#2}\xy (0,0)*+{} %
  \ar@{u(->}^-{\smash{#1}} (9,0)*+{}\endxy{#3}}

\def\overarr#1:#2->#3.{{#2}\xy (0,0)*+{} %
  \ar^-{\smash{#1}} (7,0)*+{}\endxy{#3}}

\def\spanarr#1:#2->#3.{{#1\colon #2}\xy (0,0)*+{} %
  \ar|-{\object@{|}} (9,0)*+{}\endxy{#3}}

\def\inc#1:#2->#3.{{#1\colon #2}\xy (0,0)*+{} %
  \ar@{u(->} (9,0)*+{}\endxy{#3}}

\def\nattrans#1:#2->#3.{{#1\colon #2}\xy (0,0)*+{} %
  \ar^-{.} (7,0)*+{}\endxy{#3}}

\def\iso#1:#2->#3.{{#1\colon #2}\xy (0,0)*+{} %
  \ar^-{\smash{\simeq}} (7,0)*+{}\endxy{#3}}

\def\equiv#1:#2->#3.{{#1\colon #2}\xy (0,0)*+{} %
  \ar^-{\smash{\sim}} (7,0)*+{}\endxy{#3}}

\def\anoniso#1->#2.{{#1}\xy (0,0)*+{} %
  \ar^-{\simeq} (7,0)*+{}\endxy{#2}}

\def\twocell#1:#2->#3.{{#1\colon #2}\xy (0,0)*+{} %
  \ar@{=>} (7,0)*+{}\endxy{#3}}

\def\isotwocell#1:#2->#3.{{#1\colon #2}\xy (0,0)*+{} %
  \ar@{=>}^-{\sim} (7,0)*+{}\endxy{#3}}

\def\funcat[#1,#2]{[{#1},{#2}]}
\def\cyl[#1,#2]{\Cylinder\pair<{#1};{#2}>}

\def\tDelta{t\Delta}
\def\Simp{\category{Simp}}
\def\Set{\category{Set}}
\def\SSimp{\category{SSimp}}
\def\Cat{\category{Cat}}

\def\Strat{\category{Strat}}
\def\Simplex{\category{Simplex}}
\def\Double{\category{Double}}
\def\Conn{\category{Conn}}
\def\Twocat{\category{2-Cat}}
\def\Mon{\category{Mon}}
\def\Pt{\category{Pt}}
\def\Ord{\category{Ord}}

\def\GrSet{\category{GrSet}}
\def\Parity{\category{Parity}}
\def\SParity{\category{SParity}}
\def\InfCatCat{\ECat{(\InfCat)}}

\def\Monoidal{\twocat{Mon}}
\def\Category{\twocat{Cat}}
\def\Cylinder{\twocat{Cyl}}

\def\inf{\ensuremath\omega}
\def\InfCat{\ensuremath{\ECat\inf}}
\def\nCat{\ensuremath{\ECat{n}}}

\def\Precomp{\category{Pcs}}
\def\Comp{\category{Cs}}

\def\aDelta{\Delta_{\mathord{+}}}
\def\fDelta{\Delta_f}
\def\afDelta{\Delta_{f\mathord{+}}}

\def\face{\delta}
\def\vertex{\varepsilon}
\def\degen{\sigma}
\def\tdegen{\varsigma}
\def\thop{\varrho}

\def\cface{\partial}
\def\cdegen{\epsilon}

\def\stratd#1{\pair<\mathbb{#1};t\mathbb{#1}>}
\def\strat#1{\pair<#1;t{#1}>}

\let\anytens=\circledcirc
\let\pretens=\boxtimes
\def\soplus{+}

\let\dsum=\aoplus

\def\stratc#1{{#1}^{\ast}}
\let\simpc=\stratc

\let\domexists=\exists
\def\exists{\domexists\,}

\def\shufflel#1{\alpha_{#1}}
\def\shuffler#1{\beta_{#1}}
\def\shuffle#1{\pair<\shufflel{#1};\shuffler{#1}>}
\def\topop#1{\gamma_{#1}}

\def\abpair{\pair<\alpha;\beta>}

\def\pocorner{\hbox to 10pt{{\vrule height10pt depth0pt width0.5pt}%
    \vbox to 10pt{{\hrule height0.5pt width9.5pt depth0pt}\vfill}}}
\def\poexcursion{\save[]-<12pt,-12pt>*{\pocorner}\restore}
\def\pbcorner{\vbox to 0pt{\kern 5pt\hbox to 0pt{\kern 5pt%
      \vbox{{\hrule height0.5pt width9.5pt depth0pt}}%
      {\vrule height10pt depth0pt width0.5pt}\hss}\vss}}
\def\pbexcursion{\save[]-<-2pt,2pt>*{\pbcorner}\restore}

\title{Complicial Sets}
\author[Verity]{Dominic Verity}
\address{
  Centre of Australian Category Theory \\
  Macquarie University \\
  NSW 2109 \\
  Australia
}
\email{domv@ics.mq.edu.au}
\date{October 15, 2004}

\begin{document}

\begin{abstract}
  The primary purpose of this work is to characterise strict
  \inf-categories as simplicial sets with structure. We prove the
  Street-Roberts conjecture which states that they are exactly the
  ``complicial sets'' defined and named by John Roberts in his
  handwritten notes of that title \cite{Roberts:1978:Complicial}.
\end{abstract}

\maketitle
\tableofcontents
\newpage


\section{Introduction}

\subsection{Historical Background}

This work presents a proof of a result that has sometimes been
referred to as the Street-Roberts conjecture
(after~\cite{Street:1987:Oriental}). This postulates an equivalence
between the category of {\em strict\/} \inf-categories and a category
of structures called {\em complicial sets\/} which are certain kinds
of enhanced simplicial sets originally studied by
Roberts~\cite{Roberts:1978:Complicial}.

The genesis of this work dates back to the mid-1970 and Roberts' work
on non-abelian cohomology. His original interest in this topic grew
from his conviction that (strict) \inf-categories were the appropriate
algebraic structures within which to value such theories
\cite{Roberts:1977:Complicial}. This led him to define complicial sets
to be simplicial sets with distinguished elements, which he originally
referred to as ``neutral'' then later as ``hollow'' but for which we
prefer the term {\em thin\/} (after~\cite{Dakin:1977:PhD}), satisfying
some natural conditions related to those that characterise {\em Kan
  complexes\/} in the homotopy theory of simplicial sets
\cite{GabrielZisman:1967:CFHT}.

In particular, his conditions include a certain kind of {\em unique
  thin horn filler\/} condition. These closely resemble the horn
fillers of Kan, with the notable exceptions that only certain {\em
  admissible\/} horns (for which a specified class of faces are
required to be thin) are assumed to have fillers, that those fillers
are all themselves thin and that all such fillers are assumed to be
unique.
  
Roberts was motivated in this definition by his observation that it
should be possible to naturally generalise the classical nerve
constructions of Algebraic Topology to provide a functor from the
category of (strict) \inf-categories to the category of simplicial
sets which, in a suitable sense, encapsulated a natural notion of {\em
  higher non-abelian cocycle}. While the nerve construction on
groups and partially ordered sets, and their common generalisation to
categories, is well known and easily described to students of
Algebraic Topology, the same cannot be said of its generalisation to
\inf-categories. Indeed its very definition poses substantial
technical challenges which eluded Roberts at the time.

Even without an explicit construction of this nerve functor, Roberts
set about studying those simplicial sets that would occur as the nerve
of some \inf-category. He observed that this functor would not provide
a fully-faithful representation and suggested rectifying this failing
by introducing thin elements into his study. He then set to
characterising those augmented simplicial sets which would arise in
the replete image of the postulated fully-faithful functor and thus
complicial sets were born. 

One might forgive Roberts if his identification of this category were
to prove deficient in some way, after all he was working without a
nerve functor and at a time when the theory of strict \inf-categories
had advanced little further than fundamental definitions.  It is
therefore a great tribute to Roberts' insight that our main theorem
here (theorem~\ref{street.roberts.thm} of
section~\ref{street.nerve.sec}) establishes his sought-for equivalence
under a definition of complicial set which only differs from his in
manner of expression. He provided a precise conjecture in 1978, our
only contribution has been to prove it as stated!

The next big advance in this story came with Street's papers on
orientals \cite{Street:1987:Oriental} and parity complexes
\cite{Street:1991:Parity}. He was introduced to this study,
during a visit Roberts made to Sydney in the (southern hemisphere) summer
of 1977-8, and was captivated by it. While he was quickly able to
establish the sought for equivalence at dimension 2, it became 
apparent to him that establishing the general equivalence would be
quite a difficult problem. Consequently he decided to concentrate on
providing a rigorous definition of the nerve of an \inf-category.

After a couple of false starts to that end, he soon realised that the
crux of the matter was to define the free \inf-category $\oriental_n$
on an $n$-simplex. However, in order to make this insight precise it
was necessary to define exactly the kinds of combinatorial structures
that might give rise to free \inf-categories and even to explicate the
sense in which freeness itself should be interpreted. Furthermore,
once these concepts had been defined it would still be necessary to
define exactly how an $n$-simplex might be considered to be such a
structure.

The first of these questions elicited the introduction of certain
kinds of inductively defined combinatorial structures called {\em
  \inf-computads}, which provide the desired definitions but are often
inconvenient to calculate with. To facilitate the work we are engaged
in here, Street later introduced a restricted form of \inf-computad,
called a {\em parity complex}. These satisfy some very strong
loop-freeness conditions designed to ensure that we may describe the
cells of the associated free \inf-category as (pairs of) subsets of
the parity complex itself.

To understand how one might render an $n$-simplex as a parity complex,
Street started by observing that $\oriental_n$ should play the part of
some form of non-abelian $n$-cocycle in the sought after cohomology
theory. This led him to conjecture that we could take insight from
abelian cohomology and build a consistently oriented parity complex
whose elements were the faces of a standard $n$-simplex, each of which
should be oriented as a cell to ``map from odd numbered faces to even
numbered ones''. Under this definition, he was able to demonstrate that
the resulting structure satisfied the strong loop-freeness conditions
required of a parity complex and thereby provide a completely explicit
description of the free \inf-category $\oriental_n$.

Street completed the construction of his \inf-categorical nerve by
enriching his orientals to a functor from the {\em category of
  simplicial operators\/} $\Delta$ and applying Kan's construction
\cite{Kan:1958:Adjoint}. His paper \cite{Street:1987:Oriental} goes on
to lift this functor to map \inf-categories to simplicial sets with
thin elements (which we call stratified sets in the sequel) and to
formalise Roberts' original conjecture in this context. He
identifies a family of admissible horns which {\em strictly\/}
contains the class introduced by Roberts and in
\cite{Street:1988:Fillers} he demonstrates that within the nerve of an
\inf-category such horns do indeed have unique thin fillers. It is
worth noting, however, that while Street's complicial set definition
appears, at first sight, to be stronger than Roberts', our proof here
demonstrates that they are in fact equipotent.

My own contribution to this story began in 1991 when I first read
\cite{Street:1987:Oriental} and immediately became captivated by this
problem. In particular, I had been searching for an approach to
defining structures which we now call {\em weak \inf-categories\/} and
became convinced that Street's parenthetical remarks regarding a
potential simplicial definition held great merit. In ignorance of
\cite{Street:1988:Fillers}, I set about proving that in any nerve
Street's admissible horns had unique fillers and in a very short time,
and to my great pleasure, I succeeded in constructing a decomposition
of each oriental which established this result.

Spurred on by my initial success, I decided to turn my hand to proving
Street's version of Roberts' original conjecture and soon succeeded in
showing that Street's nerve functor was fully faithful. While this
result in itself was a clear contribution to the then extant state of
knowledge in this area, I was not satisfied since the method used to
do so had originally been conceived as a proof of the complete
Street-Roberts conjecture. My initial reaction to this block was to
cast around for a categorical abstraction which might strengthen my
result, but to no avail. Somewhat disheartened by this, I wrote notes
on my proof to date and circulated this handwritten manuscript to a
few interested parties (including Street) before settling down to
complete my PhD work on enriched category theory.

Returning to this problem in 1993, I realised that an argument based
upon the decalage construction would close the gap in my proof. In
1994 I gave talks on the resulting proof to the Sydney (now
Australian) Category Seminar, at the University College of North Wales
in Bangor, at the Mathematical Sciences Research Institute in Berkley
and to a Peripatetic Seminar on Sheaves and Logic at the Newton
Institute in Cambridge UK.

Unfortunately, however, career events overtook me before I had time to
commit the final proof to paper. In early 1995 I followed up a series
of consultancy engagements in the financial markets by accepting a
full time role in investment banking. I spent the subsequent 5 years
as a derivative securities analyst and trader manager before returning
to academia in mid 2000.

On returning to this field, two things struck me immediately. Firstly,
and most pleasingly, an area of study that had given so much joy to a
small band of enthusiasts had grown into a dynamic area of wide
debate and interest, driven by an influx of new ideas and
approaches. Secondly, and much to my relief, nobody appeared to have
provided a proof of the Street-Roberts result upon which I had devoted
so much time prior to my sojourn into the business world. 

The proof presented herein is fundamentally no different to the one
that I spoke on in 1994. The primary innovation on that original work
has been the adoption of a {\em (lax) Gray tensor product\/} of
complicial sets as a unifying organisational and constructional tool.

I should like to dedicate this work to three groups of people. Firstly
of course to my wife, Sally, and children, Lottie and Florrie, who
have put up with a husband and father incessantly crouched over the
word processor. They have never once questioned the importance to me
of completing this work and have sustained me in body and soul over
the 12 years it has been in gestation. Secondly I would like to thank
Ross Street, who has been the most formative influence in my
development as a mathematician. It was his friendship, support and
inspiration which convinced me to return to academe and to work which
I had long since convinced myself I would never commit to paper. In
the end it was his quiet reminders, to the effect that I owed it to
the community to write up my ideas, that spurred me to find the time
in a busy schedule to write these scratchings. Finally I should like
to thank the staff of the Postgraduate Professional Development
Program of the Division of Information and Communication Sciences at
Macquarie University, of which I am the academic director.  Without
their support and hard work, filling in for me while I wrote this
work, I would never have been able to find the silence and space to
organise these ideas.

\subsection{Relationships to Other Work}

While this work predominantly interests itself in a study of strict
\inf-categorical structures, it is nevertheless squarely motivated by
a broader program to define and study weakened \inf-categorical
notions. In particular, its purpose from the very start has been to
act as a first step toward providing a full account of the theory of
Street's simplicial weak \inf-category notion, which first appeared in
sketch form in his orientals paper~\cite{Street:1987:Oriental} and
upon which he later elaborated in~\cite{Street:2003:WomCats}. 

Indeed, much of the work on (lax) Gray tensor product of complicial
sets given here routinely generalises to Street's category of weak
complicial sets. In particular, using these generalisations we may
show that the category of weak complicial sets supports biclosed
structures generalising those on the category of bicategories and
homomorphisms as discussed in \cite{Street:1980:FibBicat}. These allow
us to enrich the category of weak complicial sets over itself in 
natural ways and opens the possibility of providing a coherence result
for Street's weak \inf-categories along the lines of that described
for tricategories in~\cite{Gordon:1995:Tricats}. 

Notice however that, just as in the bicategorical and tricategorical
cases, the corresponding ``monoidal structure'' on weakly complicial
sets is only {\em weakly\/} adjoint (in some suitable sense) to this
biclosed structure. As a result it is only weakly coherent, making it
somewhat inconvenient to calculate with. One convenience provides by
the weak complicial approach, however, is that we can always perform
all required calculations in a bigger category of pre-complicial sets,
on which the corresponding tensor is actually part of a genuine
monoidal biclosed structure.

The reader will, of course, be aware that over the past ten years a
wide variety of weak \inf-categorical notions have been proposed by
various authors. Notable amongst these are those of
Joyal~\cite{Joyal:1997:ThetaCats}, Batanin~\cite{Batanin:1998:WomCats}
and Baez and Dolan~\cite{Baez:1998:HDA}, although more may be found in
the literature or in Leinster's survey of the current state of the
definitional art~\cite{Leinster:2002:Survey}. Many of these drew their
initial inspiration from Street's parenthetical remarks
in~\cite{Leinster:2002:Survey} although none of them, except for the
summary given in Leinster~\cite{Leinster:2002:Survey} and Street's own
account~\cite{Street:2003:WomCats}, expand directly on the purely
simplicial approach analysed here. This work starts from a point of
view neatly summed up by Street's comment
in~\cite{Street:2003:WomCats} that a simplicial formulation bears the
distinct advantage that:

\begin{quotation}
  \em Simplicial sets are lovely objects about which algebraic
  topologists know a lot. If something is described as a simplicial
  set, it is ready to be absorbed into topology. Or, in other words,
  no matter which definition of weak \inf-category eventually becomes
  dominant, it will be valuable to know its simplicial nerve.
\end{quotation}

In many ways, our work here also parallels work of Ronald Brown and
his coworkers, which established analogous results for the somewhat
simpler groupoidal case.  While our work and theirs share at least one
common motivation, to develop a coherent and complete theory of higher
dimensional cohomology, their primary interest is not the development
of an encompassing theory of weak \inf-categories but rather the
explication of a theory within which to make explicit calculations of
homotopic invariants.  To that end a study of groupoids rather than
categories is adequate for most purposes.

The complicial sets discussed here are most closely related to their
simplicial T-complexes, which were first discussed in Dakin's
thesis~\cite{Dakin:1977:PhD}, where the term ``thin'' was first
coined, and later developed and popularised in the work of Brown,
Higgins and others. Since this work was groupoid oriented, the
definitions involved are somewhat simpler than those discussed here
and a proof of the equivalence between simplicial T-complexes and
\inf-groupoids is far more easily attained. Indeed, a proof of this
equivalence in the rank 2 case dates to Dakin's thesis in 1977.

While the simplicial approach is important to their work, it should be
pointed out that in general their preference has been to work in a
cubical context. This allows them to more easily build generalised van
Kampen Theorems, for computing homotpical invariants, and to define
monoidal closed structures for analysing homotopy classes of maps.  A
comprehensive and up-to-date account of this work may be found in
Brown's excellent survey article~\cite{Brown:2004:Survey}.  In this
cubical program the closest result to the one presented here is that
of Al-Agl, Brown and Steiner~\cite{AlAgl:2002:Cubical} which
establishes an equivalence between the category of (strict)
\inf-categories and a category of multiple (cubical) categories with
connection.

While this cubical approach is attractive for the reasons discussed
above, it is our conviction that it does not provide quite such a
convenient context within which to develop a theory of weak
\inf-categories.  Nevertheless cubical calculations, and most
particularly those involved in a theory of tensor products, remain
vital even in the simplicial context. The theory discussed here,
therefore, appears to combine the best of both worlds by allowing for
a theory with all of the beauty, elegance and economy of the
simplicial approach while at the same time providing a simple and
highly explicit description of the cubical Gray tensor product.

\subsection{Overview and Structure}

This work is almost exclusively devoted to proving the Street-Roberts
conjecture, in the form presented in~\cite{Street:1987:Oriental}. On
the way we substantially develop Roberts' theory of complicial
sets~\cite{Roberts:1977:Complicial} itself and make some contributions
to Street's theory of parity complexes~\cite{Street:1991:Parity}.  In
particular, we study a new monoidal closed structure on the category
of complicial sets which we show to be the appropriate generalisation
of the (lax) Gray tensor product of 2-categories to this context.
Under the equivalence conjectured by Street
\cite{Street:1987:Oriental} and Roberts
\cite{Roberts:1977:Complicial}, which we prove here, this tensor
product coincides with those of Crans~\cite{Crans:1995:PhD},
Steiner~\cite{Steiner:1991:Tensor} and others.

From the outset, it has been designed to be as self contained as
possible. While much of the material covered in its first five
sections is classical in nature, it has been presented here in order
to fix our notation for the sequel and to aggregate together a number
of familiar (and not so familiar) pieces of Algebraic Topology and
Category Theory from a diverse range of sources.  Given the influx
into this field of mathematicians from diverse backgrounds, it was
felt prudent that no prerequisite assumptions be made. Most
particularly, it was recognised that some readers might not be fully
conversant with certain of the more abstract aspects of general
(enriched and internal) category theory and higher category theory.

In particular, where categorical abstractions are concerned we assume
little more than that the reader should have a general grounding in
basic categorical concepts, such as functor, natural transformation,
limit, colimit, comma category and 2-category etcetera, all of which
may be gleaned from Mac~Lane's book \cite{Maclane:1971:CWM}.

Congruent with this philosophy, sections~\ref{simp.sect}
to~\ref{decintro.sect} are all contextual in nature, providing fairly
standard presentations of traditional material.
Section~\ref{simp.sect} consists of a brief introduction to the theory
of simplicial sets, up to and including the theory of shuffles.  While
this is {\bf not} intended to provide an exhaustive treatment of the
simplicial algebra necessary to read this work, it should provide
most of the necessary background and adequate pointers to the
available literature.

At some points in the sequel our arguments are substantially
simplified by couching them in more abstract categorical terms, along
the lines described in Kelly's book~\cite{Kelly:1982:ECT}. In
particular, our construction of the monoidal biclosed Gray tensor
structure on complicial sets relies on Day's reflection theorem for
monoidal biclosed categories~\cite{Day:1970:ClosedFunc} and many of
our later constructions and calculations are couched in terms of left
exact theories and their coalgebras. While a thorough reading of
Kelly's book would handsomely repay the effort involved, the results
of greatest interest here are collected together in
section~\ref{cat.back.sect}; the reader should refer to the cited
literature for detailed proofs.

Section~\ref{sect.ncats} rehearses the basic definitions in the theory
of (internal) categories, double categories and \inf-categories. We
also remind the reader of the relationship between double categories
and 2-categories, by discussing the Spencer's recognition principle
for those double categories that arise as the double categories of
pasting squares in some 2-category (see Brown and
Mosa~\cite{Brown:1999:Connections}).

Section~\ref{decintro.sect}, the last of these contextual
introductions, provides a account of the classical simplicial decalage
construction and the method of simplicial reconstruction.  We review
those parts of the theory of (co)monads which were developed in order
to provide a general way of constructing functors into categories of
simplicial structures. In this context, it is a classical result that
the category of simplicial sets supports a canonical comonad, called
the decalage comonad, which is in some sense generic for this
construction.  Again, this material will be very familiar to Algebraic
Topologists and Category Theorists but may be less familiar to others
and, indeed, the 2-categorical presentation we give here may be
considered to be somewhat non-standard.

From section~\ref{strat.filt.sect} we concentrate on developing the
theory of complicial sets and for much of this work we study these as
novel structures in their own right. Only much later, in
section~\ref{street.nerve.sec}, do we ``tie the knot'' by relating our
constructions back to the traditional theory of (strict)
\inf-categories. In order to do so we contribute to the theory of parity
complexes and provide a deepened analysis of Street's nerve functor.

One of the attractions of complicial sets as a foundation for
\inf-category theory, of both the strong and weak variety, is that
they build upon the familiar theory of simplicial sets.  However it is
too much to hope that simplicial sets themselves are enough,
especially since Street's canonical nerve construction does not
provide us with a full representation of \inf-categories as simplicial
sets. The issue here is that this nerve does not record enough
information about the identities in our \inf-categories, a deficiency
we rectify by storing this missing data using a structure dubbed {\em
  hollowness\/} by Street but later renamed {\em stratification}. We
examine the theory of such {\em stratified simplicial sets\/} in
section~\ref{sect.strat} and later, in section~\ref{sect.filt.semi},
we discuss {\em filtered semi-simplicial sets\/} which may be used to
provide an alternative foundation to our work and which we apply at a
strategically important point in the sequel.

In section~\ref{precomp.sec} we introduce a class of stratified sets
which, for reasons that shall become apparent, we choose to call {\em
  pre-complicial sets}. Pre-complicial sets satisfy a relatively weak
constraint, which allows us to make inferences about the thinness of
simplices that are related by virtue of being faces on a suitable
higher dimensional thin simplex. Our primary result in this section is
that the full subcategory of pre-complicial sets admits a
(non-symmetric) biclosed monoidal structure, which will become our
main tool in much of what follows.

In section~\ref{comp.sec}, we finally define Roberts' {\em complicial
  sets\/} to be pre-complicial sets in which all {\em admissible
  horns\/} have unique thin fillers.  Most of the effort in this
section is devoted to extending the results of
section~\ref{precomp.sec} to the reflective full subcategory of
complicial sets. In particular we show that we may reflect our
biclosed monoidal structure on pre-complicial sets down onto the
category of complicial sets.

In section~\ref{comp.pathcat.sec}, we use the biclosed monoidal
structure on complicial sets to build internal and enriched categories
from complicial sets. Specifically, we present a construction which
allows us to derive a category of paths in a complicial set whose
homsets are themselves complicial sets. Later we show that, in a
suitable sense, we may iterate this construction to provide an
explicit equivalence inverse to Street's nerve functor.

In section~\ref{pathecat.ff.sec} we build analogues of the traditional
simplicial decalage construction on the categories of complicial sets
and complicially enriched categories. Using an argument which deploys a
certain complicial double category with connection derived in the
previous section, we show that these constructions correspond to each
other under the action of our path category functor. Applying
simplicial reconstruction to these then establishes the fact that the
path category construction provides us with a fully faithful
representation of complicial sets as complicially enriched categories.

Finally, in section~\ref{street.nerve.sec} we quickly review Street's
work on parity complexes \cite{Street:1991:Parity},
\cite{Street:1994:Parity} and \inf-categorical nerves
\cite{Street:1987:Oriental}, \cite{Street:1988:Fillers} and establish
its connection to our work on the theory of complicial sets.  In
particular, we construct a canonical isomorphism between the free
\inf-category on the parity complex product of a pair of simplices and
the \inf-category obtained by reflecting the tensor product of
corresponding standard simplices. This confirms our intuition that the
complicial tensor generalises the (lax) Gray tensor product of
2-categories~\cite{Gray:1974:FormalCat} and that, under Street's nerve
construction, it coincides with the \inf-categorical tensor products
of Crans~\cite{Crans:1995:PhD} and Steiner~\cite{Steiner:1991:Tensor}.

Once we have done this we can partially free ourselves from the
complicial world and discuss the all important relationship between
Street's nerve construction and our own path category construction.
Finally, we use this to provide a quick and easy proof that Street's
nerve provides an {\bf equivalence} between the categories of
\inf-categories and complicial sets, as originally conjectured in
\cite{Street:1987:Oriental}.

\nocite{Brown:1981:CC}



\section{Simplicial Operators and Simplicial Sets}
\label{simp.sect}

\subsection{Simplicial Operators}

\begin{defn}[The category $\aDelta$]
  Let $\aDelta$ denote the skeletal category of finite
  ordinals and order preserving functions. In other words, $\aDelta$
  has:
\begin{itemize}
\item {\bf Objects} ordered sets $[n] = \{ 0 < 1 < \cdots < n \}$ one
  for each $n\geq -1$,
\item {\bf Maps} $\arrow\alpha:[n]->[m].$ which are order
  preserving functions from the ordered set $[n]$ to the ordered set
  $[m]$,
\item {\bf Composition} simply that of functions, which is well
  defined since the composite of order preserving maps is again order
  preserving.
\end{itemize}
We use the notation $\id_{[n]}$ to denote the identity function on
$[n]$. The maps of $\aDelta$ are often referred to as {\em simplicial
  operators}, for reasons which will become clear.
\end{defn}

\begin{notation}[The Topologist's $\Delta$]
  The category $\aDelta$ is known as the {\em Algebraist's\/} $\Delta$
  since, as recalled later on, this category ``classifies'' the
  algebraic theory of {\em monoids}.  However, for much of this work
  we will be interested in studying simplicial operators from the
  topological perspective. To this end, we will often restrict our
  attention to $\Delta$ the full subcategory of $\aDelta$ whose
  objects are the {\bf non-zero} ordinals $[n]$ for $n\geq 0$.
  Following the usual tradition, we will usually refer to this
  category as the {\em Topologist's\/} $\Delta$.

  Under the topological interpretation of $\Delta$, the object $[n]$
  is considered to be a combinatorial rendering of the standard
  $n$-dimensional simplex. This explains our rather peculiar, but
  nonetheless entirely standard, use of the notation $[n]$ to denote
  the $(n+1)^{\text{th}}$ ordinal.
\end{notation}

\begin{notation}[faces and degeneracies] The following standard notation
  and nomenclature will be used throughout:
\begin{itemize}
\item The injective maps in $\aDelta$ are referred to as {\em face
    operators}.
\item For each $n\in\mathbb{N}$ and $j\in[n]$ define the simplicial
  operator $\arrow\face^n_j:[n-1]->[n].$ by
\begin{displaymath}
        \face^n_j(i) = 
        \left\{
        \begin{array}{cl}
                i & \text{if $i<j$,} \\
                i+1 & \text{otherwise.}
        \end{array}
        \right.
\end{displaymath}
This is called the {\em $j^{\text{\em th}}$ elementary face operator
  of $[n]$}.
\item The surjective maps in $\aDelta$ are referred to as {\em
    degeneracy operators}.
\item For each $n\in\mathbb{N}$ and $j\in[n]$ define the simplicial
  operator $\arrow\degen^n_j:[n+1]->[n].$ by
\begin{displaymath}
        \degen^n_j(i) = 
        \left\{
        \begin{array}{cl}
                i & \text{if $i\leq j$,} \\
                i-1 & \text{otherwise.}
        \end{array}
        \right.
\end{displaymath}
This is called the {\em $j^{\text{\em th}}$ elementary degeneracy
  operator of $[n]$}.
\item For each $n\in\mathbb{N}$ and $i\in[n]$ the operator
  $\arrow\vertex^n_i:[0]->[n].$ given by $\vertex^n_i(0)=i$ is
  called the {\em $i^{\text{th}}$ vertex operator of $[n]$}.
\item For each $n\in\mathbb{N}$ we use the notation $\eta^n$ to denote
  the unique operator from $[n]$ to $[0]$.
\end{itemize}
Unless doing so would introduce an ambiguity, we will tend to reduce
notational clutter by dropping the superscripts of these elementary
operators.
\end{notation}

\begin{obs}[the simplicial identities]\label{simplicial.idents} The
  following classical relationships hold in $\aDelta$ and are sufficient
  to fully characterise equalities between composites of elementary
  face and degeneracy operators in $\aDelta$:
  \begin{itemize}
  \item for any pair $j<i\in [n+1]$ we have
    $\face^{n+1}_i\circ\face^n_j=\face^{n+1}_j\circ\face^n_{i-1}$, and
  \item for any pair $j\leq i\in[n-1]$ we have
    $\degen^{n-1}_i\circ\degen^n_j=\degen^{n-1}_j\circ\degen^n_{i+1}$.
  \item for all $j\in[n]$ and $i\in[n-1]$ we have
    \begin{displaymath}
      \degen_i^{n-1}\circ\face_j^n = \left\{
        \begin{array}{lp{2in}}
          \face^{n-1}_j\circ\degen^{n-2}_{i-1} & if $j<i$, \\
          \id_{[n-1]} & if $j=i$ or $j=i+1$, \\
          \face^{n-1}_{j-1}\circ\degen^{n-2}_i & if $j>i+1$.
        \end{array}\right.
    \end{displaymath}
  \end{itemize}
\end{obs}

\begin{notation}[partition operators]\label{part.oper.def}
  We say that a pair $p,q\in\mathbb{N}$ is a {\em partition\/} of
  $n\in\mathbb{N}$ if $p+q=n$. For each such partition we have:
  \begin{itemize}
  \item face operators $\arrow\partinj^{p,q}_1:[p]->[n].$ given by
    $\partinj^{p,q}_1(i)=i$ and $\arrow\partinj^{p,q}_2:[q]->[n].$
    given by $\partinj^{p,q}_2(j)=j+p$, and
  \item degeneracy operators $\arrow\partproj^{p,q}_1:[n]->[p].$ given
    by 
    \begin{equation*}
      \partproj^{p,q}_1(i) =
      \begin{cases}
        i & \text{when $i\leq p$ and} \\
        p & \text{when $i>p$.}
      \end{cases}
    \end{equation*}
    and $\arrow\partproj^{p,q}_2:[n]->[q].$ given by:
    \begin{equation*}
      \partproj^{p,q}_2(i) =
      \begin{cases}
        0 & \text{when $i< p$ and} \\
        i - p & \text{when $i\geq p$.}
      \end{cases}
    \end{equation*}
  \end{itemize}

  We call these {\em partition operators\/} and, as is easily
  verified, they satisfy the following {\em partition identities}:
  \begin{equation}\label{part.ident}
  \begin{gathered}
    \partproj^{p,q}_1\circ\partinj^{p,q}_1 = \id_{[p]}\mkern80mu 
    \partproj^{p,q}_2\circ\partinj^{p,q}_2 = \id_{[q]} \\
    \begin{aligned}
      \partinj^{p+q,r}_1\circ\partinj^{p,q}_1 & {} =
      \partinj^{p,q+r}_1 \mkern30mu &
      \partinj^{p+q,r}_1\circ\partinj^{p,q}_2 & {} = 
      \partinj^{p,q+r}_2\circ\partinj^{q,r}_1 & \mkern30mu
      \partinj^{p,q+r}_2\circ\partinj^{q,r}_2 & {} =
      \partinj^{p+q,r}_2 \\
      \partproj^{p,q}_1\circ\partproj^{p+q,r}_1 & {} =
      \partproj^{p,q+r}_1  &
      \partproj^{p,q}_2\circ\partproj^{p+q,r}_1 & {} = 
      \partproj^{q,r}_1\circ\partproj^{p,q+r}_2 &
      \partproj^{q,r}_2\circ\partproj^{p,q+r}_2 & {} =
      \partproj^{p+q,r}_2  
    \end{aligned}
  \end{gathered}
\end{equation}
\end{notation}

\begin{obs}[duals of simplicial operators]\label{simp.op.dual} There
  exists a canonical functor $(\mathord{-})^{\mathord{\circ}}$
  from $\aDelta$ to itself which ``maps each ordinal to its dual as an
  ordered set''. Explicitly, $[n]^{\mathord{\circ}}=[n]$ and if
  $\arrow\alpha:[n]->[m].$ is a simplicial operator then, for each
  $i\in[n]$, $\alpha^{\mathord{\circ}}(i)=m-\alpha(n-i)$. Clearly this
  dual functor is {\em strictly involutive\/} in the sense that the
  diagram
  \begin{displaymath}
    \xymatrix@=2.5em{
      {\aDelta}\ar[rr]^{\textstyle\text{id}_{\aDelta}}
      \ar[dr]_{\textstyle(\mathord{-})^\circ} & 
      \save[]-<0em,1.25em>*{\textstyle ||}\restore &
      {\aDelta} \\
      & {\aDelta}\ar[ur]_{\textstyle(\mathord{-})^\circ} & }
  \end{displaymath}
  of functors commutes (on the nose). In other words, the functor
  obtained by composing $(\mathord{-})^\circ$ with itself is the
  identity on $\aDelta$. Notice also that the identities 
  \begin{equation}
    \begin{gathered}
      (\face^n_i)^\circ=\face^n_{n-i} \mkern40mu
      (\degen^n_i)^\circ=\degen^n_{n-i} \mkern40mu 
      (\vertex^n_i)^\circ = \vertex^n_{n-i}\\
      (\partinj^{p,q}_1)^\circ=\partinj^{q,p}_2 \mkern30mu
      (\partinj^{p,q}_2)^\circ=\partinj^{q,p}_1 \mkern40mu
      (\partproj^{p,q}_1)^\circ=\partproj^{q,p}_2 \mkern30mu
      (\partproj^{p,q}_2)^\circ=\partproj^{q,p}_1
    \end{gathered}\label{elem.dual.disp}
  \end{equation}
  hold between elementary operators and their duals.
\end{obs}

\begin{obs}[a useful characterisation of face and degeneracy
  operators]\label{fd.char} The following facts about simplicial
  operators $\arrow\alpha:[n]->[m].$ are sometimes of use:
\begin{enumerate}[(i)]
\item\label{degen.char} $\alpha$ is a degeneracy operator iff
  $\alpha(0)=0$, $\alpha(n)=m$ and for all $j\leq j'\in[n]$ we have
  $\alpha(j')-\alpha(j)\leq j'-j$.
\item\label{face.char} $\alpha$ is a face operator iff for all
  $j\leq j'\in[n]$ we have $\alpha(j')-\alpha(j)\geq j'-j$.
\item\label{face.char.2} $\alpha$ is a face operator iff there exists
  a simplicial operator $\arrow\tilde\alpha:[n]->[m-n].$ such that
  $\alpha(j)=j + \tilde\alpha(j)$ for all $j\in[n]$.
\end{enumerate}
\end{obs}

\begin{obs}[face-degeneracy factorisation]\label{face.degen.fac} Every
  simplicial operator $\arrow\alpha:[n]->[m].$ factors {\bf uniquely}
  into a composite $\alpha^f\circ\alpha^d$ where $\alpha^f$ is a face
  operator and $\alpha^d$ is a degeneracy. Furthermore, let
  $\im(\alpha)$ denote the subset of $[m]$ given by:
  \begin{displaymath}
    \im(\alpha) = \left\{ i\in[m] \mathrel{|} (\exists j\in[n])
      \alpha(j) = i \right\}
  \end{displaymath}
  It is a trivial, but nonetheless useful, fact that one simplicial
  operator $\arrow\alpha:[m]->[n].$ factors though another
  $\arrow\beta:[r]->[n].$, that is to say there is some
  $\arrow\gamma:[m]->[r].$ with $\beta\circ\gamma=\alpha$, if
  and only if $\im(\alpha)\subseteq \im(\beta)$.
\end{obs}      

\subsection{The Algebraist's $\Delta$ and 2-Categories}\label{subsect.2cat}
We quickly review the theory of 2-categories, which will be useful for
expressing some of the ``meta-theory'' developed in the remainder of
this section. Later on, in section~\ref{sect.ncats}, we will take a
second bite at the 2-category cherry and consider them ``in the
small'' when we review the algebraic theory of \inf-categories.

\begin{obs}[the cartesian closed category $\Category$]
  We will let $\Category$ denote the (huge) category of all (large)
  categories and functors between them. 
  
  This category is cartesian closed, where the ``function space'' from
  a category $\mathcal{C}$ to a category $\mathcal{D}$ is the {\em
    functor category\/} $\funcat[\mathcal{C},\mathcal{D}]$ which has:
  \begin{itemize}
  \item {\bf Objects} functors from $\mathcal{C}$ to $\mathcal{D}$,
  \item {\bf Maps} natural transformations between such functors,
  \item {\bf Composition} the usual, point-wise, composition of natural
    transformations.
  \end{itemize}

  By definition, this function space is characterised by the 
  adjunction 
  \begin{displaymath}
    \let\labelstyle=\textstyle
    \xymatrix@R=1ex@C=14em{
      {\Category}\ar[r]^{\bot}_{\funcat[\mathcal{C},{*}]} &
      {\Category}\ar@/_3ex/[l]_{{-}\times\mathcal{C}}}
  \end{displaymath}
  where $\times$ denotes the cartesian product of categories. In other
  words, there exists a natural bijection between functors $\arrow
  \ffunc:\mathcal{B}\times\mathcal{C}->\mathcal{D}.$ and
  $\arrow\hat{\ffunc}:\mathcal{B}->\funcat[\mathcal{C},\mathcal{D}].$.
  It follows (see \cite{Kelly:1982:ECT}) that $\Category$ is a rich
  ``universe'' over which we might enrich the homsets of other
  categories. This observation leads to the following definition:
\end{obs}

\begin{defn}[2-category]\label{2cat.defn.1} A {\em 2-category\/} is a category
  enriched in the cartesian closed category $\Category$.
  Correspondingly, a {\em 2-functor\/} is simply an
  $\Category$-enriched functor between 2-categories.
\end{defn}

The canonical reference for the theory of enriched categories is, of
course, Kelly's comprehensive book \cite{Kelly:1982:ECT} - to which we
recommend the reader.

\def\tC{\tcat{C}} 
  
\begin{obs}[2-categories explicitly]\label{2cat.expl}
  Fundamentally, a 2-category is simply a category $\tC$ in which each
  homset $\tC(C,D)$ is itself a category (object of $\Category$). This
  immediately implies that our 2-category contains 3 sorts of objects:
  \begin{itemize}
  \item {\bf 0-cells} which are the objects of $\tC$ (for which we use
    symbols $A,B,C,...$),
  \item {\bf 1-cells} which are the objects of the homsets $\tC(A,B)$
    (for which we use symbols $F,G,H,...$),
  \item {\bf 2-cells} which are the arrows of the homsets $\tC(A,B)$
    (for which we use Greek symbols $\lambda,\mu,\nu,...$).
  \end{itemize}
  
  Commensurate with the two layers of ``category-ness'' involved in
  $\tC$ we get two distinct category structures on these cells:
  \begin{itemize}
  \item {\bf horizontal} with objects which are 0-cells, arrows which
    are 1-cells {\bf and\/} 2-cells and compositional structure given by the
    quadruple $\quadruple<\circ;\dom_0;\cod_0;\id_0>$. Often we write
    this composition simply as juxtaposition.
  \item {\bf vertical} with objects which are 1-cells, arrows which
    are 2-cells and compositional structure given by
    $\quadruple<\cdot;\dom_1;\cod_1;\id_1>$.
  \end{itemize}
  In diagrams and running text we will tend to use single arrows
  $\arrow F:A->B.$ to denote 1-cells (with $A=\dom_0(F)$ and
  $B=\cod_0(F)$) and double arrows
  $\twocell:\lambda:F->G.$ to denote 2-cells (with $F=\dom_1(\lambda)$
  and $G=\cod_1(\lambda)$).

  These category structures must satisfy a number of compatibility
  conditions which bind them together intimately. Most important
  amongst these are:
  \begin{itemize}
  \item {\bf globularity:} for all 2-cells $\lambda\in\tC$ we have
    $\dom_0(\dom_1(\lambda))=\dom_0(\cod_1(\lambda))$ and
    $\cod_0(\dom_1(\lambda))=\cod_0(\cod_1(\lambda))$.
  \item {\bf middle four interchange:} if we have 2-cells
    $\lambda, \lambda', \mu$ and $\mu'$ then the equality
    \begin{math}
      (\mu'\circ\lambda')\cdot(\mu\circ\lambda) = 
      (\mu'\cdot\mu)\circ(\lambda'\cdot\lambda)
    \end{math}
    holds whenever the various composites involved are defined.
    Notice that the composite on the left of this equality is defined
    iff that on its right is defined.
  \end{itemize}

  The globularity condition on a 2-cell $\lambda$ implies that we may
  picture it as follows:
  \begin{displaymath}
    \xymatrix@R=4ex@C=8em{
      {A}\ar@/^3ex/[rr]^{\textstyle F}_{}="one"
      \ar@/_3ex/[rr]_{\textstyle G}^{}="two" & 
      \ar@{=>}^{\textstyle\lambda} "one"; "two" & {B} }
  \end{displaymath}
  where $F=\dom_1(\lambda)$, $G=\cod_1(\lambda)$, $A=\dom_0(\lambda)$ and
  $B=\cod_0(\lambda)$.
  
  Finally, a 2-functor $\arrow\mathcal{F}:\tcat{C}->\tcat{D}.$ may be
  thought of explicitly as a map which maps the cells of $\tcat{C}$ to
  those of $\tcat{D}$ in a way which acts functorially on both
  horizontal {\em and\/} vertical category structures
\end{obs}

\begin{obs}[whiskering]\label{whiskering}  
  We commonly identify 0-cells with their corresponding identity
  1-cells (under $\id_0$), and 1-cells in turn with their identity
  2-cells (under $\id_1$). For instance, we will often use this
  convention, and the one that allows us to replace $\circ$ by
  juxtaposition, to write things like
  \begin{itemize}
  \item $G\lambda$ for $(\id_1(G))\circ\lambda$ and 
  \item $\mu F$ for $\mu\circ(\id_1(F))$.
  \end{itemize}
  These two derived operations are so useful that they have been
  respectively dubbed left and right {\em whiskering\/} by
  Street~\cite{Street:1996:CatStruct}. Indeed it is possible (and quite
  informative) to re-cast the theory of 2-categories purely in terms
  of these whiskering operations and vertical composition alone.  That
  this is the case follows directly from the observation that the
  identity and middle four interchange rules imply that for
  horizontally composable 2-cells $\twocell\lambda:F->F'.$ and
  $\twocell\mu:G->G'.$ we have
  \begin{displaymath}
    (G'\lambda)\cdot(\mu F) = \mu\circ\lambda = 
    (\mu F')\cdot(G\lambda)
  \end{displaymath}
\end{obs}

\begin{obs}[$\Category$ as a 2-category]\label{cat.as.2cat}
  We know, from Kelly~\cite{Kelly:1982:ECT}, that me may immediately
  enrich the cartesian closed category $\Category$ over itself, by
  making $\funcat[\mathcal{B},\mathcal{C}]$ the enriched homset
  between categories $\mathcal{B}$ and $\mathcal{C}$.  This 2-category
  structure is often best described explicitly in terms of whiskering
  operations, it has:
  \begin{itemize}
  \item {\bf 0-cells} (large) categories, {\bf 1-cells} functors
    between these and {\bf 2-cells} natural transformations between
    those.
  \item {\bf vertical} composition of 2-cells given by the usual
    ``point-wise'' composite of natural transformations.
  \item {\bf left whiskering} of a natural transformation
    $\twocell\lambda:\ffunc->\ffunc'.\arrow:\mathcal{B}->\mathcal{C}.$ by a
    functor $\arrow \gfunc:\mathcal{C}->\mathcal{D}.$ is formed by applying
    $\gfunc$ ``point-wise'' to the components of $\lambda$, that is
    $(\gfunc\lambda)_b\defeq \gfunc(\lambda_b)$ for each $b\in\obj(\mathcal{B})$.
  \item {\bf right whiskering} of a natural transformation
    $\twocell\mu:\gfunc->\gfunc'.\arrow:\mathcal{C}->\mathcal{D}.$ by
    a functor $\arrow \ffunc:\mathcal{B}->\mathcal{C}.$ is obtained by
    re-indexing the components of $\mu$ using the action of the
    functor $\ffunc$ on objects, that is $(\mu \ffunc)_b\defeq
    \mu_{\ffunc(b)}$ for each $b\in\obj(\mathcal{B})$.
  \end{itemize}
\end{obs}

\begin{obs}[$\aDelta$ as a 2-category]\label{2cat.delta}
  Since each object of $\aDelta$ is an ordered set it follows that each
  of its homsets $\aDelta([n],[m])$ possesses a natural partial order
  $\leq$, given by
  \begin{displaymath}
    \alpha\leq\beta
    \mathrel{\text{ if and only if }}
    (\forall j\in[n])(\alpha(j)\leq\beta(j))
  \end{displaymath}
  Furthermore, since each map in $\aDelta$ is order preserving, we also
  know that composition preserves these partial orders, in the sense
  that if $\alpha\leq\beta$ in $\aDelta([n],[m])$ and $\alpha'\leq\beta'$
  in $\aDelta([m],[r])$ then $(\alpha'\circ\alpha) \leq (\beta'\circ\beta)$ 
  in $\aDelta([n],[r])$.
  
  It follows that, under these partial orders on its homsets, $\aDelta$
  becomes a partial order enriched category or, in other words, a
  2-category each homset of which is a partial order.
  
  In fact, this is simply the full sub-2-category of $\Category$ on
  those 0-cells obtained by considering each totally ordered set $[n]$
  ($n\in\mathbb{N}$) as a category in the usual way. That is to say,
  think of $[n]$ as a category with objects $\{0,1,...,n\}$ and a unique
  arrow $\overarr:i->j.$ for each pair of integers with $i\leq j$.
  
  It is worth noting that the following useful inequalities (2-cells)
  hold between (composites of) face and degeneracy operators:
  \begin{itemize}
  \item for all $n\geq 1$ and $j\in[n-1]$ we have
    \begin{math}
      \face_j^n\circ\degen_j^{n-1} > \id_{[n]}
    \end{math} and
    \begin{math}
      \face_{j+1}^n\circ\degen_j^{n-1} < \id_{[n]}
    \end{math},
  \item if $n\geq 1$ and $j\leq k\in[n]$ then
    \begin{math}
      \face^n_k \leq \face^n_j
    \end{math} and
    \begin{math}
      \degen^{n-1}_j \leq \degen^{n-1}_k 
    \end{math}.
  \end{itemize}
\end{obs}

\begin{obs}[adjoints in the 2-category $\aDelta$]\label{simp.op.adj} In
  the sequel we will have occasion to consider adjoint pairs
  $\alpha\adjoint\beta$ of simplicial operators. These can be defined
  and understood abstractly in the 2-category $\aDelta$, in terms of
  the general theory of adjoint pairs in 2-categories as expounded
  (for instance) by Kelly and Street in their classic review
  article~\cite{KellyStreet:1974:2Cats} or by Street
  in~\cite{Street:1980:FibBicat}.
  
  However, here we prefer to directly exploit the identification of
  the 2-category $\aDelta$ as a full sub-2-category of $\Category$ and
  utilise the traditional definition of adjunction.  In such terms, a
  simplicial operator $\arrow\alpha:[n]->[m].$ is left adjoint to
  $\arrow\beta:[m]->[n].$ if for all $i\in[n]$ and $j\in[m]$ we have
  $\alpha(i)\leq j$ iff $i\leq\beta(j)$. Equivalently, this condition
  holds if the inequalities $\id_{[n]}\leq\beta\circ\alpha$
  and $\alpha\circ\beta\leq\id_{[m]}$ hold.
  
  Of course, the usual properties of adjunctions hold for those in
  $\aDelta$, in particular the left (resp. right) adjoint to an
  operator $\alpha$ (if such a thing exists) is unique.  Furthermore,
  the following classical observations characterise adjunctions of
  simplicial operators:
  \begin{enumerate}[(i)]
  \item\label{so.adj.left} $\alpha$ has a (necessarily unique) {\em
      left\/} adjoint $\arrow\alpha^l:[m]->[n].$ iff
    $\alpha(n)=m$, in which case $\alpha^l(i)=\min\{j\in[n]\mid
    i\leq\alpha(j)\}$.
  \item\label{so.adj.right} Dually, $\alpha$ has a (necessarily
    unique) {\em right\/} adjoint $\arrow\alpha^r:[m]->[n].$
    iff $\alpha(0)=0$, in which case $\alpha^r(i)=\max\{j\in[n]\mid
    i\geq\alpha(j)\}$.
  \item It follows that if $\alpha$ is a degeneracy operator
    then it has both left and right adjoints.
  \item\label{so.adj.rinv} If $\alpha\adjoint\beta$ is an adjoint pair
    and either $\alpha$ is a degeneracy operator or $\beta$ is a face
    operator then $\alpha\circ\beta = \id_{[m]}$.
  \item\label{so.adj.linv} Dually, if $\alpha\adjoint\beta$ is an
    adjoint pair and either $\alpha$ is a face operator or $\beta$ is
    a degeneracy operator then $\beta\circ\alpha = \id_{[n]}$.
  \item\label{so.adj.elem} The simplicial identities and inequalities
    of observation~\ref{simplicial.idents} demonstrate that
    $\degen_j^{n-1}$ has {\bf right adjoint} $\face_j^n$ and {\bf left
      adjoint} $\face_{j+1}^n$.
  \end{enumerate}
  Finally, notice that the dual functor $(\mathord{-})^\circ$ may be naturally
  extended to a 2-functor which is {\bf contravariant on 2-cells}. By
  this we mean to say that if $\alpha$ and $\beta$ are simplicial
  operators with $\alpha\leq\beta$ then $\beta^{\mathord{\circ}}\leq
  \alpha^{\mathord{\circ}}$. It follows that $(\mathord{-})^\circ$
  carries left (resp. right) adjoints to right (resp. left) adjoints.
  In other words, $(\alpha^\circ)^r = (\alpha^l)^\circ$ and
  $(\alpha^\circ)^l = (\alpha^r)^\circ$ when these adjoints exist.
\end{obs}

\subsection{The Algebraist's $\Delta$ and Monoidal Categories}

In this subsection, we recall a few standard results regarding the
algebraic content of $\aDelta$. For proofs of these results, and the
details of many other interesting facets of the algebra of $\aDelta$,
we refer the reader to Mac~Lane~\cite{Maclane:1971:CWM}.

\begin{obs}[$\aDelta$ by generators and relations]\label{delta.gr} We
  can go further than we did in the last subsection and factor every 
  face (resp.\ degeneracy) operator as
  a (canonical) composite of elementary face (resp.\ degeneracy)
  operators. It follows that the elementary operators generate
  $\aDelta$, but more is true. 
  
  We said above that the simplicial identities served to fully
  characterise the compositional relationships between elementary
  operators. This can be made precise by observing that $\aDelta$ may
  be {\em presented\/} in terms of the elementary operators as its
  generators and the simplicial identities as its relations.
  
  In other words, in order to define a functor $\arrow
  \ladj:\aDelta->\mathcal{E}.$ it suffices to specify its action on objects
  and elementary operators and then check that the images of these
  elementary operators satisfy the simplicial identities in
  $\mathcal{E}$. It also implies that we may test the naturality of a
  family of maps between two functors $\arrow
  \ladj,\ladj':\aDelta->\mathcal{E}.$ by testing it with respect to the
  elementary operators alone.
\end{obs}

\begin{recall}[monoidal categories] Recall (from say Mac~Lane
  \cite{Maclane:1971:CWM}) that a {\em monoidal category\/}
  $\triple<\mathcal{E};\otimes;I>$ consists of a category
  $\mathcal{E}$ equipped with a bifunctor $\arrow\otimes:
  \mathcal{E}\times\mathcal{E}->\mathcal{E}.$ called a {\em tensor
    product}, an object $I\in\mathcal{E}$, and natural isomorphisms:
  \begin{displaymath}
    \begin{array}{rclp{1in}}
      X\otimes(Y\otimes Z) & \stackrel{\alpha_{X,Y,Z}}\cong &
      (X\otimes Y)\otimes Z & associativity \\
      X\otimes I & \stackrel{\iota^r_X}\cong & X & right identity \\
      I\otimes X & \stackrel{\iota^l_X}\cong & X & left identity
    \end{array}
  \end{displaymath}
  These isomorphisms must also satisfy a number of {\em coherence
    conditions}, two of which relate the identity isomorphisms to the
  associativity and one, called the {\em pentagon condition}, which
  relates the various associativities involved in a four fold tensor
  product (see Mac~Lane \cite{Maclane:1971:CWM} for greater detail).
  
  A {\em monoidal functor\/} $\arrow \ladj:{\triple<
    \mathcal{E};\otimes;I>}->{\triple<\mathcal{E}';
    \otimes';I'>}.$ consists of a functor $\arrow \ladj:
  \mathcal{E}->\mathcal{E}'.$ of underlying categories and a natural
  family of morphisms:
  \begin{displaymath}
    \begin{array}{rclp{2in}}
      \ladj(X)\otimes' \ladj(Y) & \stackrel{\mu_{X,Y}}\rightarrow & \ladj(X\otimes Y) & 
      product comparison \\ 
      I' & \stackrel{\nu}\rightarrow & \ladj(I) & identity comparison 
    \end{array}
  \end{displaymath}
  These morphisms must also satisfy coherence conditions, the first
  relating the identity comparison to right and left identity
  isomorphisms and the second relating the product comparisons to
  associativity isomorphisms (again see Mac~Lane \cite{Maclane:1971:CWM} for
  greater detail). We say that a monoidal functor is {\em strong\/}
  (resp.\ {\em strict}) if each of these comparison maps is in fact an
  isomorphism (resp.\ an identity) in $\mathcal{E}'$.

  A {\em monoidal natural transformation\/} 
  $\twocell \eta:{\ladj}->{\ladj'}.$ between monoidal functors $\ladj$ and $\ladj'$
  is simply a natural transformation the components of which commute with the
  product and identity comparisons on its domain and codomain. In
  other words, we must have $\eta_I\cdot\nu=\nu'$ and $\eta_{X\otimes
    Y}\cdot\mu_{X,Y}=\mu'_{X,Y}\cdot(\eta_X\otimes\eta_Y)$. 
  
  Notice that, we may extend the usual composites of functors and
  natural transformations to provide corresponding composites of
  monoidal functors and monoidal natural transformations. To do this
  all we need worry about is how to define the comparison maps for the
  composite of a pair of monoidal functors $\arrow \ladj:
  {\triple<\mathcal{E};\otimes;I>}->{\triple<\mathcal{E}';\otimes';I'>}.$
  and $\arrow \ladj': {\triple<\mathcal{E}';\otimes';I'>}->
  {\triple<\mathcal{E}'';\otimes'';I''>}.$. However, a moment's reflection
  reveals that the most natural candidates for these are the composites
  \begin{displaymath}
    \xymatrix@R=2ex@C=5em{
      {\ladj'(\ladj(X))\otimes'' \ladj'(\ladj(Y))}
      \ar[r]^{\mu'_{\ladj(X),\ladj(Y)}} &
      {\ladj'(\ladj(X)\otimes' \ladj(Y))} \ar[r]^{\ladj'(\mu_{X,Y})} &
      {\ladj'(\ladj(X\otimes Y))} }
  \end{displaymath}
  for tensor products and
  \begin{displaymath}
    \xymatrix@R=2ex@C=5em{
      {I''} \ar[r]^{\nu'} &
      {\ladj'(I')} \ar[r]^{\ladj'(\nu)} &
      {\ladj'(\ladj(I))} }
  \end{displaymath}
  for identities.
  Verifying the various coherence conditions for
  this structure, as well as checking that it is well behaved with
  respect to horizontal composition of monoidal natural
  transformations, is a matter of routine calculation. In this way we
  lift the structure of the 2-category $\Category$ (of categories,
  functors and natural transformations) to get a 2-category
  $\Monoidal$ with 0-cells monoidal categories, 1-cells monoidal
  functors between them and 2-cells that are monoidal natural
  transformations between those.
\end{recall}
 
\begin{recall}[monoids]\label{monoid.stuff}
  As the name implies, monoidal categories bear a strong relationship
  to {\em monoids}, that is algebraic structures $\triple<X;*;e>$
  where $*$ is an associative binary operation on $X$ for which
  $e\in X$ is a two sided identity. This relationship finds its
  clearest expression in two interesting observations:\vspace{2ex}
  
  \noindent{\bf External:} Monoidal categories are a convenient
  structure within which to interpret the monoid concept. If
  $\triple<\mathcal{E}; \otimes;I>$ is a monoidal category then a
  monoid $\triple<X;*;e>$ in $\mathcal{E}$ consists of an underlying
  object $X\in\mathcal{E}$, a {\em multiplication\/} map $\arrow
  *:X\otimes X->X.$ and a {\em unit\/} map $\arrow e:I->X.$. This data
  must satisfy the diagrammatic conditions
  \begin{displaymath}
    \def\labelstyle{\textstyle}
    \xymatrix@R=4ex@C=5em{
      {X\otimes(X\otimes X)}\ar[r]_{\scriptstyle\cong}^{\alpha_{X,X,X}}
      \ar[d]_{X\otimes *} &
      {(X\otimes X)\otimes X}\ar[r]^{*\otimes X} &
      {X\otimes X} \ar[d]^{*} \\
      {X\otimes X}\ar[rr]_{*} & & {X} \\ 
      {X\otimes I} \ar[rd]^{\scriptstyle\cong}_{\iota^r_X}
      \ar[r]^{X\otimes e} &
      {X\otimes X} \ar[d]_{*} &
      {I\otimes X} \ar[ld]_{\scriptstyle\cong}^{\iota^l_X}
      \ar[l]_{e\otimes X} \\
      & {X} & }
  \end{displaymath}
  which are simply the usual associativity and identity conditions for
  a monoid in disguise.
  
  Many common algebraic structures turn out to be monoids in a
  suitable monoidal category. For instance, rings are really no more
  than monoids within the category of abelian groups equipped with the
  usual tensor product of groups as its monoidal structure.
    
  We form a category $\Mon\triple<\mathcal{E}; \otimes;I>$ of
  monoids in $\triple<\mathcal{E}; \otimes;I>$, in which a {\em
    monoid morphism\/} $f$ from $\triple<X;*;e>$ to
  $\triple<X';*';e'>$ is simply a map of the underlying objects
  $\arrow f:X->X'.$ in $\mathcal{E}$ which commutes with the monoid
  structures on its domain and codomain, or in other words for which
  $f\circ e = e'$ and $f\circ * = *'\circ(f\otimes f)$.  \vspace{2ex}
  
  \noindent{\bf Internal:} Monoidal categories in which each of the
  structural isomorphisms $\alpha_{X,Y,Z}$, $\iota^r_X$ and
  $\iota^l_X$ are actually {\em identities} are said to be {\em
    strict}. In other words, a strict monoidal category
  $\triple<\mathcal{C};\otimes;I>$ consists of a category equipped
  with a tensor structure for which the usual associativity and
  identity conditions $\lambda\otimes(\mu\otimes\nu) =
  (\lambda\otimes\mu)\otimes\nu$, $\lambda\otimes\id_I = \lambda$ and
  $\id_I\otimes\lambda = \lambda$ hold ``on the nose'' for all arrows
  $\lambda,\mu,\nu\in\mathcal{C}$.  Equivalently, such a structure is
  no more nor less than a monoid in the (huge) monoidal category
  $\triple<\Category;\times;1>$.
  
  Notice that, if $\tcat{C}$ is a 2-category then each endo-category
  $\tcat{C}(A,A)$ on a 0-cell $A\in\tcat{C}$ is the underlying
  category of a strict monoidal category, with tensor given by
  horizontal composition $\circ$ and identity $\id_0(A)$.  Conversely,
  each strict monoidal category $\triple<\mathcal{E};\otimes;I>$ gives
  rise to a 2-category $\Omega\mathcal{E}$, called its {\em
    (one-point) suspension}, with a single 0-cell $\ast$, hom-category
  $\Omega\mathcal{E}(\ast,\ast) = \mathcal{E}$, identity
  $\id_0(\ast)=I$, and horizontal composition $\otimes$.
\end{recall}

\begin{obs}[relating monoids and monoidal functors]
  It is worth observing that monoids and monoid maps bear a close
  relationship to monoidal functors and monoidal natural
  transformation. To be precise, it is clear that the one object, one
  arrow ``terminal'' category $\mathcal{I}$ admits a unique strict
  monoidal structure, with respect to which we may consider monoidal
  functors $\arrow \gfunc:\mathcal{I}->{\triple<\mathcal{E}; \otimes;I>}.$.  
  On examining the data for such monoidal functor, it becomes clear
  that it amounts to no more nor less than an object $\gfunc(1)$, an
  identity comparison $\arrow\nu:I->\gfunc(1).$ and a product comparison
  $\arrow\mu_{1,1}:\gfunc(1)\otimes \gfunc(1) -> \gfunc(1).$ satisfying coherence
  conditions that are identical to the axioms required of a monoid
  $\triple<\gfunc(1);\mu_{1,1};\nu>$ in $\mathcal{E}$.  Furthermore, the
  data and conditions for a monoidal natural transformation
  $\twocell\eta: \gfunc -> \gfunc'.$ correspond exactly
  to those for a monoid map between the monoids corresponding to
  $\gfunc$ and $\gfunc'$. In other words, the category
  $\Mon\triple<\mathcal{E}; \otimes;I>$ and the hom-category
  $\Monoidal(\mathcal{I}, \triple<\mathcal{E}; \otimes;I>)$ are
  canonically equivalent.
    
  This observation has many useful consequences, for instance it
  implies that a monoidal functor $\arrow \ladj:
  {\triple<\mathcal{E};\otimes;I>}->{\triple<\mathcal{E}';
    \otimes';I'>}.$ may be lifted to a functor between categories of
  monoids $\arrow\Mon(\ladj): {\Mon\triple<\mathcal{E};\otimes;I>
  }->{\Mon\triple< \mathcal{E}';\otimes'; I'>}.$.  To be precise, this
  is simply the post composition functor
  $\Monoidal(\mathcal{I},\ladj)$ between corresponding hom-categories
  of the 2-category $\Monoidal$.
\end{obs}

\begin{obs}[$\aDelta$ as a strict monoidal category]
  In order to see why $\aDelta$ might be of use to Algebraists, it is
  first necessary to observe that ordinal addition of its objects
  extends to a functor $\arrow\dsum:\aDelta\times\aDelta->\aDelta.$
  where $[n]\dsum[m]\defeq [n+m+1]$ and the {\em direct sum\/}
  $\beta\dsum\alpha$ of operators $\arrow \alpha:[n]->[n'].$ and
  $\arrow\beta:[m]->[m'].$ is defined by
  \begin{displaymath}
    \beta\dsum\alpha(j) = \left\{
      \begin{array}{lp{10em}}
        \alpha(j) & if $j\leq n$, \\
        n' + 1 + \beta(j - n - 1) & if $j>n$.
      \end{array}\right.
  \end{displaymath}
  It is easily seen that the triple $\triple<\aDelta;\dsum;[-1]>$
  satisfies the strict associativity and identity conditions discussed
  in recollection~\ref{monoid.stuff} and thus it is a strict monoidal
  category. Notice also that direct sums interact gracefully with the involution
  $\arrow(\mathord{-})^{\mathord{\circ}}:\Delta ->\Delta.$, in
  particular if $\alpha$ and $\beta$ are simplicial operators then we
  have $(\alpha\dsum\beta)^\circ = \beta^\circ\dsum\alpha^\circ$.
\end{obs}

\begin{obs}[generators for $\aDelta$ as a strict
  monoidal category]\label{del.gen.rel} The strict monoidal category
  $\triple<\aDelta;\dsum;[-1]>$ contains a canonical monoid with
  underlying object $[0]$, identity $\arrow\face_0^0:[-1]->[0].$ and
  multiplication $\arrow\degen^0_0:[0]\dsum[0]=[1]->[0].$, for which
  the associativity and unit conditions may be trivially verified.
  
  Notice that all elementary face and degeneracy operators may be
  obtained as iterated direct sums of the structural components of the
  monoid $\triple<[0];\degen^0_0;\face_0^0>$ as follows:
  \begin{itemize}
  \item each object $[n]$ may be obtained from the underlying object
    $[0]$ as an iterated direct sum $[0]^{\dsum (n+1)}$,
  \item each elementary face operator $\face^n_i$ may be obtained from
    the unit map $\face_0^0$ as an iterated direct sum $[0]^{\dsum
      (n-i)}\dsum\face_0^0\dsum [0]^{\dsum i}$, and
  \item each elementary degeneracy operator $\degen^n_i$ may be
    obtained from the multiplication map $\degen^0_0$ as an iterated
    direct sum $[0]^{\dsum (n - i)}\dsum \degen^0_0\dsum [0]^{\dsum
      i}$.
  \end{itemize}
  Of course we know, by observation~\ref{delta.gr}, that every
  simplicial operator may be expressed as a composite of elementary
  face and degeneracy operators and it follows that every operator may
  be obtained as a composite of iterated direct sums of the structural
  components of $\triple<[0];\degen^0_0;\face_0^0>$
  
  In fact much more is true for this monoid, as demonstrated by the
  following classical result which states that it actually {\em
    freely\/} generates $\triple<\aDelta;\dsum;[-1]>$ in some suitable
  sense made precise in the statement of the lemma.
\end{obs}

\begin{lemma}[the universal property of $\aDelta$]\label{univ.delta}
  The monoid $\triple<[0];\degen^0_0;\face_0^0>$ in
  $\triple<\aDelta;\dsum;[-1]>$ is universal, in the sense that for
  any other monoid $\triple<X;*;e>$ in a strict monoidal category
  $\triple<\mathcal{E};\otimes;I>$ there exists a unique strict
  monoidal functor $\arrow \tilde{X}:
  {\triple<\aDelta;\dsum;[-1]>}->{\triple<\mathcal{E};\otimes;I>}.$
  such that the monoid
  $$\Mon(\tilde{X})\triple<[0];\degen^0_0;\face_0^0>=
  \triple<\tilde{X}([0]);\tilde{X}(\degen^0_0);\tilde{X}(\face_0^0)>$$
  is equal to $\triple<X;*;e>$.
\end{lemma}

\begin{proof} Standard, see Mac~Lane \cite{Maclane:1971:CWM}. It is worth
  observing, however, that observation~\ref{del.gen.rel} immediately
  implies that the functor $\tilde{X}$ is completely determined by
  $\tilde{X}([n])= X^{\otimes(n+1)}$, $\tilde{X}(\face^n_i) =
  X^{\otimes(n-i)}\otimes e \otimes X^{\otimes i}$ and
  $\tilde{X}(\degen^n_i)  = X^{\otimes(n-i)}\otimes * \otimes
  X^{\otimes i}$.
\end{proof}

\subsection{Simplicial Sets}

\begin{defn}[simplicial sets and simplicial maps] The category
  $\Simp$ of {\em simplicial sets\/} and {\em simplicial
    maps\/} between them is simply the functor category
  $\funcat[\op\Delta,\Set]$, where $\Set$ denotes the (large) category
  of all (small) sets and functions between them.
\end{defn}

\begin{obs}[simplicial sets as partial right actions]
  \label{par.right} In practise, it is easier to think of a
  simplicial set as a single set endowed with a partially defined
  right action of the simplicial operators.
  
  To be more precise, this description presents a simplicial set as a
  triple $\triple<X;\mathord{\dim};\cdot>$ where:
\begin{itemize}
\item $X$ is a (small) set and $\arrow\dim : X->\mathbb{N}.$
  is a function, which assigns to each element of $X$ a {\em
    dimension},
\item $x\cdot\alpha\in X$ is defined for any element $x\in X$ and
  simplicial operator $\arrow\alpha:[m]->[n].$ for which $\dim(x) =
  n$, in which case we say that $x$ and $\alpha$ are {\em compatible},
  and we have $\dim(x\cdot\alpha)=m$, and
\item this action satisfies the equations $x\cdot \id_{[n]} = x$ and
  $(x\cdot\alpha)\cdot\beta = x\cdot(\alpha\circ\beta)$ whenever they
  are well defined (i.e.\ when $\arrow\alpha:[m]->[n].$,
  $\arrow\beta:[r]->[m].$ and $\dim(x)=n$).
\end{itemize}

In time honoured fashion, we will usually ``overload'' $\dim$ and
$\cdot$, using them to denote the dimension and action functions of
whichever simplicial set we happen to be discussing at any given
time. This substantially simplifies our notation and obviates the need
to explicitly name these functions when introducing a new simplicial
set.

Furthermore, we will often refer to the elements of simplicial set $X$
as its {\em simplices} and say that $y\in X$ is a {\em face\/} of a
simplex $x\in X$ if there is some simplicial operator $\alpha$ such
that $x\cdot\alpha=y$.

We will often use the notation $X_n$ to denote the set of {\em
  $n$-simplices\/} $\{x\in X\mid\dim(x)=n\}$ of $X$ and if
$\arrow\alpha:[m]->[n].$ is a simplicial operator then we adopt the
notation $\arrow X_\alpha:X_n->X_m.$ for the
function which maps $x\in X_n$ to $X_\alpha(x)\defeq
x\cdot\alpha\in X_m$.
\end{obs}

\begin{obs}[simplicial maps as action preserving functions] When
  simplicial sets are expressed in this way it is most natural to
  consider simplicial maps to be action preserving functions. In other
  words, a simplicial map
  $\arrow f:\triple<X;\dim;\cdot>->\triple<Y;\dim;\cdot>.$
  consists of a function $\arrow f: X-> Y.$ of sets satisfying:
\begin{itemize}
\item $\dim(f(x))=\dim(x)$ for each $x\in X$, and
\item $f(x\cdot\alpha)=f(x)\cdot\alpha$ for each $x\in X$ and
  $\alpha\in\Delta$ with $[\dim(x)]=\cod(\alpha)$.
\end{itemize}
\end{obs}

\begin{obs}[simplicial subsets]
A subset $Y$ of a simplicial set $X$ is called a {\em simplicial
  subset}, denoted by $Y\subseteq_s X$, if it is closed in $X$ under
the action of $\Delta$. If $Y\subseteq_s X$ then $Y$ becomes a
simplicial set, by inheriting the action structure of $X$, and the
inclusion becomes a simplicial map $\overinc\subseteq_s:Y->X.$.

Notice that intersections and unions of simplicial subsets of $X$, as
mere subsets, are again closed in $X$ under the action of $\Delta$,
in other words they are themselves simplicial subsets of $X$.

Suppose that $A$ is a subset of a simplicial set $X$ then
$\simpc{A}$, its {\em simplicial closure\/} in $X$, is the smallest
simplicial subset of $X$ that contains $A$. Explicitly:
\[\simpc{A}=\{x\in X\mid (\exists \alpha\in\Delta, y\in A)\text{
  s.t. } x = y\cdot\alpha\}\]
\end{obs}

\begin{obs}[why {\em simplicial\/} sets?] The books by Gabriel and Zisman
  \cite{GabrielZisman:1967:CFHT} and May \cite{May:1967:Simp}
  are considered the canonical references for the theory of simplicial
  sets. In particular, they explain how these structures may be used
  to provide combinatorial representations of topological structures
  and expounds at some length upon their homotopy theory.
  
  For our purposes here, it is worth simply observing the following:
\begin{itemize}
\item If $x$ is an $n$-dimensional element of a simplicial set $X$,
  then $x\cdot\vertex^n_i$ is thought of as its $i^{\text{th}}$ {\em
    vertex\/} and $x\cdot\face^n_i$ is thought of as its unique
  $(n-1)$-dimensional face {\bf not} containing that vertex.
\item The simplicial identity for elementary face maps ensures that
  the various faces of a simplex agree appropriately at their
  boundaries.
\end{itemize}
\end{obs}

\begin{defn}[degenerate simplices]\label{defn.degen} 
  We say that an $r$-simplex $x$ of a simplicial set $X$ is {\em
    degenerate at $k$\/} if $x=x'\cdot\degen^{r-1}_k$ for some
  $(r-1)$-simplex $x'$.  Of course $\degen^{r-1}_k$ has right inverses
  $\face^r_k$ and $\face^r_{k+1}$, therefore such an $x'$ is unique,
  if it exists, and $x'=x\cdot\face^r_k = x\cdot\face^r_{k+1}$.
  
  We also say that a simplex is {\em degenerate\/} if it is degenerate
  at $k$ for some $k$. Equivalently, $x\in X$ is degenerate iff there
  is a simplicial operator $\alpha$ and a compatible simplex $x'\in X$
  with $\dim(x')<\dim(x)$ and $x=x'\cdot\alpha$.  We often use the
  notation $\tilde{X}$ to denote the set of those simplices of $X$
  which are {\bf not} degenerate.
\end{defn}

\begin{obs}\label{degen.obs} Let $x$ be an simplex in a simplicial set
  $X$, then
\begin{enumerate}[(a)] 
\item\label{degen.obs.a} $x$ is degenerate iff there exists some $k$
  such that $x$ is degenerate at $k$.
\item\label{degen.obs.b} $x$ is degenerate at $k$ iff $x=x'\cdot\beta$
  for a simplex $x'\in X$ and a simplicial operator $\beta$ with
  $\beta(k)=\beta(k+1)$.
\item\label{degen.obs.c} if $x$ is degenerate at $k$ then
  $x\cdot\alpha$ is degenerate for any face operator $\alpha$ with
  $k,k+1\in\im(\alpha)$.
\end{enumerate}
\end{obs}

\begin{lemma}[the Eilenberg-Zilber lemma]\label{lemma.ez} If $x$ is a
  simplex of a simplicial set $X$ then there exists a {\bf unique}
  pair $\pair<\bar{x};\alpha^x>$ such that $\bar{x}\in X$ is
  non-degenerate, $\alpha^x$ is a degeneracy operator and
  $x=\bar{x}\cdot\alpha^x$. We call this unique pair
  $\pair<\bar{x};\alpha^x>$ the {\em EZ-decomposition} of $x$.
\end{lemma}

\begin{proof} See the proof in \cite{GabrielZisman:1967:CFHT}.
\end{proof}

\begin{obs}[limits and colimits in $\Simp$]\label{simplicial.(co)limits}
  As a category of functors into the (small) complete and cocomplete
  category $\Set$, $\Simp$ also has all (small) limits and colimits.
  These are formed in the usual way for such categories, that is to
  say {\em point-wise} in $\Set$.
  
  It is useful to expand upon this definition in terms of the partial
  action description of simplicial sets given above. To this end,
  suppose that $\mathbb{C}$ is a small category and that
  $\arrow \dgm:\mathbb{C}->\Simp.$ is a functor (which we call a
  {\em diagram\/} of $\mathbb{C}$ in $\Simp$):\vspace{1.5ex}

  \noindent {\bf Limits} The elements of $\lim(\dgm)$ are families
  $\{x_i\in \dgm(i)\}_{i\in\obj(\mathbb{C})}$ satisfying the
  conditions that:
  \begin{itemize}
  \item for each pair $i,j\in\obj(\mathbb{C})$ we have
    $\dim(x_i)=\dim(x_j)$,
  \item if $\arrow f: i-> j.$ is an arrow in $\mathbb{C}$ then
    $\dgm(f)(x_i)=x_j$.
  \end{itemize}
  The dimension of $\boldsymbol{x}$ is the common dimension of its
  components $x_i$, and the action of an operator $\alpha$ is given by
  \begin{math}
    \{x_i\}_{i\in\obj(\mathbb{C})}\cdot\alpha =
    \{x_i\cdot\alpha\}_{i\in\obj(\mathbb{C})}
  \end{math}.\vspace{1.5ex}

  \noindent{\bf Colimits} We form the colimit $\colim(\dgm)$ in $\Set$, whose
  elements are equivalence classes of pairs $\pair<i;x>$ with
  $i\in\obj( \mathbb{C})$ and $x\in \dgm(i)$ under the equivalence
  relation $\sim$ generated by the relation
  $\pair<i;x>\sim_1\pair<j;y>$ which holds when there exists an arrow
  $\arrow f:i->j.$ in $\mathbb{C}$ with $\dgm(f)(x) = y$.  Following
  tradition, we use the notation $[i,x]$ to denote the equivalence
  class of a pair $\pair<i;x>$ under this equivalence relation.
  
  Notice that if $\arrow f: i-> j.$ witnesses that
  $\pair<i;x>\sim_1\pair<j;y>$, then $\dim(y) = \dim(\dgm(f)(x)) =
  \dim(x)$ and if $\alpha$ is an operator with codomain $[\dim(x)]$
  then $y\cdot\alpha = f(x)\cdot\alpha = f(x\cdot\alpha)$. It follows
  that the operations $\dim([i,x]) = \dim(x)$ and ${[i,x]}\cdot\alpha
  = [i,x\cdot\alpha]$ on the equivalence classes of $\sim$ are well
  defined and make $\colim(\dgm)$ into a simplicial set, the colimit
  of $\dgm$ in $\Simp$.
\end{obs}

\begin{obs}\label{union.qua.widepo}
  The fact that limits and colimits in $\Simp$ are constructed
  point-wise in $\Set$ implies that we may immediately ``lift'' many
  of their properties from there. For instance the following useful
  result, which holds in $\Set$ and thus immediately lifts to $\Simp$,
  allows us to describe a union of simplicial subsets as a certain
  kind of colimit, sometimes called a {\em wide pushout}, whose
  vertices are (intersections of) those subsets.  To be precise, if
  $Y$ is a simplicial set and $Y^{(i)}\subseteq Y$ (for $i=1,2,...,r$)
  is a family of simplicial subsets with $Y=\bigcup_l Y^{(l)}$ then
  the cocone of inclusion maps $\overinc:Y^{(i)}->Y.$ and
  $\overinc:Y^{(i)}\cap Y^{(j)} ->Y.$ ($i<j$) presents $Y$ as the
  colimit of a diagram $\dgm_{\{Y^{(i)}\}}$ in $\Simp$ consisting of
  the stratified sets $Y^{(i)}$ and $Y^{(i)}\cap Y^{(j)}$ and
  inclusions $\overinc:Y^{(i)}\cap Y^{(j)}-> Y^{(i)}.$ and
  $\overinc:Y^{(i)}\cap Y^{(j)}-> Y^{(j)}.$ ($i<j$) between them:
  \begin{equation}\label{union.qua.widepo.diag}
    \xymatrix@R=0.4em@C=0.6em{
      && {Y^{(i)}\cap Y^{(j)}}\ar@{u(->}[ddll]\ar@{u(-->}[dddd]
      \ar@{u(->}[ddrr] && \\
      && && \\
      {Y^{(i)}}\ar@{u(-->}[ddrr] && && {Y^{(j)}}\ar@{u(-->}[ddll] \\
      && && \\
      && {Y} && }
  \end{equation}
\end{obs}

\begin{obs}[connected components]
  \label{conn.cpt.simp}
  For each set $S\in\Set$ we have a corresponding {\em discrete\/}
  simplicial set $\dis(S)$ defined by $\dis(S)_n = S$ ($\forall
  n\in\mathbb{N}$) and $\dis(S)_{\alpha} = \id_S$ (for all simplicial
  operators $\alpha\in\Delta$). This construction is clearly provides
  us with a functor $\arrow\dis:\Set->\Simp.$ which is easily shown
  to be fully faithful. It is also clear that a simplicial set is in
  the replete image of this functor iff each of its simplices of
  dimension $>0$ is degenerate.
  
  We also have an adjunction $\cpt_0\dashv\dis$ and the set
  $\cpt_0(X)$ derived from a simplicial set $X$ is known as its {\em
    set of connected components\/}, which may be constructed
  explicitly using the coequaliser
  \begin{displaymath}
    \xymatrix@R=2ex@C=8em{
      {X_1}\ar@<1ex>[r]^{\textstyle X_{\vertex^1_0}}
      \ar@<-1ex>[r]_{\textstyle X_{\vertex^1_1}} & {X_0} 
      \ar@{->>}[r]^{\textstyle q_X} & {\cpt_0(X)}
      }
  \end{displaymath}
  in $\Set$. In other words, we start with the set of 0-simplices of
  $X$ and form $\cpt_0(X)$ from it by identifying any pair of
  0-simplices which are the 0-faces of some 1-simplex.
\end{obs}

\begin{obs}[standard simplices in $\Simp$]\label{standard.simplices}
  \label{yoneda}
  The classical Yoneda functor $\arrow\Delta:\Delta->\Simp.$ carries
  each object $[n]$ of $\Delta$ to the {\em representable\/}
  simplicial set $\Delta[n]\defeq\Delta({-},[n])$. This is called the
  {\em standard $n$-simplex\/} and is given explicitly by:
\begin{itemize}
\item $\Delta[n]_m = \Delta([m],[n])$ the set of simplicial
  operators $\arrow\gamma:[m]->[n].$,
\item if $\arrow\beta:[r]->[m].$ is a simplicial operator that is
  compatible with $\gamma\in\Delta[n]_m$ then $\beta$ and $\gamma$
  are composable and we may define $\gamma\cdot\beta =
  \gamma\circ\beta\in\Delta[n]_r$.
\end{itemize}
Furthermore, it carries the operator $\arrow\alpha:[n]->[m].$ to a
simplicial map $\arrow\Delta(\alpha):\Delta[n]-> \Delta[m].$ given by
post-composition $\Delta(\alpha)(\gamma)=\alpha\circ\gamma$.

Given an $n$-simplex $x\in X$ in a simplicial set $X$, we
adopt the notation $\arrow\yoneda{x}:\Delta[n]->X.$ for the simplicial
map defined by $\yoneda{x}(\gamma)\defeq x\cdot\gamma$. Of course,
Yoneda's lemma applied in the simplicial set context states that the
function $\arrow\yoneda{.}:X_n->\Simp\pair<\Delta[n];X>.$ is inverse
to the evaluation function which maps $f\in\Simp\pair<\Delta[n];X>$ to
the $n$-simplex $f(\id_{[n]})$ in $X$.
\end{obs}

\begin{obs}[a couple of observations regarding {$\Delta[n]$}] 
  \label{stand.simp.obs} 
  An $r$-simplex $\alpha$ of $\Delta[n]$ is non-degenerate if and only
  if, as a simplicial operator $\arrow\alpha:[r]->[m].$, it is
  injective (a face operator). It follows that, if $r>n$ then no such
  simplicial operator can be injective and so all of the $r$-simplices
  of $\Delta[n]$ must be degenerate. Consequently, the only
  non-degenerate $n$-simplex of $\Delta[n]$ is the identity
  $\arrow\id_{[n]}:[n]->[n].$.
\end{obs}

\begin{notation}[the boundaries of standard simplices] 
  \label{simplex.bound}
  The {\em boundary\/} $\boundary\Delta[n]$ of the standard
  $n$-simplex $\Delta[n]$ is the simplicial subset generated by the
  set of all the $(n-1)$-faces $\{\face^n_i\mid i=0,\dots,n\}$ of
  $\Delta[n]$. Equivalently, $\boundary\Delta[n]$ is the largest
  simplicial subset of $\Delta[n]$ which does not contain $\id_{[n]}$.
  A simplex $\alpha$ of $\Delta[n]$ is in $\boundary\Delta[n]$ iff it
  is {\bf not} a degeneracy operator.
  
  We can express $\boundary\Delta[n]$ as the union
  $\bigcup_{i=0}^n\simpc{\{\face^n_i\}}$ of the simplicial subsets
  $\simpc{\{\face^n_i\}}\subseteq \boundary\Delta[n]$ which are the
  images of the simplicial injections $\inc\Delta(\face^n_i):
  \Delta[n-1]->\Delta[n].$ and are thus isomorphic to $\Delta[n-1]$.
  Also, if $i<j$ then we have $\alpha\in\simpc{\{\face^n_i\}}\cap
  \simpc{\{ \face^n_j\}}$ iff $\alpha$ factors through both
  $\face^n_i$ and $\face^n_j$ iff $i,j\notin\im(\alpha)$ iff
  $\alpha\in\simpc{\{\face^n_j\circ \face^{n-1}_i\}}$ therefore
  $\simpc{\{\face^n_i\}} \cap \simpc{\{\face^n_j\}} =
  \simpc{\{\face^n_j\circ \face^{n-1}_i\}}$ which is the image of the
  simplicial injection $\inc\Delta(\face^n_j\circ \face^{n-1}_i)=
  \Delta(\face^n_i\circ \face^{n-1}_{j-1}):\Delta[n-2] ->\Delta[n].$
  and is thus isomorphic to $\Delta[n-2]$.  Furthermore, it is also
  clear that under these isomorphisms the inclusions of
  $\simpc{\{\face^n_i\}} \cap \stratc{\{\face^n_j\}}$ into
  $\simpc{\{\face^n_i\}}$ and $\simpc{\{\face^n_j\}}$ correspond to
  the simplicial injections $\inc\Delta(\face^{n-1}_{j-1}):
  \Delta[n-2]->\Delta[n-1].$ and $\inc\Delta(\face^{n-1}_i):
  \Delta[n-2]->\Delta[n-2].$ respectively.
  
  Applying observation~\ref{union.qua.widepo} and the isomorphisms of
  the last paragraph we see that $\boundary\Delta[n]$ may be expressed
  as the wide pushout
  \begin{equation*}
    \xymatrix@R=0.5em@C=1.5em{
      && {\Delta[n-2]}\ar@{u(->}[ddll]_{\Delta(\face^{n-1}_{j-1})}
      \ar@{u(-->}[dddd]\ar@{u(->}[ddrr]^{\Delta(\face^{n-1}_i)} && \\
      && && \\
      {\Delta[n-1]}\ar@{u(-->}[ddrr]_{\Delta(\face^n_i)} && 
      && {\Delta[n-1]}\ar@{u(-->}[ddll]^{\Delta(\face^n_j)} \\
      && && \\
      && {\boundary\Delta[n]} && }
  \end{equation*}
  of standard simplices.  Applying Yoneda's lemma, and the colimiting
  property of this wide pushout, we see that a simplicial map $\arrow
  f:\boundary\Delta[n]->X.$ corresponds to a family of
  $(n-1)$-simplices $x_0,x_1,...,x_n\in X_{n-1}$ satisfying the
  simplicial identities $x_i\cdot\face^{n-1}_j =
  x_j\cdot\face^{n-1}_{i-1}$ for $j<i$. This data is often referred to
  as a {\em $(n-1)$-dimensional cycle in $X$\/} and we say that such a
  cycle is the {\em boundary\/} of an $n$-simplex $x\in X$ if and only
  if $x_i=x\cdot\face^n_i$ for each $i=0,1,...,n$.
\end{notation}

\begin{notation}[horns] 
  \label{simplex.horn}
  For each $n\in\mathbb{N}$ and $k=0,\dots,n$, the simplicial set
  $\Lambda_k[n]$, called the {\em standard $(n-1)$-dimensional
    $k$-horn}, is the smallest simplicial subset of $\Delta[n]$
  containing the set of $(n-1)$-simplices $\{\face^n_i\mid
  i=0,\dots,k-1,k+1,\dots,n\}$. We say that such a horn is an {\em
  inner horn\/} whenever $0<k<n$ otherwise we say that it is an {\em outer
  horn}.
  
  Arguing as in notation~\ref{simplex.bound}, we see that a simplicial
  map $\arrow f:\Lambda_k[n]-> X.$ corresponds to a family of
  $(n-1)$-simplices $x_0,x_1,...,x_{k-1},x_{k+1},...,x_n\in X_{n-1}$
  satisfying the simplicial identities $x_i\cdot\face^{n-1}_j =
  x_j\cdot\face^{n-1}_{i-1}$ for $j<i$. This data is often referred to
  as a {\em $(n-1)$-dimensional $k$-horn in $X$}. We say that an
  $n$-simplex $x\in X$ {\em fills\/} such a horn if and only if
  $x_i=x\cdot\face^n_i$ for each $i=0,...,k-1,k+1,...,n$.
\end{notation}

\begin{obs}[the dual of a simplicial set]\label{simp.dual} 
  The {\em dual\/} $X^{\mathord{\circ}}$ of a simplicial set $X$ is
  obtained by pre-composing with the functor
  $\arrow(\mathord{-})^{\mathord{\circ}}:\Delta->\Delta.$. In
  action terms, $X^{\mathord{\circ}}$ has the same set of simplices as
  $X$ but an action $*$ given by
  $x*\alpha=x\cdot\alpha^{\mathord{\circ}}$
  
  Notice that the dualising functor
  $\arrow(\mathord{-})^\circ:\Simp->\Simp.$ is again strictly
  involutive, in the sense that its composite with itself is the
  identity on $\Simp$. This is an immediate consequence of the fact
  that $(\mathord{-})^\circ$ is strictly involutive on $\Delta$.
\end{obs}

\subsection{Semi-Simplicial Sets}

\begin{defn}[semi-simplicial sets] 
  Let $\fDelta$ (resp.\ $\afDelta$) denote the subcategory of face
  operators in $\Delta$ (resp.\ $\aDelta$). The category $\SSimp$ of
  {\em semi-simplicial sets\/} and {\em semi-simplicial 
    maps\/} between them is simply the functor category
  $\funcat[\op\fDelta,\Set]$.
\end{defn}
  
\begin{notation}[the category of pointed objects]
  We may also characterise the monoidal category of face operators
  $\triple<\afDelta; \dsum;[-1]>$ in the spirit of
  lemma~\ref{univ.delta}. However, to do so, we must first introduce
  the simple notion of a {\em pointed object\/} $\pair<X;e>$ in a
  monoidal category $\triple<\mathcal{E}; \otimes; I>$, which consists
  of an object $X\in\mathcal{E}$ equipped with a map $\arrow e:I->X.$
  called a {\em point}. A map of pointed objects $\arrow f:
  {\pair<X;e>}->{\pair<X';e'>}.$ is simply a map of underlying objects
  $\arrow f:X->X'.$ which ``preserves the point'' in the sense that
  $f\circ e = e'$; in other words the category of pointed objects
  $\Pt\triple<\mathcal{E}; \otimes; I>$ is simply the comma category
  $I\downarrow\mathcal{E}$.
  
  If $\arrow \ladj:{\triple<\mathcal{E}; \otimes; I>}->
  {\triple<\mathcal{E}'; \otimes'; I'>}.$ is a monoidal functor then
  we may define a functor $\Pt(\ladj)$ between corresponding
  categories of points by mapping a pointed object $\pair<X;e>$ in
  $\triple<\mathcal{E}; \otimes; I>$ to
  $\pair<\ladj(X);\ladj(e)\circ\nu>$ in $\triple<\mathcal{E}';
  \otimes'; I'>$ and a point preserving map $\arrow f:
  {\pair<X;e>}->{\pair<X';e'>}.$ to $\ladj(f)$.
\end{notation}

\begin{lemma}[the universal property of $\afDelta$]\label{univ.delta.f}
  The pointed object $\pair<[0];\face_0^0>$ in $\triple<\afDelta;
  \dsum;[-1]>$ is universal, in the sense that if $\pair<X;e>$ is a
  pointed object in a monoidal category $\triple<\mathcal{E};
  \otimes;I>$ then there exists an unique strict monoidal functor
  $\arrow \tilde{X}: {\triple<\afDelta; \dsum;[-1]>}->
  {\triple<\mathcal{E};\otimes;I>}.$ such that the pointed object
  $\Pt(\tilde{X})\pair<[0];\face_0^0>=\pair<\tilde{X}([0]);
  \tilde{X}(\face_0^0)>$ is equal to $\pair<X;e>$.
\end{lemma}

\begin{proof} Standard, see Mac~Lane \cite{Maclane:1971:CWM}. It is worth
  remarking that observation~\ref{del.gen.rel} immediately implies
  that the functor $\tilde{X}$ is completely determined by
  $\tilde{X}([n]) = X^{\otimes(n+1)}$ and $\tilde{X}(\face^n_i) =
  X^{\otimes(n-i)}\otimes e \otimes X^{\otimes i}$.
\end{proof}

\begin{obs}[the free simplicial set generated by a semi-simplicial set]
  In terms of the partial action description, a semi-simplicial set
  $\triple<X;\dim;\cdot>$ is merely a set $X$ equipped with a
  dimension function ``$\dim$'' and a partial right action ``$\cdot$''
  by simplicial face operators. Under this interpretation, a
  semi-simplicial map $\arrow f:\triple<X;\dim;\cdot>->\triple<
  Y;\dim;\cdot>.$ is simply a function from $X$ to $Y$ which preserves
  dimensions and right actions by face operators.
  
  Clearly then, every simplicial set (resp.\ map) is also a
  semi-simplicial set (resp.\ map), giving rise to a canonical
  forgetful functor $\arrow\forget:\Simp->\SSimp.$. Since limits and
  colimits are constructed point-wise in these categories it follows
  that this functor preserves them and thus has both left and right
  adjoints.  
\end{obs}

\subsection{Analysing Products of Simplicial Sets - the Theory of Shuffles}

The following observations, regarding products of standard simplices,
are classical and are applied extensively in proving the ``extension''
lemmas of sections~\ref{precomp.sec} and~\ref{comp.sec}.  

\begin{obs}[nerves of partially ordered sets]
  \label{nerves.po}
  Let $\Ord$ denote the category of all partially ordered sets and
  order preserving maps between them. We know that $\Delta$ is a full
  subcategory of $\Ord$ and it follows that we may extend the Yoneda
  functor $\arrow\Delta:\Delta->\Simp.$ to a functor
  $\arrow\Delta:\Ord->\Simp.$ which carries the partially ordered set
  $P$ to the simplicial set $\Delta(P)$ given by ``homming out'' of
  the subcategory $\Delta$ and into $P$. In other words, $\Delta(P)$
  has $n$-simplices which are order preserving maps $\arrow x:[n]->P.$
  and the simplicial operator $\arrow\alpha:[m]->[n].$ acts on such a
  simplex by pre-composition $x\cdot\alpha=\arrow
  x\circ\alpha:[m]->P.$.  If $\arrow f:P->Q.$ is an order preserving
  map then $\Delta(f)$ is the simplicial map obtained by
  post-composing each simplex of $\Delta(P)$ by $f$ to obtain a
  simplex of $\Delta(Q)$. We call $\Delta(P)$ the {\em nerve\/} of the
  partially ordered set $P$.
  
  We may construct a left adjoint to the nerve functor
  $\arrow\Delta:\Ord->\Simp.$ by left Kan extending the inclusion
  functor $\overinc:\Delta->\Ord.$ along the Yoneda embedding
  $\arrow\Delta:\Delta->\Simp.$ (cf.\ Kelly~\cite{Kelly:1982:ECT} for
  instance). It follows that the nerve construction preserves all
  limits and in particular that it preserves the product
  $[n]\times[m]$ giving a canonical isomorphism
  $\Delta([n]\times[m])\cong\Delta[n]\times\Delta[m]$ under which we
  shall usually identify these simplicial sets.
  
  More explicitly, recall from observations \ref{standard.simplices}
  and \ref{simplicial.(co)limits} that the $r$-simplices of
  $\Delta[n]\times\Delta[m]$ are pairs $\pair<\alpha;\beta>$ of
  simplicial operators $\arrow\alpha:[r]->[n].$ and
  $\arrow\beta:[r]->[m].$ and that these are acted upon by right
  composition $\pair<\alpha;\beta>\cdot\gamma =
  \pair<\alpha\circ\gamma; \beta\circ\gamma>$. The corresponding
  $r$-simplex in $\Delta([n]\times[m])$ is the unique order preserving
  map $\overarr:[r]->[n]\times[m].$ induced by the universal property
  of $[n]\times[m]$ applied to the pair $\pair<\alpha;\beta>$. As is
  traditional, we shall blur our notation a little by using
  $\pair<\alpha;\beta>$ to denote both the pair of simplicial
  operators which represent a simplex in $\Delta[n]\times\Delta[m]$
  and the single order preserving map $\arrow\pair<\alpha;\beta>:[r]->
  [n]\times[m].$, given by $\pair<\alpha;\beta>(i)=\pair<\alpha(i);
  \beta(i)>$, which represents the corresponding simplex in
  $\Delta([n]\times[m])$. 
  
  In general, our basic notational conventions for simplicial
  operators extend in a natural way to partially ordered sets and
  maps. For instance, if $\arrow f:P->Q.$ is an arrow of $\Ord$ then
  we use the notation $\im(f)$ to denote its image $\{f(p)\in Q\mid
  p\in P\}$.  Consequently, under the convention of the last
  paragraph, it follows that $\im\pair<\alpha;\beta>$ denotes the
  subset $\{\pair<\alpha(i);\beta(i)>\mid i\in[r]\}\subseteq
  [n]\times[m]$ which we think of as being the set of vertices of the
  $r$-simplex $\pair<\alpha;\beta>$ in $\Delta[n]\times\Delta[m]$. We
  leave other such mild generalisations to the imagination of the
  reader.
\end{obs}

\begin{defn}[shuffles] The non-degenerate $(n+m)$-simplices of
  $\Delta[n]\times\Delta[m]$ are called {\em shuffles}. Exploiting the
  identification of $\Delta[n]\times\Delta[m]$ with the nerve
  $\Delta([n]\times[m])$, as discussed in the last observation, we may
  think of a shuffle as a (strict) path of maximal length $(n+m)$ in
  the ordered set $[n]\times[m]$ as depicted in figure~\ref{pic.zero}.
\end{defn}

\begin{lemma}[properties of products and shuffles]\label{shuffles} 
  The following properties hold for the simplicial set
  $\Delta[n]\times\Delta[m]$:
\begin{enumerate}[(1)]
\item\label{prod.simp.1} Integer addition gives a strictly
  order preserving surjective map $\arrow
  {+}:[n]\times[m]->[n+m].$, to which we may apply the nerve functor of
  observation~\ref{nerves.po} to give a simplicial map
  $\arrow\soplus:\Delta[n]\times\Delta[m]-> \Delta[n+m].$ which
  carries a simplex $\pair<\alpha;\beta>$ to the simplex
  $\alpha+\beta$ given point-wise by $(\alpha+\beta)(j)=\alpha(j)+
  \beta(j)$.
\item\label{prod.simp.2} A simplex $\pair<\alpha;\beta>$ is degenerate
  in $\Delta[n]\times\Delta[m]$ if and only if $\alpha\soplus\beta$ is
  degenerate in $\Delta[n+m]$.  It follows that any simplex of
  $\Delta[n]\times\Delta[m]$ of dimension greater than $(n+m)$ is
  degenerate.
\item\label{prod.simp.3} The projection map that takes an $(n+m)$-simplex
  $\pair<\alpha;\beta>$ of $\Delta[n]\times\Delta[m]$ to the operator
  $\arrow\alpha:[n+m] ->[n].$ establishes a bijection between
  the set of shuffles of $\Delta[n]\times\Delta[m]$ and the set of
  degeneracy operators $\arrow\alpha:[n+m]->[n].$.
\item\label{prod.simp.4} There exists a bijection between the set of
  degeneracy operators $\arrow\alpha:[n+m]->[n].$ and the set
  of simplicial operators $\arrow\gamma:[n]->[m].$ with
  $\gamma(n)=m$.
\item\label{prod.simp.5a} Suppose that $\arrow\gamma:[n]->[m].$
  is a simplicial operator with $\gamma(n) = m$ and that
  $\shuffle\gamma$ is the shuffle it corresponds to under the
  bijections of (\ref{prod.simp.3}) and (\ref{prod.simp.4}).  Let
  $\pair<\alpha;\beta>$ be an arbitrary $r$-simplex of
  $\Delta[n]\times\Delta[m]$, then the following propositions hold:
  \begin{enumerate}
    \item if $\pair<\alpha;\beta>$ is an $r$-dimensional face of 
      $\shuffle\gamma$ then it is the face obtained by applying 
      the operator $\arrow\alpha\soplus\beta:[r]->[n+m].$ to the given 
      shuffle, or in other words:
      \begin{equation*}
        \pair<\alpha;\beta>=
        \shuffle\gamma\cdot(\alpha\soplus\beta) = 
        \pair<\shufflel\gamma\circ(\alpha\soplus\beta);
        \shuffler\gamma\circ(\alpha\soplus\beta)>
      \end{equation*} 
    \item for an $l\in[r]$ we have $\alpha(l)=\shufflel\gamma\circ(
      \alpha\soplus\beta)(l)$ and $\beta(l)=\shuffler\gamma\circ(
      \alpha\soplus\beta)(l)$ if and only if the inequalities
      \begin{equation}\label{shuff.face.ineq}
        \gamma(\alpha(l)-1)\leq\beta(l)\leq\gamma(\alpha(l))
      \end{equation}
      hold (in which we take $\gamma(-1)\defeq 0$ where necessary),
    \item consequently, $\pair<\alpha;\beta>$ is a face of
      $\shuffle\gamma$ if and only if the inequalities
      of~(\ref{shuff.face.ineq}) hold for all $l\in[r]$.
    \end{enumerate}
\item\label{prod.simp.5} Every simplex in $\Delta[n]\times\Delta[m]$
  is a face of some shuffle.
\end{enumerate}
\end{lemma}

\begin{proof} 
  These results are fundamentally classical in nature, and we leave
  their detailed verification up to the reader. However we sketch the
  main plot points.
  
  Part (\ref{prod.simp.1}) is clear from the statement and it is
  easily verified that part (\ref{prod.simp.2}) follows from the
  observation that $\arrow {+}:[n]\times[m]->[n+m].$ is strictly order
  preserving. Part (\ref{prod.simp.3}) is a direct consequence of the
  fact that a simplex $\pair<\alpha;\beta>$ is a shuffle if and only
  if $\alpha\soplus\beta = \id_{[n+m]}$, which itself follows directly
  from (\ref{prod.simp.1}) and (\ref{prod.simp.2}).
  
  Part (\ref{prod.simp.4}) is more interesting, the intuition here is
  that a shuffle $\pair<\alpha;\beta>$, as depicted in
  figure~\ref{pic.zero}, may be fully specified by providing a list of
  ``plateau'' levels $\gamma(0),\gamma(1),\dots,\gamma(n-1)$, one for
  each column in the figure (as marked).
  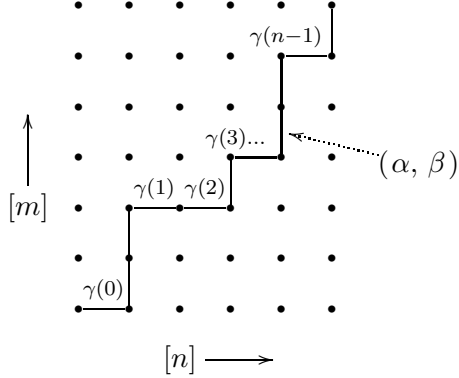
\begin{figure}[t]
    \begin{displaymath}
      \def\vertchar{\scriptscriptstyle\bullet}
      \xymatrix@!0@C=1.75em@R=1.75em{
        & *[o]{\vertchar} & *[o]{\vertchar} & *[o]{\vertchar} & *[o]{\vertchar} & 
        *[o]{\vertchar} & *[o]{\vertchar} \\ 
        & *[o]{\vertchar} & *[o]{\vertchar} & *[o]{\vertchar}
        & *[o]{\vertchar} & *[o]{\vertchar} & *[o]{\vertchar} \\ 
        & *[o]{\vertchar} & *[o]{\vertchar} & *[o]{\vertchar} & *[o]{\vertchar} & 
        *[o]{\vertchar} & *[o]{\vertchar} \\ 
        & *[o]{\vertchar} & *[o]{\vertchar} & *[o]{\vertchar} & *[o]{\vertchar} & 
        *[o]{\vertchar} & *[o]{\vertchar} \\ 
        [m]\ar[uu] & *[o]{\vertchar} & *[o]{\vertchar} & *[o]{\vertchar} & *[o]{\vertchar} & 
        *[o]{\vertchar} & *[o]{\vertchar} \\ 
        & *[o]{\vertchar} & *[o]{\vertchar} & *[o]{\vertchar} & *[o]{\vertchar} &
        *[o]{\vertchar} & *[o]{\vertchar} \\ 
        & *[o]{\vertchar}
        \ar@{-}'[r]^{\gamma(0)}'[ruu]'[ruur]^{\gamma(1)}'[ruurr]^{\gamma(2)}'[ruurru]'[ruururr]^<{\gamma(3)...}
        '[ruururruu]_(0.25){}="label"'[ruururruur]^<{\gamma(n-1)}[ruururruuru] &
        *[o]{\vertchar} &
        *[o]{\vertchar} & *[o]{\vertchar} & 
        *[o]{\vertchar} & *[o]{\vertchar}
        \save []+<3em,5em>*{\pair<\alpha;\beta>}\ar@{.>}"label"\restore \\
        & & & [n] \ar[rr] & & & 
        }
    \end{displaymath}
    \caption{A shuffle $\pair<\alpha;\beta>$ in $\Delta[n]\times\Delta[m]$}
    \label{pic.zero}
  \end{figure}
  A moment's reflection reveals that $\gamma(i)$, the height of the
  ``plateau'' in the $i^{\text{th}}$ column, may be given by the
  formula:
  \begin{equation}\label{plat.lev}
    \gamma(i)=\max\left\{j\in[n+m]\mid \alpha(j)\leq i\right\} - i
  \end{equation}
  Of course, from observation~\ref{simp.op.adj} we know that the
  degeneracy operator $\alpha$ has a right adjoint $\alpha^r$ and that
  $\alpha^r(i)=\max\{j\in[n+m]\mid\alpha(j)\leq i\}$, so we may
  re-express this as $\gamma(i) = \alpha^r(i) - i$.  Furthermore,
  observation~\ref{simp.op.adj}(\ref{so.adj.rinv}) reveals that
  $\alpha^r$ is a face operator and so we can apply
  observation~\ref{fd.char}(\ref{face.char.2}) to construct
  $\arrow\gamma:[n]->[m].$ as the unique simplicial operator
  satisfying the equation:
  \begin{equation*}
    \id_{[n]}\soplus\gamma = \alpha^r
  \end{equation*}
  It is now a routine matter to apply observations~\ref{fd.char}
  and~\ref{simp.op.adj} and show that this construction provides the
  bijection postulated in the statement of (\ref{prod.simp.4}).
  
  Given the intuition we developed in the last paragraph, regarding
  the relationship between shuffles $\pair<\alpha;\beta>$ and
  operators $\arrow\gamma:[n]->[m].$, the meaning of
  inequality~(\ref{shuff.face.ineq}) in the statement of
  part~(\ref{prod.simp.5a}) should now be clear. It simply states that
  in the $l^{\text{th}}$ vertex of $\pair<\alpha;\beta>$ is a vertex
  of the shuffle $\shuffle\gamma$ if and only if its vertical ordinate
  $\beta(l)$ lies between the plateaus of $\shuffle\gamma$ immediately
  to either side of its horizontal ordinate $\alpha(l)$ (the levels of
  which we know to be $\gamma(\alpha(l)-1)$ and $\gamma(\alpha(l))$
  respectively). Given this observation, it is a matter of
  straightforward calculation, using the results we have already
  established, to show that the results of part~(\ref{prod.simp.5a})
  also hold.
  
  Finally, to demonstrate that part~(\ref{prod.simp.5}) holds, we
  start with the $r$-simplex $\pair<\alpha;\beta>$ and define a
  simplicial operator $\arrow\gamma:[n]->[m].$ by
  \begin{equation}\label{gamma.ab.def}
    \gamma(i) = \begin{cases}
      \min\{\beta(l)\mid l\in[r]\wedge\alpha(l)>i\} & \text{if
        $\alpha(r)>i$,} \\
      m & \text{otherwise.}
    \end{cases}
  \end{equation}
  for which it is a matter of routine verification to demonstrate that
  the condition of part~(\ref{prod.simp.5a})(c) holds. It follows that
  with $\gamma$ defined in this way we can infer that our simplex
  $\pair<\alpha;\beta>$ is a face of the shuffle $\shuffle\gamma$.
\end{proof}

\begin{notation}\label{shuffles.notation} We will adopt the following
  notational conventions with respect to the results embodied in the
  previous lemma:
  \begin{itemize}
  \item Rather than consider simplicial operators
    $\arrow\gamma:[n]->[m].$ with $\gamma(n)=m$, which were
    useful in proving the results of the lemma, we will instead work
    with the set of {\em all\/} simplicial operators
    $\arrow\gamma:[n-1]->[m].$ which, whenever necessary, we
    implicitly extend to $[n]$ by setting $\gamma(n)=m$.
  \item The notation $\shuffle\gamma$ will be used
    to denote the shuffle associated with a simplicial operator
    $\arrow\gamma:[n-1]->[m].$ (as in
    lemma~\ref{shuffles}(\ref{prod.simp.5a})).
  \item Given an $r$-simplex $\pair<\alpha;\beta>$ in
    $\Delta[n]\times\Delta[m]$, the notation
    $\arrow\topop{\shuffle{}}:[n-1]->[m].$ will be used to denote the
    simplicial operator defined in display~\eqref{gamma.ab.def} of the
    proof of lemma~\ref{shuffles}(\ref{prod.simp.5}).
  \end{itemize}
\end{notation}

\begin{obs}[A linear ordering of shuffles] 
  The set of simplicial operators from $[n-1]$ to $[m]$ and
  thus, by the bijection of lemma~\ref{shuffles}, the set of shuffles
  in $\Delta[n]\times\Delta[m]$ may be linearly ordered by the
  ``lexicographic ordering'' relation $\triangleleft$, defined by
  \begin{alignat*}{1}
    \gamma\triangleleft\tau \iff (\exists i\in[n-1])(
    &\gamma(0)=\tau(0)\wedge\gamma(1)=\tau(1)\wedge\dots\wedge \\
    &\gamma(i-1)=\tau(i-1)\wedge\gamma(i)<\tau(i))
  \end{alignat*}
  which extends the (strict) point-wise ordering of
  observation~\ref{2cat.delta}, in other words whenever $\gamma<\tau$
  (ie whenever $\gamma\leq\tau$ and $\gamma\neq\tau$) then we also
  have $\gamma\triangleleft\tau$. In particular, the shuffle
  $\pair<\partproj^{n,m}_1; \partproj^{n,m}_2>$ of partition operators
  (notation~\ref{part.oper.def}) is minimal under the linear ordering
  $\triangleleft$ and $\pair<\partproj^{m,n}_2;\partproj^{m,n}_1>$ is
  maximal under the same relation.  We also often use the notation
  $\trianglelefteq$ to denote the non-strict version of
  $\triangleleft$, in other words $\gamma\trianglelefteq\tau$ if and
  only if $\gamma\triangleleft\tau$ or $\gamma=\tau$.
  
  Let $\#(n,m)$ denote the cardinality of the set of shuffles of
  $\Delta[n]\times\Delta[m]$. We know, from the fact that
  $\triangleleft$ is a linear ordering, that there exists a unique
  enumeration of these shuffles (respectively the simplicial operators
  from $[n-1]$ to $[m]$) $\{\shuffle{i}\}_{i=1}^{\#(n,m)}$
  (respectively $\{\topop{i}\}_{i=1}^{\#(n,m)}$) which is compatible
  with $\triangleleft$ in the sense that $\shuffle{i}\triangleleft
  \shuffle{i'}$ (respectively $\topop{i}\triangleleft\topop{i'}$) if
  and only if $i<i'$. Of course, since $\triangleleft$ extends the
  point-wise order $<$, it follows that in order to demonstrate that
  $i<i'$ it is sufficient to show that $\gamma_i<\gamma_{i'}$.
\end{obs}

\begin{obs}\label{shuffles.maxl} Notice that if $\pair<\alpha;\beta>$
  is an $r$-simplex of $\Delta[n]\times\Delta[m]$ then the associated
  operator $\arrow\gamma:[n-1]->[m].$ of display~\eqref{gamma.ab.def}
  in the proof of lemma~\ref{shuffles}(\ref{prod.simp.5}) is the upper
  bound, in the point-wise ordering $\leq$, of the set of operators
  $\arrow\tau:[n-1]->[m].$ for which $\pair<\alpha;\beta>$ is a face
  of $\shuffle\tau$.
  
  To prove this suppose that $\arrow\tau:[n-1]->[m].$ is such an
  operator and observe that the definition of $\gamma$ tells us that
  for each $i\in[n-1]$ either $m=\gamma(i)$ or there exists an
  $l\in[r]$ with $\alpha(l)> i$ (equivalently $\alpha(l)-1\geq i$) and
  $\gamma(i)=\beta(l)$. In the first case we know that $\tau(i)\leq
  m=\gamma(i)$, since $m$ is the largest integer the codomain of
  $\tau$ and in the second case we may apply
  lemma~\ref{shuffles}(\ref{prod.simp.5a}) to show that
  $\tau(\alpha(i)-1)\leq\beta(i)$ which we combine with the
  (in)equalities of the last sentence and order preservation by $\tau$
  to give $\tau(i)\leq\tau(\alpha(l)-1) \leq\beta(l)= \gamma(i)$ as
  required. So, by the point-wise definition of the ordering $\leq$ on
  simplicial operators, it follows that $\gamma\geq\tau$.
  
  In particular, applying the last observation we see that if $w$ is
  the unique integer with $\gamma=\topop{w}$ then it may be
  characterised as being the largest integer for which $\shuffle{}$ is
  a face of the shuffle $\shuffle{w}$.
\end{obs}



\section{Some Categorical Background}
\label{cat.back.sect}

\subsection{Reflective Full Sub-Categories}\label{append.reflective}

We will often be interested in a context in which we wish to study a
reflective full subcategory of some ``carrier'' category. In these
cases it is often inconvenient to work in the subcategory alone.
Instead, it is usually better to make explicit calculations in the
carrier category and then to reflect the results of these exploits
into the subcategory.

\begin{obs}[reflective sub-categories]\label{refl.full.subcat}
  A {\em reflective full subcategory\/} $\mathcal{E}'$ of
  $\mathcal{E}$ is simply a full subcategory which is {\em replete},
  in the sense that it is closed under isomorphisms in $\mathcal{E}$,
  and for which the inclusion $\inc \incl:\mathcal{E}'->
  \mathcal{E}.$ has a left adjoint $\arrow
  \ladj:\mathcal{E}->\mathcal{E}'.$.
  
  We follow Kelly~\cite{Kelly:1982:ECT} (page~53) by selecting this
  adjunction so that its counit is the identity and concentrating our
  attention on the pair $\pair<\refl;\eta>$ consisting of the
  endo-functor $\arrow \refl=\incl\ladj:\mathcal{E}->
  \mathcal{E}.$ and the unit $\twocell \eta:\id_{\mathcal{E}}->\refl.$ of
  the adjunction $\ladj\dashv\incl$. It is easily seen
  that this pair satisfies the identities $\refl\refl=\refl$ and $\eta
  \refl=\refl\eta=\id_\refl$ (see observations~\ref{whiskering}
  and~\ref{cat.as.2cat} for more regarding the 2-categorical notation
  here) making it into an {\em idempotent monad\/} on $\mathcal{E}$.
  
  Indeed, any such idempotent monad $\pair<\refl;\eta>$ on $\mathcal{E}$
  gives rise to a reflective full subcategory $\mathcal{E}_\refl$ on
  those objects $A\in\mathcal{E}$ which are isomorphic to $\refl(X)$ for
  some $X\in\mathcal{E}$; we might say that $\mathcal{E}_\refl$ is the
  {\em replete image\/} of $\refl$. Notice that an object
  $A\in\mathcal{E}$ is in $\mathcal{E}_\refl$ if and only if the
  corresponding component of the unit $\eta_A$ is an isomorphism. It
  may be easily demonstrated that the required left adjoint to the
  inclusion $\overinc:\mathcal{E}_\refl->\mathcal{E}.$ is provided by $\refl$
  itself.  In this way we get a 1-1 correspondence between idempotent
  monads on $\mathcal{E}$ and its reflective full sub-categories.
  
  We will usually find it easier to identify these two concepts and
  work exclusively with idempotent monads $\pair<\refl;\eta>$.
  Consequently, we refer to the $\refl$ in such a structure as a {\em
    reflector\/} and the natural transformation $\eta$ as its {\em
    unit}.  Usually if we say that $\mathcal{E}_\refl$ is a reflective
  full subcategory of $\mathcal{E}$ we take it as assumed that we
  have an associated idempotent monad $\pair<\refl;\eta>$ which defines
  $\mathcal{E}_\refl$ in the way given above.

  In general, we will use capitals from the end of the alphabet
  $X,Y,Z...$ to denote objects in $\mathcal{E}$ and those at the
  beginning $A,B,C...$ to denote objects in the subcategory
  $\mathcal{E}_\refl$.
\end{obs}

For the next few observations we work within the context of a fixed
category $\mathcal{E}$ and a chosen reflective full subcategory
$\mathcal{E}_\refl$.

\begin{defn}[L-invertible arrows]\label{l-inv.def} 
  An arrow $\arrow f:X->Y.$ in $\mathcal{E}$ is said to be {\em
    L-invertible\/} if and only if its image under the reflector $\refl$
  is an isomorphism.
\end{defn}
  
\begin{obs}[orthogonality] 
  Equivalently, by Yoneda's lemma, an arrow $f$ is L-invertible if
  and only if it is {\em orthogonal\/} to each object
  $A\in\mathcal{E}_\refl$, in the sense that if $\arrow g: X->A.$ is any
  map with $A\in\mathcal{E}_\refl$ then it admits a unique extension
  $\arrow \bar{g}:Y->A.$ along $f$ (i.e. $\bar{g}\circ f = g$):
  \begin {displaymath}
    \let\labelstyle=\textstyle
    \xymatrix@R=3em@C=3em{ {X}\ar[rr]^f\ar[dr]_g & & {Y}
      \ar@{-->}[dl]^{\exists!\bar{g}} \\
      & {A} & }
  \end{displaymath}
  
  We often define reflective full sub-categories in terms of
  orthogonality conditions. In these cases, we start with a set $S$ of
  arrows $\mathcal{E}$, sometimes called a {\em regulus}, and define
  $\mathcal{E}_{\perp S}$ to be the full subcategory of those objects
  in $\mathcal{E}$ which are orthogonal to all of the arrows in $S$.
  Usually, proving that such an $\mathcal{E}_{\perp S}$ is actually
  reflective in $\mathcal{E}$ is a matter of applying some form of
  completion process expressed as a transfinite construction, as
  described in Kelly~\cite{Kelly:1980:Trans}.
  
  In general, if we are given a set of L-invertible arrows $S$ then
  we say that it is {\em adequate to detect\/} objects of
  $\mathcal{E}_\refl$ if we can show that any object in $\mathcal{E}$
  which is orthogonal to all arrows in $S$ is an object of
  $\mathcal{E}_\refl$. Observe that this condition holds iff
  $\mathcal{E}_{\perp S}=\mathcal{E}_\refl$.
\end{obs}

\begin{obs}[fundamental properties of reflective full sub-categories 
  and L-invertible arrows]\label{l-inv} The following observations
  follow directly by standard elementary categorical arguments, which
  we leave up to the reader:
  \begin{enumerate}[1)]
  \item\label{l-inv.limits} If $\arrow \dgm:\mathbb{C}->\mathcal{E}_\refl.$
    is a diagram with limit $\lim(\dgm)$ in $\mathcal{E}$ then $\lim(\dgm)$
    is actually in the subcategory $\mathcal{E}_\refl$.  We simply say
    that {\em $\mathcal{E}_\refl$ is closed in $\mathcal{E}$ under all
      limits which exist there}.
  \item\label{l-inv.colimits} If $\arrow \dgm:\mathbb{C}->\mathcal{E}_\refl.$ is
    a diagram then we may form its colimit in $\mathcal{E}_\refl$ by
    constructing $\colim(\dgm)$ in $\mathcal{E}$ and reflecting it into
    $\mathcal{E}_\refl$ using $\refl$. The corresponding colimit cocone in
    $\mathcal{E}_\refl$ is obtained by composing the colimit cocone in
    $\mathcal{E}$ with the unit map
    $\arrow\eta_{\colim(\dgm)}:\colim(\dgm)->\refl(\colim(\dgm)).$.
  \item\label{l-inv.unit} The condition $\refl\eta=\id_{\refl}$ satisfied by
    our idempotent monad $\pair<\refl;\eta>$ may be usefully restated as
    saying that for each $X\in\mathcal{E}$ the component
    $\arrow\eta_X:X->\refl(X).$ of the unit $\eta$ is L-invertible.
  \item\label{l-inv.inc} If $\mathcal{E}_{\refl'}$ is another reflective
    full subcategory of $\mathcal{E}$ then $\mathcal{E}_\refl$ is
    contained in $\mathcal{E}_{\refl'}$ if and only if every
    $\refl'$-invertible arrow is L-invertible.
  \item\label{l-inv.comp} The class of L-invertible arrows is closed
    under composition, right cancellation of epimorphisms /
    L-invertibles and left cancellation of L-invertibles. In other
    words:
    \begin{itemize}
    \item if $\arrow f:X->Y.$ and $\arrow g: Y->Z.$ are both
      L-invertible then so is their composite $\arrow g\circ f:
      X->Z.$,
    \item if $\arrow g\circ f:X->Z.$ is L-invertible and $\arrow f:
      X-> Y.$ is either epimorphic or L-invertible 
      then the factor $\arrow g: Y->Z.$ is also L-invertible, and
    \item if $\arrow g\circ f:X->Z.$ is L-invertible and either 
      $\arrow g: Y-> Z.$ is L-invertible or $\arrow f:X->Y.$ is epimorphic
      then the factor $\arrow f:X->Y.$ is also L-invertible.
    \end{itemize}
  \item\label{l-inv.retr} Consequently, if the arrows $\arrow m:X->Y.$
    and $\arrow e:Y->X.$ form an inclusion / retraction pair (that is
    $e\circ m=\id_X$) then $m$ is L-invertible if and only if $e$ is
    L-invertible.
  \item\label{l-inv.colim} The class of L-invertible arrows is closed
    under colimits. More explicitly, if the transformation of diagrams
    $\arrow \gamma\colon \dgm\Rightarrow \dgm':
    \mathbb{C}->\mathcal{E}.$ has the property that each of its
    components $\gamma_c$ ($c\in\obj(\mathbb{C})$) is L-invertible,
    then the induced map $\arrow \bar{\gamma}:
    \colim(\dgm)->\colim(\dgm').$ is also L-invertible.
  \item\label{l-inv.po} The class of L-invertible arrows is closed
    under pushout. Explicitly, if the square
    \begin{displaymath}
      \let\labelstyle=\textstyle
      \xymatrix@R=3em@C=3em{
        {X}\ar[r]^f\ar[d]_h & {Y}\ar[d]^k \\
        {Z}\ar[r]_g & \poexcursion {W}
      }
  \end{displaymath}
  is a pushout in $\mathcal{E}$ and $f$ is L-invertible then so is
  $g$.
  \end{enumerate}
\end{obs}

\begin{obs}[L-almost concepts]\label{l-almost}
  We are often interested in proving categorical results about
  $\mathcal{E}_\refl$ by making calculations within $\mathcal{E}$.
  This usually allows us to exploit simpler explicit descriptions of
  objects in $\mathcal{E}$ than are available to us for the
  corresponding objects in $\mathcal{E}_\refl$. However, in order to
  do this we will need to generalise certain categorical concepts.
  Broadly speaking if we are given a categorical property ``X'' of a
  given structure then we say that such a structure in $\mathcal{E}$
  {\em L-almost has property ``X''\/} if applying $\refl$ to it gives
  us a structure which has property ``X'' in $\mathcal{E}_\refl$.  For
  instance an arrow $\arrow f:X->Y.$ is L-almost an epimorphism if the
  arrow $\arrow \refl(f):\refl(X)->\refl(Y).$ is an epimorphism in
  $\mathcal{E}_\refl$.
  
  In particular, if $\arrow \dgm:\mathbb{C}->\mathcal{E}.$ is a
  diagram then we say that the cocone $\nattrans i:\dgm->X.$ in
  $\mathcal{E}$ is {\em L-almost colimiting\/} iff the cone $\nattrans
  \refl(i):\refl(\dgm({-})) ->\refl(X).$ is a colimiting cocone in
  $\mathcal{E}_\refl$. Notice that if the colimit $\colim(\dgm)$
  exists in $\mathcal{E}$ then we know that our cocone induces a
  unique arrow $\arrow \hat{i}:\colim(\dgm)->X.$ and that this is
  L-invertible if and only if our cocone is L-almost colimiting.
  
  Most abstract properties of L-almost constructions may be proved
  simply by reflecting them into the subcategory $\mathcal{E}_\refl$ and
  then exploiting the corresponding traditional categorical theorem in
  that context.
\end{obs}

\begin{obs}[partial adjoints]\label{partial.adj}
  If $\arrow \ladj:\mathcal{E}->\mathcal{F}.$ is a functor then, following
  Kelly~\cite{Kelly:1982:ECT} (page~50), we say that $\ladj$ has a {\em
    partial right adjoint\/} $\forget$ defined on a full subcategory
  $\mathcal{F}'$ of $\mathcal{F}$ if there exists a functor $\arrow
  \forget:\mathcal{F}'-> \mathcal{E}.$ and a family of isomorphisms
  \begin{displaymath}
    \mathcal{F}(\ladj(X),A) \cong \mathcal{E}(X,\forget(A))
  \end{displaymath}
  which is natural in $X\in\mathcal{E}$ and $A\in\mathcal{F}'$.
  
  Suppose now that we are given reflective full sub-categories
  $\mathcal{E}_\refl$ and $\mathcal{F}_{\refl'}$ of $\mathcal{E}$ and
  $\mathcal{F}$ respectively. In such cases, we are predominantly
  interested in those partial right adjoints which are defined on
  $\mathcal{F}_{\refl'}$ {\bf and} map into $\mathcal{E}_\refl$. The
  following results hold by elementary arguments, in which we assume
  that $\arrow \forget:\mathcal{F}_{\refl'}->\mathcal{E}.$ is a
  partial right adjoint to $\arrow \ladj:\mathcal{E}->\mathcal{F}.$:
  \begin{enumerate}
  \item\label{partial.adj.a} the image of $\forget$ is contained in
    $\mathcal{E}_\refl$ if and only if $\ladj$ carries all L-invertible
    arrows in $\mathcal{E}$ to $\refl'$-invertible arrows in
    $\mathcal{F}$.
  \item\label{partial.adj.b} if $S$ is a set of L-invertible arrows
    which is adequate to detect objects of $\mathcal{E}_\refl$ and $\ladj$
    maps all elements of $S$ to $\refl'$-invertible arrows in
    $\mathcal{F}$ then the image of $\forget$ is contained in
    $\mathcal{E}_\refl$.
  \item\label{partial.adj.c} if the image of $\forget$ is contained in
    $\mathcal{E}_\refl$ then it is (genuinely) right adjoint to the
    composite $\arrow \refl'\ladj:\mathcal{E}_\refl->\mathcal{F}_{\refl'}.$.
  \item\label{partial.adj.d} under the same conditions the functor $\ladj$
    preserves almost colimits, in the sense that if $\arrow
    \dgm:\mathbb{C}->\mathcal{E}.$ is a diagram and the cocone $\nattrans
    i:\dgm->X.$ is L-almost colimiting in $\mathcal{E}$ then the cocone
    $\nattrans \ladj i:\ladj\dgm->\ladj X.$ is $\refl'$-almost colimiting in $\mathcal{F}$.
  \end{enumerate}
  The proofs of these propositions are a matter of routine abstract
  category theory and are left up to the reader. When speaking of
  partial right adjoints we will distinguish those that satisfy the
  above property simply by notationally restricting their codomains.
  In other words, we will simply say that $\arrow
  \forget:\mathcal{F}_{\refl'}->\mathcal{E}_\refl.$ is a partial right
  adjoint to $\arrow \ladj:\mathcal{E}->\mathcal{F}.$.
\end{obs}

We'll also make heavy use of the following, slightly non-standard,
version of Day's reflection theorem~\cite{Day:1970:ClosedFunc} for
monoidal biclosed categories:

\begin{thm}[Day's reflection theorem]\label{day.refl} 
  Suppose that $\triple<\mathcal{E};\otimes;I>$ is a monoidal category
  and that $\mathcal{E}_\refl$ is a reflective full subcategory of
  $\mathcal{E}$ with associated idempotent monad $\pair<\refl;\eta>$.
  Also assume that for objects $X,Y\in\mathcal{E}$ the right and
  left tensoring functors $\arrow {-}\otimes Y:
  \mathcal{E}->\mathcal{E}.$ and $\arrow X\otimes{-}:
  \mathcal{E}->\mathcal{E}.$ have partial right adjoints $\arrow
  {*}\Leftarrow Y:\mathcal{E}_\refl->\mathcal{E}_\refl.$ and $\arrow
  X\Rightarrow{*}:\mathcal{E}_\refl->\mathcal{E}_\refl.$ respectively.
  
  As usual, we may enrich the family of functors ${*}\Leftarrow Y$
  (resp.\ $X\Rightarrow {*}$) with a canonical functorial structure in
  the variable $Y$ (resp.\ $X$), thereby forming them into a
  bi-functor $\arrow\Leftarrow: \mathcal{E}_\refl\times\op{\mathcal{E}}->
  \mathcal{E}_\refl.$ (resp.\ $\arrow\Rightarrow:\op{\mathcal{E}}\times
  \mathcal{E}_\refl-> \mathcal{E}_\refl.$).
  
  Under these assumptions, the following important results hold:
  \begin{enumerate}
  \item\label{day.refl.1} We may reflect the monoidal structure
    $\triple<\mathcal{E};\otimes;I>$ onto $\mathcal{E}_\refl$, making it
    into a monoidal biclosed category with:
    \begin{itemize}
    \item tensor product $A\otimes_\refl B\defeq \refl(A\otimes B)$ and
    \item left and right closures $A\Rightarrow{*}$ and 
      ${*}\Leftarrow B$.
    \end{itemize}
  \item\label{day.refl.2} An arrow $\arrow f:X->Y.$ in $\mathcal{E}$
    is L-invertible if and only if $\arrow A\Leftarrow f:A\Leftarrow
    Y-> A\Leftarrow X.$ (or dually $\arrow f\Rightarrow A:Y\Rightarrow
    A-> X\Rightarrow A.$) is an isomorphism for each
    $A\in\mathcal{E}_\refl$.
  \item\label{day.refl.3} The right and left tensoring functors
    $\arrow {-}\otimes Y, X\otimes{-}:\mathcal{E}->\mathcal{E}.$ both
    preserve L-almost colimits.
  \item\label{day.refl.4} If $\arrow \dgm:\mathbb{C}->\mathcal{E}.$ is a
    diagram then a cocone $\nattrans i:\dgm->X.$ is L-almost colimiting
    if and only if for each $A\in\mathcal{E}_\refl$ the contravariant
    functor $\arrow A\Leftarrow{*}:\op{\mathcal{E}}->\mathcal{E}_\refl.$
    (or dually $\arrow {*}\Rightarrow
    A:\op{\mathcal{E}}->\mathcal{E}_\refl.$) maps it to a limit cone in
    $\mathcal{E}_\refl$.
  \end{enumerate}
\end{thm}

\begin{proof}
  Our proof of (\ref{day.refl.1}) follows Day's original
  argument of \cite{Day:1970:ClosedFunc}, which revolves around the
  observation that if $X$ and $Y$ are objects in $\mathcal{E}$ then
  the arrows $\arrow \eta_X\otimes Y: X\otimes Y->\refl(X)\otimes Y.$,
  $\arrow X\otimes \eta_Y: X\otimes Y->X\otimes \refl(Y).$ and $\arrow
  \eta_X\otimes\eta_Y:X\otimes Y->\refl(X)\otimes\refl(Y).$ are all
  L-invertible. This follows from the observation that each component
  of the unit $\eta$ is L-invertible (cf.\ 
  observation~\ref{l-inv}(\ref{l-inv.unit})), that the left and right
  tensoring functors $X\otimes{-}$ and ${-}\otimes Y$ both preserve
  L-invertibility (by
  observation~\ref{partial.adj}(\ref{partial.adj.a}) and the
  assumption that they have appropriate partial right adjoints) and
  that L-invertible arrows are closed under composition.
  
  This immediately allows us to enrich the bi-functor $\otimes_\refl$
  with a monoidal structure which makes $\refl$ into a strong monoidal
  functor $\arrow \refl:\triple<\mathcal{E};\otimes;I>->
  \triple<\mathcal{E}_\refl;\otimes_\refl;\refl(I)>.$.  For instance,
  we may then construct an associativity isomorphism for
  $\otimes_\refl$ as the composite
  \begin{displaymath}
    \xymatrix@R=0.4ex@C=0.75em{
      {A\otimes_\refl(B\otimes_\refl C)}\ar@{=}[r] &
      {\refl(A\otimes \refl(B \otimes C))}
      \ar[rr]|\cong^{\refl(A\otimes\eta_{(B\otimes C)})^{-1}} &
      {\kern6em} &
      {\refl(A\otimes(B\otimes C))} \\
      & {\phantom{\refl(A\otimes \refl(B \otimes C))}}
      \ar[rr]^{\refl(\alpha_{A,B,C})}|\cong &&
      {\refl((A\otimes B)\otimes C)} \\
      &  {\phantom{\refl(A\otimes \refl(B \otimes C))}}
      \ar[rr]|{\cong}^{\refl(\eta_{(A\otimes B)}\otimes C)} &&
      {\refl(\refl(A\otimes B)\otimes C)} 
      \ar@{=}[r] &
      {(A\otimes_\refl B)\otimes_\refl C}
    } 
  \end{displaymath}
  in which the isomorphisms in the first and last line are a
  consequence of the L-invertibility result of the previous paragraph.
  To enrich $\arrow\refl:\mathcal{E}->\mathcal{E}_\refl.$ to a strong
  monoidal functor we have the isomorphism $\arrow \refl(\eta_X
  \otimes \eta_Y):\refl(X\otimes Y)->\refl(\refl(X)\otimes\refl(Y))
  \defeq \refl(X)\otimes_\refl\refl(Y).$ obtained by applying $\refl$
  to the L-invertible arrow $\eta_X\otimes\eta_Y$.  Furthermore, the
  fact that the bi-functors $\Leftarrow$ and $\Rightarrow$ provide us
  with right and left closures for this structure follows directly
  from observation~\ref{partial.adj}(\ref{partial.adj.c}).
  
  To prove (\ref{day.refl.2}) observe that if $A\in\mathcal{E}_\refl$ then
  $A\Leftarrow f$ is defined to be the unique arrow making the diagram
  \begin{displaymath}
    \let\labelstyle=\textstyle
    \xymatrix@R=4ex@C=6em{
      {\mathcal{E}(Z\otimes Y, A)}
      \ar[r]^{\mathcal{E}(Z\otimes f, A)}
      \ar[d]_{\cong} &
      {\mathcal{E}(Z\otimes X, A)}
      \ar[d]^{\cong} \\
      {\mathcal{E}(Z, A\Leftarrow Y)}
      \ar[r]_{\mathcal{E}(Z,A\Leftarrow f)} &
      {\mathcal{E}(Z, A\Leftarrow X)}
    }
  \end{displaymath}
  commute for each $Z\in\mathcal{E}$. By Yoneda's lemma, $A\Leftarrow
  f$ is an isomorphism iff the map at the bottom of this diagram is an
  isomorphism for each $Z$ and, since the verticals in the
  diagram are isomorphisms, it follows that this happens iff the map
  at the top of the diagram is an isomorphism. However, by definition,
  this latter condition simply states that $Z\otimes f$ is orthogonal
  to $A$ for all $Z\in\mathcal{E}$. 
  
  So if $A\Leftarrow f$ is an isomorphism for each $A\in\mathcal{E}_\refl$
  then it follows, since $f$ is isomorphic to $I\otimes f$, that $f$
  is orthogonal to each $A\in\mathcal{E}_\refl$ and thus that $f$ is
  L-invertible.  Conversely, if $f$ is L-invertible then we know
  that $Z\otimes{-}$ preserves L-invertibility, so $Z\otimes f$
  is also L-invertible for each $Z\in\mathcal{E}$ and thus
  orthogonal to any $A\in\mathcal{E}_\refl$, thereby demonstrating that
  $A\Leftarrow f$ is an isomorphism for each $A\in\mathcal{E}_\refl$ by
  the argument of the last paragraph.
  
  Part (\ref{day.refl.3}) is simply a restatement of
  observation~\ref{partial.adj}(\ref{partial.adj.d}) so all that
  remains is to prove part (\ref{day.refl.4}). To do so we simply
  apply the reflector $\refl$, and the natural family of isomorphisms
  $\iso A\Leftarrow\eta_X:A\Leftarrow \refl(X)->A\Leftarrow X.$
  arising from part (\ref{day.refl.2}) of this theorem (and the fact
  that each $\eta_X$ is L-invertible), to reflect the desired result
  regarding L-almost colimits into a corresponding characterisation of
  (genuine) colimits in the biclosed monoidal category
  $\mathcal{E}_\refl$. This latter result follows by an elementary
  argument in the theory of biclosed monoidal categories, which we
  leave up to the reader.
\end{proof}

\subsection{LFP Categories and LE-Theories}\label{lfp.le.subsect}

In this subsection we simply summarise those parts of the theory of
locally finitely presentable categories and left exact theories that
we apply in the subsequent text. For greater detail we refer the
reader to any number of standard presentations of this material,
including Gabriel and Ulmer's book~\cite{GabrielUlmer:1971:LPK},
Kennison's work on limit preserving functors~\cite{Kennison:1968:LPF},
Barr and Wells' textbook~\cite{BarrWells:1985:TTT} and Kelly's book on
enriched category theory~\cite{Kelly:1982:ECT}.

\begin{defn}[weighted limits and colimits]
  \label{weighted.lims}
  Recall, from Kelly~\cite{Kelly:1982:ECT}, that the colimit of a
  diagram $\arrow \dgm:\mathbb{C}->\mathcal{E}.$ {\em weighted\/} or
  {\em indexed\/} by $\arrow W:\op{\mathbb{C}}->\Set.$ is an object
  $\colim(W,\dgm)$ of $\mathcal{E}$ which represents the functor
  $\functor{C}_{W,\dgm}(E)=\funcat[\op{\mathbb{C}},\Set](W,
  \mathcal{E}(\dgm({*}),E))$ from $\mathcal{E}$ to $\Set$. The unit of
  this representation is a natural transformation $\arrow\iota:W->
  \mathcal{E}(\dgm({*}),\colim(W,\dgm)).$ called its {\em colimiting
    cylinder}. In this context, the traditional colimit of $\dgm$ is
  referred to as its {\em conical colimit\/} and it may be obtained by
  weighting our colimit by the functor which maps each object of
  $\mathbb{C}$ to a singleton set.
  
  This weighted colimit construction (wherever it exists) is
  covariantly functorial in each argument. To be precise, any natural
  transformation of weights $\arrow\gamma:W->V.\in\funcat[
  \op{\mathbb{C}},\Set]$ induces an unique arrow
  $\arrow\colim(\gamma,\dgm):\colim(W,\dgm)-> \colim(V,\dgm).$ in
  $\mathcal{E}$ which represents the family of pre-composition
  functions $\arrow{-}\circ\gamma:\funcat[ \op{\mathbb{C}},\Set](
  V,\mathcal{E}( \dgm({*}),E)) ->\funcat[\op{\mathbb{C}},\Set](
  W,\mathcal{E}(\dgm({*}),E)).$.  Similarly, if
  $\arrow\phi:\dgm->\dgm'.$ is a natural transformation of diagrams
  then the induced arrow $\arrow\colim(W,\phi):\colim(W, \dgm)->
  \colim(W, \dgm').$ represents the family of post-compositions
  $\arrow \mathcal{E}(\phi_{{*}},E)\circ{-}:
  \funcat[\op{\mathbb{C}},\Set](W,\mathcal{E}(\dgm({*}),E))->
  \funcat[\op{\mathbb{C}},\Set](W,\mathcal{E}(\dgm'({*}),E)).$.
  
  Dually the limit $\lim(U,\dgm)\in\mathcal{E}$ of our diagram
  weighted by $\arrow U:\mathbb{C}->\Set.$ represents the functor
  $\functor{L}_{U,\dgm}(E)=\funcat[\mathbb{C},
  \Set](U,\mathcal{E}(E,\dgm({*})))$. It is worth remarking that under
  the duality involved here $\lim(U,\dgm)$ becomes contravariantly
  functorial in its first variable and covariantly so in its second.
  
  Let $\arrow\Yoneda:\mathbb{C}->\funcat[\op{\mathbb{C}},\Set].$
  denote the usual Yoneda functor, which maps an object
  $c\in\mathbb{C}$ to the representable $\Yoneda_c=\mathbb{C}({-},c)$,
  and observe that Yoneda's lemma immediately implies that
  $\colim(W,\Yoneda)\cong W$ and that $\colim(\Yoneda_c,\dgm)\cong
  \dgm(c)$.  It is also worth recalling that from the weight $W$ we
  may construct a small comma category $\groth(W)=\Yoneda
  {\downarrow}W$, called the {\em Grothendieck category\/} of $W$. By
  Yoneda's lemma we can describe $\groth(W)$ as a category whose
  objects are pairs $\pair<c;x>$, where $c$ is an object of
  $\mathbb{C}$ and $x\in W(c)$, and in which an arrow $\arrow
  f:\pair<c;x>->\pair<c';x'>.$ is simply an arrow $\arrow f:c->c'.$ of
  $\mathbb{C}$ such that $W(f')(x')=x$.  Composing the obvious
  projection functor $\arrow p:\groth(W)-> \mathbb{C}.$ with the
  diagram $\arrow\dgm:\mathbb{C}->\mathcal{E}.$ and taking the conical
  colimit of this, if it exists, we obtain an object which represents
  the functor $\functor{C}_{W,\dgm}$ of the last paragraph.  In other
  words the weighted limit $\colim(W,\dgm)$ exists if and only if the
  conical colimit $\colim_{\groth(W)}(\dgm\circ p)$ does and in that
  case these colimits are isomorphic.
\end{defn}

\begin{defn}[filtered colimits]
  \label{filtered.colim}
  We say that a small, non-empty category $\mathbb{D}$ is {\em
    filtered\/} iff
  \begin{itemize}
  \item given any two objects $d,d'\in\mathbb{D}$ there is a third
    object $\bar{d}\in\mathbb{D}$ to which they both map, and
  \item given any parallel pair of arrows $\arrow f,g:d->d'.$ in
    $\mathbb{D}$ there exists a $\arrow h:d'->\bar{d}.$ with
    $h\circ f= h\circ g$.
  \end{itemize}
  A {\em filtered colimit\/} is a (conical) colimit on a diagram whose
  indexing category is filtered.
  
  Most importantly, in the category $\Set$ filtered colimits bear a
  simple description which allows us to easily demonstrate that they
  commute with all finite limits.
  
  We say that a functor $\arrow \ladj:\mathcal{E}->\mathcal{F}.$ is {\em
    finitely accessible\/} iff it preserves all the (small) filtered
  colimits which exist in $\mathcal{E}$.
\end{defn}
 
\begin{defn}[locally finitely presentable categories]
  \label{lfp-cat.defn}
  We say that a category $\mathcal{E}$ is {\em locally finitely
    presentable\/} (or simply an {\em LFP-category}) iff it is
  (equivalent to a) reflective full subcategory of some functor
  category $\funcat[\op{\mathbb{C}},\Set]$ for which $\mathbb{C}$ is a
  small category and the reflector $\refl$ associated with $\mathcal{E}$
  is finitely accessible.
  
  In general, we say that an object $E$ in some category $\mathcal{E}$
  is {\em finitely presentable\/} (or simply {\em FP}) iff the
  associated representable functor $\arrow \mathcal{E}\pair<A;{-}>:
  \mathcal{E}->\Set.$ preserves all (small) filtered colimits. Using
  this concept, we may provide a simple ``intrinsic'' characterisation
  of LFP-categories as being those locally small and small cocomplete
  categories $\mathcal{E}$ which possess a small dense set of finitely
  presentable objects.
\end{defn}

\begin{obs}[reflective full sub-categories of LFP-categories]
  \label{full.refl.LFP}
  If $\pair<\refl;\eta>$ is an idempotent monad on an LFP-category
  $\mathcal{E}$ for which $\refl$ is finitely accessible then it is a
  routine matter to demonstrate that the corresponding reflective full
  subcategory $\mathcal{E}_\refl$ is also an LFP-category.
  
  In this context, it is quite common to specify such reflective full
  sub-categories by supplying a defining {\em FP-regulus\/} $T$, which
  is simply a (essentially small) set of arrows in $\mathcal{E}$ all
  of whose domains and codomains are FP-objects.  Applying
  theorem~10.2 of~\cite{Kelly:1980:Trans}, we see that the full
  subcategory $\mathcal{E}_{\perp T}$ of objects orthogonal to $T$ is
  reflective in $\mathcal{E}$. Furthermore, an immediate consequence
  of our assumption of finitely presentable domains and codomains is
  that $\mathcal{E}_{\perp T}$ is closed in $\mathcal{E}$ under
  filtered colimits. It follows that the reflector associated with
  $\mathcal{E}_{\perp T}$ must preserve filtered colimits and that we
  may apply the result of the previous paragraph to demonstrate that
  this subcategory is also LFP.
\end{obs}

\begin{defn}[left exact theories]
  \label{le-th.defn}
  LFP-categories may be presented explicitly in terms of a kind of
  categorical theory $\mathbb{T}=\pair<\mathbb{C};T>$ called a {\em
    left exact theory\/} (or simply an {\em LE-theory}) which consists
  of a small category $\mathbb{C}$ and an FP-regulus $T$ in the
  presheaf category $\funcat[\op{\mathbb{C}},\Set]$. We use the
  notation $\Al{\mathbb{T}}$ to denote the reflective full-subcategory
  $\funcat[\op{\mathbb{C}},\Set]_{\perp T}$ of objects orthogonal to
  the arrows in the FP-regulus $T$, call its objects  {\em
  $\mathbb{T}$-algebras\/} and use the notation $\refl_{\mathbb{T}}$
  for the corresponding finitely accessible reflector. 
  
  We may generalise the $\mathbb{T}$-algebra notion representably to
  any category $\mathcal{E}$, by defining such algebras to be those
  contravariant functors $\arrow A:\op{\mathbb{C}}->\mathcal{E}.$
  satisfying the condition that for each object $E\in\mathcal{E}$ the
  functor $\mathcal{E}(E,A(*))\in\funcat[\op{\mathbb{C}},\Set]$ is
  actually a $\mathbb{T}$-algebra.  We use the notation
  $\Alg<\mathbb{T};\mathcal{E}>$ to denote the full subcategory of
  the functor category $\funcat[\op{\mathbb{C}},\mathcal{E}]$ of those
  functors which are $\mathbb{T}$-algebras in $\mathbb{C}$, and refer
  to the arrows of $\Alg<\mathbb{T};\mathcal{E}>$ as {\em
    $\mathbb{T}$-algebra morphisms}.
  
  However, we are usually only interested in considering
  $\mathbb{T}$-algebras in categories $\mathcal{E}$ which possess all
  finite limits. These are called {\em LE-categories} and functors
  between them which preserve all finite limits are called {\em
    LE-functors}. Notice that when $W$ is a finitely presented object
  of $\funcat[\op{\mathbb{C}},\Set]$ we may express it as a finite
  colimit of representables and consequently show that $\lim(W,A)$ may
  be constructed as a finite limit in $\mathcal{E}$.  So if
  $\mathcal{E}$ is an LE-category then for each $\arrow\gamma:W->V.$
  in the regulus $T$ we may form the arrow $\arrow\lim(\gamma,A):
  \lim(V,A)-> \lim(W,A).$ in $\mathcal{E}$ representing the
  pre-composition $\arrow {-}\circ\gamma:\funcat[\op{\mathbb{C}},
  \Set](V,\mathcal{E}(E,A(*))->\funcat[\op{\mathbb{C}},
  \Set](W,\mathcal{E}(E,A(*)).$ for each object $E\in\mathcal{E}$
  and therefore we see, by Yoneda's lemma, that $\mathcal{E}(E,A(*))$
  is orthogonal to $\gamma$ for all objects $E\in\mathcal{E}$ if and
  only if $\lim(\gamma,A)$ is an isomorphism. It follows that in the
  context of LE-categories the $\mathbb{T}$-algebra property may be
  re-expressed as stating that the arrow
  $\arrow\lim(\gamma,A):\lim(V,A)->\lim(W,A).$ is an isomorphism for
  each arrow $\arrow \gamma:W->V.$ in the FP-regulus $T$.
 
  Having expressed our $\mathbb{T}$-algebra property in terms of
  finite limits we immediately see that if $\arrow
  \ladj:\mathcal{E}->\mathcal{F}.$ is an LE-functor then the
  post-composition functor $\arrow\ladj\circ{-}:
  \funcat[\op{\mathbb{C}},\mathcal{E}] -> \funcat[\op{\mathbb{C}},
  \mathcal{F}].$ carries $\mathbb{T}$-algebras in $\mathcal{E}$ to
  $\mathbb{T}$-algebras in $\mathcal{F}$. Therefore it restricts to a
  functor from $\Alg<\mathbb{T};\mathcal{E}>$ to
  $\Alg<\mathbb{T};\mathcal{F}>$ which we will sometimes call
  $\Alg<\mathbb{T};\ladj>$, although in many cases it is convenient to
  overload our notation and simply think of $\ladj$ itself as having
  been extended or lifted to a functor $\arrow
  \ladj:\Alg<\mathbb{T};\mathcal{E}>
  ->\Alg<\mathbb{T};\mathcal{F}>.$. We often say that
  $\Alg<\mathbb{T};\ladj>$ is obtained by {\em applying $\ladj$
    point-wise to the $\mathbb{T}$-algebras in its domain}.
\end{defn}  

\begin{obs}[more about LE-theories]
  Sometimes it is more convenient to present LE-theories in terms of
  {\em finite sketches\/} rather than FP-reguli. A finite sketch on a
  small category $\mathbb{C}$ consists of a family of diagrams $\arrow
  \dgm_i:\mathbb{D}_i->\mathbb{C}.$ ($i\in I$), for which each
  $\mathbb{D}_i$ is a finite category, and a family of cocones
  $\nattrans\iota_i:\dgm_i->c_i.$ ($i\in I$). We say that a functor
  $\arrow A:\op{\mathbb{C}}->\Set.$ is an algebra for such a sketch if
  it carries each of the specified cocones to a limiting cone in
  $\Set$.
  
  From a sketch we may construct an FP-regulus whose members are
  obtained by applying the Yoneda functor $\arrow \Yoneda:\mathbb{C}->
  \funcat[\op{\mathbb{C}}, \Set].$ to each of its diagrams and
  cocones, taking the colimit of the former to obtain finitely
  presented objects $\colim_{\mathbb{D}_i}(\Yoneda\circ \dgm_i)\in
  \funcat[\op{\mathbb{C}},\Set]$ and forming the arrows $\arrow
  f_i:\colim(\Yoneda\circ\dgm_i)->\Yoneda_{c_i}.$ induced by the
  latter. A trivial application of Yoneda's lemma now demonstrates
  that the algebras for our sketch and its associated FP-regulus
  coincide.
  
  These definitions are motivated by the observation that a
  $\mathbb{T}$-algebra in $\Set$ for an LE-theory is really
  no-more-nor-less than a family of basic underlying sets equipped
  with a family of operations whose domains and codomains are finite
  limits of those basic sets and for which a specified class of
  finitary equational conditions must hold. In particular these
  operations and equational conditions are coded into the structure of
  $\mathbb{C}$ while the FP-regulus $T$ (or sketch) specifies
  precisely which limits constitute the domains and codomains of each
  operation. We often think of each of these operations as being
  partially defined upon equationally specified subsets of products of
  basic sets.
  
  We refer the reader to~\cite{BarrWells:1985:TTT} for a thorough
  exposition of theory of LE-theories, including a proof of Kennison's
  theorem which demonstrates that $\Alg<\mathbb{T};\mathcal{E}>$ is a
  reflective full subcategory of $\funcat[\op{\mathbb{C}},
  \mathcal{E}]$ (for all sufficiently cocomplete LE-categories
  $\mathcal{E}$). It is also straightforward to show that
  $\Alg<\mathbb{T};\mathcal{E}>$ is closed in
  $\funcat[\op{\mathbb{C}},\mathcal{E}]$ under all the limits and
  filtered colimits which exist in there.  Combining these two
  results, it follows that $\Alg<\mathbb{T};\mathcal{E}>$ is an
  LFP-category for any LFP-category $\mathcal{E}$.
\end{obs}

\begin{obs}[tensor products of LE-theories]
  \label{le-th.tensor}
  Thus far we know that if $\mathbb{T}$ and $\mathbb{S}=\pair<
  \mathbb{D};S>$ are LE-theories and $\mathcal{E}$ is an LFP-category
  then the category $\Alg<\mathbb{S};{\Alg<\mathbb{T}; \mathcal{E}>}>$
  is also LFP and may, therefore, be presented as a category of
  algebras for some LE-theory in $\mathcal{E}$. However, for many
  purposes it is useful to have an explicit presentation of this
  theory in terms of $\mathbb{T}$ and $\mathbb{S}$ themselves. To this
  end we define their {\em tensor product\/}
  $\mathbb{T}\otimes\mathbb{S}$ to be the LE-theory with underlying
  category $\mathbb{C}\times\mathbb{D}$ and FP-regulus $T\otimes S$
  in $\funcat[\op{(\mathbb{C}\times \mathbb{D})},\Set]$ given by
  \begin{equation*}
    \begin{split}
      T\otimes S = {} &\{ \arrow \gamma\circ\op\pi_{\mathbb{C}}:
      W\circ\op\pi_{\mathbb{C}}->V\circ\op\pi_{\mathbb{C}}.\mid
      \arrow\gamma:W->V.\in T\}\cup{}\\
      &\{ \arrow \gamma'\circ\op\pi_{\mathbb{D}}:
      W'\circ\op\pi_{\mathbb{D}}->V'\circ\op\pi_{\mathbb{D}}.\mid
      \arrow\gamma':W'->V'.\in S\} \\
    \end{split}
  \end{equation*}
  where $\arrow\pi_{\mathbb{C}}:\mathbb{C}\times\mathbb{D}->
  \mathbb{C}.$ and $\arrow\pi_{\mathbb{D}}:\mathbb{C}\times\mathbb{D}->
  \mathbb{D}.$ are the canonical projection functors.
  
  This definition is designed to ensure that a
  $\mathbb{T}\otimes\mathbb{S}$-algebra in an LE-category
  $\mathcal{E}$ is precisely a functor $\arrow
  A:\mathbb{C}\times\mathbb{D}->\mathcal{E}.$ satisfying the
  condition that for each object $c\in\mathbb{C}$ the functor $\arrow
  A\pair<e;{-}>:\mathbb{D}->\mathcal{E}.$ is a $\mathbb{S}$-algebra
  and for each object $e'\in\mathbb{D}$ the functor $\arrow
  A\pair<{-};d>: \mathbb{C}->\mathcal{E}.$ is a $\mathbb{T}$-algebra.
  Furthermore, it is clear that the canonical isomorphisms
  \begin{displaymath}
    \funcat[\op{\mathbb{C}},{\funcat[\op{\mathbb{D}},\mathcal{E}]}]\cong
    \funcat[\op{(\mathbb{C}\times\mathbb{D})},\mathcal{E}]]\cong
    \funcat[\op{\mathbb{D}},{\funcat[\op{\mathbb{C}},\mathcal{E}]}]
  \end{displaymath}
  of functor categories restrict to isomorphisms
  \begin{displaymath}
    \Alg<\mathbb{T};{\Alg<\mathbb{S};\mathcal{E}>}>\cong
    \Alg<\mathbb{T}\otimes\mathbb{S};\mathcal{E}>\cong
    \Alg<\mathbb{S};{\Alg<\mathbb{T};\mathcal{E}>}>
  \end{displaymath}
  of categories of algebras which are natural in $\mathcal{E}$. 
\end{obs} 

\begin{obs}[Kan's construction - right adjoint functors from L-almost 
  coalgebras]\label{func.from.coalg} 
  Of course, we say that a functor $\arrow C:\mathbb{C}->
  \mathcal{E}.$ is a {\em $\mathbb{T}$-coalgebra\/} if its dual
  is a $\mathbb{T}$-algebra in $\op{\mathcal{E}}$. Taking
  the dual of our characterisation of $\mathbb{T}$-algebras in
  LE-categories, we see that we can characterise
  $\mathbb{T}$-coalgebras as being those functors $\arrow C:\mathbb{C}
  ->\mathcal{E}.$ such that for each arrow $\arrow\gamma:W->V.$ in
  the FP-regulus $T$ the induced arrow $\arrow\colim(\gamma,C):
  \colim(W,C)->\colim(V,C).$ is an isomorphism.
  
  Our primary interest in $\mathbb{T}$-coalgebras here is that each
  one gives rise to a right adjoint functor into the category of
  $\mathbb{T}$-algebras, by an application of (a variant of) {\em
    Kan's construction\/} \cite{Kan:1958:Adjoint}.  However, in the
  sequel we will actually be interested in taking the approach
  described in observation~\ref{l-almost}, that is we will work in a
  category $\mathcal{E}$ equipped with an idempotent monad
  $\pair<\refl;\eta>$ and use L-almost coalgebras to construct right
  adjoint functors with domain the reflective full subcategory
  $\mathcal{E}_\refl$.

  Following the convention established in observation~\ref{l-almost},
  we now define an L-almost $\mathbb{T}$-coalgebra in $\mathcal{E}$ to
  be a functor $\arrow C:\mathbb{C}->\mathcal{E}.$ for which the
  composite $\arrow \refl\circ C:\mathbb{C}->\mathcal{E}_\refl.$ is a
  genuine $\mathbb{T}$-coalgebra in $\mathcal{E}_\refl$. Using our
  characterisation of coalgebras in terms of colimits in
  $\mathcal{E}_\refl$ it follows that this is equivalent to saying
  that for each $\arrow \gamma:W->V.$ in the FP-regulus $T$ the
  induced arrow $\arrow \colim(\gamma,C): \colim(W,C)->\colim(V,C).$
  in $\mathcal{E}$ is L-invertible. We use the notation
  $\CoAl{\mathbb{T}}_\refl(\mathcal{E})$ for the full subcategory of
  L-almost $\mathbb{T}$-coalgebras in $\funcat[\mathbb{C},
  \mathcal{E}]$. Again it follows that if the functor
  $\arrow\ffunc:\mathcal{E}-> \mathcal{F}.$ preserves L-almost finite
  colimits then it may be applied point-wise to provide a functor
  $\arrow\CoAl{\mathbb{T}}_\refl(\ffunc): \CoAl{ \mathbb{T}
  }_\refl(\mathcal{E})->\CoAl{\mathbb{T}}_\refl( \mathcal{F}).$.
  
  Any L-almost $\mathbb{T}$-coalgebra $C$ gives rise to an adjoint
  pair
  \begin{equation*}
    \let\labelstyle=\textstyle
    \xymatrix@R=1ex@C=14em{
      {\mathcal{E}_\refl}\ar[r]^{\bot}_{\homout{C}} &
      {\Al{\mathbb{T}}}\ar@/_3ex/[l]_{\kanladj{C}}
    }
  \end{equation*}
  in which $\homout{C}$ is the {\em Kan functor\/} constructed by
  ``homming out'' of $C$. In other words, if $E$ is an object of
  $\mathcal{E}_\refl$ then we define
  $\homout{C}(E)=\arrow\mathcal{E}(C(-),E):\op{\mathbb{C}} ->\Set.$
  and if $\arrow f:E->E'.$ is an arrow of $\mathcal{E}$ then
  $\homout{C}(f)$ is the natural transformation with components
  $\homout{C}(f)_c= \mathcal{E}(C(c),f)=f\circ{-}$ for
  $c\in\mathbb{C}$. That $\homout{C}(E)$ is actually a
  $\mathbb{T}$-algebra follows from the fact that we have the
  reflection isomorphism $\mathcal{E}(C({-}),E)\cong\mathcal{E}_\refl(
  \refl\circ C({-}),E)$ for each $E\in\mathcal{E}_\refl$ and by
  definition we know that $\refl\circ C$ is a $\mathbb{T}$-coalgebra
  in $\mathcal{E}_\refl$ iff the functor on the right of this
  isomorphism is a $\mathbb{T}$-algebra (in $\Set$) for each object
  $E\in\mathcal{E}_\refl$.  
  
  This construction is clearly contravariantly functorial in the
  L-almost coalgebra $C$ so we get a functor $\arrow\homout{*}:
  \op{\CoAl{\mathbb{T}}_\refl(\mathcal{E})} ->
  \funcat[\mathcal{E}_\refl, \Al{\mathbb{T}}].$. If a functor $\arrow
  \gfunc:\mathcal{E}_\refl->\Al{\mathbb{T}}.$ is isomorphic to the Kan
  functor on some L-almost coalgebra $C$ then we say that $\gfunc$ is
  {\em representable\/} or that it is {\em represented by\/} $C$.
  
  Each $\mathbb{T}$-algebra $A$ can, of course, be considered to be a
  weight on $\mathbb{C}$ and so we may define $\kanladj{C}(A)$ to be
  the colimit $\colim(A,\refl\circ C)$ in $\mathcal{E}_\refl$. That
  this does indeed provide a left adjoint to $\homout{C}$ follows from
  the calculation:
  \begin{equation*}
    \begin{aligned}
      \mathcal{E}_\refl(\colim(A,\refl\circ C), E) &
      {} \cong \funcat[\op{\mathbb{C}},\Set](A,\mathcal{E}_\refl( 
      \refl\circ C({-}),E) && \text{definition of weighted colimit,} \\
      & {} \cong \funcat[\op{\mathbb{C}},\Set](A,\mathcal{E}(C({-}),E))
      && \text{reflection isomorphism,} \\
      & {}\cong \Al{\mathbb{T}}(A,\kan{C}(E))
      && \text{definition of $\kan{C}$}
    \end{aligned}
  \end{equation*}
  Finally, we say that an L-almost coalgebra $C$ in $\mathcal{E}$ is
  {\em finitely presented\/} iff it maps each object of $\mathbb{C}$
  to a finitely presented object in $\mathcal{E}$. We will sometimes
  have use for the routine observation that if the reflector $\refl$
  is a finitely accessible functor and the L-almost coalgebra $C$ is
  finitely presented then $\homout{C}$ also preserves all filtered
  colimits. The proof of this result is a routine matter of elementary
  category theory, which we leave up to reader.
\end{obs}  

\begin{obs}[An internal version of Kan's construction]
  \label{internal.kan}
  For the remainder of this subsection, we'll consider a context in
  which $\mathcal{E}$ is an LFP-category that comes equipped with a
  finitely accessible idempotent monad $\pair<\refl;\eta>$ and the
  monoidal and partial closures described in our version of Day's
  reflection theorem (theorem~\ref{day.refl}). In such situations we
  will be interested in studying $\mathbb{T}$-algebras in the LFP
  biclosed reflective full subcategory $\mathcal{E}_\refl$ and
  constructing functors $\overarr:\mathcal{E}_\refl->\Alg<\mathbb{T};
  \mathcal{E}_\refl>.$ using L-almost $\mathbb{T}$-coalgebras in
  $\mathcal{E}$.
  
  Using the tensor product and partial closures on $\mathcal{E}$ we
  may internalise the Kan's construction to derive right adjoint
  functors from $\mathcal{E}_\refl$ to
  $\Alg<\mathbb{T};\mathcal{E}_\refl>$. To see that this is the case
  recall from theorem~\ref{day.refl} that for each object
  $E\in\mathcal{E}_\refl$ the contravariant left closure functor
  $\arrow{*}\Rightarrow E: \op{\mathcal{E}}-> \mathcal{E}_\refl.$
  inverts L-invertible arrows and carries L-almost colimits in
  $\mathcal{E}$ to limits in $\mathcal{E}_\refl$. It follows
  therefore, from our colimit and limit characterisations, that
  post-composition by ${*}\Rightarrow E$ carries L-almost
  $\mathbb{T}$-coalgebras in $\mathcal{E}$ to $\mathbb{T}$-algebras in
  $\mathcal{E}_\refl$ and thus provides a functor $\overarr:\op{
    \CoAl{\mathbb{T}}_\refl(\mathcal{E})}->\Al{\mathbb{T}}(
  \mathcal{E}_\refl).$. The functoriality of this construction in
  $E\in\mathcal{E}_\refl$ is clear and we have therefore succeeded in
  constructing a functor $\arrow\ihoml{*}:\op{\CoAl{\mathbb{T}
    }_\refl(\mathcal{E})} ->\funcat[\mathcal{E}_\refl,\Al{
    \mathbb{T}}(\mathcal{E}_\refl)].$ which maps an L-almost
  $\mathbb{T}$-coalgebra $\arrow C:\mathbb{C}->\mathcal{E}.$ to the
  Kan functor $\arrow \ihoml{C}:\mathcal{E}_\refl->\Alg<\mathbb{T};
  \mathcal{E}_\refl>.$ obtained by ``internally homming out of $C$''
  or explicitly:
  \begin{equation*}
    \begin{aligned}
      \ihoml{C}(E)&{}\defeq C({-})\Rightarrow E &&
      \text{for each object $E\in\mathcal{E}_\refl$, and} \\
      \ihoml{C}(f)_c&{}\defeq C(c)\Rightarrow f && \text{for each
        arrow $f\in\mathcal{E}_\refl$ and object $c\in\mathbb{C}$.}
    \end{aligned}
  \end{equation*}
  Furthermore, if $\arrow \alpha:C->C'.$ is an L-almost
  $\mathbb{T}$-coalgebra map then it is mapped to the natural
  transformation $\arrow\ihoml{\alpha}:\ihoml{C'}->\ihoml{C}.$ whose
  component at the object $E\in\mathcal{E}_\refl$ is the
  $\mathbb{T}$-algebra morphism $(\ihoml{\alpha})_E\defeq\Alg<
  \mathbb{T};{*}\Rightarrow E>(\alpha)$ which consists of the family
  of arrows $\arrow\alpha_c\Rightarrow E:C'(c)\Rightarrow E->
  C(c)\Rightarrow E.$ indexed by objects $c\in\mathbb{C}$.
  
  If $A$ is a $\mathbb{T}$-algebra in $\mathcal{E}_\refl$ then we may
  tensor it with the L-almost $\mathbb{T}$-coalgebra and compose the
  result with $\refl$ to give a functor $\arrow\refl(A({-})\otimes
  C({*})):\op{\mathbb{C}}\times\mathbb{C}->\mathcal{E}_\refl.$. Now we
  may define $\kanladj{C}\lhv(A)$, the value of the left adjoint to
  $\ihoml{C}$ on the algebra $A$ to be a certain kind of colimit of
  this diagram in $\mathcal{E}_\refl$ called its {\em coend\/} by Kelly
  in~\cite{Kelly:1982:ECT} and denoted by $\int^{c\in\mathbb{C}}(\refl(
  A(c)\otimes C(c)))\in\mathcal{E}_\refl$.
  
  Notice also that if we know that the tensor product of any pair of
  finitely presentable objects in $\mathcal{E}$ is always finitely
  presentable then we may show that $C(c)\Rightarrow {-}$ preserves
  all filtered colimits in $\mathcal{E}_\refl$ whenever $C(c)$ is
  finitely presentable. Furthermore, we also know that
  $\Alg<\mathbb{T}; \mathcal{E}_\refl>$ is closed in
  $\funcat[\mathbb{C}, \mathcal{E}_\refl]$ under filtered limits which
  are calculated point-wise in there. Combining these two facts, it
  follows therefore that the functor $\arrow \ihoml{C}:
  \mathcal{E}_\refl-> \Alg<\mathbb{T};\mathcal{E}_\refl>.$ preserves
  all filtered colimits whenever $C$ is a finitely presented L-almost
  coalgebra.

  Observe that if $I$ is the identity object of the monoidal structure
  on $\mathcal{E}$ then, since we generalised the $\mathbb{T}$-algebra
  notion to arbitrary categories representably, we know that
  post-composition by the representable functor
  $\arrow\mathcal{E}_\refl(\refl(I),-)\cong \mathcal{E}(I,-):
  \mathcal{E}_\refl->\Set.$ gives a functor $\arrow\Al{\mathbb{T}}(
  \mathcal{E}(I,{-})):\Al{\mathbb{T}}(\mathcal{E}_\refl)->
  \Al{\mathbb{T}}.$. By partial adjointness we also have natural
  isomorphisms $\mathcal{E}(I,C(-)\Rightarrow
  E)\cong\mathcal{E}(C(-)\otimes I,E)\cong\mathcal{E}(C(-), E)$ for
  each object $E\in\mathcal{E}_\refl$ and we may use these to
  construct an isomorphism $\mathcal{E}(I,{-})\circ
  \ihoml{C}\cong\homout{C}$ in
  $\funcat[\mathcal{E}_\refl,\Al{\mathbb{T}}]$ which is natural in the
  L-almost $\mathbb{T}$-coalgebra $C$.
  
  Finally, we could equally well have used the partial right closure
  functor $E\Leftarrow{*}$ in the calculations above and thereby
  obtain another functor $\arrow \ihomr{*}: \op{\CoAl{\mathbb{T}
    }_\refl(\mathcal{E})} ->\funcat[\mathcal{E}_\refl, \Alg<
  \mathbb{T};\mathcal{E}_\refl>].$.  In general, we will distinguish
  these constructions in our text by referring to the former as the
  {\em left handed\/} construction and the latter as the {\em right
    handed\/} version.  If our monoidal structure on $\mathcal{E}$ is
  augmented by a suitable duality operation
  $\arrow({-})^\circ:\mathcal{E}-> \mathcal{E}.$ then we may relate
  these constructions as (appropriately defined) duals of each other.
\end{obs}

\begin{obs}[tensor products of L-almost coalgebras]
  \label{coalg.tensor}
  Assume again that $\mathcal{E}$ is as in the last observation, and
  suppose that $C\in\CoAl{\mathbb{T}}_\refl(\mathcal{E})$ and
  $D\in\CoAl{\mathbb{S}}_\refl(\mathcal{E})$ then the functor $\arrow
  C\otimes D: \mathbb{C}\times\mathbb{D}->\mathcal{E}.$ defined by
  $(C\otimes D)\pair<c;d>\defeq C(e)\otimes D(d)$ is an L-almost
  $\mathbb{T} \otimes \mathbb{S}$-coalgebra called the {\em tensor
    product\/} of $C$ and $D$.
  
  To show that this is indeed an L-almost coalgebra consider an object
  $c\in\mathbb{C}$ and observe that $C(c)\otimes{-}$ preserves
  L-almost colimits, since it has a partial right adjoint
  $C(c)\Rightarrow{*}$, and thus may be applied point-wise to the
  L-almost $\mathbb{S}$-coalgebra $D$ to show that $(C\otimes
  D)\pair<c;{-}>= C(c)\otimes D({-})$ is also an L-almost
  $\mathbb{S}$-coalgebra.  By a dual argument $(C\otimes
  D)\pair<{-};d>=C({-})\otimes D(d)$ is an L-almost
  $\mathbb{T}$-coalgebra for each object $d\in\mathbb{D}$.  Consulting
  observation~\ref{le-th.tensor} it is clear that these two
  observations are enough to demonstrate that $C\otimes D$ is an
  L-almost $\mathbb{T}\otimes\mathbb{S}$-coalgebra as required.
  
  This definition is motivated by the observation that the monoidal
  and partial closure structures on $\mathcal{E}$ provide canonical
  isomorphisms
  \begin{align*}
    D(d)\Rightarrow(C(c)\Rightarrow E)&\mkern10mu\cong\mkern10mu
    (C(c)\otimes D(d))\Rightarrow E \\
    (E\Leftarrow C(c))\Leftarrow D(d)&\mkern10mu\cong\mkern10mu
    E\Leftarrow(D(d)\otimes C(c))
  \end{align*}
  which are natural in $E\in\mathcal{E}_\refl$, $c\in\mathbb{C}$ and
  $d\in\mathbb{D}$ and thus may be collected together to provide
  natural isomorphisms
  \begin{equation}\label{coalg.tensor.diag}
    \let\labelstyle=\textstyle
    \xymatrix@=2.25em{
      {\mathcal{E}_\refl}\ar[r]^<>(0.5){\ihoml{C}}
      \ar[d]_{\ihoml{C\otimes D}} &
      {\Alg<\mathbb{T};\mathcal{E}_\refl>}
      \ar[d]^{\Alg<\mathbb{T};\ihoml{D}>} \\
      {\Alg<\mathbb{T}\otimes\mathbb{S};\mathcal{E}_\refl>} 
      \ar[r]_<>(0.5){\cong} &
      {\Alg<\mathbb{T};{\Alg<\mathbb{S};\mathcal{E}_\refl>}>}
      \ar@{}[ul]|{\textstyle\cong}
    }
    \xymatrix@=2.25em{
      {\mathcal{E}_\refl}\ar[r]^<>(0.5){\ihomr{C}}
      \ar[d]_{\ihomr{D\otimes C}} &
      {\Alg<\mathbb{T};\mathcal{E}_\refl>}
      \ar[d]^{\Alg<\mathbb{T};\ihomr{D}>} \\
      {\Alg<\mathbb{S}\otimes\mathbb{T};\mathcal{E}_\refl>} 
      \ar[r]_<>(0.5){\cong} &
      {\Alg<\mathbb{T};{\Alg<\mathbb{S};\mathcal{E}_\refl>}>}
      \ar@{}[ul]|{\textstyle\cong}
    }
  \end{equation}
  respectively. Sometimes it is useful to partially externalise this
  result and use the partial adjunction isomorphisms
  \begin{equation*}
    \mathcal{E}_\refl(D(d),C(c)\Rightarrow E)\cong
    \mathcal{E}(C(c)\otimes D(d), E)\mkern40mu
    \mathcal{E}_\refl(D(d),E\Leftarrow C(c))\cong
    \mathcal{E}(D(d)\otimes C(c),E)
  \end{equation*}
  which are natural in $E\in\mathcal{E}_\refl$, $c\in\mathbb{C}$ and
  $d\in\mathbb{D}$ to construct natural isomorphisms
  \begin{equation}\label{coalg.tensor.diag2}
    \let\labelstyle=\textstyle
    \xymatrix@=2.25em{
      {\mathcal{E}_\refl}\ar[r]^<>(0.5){\ihoml{C}}
      \ar[d]_{\homout{C\otimes D}} &
      {\Alg<\mathbb{T};\mathcal{E}_\refl>}
      \ar[d]^{\Alg<\mathbb{T};\homout{D}>} \\
      {\Al{\mathbb{T}\otimes\mathbb{S}}} 
      \ar[r]_<>(0.5){\cong} &
      {\Al{\mathbb{T}}(\Al{\mathbb{S}})}
      \ar@{}[ul]|{\textstyle\cong}
    }\mkern15mu
    \xymatrix@=2.25em{
      {\mathcal{E}_\refl}\ar[r]^<>(0.5){\ihomr{C}}
      \ar[d]_{\homout{D\otimes C}} &
      {\Alg<\mathbb{T};\mathcal{E}_\refl>}
      \ar[d]^{\Alg<\mathbb{T};\homout{D}>} \\
      {\Al{\mathbb{S}\otimes\mathbb{T}}} 
      \ar[r]_<>(0.5){\cong} &
      {\Al{\mathbb{T}}(\Al{\mathbb{S}})}
      \ar@{}[ul]|{\textstyle\cong}
    }
  \end{equation}
  respectively.
\end{obs}
  
\begin{obs}
  \label{conc.tensor.coalg}
  Consider the concrete case where $\mathbb{S}$ is an LE-theory,
  $\mathcal{E}$ is the presheaf category $\funcat[\op{\mathbb{D}},
  \Set]$ and $\refl$ is the reflector for which $\mathcal{E}_\refl$ is
  the category of algebras $\Al{\mathbb{S}}$.
  
  Applying Yoneda's lemma as in definition~\ref{weighted.lims}, we
  know that if $\arrow \gamma:W->V.$ is any natural transformation of
  weights then $\arrow\colim(\gamma,\Yoneda):\colim(W,\Yoneda)->
  \colim(V,\Yoneda).$ is actually isomorphic to $\gamma$ itself, so it
  follows immediately, by our colimit characterisation of L-almost
  coalgebras, that the Yoneda functor $\arrow\Yoneda:\mathbb{D}->
  \funcat[\op{\mathbb{D}}, \Set].$ is an L-almost
  $\mathbb{S}$-coalgebra. Furthermore, applying Yoneda's lemma once
  again we see that the associated functor $\homout{\Yoneda}$ is
  actually the identity functor on $\Al{\mathbb{S}}$.  
  
  Consequently if $C$ is any L-almost $\mathbb{T}$-coalgebra in
  $\funcat[\op{ \mathbb{D}},\Set]$ then we may apply this result in
  the context of the essentially commutative squares of
  display~(\ref{coalg.tensor.diag2}), reducing their right hand
  verticals to identities and yielding essentially commutative
  triangles:
  \begin{equation*}
    \let\labelstyle=\textstyle
    \xymatrix@C=0.6em@R=2.6em{
      & {\Al{\mathbb{S}}}
      \ar[dl]_{\homout{C\otimes\Yoneda}}^<>(0.6){}="1"
      \ar[dr]^{\ihoml{C}}_<>(0.6){}="2" & \\
      {\Al{\mathbb{T}\otimes\mathbb{S}}}
      \ar[rr]_<>(0.5){\cong} && 
      {\Al{\mathbb{T}}(\Al{\mathbb{S}})}
      \ar@{}"1";"2"|{\textstyle\cong}
    }\mkern30mu
    \xymatrix@C=0.6em@R=2.6em{
      & {\Al{\mathbb{S}}}
      \ar[dl]_{\homout{\Yoneda\otimes C}}^<>(0.6){}="1"
      \ar[dr]^{\ihomr{C}}_<>(0.6){}="2" & \\
      {\Al{\mathbb{S}\otimes\mathbb{T}}}
      \ar[rr]_<>(0.5){\cong} && 
      {\Al{\mathbb{T}}(\Al{\mathbb{S}})}
      \ar@{}"1";"2"|{\textstyle\cong}
    }
  \end{equation*}
  In practise the left hand diagonals of these triangles provide
  conveniently explicit external presentations of the internal
  constructions on their right hand diagonals.
\end{obs}



\section{Double Categories, 2-Categories and n-Categories}\label{sect.ncats}

\subsection{Categories in the Small}

\begin{obs}
  In subsection~\ref{subsect.2cat} we introduced 2-categories ``in the
  gros sense'', that is to say our definitions were selected to make
  them a little more amenable to expressing some of the meta-theory in
  this work.  In this section we re-introduce categories,
  2-categories and their higher dimensional relatives the
  \inf-categories ``in the petit sense''.  That is to say, our purpose
  here is to present them as objects of algebraic study in their own
  right rather than as contexts within which to study the properties
  of other mathematical entities.
  
  Notice here that the distinction we make between the gros and the
  petit isn't really a matter of set-theoretic size. On the one hand,
  it is true that all of our petit categories and $n$-categories will
  be (defined to be) set-theoretically small and that in our
  meta-theory the gros categories and 2-categories we manipulate will
  be set-theoretically large or huge respectively. However, the
  important distinction here really is one of intention, petit
  categories are inhabitants of our theory while our gros categories
  are denizens of our meta-theory.
\end{obs}

\begin{block}
\def\alpha{p} 
\def\beta{q}
\def\gamma{r}
\begin{defn}[categories in the small]\label{cat.two.sort}
  A (small) category $\sc{C}$ is a sextuple $\cattuple{C}$ consisting of
  \begin{itemize}
  \item a set $\scobj{C}$ whose elements we think of as being the
    {\em objects\/} of $\sc{C}$,
  \item a set $\scarr{C}$ whose elements which we think of as being the
    {\em arrows\/} of $\sc{C}$,
  \item a function $\arrow \idt:\scobj{C}->\scarr{C}.$ which we think of
    as mapping an object to the {\em identity\/} arrow on that object,
  \item two functions $\arrow \s,\t:\scarr{C}->\scobj{C}.$ mapping
    each arrow to its {\em source\/} and {\em target\/} objects
    respectively.  Collectively, these functions must satisfy the
    condition that $\idt$ is a right compositional inverse for both $\s$
    and $\t$ (that is $\t\circ\idt=\s\circ\idt=\id_{\scobj{C}}$). In other words
    for any object $c\in \scobj{C}$ the source and target of the
    corresponding identity arrow $\idt(c)$ must be $c$ itself. We adopt
    the notation $\s'\defeq \idt\circ\s$ and $\t'\defeq\idt\circ\t$ to
    denote the idempotent functions which map each arrow to the
    identity arrows associated with its source and target objects
    respectively.
  \item a partially defined binary composition operation $\comp$ on
    $\scarr{C}$, such that $\beta\comp\alpha\in\scarr{C}$ is defined
    for each pair of arrows $\alpha,\beta\in\scarr{C}$ with
    $\s(\beta)=\t(\alpha)$ (we say that such a pair is {\em
      compositionally compatible}), and
  \item this data should satisfy the equations
    \begin{displaymath}
      \s(\beta\comp\alpha) = \s(\alpha) \mkern20mu 
      \t(\beta\comp\alpha) = \t(\beta) \mkern20mu 
      \t'(\alpha)\comp\alpha = \alpha \mkern20mu 
      \alpha\comp\s'(\alpha) =  \alpha \mkern20mu 
      \gamma\comp(\beta\comp\alpha) = (\gamma\comp\beta)\comp\alpha 
    \end{displaymath}
    whenever the various composites involved are well defined.
  \end{itemize}
  
  When speaking of more than one category, we will tend to follow
  traditional practise by using the symbols $\comp$,$\idt$,$\s$ and $\t$
  ``polymorphically'' and rely on context to disambiguate the
  particular instance of these operators being applied in any given
  case.
  
  A functor $\arrow f:\sc{C}->\sc{C}'.$ consists of a function $\arrow
  f:\scarr{C}->\scarr{C}'.$ between sets of arrows, which must satisfy
  the conditions
  \begin{displaymath}
    \s'(f(\alpha)) =  f(\s'(\alpha)) 
    \t'(f(\alpha)) = f(\t'(\alpha)) \mkern20mu 
    f(\beta\comp\alpha) = f(\beta)\comp f(\alpha)
  \end{displaymath}
  wherever the composites involved are well defined. Of course, the
  composite of (the underlying functions of) two functors is itself a
  functor, and so the collection of all categories and functors itself
  possesses the structure of a (large) category, which we shall call
  $\Cat$.
  
  Notice that a functor $\arrow f:\sc{C}->\sc{C}'.$ gives rise to a
  map $\arrow f_o\defeq\s\circ f\circ \idt:\scobj{C}->\scobj{C}'.$ of
  sets of objects which is the unique function for which any one (and
  thus all) of the equalities $\idt\circ f_o=f\circ\idt$,
  $f_o\circ\s=\s\circ f$ and $f_o\circ\t=\t\circ f$ hold.
\end{defn}

\begin{obs}[single-sorted categories]
  \label{cat.one.sort}
  The category notion given in the last definition is often referred
  to as the {\em two-sorted\/} theory of categories, since each such
  category possesses both a sort of arrows and a separately specified
  sort of objects. In order to avoid definitional burden, some authors
  have chosen not to carry around this separate sort of objects and
  have instead relied upon a {\em single-sorted\/} presentation.  In
  this definition we don't distinguish objects and arrows, instead we
  identify each object with its associated identity arrow.
  
  Under such a presentation, a category becomes a quadruple
  $\quadruple<\scarr{C};\comp;\s';\t'>$ where $\s'$ and $\t'$ are
  endo-functions on $\comp$ satisfying the ``absorption properties''
  $\s'\circ\s'=\t'\circ\s'=\s'$ and $\t'\circ\t'=\s'\circ\t'=\t'$ and
  a pair of arrows $\alpha$ and $\beta$ in $\scarr{C}$ is composable
  iff $\t'(\alpha)=\s'(\beta)$. No information is lost in the process
  and we may regain the sort of objects $\scobj{C}$ simply by
  splitting the idempotent $\arrow \s':\scarr{C}->\scarr{C}.$. That is
  we let $\scobj{C}=\{\alpha\in \scarr{C}| \s'(\alpha)=\alpha\}$, the set of
  {\em identity\/} arrows for the composition operation $\comp$.
  
  The absorption properties of $\s'$ and $\t'$ ensure that we may
  factor them simultaneously through the inclusion $\idt$ of $\scobj{C}$
  into $\scarr{C}$, giving us functions $\arrow
  \s,\t:\scarr{C}->\scobj{C}.$ for which $\s'=\idt\circ\s$,
  $\s\circ\idt=\id_{\scobj{C}}$, $\t'=\idt\circ\t$ and
  $\t\circ\idt=\id_{\scobj{C}}$. This data provides us with the
  underlying structure of a two-sorted category, with $\idt$ as the
  identities function and $\s$ and $\t$ as source and target functions
  respectively.
  
  Given this discussion, we will sometimes take the liberty of passing
  without comment between single-sorted and two-sorted descriptions
  of categories. However in most cases we will prefer the two-sorted
  definition, simply because it will be convenient in many proofs to
  carry around an explicitly chosen presentation of each set of
  objects.
\end{obs}
\end{block}

\begin{obs}[categories as algebras for an LE-theory]
  \label{le-th.cat}
  The theory of categories may be expressed as an LE-theory
  $\mathbb{T}_{\Cat}$ which has underlying category $\Delta$ and the
  set of horn inclusions $\overinc\subseteq_s:\Lambda_k[n]->
  \Delta[n].\in\Simp=\funcat[\op\Delta,\Set]$ ($n=2,3,...$,
  $k=1,...,n-1$) as its FP-regulus.  We let
  $\arrow\refl_{\Cat}:\Simp->\Simp.$ denote the reflector associated
  with this theory.
  
  Alternatively, it is also common to present this LE-theory in terms
  of a finite sketch with pushout-like cocones of the form
  \begin{equation}\label{cat.th.sk}
    \xymatrix@=2.5em{
      {[0]}\ar[r]^{\vertex^p_p}
      \ar[d]_{\vertex^q_0}\ar@{-->}[dr] &
      {[p]}\ar@{-->}[d]^{\partinj^{p,q}_1} \\
      {[q]}\ar@{-->}[r]_-{\partinj^{p,q}_2} & {[p+q]}
    }
  \end{equation}
  in $\Delta$ for each pair of integers $p,q\geq 1$. The intuition
  behind this presentation is that if $A\in\Simp=\funcat[\op{\Delta},
  \Set]$ is a $\mathbb{T}_{\Cat}$-algebra then we would like to think
  of its set of $n$-simplices $A_n$ as being the set of composable
  paths of arrows of length $n$ in some category. Were this the case,
  then we should be able to form a path in $A_{p+q}$ by adjoining a
  pair of compatible paths in $A_p$ and $A_q$, which is precisely what
  the algebra condition with respect to the cocone in
  display~\eqref{cat.th.sk} states.

  If $A\in\Simp$ is a $\mathbb{T}_{\Cat}$-algebra we can see that it
  gives rise to a category, in the sense of
  definition~\ref{cat.two.sort}, with
  \begin{itemize}
  \item\label{pathcat.arrows} Sets of objects $A_0$ and arrows
    $A_1$, the sets of 0- and 1-dimensional simplices in $A$
    respectively.
  \item\label{pathcat.st} Identity, source and target functions $\idt
    \defeq A_{\eta^1}$, $\s \defeq A_{\vertex^1_0}$ and $\t \defeq
    A_{\vertex^1_1}$, for which we may apply the functoriality of $A$,
    along with elementary equations in $\Delta$, to demonstrate the
    required identity rules $\s\circ\idt=\t\circ\idt=\id_{A_0}$.
  \item\label{pathcat.comp} From notation~\ref{simplex.horn} we know
    that a simplicial map $\arrow f:\Lambda_1[2]->A.$ corresponds to a
    composable pair of arrows in $A_1$ and since $A$ is a
    $\mathbb{T}_{\Cat}$-algebra it is perpendicular to
    $\overinc\subseteq_s:\Lambda_1[2]->\Delta[2].$ and it follows that
    the function $\arrow\pair<{-}\cdot\face^2_0;{-}\cdot\face^2_2>:
    A_2->{A_1\pb{\s}{\t}A_1}.$ is an isomorphism with inverse which
    we'll denote by $w$.  In other words, if $p$ and $q$ are a
    composable pair of arrows then there exists a (unique) 2-simplex
    $w\pair<p;q>\in A_2$ with $p=w\pair<p;q>\cdot\face^2_2$ and
    $q=w\pair<p;q>\cdot\face^2_0$ and we define their composite
    $q\comp p$ to be the $3^{\text{rd}}$ 1-dimensional face
    $w\pair<p;q>\cdot\face^2_1$ of $w\pair<p;q>$. We say that
    $w\pair<p;q>$ {\em witnesses\/} the composite $q\comp p$.
  \item\label{pathcat.assoc} The associativity of this composition
    operation can be established by taking a composable sequence
    $p,q,r$ of 1-simplices and exploiting orthogonality to
    $\overarr\subseteq_s:\Lambda_1[2]->\Delta[2].$ three times to
    build a 2-dimensional 1-horn which we fill to a 3-simplex using
    orthogonality to $\overarr\subseteq_s:\Lambda_1[3]-> \Delta[3].$.
    The four 2-dimensional faces of this 3-simplex witness the various
    composites involved in the two sides of the associativity
    condition, viz $r\comp(q\comp p)$ and $(r\comp q)\comp p$, however
    its 1-dimensional faces corresponding to these two composites
    coincide, and so the composites themselves must be equal as
    required.
  \end{itemize}

  If $\arrow f:A->A'.$ is a simplicial map between
  $\mathbb{T}_{\Cat}$-algebras then it carries the 2-simplex
  witnessing a composite $q\comp p$ in $A$ to the 2-simplex witnessing
  the composite $f(q)\comp f(p)$ in $A'$ and so the action of $f$ is
  functorial with respect to these category structures.  It follows
  that this construction provides us with a functor from
  ${\Alg<\mathbb{T}_{\Cat};\Set>}$ to $\Cat.$ and it is a matter of
  routine verification to demonstrate that this is an equivalence. It
  follows that the LE-theory $\mathbb{T}_{\Cat}$ does indeed provide
  us with a suitable presentation of the theory of (small) categories.
  Hence forth, we shall move backward and forward between these
  presentations of the theory of categories without comment.
\end{obs}

\begin{obs}[relativising our theory]
  Generally, all of the theories introduced in this work, including
  the all important theory of complicial sets introduced later on, may
  be expressed as LE-theories. Furthermore, all of the theorems that
  we prove about these structures may be lifted {\em representably\/}
  (that is to say via an application of Yoneda's lemma) to an internal
  context where we take models of our LE-theories in an arbitrary
  LE-category.
  
  Indeed, there will many places in the proofs that follow in which we
  will prove a theorem in $\Set$ and then apply it to algebras
  internal to some explicitly described LFP-category $\mathcal{E}$. In
  all such cases, however, it will be clear that the results we are
  applying immediately lift ``point-wise'' to corresponding results on
  structures defined in terms of the concrete algebras in
  $\mathcal{E}$.
  
  Consequently, in future we will silently assume that the reader is
  aware of this process of {\em relativisation\/} and make no further
  comment when lifting results to an internal context. Furthermore, we
  will generally adopt the notational convention that if
  $\category{C}$ denotes the category of algebras of such a theory
  then $\category{C}(\mathcal{E})$ is used to denote the corresponding
  category of models in $\mathcal{E}$. As before, if $\arrow
  F:\mathcal{E}->\mathcal{E}'.$ is an LE-functor then the
  corresponding functor of categories of algebras may sometimes be
  denoted by $\category{C}(F)$, but is most likely to simply be
  referred to by ``overloading'' the symbol $F$ and saying that it
  lifts to a LE-functor $\arrow F:\category{C}(\mathcal{E})->
  \category{C}(\mathcal{E}').$
\end{obs}

\vfill

\subsection{Double Categories}

\begin{defn}[double categories and double functors]\label{doub.defn}
  A (small) double category is a model for the finite limit theory of
  categories in the category $\Cat$ of (small) categories and
  functors. We define the category $\Double$ to be $\Cat(\Cat)$ and
  call its arrows {\em double functors}.
  
  Alternatively we know, by observations~\ref{le-th.cat}
  and~\ref{le-th.tensor}, that a double category $\mathbb{D}$ may be
  presented more symmetrically as an algebra of the LE-theory
  $\mathbb{T}_{\Cat}\otimes\mathbb{T}_{\Cat}$, that is a functor
  $\arrow \mathbb{D}:\op\Delta\times\op\Delta->\Set.$ such that each
  of the functors $\arrow \mathbb{D}\pair<[n];{-}>,
  \mathbb{D}\pair<{-};[m]>: \op\Delta ->\Set.$ are categories
  (presented as a $\mathbb{T}_{\Cat}$-algebras) for all
  $[n],[m]\in\Delta$.
  
  In more economical and explicit terms, we may present our double
  category $\dc{D}$ as having:
  \begin{itemize}
  \item a set $\dcsq{D}$ whose elements are called {\em squares\/}, sets
    $\dcvarr{D}$ and $\dcharr{D}$ whose elements are called {\em vertical\/} and
    {\em horizontal\/} arrows respectively and a set $\dcobj{D}$ of {\em
      objects}.
  \item four category structures 
    \begin{itemize}
    \item $\dcattup<\dcharr{D};\dcobj{D};h>$ the {\em horizontal 
        category of arrows}, denoted $\arr_h(\mathbb{D})$,
    \item $\dcattup<\dcvarr{D};\dcobj{D};v>$ the {\em vertical 
        category of arrows}, denoted $\arr_v(\mathbb{D})$,
    \item $\dcattup<\dcsq{D};\dcvarr{D};h>$ the {\em horizontal
        category of squares}, denoted $\sq_h(\mathbb{D})$, and
    \item $\dcattup<\dcsq{D};\dcharr{D};v>$ the {\em vertical
        category of squares}, denoted $\sq_v(\mathbb{D})$.
    \end{itemize}
  \end{itemize}
  This data must satisfy the following rules:
  \begin{enumerate}[(i)]
  \item the vertical structural functions $\idt_v$, $\s_v$ and $\t_v$ of
    $\sq_v(\mathbb{D})$ are functors with respect to the horizontal
    category structures $\arr_h(\mathbb{D})$ and $\sq_h(\mathbb{D})$,
  \item the horizontal structural functions $\idt_h$, $\s_h$ and $\t_h$ of
    $\sq_h(\mathbb{D})$ are functors with respect to the vertical
    category structures $\arr_v(\mathbb{D})$ and $\sq_v(\mathbb{D})$,
  \item the middle four interchange rule holds, that is for any
    quadruple of squares $\alpha,\beta,\gamma,\delta$ we have
    $(\delta\comp_h\gamma)\comp_v(\beta\comp_h\alpha) =
    (\delta\comp_v\beta)\comp_h(\gamma\comp_v\alpha)$ where the
    composite on the left is well defined iff that on the right is.
  \end{enumerate} 
  Notice that the middle four interchange rule is equivalent to saying
  that the horizontal composite map $\comp_h$ of $\sq_h(\mathbb{D})$
  is functorial with respect to the vertical category structure on
  $\sq_v(\mathbb{D})$ and that of its pullback
  $\sq_v(\mathbb{D})\pb{\s_h}{\t_h}\sq_v(\mathbb{D})$. Dually, this
  rule may be interpreted as postulating the functoriality of vertical
  composition $\comp_v$ with respect to horizontal category
  structures.
  
  Under this presentation, a double functor $\arrow
  f:\mathbb{D}->\mathbb{D}'.$ is simply a function $\arrow
  f:\dcsq{D}->\dcsq{D}'.$ which is functorial with respect to both
  the vertical and the horizontal categories of squares of
  $\mathbb{D}$ and $\mathbb{D}'$. In other words, $f$ is the
  underlying function of two functors $\arrow
  f:\sq_v(\mathbb{D})->\sq_v(\mathbb{D}').$ and $\arrow
  f:\sq_h(\mathbb{D})->\sq_h(\mathbb{D}').$.
  
  Given a double category $\mathbb{D}$ presented as a
  $\mathbb{T}_{\Cat}\otimes\mathbb{T}_{\Cat}$-algebra its
  corresponding explicit presentation may be obtained by:
  \begin{itemize}
  \item defining the various sets of squares and arrows in
    $\mathbb{D}$ to be:
    \begin{displaymath}
        \dcsq{D} \defeq \dc{D}\pair<[1];[1]> \mkern10mu 
        \dcobj{D} \defeq  \dc{D}\pair<[0];[0]> \mkern10mu 
        \dcvarr{D} \defeq  \dc{D}\pair<[1];[0]> \mkern10mu 
        \dcharr{D} \defeq \dc{D}\pair<[0];[1]>
     \end{displaymath}
  \item letting the various categories $\arr_h(\mathbb{D})$,
    $\arr_v(\mathbb{D})$, $\sq_h(\mathbb{D})$ and $\sq_v(\mathbb{D})$
    be those which correspond to the $\mathbb{T}_{\Cat}$-algebras
    $D\pair<[0];{-}>$, $D\pair<{-};[0]>$, $D\pair<[1];{-}>$ and
    $D\pair<{-};[1]>$ respectively.
  \end{itemize}
  
  The isomorphisms of categories of algebras given at the end of
  observation~\ref{le-th.tensor}, provide us with two distinct
  presentations of each ``symmetrised'' double category $\mathbb{D}$
  as an internal category in $\Cat$, viz:
  \begin{itemize}
  \item its {\em horizontal presentation\/} $\sextuple<
    \sq_v(\mathbb{D}); \arr_v(\mathbb{D}); \comp_h; \idt_h; \s_h; \t_h>$, an
    internal category of vertical categories, and
  \item its {\em vertical presentation\/} $\sextuple<
    \sq_h(\mathbb{D}); \arr_h(\mathbb{D}); \comp_v; \idt_v; \s_v; \t_v>$, an
    internal category of horizontal categories.
  \end{itemize}
\end{defn}

\begin{obs}[double categories in pictures]\label{dcat.piccy} 
  A square $\alpha$ of a double category $\mathbb{D}$ may literally be
  pictured as a square thus:
  \begin{displaymath}
    \let\labelstyle=\textstyle
    
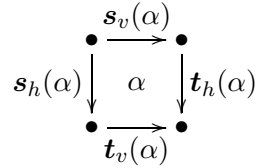
\xymatrix@=2em{
      \bullet \ar[r]^{\s_v(\alpha)} 
      \ar[d]_{\s_h(\alpha)} 
      \ar@{}[dr]|{\alpha} & 
      \bullet \ar[d]^{\t_h(\alpha)} \\
      \bullet \ar[r]_{\t_v(\alpha)} & \bullet }
  \end{displaymath}
  Here, the vertical (resp. horizontal) lines represent the vertical
  (resp. horizontal) arrows which constitute the horizontal (resp.
  vertical) domain and codomain of $\alpha$.  It is a consequence of
  the functoriality of horizontal and vertical domain and codomain
  functions that adjacent pairs of sides in this diagram ``meet up''
  at a corner, each of which is an object.  Furthermore, two squares
  are horizontally (resp.\ vertically) composable if they abut at a
  common vertical (resp.\ horizontal) arrow.  Finally, we can picture
  the middle four interchange rule in terms of the two possible ways
  of composing four squares, as shown in figure~\ref{doub.mid.four}.
  \begin{figure}[h]
    \begin{displaymath}
      \xymatrix@R=1.2em@C=1.2em{
        \bullet\ar[r]\ar[d]\ar@{}[dr]|{\textstyle\alpha} & 
        \bullet\ar[d] & 
        \bullet\ar[r]\ar[d]\ar@{}[dr]|{\textstyle\beta} & 
        \bullet\ar[d]^{}="one" & &
        \bullet\ar[rr]\ar[d]_{}="three" & 
        \ar@{}[d]|{\beta\comp_h\alpha} & \bullet\ar[d] \\
        \bullet\ar[r] & \bullet & \bullet\ar[r] & \bullet & &
        \bullet\ar[rr] & & \bullet \\
        \bullet\ar[r]\ar[d]\ar@{}[dr]|{\textstyle\gamma} & 
        \bullet\ar[d] & 
        \bullet\ar[r]\ar[d]\ar@{}[dr]|{\textstyle\delta} & 
        \bullet\ar[d]^{}="two" & &
        \bullet\ar[rr]\ar[d]_{}="four" & 
        \ar@{}[d]|{\delta\comp_h\gamma} & \bullet\ar[d] \\
        \bullet\ar[r]_{}="five" & \bullet & \bullet\ar[r]_{}="six" & 
        \bullet & & \bullet\ar[rr]_{}="nine" & & \bullet \\
        & & & & & & & \\
        \bullet\ar[dd]\ar[r]^{}="seven" & \bullet\ar[dd] &
        \bullet\ar[dd]\ar[r]^{}="eight" & 
        \bullet\ar[dd]^{}="eleven" & &
        \bullet\ar[dd]_{}="twelve"\ar[rr]^{}="ten" & & 
        \bullet\ar[dd] \\
        \ar@{}[r]|{\gamma\comp_v\alpha} & & 
        \ar@{}[r]|{\delta\comp_v\beta} & & &
        \ar@{}[rr]|{\begin{array}{l}
            \scriptstyle(\delta\comp_h\gamma)\comp_v(\beta\comp_h\alpha) = \\
            \scriptstyle(\delta\comp_v\beta)\comp_h(\gamma\comp_v\alpha)
          \end{array}} & {{\vrule height0em depth0em width3em}
          {\vrule height3em depth0em width0em}} & \\
        \bullet\ar[r] & \bullet & \bullet\ar[r] & \bullet & &
        \bullet\ar[rr] & & \bullet
        \ar@{|->} "one" ; "three"
        \ar@{|->} "two" ; "four"
        \ar@{|->} "five" ; "seven"
        \ar@{|->} "six" ; "eight"
        \ar@{|->} "nine" ; "ten"
        \ar@{|->} "eleven" ; "twelve"
      }
    \end{displaymath}
    \caption{The middle four interchange rule for double categories.}
    \label{doub.mid.four}
  \end{figure}
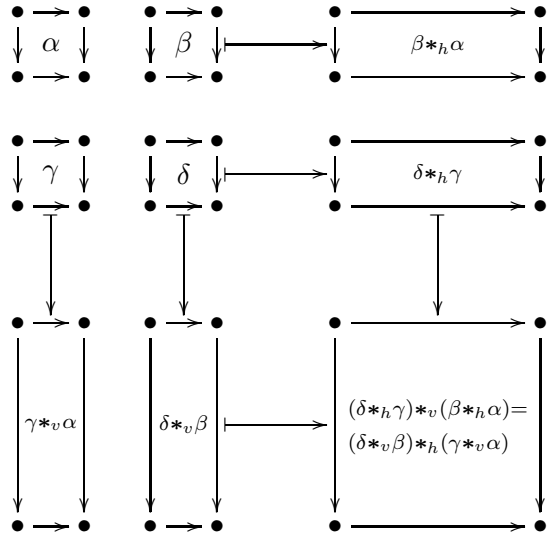
\end{obs}

\subsection{2-Categories and Double Categories with Connections}
\label{sect.dcat.with.conn}

Now we are ready to give another definition of 2-categories, which is
easily seen to be equivalent to that given in
subsection~\ref{subsect.2cat}:

\begin{defn}[2-categories in the small]\label{2cat.qua.doub}
  A category $\mathbb{C}$ is said to be {\em discrete\/} iff each of
  its arrows is an identity on some object. Equivalently $\mathbb{C}$
  is discrete iff its identity function
  $\arrow\idt:\scobj{C}->\scarr{C}.$ is an isomorphism.
  
  A 2-category is a double category $\mathbb{D}$ whose category of
  vertical arrows $\arr_v(\mathbb{D})$ is discrete and a 2-functor is
  simply a double functor between 2-categories. We use the notation
  $\Twocat$ for the full subcategory of 2-categories in $\Double$.
\end{defn}

\begin{obs}[The relationship to definition~\ref{2cat.defn.1}]
  Under our new definition, a 2-category is simply a double category
  in which each vertical arrow is an identity. In other words, if we
  depict identity arrows as equalities then we might depict the
  squares of a 2-category thus:
  \begin{displaymath}
    \xymatrix@=2em{
      \bullet\ar[r]^p\ar@{=}[d] & \bullet\ar@{=}[d]
      \ar@{}[dl]|{\textstyle\alpha} \\
      \bullet\ar[r]_{p'} & \bullet }
  \end{displaymath}
  Contracting these equalities to a single point and depicting our
  square as a double arrow from its vertical domain to its vertical
  codomain we simply regain the ``globular'' 2-cell picture of
  observation~\ref{2cat.expl}. 
  
  Consequently, we will adopt the usual nomenclature for the
  structural components of 2-categories presented in this way,
  referring to their squares as 2-cells, horizontal arrows as 1-cells
  and objects as 0-cells.
\end{obs}

We've seen that 2-categories may be considered simply to be special
kinds of double categories. There is, however, a slightly less trivial
construction by which we can build a double category out of a
2-category and which we shall need to apply in the sequel.

\begin{obs}[the double category of pasting squares in a
  2-category]\label{dcat.psq} If $\mathbb{C}$ is a 2-category then we
  may define a double category $\squares(\mathbb{C})$ with
  squares consisting of quintuples $\quintuple<\psi;d_h;c_h;d_v;c_v>$
  where $c_h$, $d_h$, $c_v$ and $d_v$ are all 1-cells (in
  $\mathbb{C}$), $\psi$ is a 2-cell and this information fits together
  into a pasting square:
  \begin{displaymath}
    \let\labelstyle=\textstyle
    \xymatrix@=2.5em{
      \bullet \ar[r]^{d_v} \ar[d]_{d_h} & 
      \bullet \save []-<1.2em,1.2em>="test2" \restore
      \ar[d]^{c_h} \\
      \bullet \save []+<1.2em,1.2em>="test" \restore 
      \ar[r]_{c_v} & \bullet 
      \ar@{{}=>}_{\psi} "test2";"test"
    }
  \end{displaymath}
  The operations of horizontally and vertically pasting together such
  squares (see Kelly and Street~\cite{KellyStreet:1974:2Cats} or
  Power~\cite{Power:1988:Paste}) suffices to provide us with natural
  horizontal and vertical composition operations. Furthermore, it is
  an immediate consequence of the theory of pasting that a 2-functor
  $\arrow f:\mathbb{C} ->\mathbb{C}'.$ gives rise to a double functor
  $\arrow \squares(f): \squares(\mathbb{C})-> \squares(\mathbb{C}').$,
  which acts ``point-wise'' on pasting squares as in the next diagram:
  \begin{displaymath}
    \let\labelstyle=\textstyle
    \xymatrix@=2.6em{
      \bullet \ar[r]^{d_v} \ar[d]_{d_h} & 
      \bullet \save []-<1.2em,1.2em>="test2" \restore
      \ar[d]^{c_h}="test6" & &
      \bullet \ar[r]^{f(d_v)} \ar[d]_{f(d_h)}="test5" & 
      \bullet \save []-<0.8em,1.8em>="test3" \restore
      \ar[d]^{f(c_h)} \\
      \bullet \save []+<1.2em,1.2em>="test" \restore 
      \ar[r]_{c_v} & \bullet & &
      \bullet \save []+<1.8em,0.8em>="test4" \restore 
      \ar[r]_{f(c_v)} & \bullet
      \ar@{{}=>}_{\psi} "test2";"test"
      \ar@{{}=>}_(0.3){f(\psi)} "test3" ; "test4"
      \ar@{|->} "test6" ; "test5"
    }
  \end{displaymath}  
\end{obs}

\begin{obs}[duals of double categories and 2-categories]\label{double.duals}
  A double category $\mathbb{D}$ admits two basic kinds of duality:
  \begin{itemize}
  \item The {\em horizontal\/} and {\em vertical\/} duals
    $\hop{\mathbb{D}}$ and $\vop{\mathbb{D}}$ which are obtained by
    applying the usual categorical duality to its horizontal and
    vertical category structures respectively. We will often use
    $\hvop{\mathbb{D}}$ to denote the result of applying both of these
    dualities.
  \item The {\em reflection\/} dual $\refld{\mathbb{D}}$ which is
    obtained by swapping the roles of these horizontal and vertical
    category structures.
  \end{itemize}
  Geometrically the first two dualities correspond to reversing the
  orientation of squares in one dimension or the other. The second
  corresponds to flipping or reflecting squares through their leading
  diagonals.
  
  Of course, these dualities are all canonically functorial and
  involutive, in the sense that apply any one of them twice returns us
  to the double category we started with, and trivially interrelated
  via the (inter-derivable) identities $\hop{(\refld{\mathbb{D}})} =
  \refld{(\vop{\mathbb{D}})}$ and $\vop{(\refld{\mathbb{D}})} =
  \refld{(\hop{\mathbb{D}})}$.

  2-categorists often use the notation $\op{\mathbb{C}}$ for the
  horizontal dual and $\co{\mathbb{C}}$ for the vertical dual of a
  2-category $\mathbb{C}$. They also tend to use the notation
  $\coop{\mathbb{C}}$ to denote the 2-category obtained by applying
  both of these dualities.
\end{obs}

\begin{obs}[The 2-category of globs in a double category]\label{globs.defn}
  For a double category $\mathbb{D}$, let $\glob(\mathbb{D})$ be the
  sub-double category of squares $\psi\in\mathbb{D}$ such that
  $\t_h(\psi)$ and $\s_h(\psi)$ are both identities in
  $\arr_v(\dc{D})$.  From this definition it is clear that
  $\glob(\mathbb{D})$ is actually a 2-category, which we call the {\em
    2-category of globs\/} in $\mathbb{D}$. 
  
  It is easily shown that a double functor $\arrow
  f:\mathbb{D}->\mathbb{D}'.$ maps squares in $\glob(\mathbb{D})$ into
  $\glob(\mathbb{D}')$ and so restricts to a 2-functor
  $\arrow\glob(f):\glob(\mathbb{D})-> \glob(\mathbb{D}').$. It follows
  that $\glob$ may be extended to a functor
  $\arrow\glob:\Double->\Twocat.$.
  
  Notice that the 2-category of globs construction allows us to
  re-construct a 2-category from its double category of squares, in
  the sense that there exists a canonical family of isomorphisms
  $\mathbb{C}\cong\glob(\squares(\mathbb{C}))$ which is natural in the
  2-category $\mathbb{C}$.
  
  It follows that we may immediately infer that
  $\arrow\squares:\Twocat->\Double.$, the pasting squares functor, is
  faithful. However, it is {\bf not} the case that this functor is also
  full, a fact which immediately begs two related questions:
  \begin{itemize}
  \item What extra structure might double categories of squares be
    asked to carry to restrict the maps between them and thus make the
    squares functor full?
  \item How can we characterise those double categories that arise by
    applying the squares functor to some 2-category?
  \end{itemize}
  The remainder of this subsection recalls the answers to these
  questions, in a form originally due to Chris Spencer and presented
  in Brown and Mosa~\cite{Brown:1999:Connections}.
\end{obs}

\begin{defn}[thin squares]\label{strat.double.cat} 
  A {\em double category with thinness\/} is a pair $\stratd{D}$
  consisting of a double category $\mathbb{D}$ and a specified
  sub-double category $t\mathbb{D}\subseteq \mathbb{D}$ which contains
  all of the horizontal and vertical arrows of $\mathbb{D}$ (but which
  does not, in general, contain all of the squares of $\mathbb{D}$).
  It should come as no surprise that the squares of the {\em thinness
  structure\/} $t\mathbb{D}$ are said to be {\em thin}.
  
  We form a category $\Double_T$ whose objects are double categories
  with thinness and whose arrows $\arrow f:\stratd{D}-> \stratd{D'}.$
  are double functors $\arrow f:\mathbb{D}-> \mathbb{D}'.$ which {\em
  preserve thinness\/} in the sense that they satisfy the condition 
  $t\mathbb{D}\subseteq f^{-1}(t\mathbb{D}')$. 
  
  Clearly, in order to provide a double category $\mathbb{D}$ with a
  thinness structure it is enough to specify a subset $t\dcsq{D}$ of
  the set of squares of $\mathbb{D}$ which:
  \begin{itemize}
  \item contains all of the horizontal and vertical identity squares
    of $\mathbb{D}$, and
  \item is closed in $\mathbb{D}$ under horizontal and vertical
    composition of squares.
  \end{itemize}
\end{defn}

\begin{obs}
  The double category of pasting squares in a 2-category $\mathbb{C}$
  admits a canonical thinness structure under which we declare a
  square to be thin if its 2-cellular component is an identity in
  $\mathbb{C}$. If $\arrow f:\mathbb{C}->\mathbb{C}'.$ is a 2-functor
  then it is clear that the associated double functor $\arrow
  \squares(f): \squares(\mathbb{C})-> \squares(\mathbb{C}').$
  preserves thinness with respect to these canonical structures. In
  other words, the pasting squares construction admits a canonical
  lifting to a functor $\arrow \squares:\Twocat->\Double_T.$.
\end{obs}

\begin{lemma}[Brown and Mosa~\cite{Brown:1999:Connections}]
  \label{dcat.with.con} 
  The functor $\arrow \squares:\Twocat->\Double_T.$ is fully faithful,
  furthermore a double category with thinness $\stratd{D}$ is in its
  replete image if and only if it satisfies the following conditions:
  \begin{enumerate}[(i)]
  \item\label{con.one} If $\alpha$ is a thin square such that
    $\s_h(\alpha)$ and $\t_h(\alpha)$ are identities in the category
    of vertical arrows $\arr_v(\mathbb{D})$ then $\alpha$ is a
    vertical identity, diagrammatically:
    \begin{displaymath}
      \text{if }\mkern10mu\vcenter{
        \xymatrix@R=2em@C=2em{
          \bullet\ar[r]^p\ar@{=}[d] & \bullet\ar@{=}[d]
          \ar@{}[dl]|{\textstyle\alpha} \\
          \bullet\ar[r]_{p'} & \bullet }}
      \text{ is thin then $\alpha=\idt_v(p)=\idt_v(p')$ (so $p=p'$)}
    \end{displaymath}
  \item\label{con.two} If $\alpha$ is a thin square such that
    $\s_v(\alpha)$ and $\t_v(\alpha)$ are identities in the category of
    horizontal arrows $\arr_h(\mathbb{D})$ then $\alpha$ is a horizontal
    identity, diagrammatically:
    \begin{displaymath}
      \text{if }\mkern10mu\vcenter{
        \xymatrix@R=2em@C=2em{
          \bullet\ar[d]_q\ar@{=}[r] & \bullet\ar[d]^{q'}
          \ar@{}[dl]|{\textstyle\alpha} \\
          \bullet\ar@{=}[r] & \bullet }}
      \text{ is thin then $\alpha=\idt_h(q)=\idt_h(q')$ (so $q=q'$)}
    \end{displaymath}
  \item\label{con.three} For each horizontal arrow $p\in \arr_h(\dc{D})$
    there exists a thin square $\mu_p$ such that $\t_v(\mu_p)=p$,
    $\s_v(\mu_p)$ is an identity in $\arr_h(\dc{D})$ and $\s_h(\mu_p)$ is
    an identity in $\arr_v(\dc{D})$, diagrammatically:
    \begin{displaymath}
      \text{for each }
      \vcenter{\xymatrix@=2em{\bullet\ar[r]^p & \bullet}}
      \text{ we have a thin }
      \vcenter{\xymatrix@=2em{
        \bullet\ar@{=}[d]\ar@{=}[r] &
        \bullet\ar[d]^{\hat{p}}\ar@{}[dl]|{\textstyle\mu_p} \\
        \bullet\ar[r]_p & \bullet }}
    \end{displaymath}
  \item\label{con.four} For each vertical arrow $q\in\arr_v(\dc{D})$
    there exists a thin square $\nu_q$ such that $\s_h(\nu_q)=q$,
    $\t_h(\nu_q)$ is an identity in $\arr_v(\dc{D})$ and $\t_v(\nu_q)$
    is an identity in $\arr_h(\dc{D})$, diagrammatically:
    \begin{displaymath}
      \text{for each }
      \vcenter{\xymatrix@=2em{\bullet\ar[d]_q \\ \bullet}}\
      \text{ we have a thin }
      \vcenter{\xymatrix@=2em{
        \bullet\ar[d]_q\ar[r]^{\tilde{q}} &
        \bullet\ar@{=}[d]\ar@{}[dl]|{\textstyle\nu_q}\\
        \bullet\ar@{=}[r] & \bullet }}
    \end{displaymath}
  \end{enumerate}
  Note that the various squares $\mu_p$ and $\nu_q$ are called {\em
    connection squares\/} or simply {\em connections\/} and so a
  double category with thinness which satisfies the conditions above
  is often referred to as a {\em double category with connections\/}
  and we use $\Conn$ to denote the full subcategory of such things in
  $\Double_T$.
\end{lemma}

\begin{proof}(sketch, see Brown and Mosa~\cite{Brown:1999:Connections} for
  more detail) Notice that, as pointed out in
  observation~\ref{dcat.psq}, the 2-cells of
  $\glob(\squares(\mathbb{C}))$ are pretty much no more nor less than
  2-cells in $\mathbb{C}$.  A routine calculation shows that this
  correspondence extends to a (natural) isomorphism of 2-categories
  $\mathbb{C}\cong \glob(\squares(\mathbb{C}))$.
  
  It remains to prove that for each double category with connections
  $\stratd{D}$ we have a natural isomorphism $\stratd{D}\cong
  \squares(\glob{\stratd{D}})$ in $\Double_T$. We do this in a number
  of steps:
  \begin{enumerate}[1.]
  \item Our conditions set up a bijection between the horizontal and
    vertical arrows of $\mathbb{D}$:
    \begin{enumerate}[(i)]
    \item if $p$ is a horizontal arrow then condition~(\ref{con.four})
      provides us with a thin square $\mu_p$ and, as a result, a
      vertical arrow $\hat{p}=\t_h(\mu_p)$,
    \item if $q$ is a vertical arrow then condition~(\ref{con.three})
      provides us with a thin square $\nu_q$ and a horizontal arrow
      $\tilde{q}=\s_v(\nu_q)$,
    \item the horizontal arrow $\widetilde{\widehat{p}}$ is equal to
      $p$, which fact may be demonstrated by forming the horizontal
      composite $\nu_{\hat{p}}\comp_h\mu_p$ and appealing to
      condition~(\ref{con.one}).
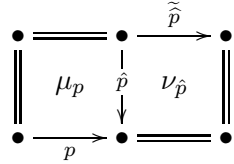
      \begin{displaymath}
        \xymatrix@=2.5em{
          \bullet\ar@{=}[r]\ar@{=}[d]\ar@{}[dr]|{\textstyle\mu_p} & 
          \bullet\ar[r]^{\widetilde{\widehat{p}}}
          \ar[d]|(0.45)*+<3pt>{\scriptstyle\hat{p}}
          \ar@{}[dr]|{\textstyle\nu_{\hat{p}}} & \bullet\ar@{=}[d] \\
          \bullet\ar[r]_p & \bullet\ar@{=}[r] & \bullet }
      \end{displaymath}
    \item the vertical arrow $\widehat{\widetilde{q}}$ is equal to
      $q$, which we demonstrate by forming the vertical composite
      $\nu_q\comp_v\mu_{\tilde{q}}$ and appealing to
      condition~(\ref{con.two}).
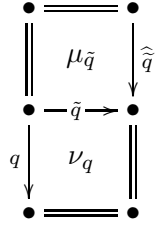
      \begin{displaymath}
        \xymatrix@=2.5em{
          \bullet\ar@{=}[r]\ar@{=}[d]\ar@{}[dr]|{\textstyle\mu_{\tilde{q}}} &
          \bullet\ar[d]^{\widehat{\widetilde{q}}} \\
          \bullet\ar[d]_q\ar[r]|(0.45)*+<3pt>{\scriptstyle\tilde{q}}
          \ar@{}[dr]|{\textstyle\nu_q} &
          \bullet\ar@{=}[d] \\
          \bullet\ar@{=}[r] & \bullet }
      \end{displaymath}
    \item If $p$ and $q$ are a composable pair of vertical arrows then
      we have $\widetilde{q\comp_v p} = \tilde{q}\comp_h\tilde{p}$,
      which we can demonstrate by forming the composite $\nu_{q\comp_v
        p}\comp_h (\mu_{\tilde{q}}\comp_v p) \comp_h\mu_{\tilde{p}}$
      appealing to condition~(\ref{con.one}).
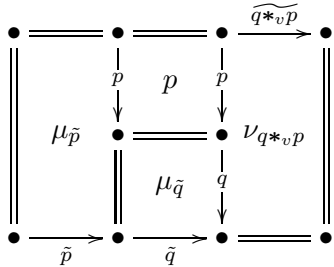
      \begin{displaymath}
        \xymatrix@=2.5em{
          \bullet\ar@{=}[r]\ar@{=}[dd] &
          \bullet\ar@{=}[r]\ar[d]|(0.45)*+<3pt>{\scriptstyle{p}}
          \ar@{}[dr]|{\textstyle{p}} &
          \bullet\ar[r]^{\widetilde{q\comp_v p}}
          \ar[d]|(0.45)*+<3pt>{\scriptstyle{p}} & \bullet\ar@{=}[dd] \\
          {}\ar@{}[r]|{\textstyle\mu_{\tilde{p}}} & 
          \bullet\ar@{=}[r]\ar@{=}[d]
          \ar@{}[dr]|{\textstyle\mu_{\tilde{q}}} &
          \bullet\ar[d]|(0.45)*+<3pt>{\scriptstyle{q}}
          \ar@{}[r]|{\textstyle\nu_{q\comp_v p}} & \\
          \bullet\ar[r]_{\tilde{p}} &
          \bullet\ar[r]_{\tilde{q}} &
          \bullet\ar@{=}[r] & \bullet }
      \end{displaymath}
    \end{enumerate}
    These results provide us with a natural isomorphism between the
    horizontal category $\arr_h(\mathbb{D})$ and the vertical category
    $\arr_v(\mathbb{D})$. Consequently, from here on we identify these
    categories and simplify our notation by dropping explicit use of
    the accents $\sim$ and $\wedge$ to pass to and fro between
    horizontal and vertical arrows.

  \item Each square $\quintuple<\psi;d_h;c_h;d_v;c_v>$ in
    $\squares(\glob \stratd{D})$ gives rise to a square
    $\psi^{\mathord{\circ}}$ in $\mathbb{D}$ of the form
    
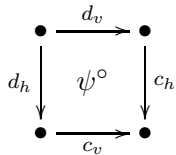
\begin{displaymath}
      \xymatrix@=2.5em{
        \bullet\ar[r]^{d_v}\ar[d]_{d_h}
        \ar@{}[dr]|{\textstyle\psi^{\mathord{\circ}}} &
        \bullet\ar[d]^{c_h} \\
        \bullet\ar[r]_{c_v} & \bullet }
    \end{displaymath}
    obtained by taking the composite $(c_v\comp_h\nu_{d_h})
    \comp_v\psi\comp_v(\mu_{c_h}\comp_h d_v)$, or diagrammatically:
    \begin{displaymath}
      \xymatrix@=2.5em{
        \bullet\ar@{=}[d]\ar[r]^{d_v}\ar@{}[dr]|{\textstyle d_v} & 
        \bullet\ar@{=}[d]\ar@{=}[r]\ar@{}[dr]|{\textstyle\mu_{c_h}} &
        \bullet\ar[d]^{c_h} \\
        \bullet\ar[r]|(0.45)*+<3pt>{\scriptstyle d_v}\ar@{=}[d] &
        \bullet\ar[r]|(0.45)*+<3pt>{\scriptstyle c_h}
        \ar@{}[d]|{\textstyle\psi} &
        \bullet\ar@{=}[d] \\
        \bullet\ar[r]|(0.45)*+<3pt>{\scriptstyle d_h}\ar[d]_{d_h}
        \ar@{}[dr]|{\textstyle\nu_{d_h}} &
        \bullet\ar@{=}[d]\ar[r]|(0.45)*+<3pt>{\scriptstyle c_v}
        \ar@{}[dr]|{\textstyle c_v} &
        \bullet\ar@{=}[d] \\
        \bullet\ar@{=}[r] & \bullet\ar[r]_{c_v} & \bullet }
    \end{displaymath}
    
  \item\label{brown.double.pr.3} A square $\lambda\in\mathbb{D}$ gives
    rise to a pasting square $\quintuple<\lambda^*;q;q';p;p'>$ in
    $\glob\stratd{D}$ where $p=\s_v(\lambda)$, $p'=\t_v(\lambda)$,
    $q=\s_h(\lambda)$, $q'=\t_h(\lambda)$ and $\lambda^*$ is the
    square given by the composite
    $\nu_{q'}\comp_h\lambda\comp_h\mu_{q}$, or
    diagrammatically:
    \begin{displaymath}
      \xymatrix@=2.5em{
        \bullet\ar@{=}[r]\ar@{=}[d]\ar@{}[dr]|*{\textstyle\mu_{q}} &
        \bullet\ar[r]^{p}\ar[d]|(0.45)*+<3pt>{\scriptstyle q}
        \ar@{}[dr]|{\textstyle\lambda} &
        \bullet\ar[r]^{q'}
        \ar[d]|(0.45)*+<3pt>{\scriptstyle q'}\ar@{}[dr]|{\textstyle{\nu_{q'}}} &
        \bullet\ar@{=}[d] \\
        \bullet\ar[r]_{q} & 
        \bullet\ar[r]_{p'} &
        \bullet\ar@{=}[r] & \bullet }
    \end{displaymath}
  \item It is a matter of routine calculation, involving heavy use of
    the middle four rule and conditions (\ref{con.one}) and
    (\ref{con.two}), to demonstrate that these two operations are
    mutually inverse, respect thinness and preserve the composition
    structures of $\stratd{D}$ and $\squares(\glob\stratd{D})$. It
    follows that they give rise to a (natural) isomorphism between
    these two double categories with connections, as required.\qedhere
  \end{enumerate}
\end{proof}

\begin{obs}[duals and double categories of pasting squares]
  \label{duals.and.psquares}
  Later on we will need to understand the relationship between the
  various dualities of observation~\ref{double.duals} and the
  constructions of the last lemma. For the horizontal and vertical
  dualities things are simple since it is easily seen that if
  $\mathbb{D}$ is any double category then $\glob(\hop{\mathbb{D}}) =
  \op{\glob(\mathbb{D})}$ and $\glob(\vop{\mathbb{D}}) =
  \co{\glob(\mathbb{D})}$. In general there is, however, no reason to
  believe that any simple relationship should hold between the
  2-categories $\glob(\mathbb{D})$ and $\glob(\refld{\mathbb{D}})$.
  
  One useful exception, which we shall apply later on, is the case
  when our double category is of the form $\squares(\mathbb{C})$ for
  some 2-category $\mathbb{C}$.  In this case, it is easy to construct
  a canonical family of isomorphisms $\co{\mathbb{C}}\cong\glob(\refld{
    \squares( \mathbb{C})})$ which is natural in the 2-category
  $\mathbb{C}$. Equivalently, we can exploit the result of the last
  lemma to express this as a family of isomorphisms
  $\co{\glob{\stratd{D}}}\cong\glob\strat{\refld{\mathbb{D}}}$ which is
  natural in the double category with connections $\stratd{D}\in\Conn$.
\end{obs}

\subsection{$n$-Categories and \inf-Categories}

One exception to preferring many sorted definitions over single sorted
ones is that of \inf-categories, which appears in the following form
in Street~\cite{Street:1987:Oriental}:

\begin{defn}[\inf-categories]\label{omegacat}
  An {\em \inf-category\/} $\mathbb{C}$ consists of the following data
  \begin{itemize}
  \item a set of elements called {\em cells\/} and
  \item a family of (one-sorted) category structures
    $\triple<\comp_i;\s_i;\t_i>$ one for each $i\in\mathbb{N}$,
  \end{itemize}
  which must satisfy the following conditions
  \begin{enumerate}[(i)]
  \item\label{2cat.cond} if $i<j\in\mathbb{N}$ then category
    structures $\triple<\comp_i;\s_i;\t_i>$ and
    $\triple<\comp_j;\s_j;\t_j>$ constitute the horizontal and
    vertical category structures (respectively) of a 2-category
    structure on $\mathbb{C}$ and
  \item for each cell $\alpha$ there exists an $i\in\mathbb{N}$ such
    that $\alpha$ is an identity for the category structure
    $\triple<\comp_i;\s_i;\t_i>$.
  \end{enumerate}
  If $\mathbb{C}$ and $\mathbb{C}'$ are \inf-categories then an {\em
    \inf-functor\/} is simply a function $\arrow
  f:\mathbb{C}->\mathbb{C}'.$ which is a functor with respect to the
  $i^{\text{th}}$ category structure on $\mathbb{C}$ and $\mathbb{C}'$
  for each $i\in\mathbb{N}$. As usual we use the notation $\InfCat$ to
  denote the category of all \inf-categories and \inf-functors between
  them.
\end{defn}

\begin{obs}[$n$-cells and $n$-categories]
  \label{ncells.ncats}
  A cell in an \inf-category $\mathbb{C}$ is said to be an {\em
    $n$-cell\/} (for $n\in\mathbb{N}$) if it is an identity for the
  category structure $\triple<\comp_n;\s_n;\t_n>$.  Furthermore, we
  say that our cell is a {\em non-trivial $n$-cell\/} if it is an
  $n$-cell but it isn't an $(n-1)$-cell.  We will generally use the
  notation $\mathbb{C}_n\subseteq \mathbb{C}$ to denote the set of
  $n$-cells in $\mathbb{C}$.  Condition~(\ref{2cat.cond}) of
  definition~\ref{omegacat} implies that if $m\geq n$ then any
  $n$-cell is also an $m$-cell (that is $\mathbb{C}_n\subseteq
  \mathbb{C}_m$).

  An {\em $n$-category\/} is an \inf-category in which every cell is
  an $n$-cell and $\nCat$ denotes the full subcategory of
  $n$-categories in $\InfCat$. From the previous paragraph, we know
  that if $m\geq n$ then the category structure
  $\triple<\comp_m;\s_m;\t_m>$ is {\em discrete}, in the sense that
  under $\comp_m$ cells may only compose with themselves. It follows
  that our latter definition of 2-category is entirely consistent with
  the one which we originally gave in definition~\ref{2cat.qua.doub},
  since all of its category structures aside from $\comp_0$ and
  $\comp_1$ are completely trivial.
  
  Observe that the sources, targets and composites of $n$-cells are
  again $n$-cells, thus it follows that the \inf-category structure of
  $\mathbb{C}$ restricts to provide us with an $n$-category structure
  on $\mathbb{C}_n$. Following Street~\cite{Street:1987:Oriental}, we
  adopt the notation $\Sup_n(\mathbb{C})$ for this sub-$n$-category
  which we'll refer to as the {\em $n$-superstructure of
    $\mathbb{C}$}.
  
  Of course, since an \inf-functor $\arrow f:\mathbb{C}->\mathbb{D}.$
  preserves sources and targets it is clear that it maps $n$-cells to
  $n$-cells and thus restricts to an \inf-functor $\arrow f:
  \Sup_n(\mathbb{C})->\Sup_n(\mathbb{D}).$. This makes $\Sup_n(\cdot)$
  into the underlying functor of an idempotent comonad on $\InfCat$,
  the full image of which is the subcategory of $n$-categories.
  
  The contrapositive of this ``dimension preservation'' property of
  the \inf-functor $f$, the observation that if $f(c)$ is not an
  $n$-cell in $\mathbb{D}$ then $c$ itself is not an $n$-cell in
  $\mathbb{C}$, is also, somewhat surprisingly, of use in its own
  right. In particular, it is often used to demonstrate the
  non-triviality of an $n$-cell in $\mathbb{C}$ by constructing a
  functor with that domain which maps our cell to a demonstrably
  non-trivial $n$-cell in its codomain.
  
  The $n$-superstructure functor has a right adjoint
  $\arrow\Cosup_n:\InfCat->\InfCat.$, we call the \inf-category
  $\Cosup_n(\mathbb{C})$ the {\em $n$-cosuperstructure\/} of the
  \inf-category $\mathbb{C}$. This may be constructed as the
  $(n+1)$-category whose cells are pairs $\pair<c_s;c_t>$ of $n$-cells
  in $\mathbb{C}$ satisfying the ``globularity'' condition that
  $\s_{n-1}(c_s)=\s_{n-1}(c_t)$ and $\t_{n-1}(c_s)=\t_{n-1}(c_t)$.
  The category structure $\triple<\comp_m;\s_m;\t_m>$ on
  $\Cosup_n(\mathbb{C})$ is given point-wise on components when $m<n$,
  is discrete when $m>n$ and is defined by
  \begin{equation*}
    \s_n\pair<c_s;c_t>=\pair<c_s;c_s>\mkern60mu
    \t_n\pair<c_s;c_t>=\pair<c_t;c_t>\mkern60mu
    \pair<c'_s;c'_t>\comp_m\pair<c_s;c_t>=\pair<c'_s;c_t>
  \end{equation*}
  when $m=n$. By construction, it is clear that we have a canonical
  isomorphism of \inf-categories $\overarr:\Sup_n(\mathbb{C})->
  \Sup_n({\Cosup_n(\mathbb{C})}).$, which carries an $n$-cell
  $c\in\mathbb{C}$ to the $n$-cell $\pair<c;c>\in\Cosup_n(\mathbb{C})$.
  Furthermore, the remaining structure of $\Cosup_n(\mathbb{C})$ is
  completely determined by the fact that each pair of its $n$-cells
  which satisfy the globularity condition required by the $n$-source
  and $n$-target of an $(n+1)$-cell do indeed bound a unique such
  cell.  The adjoint transpose $\arrow \hat{f}:\mathbb{C}->\Cosup_n(
  \mathbb{D}).$ of an \inf-functor $\arrow f:\Sup_n(\mathbb{C})->
  \mathbb{D}.$ is given by $\hat{f}(c)=\pair<f(\s_n(c));f(\t_n(c))>$,
  which is well defined since each $c\in\mathbb{C}$ has $\s_n(c),
  \t_n(c)\in\Sup_n(\mathbb{C})$ and is easily seen to be appropriately
  functorial.
  
  As a left adjoint $\Sup_n(\cdot)$ preserves all (small) colimits in
  $\InfCat$.  On the other hand, (small) limits in $\InfCat$ are
  constructed by taking the corresponding limit of underlying sets of
  cells in $\Set$ and defining category structures on this set
  point-wise, therefore it is also clear that $\Sup_n(\cdot)$
  preserves all (small) limits. Of course, these preservation
  properties imply that $\Sup_n(\cdot)$ preserves any structure or
  property which can be described in terms of limits and colimits,
  such as (regular) monomorphisms and (regular) epimorphisms.
\end{obs}

\begin{obs}[\inf-categories as enriched categories]\label{infcat.enr}
  Recall, from Street~\cite{Street:1987:Oriental}, that the category
  $\InfCat$ is cartesian closed and equivalent to the category
  $\InfCatCat$ of \inf-category enriched categories (in the sense of
  Kelly~\cite{Kelly:1982:ECT}).  The functor
  $\equiv\pathecat:\InfCat-> \InfCatCat.$ carries an \inf-category
  $\mathbb{C}$ to the \inf-category enriched category
  $\pathecat(\mathbb{C})$ with:
  \begin{itemize}
  \item objects $u$ which are the 0-cells of $\mathbb{C}$,
  \item hom-\inf-categories $\pathecat(\mathbb{C})(u,v)$ with set of
    cells $\{c\in\mathbb{C}\mid \s_0(c)=u\text{ and }\t_0(c)=v\}$ and
    compositions which are restrictions of the compositions
    $\comp_1,\comp_2,...$ of $\mathbb{C}$ to this set,
  \item composition \inf-functor $\arrow\circ:\pathecat(\mathbb{C})(v,w)
    \times\pathecat(\mathbb{C})(u,v)->\pathecat(\mathbb{C})(u,w).$
    given by restricting the composition $\comp_0$ of $\mathbb{C}$,
    and corresponding identities $\id_u=\idt_0(u)\in\pathecat(\mathbb{C})(u,u)$.
  \end{itemize}
  Notice also that this functor restricts to an equivalence
  $\equiv\pathecat:\ECat{(n+1)}->\ECat{(\nCat)}.$ for each
  $n\in\mathbb{N}$.

  It is sometimes instructive to think of \inf-categories as being
  {\em oriented combinatorial CW-complexes\/} or {\em globular
    spaces\/} and under this analogy it is clearly natural to think of
  this functor as being a form of {\em path category
  construction}. 
\end{obs}



\section{An Introduction to the Decalage Construction}
\label{decintro.sect}

\subsection{Nerves and Decalage}\label{nerves.dec}

First we introduce the nerve construction and the decalage comonad.
We also make precise the sense in which the latter is generic for the
former.

\begin{obs}[comonads and left actions of $\op\aDelta$]
  \label{comonad.lact} First recall, from
  say~\cite{BarrWells:1985:TTT}, that a {\em comonad\/}
  $\triple<\gfunc;\nu;\varepsilon>$ on a category $\mathcal{C}$ is
  no-more-nor-less than a monoid in the strict monoidal category
  $\triple<\op{\funcat[\mathcal{C},\mathcal{C}]}; \circ;
  \id_{\mathcal{C}}>$, where $\circ$ denotes the {\em horizontal\/}
  composition of the 2-category $\Category$ and $\id_{\mathcal{C}}$ denotes
  the identity functor on $\mathcal{C}$. So we know, by
  lemma~\ref{univ.delta}, that such a comonad corresponds to a unique,
  strictly monoidal, functor:
  \begin{displaymath}
    \xymatrix@R=2ex@C=10em{
      {\triple<\op\aDelta;\oplus;[-1]>}\ar[r]^{\tilde{\gfunc}} & 
      {\triple<\funcat[\mathcal{C},\mathcal{C}];\circ;\id_{\mathcal{C}}>}}
  \end{displaymath}
  However, just as in classical monoid theory, any such representation
  corresponds, under the adjunction
  $-\times\mathcal{C}\dashv\funcat[\mathcal{C},\ast]$, to a unique
  left action of the monoid $\op\aDelta$ on the category
  $\mathcal{C}$. That is to say, it uniquely determines 
  a functor
  \begin{displaymath}
    \let\labelstyle=\textstyle
    \xymatrix@R=2ex@C=10em{
      {\op\aDelta\times\mathcal{C}}\ar[r]^{\decact_\gfunc} & {\mathcal{C}}}
  \end{displaymath}
  which we shall usually write in infix form.  This may be given
  explicitly by $\alpha\decact_\gfunc f \defeq (\tilde{\gfunc}\alpha)f$, for all
  simplicial operators $\alpha\in\op\aDelta$ and arrows
  $f\in\mathcal{C}$, and it satisfies the usual axioms for a monoid
  left action of the arrows of $\op\aDelta$ on those of $\mathcal{C}$:
  \begin{displaymath}
    \begin{array}{rcll}
      \id_{[-1]}\decact_\gfunc f & = & f & \text{for all arrows
        $f\in\mathcal{C}$,} \\
      \beta\decact_\gfunc (\alpha\decact_\gfunc f) & = &
      (\beta\oplus\alpha)\decact_\gfunc f \mkern5mu & 
      \text{for all arrows
        $f\in\mathcal{C}$ and $\alpha,\beta\in\op\aDelta$.}
    \end{array}
  \end{displaymath}
  Classically, it is also common to go a step further and transpose
  such actions once more, this time under the adjunction
  $\op\aDelta\times\mathord{-}\dashv\funcat[\op\aDelta, \ast]$, to give
  a functor
  \begin{equation}\label{trans.action}
    \let\labelstyle=\textstyle
    \xymatrix@R=2ex@C=10em{
      {\mathcal{C}}\ar[r]^{\hat{\gfunc}} &       
      {\funcat[\op\aDelta,\mathcal{C}]}
    }
  \end{equation}
  satisfying rules which are transposes of those satisfied by
  $\decact_\gfunc$. Explicitly, this latter functor maps an object
  $C\in\mathcal{C}$ to a functor $\hat{\gfunc}(C)\in\funcat[\op\aDelta,
  \mathcal{C}]$, which we might call an {\em augmented simplicial
    object in $\mathcal{C}$}, determined by 
    \begin{gather*}
      \hat{\gfunc}(C)([n]) = \gfunc^{n+1}(C) \\ 
      \hat{\gfunc}(C)(\face^n_i)
      =\gfunc^{n-i}(\varepsilon_{\gfunc^i(C)})
      \mkern80mu
      \hat{\gfunc}(C)(\degen^n_i) = \gfunc^{n-i}(\nu_{\gfunc^i(C)})
    \end{gather*}
    and maps an arrow $\arrow f: C-> D.\in\mathcal{C}$ to an {\em
      augmented simplicial map\/} $\hat{\gfunc}(f)$ with components:
  \begin{eqnarray*}
    \hat{\gfunc}(f)_n & = & \gfunc^{n+1}(f)
  \end{eqnarray*}
\end{obs}
  
\begin{defn}[nerves of comonads]\label{defn.comonad.nerve} 
  Suppose that we are given a category $\mathcal{C}$ equipped with a
  functor $\arrow\cpt:\mathcal{C}->\Set.$ and a comonad
  $\triple<\gfunc;\nu;\varepsilon>$ on $\mathcal{C}$. As in
  observation~\ref{comonad.lact} we may transpose the
  $\op\aDelta$-action corresponding our comonad to give the functor
  $\hat{\gfunc}$ of display~(\ref{trans.action}). This we compose with
  the functor
  \begin{displaymath}
    \let\labelstyle=\textstyle
    \xymatrix@R=2ex@C=10em{
      {\funcat[\op\aDelta,\mathcal{C}]}
      \ar[r]^{\funcat[i,\cpt]} &       
      {\funcat[\op\Delta,\Set]
        \hbox to 0pt{$\mathsurround=0pt {}\defeq\Simp$\hss}}
    }
  \end{displaymath}
  which post-composes $\cpt$ onto each $\arrow
  \ffunc:\op\aDelta->\mathcal{C}.$ and restricts the domain of the result
  to $\op\Delta$ by pre-composing it with the inclusion $\inc
  i:\op\Delta->\op\aDelta.$. We call this composite
  $\nerv_{\gfunc,\cpt}\defeq\funcat[i,\cpt]\circ\hat{\gfunc}$ the {\em nerve functor
    associated with $\triple<\gfunc;\nu;\varepsilon>$ under $\arrow\cpt:
    \mathcal{C}->\Set.$}.
\end{defn}

\begin{defn}[the decalage construction]\label{defn.dec}
  Using the properties of $\aDelta$ as a strict monoidal category and
  the fact that $\Simp\defeq\funcat[\op\Delta,\Set]$, we may define a
  canonical left action of $\op\aDelta$ on $\Simp$ as follows:
  \begin{equation}\label{dec.lact}
    \let\labelstyle=\textstyle
    \xymatrix@R=1.8ex@C=10em{
      {\op\aDelta\times\Simp} \ar[r]^{\decact} & {\Simp} \\
      {([n],X)}\ar@{|->}[r]\ar@<1.5ex>[dd]^{\textstyle f}="one" & 
      {X\circ(\mathord{-}\dsum[n])} 
      \ar[dd]^{\textstyle f\circ(\mathord{-}\dsum\alpha)}_{}="two" \\
      {\hbox{\vrule width0pt height2ex}} & \\
      {([m],Y)}\ar@<1.5ex>[uu]^{\textstyle\alpha}_{}
      \ar@{|->}[r] & {Y\circ(\mathord{-}\dsum[m])}
      \ar@{|->}"one";"two"}
  \end{equation}
  Showing that this is functorial and that it satisfies the axioms
  required of a left action is a matter of the routine application of
  the properties of $\dsum$ as the tensor of a strict monoidal
  category and of $\circ$ as the horizontal composition of the
  2-category $\Category$.  Appealing to
  observation~\ref{comonad.lact}, we know that this left action
  corresponds to a unique comonad $\triple<\Dec;\cdegen;\cface>$ on
  $\Simp$, called the {\em decalage comonad}, which has:
  \begin{itemize}
  \item {\bf underlying functor} given by the ``section'': $$\arrow\Dec\defeq
    ([0]\decact\mathord{-}):\Simp->\Simp.$$
  \item {\bf comultiplication} the natural transformation:
    $$\arrow\cdegen\defeq(\degen^0_0\decact\mathord{-}):\Dec=
    ([0]\decact\mathord{-})->([1]\decact\mathord{-})=\Dec\circ\Dec.$$
  \item {\bf counit} the natural transformation:
    $$\arrow\cface\defeq(\face^0_0\decact\mathord{-}):
    \Dec=([0]\decact\mathord{-})->([-1]\decact\mathord{-})=
    \id_{\Simp}.$$
  \end{itemize}
\end{defn}

In a precise sense, the pair consisting of the decalage comonad and
the connected components functor constitutes an ``identity'' for the
nerve construction:

\begin{lemma}[the nerve functor associated with the decalage comonad]%
  \label{nerve.dec}%
  The nerve functor $\Ndec$ associated with the decalage comonad
  $\triple<\Dec;\cdegen;\cface>$ under the connected components
  functor $\arrow\cpt_0:\Simp->\Set.$ of observation~\ref{conn.cpt.simp} is
  (isomorphic to) the identity functor on $\Simp$.
\end{lemma}

\begin{proof}
  We know, from the constructions in
  definitions~\ref{defn.comonad.nerve} and~\ref{defn.dec}, that if
  $X$ is a simplicial set and $n\in\mathbb{N}$ then
  $\Ndec(X)_{n}\defeq\cpt_0([n]\decact X)$, where $\decact$ is the
  canonical left action defined in display~(\ref{dec.lact}).  Of
  course we also know, from observation~\ref{conn.cpt.simp}, that
  $\cpt_0([n]\decact X)$ may be constructed as the coequaliser of the
  pair
  \begin{displaymath}
    \let\labelstyle=\textstyle
    \xymatrix@R=2ex@C=8em{
      {([n]\decact X)_{1}}
      \ar@<1ex>[r]^{{-}\cdot\face^1_0}
      \ar@<-1ex>[r]_{{-}\cdot\face^1_1} & 
      {([n]\decact X)_{0}} 
    }
  \end{displaymath}
  which is easily seen, from the definition of $[n]\decact X$ in
  display~(\ref{dec.lact}), to be equal to the pair:
  \begin{displaymath}
    \let\labelstyle=\textstyle
    \xymatrix@R=2ex@C=8em{
      {X_{n+2}}
      \ar@<1ex>[r]^{{-}\cdot\face^{n+2}_{n+1}}
      \ar@<-1ex>[r]_{{-}\cdot\face^{n+2}_{n+2}} & 
      {X_{n+1}} 
    }
  \end{displaymath}
  Notice now that a direct application of the simplicial identities of
  observation~\ref{simplicial.idents} demonstrates that both squares
  in the following diagram of simplicial operators commute
  \begin{displaymath}
    \xymatrix@R=2em@C=3em{
      {[n]}\ar[r]^{\face^{n+1}_{n+1}}\ar[d]_{\face^{n+1}_{n+1}} &
      {[n+1]}\ar[r]^{\degen^n_n}\ar[d]^{\face^{n+2}_{n+2}} &
      {[n]}\ar[d]^{\face^{n+1}_{n+1}} \\
      {[n+1]}\ar[r]_{\face^{n+2}_{n+1}} &
      {[n+2]}\ar[r]_{\degen^{n+1}_{n}} & {[n+1]}
    }
  \end{displaymath}
  and that the composites of the upper and lower horizontal pairs of
  operators evaluate to identities. Applying the simplicial set $X$ to
  this diagram we get a corresponding diagram in $\Set$ satisfying
  these commutativity conditions. This demonstrates that the following
  diagram satisfies the conditions required of it to be a {\em
    split coequaliser\/} (cf.~\cite{BarrWells:1985:TTT})
  \begin{displaymath}
    \xymatrix@R=2ex@C=8em{
      {X_{n+2}}
      \ar@<-2ex>[r]_{{-}\cdot\face^{n+2}_{n+2}}
      \ar[r]|{\mskip4mu{-}\cdot\face^{n+2}_{n+1}\mskip4mu} & 
      {X_{n+1}} 
      \ar@<-2ex>[l]_{{-}\cdot\degen^{n+1}_n}
      \ar@<-1ex>[r]_{{-}\cdot\face^{n+1}_{n+1}} &
      {X_{n}}
      \ar@<-1ex>[l]_{{-}\cdot\degen^n_n}
    }
  \end{displaymath}
  from which it follows that $\Ndec(X)_{n}\defeq\cpt_0([n]\decact
  X)\cong X_{n}$. It is now a routine matter to demonstrate that
  these isomorphisms are natural in $[n]\in\op\Delta$, making them the
  components of a simplicial isomorphism $\Ndec(X)\cong X$, and in
  $X\in\Simp$, thus making the whole collection into a natural
  isomorphism $\Ndec\cong\id_{\Simp}$ as required.
\end{proof}

\subsection{Comonad Transformations and Simplicial Reconstruction}

\begin{obs}[the 2-category $\Cylinder$ of cylinders]
  As a step towards studying maps of comonads, we may define a certain
  2-category $\Cylinder$ whose 0-cells are arbitrary functors and
  whose hom-category $\cyl[\ffunc,\gfunc]$, from one functor $\arrow
  \ffunc:\mathcal{C}->\mathcal{C}'.$ to another $\arrow
  \gfunc:\mathcal{D}->\mathcal{D}'.$, is defined to be the {\bf comma}
  category $\funcat[\mathcal{C},\gfunc]\downarrow
  \funcat[\ffunc,\mathcal{D}']$ (cf.~\cite{BarrWells:1985:TTT}
  or~\cite{Maclane:1971:CWM}) from the post-composite
  $\arrow\funcat[\mathcal{C},\gfunc]=\gfunc\circ{-}:
  \funcat[\mathcal{C},\mathcal{D}]-> \funcat[\mathcal{C},
  \mathcal{D}'].$ to the pre-composite
  $\arrow\funcat[\ffunc, \mathcal{D}']={-}\circ
  \ffunc:\funcat[\mathcal{C'},\mathcal{D}']->
  \funcat[\mathcal{C},\mathcal{D}'].$.  
  
  More explicitly, this hom-category has:
  \begin{itemize}
  \item {\bf objects (1-cells)}, called {\em squares},
    $\triple<\kfunc;\kfunc';\kappa>$ which consist of a pair of
    functors $\arrow \kfunc:\mathcal{C}->\mathcal{D}.$ and $\arrow
    \kfunc':\mathcal{C}'->\mathcal{D}'.$ and a natural transformation
    $\twocell \kappa:\gfunc\circ \kfunc->\kfunc'\circ \ffunc.$; as in
    observation~\ref{dcat.psq} this data is depicted as a {\em pasting
      square\/}
    \begin{displaymath}
      \let\labelstyle=\textstyle
      \xymatrix@=2.5em{
        {\mathcal{C}}\ar[d]_{\ffunc}\ar[r]^{\kfunc} &
        {\mathcal{D}\save []-<1.25em,1.25em>="three"\restore}
        \ar[d]^{\gfunc} \\
        {\mathcal{C}'\save []+<1.25em,1.25em>="one"\restore}
        \ar[r]_{\kfunc'} & 
        {\mathcal{D}'}
        \ar@{=>}_{\textstyle\kappa} "three"; "one"
      }
    \end{displaymath}
  \item {\bf arrows (2-cells)}, called {\em cylinders},
    $\arrow\pair<\gamma;\gamma'>:
    \triple<\kfunc;\kfunc';\kappa>->\triple<\bar{\kfunc};\bar{\kfunc}';\bar\kappa>.$
    consisting of a pair of natural transformations
    $\twocell\gamma:\kfunc->\bar{\kfunc}.$ and $\twocell\gamma':\kfunc'->\bar{\kfunc}'.$
    which satisfy the equation
    $\bar\kappa\cdot(\gfunc\gamma)=(\gamma' \ffunc)\cdot\kappa$, and
  \item {\bf vertical composition} of cylinders calculated by applying
    vertical composition point-wise to the components of each
    cylinder.
  \end{itemize}
  Just as in observation~\ref{dcat.psq} the horizontal composite of a
  pair of 1-cells in $\Cylinder$ is constructed by horizontal
  ``pasting'' of squares or more explicitly:
  $$\triple<\lfunc;\lfunc';\lambda>\circ \triple<\kfunc;\kfunc'; \kappa> \defeq
  \triple<\kfunc\circ \lfunc; \kfunc'\circ \lfunc'; (\lfunc'\kappa)\cdot(\lambda \kfunc)>$$
  Finally, horizontal composition of cylinders is a trivial matter of
  point-wise composition of their underlying natural transformations.
  
  Of course, pasting squares are familiar to us from our
  considerations of double categories in
  subsection~\ref{sect.dcat.with.conn}.  However, it is worth
  stressing that in the earlier context our pasting squares comprised
  dimension 2 of our double categories whereas here we've shifted
  dimensions and they now inhabit dimension 1 of $\Cylinder$.
\end{obs}

\begin{obs}[cotensors and tensors in $\Cylinder$]\label{cyl.cotens}
  In \cite{Kelly:1982:ECT} Kelly discusses a primitive kind of
  enriched weighted colimit construction, called a {\em cotensor},
  which generalises the notion of an iterated product in un-enriched
  category theory. In the 2-categorical case, if $\twocat{A}$ is a
  2-category, $\mathcal{E}$ is a category in $\Category$ and $A'$ is a
  0-cell in $\twocat{A}$ then the cotensor $\mathcal{E}\cotens
  A'\in\twocat{A}$ of $A'$ by $\mathcal{C}$, if it exists, is the
  (essentially) unique 0-cell for which there is an isomorphism of
  categories
  \begin{displaymath}
    \funcat[\mathcal{E},\twocat{A}(A,A')] \cong 
    \twocat{A}\pair<A; \mathcal{E}\cotens A'>
  \end{displaymath}
  which is 2-natural in $A\in\twocat{A}$. For instance, it is a
  trivial consequence of the way the we defined the enrichment of
  $\Category$ over itself that the cotensor of a category
  $\mathcal{D}\in\Category$ by $\mathcal{E}$ is simply (isomorphic to)
  the functor category $\funcat[\mathcal{E},\mathcal{C}]$.
  
  We might guess that the cotensor of a 0-cell $\arrow \gfunc:\mathcal{D}->
  \mathcal{D}'.\in\Cylinder$ in our 2-category of cylinders by the
  category $\mathcal{E}$ is simply the post-composition functor
  $\arrow \funcat[\mathcal{E},\gfunc]: \funcat[\mathcal{E},\mathcal{D}]->
  \funcat[\mathcal{E},\mathcal{D}'].$. To prove this first recall, by
  the definition of $\cyl[\ffunc,\gfunc]$ as a comma category, that we have a
  square
  \begin{equation}
    \let\labelstyle=\textstyle
    \xymatrix@R=2.25em@C=6em{
      & {\cyl[\ffunc,\gfunc]}\ar[r]^{\Pi}\ar[d]_{\Pi'} &
      {\funcat[\mathcal{C},\mathcal{D}]}
      \ar[d]^{\funcat[\mathcal{C},\gfunc]}_{}="two" & \\
      & {\funcat[\mathcal{C}',\mathcal{D}']}
      \ar[r]_{\funcat[\ffunc,\mathcal{D}']}^{}="one" &
      {\funcat[\mathcal{C},\mathcal{D}']}
      \ar@/^1.5ex/@<1.5ex>@{<=}^{\pi} "one";"two" &
    }
  \end{equation}
  in which $\Pi$ (resp.\ $\Pi'$) is the functor which projects cylinders
  onto their first (resp. second) component and $\pi$ is the natural
  transformation with $\pi_{\triple<\kfunc;\kfunc';\kappa>} \defeq \kappa$ for
  each square $\arrow\triple<\kfunc;\kfunc';\kappa>:\ffunc->\gfunc.$. This square is, in
  fact, 2-universal (in the strict enriched sense of
  Kelly~\cite{Kelly:1982:ECT}) amongst all 2-cones of this form over the
  diagram consisting of the functors $\funcat[\mathcal{C},\gfunc]$ and
  $\funcat[\ffunc,\mathcal{D}']$. Squares with this property are usually
  called {\em comma squares\/} by 2-categorists.
  
  Taking two copies of this comma square, one for $\cyl[\ffunc,\gfunc]$ (to
  which we apply the 2-functor $\funcat[\mathcal{E},\ast]$) and the
  other for $\cyl[\ffunc,{\funcat[\mathcal{E},\gfunc]}]$, we get a diagram
  \begin{equation}\label{cotens.big.diag}
    \xy
      \xymatrix"*"@R=2.5em@C=4em{
        {\cyl[\ffunc,{\funcat[\mathcal{E},\gfunc]}]}\ar[r]^{\Pi}
        \ar[d]_(0.6){\Pi'} &
        {\funcat[\mathcal{C},{\funcat[\mathcal{E},\mathcal{D}]}]}
        \ar[d]^{\funcat[\mathcal{C},{\funcat[\mathcal{E},\gfunc]}]}_{}="two" \\
        {\funcat[\mathcal{C}',{\funcat[\mathcal{E},\mathcal{D}']}]}
        \ar[r]_{\funcat[\ffunc,{\funcat[\mathcal{E},\mathcal{D}']}]}^{}="one" &
        {\funcat[\mathcal{C},{\funcat[\mathcal{E},\mathcal{D}']}]}
        \ar@/^1.5ex/@<1.5ex>@{<=}^{\pi} "one";"two" 
      }
      \POS+<-18em,4.5em>
      \xymatrix@R=2.5em@C=4em{
        {\funcat[\mathcal{E},{\cyl[\ffunc,\gfunc]}]}
        \ar[r]^{\funcat[\mathcal{E},\Pi]}
        \ar[d]_{\funcat[\mathcal{E},\Pi']} 
        \ar@{.>}["*"] &
        {\funcat[\mathcal{E},{\funcat[\mathcal{C},\mathcal{D}]}]}
        \ar@/^0.1pt/ [d] |!{"1,1";"*1,1"}{\hole} _{}="four"%
          ^(0.4){\funcat[\mathcal{E},{\funcat[\mathcal{C},\gfunc]}]}\relax
        \ar@{-->}["*"]^{\scriptstyle\cong} \\
        {\funcat[\mathcal{E},{\funcat[\mathcal{C}',\mathcal{D}']}]}
        \ar[r]_(0.6){\funcat[\mathcal{E},{\funcat[\ffunc,\mathcal{D}']}]}^{}="three" 
        \ar@{-->}["*"]_{\scriptstyle\cong} &
        {\funcat[\mathcal{E},{\funcat[\mathcal{C},\mathcal{D}']}]}
        \ar@{-->}["*"]^{\scriptstyle\cong}|!{"*1,1";"*2,1"}\hole
        \ar@/^1.5ex/@<1.5ex>@{<=}^(.3){\funcat[\mathcal{E},\pi]} "three";"four" 
      }
    \endxy
  \end{equation}
  in which the dashed diagonal arrows are canonical isomorphisms
  derived from the cartesian closed structure on $\Category$ and the
  skewed rectangles that they bound commute by the naturality of those
  constructions. The existence of the last of these diagonal arrows,
  depicted with a dotted line, follows from the 2-universal property
  of the comma square at the bottom right-hand of the diagram and it
  is the unique arrow making the remaining skewed rectangles commute.
  However, the square in the upper left hand of the diagram is also a
  comma square since it is constructed by applying
  $\funcat[\mathcal{E},\ast]$, which as a right 2-adjoint preserves
  all 2-limits, to the comma square associated with $\cyl[\ffunc,\gfunc]$ and so
  it follows that our dotted arrow is also an isomorphism as required.
  Proving the appropriate 2-naturality of this isomorphism in $\ffunc$, $\gfunc$
  and $\mathcal{E}$ is now a matter of routine verification.
  
  For the sake of symmetry, we might also mention the dual concept
  which generalises iterated sums. The {\em tensor} of a 0-cell $A$ in
  a 2-category $\twocat{A}$ by a category $\mathcal{E}$, if it exists,
  is the (essentially) unique 0-cell $\mathcal{E}\circledast A$ for
  which there is an isomorphism of categories
  \begin{displaymath}
    \funcat[\mathcal{E},\twocat{A}(A,A')] \cong
    \twocat{A}\pair<\mathcal{E}\circledast A;A'>
  \end{displaymath}
  which is 2-natural in $A'\in\twocat{A}$. Again, it is practically a
  tautology that the category $\mathcal{E}\times\mathcal{C}$ is the
  tensor in $\Category$ of $\mathcal{C}$ by $\mathcal{E}$.  It is
  also the case that the tensor of a 0-cell $\arrow \ffunc:
  \mathcal{C}->\mathcal{C}'.\in\Cylinder$ in our 2-category of
  cylinders by the category $\mathcal{E}$ is simply the functor
  $\arrow\mathcal{E}\times \ffunc: \mathcal{E}\times\mathcal{C}->
  \mathcal{E}\times\mathcal{C}'.$, which applies $\ffunc$ ``in the second
  variable''. A proof of this fact may be constructed along the
  lines of that given for cotensors in $\Cylinder$ above and is left
  as an exercise for the reader. \qed
\end{obs}

\begin{obs}[comonad transformations]\label{comonad.trans}  
  Since the comonad notion is entirely equational, we may follow
  Street~\cite{Street:1972:Monads} and interpret it in any 2-category;
  in particular we have an interest in studying comonads in
  $\Cylinder$. If $\arrow \ffunc: \mathcal{C}->\mathcal{C}'.$ is a functor
  then a comonad on the 0-cell $\ffunc$ in $\Cylinder$ is simply a monoid
  in the (dual) endo-category $\op{\cyl[\ffunc,\ffunc]}$, which is a strict
  monoidal category under the horizontal composition of $\Cylinder$ in
  the usual way. Now, since horizontal and vertical composition of
  cylinders is calculated point-wise, we know that the function which
  projects cylinders onto their first (respectively second) component
  provides us with a 2-functor $\arrow \Pi:\Cylinder->\Category.$
  (resp.\ $\arrow \Pi':\Cylinder->\Category.$). It follows that this
  2-functor restricts to a strict monoidal functor from
  $\op{\cyl[\ffunc,\ffunc]}$ to $\op{\funcat[\mathcal{C},\mathcal{C}]}$ (resp.\ 
  $\op{\funcat[\mathcal{C}',\mathcal{C}']}$) and, consequently, that
  this carries our monoid to a comonad $\triple<\gfunc;\nu;\varepsilon>$ on
  $\mathcal{C}$ (resp.\ $\triple<\gfunc';\nu';\varepsilon'>$ on
  $\mathcal{C}'$).  The remaining data encapsulated in our monoid in
  $\cyl[\ffunc,\ffunc]$ simply consists of a natural transformation
  $\twocell\gamma:\ffunc\circ \gfunc-> \gfunc'\circ \ffunc.$ making $\triple<\gfunc;\gfunc';\gamma>$
  into an endo-square (endo-1-cell in $\Cylinder$) on $\ffunc$ and
  satisfying the cylinder conditions:
  \begin{itemize}
  \item the pair $\pair<\nu;\nu'>$ is a cylinder from
    $\triple<\gfunc;\gfunc';\gamma>$ to $\triple<\gfunc;\gfunc';\gamma>\circ
    \triple<\gfunc;\gfunc';\gamma>$, in other words we must have
    $(\gfunc'\gamma)\cdot(\gamma \gfunc)\cdot(\ffunc\nu)=(\nu' \ffunc)\cdot\gamma$,
    and
  \item the pair $\pair<\varepsilon;\varepsilon'>$ is a cylinder
    from $\triple<\gfunc;\gfunc';\gamma>$ to the identity square
    $\triple<\id_{\mathcal{C}};\id_{\mathcal{C}};\id_\ffunc>$, equivalently
    we require that $\ffunc\varepsilon = (\varepsilon' \ffunc)\cdot\gamma$.
  \end{itemize}
  The pair $\pair<\ffunc;\gamma>$ of a functor $\ffunc$ and a natural
  transformation $\gamma$ satisfying these conditions is known as a
  {\em comonad transformation\/} from
  $\triple<\gfunc;\nu;\varepsilon>$ to
  $\triple<\gfunc';\nu';\varepsilon'>$. Such transformations are said to
  be {\em strong\/} if $\gamma$ is invertible.
  
  Of course we know, from lemma~\ref{univ.delta}, that our comonad
  transformation $\triple<\gfunc;\gfunc';\gamma>$, as a monoid in
  $\op{\cyl[\ffunc,\ffunc]}$, corresponds to a unique strict monoidal
  functor
  \begin{displaymath}
    \let\labelstyle=\textstyle
    \xymatrix@R=2ex@C=10em{
      {\triple<\op\aDelta;\oplus;[-1]>}
      \ar[r]^<>(0.5){\suptilde{\triple<\gfunc;\gfunc';\gamma>}} &
      {\triple<\cyl[\ffunc,\ffunc];\circ;\triple<\id_{\mathcal{C}};
        \id_{\mathcal{C}'};\id_\ffunc>>}
    }
  \end{displaymath}
  and by applying the isomorphism $ \funcat[\op\aDelta,{\cyl[\ffunc,\ffunc]}]
  \cong \cyl[\ffunc,{\funcat[\op\aDelta,\ffunc]}]$ of
  observation~\ref{cyl.cotens} it follows that this, in turn,
  corresponds to a square
  \begin{displaymath}
    \let\labelstyle=\textstyle
    \xymatrix@R=2ex@C=10em{
      {\ffunc}\ar[r]^<>(0.5){\suphat{\triple<\gfunc;\gfunc';\gamma>}} &
      {\funcat[\op\aDelta,\ffunc]}
    }
  \end{displaymath}
  which, just as in observation~\ref{comonad.lact}, satisfies rules
  that are transposes of those for a left action of $\op\aDelta$ on
  $\ffunc$. Tracing through the constructions of lemma~\ref{univ.delta} and
  observations~\ref{comonad.lact} and~\ref{cyl.cotens} it is easily
  seen that this square is of the form
  \begin{equation}\label{square.trans}
    \let\labelstyle=\textstyle
    \xymatrix@R=2.5em@C=8em{
      {\mathcal{C}}\ar[r]^{\hat{\gfunc}}_{}="two"\ar[d]_{\ffunc} &
      {\funcat[\op\aDelta,\mathcal{C}]}
      \ar[d]^{\funcat[\op\aDelta,\ffunc]} \\
      {\mathcal{C}'}\ar[r]_{\hat{\gfunc}'}^{}="one" &
      {\funcat[\op\aDelta,\mathcal{C}']}    
      \ar@{<=} "one"+<0em,1.2em>;"two"-<0em,1.2em> ^{\hat\gamma}
    }
  \end{equation}
  in which the horizontal functors $\hat{\gfunc}$ and $\hat{\gfunc}'$ are those
  which correspond to the comonads $\triple<\gfunc;\nu;\varepsilon>$ and
  $\triple<\gfunc';\nu';\varepsilon'>$ (respectively) as in
  observation~\ref{comonad.lact}. At the same time it is easily
  demonstrated that, for objects $C\in\mathcal{C}$ and
  $[n]\in\op\aDelta$, the component $(\hat\gamma_C)_{n}$  of the
  natural transformation in square~(\ref{square.trans}) is equal to
  $\gamma^{(n+1)}_C$, where $\gamma^{(n+1)}$ is the natural
  transformation in the square:
  \begin{displaymath}
    \triple<\gfunc^{n+1};(\gfunc')^{n+1};\gamma^{(n+1)}> \defeq 
    \suptilde{\triple<\gfunc;\gfunc';\gamma>}([n])
  \end{displaymath}
  Since this latter square is simply equal to the $(n+1)$-fold
  iterated horizontal power of the square $\triple<\gfunc;\gfunc';\gamma>$ in
  $\Cylinder$, it follows that our original transformation of comonads
  $\pair<\ffunc;\gamma>$ is {\bf strong} if and only if the corresponding
  natural transformation $\hat\gamma$ of display~(\ref{square.trans})
  is an isomorphism. \qed
\end{obs}  

\begin{lemma}\label{nerves.and.comonad.tr}
  Suppose that we are given a categories $\mathcal{C}$ and
  $\mathcal{C}'$, along with
  \begin{itemize}
  \item comonads $\triple<\gfunc;\nu;\varepsilon>$ on $\mathcal{C}$ and
    $\triple<\gfunc';\nu';\varepsilon'>$ on $\mathcal{C}'$,
  \item a strong comonad transformation $\pair<\ffunc;\gamma>$ from the
    $\triple<\gfunc;\nu;\varepsilon>$ to $\triple<\gfunc';\nu';\varepsilon'>$,
    and
  \item functors $\arrow\cpt:\mathcal{C}->\Set.$ and
    $\arrow\cpt':\mathcal{C}'->\Set.$ along with an isomorphism
    \begin{displaymath}
      \let\labelstyle=\textstyle
      \xymatrix@R=2em@C=2em{
        {\mathcal{C}}\ar[rr]^{\ffunc}\ar[dr]_(0.4){\cpt}^{}="one" & &
        {\mathcal{C}'}\ar[dl]^(0.4){\cpt'}_{}="two" \\
        & {\Set} & 
        \ar@{<=}_{\scriptstyle\simeq}^{\alpha}%
          "two"+<-1em,0.4em>;"one"+<1em,0.4em>
      }
    \end{displaymath}
  \end{itemize}
  then we may construct a natural isomorphism
  \begin{displaymath}
    \let\labelstyle=\textstyle
    \xymatrix@R=2em@C=2em{
      {\mathcal{C}}\ar[rr]^{\ffunc}\ar[dr]_(0.4){\nerv_{\gfunc,\cpt}}^{}="one" & &
      {\mathcal{C}'}\ar[dl]^(0.4){\nerv_{\gfunc',\cpt'}}_{}="two" \\
      & {\Simp} & 
      \ar@{<=}_{\scriptstyle\simeq}^{\scriptstyle \nerv_{\gamma,\alpha}}%
      "two"+<-1em,0.5em>;"one"+<1em,0.5em>
    }
  \end{displaymath}
  in which the diagonals are the nerve functors associated with our
  comonads as in definition~\ref{defn.comonad.nerve}.
\end{lemma}

\begin{proof}
  Consider the following diagram
  \begin{displaymath}
    \let\labelstyle=\textstyle
    \xymatrix@R=1.5em@C=4em{
      {\mathcal{C}}\ar[rr]^<>(0.5){\hat{\gfunc}}_{}="two"\ar[dd]_{\ffunc} &&
      {\funcat[\op\aDelta,\mathcal{C}]}
      \ar[dd]_{\scriptstyle\funcat[\op\aDelta,\ffunc]}
      \ar[dr]^{\funcat[\op\aDelta,\cpt]}_{}="four" && \\
      &&& {\funcat[\op\aDelta,\Set]}\ar[r]^<>(0.5){\funcat[i,\Set]} &
      {\Simp} \\
      {\mathcal{C}'}\ar[rr]_<>(0.5){\hat{\gfunc'}}^{}="one" && 
      {\funcat[\op\aDelta,\mathcal{C}']}
      \ar[ur]_{\funcat[\op\aDelta,\cpt']}^{}="three" &&
      \ar@{<=} "one"+<-1em,2em>;"two"+<-1em,-2em> ^{\hat\gamma} _\cong
      \ar@{<=} "three"+<-0.6em,0.8em>;"four"+<-0.6em,-0.8em>%
        ^{\scriptstyle\funcat[\op\aDelta,\alpha]} _\cong
    }
  \end{displaymath}
  in which 
  \begin{itemize}
  \item the left hand square is that derived from the strong
    comonad transformation $\pair<\ffunc;\gamma>$ as in
    observation~\ref{comonad.trans},
  \item the upper composite of functors (from $\mathcal{C}$ to
    $\Simp$) is, by definition, the nerve functor associated with
    $\triple<\gfunc;\nu;\varepsilon>$ under $\cpt$.
  \item the lower composite of functors (from $\mathcal{C}'$ to
    $\Simp$) is, by definition, the nerve functor associated with
    $\triple<\gfunc';\nu';\varepsilon'>$ under $\cpt'$.
  \end{itemize}
  It follows that the natural transformation obtained by pasting this
  diagram provides us with the isomorphism required in the statement
  of the lemma.
\end{proof}

The following formalises a common application of the theory of
decalage:

\begin{cor}[simplicial reconstruction]\label{simplicial.recons.cor}
  Suppose that we are give a category $\mathcal{C}$, along with
  \begin{itemize}
  \item a comonad $\triple<\gfunc;\nu;\varepsilon>$ on $\mathcal{C}$,
  \item a strong comonad transformation $\pair<\ffunc;\gamma>$ from the
    decalage comonad $\triple<\Dec;\cdegen;\cface>$ on $\Simp$ to
    $\triple<\gfunc;\nu;\varepsilon>$, and
  \item a functor $\arrow\cpt:\mathcal{C}->\Set.$ along with an
    isomorphism
    \begin{displaymath}
      \let\labelstyle=\textstyle
      \xymatrix@R=2em@C=2em{
        {\Simp}\ar[rr]^<>(0.5){\ffunc}\ar[dr]_(0.4){\cpt_0}^{}="one" & &
        {\mathcal{C}}\ar[dl]^(0.4){\cpt}_{}="two" \\
        & {\Set} & 
        \ar@{<=}_{\scriptstyle\simeq}^{\alpha}%
          "two"+<-1em,0.4em>;"one"+<1em,0.4em>
      }
    \end{displaymath}
  \end{itemize}
  then the nerve functor $\arrow \nerv_{\gfunc,\cpt}:\mathcal{C}->\Simp.$
  associated with $\triple<\gfunc;\nu;\varepsilon>$ under $\cpt$ is a
  left pseudo-inverse of $\arrow \ffunc:\Simp->\mathcal{C}.$. In other
  words, there exists a natural isomorphism $\nerv_{\gfunc,\cpt}\circ \ffunc\cong
  \id_{\Simp}$.
\end{cor}

\begin{proof}
  Simply apply the previous lemma to the data given in the statement
  to obtain a natural isomorphism
  \begin{displaymath}
    \let\labelstyle=\textstyle
    \xymatrix@R=2em@C=2em{
      {\Simp}\ar[rr]^<>(0.5){\ffunc}\ar[dr]_(0.4){\Ndec}^{}="one" & &
      {\mathcal{C}}\ar[dl]^(0.4){\nerv_{\gfunc,\cpt}}_{}="two" \\
      & {\Simp} & 
      \ar@{<=}_{\scriptstyle\simeq}%
      "two"+<-1em,0.4em>;"one"+<1em,0.4em>
    }
  \end{displaymath}
  and observe, by lemma~\ref{nerve.dec}, that the nerve functor
  $\Ndec$ is isomorphic to the identity functor on $\Simp$.
\end{proof}

In other words, this lemma simply shows that under the conditions of
its statement we may use the comonad $\triple<\gfunc;\nu;\varepsilon>$ to
reconstruct each simplicial set $X$ from the object $\ffunc(X)\in\mathcal{C}$.

\begin{obs}[semi-simplicial versions]
  In this paper we will apply a version of the simplicial
  reconstruction lemma above, but it will suffice for us to
  reconstruct semi-simplicial structures only. Of course, all of the
  theory in this appendix may be restricted to consider only
  semi-simplicial sets, copointed endo-functors (rather then
  comonads), their transformations (appropriately defined) and the
  corresponding subcategory $\afDelta$ (in the place of
  $\aDelta$). Proofs of these semi-simplicial results remain
  substantially unchanged from their simplicial counterparts discussed
  above, so we leave the (trivial) details up to the reader.
\end{obs}



\section{Stratifications and Filterings of Simplicial Sets}
\label{strat.filt.sect}

\subsection{Stratified Simplicial Sets}\label{sect.strat}

\begin{defn}[stratified simplicial sets]\label{strat.defn}
  A {\em stratified simplicial set\/} (or sometimes just a {\em
    stratified set}) is a pair $\strat{X}$ where $X$ is a simplicial
  set and $tX$ is a subset\footnote{Note that $tX$ is merely a subset
    of $X$, {\bf not} a simplicial subset, in general it will not be
    closed in $X$ under the action of $\Delta$.} of its simplices,
  which are said to be {\em thin}, satisfying the conditions that
  \begin{itemize}
  \item no 0-simplex of $X$ is in $tX$,
  \item all of the degenerate simplices of $X$ are in $tX$.
  \end{itemize}
  
  A {\em stratified simplicial map\/} (or sometimes just a {\em
    stratified map}) $\arrow f:\strat{X}->\strat{Y}.$ is a
  simplicial map $\arrow f: X-> Y.$ such that $tX\subseteq
  f^{-1}(tY)$. In other words, a stratified simplicial map is a
  simplicial map that {\em preserves thinness}.

  The composite of two stratified simplicial maps again preserves
  thinness, so we may collect stratified simplicial sets and maps
  together into a category $\Strat$.
  
  We will often drop explicit mention of the set of thin simplices
  $tX$ and elect to notationally confuse a stratified set with its
  underlying simplicial set by simply declaring that $X$ (rather than
  $\strat{X}$) is a stratified set. In such cases, the set of thin
  simplices associated with the stratified sets $X,Y,Z,\dots$ will
  always be denoted by $tX,tY,tZ,\dots$ respectively.
\end{defn}

\begin{obs}[unity and identity of opposites]
  The forgetful functor $\arrow \forget:\Strat->\Simp.$,
  which takes $\strat{X}$ to $X$, admits both right {\bf and} left
  adjoints.  The right adjoint takes a simplicial set $X$ to
  $\radj(X)=\pair<X;X\setminus X_0>$, in which all simplices
  of dimension $>0$ are thin, whereas the left adjoint takes $X$ to
  $\refl(X)=\pair<X;X_d>$ where $X_d=\{x\in X\mid x\text{ is
    degenerate}\}$, making as few simplices thin as possible. Both of
  the functors $\radj$, $\refl$ are fully faithful, so
  this sequence of adjoints is a {\em unity and identity of
  opposites\/} (a situation studied in general by F.\ W.\ Lawvere).
  
  In general, we will identify $\Simp$ with its image under $\refl$, that
  is the full subcategory of $\Strat$ of those stratified sets such
  that every thin simplex is degenerate. In particular, we will
  assume, without comment, that the standard simplicial sets such as
  $\Delta[n]$ and $\Lambda_k[n]$ also live in $\Strat$.
\end{obs}

\begin{obs}[duals of stratified simplicial sets]\label{strat.dual} The dual
  $\strat{X}^\circ$ of a stratified simplicial set $\strat{X}$ is the
  pair $\strat{X^\circ}$ where $X^\circ$ is the dual simplicial set
  introduced in observation~\ref{simp.dual} and $tX^\circ\defeq
  tX$.\footnote{It is often quite useful to distinguish $tX$ and
    $tX^\circ$ notationally, even though they are identical as sets.}
  Quite clearly this pair constitutes a well defined stratified
  simplicial set.
  
  If $\arrow f:\strat{X}->\strat{Y}.$ is a stratified simplicial map
  then the dual simplicial map $\arrow f^\circ: X^\circ-> Y^\circ.$
  carries $tX^\circ$ into $tY^\circ$. It follows that $f^\circ$ is a
  stratified map from $\strat{X}^\circ$ to $\strat{Y}^\circ$ and that
  we may extend our dual operation to a functor
  $\arrow(\mathord{-})^\circ:\Strat->\Strat.$ which makes the diagram
  of functors
  \begin{displaymath}
    \xymatrix@C=4em@R=2em{
      {\Strat}\ar[r]^{\textstyle(\mathord{-})^\circ}
      \ar[d]_{\textstyle \forget} & {\Strat}\ar[d]^{\textstyle \forget} \\
      {\Simp}\ar[r]_{\textstyle(\mathord{-})^\circ} & {\Simp}}
  \end{displaymath}
  commute.
  Yet again this dual functor is strictly involutive, in that its
  composite with itself is the identity functor on $\Strat$.
\end{obs}

\begin{notation}[stratified simplicial subsets]
  A pair $\strat{Y}$ is a {\em stratified simplicial subset\/} (or
  just a {\em stratified subset}) of a stratified set $\strat{X}$,
  denoted $\strat{Y}\subseteq_s\strat{X}$, if:
  \begin{itemize}
  \item $Y$ is a simplicial subset of $X$,
  \item $tY$ contains all of the degenerate simplices of $Y$, and
  \item $tY$ is a subset of $tX$.
  \end{itemize}
  It follows from this definition that $\strat{Y}$ inherits a
  stratified simplicial structure from $\strat{X}$ and that the
  inclusion $Y\subseteq_s X$ becomes a stratified map
  $\overinc{\subseteq_s}:\strat{Y}->\strat{X}.$.
  
  Intersections and unions of stratified subsets
  $\strat{Y_i}\subseteq_s\strat{X}$ are given by the formulae:
  \begin{equation*}
    \bigcup_i\strat{Y_i} \defeq \pair<\bigcup_i Y_i;
    \bigcup_i tY_i> \mkern40mu
    \bigcap_i\strat{Y_i} \defeq \pair<\bigcap_i Y_i;
    \bigcap_i tY_i> 
  \end{equation*}
\end{notation}

\begin{notation}[regular, entire and inclusive stratified maps]
  \label{reg.entire.maps}\label{strat.isos} 
  We say that a stratified map $\arrow f:X->Y.$ is
  \begin{itemize}
  \item {\bf regular} if it {\em reflects thin simplices}, 
    which means that whenever $f(x)$ is thin in $Y$ it follows that
    $x$ is thin in $X$. In other words, $f$ is regular iff
    $f^{-1}(tY)=tX$. 
  \item {\bf entire} if it is surjective on simplices, that is to say
    if its underlying simplicial map is surjective. 
  \item {\bf an inclusion} if it is injective on simplices, that is to
    say if its underlying simplicial map is injective.
  \end{itemize}
  Sometimes we will adopt notation $\mathcal{R}$, $\mathcal{E}$ and
  $\mathcal{I}$ to denote the classes of regular, entire and inclusive
  stratified maps respectively.

  Notice that a stratified map is entire (resp.\ an inclusion) if and
  only if it is an epimorphism (resp.\ a monomorphism) in the usual
  categorical sense. Furthermore, observe that our stratified map
  $\arrow f:X->Y.$ is a stratified isomorphism (that is an isomorphism
  in the category $\Strat$) if and only if it is both bijective on
  simplices and regular. In other words, the class of stratified
  isomorphisms coincides with the intersection $\mathcal{E}\cap
  \mathcal{I}\cap\mathcal{R}$ of the three classes defined above.
\end{notation}

\begin{obs}[stratified images and inverse images]\label{strat.im}
  If $\arrow f:X->Y.$ is a stratified simplicial
  map then define
  \begin{itemize}
  \item the {\em image\/} of a stratified subset
    $X'\subseteq_s X$ under $f$ to be the stratified
    subset given by $f(X')\defeq\pair<f(X');f(tX')>$, and
  \item the {\em inverse image\/} of a stratified subset 
    $Y'\subseteq_s Y$ under $f$ to be the stratified
    subset given by $f^{-1}(Y')\defeq\pair<f^{-1}(Y');f^{-1}(tY')>$.
  \end{itemize}
  We leave it up to the reader to verify these definitions do indeed
  specify well defined stratified simplicial subsets. Observe that:
  \begin{itemize}
  \item for any stratified map $\arrow f:X->Y.$ the inverse image
    $f^{-1}(Y')$ of a regular subset $Y'\subseteq_r Y$ is a regular
    subset of $X$,
  \item for any {\bf regular} stratified map $\arrow
    f:X->Y.$ the image $f(X')$ of a regular subset $X'\subseteq_r X$
    is again a regular subset of $Y$.
  \end{itemize}
\end{obs}

\begin{notation}[regular and entire stratified subsets]
  \label{substrat.reg}
  We say that $X\subseteq_s Y$ is:
  \begin{itemize}
  \item {\em regular}, and use the notation $X\subseteq_r Y$, if the
    associated inclusion map $\overinc\subseteq_s:X->Y.$ is regular or
    equivalently when $tX=X\cap tY$.
  \item {\em entire}, and use the notation $X\subseteq_e Y$, if the
    associated inclusion map $\overinc\subseteq_s:X->Y.$ is entire or
    equivalently when $X$ and $Y$ have the same underlying simplicial
    set.
  \end{itemize}
  We define $\stratc{S}$, the regular stratified subset {\em
    generated\/} by a subset $S\subseteq X$, to be the smallest
  regular stratified subset of $X$ containing $S$ and clearly
  $\stratc{A} = \pair<\simpc{A};tX\cap\simpc{A}>$.
  
  Notice that the classes $\mathcal{R}$, $\mathcal{E}$ and
  $\mathcal{I}$ provide us with two distinct factorisation systems on
  $\Strat$ (in the sense of Freyd and Kelly~\cite{FreydKelly:1972:Fact}):
  \begin{itemize}
  \item $\pair<\mathcal{E};\mathcal{R}\cap\mathcal{I}>$ which
    factorises a stratified map $\arrow f:X->Y.$ as $\overarr
    f_e:X->\im_r(f).\overinc\subseteq_r:->Y.$, where
    $\im_r(f)\subseteq_r Y$ is the {\em regular image\/} of $f$, which
    is the regular stratified subset of $Y$ generated by $\{f(x)\in Y|
    x\in X\}$.
  \item $\pair<\mathcal{E}\cap\mathcal{I};\mathcal{R}>$ which
    factorises a stratified map $\arrow f:X->Y.$ as $\overinc
    \subseteq_e:X->\coim_e(f).\overarr f_r:->Y.$, where $X\subseteq_e
    \coim_e(f)$ is the {\em entire coimage\/} of $f$, which has the
    same underlying simplicial set as $X$ and thin simplices those
    $x\in X$ with $f(x)$ thin in $Y$.
  \end{itemize}
\end{notation}

\begin{notation}[lifting and extension]
  Suppose that $\arrow f:X->Y.$ is a stratified map:
  \begin{itemize}
  \item if $Y'$ is a regular subset of $Y$ then we say that $f$ {\em
      lifts\/} to a stratified map with codomain $Y'$ iff the
    simplex $f(x)$ is in $Y'$ for all $x\in X$, in which case we can
    factor $f$ through $Y'\subseteq_r Y$ and consider it to be a
    stratified map $\arrow f:X->Y'.$.
  \item if $X'$ is an entire superset of $X$ (i.e.\ $X\subseteq_e X'$)
    then we say that $f$ {\em extends\/} to a stratified map with domain
    $X'$ iff the simplex $f(x)$ is thin in $Y$ for all thin simplices
    $x$ in $X'$, in which case we can factor it through $X\subseteq_e
    X'$ and consider it to be a stratified map $\arrow f:X'->Y.$.
  \end{itemize}
\end{notation}

\begin{defn}\label{thinner.simplices}
  Let $\Simplex$ denote the full subcategory of $\Strat$ consisting of
  those stratified sets whose underlying simplicial sets is a standard
  simplex $\Delta[n]$ for some $n\in\mathbb{N}$. 
  
  We'll often adopt the alphabetic convention of assuming that the
  uppercase letter used to denote a given stratified set in $\Simplex$
  matches the lower case letter used to denote the dimension of its
  underlying simplex.  In other words, under this convention the
  stratified sets $N,M,P$ and $Q$ in $\Simplex$ would have underlying
  simplicial sets $\Delta[n]$, $\Delta[m]$, $\Delta[p]$ and
  $\Delta[q]$ respectively whereas $N,N'$ and
  $N_i$ in $\Simplex$ would all have underlying simplicial set
  $\Delta[n]$.
\end{defn}

\begin{notation}[some standard stratified simplicial sets]\label{stan.strat}
  A few stratified simplicial sets of note are:
\begin{enumerate}[(i)]
\item {\bf The standard thin $n$-simplex} $\Delta[n]_t$, which is
  obtained from $\Delta[n]$ by making its only non-degenerate
  $n$-simplex $\arrow\id_{[n]}:[n]->[n].$ thin.
\item\label{admissible} {\bf The $k^{\text{th}}$ standard admissible
    $n$-simplex} $\Delta^a_k[n]$(for $n\geq 2$, $k=1,\dots,n-1$) is
  obtained from $\Delta[n]$ by making all (non-degenerate) simplices
  $\arrow\alpha:[m]->[n].$ with $k-1,k,k+1\in\im(\alpha)$ thin.  These
  are particularly important in characterising the nerves of
  $n$-categories.
\item {\bf The standard admissible $(n-1)$-dimensional
    $k$-horn} $\Lambda^a_k[n]$ (for $n\geq 2$, $k=1,\dots,n-1$) is
  obtained from $\Lambda_k[n]$ by making any simplex which is thin in
  $\Delta_k[n]$ thin in $\Lambda^a_k[n]$. In other words,
  $\Lambda^a_k[n]$ is the regular stratified subset of $\Delta^a_k[n]$
  generated by the set of simplices $\{\delta^n_i\mid i=0,\dots, k-1,
  k+1,\dots,n\}$.
\end{enumerate}
The non-degenerate simplices of $\Delta[n]$ satisfying the condition
given in example~(\ref{admissible}) are said to be {\em $k$-divided}.
It should be noted that this use of the term differs slightly from
that used by Street in~\cite{Street:1987:Oriental}, a difference we
discuss a little more in theorem~\ref{nerves.are.complicial}.
\end{notation}

\begin{obs}[admissible simplices and horns in a stratified simplicial set]
\label{admiss.horn} An $n$-simplex $x$ in a stratified
simplicial set $X$ is said to be {\em $k$-admissible\/} if the
corresponding simplicial map $\arrow\yoneda{x}:\Delta[n]-> X.$ (under
Yoneda's lemma) extends to a stratified simplicial map $\arrow\yoneda{x}
:\Delta^a_k[n]->X.$.

By the definition of $\yoneda{x}$, that is $\yoneda{x}(\alpha)\defeq
x\cdot\alpha$, this holds if and only if $x\cdot\alpha$ is thin in
$X$ for all $k$-divided simplices $\alpha\in\Delta[n]$.

Similarly, we say that an $(n-1)$-dimensional $k$-horn $\{x_i\in
X_{n-1}\mid i=0,1,...,k-1,k+1,k+2,...,n\}$ in $X$ is {\em
  admissible\/} if the corresponding simplicial map
$\arrow:\Lambda_k[n]->X.$ extends to a stratified map from
$\Lambda^a_k[n]$ to $X$. By Yoneda's lemma, this condition holds if
and only if
\begin{itemize}
\item $x_i$ is $(k-1)$-admissible in $X$ for each $i<k-1$, and
\item $x_i$ is $k$-admissible in $X$ for each $i>k+1$.
\end{itemize}
\end{obs}

\begin{obs}[duals of standard stratified simplicial sets]\label{stan.dual}
The dual operation on $\Delta$ gives rise to a canonical simplicial
isomorphism:
\begin{displaymath}
  \xymatrix@C=8em@R=0.1em{ {\Delta[n]}\ar[r]^{(\mathord{-})^\circ} &
    \Delta[n]^\circ \\
    \alpha\ar@{|->}[r] & \alpha^\circ }
\end{displaymath}
This extends to a family of stratified simplicial isomorphisms between
the standard stratified sets introduced in the last observation and
their duals $\Delta[n]_t\cong\Delta[n]_t^\circ$ and
$\Delta^a_{n-k}[n]\cong\Delta^a_k[n]^\circ$. Furthermore this latter
isomorphism restricts to an isomorphism
$\Lambda^a_{n-k}[n]\cong\Lambda^a_k[n]^\circ$ of admissible horns and
their duals.
\end{obs}

\begin{obs}[$\Strat$ as an LFP-category]\label{deltat.dense}
  The small full subcategory $\tDelta$ on the set of objects
  $\{\Delta[n],\Delta[n]_t\mid n=0,1,\dots\}$ in $\Strat$ gives rise
  to a right adjoint functor 
  \begin{equation}\label{inc.strat.pres}
    \let\labelstyle=\textstyle
    \xymatrix@R=1ex@C=14em{
      {\Strat}\ar[r]^{\bot} &
      {[\op{\tDelta},\Set]}\ar@{.>}@/_3ex/[l] }
  \end{equation}
  by Kan's construction (observation~\ref{func.from.coalg}).
  
  Applying Yoneda's lemma, we may describe $\tDelta$ as an extension
  of $\Delta$ in the way originally proposed by Street
  in~\cite{Street:1987:Oriental}. Under this description $\tDelta$ has
  two families of objects, the original family of finite ordinals
  $\{[n]\mid n\geq 0\}$ of $\Delta$ and a family of copies of the
  non-zero ordinals $\{[n]_t\mid n>0\}$. In other words $[n]$ and
  $[n]_t$ are distinct copies of the finite ordinal with $n+1$
  elements, where the former copy represents the standard $n$-simplex
  and the latter is a model for the standard thin $n$-simplex. In the
  sequel, we will sometimes also use the notation $[n]_?$ if we wish
  to refer to either of the objects $[n]$ or $[n]_t$. The arrows of
  $\tDelta$ are simply order preserving maps subject to a restriction
  which excludes all non-identity face operators (injective maps)
  whose domains are of the form $[n]_t$.
  
  This category may be presented in terms of generators which include
  the elementary face and degeneracy operators of $\Delta$, the
  elementary degeneracies $\arrow \tdegen^n_i:[n+1]_t->[n].$ and the
  operator $\arrow \thop^n:[n]->[n]_t.$ that is the unique order
  isomorphism between the ordinal $[n]$ and its copy $[n]_t$. The
  relations that must hold between these operators are those of
  $\Delta$ plus the extra family of relations
  $\degen^n_i=\tdegen^n_i\circ \thop^{n+1}$ (for $n\in\mathbb{N}$). 
  
  In line with our notation for stratified standard simplices, we will
  let $\arrow\Delta:\tDelta->\Strat.$ denote the Yoneda functor which
  embeds $\tDelta$ (under this alternative presentation) into
  $\Strat$. In particular, this means that $\arrow\Delta(\thop^n):
  \Delta[n]->\Delta[n]_t.$ and $\arrow\Delta(\tdegen^n_i):
  \Delta[n+1]_t->\Delta[n].$ denote the stratified maps whose
  underlying simplicial maps are the identity on $\Delta[n]$ and the
  map $\arrow\Delta(\degen^n_i):\Delta[n+1] ->\Delta[n].$
  respectively.
  
  Using a Yoneda argument we may demonstrate that the right adjoint
  functor in display~(\ref{inc.strat.pres}) is fully faithful, or in
  other words that $\tDelta$ is embedded in $\Strat$ as a small {\em
    dense\/} subcategory. Consequently, we will identify $\Strat$ with
  the reflective full subcategory of $\funcat[\op\tDelta,\Set]$ which
  is the replete image of this right adjoint functor and use the
  notation $\arrow\refl_s:\funcat[ \op\tDelta,\Set]->\funcat[
  \op\tDelta,\Set].$ for the corresponding reflector.  Notice also
  that an object $F$ of $\funcat[\op\tDelta,\Set]$ is in the full
  subcategory $\Strat$ if and only if the function $\arrow
  F(\thop^n):F([n]_t)-> F([n]).$ is injective for each $n>0$.  This
  condition may equivalently be expressed as stating that the pullback
  of each $F(\thop^n)$ along itself should be the identity on
  $F([n]_t)$ which in turn may be written as a perpendicularity
  condition with respect to certain maps between finitely presented
  objects in $\funcat[\op\tDelta,\Set]$. It follows, therefore, that
  $\Strat$ may be described as the LFP-category of algebras of an
  LE-theory $\mathbb{T}_{\Strat}$ with underlying category $\tDelta$
  and that the associated reflector $\refl_s$ is finitely accessible
  (see subsection~\ref{lfp.le.subsect}).
  
  Of course, now observation~\ref{full.refl.LFP} may be applied to
  show that $\Strat$ is closed in $\funcat[\op\tDelta,\Set]$ under
  limits and filtered colimits, with all other colimits in $\Strat$
  being constructed by forming them point-wise in
  $\funcat[\op\tDelta,\Set]$ and then reflecting them back into
  $\Strat$ using $\refl_s$. Furthermore, consulting
  definition~\ref{weighted.lims} we also see that every stratified
  simplicial set $X$ may be expressed in $\Strat$ as the canonical
  colimit $\colim(X,\Delta)$ of standard simplices and standard thin
  simplices weighted by $X$ (regarded as an object in
  $\funcat[\op\tDelta,\Set]$).  Alternatively, as discussed in the
  text accompanying definition~\ref{weighted.lims}, we may describe
  this colimit in terms of a traditional conical colimit by
  constructing the Grothendieck category $\groth(X)\defeq
  \tDelta\downarrow X$, forming the diagram $\arrow
  \dgm_X:\groth(X)->\Strat.$ by domain projection and using a Yoneda
  argument to show that the obvious cocone $\nattrans i: \dgm_X->X.$
  with components $\arrow i_f=f:\dgm_X(f)=\Delta[n]_{?}-> X.$ induces
  a canonical isomorphism $X\cong\colim_{\groth(X)}(\dgm_X)$.
  
  If $X$ is a stratified set then the category $\groth(X)$ has a {\em
    cofinal\/} set of objects consisting of those stratified maps
  $\arrow \yoneda{x}:\Delta[n]_{?}->X.$ with $[n]_{?}\in\tDelta$ which
  correspond under Yoneda's lemma to simplices $x\in X$ which are
  non-degenerate. It follows that if $X$ has only a finite set of
  non-degenerate simplices then it may be expressed as a {\bf finite}
  colimit of standard (thin) simplices. However, we know that every
  standard (thin) simplex is finitely presentable and that any finite
  colimit of FP-stratified sets is also finitely presentable, so we
  may infer that any stratified set with only a finite set of
  non-degenerate simplices is finitely presentable. Indeed the reverse
  implication also holds, but we won't use it here and we leave its
  proof up to the reader.
\end{obs}

\begin{obs}[limits and colimits in $\Strat$]
  \label{stratified.(co)limits} If $\arrow \dgm:\mathbb{C}->
  \Strat.$ is a diagram of stratified simplicial sets, where
  $\mathbb{C}$ is a small category, then the limit and colimit of $\dgm$
  may be obtained via the following explicit constructions:\vspace{1ex}

  \noindent{\bf Limits} First form the limit $\lim(\forget\circ \dgm)$ 
  in $\Simp$, as described in observation~\ref{simplicial.(co)limits},
  then raise this to a stratified simplicial set by making thin those
  simplices $\{x_i\in \dgm(i)\}_{i\in\obj(\mathbb{C})}$ for which
  $x_i$ is thin in $\dgm(i)$ for each $i\in\obj(\mathbb{C})$.\vspace{1ex}

  \noindent{\bf Colimits} Again start by forming the colimit 
  $\colim(\forget\circ \dgm)$ in $\Simp$, cf.\ 
  observation~\ref{simplicial.(co)limits}, then raise this to $\Strat$
  by making thin those simplices that can be written as an equivalence
  class $[x,i]$ with $x$ thin in $\dgm(i)$.\vspace{1ex}

Given these descriptions, it is clear that many results about limits
and colimits in $\Simp$ may be routinely raised to $\Strat$, simply by
appealing to the result on underlying simplicial sets and then
checking that stratifications are appropriately respected. In
particular, it is immediate that the result of
observation~\ref{union.qua.widepo} also holds for unions of stratified
subsets.
\end{obs}

\begin{obs}[skeleta and coskeleta of stratified sets]
  \label{(co)skeletal.strat} 
  Fixing $n\in\mathbb{N}$, we say that a stratified set $X$ is {\em
    $n$-skeletal\/} iff each simplex $x\in X$ with dimension
  $\dim(x)>n$ is degenerate.  The full subcategory of $n$-skeletal
  stratified sets is coreflective in $\Strat$, with coreflector
  $\arrow\Sk_n(\cdot):\Strat->\Strat.$ (cf.\ 
  observation~\ref{refl.full.subcat}) which maps a stratified set $X$
  to the regular stratified subset $\Sk_n(X)\subseteq_r X$, called its
  {\em $n$-skeleton}, of those simplices $x\in X$ with
  $x=x'\cdot\alpha$ for some $x'$ with $\dim(x')\leq n$.  In other
  words, this is the unique largest $n$-skeletal regular subset of
  $X$.
 
  This $n$-skeleton functor has a right adjoint
  $\arrow\Cosk_n:\Strat-> \Strat\subseteq\funcat[\op\tDelta,\Set].$
  formed by applying Kan's construction
  (observation~\ref{func.from.coalg}) to the functor
  $\arrow\Sk_n({\Delta( \cdot)}):\tDelta ->\Strat.$.  Equivalently,
  applying a Yoneda argument we see that $\Cosk_n(X)$, the {\em
    $n$-coskeleton\/} of the stratified set $\Strat$, may be described
  as the stratified set whose simplices, operator actions and
  stratification coincide with that of $X$ at and below dimension $n$
  and which is completely determined at each dimension $r$ above $n$
  by the fact that each one of its $(r-1)$-dimensional cycles is the
  boundary of a unique thin $r$-simplex.
  
  If we let $\tDelta_n$ denote the full subcategory of $\tDelta$ on
  those standard (thin) simplices of dimension $\leq n$ then it is an
  immediate consequence of the Eilenberg-Zilber lemma that $\tDelta_n$
  is dense in the full subcategory of $n$-skeletal stratified sets in
  $\Strat$. It follows that we may use the Grothendieck construction,
  just as we did in observation~\ref{deltat.dense}, to express the
  $n$-skeleton of a stratified set $X$ as a canonical colimit of
  standard (thin) simplices of dimension $\leq n$. Adopting the
  notation $\groth_n(X)$ for the comma category $\tDelta_n\downarrow
  X$ and $\dgm_X^n$ for the canonical diagram obtained by ``codomain
  projection'' then we find that each leg of the canonical cocone from
  this diagram into $X$ restricts to a map into the regular subset
  $\Sk_n(X)$ and provides us with a colimiting cocone which induces a
  canonical isomorphism $\colim(\dgm_X^n)\cong\Sk_n(X)$.
\end{obs}

\begin{obs}[connected components of stratified sets]
  \label{conn.cpt}
  Under our identification of $\Simp$ with the subcategory of
  minimally stratified sets in $\Strat$, we may extend the adjunction
  $\arrow\cpt_0\dashv\dis:\Set->\Simp.$ of
  observation~\ref{conn.cpt.simp} to one between $\Set$ and the
  category of stratified sets $\Strat$. If $X$ is a stratified set
  then its set of connected components $\cpt_0(X)$ is simply formed by
  applying the simplicial connected components functor to the
  underlying simplicial set $\forget(X)$. As before, it is clear that
  $\arrow\dis:\Set->\Strat.$ is fully faithful and that a stratified
  set is in its replete image iff it is 0-skeletal. It follows that
  $\dis$ provides an equivalence between $\Set$ and the full
  subcategory of 0-skeletal stratified sets in $\Strat$.
  
  It follows that the category of 0-skeletal stratified sets is also
  a {\bf reflective} full subcategory of $\Strat$ with corresponding
  reflector $\refl_0=\dis\circ\cpt_0$. A stratified map $\arrow
  f:X->Y.$ is $\refl_0$-invertible (cf.definition~\ref{l-inv.def}) iff
  the function $\arrow\cpt_0(f): \cpt_0(X)->\cpt_0(Y).$ is an
  isomorphism, consequently we say that such a map is {\em bijective
    on components}.  We also say that $X$ is {\em connected\/} if
  $\cpt_0(X)$ is the singleton set. From the description of $\cpt_0$
  above, it follows immediately that:
  \begin{itemize}
  \item The set of connected components of a stratified set is
    dependent only upon the structure of its underlying simplicial
    set.
  \item All of the standard simplices and horns are connected.
  \item The connected components functor $\cpt_0$ preserves products
    and it follows that products of connected sets are connected.
  \end{itemize}
  Furthermore, it is worth mentioning that any stratified map $\arrow
  f:X->Y.$ whose domain and codomain are connected is bijective on
  components, since any function between the singleton sets
  $\cpt_0(X)$ and $\cpt_0(Y)$ must be a bijection.
\end{obs}

\subsection{Superstructures and Filtered Semi-Simplicial Sets}
\label{sect.filt.semi}

\begin{defn}[superstructures of stratified sets]\label{superstr.defn}
  For each $n\in\mathbb{N}$ there exists a pair of adjoint functors
  \begin{displaymath}
    \let\labelstyle=\textstyle
    \xymatrix@R=5ex@C=18em{
      {\Strat}
      \ar@/_2.5ex/[r]_{\labelstyle \Sup_n(\cdot)}^{}="1" &
      {\Strat}\ar@/_2.5ex/[l]_{\Th_n}^{}="3"
      \ar@{}"1";"3"|{\bot} 
    }
  \end{displaymath}
  which are defined as follows:

  \begin{domsitem}
  \item {\boldmath $\Th_n(X)$} is defined to be $\pair<X;tX\cup{\{x\in
      X\mid \dim(x)>n\}}>$, in other words it is the stratified
    simplicial set obtained by making thin all of the simplices of $X$
    with dimension $>n$. It is clear that any stratified map $\arrow
    f:X->Y.$ extends to a stratified map from $\Th_n(X)$ to $\Th_n(Y)$,
    so we let $\Th_n(f)\defeq f$.
 
  \item {\boldmath $\Sup_n(X)$}, called the {\em $n$-dimensional
      superstructure\/} of $X$ is its regular stratified subset with
    underlying simplicial set given by:
    \begin{displaymath}
      \left\{x\in X\mid(\forall\arrow\alpha:[r]->[\dim(x)].\in\Delta)\; 
      r>n \Rightarrow x\cdot\alpha\in tX\right\}
    \end{displaymath}
    In other words, $\Sup_n(X)$ is the regular subset of those
    simplices of $X$ for which each $r$-dimensional face with $r>n$ is
    thin.  We may easily verify that a stratified map $\arrow f:X->Y.$
    has $\Sup_n(X)\subseteq_s f^{-1}(\Sup_n(Y))$ and thus restricts to
    a stratified map $\arrow \Sup_n(f):\Sup_n(X)->\Sup_n(Y).$.
  \end{domsitem}
  
  We say that a stratified set $X$ is {\em $n$-trivial\/} if and only
  if $X=\Sup_n(X)$ (or equivalently iff $X=\Th_n(X)$) and we use the
  notation $\Strat_n$ to denote the full subcategory of these in
  $\Strat$. In other words, $X$ is $n$-trivial if and only if all
  of its simplices of dimension greater than $n$ are thin.
\end{defn}

\begin{obs}\label{obs.super}
  The following properties of the superstructures of a stratified set
  $X$ are worthy of note.
  \begin{enumerate}
  \item\label{obs.super.a} If $m\leq n$ then $\Sup_m(X)$ is a regular
    subset of $\Sup_n(X)$.
  \item\label{obs.super.b} If $m < r$ and $x$ is an $m$-simplex of $X$
    then all of the $r$-dimensional faces of $x$ are degenerate, and
    thus thin, in $X$. It follows that $\Sup_n(X)$ contains {\bf all}
    of the simplices of $X$ with dimension $\leq n$.
  \item\label{obs.super.c} If $x$ is an $n$-simplex of $X$ then the
    only non-degenerate $n$-dimensional face of $x$ is $x$ itself, so
    it follows that this is an element of $\Sup_{n-1}(X)$ if and only if
    it is thin in $X$.
  \end{enumerate}
  In particular, the second of these observations immediately implies
  that $X$ is equal to the union $\bigcup_{n\in\mathbb{N}} \Sup_n(X)$
  of its superstructures.
\end{obs}

\begin{defn}[categories of filtered objects] 
  The partially ordered set of natural numbers $\mathbb{N}$ may be
  considered to be a category in the usual way. That is to say we
  think of it as a category with objects non-negative integers and a
  unique arrow $\overarr:i->j.$ for each pair of integers with $i\leq
  j$. If $\mathcal{C}$ is any category then we call the functor
  category $\funcat[\mathbb{N},\mathcal{C}]$ the {\em category of
    filtered objects\/} in $\mathcal{C}$.
\end{defn}

\begin{obs}\label{filtering.strat} 
  The family of superstructures of $X$ provides us with a filtered
  family of stratified sets, since we have $\Sup_i(X)\subseteq_r
  \Sup_{i+1}(X)$ for each $i\in\mathbb{N}$.  Now, if we throw away
  stratifications and forget about the actions of degeneracy operators
  we may consider this to be a filtered family of semi-simplicial
  sets. It follows that we may gather this information together to
  provide a single functor
  \begin{displaymath}
    \xymatrix@R=2ex@C=6em{
      {\Strat}\ar[r]^<>(0.5){\textstyle\Sup_{\bullet}(\cdot)} & 
      \funcat[\mathbb{N},\SSimp]
    }
  \end{displaymath}
  which represents each stratified set as a filtered semi-simplicial
  set. 
  
  It is not, in general, the case that this functor is fully
  faithful, however it is possible to find a substantial full
  subcategory of $\Strat$ on which it is. In the remainder of this
  section we establish a characterisation of this subcategory and in
  section~\ref{comp.sec} we demonstrate that it contains all of those
  stratified sets which Street~(\cite{Street:1987:Oriental} and
  \cite{Street:1988:Fillers}) identified as being candidates for those
  in the replete image of his \inf-categorical nerve functor.
\end{obs}

\begin{defn}\label{defn.predegen}  
  We say that a simplex $x$ in a stratified set $X$ is {\em
    pre-degenerate at $k$\/} (where $0\leq k < \dim(x)$) if
  $x\cdot\alpha$ is thin in $X$ for each face operator $\alpha$ for
  which $k,k+1\in\im(\alpha)$.
  
  In particular, observation~\ref{degen.obs}(\ref{degen.obs.c}) and
  the fact that in a stratified set all degenerate simplices are thin
  implies that every simplex $x\in X$ which is degenerate at $k$ is
  pre-degenerate at $k$. The converse is true for some stratified
  sets, including those that occur as the nerves of \inf-categories.
\end{defn}

\begin{defn}\label{well.tempered} 
  If $X$ is a stratified simplicial set then we say that it is {\em
    well tempered\/} iff any simplex $x\in X$ which is pre-degenerate
  at some $k$ is in fact degenerate at $k$.  In other words, $X$ is
  well tempered if its collection of thin simplices is sufficient to
  detect those simplices which are degenerate.

  We use the notation $\Strat_W$ to denote the full subcategory of
  $\Strat$ on those stratified sets which are well tempered.
\end{defn}

Our primary motivation for introducing these definitions is:

\begin{lemma}\label{ss.tp=>strat} 
  Suppose that $Y$ is a well tempered stratified set and $\arrow
  f:X->Y.$ is a semi-simplicial map which preserves thinness then $f$
  is a stratified map.
\end{lemma}

\begin{proof} 
  First recall that every simplicial operator factors as a face
  operator following a degeneracy operator and that this latter
  operator, in turn, factors as a composite of elementary degeneracy
  operators. However, we know that $f$ preserves the actions of face
  operators (since it is a semi-simplicial map) so if we can show
  that it also preserves the action of elementary degeneracy operators
  then it preserves the actions of all operators and is thus
  a simplicial map.
  
  Now the notion of pre-degeneracy is defined in terms of actions of
  {\bf face} operators and thinness, thus it follows that a map like
  $f$ which preserves both of these structures will also preserve
  pre-degeneracy. So suppose that we are given an $n$-simplex $x\in X$
  and an elementary degeneracy operator
  $\arrow\degen^n_k:[n+1]->[n].$, then we know that $x\cdot\degen^n_k$
  is pre-degenerate at $k$ and so we may infer that
  $f(x\cdot\degen^n_k)\in Y$ is also pre-degenerate at $k$. But $Y$ is
  well tempered and so it follows that $f(x\cdot\degen^n_k)$ is
  actually degenerate at $k$.
  
  Finally we may demonstrate that $f$ preserves the action of
  $\degen^n_k$ on $x$ using the following simple calculation:
  \begin{displaymath}
    \begin{array}{rcl@{\hspace{0.5cm}}p{8cm}}
      f(x\cdot\degen^n_k) & = &  (f(x\cdot\degen^n_k)\cdot
      \face^{n+1}_k)\cdot\degen^n_k & since $f(x\cdot\degen^n_k)$ is
      degenerate at $k$, \\
      & = & f((x\cdot\degen^n_k)\cdot\face^{n+1}_k)\cdot\degen^n_k &
      $f$ preserves the action of face operators, \\
      & = & f(x\cdot(\degen^n_k\circ\face^{n+1}_k))\cdot\degen^n_k &
      right action, \\
      & = & f(x)\cdot\degen^n_k & $\face^{n+1}_k$ is a right
      inverse of $\degen^n_k$.
    \end{array}
  \end{displaymath}
  Therefore, quantifying over $x$ and $\degen^n_k$, we see that $f$
  preserves the actions of all elementary face operators and so, as we
  argued in the first paragraph, it follows that $f$ is a simplicial
  map. However, we also assumed that $f$ preserved thinness and it is
  thus a stratified map as required.
\end{proof}

\begin{lemma}
  \label{well.tempered.rep}
  Restricting the domain of the superstructure functor of
  observation~\ref{filtering.strat} to the category of well tempered
  stratified sets
  \begin{displaymath}
    \xymatrix@R=2ex@C=6em{
      {\Strat_W}\ar[r]^<>(0.5){\textstyle\Sup_{\bullet}(\cdot)} & 
      \funcat[\mathbb{N},\SSimp]
    }
  \end{displaymath}
  we obtain a functor which is fully faithful.
\end{lemma}

\begin{proof}
  To show that $\Sup_{\bullet}(\cdot)$ is faithful, fix two stratified maps
  $\arrow f,g:X->Y.$ and observe that, by definition,
  $\Sup_{\bullet}(f)=\Sup_{\bullet}(g)$ if and only if $f$ and $g$
  agree when restricted to each superstructure $\Sup_n(X)$ of $X$. It
  follows that $f$ and $g$ agree on the union of these superstructures
  which, by observation~\ref{obs.super}, is equal to $X$ itself and so
  they are indeed equal on the whole of $X$ as required.
  
  To show that $\Sup_{\bullet}(\cdot)$ is full, suppose that $X$ and
  $Y$ are well tempered stratified sets and consider a filtered map
  $\arrow f_\bullet: \Sup_{\bullet}(X)->\Sup_{\bullet}(Y).$ in
  $\funcat[\mathbb{N},\SSimp]$. This is simply a family of
  semi-simplicial maps $\arrow f_i:\Sup_i(X)-> \Sup_i(Y).$ for
  $i\in\mathbb{N}$ satisfying the naturality condition that if $i<j$
  then $f_i$ and $f_j$ act identically on the subset $\Sup_i(X)
  \subseteq\Sup_j(X)$.
  
  By observation~\ref{obs.super}(\ref{obs.super.b}) above it is clear
  that we can define a map $\arrow f:X->Y.$ by $f(x)=f_{\dim(x)}(x)$
  and, indeed, the naturality property of $f_\bullet$ actually implies
  that if $n$ is an integer for which $x\in\Sup_n(X)$ then we have
  $f(x)=f_n(x)$.  In particular, it follows that the restriction of
  $f$ to $\Sup_i(X)$ is equal to $f_i$ for each $i\in\mathbb{N}$, in
  other words we can infer that $\Sup_{\bullet}(f)=f_\bullet$ as
  required so long as we can prove that $f$ is a stratified map. To
  prove this all we need do is demonstrate that this map:
  \begin{enumerate}[(a)]
  \item\label{filt.embed.pf.a} {\bf is semi-simplicial:} If $x$ is an
    $n$-simplex in $X$ and $\arrow\alpha:[m]->[n].$ is a face operator
    then we know that $m<n$ and so, by
    observation~\ref{obs.super}(\ref{obs.super.b}), we may infer that
    both $x$ and $x\cdot\alpha$ are elements of $\Sup_n(X)$. However,
    we also know that $f$ and $f_n$ coincide on $\Sup_n(X)$ and that
    $f_n$ is semi-simplicial, so it follows that $f(x\cdot\alpha)=
    f_n(x\cdot\alpha)=f_n(x)\cdot\alpha=f(x)\cdot\alpha$, that is to
    say $f$ preserves the action of the face operator $\alpha$ on
    $x$. Now, quantifying over $x$ and $\alpha$, we see that $f$ is
    semi-simplicial as required.
  \item\label{filt.embed.pf.b} {\bf preserves thinness:} If $x$ is a
    thin $n$-simplex in $X$ then we know, by
    observation~\ref{obs.super}(\ref{obs.super.c}), that it is an
    element of $\Sup_{n-1}(X)$. However, we also know that $f_{n-1}$
    and $f$ agree on $\Sup_{n-1}(X)$ and so it follows that $f(x)$ is
    an element of the codomain $\Sup_{n-1}(Y)$ of $f_{n-1}$. Applying
    observation~\ref{obs.super}(\ref{obs.super.c}) again we may infer
    that $f(x)$ is thin in $Y$ as required.
  \end{enumerate}
  Finally, we may apply lemma~\ref{ss.tp=>strat} to demonstrate that
  $f$ is a stratified map as required. 
\end{proof}



\section{Pre-Complicial Sets}\label{precomp.sec}

\subsection{Introducing Pre-Complicial Sets}

\begin{defn}[primitive t-extensions] For $n\geq 2$ and $k=1,\dots,n-1$
  let
\begin{itemize}
\item $\Delta_k^a[n]'$ denote the stratified set obtained
  from $\Delta_k^a[n]$ by making the $(n-1)$-simplices $\face^n_{k-1}$
  and $\face^n_{k+1}$ thin, and
\item $\Delta_k^a[n]''$ denote the stratified set obtained
  from $\Delta^a_k[n]'$ by making the $(n-1)$-simplex $\face^n_k$
  thin.
\end{itemize}
We say that the entire subset inclusion 
\begin{equation*}
  \xymatrix@R=1ex@C=10em{
    {\Delta^a_k[n]'}\ar@{^(->}[r]^{\textstyle\subseteq_e} &
    {\Delta^a_k[n]''}
  }
\end{equation*}
is a {\em primitive t-extension}.
\end{defn}

\begin{defn}[pre-complicial sets]\label{defn.precomp} A stratified set $A$ is
  said to be a {\em pre-complicial set\/} if it is orthogonal to all
  primitive t-extensions.  Let $\Precomp$ denote the full subcategory
  of $\Strat$ whose objects are the pre-complicial sets.
  
  Of course, we know that each of the stratified sets $\Delta^a_k[n]'$
  and $\Delta^a_k[n]''$ have only a finite set of non-degenerate
  simplices and are therefore finitely presentable in $\Strat$ (cf.\ 
  observation~\ref{deltat.dense}). It follows that the set of
  primitive t-extensions is an FP-regulus and so, by
  observation~\ref{full.refl.LFP}, that $\Precomp$ is also an
  LFP-category which is reflective in $\Strat$ (with associated
  idempotent monad $\pair<\refl_p;\eta^p>$). 
  
  We let $\mathbb{T}_{\Precomp}$ denote the LE-theory of
  pre-complicial sets, which we may construct from
  $\mathbb{T}_{\Strat}$, the LE-theory of stratified sets, by adding
  the set of primitive t-extensions to its FP-regulus.
\end{defn}

\begin{obs}[describing pre-complicial sets explicitly] By Yoneda's lemma,
  a stratified map $\overarr:\Delta^a_k[n]->A.$ corresponds to an
  $n$-simplex $a\in A$ such that $a\cdot\alpha$ is thin for each
  $k$-divided operator $\alpha$. We call such a simplex
  $k$-admissible.
  
  So, a set $A\in\Strat$ is pre-complicial iff for each $k$
  and each $k$-admissible simplex $a\in A$ if $a\cdot\face_{k-1}$ and
  $a\cdot\face_{k+1}$ are thin then $a\cdot\face_k$ is thin as well.
\end{obs}

\begin{obs}[an explicit description of the reflector $\refl_p$]
  \label{precomp.refl}
  In the sequel it will be useful to have a concrete description of
  the reflector $\arrow \refl_p:\Strat->\Strat.$ associated with
  $\Precomp$. Starting with a stratified set $X$ we inductively
  construct a sequence of subsets:
  \begin{eqnarray*}
    tX_0 & = & tX \\ 
    tX_{n+1} & = & \left\{x\in X \left| 
        \begin{array}[c]{l} 
          (\exists w\in X, k\in\mathbb{N}) \text{ s.t. $w$ is
            $k$-admissible in $\pair<X;tX_n>$,} \\
          \text{$w.\face_{k-1}, w.\face_{k+1}\in
            tX_n$ and $w\cdot\face_k=x$}
        \end{array}\right.\right\}
  \end{eqnarray*}
  Notice that a degenerate simplex $x\cdot\sigma_k$ is always
  $k$-admissible, from which it follows that $tX_0\subseteq
  tX_1\subseteq\cdots\subseteq tX_n\subseteq\cdots$, and define
  $\refl_p(X)$ to be the stratified set which has the same underlying
  simplicial set as $X$ but has $t\refl_p(X)\defeq\bigcup_{n=1}^{\infty}
  tX_n$ as its set of thin simplices. It is clear, by
  construction, that if $\arrow f:X->Y.$ is a stratified map then
  $f(tX_n)\subset tY_n$ for each $n\geq 0$ and so $f$ extends to a
  stratified map $\arrow f:\refl_p(X)->\refl_p(Y).$, making $\refl_p$ into a
  functor. It is also easily shown that $\refl_p(X)$ is a pre-complicial
  set and that the family of inclusions $\overinc\subseteq_r:X->\refl_p(X).$
  constitutes the unit $\eta^p$ for our idempotent monad.
\end{obs}

\begin{defn}[t-extensions] We call the $\refl_p$-invertible stratified
  maps {\em t-extensions}. For more on the properties of t-extensions
  please consult subsection~\ref{append.reflective}. We will also use
  the term {\em t-almost\/} as a synonym for $\refl_p$-almost (cf.\ 
  observation~\ref{l-almost}).
\end{defn}

\begin{obs}[duals of pre-complicial sets]\label{precomp.duals} We can
  extend the stratified isomorphisms of observation~\ref{stan.dual} to
  get two further families of stratified isomorphisms
  $\Delta^a_{n-k}[n]'\cong(\Delta^a_k[n]')^\circ$ and
  $\Delta^a_{n-k}[n]''\cong(\Delta^a_k[n]'')^\circ$.  Of course, the
  diagram
  \begin{displaymath}
    \xymatrix@C=8em@R=2.5em{
      \Delta^a_{n-k}[n]'\ar@{u(->}[r]^{\textstyle \subseteq_e}
      \ar[d]_{\cong} &
      \Delta^a_{n-k}[n]''\ar[d]^{\cong} \\
      (\Delta^a_k[n]')^\circ\ar@{u(->}[r]_{\subseteq_e} &
      (\Delta^a_k[n]'')^\circ }
  \end{displaymath}
  commutes for each $n$ and $1\leq k < n$, in other words the dual of
  each primitive t-extension is again (isomorphic to) a primitive
  t-extension.
  
  Now, since the dual operation on $\Strat$ is a (strict) involution,
  and in particular an equivalence of categories, it follows that a
  stratified set $A$ is orthogonal to each primitive t-extension if
  and only if its dual $A^\circ$ is orthogonal to each primitive
  t-extension dual. However, as we have seen, every primitive
  t-extension dual is in fact (isomorphic to) a primitive t-extension,
  so it follows that $A$ is pre-complicial if and only if $A^\circ$ is
  pre-complicial.
  
  Equivalently, we have demonstrated that the class of t-extensions is
  closed under the action of applying the dual functor $\arrow
  (\mathord{-})^\circ:\Strat->\Strat.$, in other words $\arrow
  f:X->Y.$ is a t-extension if and only if $\arrow
  f^\circ:X^\circ->Y^\circ.$ is a t-extension.
\end{obs} 

\begin{obs}[a simple explicit characterisation of
  t-extensions]\label{t-ext.expl.obs} It is clear from the
  construction of observation~\ref{precomp.refl} that the underlying
  simplicial map of any t-extension is an isomorphism in $\Simp$.  In
  general, we will work with t-extensions which are inclusions
  $\overinc\subseteq_e:X->Y.$ associated with an entire subset
  $X\subseteq_e Y$ (cf.\ notation~\ref{substrat.reg}).  Such an
  inclusion is a t-extension if and only if there exists a sequence of
  entire subsets $X=X^{(0)}\subseteq_e X^{(1)}
  \subseteq_e\cdots\subseteq_e X^{(n)}\subseteq_e\cdots$ of $Y$ such
  that $Y=\bigcup_{i=0}^\infty X^{(i)}$ and for each $i\geq 0$ and
  each simplex $x\in tX^{(i+1)}\smallsetminus tX^{(i)}$ there exists a
  $w\in Y$ and a $k\in\mathbb{N}$ satisfying the conditions
  \begin{itemize}
  \item $w$ is $k$-admissible in $X^{(i)}$,
  \item $w\cdot\face_{k-1}$ and $w\cdot\face_{k+1}$ are both in
    thin in $X^{(i)}$, and
  \item $x=w\cdot\face_k$
  \end{itemize}
  (in which case we say that {\em $w$ witnesses the extension of
    thinness to $x$}).
\end{obs}

\subsection{Tensor Products of Pre-Complicial Sets}

\begin{defn}\label{cleaves} 
  Suppose that $X$ and $Y$ are stratified sets then we say that a
  partition $p,q\in\mathbb{N}$ of $r\in\mathbb{N}$ {\em cleaves\/} an
  $r$-simplex $\pair<x;y>$ of $X\times Y$ if either the $p$-simplex
  $x\cdot\partinj^{p,q}_1$ is thin in $X$ or the $q$-simplex
  $y\cdot\partinj^{p,q}_2$ is thin in $Y$.
\end{defn}

\begin{defn}[the lax Gray-tensor of stratified sets]
  \label{lax.gray} We may define a bifunctor
  \begin{displaymath}
    \xymatrix@C=10pc{
      {\Strat\times\Strat} \ar[r]^{\displaystyle\otimes} & {\Strat}}
  \end{displaymath}
  called the {\em lax Gray-tensor} of stratified sets, by
  letting
  \begin{displaymath}
    X\otimes Y \stackrel{\text{def}}=
    \pair<X\times Y;t(X\otimes Y)>
  \end{displaymath}
  for each pair of stratified sets $X$ and $Y$, where:
  \begin{equation}\label{tens.thin.def}
    t(X\otimes Y) \defeq
    \left\{ \pair<x;y>\in X\times Y \left|
        \begin{array}[c]{l}
          \parbox[l]{7cm}{all partitions $p,q\in\mathbb{N}$ of 
            $r=\dim\pair<x;y>$ cleave $\pair<x;y>$.}
        \end{array}\right.\right\}
  \end{equation}
  The Gray-tensor $\arrow f\otimes g:X\otimes Y->X'\otimes Y'.$ of a
  pair of stratified maps $\arrow f:X->X'.$ and $\arrow g:Y->Y'.$ is
  simply the stratified map whose underlying simplicial map is the product
  $\arrow f\times g: X\times Y-> X'\times Y'.$, in other words
  $(f\otimes g)\pair<x;y> = \pair<f(x);g(y)>$.
\end{defn}

Before proving that this bifunctor is well defined, we first establish
the following simple lemma:

\begin{lemma}[simplices in $X\otimes Y$ with degenerate
  ordinates]\label{deg.ord} Let $X$ and $Y$ be stratified
  sets and suppose that $\pair<x;y>$ is an $r$-simplex in $X\times Y$
  and $p,q\in\mathbb{N}$ is an arbitrary partition of $r$. We have:
  \begin{enumerate}[(a)]
  \item\label{deg.xord} If $x$ is degenerate at $u$ then our partition
    cleaves $\pair<x;y>$ whenever $u < p$.
  \item\label{deg.yord} If $y$ is degenerate at $v$ then our partition
    cleaves $\pair<x;y>$ whenever $p \leq v$.
  \item\label{deg.ord.both} If there exists $u\leq v$ such that $x$ is
    degenerate at $u$ and $y$ is degenerate at $v$ then $\pair<x;y>$
    is in $t(X\otimes Y)$.
  \end{enumerate}
\end{lemma}

\begin{proof}\hfil\vspace{1ex}

\noindent(\ref{deg.xord}) We have $\dom(\partinj^{p,q}_1)=[p]$ so
$u$ and $u+1$ are both elements of the domain of this simplicial
operator, by our assumption that $0\leq u< p$, furthermore
$\partinj^{p,q}_1(u)=u$ and $\partinj^{p,q}_1(u+1)=u+1$. But $x$ is
degenerate at $u$ and it follows, by
observation~\ref{degen.obs}(\ref{degen.obs.c}), that
$x\cdot\partinj^{p,q}_1$ is degenerate, and thus thin, in $X$ so
our partition cleaves $\pair<x;y>$ as required.  \vspace{1ex}
  
\noindent(\ref{deg.yord}) The argument to establish this point is dual
to that of (\ref{deg.xord}).  We have $\dom(\partinj^{p,q}_2)
=[q]=[r-p]$ so $v'=v-p$ and $v'+1$ are both elements of the
domain of this simplicial operator, by our assumption that $p\leq v
< r$), furthermore $\partinj^{p,q}_2(v')=v'+p=v$ and
$\partinj^{p,q}_2(v')=v+1$. But $y$ is degenerate at $v$ and it follows, by
observation~\ref{degen.obs}(\ref{degen.obs.c}), that
$y\cdot\partinj^{p,q}_2$ is degenerate, and thus thin, in $X$
so our partition cleaves $\pair<x;y>$ as required.
\vspace{1ex}
  
\noindent(\ref{deg.ord.both}) Under the condition given we may infer,
by parts (\ref{deg.xord}) and (\ref{deg.yord}), that our partition
cleaves $\pair<x;y>$ whenever $p>u$ or $p\leq v$.  However, since
$u\leq v$ also holds we know that one or other of these conditions
will always apply. It follows that every partition will cleave
$\pair<x;y>$ which is therefore an element of $t(X\otimes Y)$ as
stated.
\end{proof}

\begin{proof} (that $\otimes$ is a well defined bifunctor as stated in
  definition~\ref{lax.gray}) We need to demonstrate that: \vspace{1ex}

  \noindent 1) {\bf\boldmath For any pair of stratified sets
    $X$ and $Y$ then $X\otimes Y$ is a well defined stratified set.}
  Clearly no 0-simplex of $X\times Y$ can be in $t(X\otimes Y)$ (since
  $tX$ and $tY$ contain no 0-simplices) therefore all we need show is
  that any degenerate simplex in $X\times Y$ is in $t(X\otimes Y)$. So
  suppose that $\pair<x;y>$ is a degenerate $r$-simplex in $X\times Y$
  then, by observation~\ref{degen.obs}(\ref{degen.obs.a}) and the
  point-wise definition of the simplicial structure of $X\times Y$, we
  know that this holds if and only if there is some $k$ such that both
  $x$ and $y$ are degenerate at $k$.  So $\pair<x;y>$ satisfies the
  conditions of lemma~\ref{deg.ord}(\ref{deg.ord.both}) with $u=v=k$,
  by which result we know that it is an element of $t(X\otimes Y)$ as
  required.\vspace{1ex}

  \noindent 2) {\bf\boldmath For any pair of stratified maps 
    $\arrow f:X->X'.$ and $\arrow g:Y->Y'.$ the simplicial map
    $f\otimes g$ is a well defined stratified map from $X\otimes Y$ to
    $X'\otimes Y'$.} All we need to show is that $f\otimes g$
  preserves thinness. So, suppose that the $r$-simplex $\pair<x;y>$ of
  $X\times Y$ is in $t(X\otimes Y)$ and that $p,q\in\mathbb{N}$ is an
  arbitrary partition of $r$, then we know, by defining equation
  (\ref{tens.thin.def}), that this partition cleaves $\pair<x;y>$ so
  either
  \begin{itemize}
  \item $x\cdot\partinj^{p,q}_1$ is in $tX$, in which case, since
    $f$ is a stratified map, it follows that
    $f(x)\cdot\partinj^{p,q}_1 = f(x\cdot\partinj^{p,q}_1)$ is in
    $tX'$, or
  \item $y\cdot\partinj^{p,q}_2$ is in $tY$, in which case, since $g$
    is a stratified map, it follows that $g(y)\cdot\partinj^{p,q}_2 =
    f(y\cdot\partinj^{p,q}_2)$ is in $tY'$.
  \end{itemize}
  In either case, our partition cleaves the simplex $(f\otimes
  g)\pair<x;y> = \pair<f(x);g(y)>$ in $X'\times Y'$ and therefore,
  since this partition was chosen arbitrarily, it follows that
  $(f\otimes g)\pair<x;y>$ is in $t(X'\otimes Y')$ as required.
  \vspace{1ex}

  \noindent 3) {\bf\boldmath $\otimes$ is functorial.} This 
  follows trivially from the functoriality of $\times$ as a bifunctor
  from $\Simp\times\Simp$ to $\Simp$.
\end{proof}

\begin{cor}[of lemma~\ref{deg.ord}]\label{tensor.simp.set} In a tensor
  product of standard simplices $\Delta[n]\otimes\Delta[m]$, a simplex
  $\pair<\alpha;\beta>$ is thin if and only if it satisfies the
  condition of lemma~\ref{deg.ord}(\ref{deg.ord.both}).
\end{cor}

\begin{proof}
  The ``if'' direction is simply
  lemma~\ref{deg.ord}(\ref{deg.ord.both}).  For the ``only if''
  direction, suppose that the $r$-simplex $\pair<\alpha;\beta>$ is
  thin in $\Delta[n]\otimes\Delta[m]$. Pick the partition $p,q$ of $r$
  by selecting $p\in[r]$ to be the largest integer such that
  $\alpha(0)<\alpha(1)<\cdots<\alpha(p)$ and letting $q=r-p$.  Of
  course, by our choice of $p$, the $p$-simplex
  $\alpha\circ\partinj^{p,q}_1$ is non-degenerate and thus non-thin in
  $\Delta[n]$, since this has the minimal stratification, in which
  only degenerate simplices are thin. However we know that our
  partition must cleave $\pair<\alpha;\beta>$, by the definition of
  thinness for the tensor product $\Delta[n]\otimes\Delta[m]$, and
  thus $\beta\circ\partinj^{p,q}_2$ must be thin in $\Delta[m]$. It
  follows that this simplex must be degenerate, again by the
  minimality of the stratification of $\Delta[m]$, so there must exist
  some $0\leq l<q$ such that $\beta\circ\partinj^{p,q}_2(l)=
  \beta\circ\partinj^{p,q}_2(l+1)$ equivalently, using the definition
  of $\partinj^{p,q}_2$, we see that there exists some $p\leq k< r$
  (which is related to $l$ by $k=l+p$) such that $\beta$ is degenerate
  at $k$.  Notice that this also implies that $p<r$ (since $0<q$ and
  $p+q=r$) from which we may infer, from the maximality clause in the
  definition of $p$, that $\alpha$ is degenerate at $p$ (i.e.
  $\alpha(p)=\alpha(p+1)$ otherwise $p$ would not be maximal for the
  stated property). So we have found a pair $p\leq k$ satisfying the
  condition of lemma~\ref{deg.ord}(\ref{deg.ord.both}) as required.
\end{proof}

\begin{lemma}\label{monoidal.strat}
  Given stratified sets $X$, $Y$ and $Z$ the structural isomorphisms
  of the cartesian category $\triple<\Simp;\times;\Delta[0]>$ extend to
  stratified isomorphisms:
  \begin{eqnarray}
    X\otimes(Y\otimes Z) & \cong & 
    (X\otimes Y)\otimes Z\label{strat.assoc} \\
    X\otimes\Delta[0] & \cong & X \label{strat.rident} \\
    \Delta[0]\otimes X & \cong & X\label{strat.lident}
  \end{eqnarray} 
  Collectively these provide the structural isomorphisms which make
  the triple $\triple<\Strat;\otimes;\Delta[0]>$ into a monoidal
  category.  Furthermore, the canonical isomorphism
  $(\forget(X)\times\forget(Y))^\circ \cong
  \forget(X)^\circ\times\forget(Y)^\circ$ in $\Simp$ extends to a
  stratified isomorphism:
  \begin{equation}\label{sym.strat}
    (X\otimes Y)^\circ \cong
    Y^\circ\otimes X^\circ
  \end{equation}
\end{lemma}

\begin{proof} The forgetful functor $\arrow \forget:\Strat->\Simp.$
  is faithful and $\otimes$ is defined so that $\forget$ commutes with it
  and $\times$, in other words the diagram of functors
  \begin{displaymath}
    \xymatrix@R=2em@C=4em{
      {\Strat\times\Strat}\ar[r]^-{\textstyle\otimes}
      \ar[d]_{\textstyle \forget\times \forget} & {\Strat}\ar[d]^{\textstyle \forget} \\
      {\Simp\times\Simp}\ar[r]_-{\textstyle\times} & {\Simp} }
  \end{displaymath}
  commutes. So, assuming that we can show that the various structural
  isomorphisms associated with the cartesian structure on $\Simp$ extend
  as stated to stratified isomorphisms, it follows that they will
  trivially satisfy the various conditions required of the structural
  maps of a monoidal category $\triple<\Strat;\otimes;\Delta[0]>$,
  simply because they do so in $\triple<\Simp;\times;\Delta[0]>$. We
  examine each of these isomorphisms in turn: \vspace{1ex}

  \noindent{\bf Associativity.} An $n$-simplex $\pair<x;{\pair<y;x>}>$ of
  $\forget(X)\times(\forget(Y)\times \forget(Z))$ is in the subset
  $t(X\otimes(Y\otimes Z))$ if and only if for all partitions $p,u$ of
  $n$ either
  \begin{itemize}
  \item $x\cdot\partinj_1^{p,u}$ is in $tX$, or
  \item $\pair<y;z>\cdot\partinj_2^{p,u}=
    \pair<y\cdot\partinj_2^{p,u};z\cdot\partinj_2^{p,u}>$ is in
    $t(Y\otimes Z)$.
  \end{itemize}
  Furthermore, $\pair<y\cdot\partinj_2^{p,u};
  z\cdot\partinj_2^{p,u}>$ is in $t(Y\otimes Z)$ if and only if for
  all partitions $q,r$ of $u$ either
  \begin{itemize}
  \item $(y\cdot\partinj_2^{p,u})\cdot \partinj_1^{q,r} =
    y\cdot(\partinj_2^{p,q+r}\circ \partinj_1^{q,r})$ is in $tY$,
    or
  \item $(z\cdot\partinj_2^{p,u})\cdot \partinj_2^{q,r} =
    z\cdot(\partinj_2^{p,q+r}\circ \partinj_2^{q,r})$ is in $tZ$.
  \end{itemize}
  By a dual argument, the corresponding $n$-simplex
  $\pair<{\pair<x;y>};z>$ in $(\forget(X)\times\forget(Y))\times
  \forget(Z)$ is in $t((X\otimes Y)\otimes Z)$ iff for all $p,q,r$
  with $p+q+r=n$ we have:
  \begin{itemize}
  \item $x\cdot(\partinj^{p+q,r}_1\circ\partinj^{p,q}_1)$ is in $tX$,
    or
  \item $y\cdot(\partinj^{p+q,r}_1\circ\partinj^{p,q}_2)$ is in $tY$,
    or
  \item $z\cdot\partinj^{p+q,r}_2$ is in $tZ$.
  \end{itemize}
  Applying the partition identities laid out in
  display~\eqref{part.ident} of notation~\ref{part.oper.def}, we see
  that corresponding pairs of simplices from these characterisations
  are actually equal in $X$, $Y$ and $Z$ respectively and it follows
  that $\pair<x;{\pair<y;z>}>$ is thin in $X\otimes (Y \otimes Z)$ if
  and only if $\pair<{\pair<x;y>};z>$ is thin in $(X \otimes Y)\otimes
  Z$.  In other words the associativity isomorphism of $\times$ on
  $\Simp$ preserves and reflects the stratifications associated with
  $\otimes$ and so extends to the stratified isomorphism of
  (\ref{strat.assoc}) as required.\vspace{1ex}

  \noindent{\bf Right and Left Identity.} Consider an $r$-simplex
  $\pair<x;y>$ in $X\otimes\Delta[0]$ and a partition
  $p,q$ of $r$. There are two cases:
  \begin{itemize}
  \item {\boldmath $q>0$} in which case $y\cdot\partinj^{p,q}_2$ would
    be a simplex of dimension greater than $0$ in the terminal
    simplicial set $\Delta[0]$.  However, all simplices of $\Delta[0]$
    of dimension greater than $0$ are degenerate and therefore thin in
    $\Delta[0]$. So, it follows that our partition cleaves
    $\pair<x;y>$.
  \item {\boldmath $q=0$} or equivalently $p=r$, in this case
    $x\cdot\partinj^{p,q}_1=x$ and $y\cdot\partinj^{p.q}_2$ is the
    unique 0-simplex of $\Delta[0]$ and is therefore not thin. It follows
    that our partition would cleave $\pair<x;y>$ if and only if $x$ is thin
    in $X$.
  \end{itemize}
  Quantifying over the partitions of $r$, it follows that $\pair<x;y>$
  is thin in $X\otimes\Delta[0]$ if and only if $x$ is thin in $X$. In
  other words, the right identity isomorphism $\triple<\Simp;\times;
  \Delta[0]>$ extends to the stratified isomorphism of
  (\ref{strat.rident}) and a dual argument demonstrates that the left
  identity isomorphism extends to the stratified isomorphism of
  (\ref{strat.lident}), as required.\vspace{1ex}

  \noindent{\bf Dual Symmetry.} Applying
  observation~\ref{simp.op.dual}, it is easily seen that a partition
  $p,q$ cleaves an $r$-simplex $\pair<y;x>$ in $Y^\circ\otimes
  X^\circ$ if and only if either
  \begin{itemize}
  \item $y*\partinj^{p,q}_1\defeq y\cdot (\partinj^{p,q}_1)^\circ$ is in
    $tY^\circ$ but then, by observation~\ref{simp.op.dual},
    $(\partinj^{p,q}_1)^\circ=\partinj^{q,p}_2$ and so this condition
    holds iff $y\cdot\partinj^{q,p}_2$ is in $tY$, or 
  \item $x*\partinj^{p,q}_2\defeq x\cdot (\partinj^{p,q}_2)^\circ$ is
    in $tX^\circ$ but then, by observation~\ref{simp.op.dual},
    $(\partinj^{p,q}_2)^\circ=\partinj^{q,p}_1$ and so this condition
    holds iff $x\cdot\partinj^{q,p}_1$ is in $tX$
  \end{itemize}
  which says no more nor less than that $q,p$ cleaves the simplex
  $\pair<x;y>$ of $X\otimes Y$. Quantifying this result over the
  partitions of $r$, we see that $\pair<y;x>$ is in $t(Y^\circ\otimes
  X^\circ)$ if and only if $\pair<x;y>$ is in $t(X\otimes
  Y)^\circ\defeq t(X\otimes Y)$. In other words, the dual symmetry
  isomorphism on $\Simp$ extends to the stratified isomorphism of
  (\ref{sym.strat}) as required.
\end{proof}

\begin{obs}[tensor products of regular maps]
  \label{tens.thin.refl}
  If the stratified maps $\arrow f:X->X'.$ and $\arrow g:Y->Y'.$ are both
  regular then so is their tensor product $\arrow f\otimes
  g:X\otimes Y->X'\otimes Y'.$. It follows that a tensor product of
  regular subsets $X\subseteq_r X'$ and $Y\subseteq_r Y'$ is again a
  regular subset $X\otimes Y\subseteq_r X'\otimes Y'$.
\end{obs}

\begin{proof} Suppose that the $r$-simplex $\pair<f(x);g(y)>$ is thin in
  $X'\otimes Y'$ and that $p,q$ is a partition of $r$. By the
  definition of $t(X'\otimes Y')$, this cleaves $\pair<f(x);g(y)>$,
  so either
  \begin{itemize}
  \item $f(x\cdot\partinj^{p,q}_1)=f(x)\cdot\partinj^{p,q}_1$ is
    is in $tX'$, in which case $x\cdot\partinj^{p,q}_1$ is in $tX$
    since $f$ is regular, or
  \item $g(y\cdot\partinj^{p,q}_2)=g(y)\cdot\partinj^{p,q}_2$ is
    is in $tY'$, in which case $y\cdot\partinj^{p,q}_2$ is in $tY$
    since $g$ regular.
  \end{itemize}
  In either case, we have shown that our partition also cleaves
  $\pair<x;y>\in X\times Y$ and it follows, by quantifying over such
  partitions, that $\pair<x;y>$ is in $t(X\otimes Y)$ as required.
\end{proof}

Recall that entire maps (resp.\ inclusions) are defined to be those
stratified maps whose action on underlying simplicial sets is
surjective (resp.\ injective). It therefore goes without saying, since
$\otimes$ acts as the cartesian product on underlying sets, that these
classes of stratified maps are also closed under (pre-)tensoring.

\begin{lemma}\label{tensor.Th} 
  If $X$ and $Y$ are stratified sets and $n,m\geq 0$ are fixed
  integers then we have entire subset inclusions 
  \begin{equation}\label{tensor.Th.ineq}
    \begin{aligned}
      \Th_{n+m}(X\otimes Y) &{}\subseteq_e\Th_n(X) \otimes \Th_m(Y) &&
      \text{and} \\
      \Th_n(X)\otimes\Th_m(Y) &{}\subseteq_e\Th_{\min(n,m)}( X\otimes Y)&& 
    \end{aligned}
  \end{equation}
  furthermore the argument used to establish the second of these also
  shows that:
  \begin{equation*}
    \begin{aligned}
      X\otimes\Th_m(Y) &{}\subseteq_e \Th_m(X\otimes Y) && \text{and} \\
      \Th_n(X)\otimes Y &{}\subseteq_e \Th_n(X\otimes Y) &&
    \end{aligned}
  \end{equation*}
  In particular, applying the inclusions in
  display~(\ref{tensor.Th.ineq}) it follows that $\Th_0(X\otimes
  Y)=\Th_0(X) \otimes\Th_0(Y)$.
\end{lemma}

\begin{proof}
  Of course $X\otimes Y\subseteq_e \Th_n(X)\otimes\Th_m(Y)$, so in
  order to prove the first of these inclusions all we need to show is
  that every simplex of $\Th_n(X)\otimes\Th_m(Y)$ of dimension $>n+m$
  is thin.  To this end, let $r>n+m$ and consider an $r$-simplex
  $\pair<x;y>$ in $\Th_n(X)\otimes\Th_m(Y)$. If $p,q$
  is a partition of $r$ then, since $p+q=r>n+m$, we know that
  either $p>n$ or $q>m$. In the first case
  $x\cdot\partinj^{p,q}_1$ is of dimension $p>n$, and thus it is
  thin in $\Th_n(X)$, and in the latter $y\cdot\partinj^{p,q}_2$ is of
  dimension $q>m$, which implies that it is thin in $\Th_m(Y)$;
  either way we see that our partition cleaves $\pair<x;y>$. It follows that
  $\pair<x;y>$ is thin in $\Th_n(X)\otimes \Th_m(Y)$, as required.
  
  Now we also know that $X\otimes Y\subseteq_e\Th_{\min(n,m)}(X\otimes
  Y)$, so in order to prove the second inclusion all we need to do is
  show that any $r$-simplex $\pair<x;y>$ which is thin in
  $\Th_n(X)\otimes \Th_m(Y)$ and not thin in $X\otimes Y$ is actually
  also thin in $\Th_{\min(n,m)}(X\otimes Y)$. However, for any such
  simplex it is clear that there must exist some partition
  $p,q$ of $r$ for which either $x\cdot\partinj^{p,q}_1$ is
  thin in $\Th_n(X)$ but not in $X$ or $y\cdot\partinj^{p,q}_2$ is
  thin in $\Th_m(Y)$ but not in $Y$. In the first of these cases we
  have $p>n$ and in the second we have $q>m$, but $r=p+q$ so
  in either case we have $r>\min(n,m)$. It follows, therefore, that
  any such simplex is thin in $\Th_{\min(n,m)}(X\otimes Y)$ as
  required. Everything else follows trivially.
\end{proof}

\subsection{Pre-Tensors and Preservation of t-Extensions}

\begin{obs}[some thin simplices in $X\otimes Y$]
  \label{med.cyl} While $\otimes$ provides us with a
  rich and easily defined monoidal structure on $\Strat$ it is
  unfortunately deficient in one important respect - it is not
  biclosed.  However, the remainder of this section will demonstrate
  that it may be reflected to a monoidal structure on the full
  subcategory $\Precomp$ which {\bf is} biclosed.
  
  To do this we need to analyse the thin simplices of a tensor product
  $X\otimes Y$ in some greater detail. In fact there
  are two quite distinct classes of simplices in such tensors which
  turn out to be of particular importance:\vspace{1ex}

  \noindent{\bf Mediator Simplices.} An $r$-simplex $\pair<x;y>$ of
  $X\otimes Y$ is called a {\em mediator simplex} if there exists some
  $1\leq k< r$ such that $x = x\cdot(\face_{k-1}\circ\degen_{k-1})$
  and $y = y\cdot(\face_k\circ\degen_k)$.  In this case we say that
  $k$ {\em witnesses\/} the fact that $\pair<x;y>$ is a mediator
  simplex. Notice that if this simplex is non-degenerate then the
  mediation condition implies immediately that $x$ cannot be
  degenerate at $k$ and $y$ cannot be degenerate at $k-1$.
  
  An equivalent way of stating this condition is to say that we may
  find decompositions $x=x'\cdot\alpha$ and $y=y'\cdot\beta$ (for some
  $x'\in X$, $y'\in Y$ and $\alpha,\beta\in\Delta$) which have
  $\alpha(k-1)=\alpha(k)$ and $\beta(k)=\beta(k+1)$.

  \begin{figure}
  \begin{displaymath}
    \def\vertchar{\scriptscriptstyle\bullet}
    \xymatrix@!0@C=1.5em@R=1.5em{
       & *[o]{\vertchar} & *[o]{\vertchar} & *[o]{\vertchar} & *[o]{\vertchar} & 
      *[o]{\vertchar} & *[o]{\vertchar} \\ 
       & *[o]{\vertchar} & *[o]{\vertchar} & *[o]{\vertchar}
      \save [].[d]!C*++[o][F.]\frm{}\restore
      \save [].[r]!C*++[o][F.]\frm{}\restore
      \save []+<10em,-4em>*{\text{path vertex $k$}}\ar[]\restore
      & *[o]{\vertchar} & *[o]{\vertchar} & *[o]{\vertchar} \\ 
      & *[o]{\vertchar} & *[o]{\vertchar} & *[o]{\vertchar} & *[o]{\vertchar} & 
      *[o]{\vertchar} & *[o]{\vertchar} \\ 
      & *[o]{\vertchar} & *[o]{\vertchar} & *[o]{\vertchar} & *[o]{\vertchar} & 
      *[o]{\vertchar} & *[o]{\vertchar} \\ 
      [m]\ar[uu] & *[o]{\vertchar} & *[o]{\vertchar} & *[o]{\vertchar} & *[o]{\vertchar} & 
      *[o]{\vertchar} & *[o]{\vertchar} \\ 
      & *[o]{\vertchar} & *[o]{\vertchar} & *[o]{\vertchar} & *[o]{\vertchar} & 
      *[o]{\vertchar} & *[o]{\vertchar} \\ 
      & *[o]{\vertchar}\ar@{-}'[u]'[uruu]'[uruuur]'[uruuuru]'[uruuurur]
      '[uruuururur][uruuurururr] & 
      *[o]{\vertchar} & *[o]{\vertchar} & *[o]{\vertchar} & 
      *[o]{\vertchar} & *[o]{\vertchar} \\
       & & & [n] \ar[rr] & & & 
      }
    \end{displaymath}
    \caption{A mediator simplex in $\Delta[n]\otimes\Delta[m]$.}
    \label{mediator.pic}
  \end{figure}
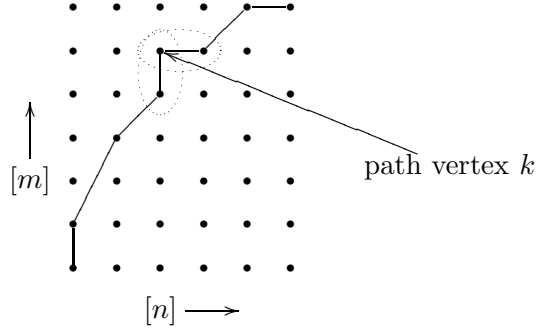 
  
  To visualise this {\em mediation condition}, consider an $r$-simplex
  $\pair<\alpha;\beta>$ in the tensor $\Delta[n]\otimes\Delta[m]$
  recall, from observation~\ref{nerves.po}, that we may consider such
  a simplex to be a single order preserving map
  $\overarr:[r]->[n]\times[m].$. We can draw this as a path in
  $[n]\times[m]$, as depicted in figure~\ref{mediator.pic}, every step
  of which is a move toward the upper-right corner of that grid.
  Pictorially the mediation condition at the highlighted path point
  $k$ states that our arrival there is achieved by a pure upward
  movement and our departure by a pure rightward one.
  
  {\em Most importantly}, all mediator simplices are thin in a tensor
  $X\otimes Y$; a fact which follows directly from
  lemma~\ref{deg.ord}(\ref{deg.ord.both}) since the mediation
  condition is simply a special case of the condition given there.

  \vspace{1ex}

  \noindent{\bf Cylinders.} An $r$-simplex $\pair<x;y>$ in
  $X\times Y$ is said to be a {\em cylinder\/} if there exists a
  partition $p,q$ of $r$ such that $x = x\cdot(
  \partinj^{p,q}_1\circ\partproj^{p,q}_1)$ and $y =
  y\cdot(\partinj^{p,q}_2\circ\partproj^{p,q}_2)$.  In this case we
  say that our partition {\em witnesses\/} the fact that $\pair<x;y>$
  is a cylinder.
  
  The fact that $\partinj^{p,q}_i$ is a right inverse of
  $\partproj^{p,q}_i$ (for $i=1,2$, see notation~\ref{part.oper.def})
  implies that our partition witnesses $\pair<x;y>$ as a cylinder if and
  only if there exists an $p$-simplex $x'\in X$ and an $q$-simplex
  $y'\in Y$ such that $x=x'\cdot\partproj^{p,q}_1$ and
  $y=y'\cdot\partproj^{p,q}_2$.

  \begin{figure}[b]
  \begin{displaymath}
    \def\vertchar{\scriptscriptstyle\bullet}
    \xymatrix@!0@C=1.5em@R=1.5em{
       & *[o]{\vertchar} & *[o]{\vertchar} & *[o]{\vertchar} & *[o]{\vertchar} & 
      *[o]{\vertchar} & *[o]{\vertchar} \\ 
       & *[o]{\vertchar} & *[o]{\vertchar} & *[o]{\vertchar}
      & *[o]{\vertchar} & *[o]{\vertchar}\ar@{.}[llll] & *[o]{\vertchar} \\ 
      & *[o]{\vertchar} & *[o]{\vertchar} & *[o]{\vertchar} & *[o]{\vertchar} & 
      *[o]{\vertchar} & *[o]{\vertchar} \\ 
      & *[o]{\vertchar} & *[o]{\vertchar} & *[o]{\vertchar} & *[o]{\vertchar} & 
      *[o]{\vertchar} & *[o]{\vertchar} \\ 
      [m]\ar[uu] & *[o]{\vertchar} & *[o]{\vertchar} & *[o]{\vertchar} & *[o]{\vertchar} & 
      *[o]{\vertchar} & *[o]{\vertchar} \\ 
      & *[o]{\vertchar}\save[].[uuuu]!C*+[F.]\frm{}\restore
      & *[o]{\vertchar} 
      \ar@{-}'[rrr][rrruuuu]\ar@{.}[l]\ar@{.}[d]
      & *[o]{\vertchar} & *[o]{\vertchar} &
      *[o]{\vertchar}\ar@{.}[d] & *[o]{\vertchar} \\ 
      & *[o]{\vertchar} &
      *[o]{\vertchar}\save[].[rrr]!C*+[F.]\frm{}\restore & 
      *[o]{\vertchar} & *[o]{\vertchar} & 
      *[o]{\vertchar} & *[o]{\vertchar} \\
       & & & [n] \ar[rr] & & & 
      }
    \end{displaymath}
    \caption{A cylinder in $\Delta[n]\otimes\Delta[m]$.}
    \label{cylinder.pic}
  \end{figure}
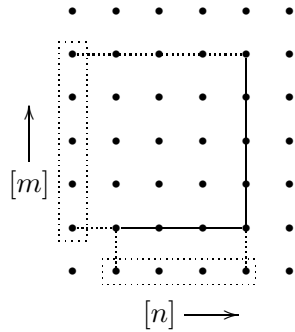
  
  To understand what a cylinder ``looks'' like we can again consider
  the cylinders of $\Delta[n]\otimes\Delta[m]$. With reference to
  figure~\ref{cylinder.pic} we see that such a cylinder
  $\pair<\alpha;\beta>$ corresponds to a path consisting of a
  (possibly empty) series of pure rightward moves followed by a
  (possibly empty) series of pure upward moves.  The projections of
  this simplex onto the x- and y-axes, as marked in the figure,
  correspond to the simplices
  $\alpha\circ\partinj^{p,q}_1\in\Delta[n]_p$ and
  $\beta\circ\partinj^{p,q}_2\in\Delta[m]_q$ respectively.
  
  Now, suppose that the $r$-simplex $\pair<x;y>$ is witnessed as a
  cylinder by a partition $p,q$ and suppose that $s,t$ is any other
  partition. Observe that:
  \begin{itemize}
  \item $\partproj^{{p,q}}_1(p+1)=\partproj^{p,q}_1(p)$ and
    $x=(x\cdot\partinj^{{p,q}}_1)\cdot\partproj^{{p,q}}_1$ which is
    therefore degenerate at $p$, so applying
    lemma~\ref{deg.ord}(\ref{deg.xord}) we see that ${s,t}$ cleaves
    $\pair<x;y>$ whenever $p < s$.
  \item $\partproj^{{p,q}}_2(p-1)=\partproj^{{p,q}}_2(p)$ and
    $y=(y\cdot\partinj^{{p,q}}_2)\cdot\partproj^{{p,q}}_2$ which is
    therefore degenerate at $p-1$, so applying
    lemma~\ref{deg.ord}(\ref{deg.yord}) with we see that ${s,t}$
    cleaves $\pair<x;y>$ whenever $s\leq p-1$.
  \end{itemize}
  It follows that $\pair<x;y>$ is in $t(X\otimes Y)$ if and only if
  ${p,q}$ cleaves it.\vspace{0.5ex}
  
  Consequently, we say that an $r$-simplex $\pair<x;y>$ of
  $X\otimes Y$ is a {\em crushed cylinder\/} if it is
  witnessed as a cylinder by a partition which also cleaves it. Our
  last result simply establishes that a cylinder
  is thin in $X\otimes Y$ if and only if it is
  a crushed cylinder in there.
\end{obs}

\begin{defn}[the pre-tensor of stratified sets]
  \label{pre.tens} We may define a bifunctor
  \begin{displaymath}
    \xymatrix@C=10pc{
      {\Strat\times\Strat} \ar[r]^{\displaystyle\pretens} & {\Strat}}
  \end{displaymath}
  called the {\bf pre-tensor} of stratified sets, by
  letting
  \begin{displaymath}
    X\pretens Y \stackrel{\text{def}}=
    \pair<X\times Y; t(X\pretens Y)>
  \end{displaymath}
  where a simplex $\pair<x;y>$ is in $t(X\pretens Y)$ if
  \begin{itemize}
  \item it is degenerate, or
  \item it is a mediator simplex, or
  \item it is a crushed cylinder.
  \end{itemize}
  Given a pair of stratified maps $\arrow f:X-> X'.$ and $\arrow
  g:Y->Y'.$, the stratified map $f\pretens g$ is defined to have the
  product map $\arrow f\times g: X\times Y-> X'\times Y'.$ as its
  underlying simplicial map, in other words $(f\pretens g)\pair<x;y> =
  \pair<f(x);g(y)>$.
  
  Furthermore, by observation~\ref{med.cyl}, we know that $X \pretens
  Y$ is an entire stratified subset of $X\otimes Y$ and so we have an
  entire inclusion map:
  \begin{equation}\label{pt.incl}
    \let\labelstyle=\textstyle
    \xymatrix@C=8em{{X\pretens Y}
      \ar@{u(->}[r]^{\subseteq_e} & {X\otimes Y}}
  \end{equation}
  The canonical isomorphism $(X\otimes Y)^\circ\cong Y^\circ\otimes
  X^\circ$ restricts along these inclusions to give an isomorphism
  $(X\pretens Y)^\circ\cong Y^\circ\pretens X^\circ$.
\end{defn}

\begin{proof} (that $\pretens$ is a well defined bifunctor as stated)
  We showed earlier that $t(X\otimes Y)$ contains all degenerate
  simplices (definition~\ref{lax.gray}), all mediator simplices and
  all crushed cylinders (observation~\ref{med.cyl}); it follows that
  $t(X\pretens Y)$ is a subset of $t(X\otimes Y)$. But $t(X\otimes Y)$
  contains no 0-simplices, by the proof attached to
  definition~\ref{lax.gray}, furthermore every degenerate simplex of
  $X\times Y$ is in $t(X\pretens Y)$, by definition, so it is
  certainly true that $\pair<X\times Y;t(X\pretens Y)>$ is a well
  defined stratified subset of $X\otimes Y$.
  
  Inspecting the defining properties of mediator simplices and
  crushed cylinders, we see that both concepts are defined in terms
  of operations and properties which are preserved by the stratified
  maps $f$ and $g$. It follows that if the $r$-simplex
  $\pair<x;y>$ is a mediator simplex witnessed by $k$ (respectively a
  crushed cylinder witnessed by $p,q$) then so is $(f\pretens
  g)\pair<x;y> = \pair<f(x);g(y)>$. From this fact, it follows that
  $f\pretens g$ preserves thinness as required.
\end{proof}

\begin{obs}[why introduce $\pretens$]\label{why.pretens}
  While the tensor $\otimes$ possesses many properties that $\pretens$
  lacks, for instance it is canonically part of a genuine monoidal
  structure on $\Strat$, it is also deficient in one major respect -
  the left and right tensoring functors $X\otimes\mathord{-}$ and
  $\mathord{-}\otimes Y$ do not preserve colimits. Consequently,
  $\otimes$ cannot be the tensor of a biclosed monoidal structure on
  $\Strat$. 

  This deficiency can be traced back to the very definition of the
  thin simplices in $X\otimes Y$, which relies upon a
  property that must hold for {\bf every} partition of the dimension
  of a given simplex.
  
  On the other hand, $\pretens$ does not suffer from this problem and
  consequently, as we shall show later on, the left and right
  pre-tensors $X\pretens\mathord{-}$ and $\mathord{-}\pretens Y$ {\bf
    do} preserve the colimits of $\Strat$. Furthermore, the more
  ``local'' explicit description of the thin simplices of $X\pretens
  Y$ often makes it easy to demonstrate that a given simplicial map
  $\arrow f:X\times Y->A.$ extends to a stratified map $\arrow
  f:X\pretens Y->A.$.
  
  In this context, the next lemma provides the most cogent
  justification for introducing the pre-tensor $\pretens$. In essence
  it says that ``pre-complicial sets do not recognise the difference
  between $X\otimes Y$ and $X\pretens Y$''.
\end{obs}

\begin{lemma}[mediator simplices witness thinness extension]
  \label{mediator.admissible}
  Suppose that the $r$-simplex $\pair<x;y>$ is a mediator simplex in
  $X\otimes Y$, as witnessed by some $1\leq k < r$, then
  \begin{enumerate}
  \item\label{mediator.admissible.1} $\pair<x;y>$ is $k$-admissible in
    $X\pretens Y$, and
  \item\label{mediator.admissible.2} if $\pair<x;y>\cdot\face^r_k$ is
    thin in $X\otimes Y$ then so are the faces
    $\pair<x;y>\cdot\face^r_{k+1}$ and $\pair<x;y>\cdot\face^r_{k-1}$.
  \end{enumerate}
\end{lemma}

\begin{proof}
  To establish part~(\ref{mediator.admissible.1}) we need to show that
  $\pair<x;y>\cdot\mu$ is thin in $X\pretens Y$ for each face operator
  $\arrow \mu:[m]->[r].$ with $k-1$, $k$ and $k-1$ in $\im(\mu)$. So,
  assume that $\mu$ is such a face operator and let $l\in[m]$ be the
  unique integer with $\mu(l)=k$, for which we know that
  $\mu(l-1)=k-1$ and $\mu(l+1)=k+1$. Now it is easily demonstrated
  that $\face^r_{k-1}\circ\degen^{r-1}_{k-1}\circ\mu =
  \mu\circ\face^m_{l-1} \circ\degen^{m-1}_{l-1}$ and
  $\face^r_k\circ\degen^{r-1}_k\circ\mu =
  \mu\circ\face^m_l\circ\degen^{m-1}_l$ simply by considering cases
  and applying these equations along with the mediation condition for
  $\pair<x;y>$, we get
  \begin{equation*}
    \begin{split}
    x\cdot\mu & {} = 
    (x\cdot(\face^r_{k-1}\circ\degen^{r-1}_{k-1}))\cdot\mu =
    x\cdot(\face^r_{k-1}\circ\degen^{r-1}_{k-1}\circ\mu) \\
    & {} = x\cdot(\mu\circ\face^m_{l-1}\circ\degen^{m-1}_{l-1}) =
    (x\cdot\mu)\cdot(\face^m_{l-1}\circ\degen^{m-1}_{l-1}) \\
    y\cdot\mu & {} = 
    (y\cdot(\face^r_k\circ\degen^{r-1}_k))\cdot\mu =
    y\cdot(\face^r_k\circ\degen^{r-1}_k\circ\mu) \\
    & {} =  y\cdot(\mu\circ\face^m_l\circ\degen_l) =
    (y\cdot\mu)\cdot(\face^m_l\circ\degen_l)
    \end{split}
  \end{equation*}
  or, in other words, $l$ witnesses the fact that $\pair<x;y>\cdot\mu=
  \pair<x\cdot\mu;y\cdot\mu>$ is a mediator simplex. It follows, by
  the fact that all mediator simplices are thin in $X\pretens Y$, that
  $\pair<x;y>\cdot\mu$ is thin there as required.
  
  For part~(\ref{mediator.admissible.2}), it is only necessary to
  prove the stated result for $\pair<x;y>\cdot\face^r_{k+1}$, since
  the corresponding result for $\pair<x;y>\cdot\face^r_{k-1}$ may be
  obtained by applying the former result to the simplex $\pair<y;x>$
  in the dual $(X\otimes Y)^\circ\cong Y^\circ\otimes X^\circ$. So we
  assume that $\pair<x;y>\cdot\face^r_k=\pair<x\cdot\face^r_k;
  y\cdot\face^r_k>$ is thin in $X\otimes Y$ and demonstrate that the
  $(r-1)$-simplex $\pair<x;y>\cdot\face^r_{k+1}
  =\pair<x\cdot\face^r_{k+1};y\cdot\face^r_{k+1}>$ is also thin there
  by considering an arbitrary partition $p,q$ of $r-1$, for which
  there are two cases:
  \begin{itemize}
  \item {\boldmath $p\geq k$} applying the mediation condition we know
    that $x=x'\cdot\degen^{r-1}_{k-1}$ for some $(r-1)$-simplex $x'$
    so $x\cdot\face^r_{k+1}=x'\cdot(\degen^{r-1}_{k-1}\circ
    \face^r_{k+1})=x'\cdot(\face^{r-1}_k\circ \degen^{r-2}_{k-1})=
    (x'\cdot\face^{r-1}_k)\cdot \degen^{r-2}_{k-1}$ where the
    penultimate equality is an application of the appropriate
    simplicial identity. It follows that $x\cdot\face^r_{k+1}$ is
    degenerate at $k-1$ and so our assumption that $p\geq k$ allows us to
    apply lemma~\ref{deg.ord}(\ref{deg.xord}) to show that $p,q$
    cleaves $\pair<x;y>\cdot\face^r_{k+1}$.
  \item {\boldmath $p < k$} which implies that $\face^r_k
    \circ\partinj^{p,q}_1=\partinj^{p,q+1}=\face^r_{k+1}\circ
    \partinj^{p,q}_1$ and it follows that $(x\cdot\face^r_{k+1})
    \cdot\partinj^{p,q}_1=(x\cdot\face^r_k)\cdot\partinj^{p,q}_1$.
    Furthermore, applying the mediation condition on $y$ we see that
    $y\cdot\face^r_{k+1}=(y\cdot(\face^r_k\circ \degen^{r-1}_k))
    \cdot\face^r_{k+1}=y\cdot(\face^r_k\circ( \degen^{r-1}_k
    \circ\face^r_{k+1}))=y\cdot\face^r_k$ where the last equality
    follows by applying the appropriate simplicial identity and so, in
    particular, $(y\cdot\face^r_k)\cdot \partinj^{p,q}_2
    =(y\cdot\face^r_{k+1})\cdot\partinj^{p,q}_2$. From these two
    equalities, it clearly follows that our partition $p,q$ cleaves
    $\pair<x;y>\cdot\face^r_{k+1}$ iff it cleaves
    $\pair<x;y>\cdot\face^r_k$, which it does by the assumption that
    this latter simplex is thin in $X\otimes Y$.
  \end{itemize}
  So we have shown that every partition of $r-1$ cleaves
  $\pair<x;y>\cdot\face^r_{k+1}$ and thus that it is thin in $X\otimes
  Y$ as required.
\end{proof}

\begin{lemma}\label{(pre)tens.diff}
  If the $r$-simplex $\pair<x;y>=\pair<x'\cdot\alpha;y'\cdot\beta>$ is
  thin in $X\otimes Y$ and not thin in $X\pretens Y$ then there exists
  a $0< k\leq r$ such that $\alpha(k-1)<\alpha(k)$ and $\beta(k-1)<
  \beta(k)$.
\end{lemma}

\begin{proof}
  Define two integers $p$ and $q$ as follows
  \begin{itemize}
  \item let $p$ be the largest integer in $[r]$ such that
    $\alpha(0)<\alpha(1) < ... < \alpha(p)$ and $\beta(0) = \beta(1) =
    ... = \beta(p)$, and
  \item let $q$ be the largest integer in $[r-p]$ such that $\alpha(p)
    = \alpha(p+1) = ... = \alpha(p+q)$ and $\beta(p) < \beta(p+1) <
    ... < \beta(p+q)$,
  \end{itemize}
  and suppose, for a contradiction, that $p,q$ is a partition of $r$
  (that is $p+q=r$). Were this the case then the defining properties
  of $p$ and $q$ would imply that
  $\alpha=\alpha\circ\partinj^{p,q}_1\circ\partproj^{p,q}_1$ and
  $\beta=\beta\circ\partinj^{p,q}_2\circ\partproj^{p,q}_2$. It
  follows that $x\cdot(\partinj^{p,q}_1\circ\partproj^{p,q}_1) =
  (x'\cdot\alpha)\cdot(\partinj^{p,q}_1\circ\partproj^{p,q}_1) =
  x'\cdot(\alpha\circ\partinj^{p,q}_1\circ\partproj^{p,q}_1) =
  x'\cdot\alpha = x$ and
  $y\cdot(\partinj^{p,q}_2\circ\partproj^{p,q}_2) =
  (y'\cdot\beta)\cdot(\partinj^{p,q}_2\circ\partproj^{p,q}_2) =
  y'\cdot(\beta\circ\partinj^{p,q}_2\circ\partproj^{p,q}_2) =
  y'\cdot\beta = y$, and so our partition witnesses the fact that
  $\pair<x;y>$ is a cylinder. Of course we know, from our analysis in
  observation~\ref{med.cyl}, that a cylinder is thin in $X\otimes Y$ if
  and only if it is a crushed cylinder. But, all crushed cylinders are
  thin in $X\pretens Y$ (by definition) which contradicts our original
  assumption that $\pair<x;y>$ is {\bf not} thin in $X\pretens Y$, so
  it follows that $p+q<r$. 

  Now, let $k=p+q+1\in[r]$ for which we may show that:
  \begin{itemize}
  \item {\bf\boldmath $\alpha(k-1)<\alpha(k)$}. {\em Proof}.  Suppose,
    for a contradiction, that $\alpha(k-1)=\alpha(k)$; then either
    \begin{itemize}
    \item $\beta(k-1)=\beta(k)$, in which case $\pair<x;y>$ would be
      degenerate, which would contradict the assumption that it is not
      thin in $X\pretens Y$, or
    \item$\beta(k-1)<\beta(k)$, which would contradict the maximality
      clause in the definition of $q$ given above.
    \end{itemize}
    In either case we get a contradiction and it follows that
    $\alpha(k-1)<\alpha(k)$ as required.
    
  \item {\bf\boldmath $\beta(k-1)<\beta(k)$}. {\em Proof}. Suppose,
    for a contradiction, that $\beta(k-1)=\beta(k)$, then we already
    know that $\alpha(k-1)<\alpha(k)$ and so either:
    \begin{itemize}
    \item $q=0$, in which case $k-1=p$ and we get a contradiction of
      the maximality clause in the definition of $p$ given above.
    \item $q>0$, from which we may infer that
      $\alpha(k-2)=\alpha(k-1)$ which, when combined with the
      assumption that $\beta(k-1)=\beta(k)$, would imply that $k-1$
      witnesses that $\pair<x;y>$ is a mediator simplex, again
      contradicting the assumption that it is not thin in
      $X\pretens Y$.
    \end{itemize}
    In either case we get a contradiction and it follows that
    $\beta(k-1)<\beta(k)$ as required.\qedhere
  \end{itemize}
\end{proof}

\begin{lemma}\label{pretens.density} For each pair of stratified
  sets $X$ and $Y$ the entire inclusion map of display~(\ref{pt.incl})
  is a t-extension.
\end{lemma}

\begin{proof}
  Our intention is to prove this lemma by applying
  observation~\ref{t-ext.expl.obs}. To this end, we define a function
  $\arrow \Phi:\Delta[n]\times\Delta[m]->\mathbb{N}.$ by
  \begin{displaymath}
    \Phi\pair<\alpha;\beta>\defeq \sum_{i=0}^{\dim\pair<\alpha;\beta>}
    \alpha(i) + (m - \beta(i))
  \end{displaymath}
  and an increasing sequence of subsets $tA^{(0)}\subseteq tA^{(1)}\subseteq
  ...  \subseteq tA^{(n)}\subseteq ...$ of $t(X\otimes Y)$ by
  \begin{displaymath}
    tA^{(n)}\defeq t(X\pretens Y)\cup\left\{\pair<x;y>\in t(X\otimes Y) \left|
        \begin{array}{l}
          \exists x'\in X, y'\in Y, \alpha,\beta\in\Delta \text{ such that} \\
          x=x'\cdot\alpha\wedge y=y'\cdot\beta \wedge \Phi\pair<\alpha;\beta>\leq n
        \end{array}
    \right.\right\}
  \end{displaymath}
  and let $A^{(n)}$ be the stratified set $\pair<X\times Y;tA^{(n)}>$.
  The following observations are of note:
  
  \parhead[$tA^{(0)} = t(X\pretens Y)$] {\em Proof}.
    Suppose that the $r$-simplex $\pair<x;y>\in t(X\otimes Y)$
    possesses decompositions $x=x'\cdot\alpha$ and $y=y'\cdot\beta$
    such that $\Phi\pair<\alpha;\beta>\leq 0$. Consulting the
    definition of $\Phi$, we see that for this latter condition to
    hold each $\alpha(i)=0$ and that each $\beta(i)=m$ where
    $\cod(\beta)=[m]$.  Furthermore $r$ must be greater than $0$,
    since $t(X\otimes Y)$ does not contain any 0-simplices, therefore
    $0,1\in[r]$, $\alpha(0)=\alpha(1)$ and $\beta(0)=\beta(1)$ from
    which we may infer that
    $\pair<x;y>=\pair<x'\cdot\alpha;y'\cdot\beta>$ is degenerate and
    so it is an element of $t(X\pretens Y)$. It follows, by the
    definition of $tA^{(0)}$, that $tA^{(0)}\subseteq t(X\pretens Y)$ as
    required.
  
  \parhead[$\bigcup_{n=0}^{\infty} tA^{(n)} = t(X\otimes Y)$] {\em Proof}.
  Given an $r$-simplex $\pair<x;y>\in t(X\otimes Y)$ we know that
  $x=x\cdot\id_{[r]}$, $y=y\cdot\id_{[r]}$ and
  $\Phi\pair<\id_{[r]};\id_{[r]}> = r^2$; it follows, by the defining
  property of $tA^{({r^2})}$, that $\pair<x;y>\in tA^{({r^2})}$. In other
  words, every element of $t(X\otimes Y)$ is a member of those $tA^{(n)}$
  for which $n$ is ``sufficiently large'' and therefore $t(X\otimes
  Y)\subseteq\bigcup_{n=0}^{\infty} tA^{(n)}$ as required.
  
  Now we are in position to verify the postulates of
  observation~\ref{t-ext.expl.obs}, so suppose that $\pair<x;y>$ is an
  $r$-simplex in $tA^{({n+1})}\smallsetminus tA^{(n)}$. We know, by
  the definition of $tA^{({n+1})}$, that we may find decompositions
  $x=x'\cdot\alpha$ and $y=y'\cdot\beta$ such that
  $\Phi\pair<\alpha;\beta>\leq n+1$, furthermore $tA^{(n+1)}\subseteq
  t(X\otimes Y)$ and $t(X\pretens Y)\subseteq tA^{(n)}$ so
  $\pair<x;y>$ is thin in $X\otimes Y$ and it is {\bf not} thin in
  $X\pretens Y$, and it follows that we may apply
  lemma~\ref{(pre)tens.diff} to find a $0<k\leq r$ with
  $\alpha(k-1)<\alpha(k)$ and $\beta(k-1)<\beta(k)$.  So,
  consider the $(r+1)$-simplex $\pair<x\cdot\degen^r_{k-1};
  y\cdot\degen^r_k>$ and observe that:
  \begin{itemize}
  \item $k$ witnesses that $\pair<x\cdot\degen^r_{k-1};
    y\cdot\degen^r_k>$ is a mediator simplex, in particular it is
    $k$-admissible in $X\pretens Y$, by
    lemma~\ref{mediator.admissible}, and $t(X\pretens Y)\subseteq
    tA^{(n)}$ so it is also $k$-admissible in $A^{(n)}$.
  \item $\pair<x\cdot\degen^r_{k-1};
    y\cdot\degen^r_k>\cdot\face^{r+1}_k =
    \pair<x\cdot(\degen^r_{k-1}\circ\face^{r+1}_k);
    y\cdot(\degen^r_k\circ\face^{r+1}_k)> = \pair<x;y>$. However, we
    selected this so that it was thin in $X\otimes Y$ and we may apply
    lemma~\ref{mediator.admissible}(\ref{mediator.admissible.2}) to
    show that the faces $\pair<x\cdot\degen^r_{k-1};
    y\cdot\degen^r_k>\cdot\face^{r+1}_{k+1}$ and $\pair<x\cdot\degen^r_{k-1};
    y\cdot\degen^r_k>\cdot\face^{r+1}_{k-1}$ are also thin in there.
  \item $\pair<x\cdot\degen^r_{k-1};
    y\cdot\degen^r_k>\cdot\face^{r+1}_{k-1} =
    \pair<x\cdot(\degen^r_{k-1}\circ\face^{r+1}_{k-1});
    y\cdot(\degen^r_k\circ\face^{r+1}_{k-1})> = \pair<x;
    y\cdot(\degen^r_k\circ\face^{r+1}_{k-1})>$ and applying our
    decompositions of $x$ and $y$ this is equal to
    $\pair<x'\cdot\alpha; y'\cdot(\beta\circ\degen^r_k\circ
    \face^{r+1}_{k-1})>$. However, notice that
    $\beta\circ\degen^r_k\circ\face^{r+1}_{k-1}$ and $\beta$ only
    differ at $k-1$ where $\beta\circ\degen^r_k\circ\face^{r+1}_{k-1}(
    k-1) = \beta(k)$, so it follows that
    \begin{displaymath}
      \begin{array}{rcll}
        \Phi\pair<\alpha;
        \beta\circ\degen^r_k\circ\face^{r+1}_{k-1}> & = &
        \Phi\pair<\alpha; \beta> + \beta(k-1) - \beta(k) &
        \text{by the definition of $\Phi$,} \\
        & < & \Phi\pair<\alpha;\beta> - 1 & \text{since
          $\beta(k-1)<\beta(k)$,} \\
        & \leq & n & 
      \end{array}
    \end{displaymath}
    and so $\pair<x\cdot\degen^r_{k-1}; y\cdot\degen^r_k>\cdot
    \face^{r+1}_{k-1}$, which we've already shown to be thin in
    $X\otimes Y$, is an element of $tA^{(n)}$ (by definition).
  \item dually $\pair<x\cdot\degen^r_{k-1};
    y\cdot\degen^r_k>\cdot\face^{r+1}_{k+1} =
    \pair<x'\cdot(\alpha\circ\degen^r_{k-1}\circ\face^{r+1}_{k+1});
    y'\cdot\beta>$, but $\alpha\circ\degen^r_{k-1}\circ
    \face^{r+1}_{k+1}$ and $\alpha$ only differ at $k$ where
    $\alpha\circ\degen^r_{k-1}\circ\face^{r+1}_{k+1}(k) =
    \alpha(k-1)$. It follows that
    \begin{displaymath}
      \begin{array}{rcll}
        \Phi\pair<\alpha\circ\degen^r_{k-1}\circ\face^{r+1}_{k+1};
        \beta> & = &
        \Phi\pair<\alpha; \beta> - \alpha(k) + \alpha(k-1) &
        \text{by the definition of $\Phi$,} \\
        & \leq & \Phi\pair<\alpha;\beta> - 1 & \text{since
          $\alpha(k-1)<\alpha(k)$,} \\
        & \leq & n & 
      \end{array}
    \end{displaymath}
    and so $\pair<x\cdot\degen^r_{k-1}; y\cdot\degen^r_k>\cdot
    \face^{r+1}_{k+1}$, which we've already shown to be thin in
    $X\otimes Y$, is also an element of $tA^{(n)}$ (by definition).
  \end{itemize} 
  In summary, if $\pair<x;y>$ is an $r$-simplex of
  $tA^{({n+1})}\smallsetminus tA^{(n)}$ then the $r+1$-simplex
  $\pair<x\cdot\degen^r_{k-1}; y\cdot\degen^r_k>$ which we constructed
  witnesses the extension of thinness to $\pair<x;y>$, as required by
  observation~\ref{t-ext.expl.obs}.
  
  This completes the proof that the sequence $\{A^{(n)}\}_{i=0}^{\infty}$
  satisfies all of the conditions of observation~\ref{t-ext.expl.obs},
  from which it follows that the stratified map $\overinc\subseteq_e:
  X\pretens Y-> X\otimes Y.$ is a t-extension as required.
\end{proof}

\begin{lemma}\label{a.useful.t-ext.res}
  If $X'\subseteq_e X$ and $Y'\subseteq_e Y$ are entire stratified
  subsets then we have a commutative square
  \begin{equation}\label{a.useful.sq}
    \xymatrix@R=3em@C=6em{
      {(X\pretens Y')\cup(X'\pretens Y)}
      \ar@{u(->}[r]^-{\subseteq_e}
      \ar@{u(->}[d]_{\subseteq_e} &
      {X\pretens Y}
      \ar@{u(->}[d]^{\subseteq_e} \\
      {(X\otimes Y')\cup(X'\otimes Y)}
      \ar@{u(->}[r]_-{\subseteq_e} &
      {X\otimes Y} }
  \end{equation}
  of entire inclusions in which the upper horizontal is a stratified
  isomorphism and the lower horizontal is a t-extension.
\end{lemma}

\begin{proof}
  By lemma~\ref{pretens.density} the vertical inclusions in this
  square are t-extensions so, by
  observation~\ref{l-inv}(\ref{l-inv.comp}), it follows that its lower
  horizontal is a t-extension if and only if its upper horizontal is a
  t-extension.  So in order to establish both of the stated results it
  is enough to demonstrate that this latter map is an isomorphism,
  however since it is an entire inclusion all that remains is to
  demonstrate that it is also regular.
  
  To show that this is the case, it is clear that all we need
  do is demonstrate that any crushed cylinder $\pair<x;y>$ in
  $X\pretens Y$ is a crushed cylinder in one of $X\pretens Y'$ or
  $X'\pretens Y$. So suppose that $\pair<x;y>$ is an $r$-simplex and
  that it is witnessed as a crushed cylinder by the partition
  $p,q$ of $r$ then, by definition, we have two cases 
  \begin{itemize}
  \item $x\cdot\partinj^{p,q}_1$ is thin in $X$, so
    our partition also witnesses $\pair<x;y>$ as a crushed cylinder in
    $X\pretens Y'$, or
  \item $y\cdot\partinj^{p,q}_2$ is thin in $Y$, so
    our partition also witnesses $\pair<x;y>$ as a crushed cylinder in
    $X'\pretens Y$
  \end{itemize}
  as required.
\end{proof}

\begin{obs}\label{tensor.precomp.test}
  Of course, lemma~\ref{pretens.density} implies that if $A$ is a
  pre-complicial set then we may show that a simplicial map $\arrow
  f:X\times Y->A.$ extends to a stratified map $\arrow f:X\otimes Y->A.$
  simply by verifying that $f$ maps all {\bf mediator simplices} and
  {\bf crushed cylinders} in $X\otimes Y$ to thin simplices in $A$.
\end{obs}

\subsection{Some Other Preservation Properties}

\begin{lemma}\label{colim.pres} For each stratified set
  $Y$ the functor 
  \begin{displaymath}
    \let\labelstyle=\textstyle
    \xymatrix@C=10em{{\Strat}\ar[r]^{\mathord{-}\pretens Y} & {\Strat}}
  \end{displaymath}
  preserves all small colimits.
\end{lemma}

\begin{proof} Given a diagram $\arrow D:\mathbb{C}->\Strat.$ our
  aim is to show that the canonical stratified map
  \begin{equation}\label{comparison}
    \xymatrix@C=5em@R=1ex{
      {\colim(D(\mathord{-})\pretens Y)}
      \ar[r]^{\textstyle c_{D,Y}} &
      {\colim(D)\pretens Y}}
  \end{equation}
  is a stratified isomorphism. In fact, we know that the functor
  $\arrow \forget:\Strat->\Simp.$ preserves colimits and that, for a
  fixed stratified set $Y$, we have a commutative square of
  functors:
  \begin{displaymath}
    \let\labelstyle=\textstyle
    \xymatrix@R=1em@C=2em{
      {\Strat}\ar[rr]^{\mathord{-}\pretens Y}\ar[dd]_{\forget} & &
      {\Strat}\ar[dd]^{\forget} \\
      & {\cong} & \\
      {\Simp}\ar[rr]_{-\times Y} & &
      {\Simp} }
  \end{displaymath}
  So the underlying simplicial map of $c_{D,Y}$ is
  simply the corresponding comparison map
  \begin{displaymath}
    \xymatrix@C=5em@R=1ex{
      {\colim(\forget\circ D(-)\times Y)} \ar[r] &
      {\colim(\forget\circ D)\times Y}}
  \end{displaymath}
  in $\Simp$, which is an isomorphism since $\Simp$ is a cartesian
  closed category (from which it follows that $-\times Y$ preserve
  colimits). We can infer, using observation~\ref{strat.isos}, that
  all we need to do is demonstrate that $c_{D,Y}$ is regular. 
  
  Under the concrete presentation of the colimits in $\Strat$ given in
  observations~\ref{simplicial.(co)limits}
  and~\ref{stratified.(co)limits}, it is clear that this map bears an
  explicit description under which it carries the equivalence class
  $[i,\pair<z;y>]$ in $\colim(\dgm(-)\pretens Y)$ to the pair
  $\pair<[i,z];y>$ in $\colim(D)\pretens Y$.  So suppose that the
  non-degenerate $r$-simplex $\pair<[i,z];y>$ is thin in
  $Y\pretens\colim(D)$, our aim will be to find another
  $j\in\mathbb{C}$ and $w\in D(j)$ such that $\pair<w;y>$ is thin
  in $D(j)\pretens Y$ and $[j,w]=[i,z]$.  By
  definition~\ref{pre.tens} there are two possibilities:
  \begin{itemize}
  \item {\bf\boldmath $\pair<[i,z];y>$ is a mediator simplex as
      witnessed by some $k$}, in other words
    $[i,z]=[i,z]\cdot(\face^r_{k-1}\circ \degen^{r-1}_{k-1})$ and
    $y=y\cdot(\face^r_k\circ \degen^{r-1}_k)$. Notice that
    $[i,z]\cdot(\face^r_{k-1}\circ
    \degen^{r-1}_{k-1})=[i,z\cdot(\face^r_{k-1}\circ
    \degen^{r-1}_{k-1})]$, so let $w=z\cdot(\face^r_{k-1}\circ
    \degen^{r-1}_{k-1})$, $j=i$ and observe that the integer $k$ also
    witnesses the $r$-simplex $\pair<w;y>$ in
    $D(j)\pretens Y$ as a mediator simplex.
  \item {\bf\boldmath $\pair<[i,z];y>$ is a crushed cylinder as
      witnessed by some partition $p,q$ of $r$}, in other words
    $[i,z]=[i,z]\cdot(\partinj^{p,q}_1\circ\partproj^{p,q}_1)$,
    $y=y\cdot(\partinj^{p,q}_2\circ\partproj^{p,q}_2)$ and either
    of $[i,z]\cdot\partinj^{p,q}_1$ or $y\cdot\partinj^{p,q}_2$ is
    thin. If $[i,z]\cdot\partinj^{p,q}_1$ is thin then we can find
    an $j\in\mathbb{C}$ and a $v\in D(j)$ such that $v$ is thin
    and $[j,v]=[i,z]\cdot\partinj^{p,q}_1$, otherwise, simply let
    $j=i$ and $v=z\cdot\partinj^{p,q}_1$. Setting
    $w=v\cdot\partproj^{p,q}_1$, it is easy to show that the
    $r$-simplex $\pair<w;y>$ is a crushed cylinder in
    $D(j)\pretens Y$ and that $[i,z]=[j,w]$.
  \end{itemize}
  In either case $\pair<w;y>$ is thin in $D(j)\pretens Y$, which
  implies that $[j,\pair<w;y>]$ is thin in
  $\colim(D(\mathord{-})\pretens Y)$. Furthermore
  $c_{D,Y}[j,\pair<w;y>] = \pair<[j,w];y> = \pair<[i,z];y>$, so
  the unique simplex mapping to $\pair<[i,z];y>$ under $c_{D,Y}$ is
  thin in $\colim(D(\mathord{-})\pretens Y)$; in other words
  $c_{D,Y}$ is regular as required.
\end{proof}

\begin{lemma}[tensor products of primitive t-extensions]
  \label{prim.t.tens} 
  If $Y$ is a stratified set then in the commutative square
  \begin{equation}\label{prim.t.rtens}
    \let\labelstyle=\textstyle
    \xymatrix@R=2.5em@C=8em{
      {\Delta^a_k[n]'\pretens Y}
      \ar@{u(->}[r]^{\subseteq_e}
      \ar@{u(->}[d]_{\subseteq_e} &
      {\Delta^a_k[n]''\pretens Y}
      \ar@{u(->}[d]^{\subseteq_e} \\
      {\Delta^a_k[n]'\otimes Y}
      \ar@{u(->}[r]_{\subseteq_e} &
      {\Delta^a_k[n]''\otimes Y}
    }
  \end{equation}
  the horizontal inclusions, which may be constructed by tensoring and
  pre-tensoring the primitive t-extension $\inc
  \subseteq_e:\Delta_k^a[n]'->\Delta_k^a[n]''.$ by $Y$, are both
  t-extensions.
\end{lemma}

\begin{proof}
  Firstly observe that the vertical inclusions in the square are both
  t-extensions, by lemma~\ref{pretens.density}, and therefore applying
  observation~\ref{l-inv}(\ref{l-inv.comp}) we see that its lower
  horizontal is a t-extension if and only if its upper horizontal is a
  t-extension. We may factor the lower of these two into a pair of
  entire inclusions
  \begin{equation}\label{prim.t.tens.2}
    \xymatrix@R=1ex@C=6em{
      {\Delta^a_k[n]'\otimes Y}
      \ar@{^{(}->}[r]^-{\subseteq_e} &
      {(\Delta^a_k[n]' \otimes Y)\cup (\Delta^a_k[n]''\pretens Y)}
      \ar@{^{(}->}[r]^-{\subseteq_e} &
      {\Delta^a_k[n]''\otimes Y}
    }
  \end{equation}
  and show that their composite is a t-extension by showing that each
  of its factors is and applying the composition result of
  observation~\ref{l-inv}(\ref{l-inv.comp}).
  
  To show that the right hand inclusion here is a t-extension, we
  factor the t-extension of lemma~\ref{pretens.density} as
  \begin{equation*}
    \xymatrix@R=1ex@C=6em{
      {\Delta^a_k[n]''\pretens Y}
      \ar@{^{(}->}[r]^-{\subseteq_e} &
      {(\Delta^a_k[n]' \otimes Y)\cup (\Delta^a_k[n]''\pretens Y)}
      \ar@{^{(}->}[r]^-{\subseteq_e} &
      {\Delta^a_k[n]''\otimes Y}
    }
  \end{equation*}
  and appeal to the right cancellation result given in
  observation~\ref{l-inv}(\ref{l-inv.comp}), which we may apply in
  this case since all entire subset inclusions are epimorphisms in
  $\Strat$.
  
  Finally we show that the left hand inclusion in
  display~\eqref{prim.t.tens.2} is also a t-extension by demonstrating
  that we may apply observation~\ref{t-ext.expl.obs} to the ``single
  step'' sequence $\Delta^a_k[n]'\otimes
  Y\subseteq_e(\Delta^a_k[n]'\otimes Y)\cup (\Delta^a_k[n]''\pretens
  Y)$. By the definition of the pre-tensor product $\pretens$ we know
  that if an $r$-simplex $\pair<\alpha;y>$ is thin in
  $\Delta^a_k[n]''\pretens Y$ and not thin in $\Delta^a_k[n]'\otimes
  Y$ then it must be a cylinder witnessed by a partition $s,t$ for
  which $y'\defeq y\cdot\partinj^{s,t}_2$ is not thin in $Y$ and
  $\alpha\circ\partinj^{s,t}_1$ is thin in $\Delta^a_k[n]''$ but not
  in $\Delta^a_k[n]'$. This latter condition implies that
  $\alpha\circ\partinj^{s,t}_1=\face^n_k$ and so we have $s=n-1$ and
  $t=r-n+1$ (in particular $r\geq n-1$), Consequently consider the
  $(r+1)$-simplex $\pair<\partproj^{n,t}_1; y'\cdot\partproj^{n,t}_2>$
  and observe that this witnesses the extension of thinness from
  $\Delta^a_k[n]'\otimes Y$ to the simplex $\pair<\alpha;y>$ as
  verified in the following points:
  \begin{enumerate}[(i)]
  \item If $i\in[n]$ then a simple calculation with simplicial
    operators demonstrates that:
    \begin{equation}\label{face.calc.t-ext.tens}
      \begin{split}
        \pair<\partproj^{n,t}_1;y'\cdot\partproj^{n,t}_2>\cdot
        \face^{r+1}_i & {} =\pair<\partproj^{n,t}_1\circ\face^{r+1}_i;
        y'\cdot(\partproj^{n,t}_2\circ\face^{r+1}_i)> \\
        & {} =
        \begin{cases}
          \pair<\face^n_i\circ\partproj^{n-1,t}_1;y'\cdot\partproj^{n-1,t}_2>
          & \text{when $i<n$,}\\
          \pair<\partproj^{n,t-1}_1;y'\cdot\partproj^{n-1,t}_2>
          & \text{when $i=n$.}
        \end{cases}
      \end{split}
    \end{equation}
    In particular, the assumption that $s,t$ witnesses our simplex
    $\pair<\alpha;y>$ as a cylinder implies that $\alpha=
    \face^n_k\circ\partproj^{s,t}_1$ and $y=y'\cdot\partproj^{s,t}_2$,
    so it follows, from the calculation above, that
    $\pair<\partproj^{n,t}_1; y'\cdot\partproj^{n,t}_2>\cdot
    \face^{r+1}_k=\pair<\alpha;y>$ (since $k<n$).
  \item\label{prim.t.tens.pr2} We know that $k-1<n$ so applying
    equation~\eqref{face.calc.t-ext.tens} we see that the face
    $\pair<\partproj^{n,t}_1; y'\cdot\partproj^{n,t}_2>\cdot
    \face^{r+1}_{k-1}$ is a cylinder which is thin in
    $\Delta^a_k[n]'\otimes Y$ since $\face^n_{k-1}$ is thin in
    $\Delta^a_k[n]'$.
  \item\label{prim.t.tens.pr3} If $k+1<n$ then we can apply the same
    argument to show that the face $\pair<\partproj^{n,t}_1;
    y'\cdot\partproj^{n,t}_2> \cdot \face^{r+1}_{k+1}$ is thin in
    $\Delta^a_k[n]'\otimes Y$, leaving us the special case $k=n-1$ in
    which equation~\eqref{face.calc.t-ext.tens} identifies this face
    as the simplex
    $\pair<\partproj^{n,t-1}_1;y'\cdot\partproj^{n-1,t}_2>$.  To show
    that this face is thin in $\Delta^a_k[n]'\otimes Y$ simply observe
    that
    \begin{itemize}
    \item $\partproj^{n,t-1}_1\circ\partinj^{n-1,t}=\face^n_n$ and
      $\partproj^{n,t-1}_1\circ\partinj^{n,t-1}=\id_{[n]}$ are both
      thin in $\Delta^a_k[n]'$ so the partitions $n-1,t$ and $n,t-1$
      both cleave $\pair<\partproj^{n,t-1}_1;y'\cdot
      \partproj^{n-1,t}_2>$, and
    \item $\partproj^{n,t-1}_1$ is degenerate at $n$ and
      $y'\cdot\partproj^{n-1,t}_2$ is degenerate at $n-2$ so we may
      apply observation~\ref{deg.ord} to show that all other
      partitions also cleave that simplex, thereby completing the
      demonstration that it is thin as stated.
    \end{itemize}
  \item\label{prim.t.tens.admprf} If $\arrow\mu:[m]->[r+1].$ is a
    $k$-divided face operator and $l\in[m]$ is the unique integer for
    which $\mu(l)=k$ (and thus $\mu(l-1)=k-1$ and $\mu(l+1)=k+1$) and
    consider the face
    $\pair<\partproj^{n,t}_1;y'\cdot\partproj^{n,t}_2>\cdot\mu=
    \pair<\partproj^{n,t}_1\circ\mu;y'\cdot(\partproj^{n,t}_2\circ\mu)>$.
    We show that this is thin in $\Delta^a_k[n]'\otimes Y$ by
    selecting an arbitrary partition $p,q$ of $m$ and considering
    cases:
    \begin{itemize}
    \item {\boldmath $p\leq l$} in which case we have
      $\mu\circ\partinj^{p,q}_2(l-p)=\mu(l)=k$ and
      $\mu\circ\partinj^{p,q}(l-p+1)=\mu(l+1)=k+1$, however $k<k+1\leq
      n$ so $\partproj^{n,t}_2(k)=\partproj^{n,t}(k+1) =0$ and it
      follows that $(y'\cdot(\partproj^{n,t}_2\circ\mu))\cdot
      \partinj^{p,q}_2=y'\cdot(\partproj^{n,t}_2\circ\mu\circ
      \partinj^{p,q}_2)$ is degenerate at $l-p$ and is thus thin in
      $Y$, or
    \item {\boldmath $p\geq l+1$} which implies that $l-1, l$ and
      $l+1$ are all elements in the domain of
      $\partproj^{n,t}_1\circ\mu\circ\partinj^{p,q}_1$ which maps them
      to $k-1, k$ and $k+1$ respectively, so it follows that this
      operator is $k$-divided and thus thin in $\Delta^a_k[n]$.
    \end{itemize}
    In either case we've shown that $p,q$ cleaves the simplex
    $\pair<\partproj^{n,t}_1; y'\cdot\partproj^{n,t}_2>\cdot\mu$ so
    quantifying over all such partitions we see that this face is thin
    in $\Delta^a_k[n]\otimes Y\subseteq_e\Delta^a_k[n]'\otimes Y$.\qedhere
  \end{enumerate}
\end{proof}

\subsection{A Monoidal Biclosed Structure on Pre-Complicial sets}

\begin{cor}[of lemma~\ref{colim.pres}]\label{closure} The functor
  $\arrow\mathord{-}\pretens Y:\Strat->\Strat.$ admits a right
  adjoint, which we'll denote by $\arrow \mathord{*}\Leftarrow Y:
  \Strat->\Strat.$.  Dually, for each stratified set $X$ the functor
  $\arrow X\pretens{-}: \Strat->\Strat.$ has a right adjoint $\arrow
  X\Rightarrow{*}:\Strat->\Strat.$ given by $(X\Rightarrow{*})\defeq
  (({*})^\circ\Leftarrow X^\circ)^\circ$. 
\end{cor}

\begin{proof} We know, from observation~\ref{deltat.dense}, that the
  small full subcategory $\tDelta$ is dense in $\Strat$. Applying
  theorem 5.33 of Kelly~\cite{Kelly:1982:ECT}, it follows that any
  (small) colimit preserving functor with codomain $\Strat$ has a
  right adjoint. In particular, in view of lemma~\ref{colim.pres} this
  result applies to ${-}\pretens Y$ giving us stated right adjoint
  ${*}\Leftarrow Y$.
  
  The dual result follows directly from basic properties of the
  strictly involutive dual $({-})^\circ$ and the canonical
  isomorphisms $(X\pretens Y)^\circ\cong Y^\circ\pretens X^\circ$.
\end{proof}

\begin{obs}[presenting $A\Leftarrow Y$ and $X\Rightarrow A$
  explicitly]\label{expl.power} The result of the last corollary fully
  characterises $A\Leftarrow Y$, however it will sometimes be useful
  to assume an explicit representation of this stratified set.  Using
  Yoneda's lemma, we know that an $n$-simplex in $A\Leftarrow Y$
  corresponds to a stratified map $\overarr:\Delta[n]->A\Leftarrow Y.$
  which in turn corresponds to a stratified map
  \begin{equation}\label{right.clos.elem}
    \let\labelstyle=\textstyle
    \xymatrix@R=1ex@C=14em{
      {\Delta[n]\pretens Y}\ar[r]^{f} &
      {A}}
  \end{equation}
  via the adjunction between the functors $\mathord{-}\pretens Y$ and
  $\mathord{*}\Leftarrow Y$ of the last corollary. It follows that we
  may identify $n$-simplices in $A\Leftarrow Y$ and stratified maps
  $f$ of the form shown in display~\eqref{right.clos.elem} so, using
  the usual naturality properties of adjunctions, we may derive the
  following explicit presentation of $A\Leftarrow Y$:
  \begin{itemize}
  \item {\bf\boldmath $n$-simplices} are stratified maps of the form
    given in figure~(\ref{right.clos.elem}),
  \item {\bf simplicial action} of an operator
    $\arrow\alpha:[m]->[n].$ on an $n$-simplex $f$ given by
    $f\cdot\alpha\defeq f\circ(\Delta(\alpha)\pretens Y)$,
  \item {\bf thin simplices} those $n$-simplices $\arrow
    f:\Delta[n]\pretens Y->A.$ which extend to a stratified map
    $\arrow f:\Delta[n]_t\pretens Y->A.$.
  \end{itemize}
  Dually, we may represent the $n$-simplices of $X\Rightarrow
  A$ as stratified maps $\arrow g:X\pretens\Delta[n]_?->A.$ under the
  simplicial action given by $g\cdot\alpha\defeq
  g\circ(X\pretens\Delta(\alpha))$.
\end{obs}

\begin{obs}
  Of course, if $A$ is a pre-complicial set then we know, by
  lemma~\ref{pretens.density}, that stratified maps of the form given
  in figure~(\ref{right.clos.elem}) extend to maps of the form:
  \begin{equation}\label{right.clos.elem.2}
    \let\labelstyle=\textstyle
    \xymatrix@R=1ex@C=14em{
      {\Delta[n]\otimes Y}\ar[r]^{f} &
      {A}}
  \end{equation}
  It follows that, in this case, it is equally valid to identify the
  elements of $A\Leftarrow Y$ with stratified maps of
  the form given in figure~(\ref{right.clos.elem.2}). The remainder of
  the explicit description of $A\Leftarrow Y$ given in
  the last observation carries over to elements represented in this
  form simply by substituting $\otimes$ for $\pretens$ wherever it
  occurs. 
\end{obs}

\begin{lemma}\label{exp.inc}
  Suppose that $Y$ and $A$ are stratified sets and that $Y'\subseteq_e
  Y$ is an entire subset then we may canonically identify $A\Leftarrow
  Y$ with a regular stratified subset of $A\Leftarrow Y'$.
  Furthermore, a stratified map $\arrow g:Y'->X.$ extends to a
  stratified map with domain $Y$ iff for each $A$ the corresponding
  map $\arrow A\Leftarrow g:A\Leftarrow X->A\Leftarrow Y'.$ lifts
  to a stratified map with codomain $A\Leftarrow Y$.
\end{lemma}

\begin{proof}
  Applying the contravariant functor $A\Leftarrow{*}$ to the entire
  inclusion $\overinc\subseteq_e:Y'->Y.$ gives rise to a stratified
  map
  \begin{equation}\label{power.inc}
    \let\labelstyle=\textstyle
    \xymatrix@R=1ex@C=14em{
      {A\Leftarrow Y}
      \ar[r]^{A\Leftarrow(\subseteq_e)} &
      {A\Leftarrow Y'}}
  \end{equation}
  which, in terms of the explicit presentation given above, acts on an
  $n$-simplex $f$ of $A\Leftarrow Y'$ to restrict its domain from
  $\Delta[n]\pretens Y$ to its entire stratified subset
  $\Delta[n]\pretens Y'$. This action is clearly injective, since
  these two stratified sets share the same underlying simplicial set.
  
  This result certainly allows us to identify $A\Leftarrow Y$ with a
  stratified subset of $A\Leftarrow Y'$, but we also need to show that
  the inclusion in~(\ref{power.inc}) is regular. So suppose that
  $\arrow f:\Delta[n]\pretens Y->A.$ is an $n$-simplex in $A\Leftarrow
  Y$ which maps to a thin simplex in $A\Leftarrow Y'$. We know that
  this simply means that the restriction of $f$ to $\Delta[n]\pretens
  Y'$ extends to a stratified map with domain $\Delta[n]_t\pretens Y'$.
  Clearly it follows that $f$ extends to $(\Delta[n]\pretens
  Y)\cup(\Delta[n]_t \pretens Y')$ and, by
  lemma~\ref{a.useful.t-ext.res}, this latter stratified set is equal
  to $\Delta[n]_t\pretens Y$ thus $f$ is thin in $A\Leftarrow Y$ as
  required.

  The remainder of the observation is a matter of simple abstract
  nonsense and is left as an exercise.
\end{proof}

\begin{thm}\label{otens.t-ext.pres}\label{biclosed.precomp}
  For each stratified set $Y$ the right tensor functor
  $\arrow\mathord{-}\otimes Y: \Strat-> \Strat.$ has partial right
  adjoint $\arrow\mathord{*}\Leftarrow Y:\Precomp->\Precomp.$ (cf.\ 
  observation~\ref{partial.adj}).  Dually,
  the left tensor functor $\arrow X\otimes{-}:\Strat->\Strat.$
  associated with a stratified set $X$ has
  partial right adjoint $\arrow X\Rightarrow\mathord{*}:\Precomp->
  \Precomp.$.  It follows, by
  observation~\ref{partial.adj}(\ref{partial.adj.a}), that if $\arrow
  f:Z-> Z'.$ is a t-extension then so are the maps:
  \begin{displaymath}
    \let\labelstyle=\textstyle
    \xymatrix@R=1ex@C=10em{
      Z\otimes Y\ar[r]^{\textstyle f\otimes Y} &
      Z'\otimes Y \\
      X\otimes Z\ar[r]^{\textstyle X\otimes f} &
      X\otimes Z'}
  \end{displaymath}
  Furthermore, by Day's reflection theorem (theorem~\ref{day.refl}), we
  may reflect the monoidal structure on $\Strat$ to a monoidal
  biclosed structure on $\Precomp$ with tensor product $A\otimes_p
  B\defeq \refl_p(A\otimes B)$ (where $\refl_p$ is reflector associated with
  $\Precomp$) and left and right closures $A\Rightarrow{*}$ and
  ${*}\Leftarrow B$ respectively.
\end{thm}

\begin{proof}
  Fix a stratified set $Y$. We know (by definition) that the set of
  primitive t-extensions is adequate to detect pre-complicial sets and
  that the functor $\arrow {-}\otimes Y:\Strat->\Strat.$ maps
  primitive t-extensions to t-extensions (by lemma~\ref{prim.t.tens}).
  We may infer that it is enough to show that $\arrow{*}\Leftarrow
  Y:\Precomp->\Strat.$ is a partial right adjoint to ${-}\otimes Y$,
  since if that is the case then we may apply
  observation~\ref{partial.adj}(\ref{partial.adj.b}) to demonstrate
  that ${*}\Leftarrow A$ maps pre-complicial sets to pre-complicial
  sets. However, this latter result is a direct consequence of
  corollary~\ref{closure} via lemma~\ref{pretens.density}, as
  demonstrated by the following sequence of natural bijections
  \begin{equation}\label{otens.clos}
    \let\labelstyle=\textstyle
    \xymatrix@C=1em@R=0ex{
      & {X}\ar[r] & {A\Leftarrow Y} & \\
      \ar@{-}[rrr] & & & 
      \save[]+<7em,0ex>*{\txt{corollary~\ref{closure},}}\restore \\
      & {X\pretens Y}\ar[r] & {A} & \\
      \ar@{-}[rrr] & & & 
      \save[]+<7em,0ex>*{\txt{the t-extension
        of lemma~\ref{pretens.density} is $\perp$\\
        to the pre-complicial set $A$,}}\restore \\
      & {X\otimes Y}\ar[r] & {A} & 
      }
  \end{equation}
  The remainder of this theorem, that is the corresponding result for
  $\Rightarrow$, is an immediate dual of the above result. 
\end{proof}

\begin{notation}
  When making calculations with the tensor $\otimes$ and the partial
  closures ${*}\Leftarrow Y$ and $X\Rightarrow{*}$ we will tend to
  simplify matters by following traditional conventions with regard to
  associativities and the identity object. In particular:
  \begin{itemize}
  \item We often write unbracketed $n$-fold tensors $X_1\otimes
    X_2\otimes...\otimes X_n$ and assume that such expressions are
    interpreted via the consistent application of a rule, such as
    right associativity, under which they might be bracketed into
    binary tensors.
  \item Wherever possible, we will reduce tensor and closure
    expressions which involve the identity object $\Delta[0]$ to
    eliminate them wherever possible.  In other words, we will tend to
    identify expressions such as $X\otimes\Delta[0]$ with $X$,
    $\Delta[0]\otimes Y$ with $Y$,  $A\Leftarrow\Delta[0]$ with $A$
    and $\Delta[0]\Rightarrow B$ with $B$.
 \end{itemize}
 These rules tacitly assume that whenever we manipulate expressions
 involving composites of arrows between objects which have been
 specified in this way, we will be prepared to silently introduce
 canonical associativity and identity isomorphisms wherever necessary
 to make sense of those composites. That this is possible to do in
 general is a consequence of the coherence theorem for monoidal
 categories and we refer the reader to Mac~Lane's
 book~\cite{Maclane:1971:CWM} for more on that topic.
\end{notation}

\subsection{Superstructures of Pre-Complicial Sets}

\begin{lemma}\label{precomp.superst} 
  If $A$ is a pre-complicial set then so is its $n$-dimensional
  superstructure $\Sup_n(A)$. Equivalently, if the stratified map
  $\arrow f:X ->Y.$ is a t-extension then so is the associated
  stratified map $\arrow \Th_n(f): \Th_n(X)->\Th_n(Y).$.
\end{lemma}

\begin{proof} Since $\Th_n\adjoint\Sup_n(\cdot)$ we know, by
  observation~\ref{partial.adj}, that the formulations given in the
  statement of this lemma are equivalent and that they both follow if
  we can demonstrate that each map obtained by applying $\Th_n$ to
  a primitive t-extension is again a t-extension. So consider the
  primitive t-extension $\overinc \subseteq_e:
  \Delta^a_k[m]'->\Delta^a_k[m]''.$ and observe that we have two
  cases:
  \begin{itemize}
  \item {\boldmath $n\geq m-1$} In this case, all of the simplices of
    $\Delta^a_k[m]'$ and $\Delta^a_k[m]''$ of dimension $>n$ are
    already thin, therefore $\Th_n(\Delta^a_k[m]')=\Delta^a_k[m]'$,
    $\Th_n(\Delta^a_k[m]'') =\Delta^a_k[m]''$ and, consequently,
    $\Th_n$ maps $\overinc\subseteq_e: \Delta^a_k[m]'->\Delta^a_k[m]''.$ to
    itself.
  \item {\boldmath $n<m-1$} We know that $\Delta^a_k[m]'$ and
    $\Delta^a_k[m]''$ differ only in as much as the $(m-1)$-simplex
    $\face^m_k$ is not thin in the former and is thin in the latter.
    However, since $n<m-1$ we know that $\face^m_k$ is thin in both of
    $\Th_n(\Delta^a_k[m]')$ and $\Th_n(\Delta^a_k[m]'')$ and thus that
    these two stratified sets are identical. It follows that $\Th_n$
    maps $\arrow \subseteq_e:\Delta^a_k[m]'-> \Delta^a_k[m]''.$ to the
    identity on $\Th_n(\Delta^a_k[m]')=\Th_n(\Delta^a_k[m]'')$.
  \end{itemize}
  In either case, the resulting map is a t-extension as required.
\end{proof}

\begin{lemma}\label{precomp.triv} 
  If $A$ is an $n$-trivial pre-complicial set then $X\Rightarrow A$
  and $A\Leftarrow Y$ are also $n$-trivial for any stratified sets
  $X,Y\in\Strat$.
\end{lemma}

\begin{proof}
  The results for $X\Rightarrow A$ and $A\Leftarrow Y$ are dual so we
  only consider the former.  By definition, $\Strat_n$ is a reflective
  full subcategory of $\Strat$ with reflector $\Th_n$, so it follows
  that a stratified set is $n$-trivial if and only if it is orthogonal
  to each entire inclusion $\overinc \subseteq_e:Z->\Th_n(Z).$. Now,
  we know that $X\Rightarrow{*}$ is a partial right adjoint to
  $X\otimes{-}$ on the category of pre-complicial sets, so taking
  adjoint transposes we see that $X\Rightarrow A$ is orthogonal to
  $\overinc\subseteq_e:Z-> \Th_n(Z).$ if and only if $A$ is orthogonal
  to $\overinc\subseteq_e:X\otimes Z-> X\otimes\Th_n(Z).$. However, we
  have a sequence of entire inclusions
  \begin{equation*}
    \begin{aligned}
      \Th_n(X\otimes Z) & {} \subseteq_e \Th_n(X\otimes\Th_n(Z)) &&
      \text{applying $\Th_n$ to $X\otimes Z\subseteq_e
        X\otimes\Th_n(Z)$,}\\
      & {} \subseteq_e \Th_n(\Th_n(X\otimes Z)) &&
      \parbox{20em}{applying $\Th_n$ to $X\otimes\Th_n(Z)\subseteq_e
        \Th_n(X\otimes Z)$ of lemma~\ref{tensor.Th},} \\
      & {} \subseteq_e \Th_n(X\otimes Z) &&
      \text{$\Th_n$ is idempotent.}
    \end{aligned}
  \end{equation*}
  from which we see that $\Th_n(X\otimes Z)=\Th_n(X\otimes\Th_n(Z))$
  and thus that the inclusion of the last sentence is
  $\Th_n$-invertible. It follows that the $n$-trivial set $A$ is
  perpendicular to that inclusion for each $Z\in\Strat$ and thus that
  $X\Rightarrow A$ is $n$-trivial as required.
\end{proof}

\begin{obs}\label{precomp_n.mon.bic} 
  Let $\Precomp_n=\Precomp\cap\Strat_n$ denote the full subcategory of
  $n$-trivial pre-complicial sets in $\Precomp$. A stratified set is
  in $\Precomp_n$ iff it is orthogonal to the maps in the FP-regulus
  which is the union of the set of primitive t-extensions and the set
  $\{\overinc\subseteq_e:\Delta[m]-> \Delta[m]_t.\mid
  m\in\mathbb{N}\also m>n\}$. It follows that $\Precomp_n$ is an
  LFP-category and that it is a reflective full subcategory of
  $\Strat$.

  Now, lemma~\ref{precomp.triv} implies that $X\Rightarrow{*}$ (resp.\
  ${*}\Leftarrow Y$) restricts to a partial right adjoint to
  $X\otimes{-}$ (resp.\ ${-}\otimes Y$) on $\Precomp_n$, so it follows
  that we may apply Day's reflection theorem, just as in
  theorem~\ref{biclosed.precomp}, to reflect the monoidal structure of
  $\Strat$ to a monoidal biclosed structure on $\Precomp_n$. 
\end{obs}



\section{Complicial Sets}\label{comp.sec}

\subsection{Introducing Complicial Sets}

\begin{defn}[primitive f-extensions]\label{prim.f-ext} The class of
  {\em primitive f-extensions\/} is the union of the 
  following classes of stratified maps
  \begin{enumerate}[a)]
  \item the primitive t-extensions,
  \item\label{prim.f-ext.horn} the inclusions
    \begin{displaymath}
      \xymatrix@C=8em@R=0.5em{*++{\Lambda^a_k[n]}
        \ar@{u(->}[r]^{\textstyle \subseteq_r} & *++{\Delta^a_k[n]}}
    \end{displaymath}
    (cf. notation~\ref{stan.strat}) for $n\geq 2$ and $k=1,2,...,n-1$
    which we call {\em admissible horn extensions}), and
  \item\label{prim.f-ext.one} the unique surjection:
    \begin{displaymath}
      \xymatrix@C=8em@R=0.5em{*++{\Delta[1]_t}
        \ar[r]^{\textstyle d} & *++{\Delta[0]}}
    \end{displaymath}
  \end{enumerate}
\end{defn}

\begin{defn}[complicial set] A stratified set $X$ is
  said to be a {\em complicial set\/} if it is orthogonal to each
  primitive f-extension.  Let $\Comp$ denote the full subcategory of
  $\Strat$ whose objects are the complicial sets. Of course, since
  every primitive t-extension is a primitive f-extension it follows
  that every complicial set is pre-complicial or, in other words, that
  $\Comp$ is a full subcategory of $\Precomp$.
  
  As in definition~\ref{defn.precomp}, we know that the stratified
  sets $\Delta^a_k[n]'$, $\Delta^a_k[n]''$, $\Delta^a_k[n]$,
  $\Lambda^a_k[n]$, $\Delta[1]$ and $\Delta[0]$ are all finitely
  presentable in $\Strat$ (cf observation~\ref{deltat.dense}). It
  follows, therefore, that the set of primitive f-extensions is also
  an FP-regulus and so, by observation~\ref{full.refl.LFP}, we see
  that $\Comp$ is an LFP-category which is reflective in $\Strat$
  (with associated idempotent monad $\pair<\refl_c;\eta^c>$).

  We let $\mathbb{T}_{\Comp}$ denote the LE-theory of complicial sets,
  which we construct from $\mathbb{T}_{\Strat}$, the LE-theory of
  stratified sets, by adding the set of primitive f-extensions to its
  FP-regulus.
\end{defn}

\begin{obs}[describing complicial sets explicitly]\label{comp.expl}
  More explicitly, by Yoneda's lemma a pre-complicial set $A$
  satisfies these conditions if and only if
  \begin{enumerate}[(a)]
  \item\label{comp.expl.hf} for each $(n-1)$-dimensional $k$-horn
    ($n\geq 2$ and $0<k<n$) $\{x_i\mid i\neq k\}$ in $X$ which is
    admissible (in the sense of observation~\ref{admiss.horn}) there
    is a unique thin $n$-simplex $x\in tX$, called a {\em thin
      filler}, such that $x\cdot\face^n_i=x_i$ for $i\neq k$, and
  \item\label{comp.expl.deg} each thin 1-simplex of $X$ is degenerate.
  \end{enumerate}
\end{obs}

\begin{defn}[f-extensions] We call the $\refl_c$-invertible stratified
  maps {\em f-extensions}. For more on the properties of f-extensions
  please consult subsection~\ref{append.reflective}. We will also use
  the term {\em f-almost\/} as a synonym for $\refl_c$-almost (cf.\ 
  observation~\ref{l-almost}).
\end{defn}

\begin{obs}[duals of complicial sets]\label{comp.duals} We can
  extend the stratified isomorphisms of observation~\ref{stan.dual} to
  get a further family of stratified isomorphisms:
  \begin{equation}\label{horn.duals}
    \let\labelstyle=\textstyle
    \xymatrix@C=8em@R=0.8em{
      \Delta^a_{n-k}[n]\ar[r]^\simeq & (\Delta^a_k[n])^\circ}
  \end{equation}
  Furthermore, it is a matter of a simple calculation to demonstrate
  that the stratified subset $(\Lambda^a_k[n])^\circ\subseteq_r
  (\Delta^a_k[n])^\circ$ corresponds to the subset $\Lambda^a_{n-k}[n]
  \subseteq_r\Delta^a_{n-k}[n]$ under this isomorphism. In other
  words, we have a commutative diagram
  \begin{displaymath}
    \xymatrix@C=8em@R=2.5em{
      \Lambda^a_{n-k}[n]\ar@{u(->}[r]^{\textstyle \subseteq_r}
      \ar[d]_{\cong} &
      \Delta^a_{n-k}[n]\ar[d]^{\cong} \\
      (\Lambda^a_k[n])^\circ\ar@{u(->}[r]_{\subseteq_r} &
      (\Delta^a_k[n])^\circ }
  \end{displaymath}
  for each $n\geq 2$ and $1\leq k < n$. So the dual of each admissible
  horn extension is again (isomorphic to) an admissible horn
  extension.
  
  So, by an argument identical to that of
  observation~\ref{precomp.duals}, we see that $A$ is a complicial set
  if and only if $A^\circ$ is. Equivalently, $\arrow
  f:X->Y.$ is an f-extension if and only if $\arrow
  f^\circ:X^\circ->Y^\circ.$ is.
\end{obs} 

\subsection{Pasting Squares and Filling Lemmas}

The following lemmas are mostly a direct consequence of the properties
of f-extensions listed in observation~\ref{l-inv}. They provide us
with higher level tools for constructing a whole range of complex
f-extensions and, as such, they will be used repeatedly in the sequel.

\begin{obs}[pasting squares]\label{pasting.sq} Suppose that 
  $\inc f:X->Z.$ and $\inc g:Y->Z.$ are stratified inclusions and that
  $Z= f(X)\cup g(Y)$ then the pullback
  \begin{displaymath}
    \let\labelstyle=\textstyle
    \xymatrix@R=2em@C=4em{
      {X \pb{f}{g} Y}
      \pbexcursion\ar@{u(->}[r]^-{\pi_X}
      \ar@{u(->}[d]_{\pi_Y} &
      {X}\ar@{u(->}[d]^f \\
      {Y}\ar@{u(->}[r]_-g &
      {Z} }
  \end{displaymath}
  is also a pushout in $\Strat$; such a square is known as a {\em
    pasting square}. It follows, by
  observation~\ref{l-inv}(\ref{l-inv.po}) that if the upper horizontal
  map $\pi_X$ is an f-extension then so is the lower one $g$.  Usually
  we'll apply this result in one of the following cases:
  \begin{enumerate}[(a)]
  \item\label{pasting.sq.a} $f$ is an arbitrary stratified inclusion
    and $g$ is the inclusion associated with a subset $Y\subseteq_s
    Z$, in which case our pasting square reduces to
    \begin{displaymath}
      \let\labelstyle=\textstyle
      \xymatrix@=2em{
        {f^{-1}(Y)}
        \pbexcursion
        \ar@{u(->}[r]^-{\subseteq_s}
        \ar@{u(->}[d]_f &
        {X}\ar@{u(->}[d]^f \\
        {Y}\ar@{u(->}[r]_-{\subseteq_s} &
        {f(X)\cup Y} }
    \end{displaymath}
    and we infer that $\overinc\subseteq_s:Y->f(X)\cup Y.$
    is an f-extension from the fact that $\overinc\subseteq_s:
    f^{-1}(Y)->X.$ is an f-extension (by
    observation~\ref{l-inv}(\ref{l-inv.po})). 
  \item\label{pasting.sq.b} $f$ is the inclusion associated with
    $X\subseteq_s Z$ and $g$ is the inclusion associated with a subset
    $Y\subseteq_s Z$, in which case our pasting square reduces to
    \begin{displaymath}
      \let\labelstyle=\textstyle
      \xymatrix@=2em{
        {X\cap Y}\pbexcursion
        \ar@{u(->}[r]^-{\subseteq_s}
        \ar@{u(->}[d]_{\subseteq_s} &
        {X}\ar@{u(->}[d]^{\subseteq_s} \\
        {Y}\ar@{u(->}[r]_-{\subseteq_s} &
        {X\cup Y} }
    \end{displaymath}
    and we infer that $\overinc\subseteq_s:Y->X\cup Y.$ is
    an f-extension from the fact that $\overinc\subseteq_s: X\cap
    Y->X.$ is an f-extension (by
    observation~\ref{l-inv}(\ref{l-inv.po})).
  \end{enumerate}
\end{obs}

\begin{lemma}[pasting lemma for f-extensions]\label{pasting.lemma} 
  Suppose that $\arrow f:X->Y.$ is a stratified map and that we have
  stratified subsets $X^{(i)}\subseteq_s X$ and $Y^{(i)}\subseteq_s Y$
  with $X=\bigcup_{i=1}^r X^{(i)}$, $Y=\bigcup_{i=1}^r Y^{(i)}$ and
  $X^{(i)}\subseteq f^{-1}(Y^{(i)})$ for $i=1,\dots,r$ such that $f$
  restricts to
  \begin{enumerate}[(a)]
  \item an f-extension $\arrow f^{(i)}:X^{(i)}->Y^{(i)}.$ for each $i$ and
  \item an f-extension $\arrow f^{(i,j)}:X^{(i)}\cap X^{(j)} -> Y^{(i)}\cap
    Y^{(j)}.$ for each $i<j$
  \end{enumerate}
  then $\arrow f:X->Y.$ itself is also an f-extension.
\end{lemma}

\begin{proof} 
  From observations~\ref{union.qua.widepo}
  and~\ref{stratified.(co)limits} we know that the unions $\bigcup_l
  X^{(l)}$ and $\bigcup_l Y^{(l)}$ may be expressed as wide pushouts
  as depicted in display~(\ref{union.qua.widepo.diag}).  Furthermore,
  the various restrictions of the statement fit together into
  commuting diagrams of the form
  \begin{displaymath}
    \let\labelstyle=\textstyle
    
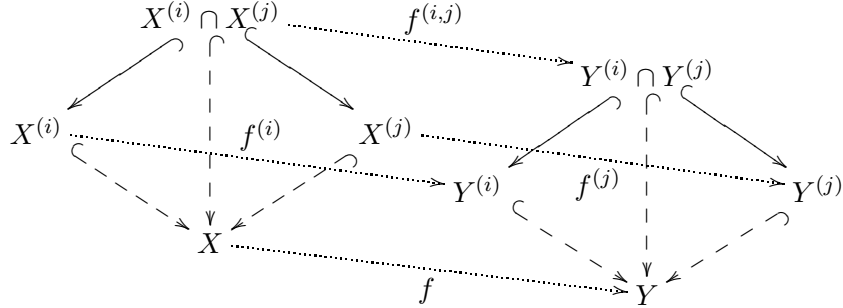
\xymatrix@R=0.5em@C=0.75em{
      && {X^{(i)}\cap X^{(j)}}\ar@{u(->}[ddll]\ar@{u(-->}[dddd]\ar@{u(->}[ddrr] 
      \ar@{.>}[rrrrrd]^{f^{(i,j)}} && & && && \\
      && && & && {Y^{(i)}\cap Y^{(j)}}\ar@{u(->}[ddll]\ar@{u(-->}[dddd]
      \ar@{u(->}[ddrr] && \\
      {X^{(i)}}\ar@{u(-->}[ddrr]\ar@{.>}[rrrrrd]^{f^{(i)}} && && 
      {X^{(j)}}\ar@{u(-->}[ddll]\ar@{.>}[rrrrrd]_{f^{(j)}} & && && \\
      && && & {Y^{(i)}}\ar@{u(-->}[ddrr] && && {Y^{(j)}}\ar@{u(-->}[ddll]\\
      && {X}\ar@{.>}[rrrrrd]_{f} && & && && \\
      && && & && {Y} &&}
  \end{displaymath}
  for each $i<j$. These demonstrate that $\arrow f:X->Y.$ is the
  unique map induced by the restrictions at the top of our diagram,
  each one of which is an f-extension (by assumption), between the
  wide pushouts at the bottom.  Thus the postulated result follows
  from the fact that the class of f-extensions is closed under
  colimits (see observation~\ref{l-inv}).
\end{proof}

\begin{cor}\label{pasting.cor} Suppose that $X$ is a
  stratified subset of $Y\in\Strat$ and that $Y^{(i)}$ ($i=1,...,r$) is a
  family of stratified subsets of $Y$ such that
  \begin{enumerate}[(a)]
  \item the inclusion $\overinc\subseteq_s:X-> X\cup Y^{(i)}.$ is an
    f-extension for each $i$ and
  \item the inclusion $\overinc\subseteq_s:X-> X\cup(Y^{(i)}\cap Y^{(j)}).$ is
    an f-extension for each $i<j$
  \end{enumerate}
  then the inclusion $\overinc\subseteq_s:X->X\cup (\bigcup_{i=1}^r
  Y^{(i)}).$ is also an f-extension.
\end{cor}

\begin{proof}
  Apply lemma~\ref{pasting.lemma} to the inclusion map
  $\overinc\subseteq_s:X->Y.$ and the families of subsets
  $\{X\}_{i=1}^r$ and $\{X\cup Y^{(i)}\}_{i=1}^r$, for which the
  conditions given in that lemma simply reduce to those given in the
  statement of this corollary.
\end{proof}

\begin{lemma}[minor filling lemma]\label{minor.fill} If $X$ is
  a regular stratified subset of $\Delta^a_k[n]$ and
  \begin{enumerate}[(a)]
  \item\label{minor.cond.i} $\face^n_{k-1}$ and $\face^n_{k+1}$ are
    elements of $X$ and
  \item\label{minor.cond.ii} if $\mu$ is a $k$-divided $r$-simplex in
    $\Delta^a_k[n]$, $l\in[r]$ is the unique integer such that
    $\mu(l)=k$ and $\mu\circ\face^r_{l}$ is an element of the subset
    $X$ then $\mu$ is also an element of $X$,
  \end{enumerate}
  then the inclusion map $\overinc\subseteq_r:X->\Delta^a_k[n].$
  is an f-extension.
\end{lemma} 

\begin{proof} By induction on $n$, using lemma~\ref{pasting.lemma}:
\vspace{1ex}
  
  \noindent{\bf n=2:} The only stratified subset satisfying the
  conditions of the lemma is $\Lambda^a_1[2]\subseteq_r \Delta^a_1[2]$
  which is a primitive f-extension. \vspace{1ex}

  \noindent{\bf the inductive case:} Suppose that if $m<n$ then the
  result is true for all subsets of $\Delta^a_l[m]$ ($l=1,\dots,m-1$)
  satisfying the conditions of the lemma.
    
  Assume that $X\subseteq_r\Delta^a_k[n]$ satisfies the
  conditions of the lemma. If $\id_{[n]}\in X$ then
  $X=\Delta^a_k[n]$ and there is nothing to prove, so assume
  that $\id_{[n]}\notin X$ from which it follows, by condition
  (\ref{minor.cond.ii}), that $\face^n_k\notin X$.  In other words,
  $X$ is a stratified subset of $\Lambda^a_k[n]$.
    
  Suppose now that $\mu$ is a non-degenerate $k$-divided
  $m$-simplex in $\Delta^a_k[n]$, that is to say a face operator
  $\arrow \mu:[m]->[n].$ with $k-1,k,k+1\in\im(\mu)$ and let
  $l\in[m]$ be the unique number such that $\mu(l)=k$, for which we
  also know that $k-1=\mu(l-1)$ and $k+1=\mu(l+1)$.  Clearly an
  element $\nu$ of $\Delta[m]$ is $l$-divided if and only if
  $\Delta(\mu)(\nu)=\mu\circ\nu\in\Delta[n]$ is $k$-divided.  It
  follows that the simplicial map $\arrow \Delta(\mu):\Delta[m]
  ->\Delta[n].$ extends to a regular stratified map from
  $\Delta^a_l[m]$ to $\Delta^a_k[n]$. So consider the inverse image
  $\Delta(\mu)^{-1}(X)\subseteq_r\Delta^a_l[m]$ of
  $X\subseteq_r\Delta^a_k[n]$, which we know to be a regular
  stratified subset by observation~\ref{strat.im}. It is
  also easily established that this inverse image satisfies the
  conditions (\ref{minor.cond.i}) and (\ref{minor.cond.ii}) in the
  statement of this lemma:
  \begin{enumerate}[(i)]
  \item Given that $\mu$ is $k$-divided with $\mu(l)=k$ and that
    $l-1\notin\im(\face^m_{l-1})$, we may infer that $k-1\notin
    \im(\mu\circ\face^m_{l-1})$ and, consequently, that
    $\mu\circ\face^m_{l-1}=\Delta(\mu)(\face^m_{l-1})$ is an
    $(m-1)$-dimensional face of $\face^n_{k-1}$ in $\Delta^a_k[n]$.
    However, by assumption~(\ref{minor.cond.i}) for $X$, we know that
    $\face^n_{k-1}$ is in $X$ as is each of its faces, including in
    particular $\Delta(\mu)(\face^m_{l-1})$. It follows that
    $\face^m_{l-1}$ is in the inverse image $\Delta(\mu)^{-1}(X)$
    as required. An identical argument demonstrates that
    $\face^m_{l+1}$ is also in $\Delta(\mu)^{-1}(X)$.
  \item Suppose that $\nu$ is an $l$-divided $r$-simplex in
    $\Delta^a_l[m]$, that $j\in[r]$ has $\nu(j)=l$ and that
    $\nu\circ\face^r_{j}$ is in $\Delta^{-1}(X)$. By
    definition, this latter condition holds iff
    $\Delta(\mu)(\nu\circ\face^r_{j})=\mu\circ(\nu\circ
    \face^r_{j})=(\mu\circ\nu)\circ\face^r_{j}$ is an element of
    $X$. Furthermore, we know that $\mu$ is $k$-divided with
    $\mu(l)=k$ and that $\nu$ is $l$-divided with $\nu(j)=l$ so it
    follows that $\mu\circ\nu$ is $k$-divided with
    $\mu\circ\nu(j)=k$. Now apply condition~(\ref{minor.cond.ii})
    for the subset $X\subseteq_r\Delta^a_k[n]$ to demonstrate
    that if $(\mu\circ\nu)\circ\face^r_{j}$ is an element of
    $X$ then so is $\mu\circ\nu=\Delta(\mu)(\nu)$.
    From this we may infer that $\nu$ is an element of
    $\Delta(\mu)^{-1}(X)$ as required.
  \end{enumerate}
  It follows, by our inductive hypothesis, that the inclusion
  $\overinc\subseteq_r:\Delta(\mu)^{-1}(X)-> \Delta^a_l[m].$ is an
  f-extension and since it features in the pasting square
  \begin{displaymath}
    \let\labelstyle=\textstyle
    \xymatrix@R=2.5em@C=4em{
      {\Delta(\mu)^{-1}(X)}
      \ar@{u(->}[r]^-{\subseteq_r}
      \ar@{u(->}[d]_{\Delta(\mu)} &
      {\Delta^a_l[m]} 
      \ar@{u(->}[d]^{\Delta(\mu)} \\
      {X} \ar@{u(->}[r]_-{\subseteq_r} &
      \poexcursion {\Delta(\mu)(\Delta^a_l[m])\cup X}}
  \end{displaymath}
  we may infer that the inclusion
  $\overinc\subseteq_r:X->\Delta(\mu) (\Delta^a_l[m])\cup X.$ is
  also an f-extension (by
  observation~\ref{pasting.sq}(\ref{pasting.sq.a})). Furthermore,
  since we know that $\arrow \Delta(\mu)
  :\Delta^a_l[m]->\Delta^a_k[n].$ is regular, we may apply
  observation~\ref{strat.im} to show that the image
  $\Delta(\mu)(\Delta^a_l[m])$ is a regular subset of
  $\Delta^a_k[n]$ and indeed that it must be $\stratc{\{\mu\}}
  \subseteq_r\Delta^a_k[n]$, the regular stratified subset generated
  by the $m$-simplex $\mu\in\Delta^a_k[n]$. To summarise, we have
  established that the inclusion map $\overinc\subseteq_r:X->
  X\cup\stratc{\{\mu\}}.$ is an f-extension for each
  $k$-divided operator $\arrow\mu:[m]->[n].$.
  
  Now, for each $i\neq k-1,k,k+1$ we know that $\face^n_i$ is
  $k$-divided, so let $Y^{(i)}\defeq\stratc{\{\face^n_i\}}$ and the
  result of the last paragraph implies that
  $\overinc\subseteq_r:X->X\cup Y^{(i)}.$ is an f-extension.
  Furthermore, if $i<j$ ($i,j\neq k-1,k,k+1$) then we have
  $Y^{(i)}\cap Y^{(j)} = \stratc{\{\face^n_i\}} \cap
  \stratc{\{\face^n_j\}} = \stratc{\{\face^n_j\circ \face^{n-1}_i\}}$
  (cf.\ notation~\ref{simplex.bound}) and
  $\face^n_j\circ\face^{n-1}_i$ is also $k$-divided thus the result of
  the last paragraph implies that $\overinc\subseteq_r:X
  ->X\cup(Y^{(i)}\cap Y^{(j)}).$ is an f-extension as well. In other
  words, the family $\{Y^{(i)}\}_{i\neq k-1,k,k+1}$ satisfies the
  conditions of corollary~\ref{pasting.cor}, from which it follows
  that the inclusion map $\overinc\subseteq_r:X->X\cup\left(
    \bigcup_{i\neq k-1,k,k+1} Y^{(i)}\right).$ is also an f-extension.
  
  However, the faces $\face^n_{k-1}$ and $\face^n_{k+1}$ are elements
  of $X\subseteq_r\Lambda^a_k[n]$ (by condition~(\ref{minor.cond.i})
  in the statement of this lemma) and the faces $\face^n_i$ ($i\neq
  k-1,k,k+1$) are the remaining non-degenerate $(n-1)$-simplices of
  $\Lambda^a_k[n]$, from which it follows that
  $X\cup\left(\bigcup_{i\neq k-1,k,k+1} Y^{(i)}\right)=\Lambda^a_k[n]$. So
  the result of the last paragraph simply demonstrated that the
  inclusion $\overinc\subseteq_r:X->\Lambda^a_k[n].$ is an
  f-extension, which we now compose with the primitive f-extension
  $\overinc\subseteq_r:\Lambda^a_k[n]->\Delta^a_k[n].$ in order to
  prove that the resulting inclusion $\overinc\subseteq_r:X->
  \Delta^a_k[n].$ is also an f-extension as stated.
\end{proof}

The following variant of this result will prove particularly useful:

\begin{cor}\label{minor.fill.cor} Suppose that $N$ is a stratified
  set in the subcategory $\Simplex$ of
  definition~\ref{thinner.simplices} and that $X\subseteq_r N$ is a
  regular subset which satisfies the conditions of
  lemma~\ref{minor.fill} and these also satisfy
  \begin{enumerate}[(a)]
  \item\label{mfc.adm} $\Delta^a_k[n]\subseteq N$, in other
    words every $k$-divided simplex of $\Delta[n]$ is in $tN$, and
  \item\label{mfc.t-ext} whenever $\mu$ is a $k$-divided $r$-simplex
    in $\Delta[n]$, $l$ is the unique integer with $\mu(l)=k$ and
    $\mu\circ\face^r_{l}\in tN\smallsetminus X$ then both of the
    faces $\mu\circ\face^r_{l-1}$ and $\mu\circ\face^r_{l+1}$ are also
    in $tN$.
  \end{enumerate}
  then the inclusion $\overinc\subseteq_r:X->N.$ is an f-extension.
\end{cor}

\begin{proof}
  Consider the stratified subsets $X\subseteq_r N$ and
  $\Delta^a_k[n]\subseteq_e N$. Since $X$ is a regular subset of $N$
  it follows that $X\cap\Delta^a_k[n]$ is a regular subset of
  $\Delta^a_k[n]$, furthermore $X\cap\Delta^a_k[n]$ satisfies the
  conditions of lemma~\ref{minor.fill} since it has the same
  underlying simplicial set as $X$ which does so by assumption.
  Consequently, we may apply that lemma to demonstrate that the
  inclusion $\overinc\subseteq_r:X \cap\Delta^a_k[n]->\Delta^a_k[n].$
  is an f-extension and, by
  observation~\ref{pasting.sq}(\ref{pasting.sq.b}), we know that this
  map features in a pasting square
  \begin{displaymath}
    \let\labelstyle=\textstyle
    \xymatrix@R=2em@C=2em{
      {X\cap\Delta^a_k[n]}
      \ar@{u(->}[r]\ar@{_{(}->}[d] & 
      \Delta^a_k[n]\ar@{u(->}[d] \\
      {X}\ar@{u(->}[r] & \poexcursion
      {X\cup\Delta^a_k[n]}
      }
  \end{displaymath}
  from which we may conclude that the inclusion $\overinc\subseteq_r:
  X->X\cup\Delta^a_k[n].$ is also an f-extension.
  
  We can also show that the inclusion $\overinc:X\cup
  \Delta^a_k[n]->N.$ is a t-extension.  To this
  end, suppose that $\nu$ is thin in $N$ and not thin in
  $X\cup\Delta^a_k[n]$. Then $\nu$ is necessarily
  non-degenerate, since otherwise it would be thin in
  $X\cup\Delta^a_k[n]$.
  
  Notice that if a simplex $\nu$ of $\Delta[n]$ has
  $k-1\notin\im(\nu)$ (resp.\ $k+1\notin\im(\nu)$) then it is a face
  of $\face^n_{k-1}$ (resp.\ $\face^n_{k+1}$) which is an element of
  $X$, by assumption~(\ref{minor.cond.i}) of lemma~\ref{minor.fill}.
  But the stratified subset $X$ is closed in $N$ under the right
  action of $\Delta$, so if it contains the simplex $\face^n_{k-1}$
  (resp. $\face^n_{k+1}$) then it must also contain all of its faces
  and it follows that $\nu\in X$. Furthermore $X$ is a {\bf regular}
  subset of $N$, so in such a case $\nu$ is thin in $X$ iff it is thin
  in $N$.
  
  So suppose, for a contradiction, that $k-1$ (resp. $k+1$) is not an
  element of $\im(\nu)$, then we could infer (by the result of the
  last paragraph) that $\nu$ is in $X$ and that it thin there,
  since we assumed that it was thin in $N$. However, this
  would contradict our assumption that $\nu$ was {\bf not} thin in
  $X\cup\Delta^a_k[n]$.
  
  It follows that both $k-1$ and $k+1$ are elements of $\im(\nu)$.
  However, $k$ is not in $\im(\nu)$, because if it were $\nu$ would be
  $k$-divided and thus thin in $X\cup \Delta^a_k[n]$. So let
  $\mu$ be the unique non-degenerate simplex of $\Delta[n]$ which has
  $\im(\mu)=\im(\nu)\cup\{k\}$ and let $l$ be the unique integer
  such that $\mu(l)=k$.
  
  By construction $\mu$ is $k$-divided, in particular it is
  $k$-admissible in $X\cup\Delta^a_k[n]$, and
  $\mu\circ\face^n_{l} = \nu$. We assumed that $\nu$ was thin in
  $N$ therefore, by assumption~(\ref{mfc.t-ext}) of this
  corollary, we can infer that both of $\mu\circ\face^n_{l-1}$ and
  $\mu\circ\face^n_{l+1}$ are also thin in $N$. But $k-1$
  (resp. $k+1$) is not an element of $\im(\mu\circ\face^n_{l-1})$
  (resp.  $\im(\mu\circ\face^n_{l+1})$) so, by our observation of a
  few paragraphs ago, both of these simplices are in $X$ and
  they are thin in there, thus they are thin in
  $X\cup\Delta^a_k[n]$.
  
  But this is precisely what we need in order to show that
  $\overinc:X\cup\Delta^a_k[n]->N.$ is a t-extension by invoking
  observation~\ref{t-ext.expl.obs}, and all t-extensions are
  f-extensions (cf.\ observation~\ref{l-inv}(\ref{l-inv.inc})).
  
  So both of $\overinc:X->X\cup\Delta^a_k[n].$ and
  $\overinc:X\cup\Delta^a_k[n]->N.$ are f-extensions and it follows
  that their composite $\overinc:X->N.$ is also an f-extension as
  required.
\end{proof}

\subsection{Tensor Products and Complicial Sets}

First we show that complicial sets are all well tempered as indicated
in section~\ref{sect.filt.semi}. Then we push on to prove a range of
particularly important results with respect to complicial sets and
their interaction with the biclosed monoidal structure introduced in
section~\ref{precomp.sec}.

\begin{lemma}\label{lemma.predegen}
  Every complicial set is well tempered (cf.
  definition~\ref{well.tempered}). In other words if $A$ is a
  complicial set then every simplex $a\in A$ which is pre-degenerate
  at $k$ is degenerate at $k$.
\end{lemma}

\begin{proof} By induction on $n=\dim(a)$:\vspace{1ex}

\noindent{\bf\boldmath $n=1$:} This is simply (equivalent to) the
axiom that states that a complicial set is orthogonal to
$\arrow d:\Delta[1]_t->\Delta[0].$, as expressed in the form
given in observation~\ref{comp.expl}(\ref{comp.expl.deg}).
\vspace{1ex}

\noindent{\bf\boldmath inductive case:} Suppose that the given result
holds for every $r$-simplex $a\in A$ with $1\geq r<n$ and each
$k\in[r]$.

Let $a\in A$ be an $n$-simplex which is pre-degenerate at $k$,
assume for the moment that $k>0$ (we'll consider the case $k=0$ later
on) and consider the $(n-1)$-dimensional faces $a\cdot\face^n_i$ of
$a$:
\begin{itemize}
\item {\bf\boldmath $i<k$} Let $\arrow \alpha:[r]->[n-1].$ be a
  face operator with $k-1,k\in\im(\alpha)$ then, since $i\leq k-1$, we
  know that $k,k+1\in\im(\face^n_i\circ\alpha)$; it follows, since $a$
  is pre-degenerate at $k$, that $(a\cdot\face^n_i)\cdot\alpha =
  a\cdot(\face^n_i\circ\alpha)$ is thin. Consequently
  $a\cdot\face^n_i$ is pre-degenerate at $k-1$ and so, by our
  inductive hypotheses, it is degenerate at $k-1$; it follows that
  $a\cdot\face^n_i = ((a\cdot\face^n_i)\cdot\face^{n-1}_k)
  \cdot\degen^{n-2}_{k-1} = a\cdot(\face^n_i\circ\face^{n-1}_k
  \circ\degen^{n-2}_{k-1})$.
  
  Applying the simplicial identities of
  observation~\ref{simplicial.idents} (twice), under the assumption
  that $i<k$, we see that $\face^n_i\circ\face^{n-1}_k
  \circ\degen^{n-2}_{k-1}= \face^n_{k+1}\circ\face^{n-1}_i\circ
  \degen^{n-2}_{k-1} = \face^n_{k+1}\circ\degen^{n-1}_k\circ\face^n_i$
  and so the equation at the end of the last paragraph becomes
  $a\cdot\face^n_i = a\cdot(\face^n_{k+1}\circ\degen^{n-1}_k
  \circ\face^n_i) = ((a\cdot\face^n_{k+1})\cdot\degen^{n-1}_k)
  \cdot\face^n_i$.
  
\item {\bf\boldmath $i>k+1$} Let $\arrow \alpha:[r]->[n-1].$ be
  a face operator with $k,k+1\in\im(\alpha)$ then, since $i>k+1$, we
  know that $k,k+1\in\im(\face^n_i\circ\alpha)$; it follows, since $a$
  is pre-degenerate at $k$, that $(a\cdot\face^n_i)\cdot\alpha =
  a\cdot(\face^n_i\circ\alpha)$ is thin. Consequently
  $a\cdot\face^n_i$ is pre-degenerate at $k$ and so, by our inductive
  hypotheses, it is degenerate at $k$; it follows that
  $a\cdot\face^n_i = ((a\cdot\face^n_i)\cdot\face^{n-1}_{k+1})
  \cdot\degen^{n-2}_k = a\cdot(\face^n_i\circ\face^{n-1}_{k+1}
  \circ\degen^{n-2}_k)$.
 
  Applying the simplicial identities of
  observation~\ref{simplicial.idents} (twice), under the assumption
  that $i>k+1$, we see that $\face^n_i\circ\face^{n-1}_{k+1}
  \circ\degen^{n-2}_k= \face^n_{k+1}\circ\face^{n-1}_{i-1}\circ
  \degen^{n-2}_k = \face^n_{k+1}\circ\degen^{n-1}_k\circ\face^n_i$ and
  so the equation at the end of the last paragraph becomes
  $a\cdot\face^n_i = a\cdot(\face^n_{k+1}\circ\degen^{n-1}_k
  \circ\face^n_i) = ((a\cdot\face^n_{k+1})\cdot\degen^{n-1}_k)
  \cdot\face^n_i$.
  
\item {\bf\boldmath $i=k+1$} Since $\degen^{n-1}_k\circ\face^n_{k+1}
  =\id_{[n-1]}$ we have
  $((a\cdot\face^n_{k+1})\cdot\degen^{n-1}_k)\cdot \face^n_{k+1} =
  (a\cdot\face^n_{k+1})\cdot(\degen^{n-1}_k\circ \face^n_{k+1}) =
  (a\cdot\face^n_{k+1})\cdot\id_{[n-1]} = a\cdot \face^n_{k+1}$.

\end{itemize}

Now we are in a position to compare the simplices $a$ and
$(a\cdot\face^n_{k+1})\cdot\degen^{n-1}_k$; the equations derived
above simply state that the corresponding horns $\{a\cdot\face^n_i\mid
i\neq k\}$ and
$\{((a\cdot\face^n_{k+1})\cdot\degen^{n-1}_k)\cdot\face^n_i\mid i\neq
k\}$ are in fact identical. Furthermore, since every $k$-divided
operator $\alpha$ has $k,k+1\in\im(\alpha)$ and $a$ is pre-degenerate
at $k$, we may infer that $a$ is $k$-admissible. In summary, $a$ and
$(a\cdot\face^n_{k+1})\cdot\degen^{n-1}_k$ are both thin fillers of
the admissible horn $\{a\cdot\face^n_i\mid i\neq k\}$.

However we know that thin fillers of admissible horns in $A$
are unique, since it is a complicial set, therefore
$a=(a\cdot\face^n_{k+1})\cdot\degen^{n-1}_k$ and thus $a$ is
degenerate at $k$.

All that remains is to prove that the result also holds when $k=0$. To
do so, observe that an $n$-simplex $a\in A$ is degenerate
(resp.  pre-degenerate) at $k$ if and only if it is degenerate (resp.
pre-degenerate) at $(n-1-k)$ when considered as a simplex in the dual
$A^\circ$. So if our simplex $a$ is pre-degenerate at $0$ as
an element of $A$ then it is pre-degenerate at $(n-1)$ when
considered as an element of the dual complicial set $A^\circ$.
Applying the result we've already proved, we may infer that $a$ is
degenerate at $(n-1)$ in $A^\circ$ and thus is degenerate at
$0$ in $A$ as required.
\end{proof}

\begin{cor}\label{prim.d.tens} For each stratified set
  $Y$ the stratified map
  \begin{equation}\label{d.tens}
    \let\labelstyle=\textstyle
    \xymatrix@R=0ex@C=10em{
      {\Delta[1]_t\otimes Y}
      \ar[r]^{d\otimes Y} &
      {Y}}
  \end{equation}
  is an f-extension.
\end{cor}

\begin{obs}\label{prim.d.tens.note}
  Before we proceed with this proof, we introduce a little notation
  for the elements of $\Delta[1]$. For $k\in[r+1]$ let
  $\arrow \rho^r_k:[r]->[1].$ be the simplicial operator given
  by
  \begin{displaymath}
    \rho^r_k(i) = \left\{
      \begin{array}{lp{1.5cm}}
        0 & if $i<k$, \\
        1 & if $i\geq k$.
      \end{array}
    \right.
  \end{displaymath}
  for which the following, easily demonstrated, identities will be of
  use in this proof:
  \begin{equation}\label{rho.equal}
    \begin{array}{rcl}
      \rho^{r+1}_{k+1}\circ\face^{r+1}_k & = & \rho^r_k \\
      \rho^{r+1}_k\circ\face^{r+1}_k & = & \rho^r_k
    \end{array}
  \end{equation}
\end{obs}
  
\begin{proof}
  We know that the map $\arrow d:\Delta[1]_t ->\Delta[0].$ has
  $\arrow \Delta(\vertex^1_0):\Delta[0] ->\Delta[1]_t.$ as a right
  inverse so, by the functoriality of $\otimes$, we also know that
  $\Delta(\vertex^1_0)\otimes Y$ is a right inverse of
  $d\otimes Y$.  We'll actually prove that $\Delta(\vertex^1_0)\otimes
  Y$ is an f-extension and then appeal to
  observation~\ref{l-inv}(\ref{l-inv.retr}) to prove that $d\otimes Y$
  is also an f-extension.
  
  Since $d\otimes Y$ is left inverse to $\Delta(\vertex^1_0)\otimes
  Y$, we may extend any stratified map $\arrow f:Y->A.$ along
  $\Delta(\vertex^1_0)\otimes Y$ by simply composing it with $d\otimes
  Y$. To complete our proof we must demonstrate that if $A$ is a
  complicial set then this extension is unique or, in other words,
  that if $\arrow g,g':\Delta[1]_t\otimes Y->A.$ are any pair of
  stratified maps with $g\circ(\Delta(\vertex^1_0)\otimes Y) =
  g'\circ(\Delta(\vertex^1_0)\otimes Y)$ then we have $g=g'$.
  
  Observe that in order to verify this latter condition it suffices to
  show that for any stratified $\arrow g:\Delta[1]_t\otimes Y->A.$
  with $A$ complicial we have $g\pair<\rho^r_k;y>=
  g\pair<\rho^r_{k+1};y>$ for each $k\in[r]$; this follows because any
  simplicial operator $\arrow \rho:[r]->[1].$ is of the form
  $\rho^r_i$ for some $i\in[r+1]$ so the stated equality (for each
  $k$) would give $g\pair<\rho;y>=g\pair<\rho^r_i;y>=
  g\pair<\rho^r_{i-1};y>= g\pair<\rho^r_{i-2};y> = \cdots =
  g\pair<\rho^r_0;y>$. However $\rho^r_0=\vertex^1_0\circ\eta^r$ so we
  have $\pair<\rho^r_0;y>=(\Delta(\vertex^1_0)\otimes Y)(y)$ and it
  follows that if $g$ and $g'$ are as in the previous paragraph we would
  have $g\pair<\rho;y>=(g\circ(\Delta(\vertex^1_0)\otimes Y))(y) =
  (g'\circ(\Delta(\vertex^1_0)\otimes Y))(y) = g'\pair<\rho;y>$ as
  required.
  
  To prove the equality $g\pair<\rho^r_k;y>= g\pair<\rho^r_{k+1};y>$
  (for a given $k\in[r]$) consider the $(r+1)$-simplex
  $\pair<\rho^{r+1}_{k+1}; y\cdot\degen^r_k>$ in $\Delta[1]_t\otimes
  Y$ which, by the equalities in (\ref{rho.equal}), has:
  \begin{equation}\label{faces.eqn}
    \begin{aligned}
      \pair<\rho^{r+1}_{k+1};y\cdot\degen^r_k>\cdot\face^{r+1}_k & {} =
      \pair<\rho^{r+1}_{k+1}\circ\face^{r+1}_k;
      y\cdot(\degen^r_k\circ\face^{r+1}_k)>  = 
      \pair<\rho^r_k; y> \\
      \pair<\rho^{r+1}_{k+1};y\cdot\degen^r_k>\cdot\face^{r+1}_{k+1} &
      {} = \pair<\rho^{r+1}_{k+1}\circ\face^{r+1}_{k+1};
      y\cdot(\degen^r_k\circ\face^{r+1}_{k+1})> =
      \pair<\rho^r_{k+1}; y>
    \end{aligned}
  \end{equation}
  Now observe that we have $(d\otimes Y)\pair<\rho^{r+1}_{k+1};y\cdot
  \degen^r_k>=y\cdot\degen^r_k$ which is degenerate, and thus
  pre-degenerate, at $k$ in $Y$.  However we know that $d$ is regular
  from which we may infer, by observation~\ref{tens.thin.refl}, that
  $d\otimes Y$ is also regular and it follows that
  $\pair<\rho^{r+1}_{k+1};y\cdot \degen^r_k>$ is pre-degenerate at $k$
  in $\Delta[1]_t\otimes Y$ (as observed in
  definition~\ref{defn.predegen}).  Finally, we see that
  $g\pair<\rho^{r+1}_{k+1}; y\cdot \degen^r_k>$ is also pre-degenerate
  at $k$ in $A$ and so, since $A$ is complicial, we may apply
  lemma~\ref{lemma.predegen} to show that this simplex is, in fact,
  degenerate at $k$. Now we may infer our desired equality since:
  \begin{displaymath}
    \begin{aligned}[b]
      g\pair<\rho^r_k;y> & {} = 
      g\pair<\rho^{r+1}_{k+1};y\cdot\degen^r_k>\cdot\face^{r+1}_k &&
      \text{from first equality in~(\ref{faces.eqn})}\\
      & {} =  
      g\pair<\rho^{r+1}_{k+1};y\cdot\degen^r_k>\cdot\face^{r+1}_{k+1}
      && \text{since $g\pair<\rho^{r+1}_{k+1};y\cdot\degen^r_k>$ is
        degenerate at $k$.} \\
      & {} = g\pair<\rho^r_{k+1};y> &&
      \text{from second equality in~(\ref{faces.eqn})}
    \end{aligned}\qedhere
  \end{displaymath}
\end{proof}

Now we present the primary combinatorial result of this section, and
indeed of this work, that being lemma~\ref{prim.f.tens}. Firstly
however, we introduce some notation and prove an associated
technical lemma.  

\begin{notation}\label{Hdelta} Suppose that ${N}$ and ${M}$
  are stratified sets in $\Simplex$ and that $\anytens$ represents any
  of the bifunctors $\times$, $\otimes$ or $\pretens$ on $\Strat$ then
  in the sequel we adopt the following notation for various useful
  regular subsets of these stratified sets and their tensor products:
\begin{enumerate}[(i)]
\item As usual, let $\boundary{N}$ and
  $\boundary{M}$ denote the ``boundaries'' of ${N}$ and
  ${M}$, that is their regular subsets consisting of those
  operators which are {\bf not} degeneracy operators.
\item\label{Hdelta.bdary} Furthermore let $\boundary({N}\anytens{M})$
  be the boundary of the tensor product ${N}\anytens{M}$, which is its
  regular subset given by the Leibniz formula:
  $$\boundary({N}\anytens{M}) \defeq (\boundary{N}\anytens{M})\cup
  ({N}\anytens\boundary{M})$$ In other words, a simplex
  $\abpair$ of ${N}\anytens{M}$ is in
  $\boundary({N}\anytens{M})$ if and only if either $\alpha$ or $\beta$
  is {\bf not} a degeneracy operator.
\item\label{Hdelta.H(N,M)} Let $H({N}\anytens{M})$ be the regular stratified
  subset of ${N}\anytens{M}$ on those simplices which do {\bf not} have
  $\pair<n;0>$ as a vertex. Recall that every simplex of
  ${N}\anytens{M}$ is a face of some shuffle and that the
  $\triangleleft$-minimal shuffle $\shuffle{1}=\pair<
  \partproj^{n,m}_1;\partproj^{n,m}_2>$ is the only one
  which has $\pair<n;0>$ as a vertex. Therefore $H({N}\anytens{M})$ is
  precisely the regular stratified subset of ${N}\anytens{M}$ generated
  by the set of shuffles $\{\shuffle{i} : 1 < i \leq \#\pair<n;m>\}$
  and
  \begin{displaymath}
    H({N}\anytens{M})\cup\stratc{\{\shuffle{1}\}} =
    {N}\anytens{M}
  \end{displaymath}
\item Let $\boundary H({N}\anytens{M})$ denote the ``boundary''
  of $H({N}\anytens{M})$
  \begin{displaymath}
    \boundary H({N}\anytens{M})\defeq
    H({N}\anytens{M})\cap 
    \boundary({N}\anytens{M})
  \end{displaymath}
  in other words, the regular stratified
  subset of $H({N}\anytens{M})$ consisting of those simplices
  $\abpair$ for which either $\alpha$ or $\beta$ is {\bf
    not} a degeneracy operator.
\item Given a natural number $i$ with $1\leq i\leq \#\pair<n;m>$ let
  $H_i({N}\anytens{M})$ denote the regular stratified subset
  given by
  \begin{displaymath}
    H_i({N}\anytens{M}) = \boundary H({N}\anytens{M})\cup \stratc{\{\shuffle{i'} \mid i'>i\}} 
  \end{displaymath}
  So $H_{\#\pair<n;m>}({N}\anytens{M})=\boundary
  H({N}\anytens{M})$,
  $H_1({N}\anytens{M})=H({N}\anytens{M})$ and
  $H_{i-1}({N}\anytens{M})$ is the regular stratified subset of
  $H({N}\anytens{M})$ obtained by adjoining the shuffle
  $\shuffle{i}$ to $H_i({N}\anytens{M})$ (in other words
  $H_{i-1}({N}\anytens{M}) = H_i({N}\anytens{M})\cup
  \stratc{\{\shuffle{i}\}}\subseteq_r H({N}\anytens{M})$).
\end{enumerate}
\end{notation}

\begin{lemma}\label{H.analysis}
  Fix $n,m\in\mathbb{N}$ and suppose that $\abpair$ is a
  non-degenerate $r$-simplex of $\Delta[n]\times\Delta[n]$ which is
  a mediator as witnessed by some $k$ with $0<k<r$, furthermore assume
  that $\abpair\cdot\face^r_k$ is in
  $H_{i-1}(\Delta[n]\times\Delta[m])$ then
  \begin{enumerate}[(a)]
  \item\label{H.analysis.1} $\abpair$ is in
    $H_{i-1}(\Delta[n]\times\Delta[m])$ and
  \item\label{H.analysis.2} $\abpair\cdot\face^r_{k-1}$
    and $\abpair\cdot\face^r_{k+1}$ are both in
    $H_i(\Delta[n]\times\Delta[m])$.
  \end{enumerate}
\end{lemma}

\begin{proof}
  Let $w$ be the largest integer for which the simple
  $\abpair\cdot\face^r_k$ is a face of the shuffle $\shuffle{w}$.  By
  observation~\ref{shuffles.maxl} we know that the associated operator
  $\arrow\topop{w}:[n-1]->[m].$ is given by the formula in
  display~\eqref{gamma.ab.def} of the proof of
  lemma~\ref{shuffles}(\ref{prod.simp.5}), which in the case of the
  particular value $\topop{w}(\alpha(k-1))$ reduces to:
  \begin{displaymath}
    \topop{w}(\alpha(k-1)) = \min\{\beta\circ\face^r_k(l)\mid l\in[r-1]
    \wedge\alpha\circ\face^r_k(l)>\alpha(k-1)\} 
  \end{displaymath}
  Now since $\abpair$ is {\em non-degenerate\/} and satisfies
  the mediation condition at $k$ we know that the (in)equalities
  \begin{equation}\label{H.anal.ineq}
    \alpha(k-1)=\alpha(k)< \alpha(k+1)\text{ and }
    \beta(k-1)<\beta(k)=\beta(k+1)
  \end{equation}
  hold. In particular, we see that $\alpha\circ\face^r_{k}(k) =
  \alpha(k+1) > \alpha(k-1)$ and $\alpha\circ \face^r_{k}(k-1) =
  \alpha(k-1) \not> \alpha(k-1)$ so it clearly follows that
  $\topop{w}(\alpha(k-1))=\beta(k+1)$. Applying these (in)equalities
  again, this time to $\beta$, we see that
  $\topop{w}(\alpha(k-1))=\beta(k+1)>\beta(k-1)\geq 0$ so $\topop{w}$
  is a non-zero operator and therefore the corresponding shuffle
  $\shuffle{w}$ is not $\pair<\partproj^{n,m}_1;\partproj^{n,m}_2>$
  and is thus an element of $H(\Delta[n]\times\Delta[m])$.
  
  Applying lemma~\ref{shuffles}(\ref{prod.simp.5a}) we have
  $\abpair\cdot \face^r_k=\shuffle{w}\cdot(\alpha\circ \face^r_k+
  \beta\circ\face^r_k)=(\shuffle{w}\cdot(\alpha+\beta)) \cdot
  \face^r_k$, where the final equality is a consequence of the fact
  that $\arrow{+}:\Delta[n]\times\Delta[m]->\Delta[n+m].$ is a
  simplicial map. It follows that $\abpair$ and the face
  $\shuffle{w}\cdot(\alpha + \beta)$ agree at vertices $j\in[r]$ with
  $j\neq k$, and so they are equal if and only if they also agree at
  vertex $k$. However, by lemma~\ref{shuffles}(\ref{prod.simp.5a}),
  this is the case if and only if the inequalities
  $\topop{w}(\alpha(k)-1)\leq \beta(k)\leq\topop{w}(\alpha(k))$ of
  display~(\ref{shuff.face.ineq}) hold, which is certainly the case
  since we may apply the equalities of the last paragraph to show that
  $\topop{w}(\alpha(k))=\topop{w}(\alpha(k-1))=\beta(k+1)=\beta(k)$.
  
  We may now dispense with the trivial case, wherein
  $\abpair\cdot\face^r_k=\pair<\alpha\circ\face^r_k;\beta\circ
  \face^r_k>$ is in the boundary $\boundary
  H(\Delta[n]\times\Delta[m])\subseteq
  H_{i-1}(\Delta[n]\times\Delta[m])$ which implies that one of
  $\alpha\circ\face^r_k$ or $\beta\circ\face^r_k$ is not surjective.
  However, consulting display~\ref{H.anal.ineq} we know that
  $\alpha(k-1)=\alpha(k)$ and $\beta(k)=\beta(k+1)$ and it follows
  easily that $\im(\alpha)=\im(\alpha\circ \face^r_k)$ and
  $\im(\beta)=\im(\beta\circ\face^r_k)$, consequently either $\alpha$
  or $\beta$ is also non-surjective and thus $\abpair$ is
  in $\boundary(\Delta[n]\times \Delta[m])$. Furthermore, as verified
  in the last paragraph, we know that $\abpair$ is a face of the
  shuffle $\shuffle{w}\in H(\Delta[n]\times\Delta[m])$ and it is
  thus an element of $\boundary H(\Delta[n]\times\Delta[m])\subseteq
  H_i(\Delta[n]\times\Delta[m])$ as are its faces
  $\abpair\cdot\face^r_{k-1}$ and $\abpair\cdot\face^r_{k+1}$ as
  required.
  
  Otherwise, if $\abpair\cdot\face^r_k$ is not in the boundary
  $\boundary H(\Delta[n]\times\Delta[m])$ then, since it is an element
  of $H_{i-1}(\Delta[n]\times\Delta[m])$, there must exist an integer
  $u>i-1$ such that it is a face of the shuffle $\shuffle{u}$. Now, by
  the maximality of $w$ it follows that $w\geq u>i-1$ and therefore
  that $\shuffle{w}$ is in $H_{i-1}(\Delta[n]\times\Delta[m])$.
  However, we've already shown that $\abpair$ is a face of the shuffle
  $\shuffle{w}$ so it follows that it is an element of
  $H_{i-1}(\Delta[n]\times\Delta[m])$ as postulated in point
  (\ref{H.analysis.1}) of the statement.
  
  Now consider the face $\abpair\cdot\face^r_{k-1}$ and observe that
  the inequality $\beta(k-1)<\beta(k)$ of display~\eqref{H.anal.ineq}
  implies that either:
  \begin{itemize}
  \item $\beta(k-1)\notin\im(\beta\circ\face^r_{k-1})$ and thus that
    $\abpair\cdot\face^r_{k-1}=\pair<\alpha\circ\face^r_{k-1};
    \beta\circ \face^r_{k-1}>$ is an element of the boundary
    $\boundary H(\Delta[n]\times\Delta[m])\subseteq
    H_i(\Delta[n]\times\Delta[m])$ as required by
    clause~(\ref{H.analysis.2}) of the statement, or
  \item $k-2\geq 0$ and $\beta(k-2)=\beta(k-1)$, in which case we let
    $l=k-1$ and discharge our remaining obligation in the last
    paragraph of this proof.
  \end{itemize}
  Arguing dually for the face $\abpair\cdot\face^r_{k+1}$ we see that
  the inequality $\alpha(k)<\alpha(k+1)$ of
  display~\eqref{H.anal.ineq} implies that either:
  \begin{itemize}
  \item $\alpha(k+1)\notin\im(\alpha\circ\face^r_{k+1})$ and thus that
    $\abpair\cdot\face^r_{k+1}=\pair<\alpha\circ\face^r_{k+1};
    \beta\circ\face^r_{k+1}>$ is an element of the boundary $\boundary
    H(\Delta[n]\times\Delta[m])\subseteq
    H_i(\Delta[n]\times\Delta[m])$ as required by
    clause~(\ref{H.analysis.2}) of the statement, or
  \item $k+2\leq r$ and $\alpha(k+1)=\alpha(k+2)$, in which case we
    let $l=k-1$ and discharge our remaining obligation in the last
    paragraph of this proof.
  \end{itemize}
  
  So in the cases which remain outstanding we have an $l\in[r]$ (with
  $0<l<r$) at which the (in)equalities
  \begin{equation}\label{H.anal.ineq.2}
    \alpha(l-1)<\alpha(l)=\alpha(l+1)\text{ and }
    \beta(l-1)=\beta(l)<\beta(l+1)
  \end{equation}
  hold and in order to establish clause~(\ref{H.analysis.2}) of the
  statement we need to prove that the face $\abpair\cdot\face^r_l$ is
  an element of $H_i(\Delta[n]\times\Delta[m])$. To that end let $v$
  be the largest integer for which $\abpair\cdot\face^r_l$ is a face
  of the shuffle $\shuffle{v}$. Notice that we've already shown that
  $\abpair$ is a face of the shuffle $\shuffle{w}$ so it follows that
  $\abpair\cdot\face^r_l$ is also a face of that shuffle and therefore
  we may apply the maximality of $v$ to show that $v\geq w$. Indeed we
  may show that $\topop{v}$ and $\topop{w}$ are actually distinct
  operators and therefore that the inequality of the last sentence is
  actually strict. To do this start by arguing just as in the first
  paragraph of this proof to show that
  $\topop{v}(\alpha(l-1))=\beta(l+1)$. To compare this with the value
  $\topop{w}(\alpha(l-1))$ we start by observing that
  $\topop{w}(\alpha(l-1))\leq\topop{w}( \alpha(l)-1)$, since the first
  inequality of display~\eqref{H.anal.ineq.2} implies that
  $\alpha(l-1)\leq\alpha(l)-1$ and $\topop{w}$ is order preserving,
  then we use the fact that $\abpair$ is a face of $\shuffle{w}$ to
  establish that the inequality $\topop{w}(\alpha(l)-1)\leq\beta(l)$
  also holds as in display~(\ref{shuff.face.ineq}) of
  lemma~\ref{shuffles}(\ref{prod.simp.5a}).  Combining these various
  (in)equalities with the second inequality of
  display~\eqref{H.anal.ineq.2} we see that
  \begin{equation*}
    \topop{w}(\alpha(l-1))\leq\topop{w}(\alpha(l)-1)\leq\beta(l)<
    \beta(l+1)=\topop{v}(\alpha(l-1)) 
  \end{equation*}
  and in particular that $\topop{w}$ and $\topop{v}$ are distinct as
  postulated. So we've succeeded in showing that $v>w$ and we already
  know that $w>i-1$ so it clearly follows that $v>i$ and therefore
  that $\shuffle{v}$ is an element of $H_i(\Delta[n]\times\Delta[m])$
  (by definition), as is its face $\abpair\cdot\face^r_l$ as required.
\end{proof}

\begin{lemma}[Major filling lemma]\label{major.fill} For fixed $n,m\in
  \mathbb{N}$ and each $i>1$ the inclusion map
  \begin{equation}\label{major.fill.inc.1}
    \xymatrix@R=1ex@C=6em{*!++{H_i(\Delta[n]\otimes\Delta[m])}
      \ar@<0.5ex>@{u(->}[r]^{\textstyle\subseteq_r} & 
      *!++{H_{i-1}(\Delta[n]\otimes\Delta[m])}}
  \end{equation}
  is an f-extension. It follows that for each $N,M\in\Simplex$ the
  horizontal inclusions in the following commutative square
  \begin{equation}\label{major.fill.sq}
    \let\labelstyle=\textstyle
    \xymatrix@R=6ex@C=6em{
      {\boundary H(N\pretens M)}\ar@{^{(}->}[r]^{\subseteq_r}
      \ar@{^{(}->}[d]_{\subseteq_e} &
      {H(N\pretens M)}
      \ar@{^{(}->}[d]^{\subseteq_e} \\
      {\boundary H(N\otimes M)}
      \ar@{^{(}->}[r]_{\subseteq_r} &
      {H(N\otimes M)}
    }
  \end{equation}
  are both f-extensions.
\end{lemma}

\begin{proof}
  Assuming, for the next few paragraphs, that the first part of the
  statement holds, it follows that the inclusion
  \begin{equation}\label{major.fill.inc.2}
    \let\labelstyle=\textstyle
    \xymatrix@R=1ex@C=6em{*!++{\boundary H(\Delta[n]\otimes\Delta[m])}
      \ar@<0.5ex>@{u(->}[r]^{\subseteq_r} & 
      *!++{H(\Delta[n]\otimes\Delta[m])}}
  \end{equation}
  is also an f-extension since we know, from notation~\ref{Hdelta},
  that we may it may be expressed as the composite of the sequence
  consisting of the inclusions of display~\eqref{major.fill.inc.1}.
  
  Notice also that in the proof of lemma~\ref{pretens.density} each
  simplex used to witness an extension of thinness is a mediator, so
  it follows by lemma~\ref{H.analysis}(\ref{H.analysis.1}) that this
  witness is in $H_i(N\otimes M)$ whenever the face to which it
  extends thinness is in there. Consequently we may easily adapt that
  proof, by intersecting each of the subsets $A^{(j)}\cap N\otimes M$
  constructed there with $H_i(N\otimes M)$, to demonstrate that each
  inclusion
  \begin{displaymath}
    \let\labelstyle=\textstyle
    \xymatrix@R=1ex@C=6em{*!++{H_i(\Delta[n]\pretens\Delta[m])}
      \ar@<0.5ex>@{u(->}[r]^{\subseteq_r} & 
      *!++{H_i(\Delta[n]\otimes\Delta[m])}}
  \end{displaymath}
  is a t-extension (and thus an f-extension). Taking the special cases
  $i=1, \#\pair<n;m>$ we see that both of the vertical maps in
  display~\eqref{major.fill.sq} are f-extensions and consequently,
  applying the composition and cancellation results of
  observation~\ref{l-inv}(\ref{l-inv.comp}), it follows that the upper
  horizontal of that square is an f-extension iff its lower horizontal
  is. In particular, it follows that we may infer that the upper
  horizontal map in the square of display~\eqref{major.fill.sq.2} below
  is an f-extension from the fact that the map in
  display~\eqref{major.fill.inc.2} is an f-extension.

  Furthermore, by the definition of $\pretens$, it is the case that
  any simplex which is thin in $N\pretens M$ and is not thin in
  $\Delta[n]\pretens\Delta[m]$ must be a crushed cylinder in
  $N\pretens M$.  However the only one of these which is not in
  $\boundary (N\pretens M)$ is
  $\pair<\partproj^{n,m}_1;\partproj^{n,m}_2>$ and this is not itself
  in $H(N\pretens M)$, so it follows that $H(N\pretens M)=\boundary
  H(N\pretens M) \cup H(\Delta[n]\pretens \Delta[m])$. Of course, we
  also know that $\boundary H(\Delta[n]\pretens\Delta[m])=\boundary
  H(N\pretens M) \cap H(\Delta[n]\pretens\Delta[m])$ so we get a
  pasting square
  \begin{equation}\label{major.fill.sq.2}
    \let\labelstyle=\textstyle
    \xymatrix@R=6ex@C=6em{
      {\boundary H(\Delta[n]\pretens\Delta[m])}
      \ar@{u(->}[r]^{\subseteq_r} \ar@{u(->}[d]_{\subseteq_e}& 
      {H(\Delta[n]\pretens\Delta[m])} \ar@{u(->}[d]^{\subseteq_e}\\
      {\boundary H(N\pretens M)}
      \ar@{u(->}[r]_{\subseteq_r} & {H(M\pretens M)}
    }
  \end{equation}
  as in observation~\ref{pasting.sq}(\ref{pasting.sq.b}). Using that
  result we may infer that the lower horizontal of this square is an
  f-extension, as postulated, from the fact we have already
  demonstrated that its upper horizontal is an
  f-extension. 
  
  It remains to prove that the inclusion in
  display~\eqref{major.fill.inc.1} an f-extension. To that end,
  start by observing that the shuffle $\shuffle{i}\in
  H_{i-1}(\Delta[n]\otimes\Delta[m])$ gives rise to a stratified
  inclusion $\inc\yoneda{\shuffle{i}}:\Delta[n+m]->
  H_{i-1}(\Delta[n]\otimes\Delta[m]).$, by Yoneda's lemma, which we
  may extend to a regular map on its entire coimage
  $\Delta[n+m]\subseteq_e P\defeq\coim_e( \yoneda{\shuffle{i}})$ as in
  notation~\ref{substrat.reg}. This provides us with a stratified
  isomorphism between $P$ and the regular stratified subset
  $\stratc{\{\shuffle{i}\}}\subseteq_r
  H_{i-1}(\Delta[n]\otimes\Delta[m])$ generated by the shuffle
  $\shuffle{i}$. Furthermore, by
  observation~\ref{pasting.sq}(\ref{pasting.sq.a}) we know that if we
  define $X\subseteq_r P$ to be the inverse image of the regular
  subset $H_i(\Delta[n]\otimes\Delta[m])\subseteq_r
  H_{i-1}(\Delta[n]\otimes\Delta[m])$ under $\yoneda{\shuffle{i}}$
  then we get a pasting square
  \begin{equation}\label{psqu.1}
    \let\labelstyle=\textstyle
    \xymatrix@R=3em@C=5em{
      {X}\ar@{u(->}[r]^{\subseteq_r}
      \ar@{u(->}[d]_{\yoneda{\shuffle{i}}} &
      {P}\ar@{u(->}[d]^{\yoneda{\shuffle{i}}} \\
      {H_i(\Delta[n]\otimes\Delta[m])}\ar@{u(->}[r]_{\subseteq_r} &
      \poexcursion {H_{i-1}(\Delta[n]\otimes\Delta[m])}}
  \end{equation}
  in $\Strat$. More explicitly, an $r$-simplex (operator)
  $\arrow\mu:[r]->[n+m].$ is thin in $P$ iff $\shuffle{i}\cdot\mu$ is
  thin in $\Delta[n]\otimes\Delta[m]$ and it is in the regular subset
  $X\subseteq_r P$ iff $\shuffle{i}\cdot\mu$ is in
  $H_i(\Delta[n]\otimes\Delta[m])$.
  
  \begin{figure}
  \begin{displaymath}
    \def\vertchar{\scriptscriptstyle\bullet}
    \xymatrix@!0@C=1.5em@R=1.5em{
       & *[o]{\vertchar} & *[o]{\vertchar} & *[o]{\vertchar} & *[o]{\vertchar} & 
      *[o]{\vertchar} & *[o]{\vertchar} \\ 
       & *[o]{\vertchar} 
       \save[]-<0.75em,0em>*{t}\restore
      & *[o]{\vertchar} & *[o]{\vertchar}
      & *[o]{\vertchar}\ar@{.}[lll]  
       \save[]+<0.6em,-0.6em>*{k}\restore
      & *[o]{\vertchar}\ar@{.}[ur] & *[o]{\vertchar} \\ 
      & *[o]{\vertchar} & *[o]{\vertchar} & *[o]{\vertchar} & *[o]{\vertchar} & 
      *[o]{\vertchar} & *[o]{\vertchar} \\ 
      & *[o]{\vertchar} & *[o]{\vertchar} & *[o]{\vertchar} & *[o]{\vertchar} & 
      *[o]{\vertchar} & *[o]{\vertchar} \\ 
      [m]\ar[uu] & *[o]{\vertchar} & *[o]{\vertchar} & *[o]{\vertchar} & *[o]{\vertchar} & 
      *[o]{\vertchar} & *[o]{\vertchar} \\ 
      & *[o]{\vertchar} & *[o]{\vertchar} & *[o]{\vertchar} & *[o]{\vertchar} &
      *[o]{\vertchar} & *[o]{\vertchar} \\ 
      & *[o]{\vertchar}\ar@{-}'[rrr]'[rrruuuuu][rrruuuuur] &
      *[o]{\vertchar} & *[o]{\vertchar} & *[o]{\vertchar} 
      \save[]-<0em, 0.6em>*{s}\restore
      & *[o]{\vertchar} & *[o]{\vertchar} \\
       & & & [n] \ar[rr] & & & 
      }
    \end{displaymath}
    \caption{The shuffle $\shuffle{i}$ in $\Delta[n]\otimes\Delta[m]$.}
    \label{shuffle.pic}
  \end{figure}
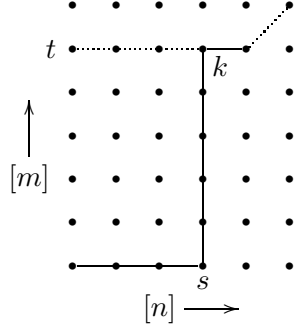
  
  It is our intention to apply corollary~\ref{minor.fill.cor} to prove
  that the upper horizontal inclusion in display~\eqref{psqu.1} is an
  f-extension, but to do so we must first select a suitable integer
  $0<k<n+m$. By the definition of the linear order $\triangleleft$, it
  is clear that the $\triangleleft$-minimal operator $\arrow
  \topop{1}:[n-1]->[m].$ identically $0$ (that is $\topop{1}(j)=0$ for
  all $j\in[n-1]$).  Since $i>1$, it follows that $\topop{i}$ is {\em
    not} identically $0$, in particular there is some $s\in[n-1]$ such
  that $\topop{i}(j)=0$ for $j<s$ and $\topop{i}(s)=t$ for some $t>0$.
  So letting $k=s+t$ and applying lemma~\ref{shuffles} we see that the
  corresponding shuffle $\shuffle{i}$ has
  \begin{equation}\label{ab.expl}
    \shufflel{i}(j) = \left\{
      \begin{array}{lp{1.2in}}
        j & for $j < s$, \\
        s & for $s\leq j\leq k$, \\
        s+1 & for $j = k+1$, \\
        ... & 
      \end{array}\right.
    \shuffler{i}(j) = \left\{
      \begin{array}{lp{1.2in}}
        0 & for $j< s$, \\
        j-s & for $s\leq j\leq k$, \\
        t & for $j = k+1$, \\
        ... & 
      \end{array}\right.
  \end{equation}
  (see figure~\ref{shuffle.pic}). Now we may check each of the
  conditions of corollary~\ref{minor.fill.cor} in turn:
  
  \vspace{1ex}
  \noindent{\bf corollary~\ref{minor.fill.cor}(\ref{mfc.adm})}
  Consulting display~\eqref{ab.expl} it is clear that $\shuffle{i}$
  satisfies the meditation condition at $k$, so applying
  lemma~\ref{mediator.admissible}(\ref{mediator.admissible.1}) we know
  that it is $k$-admissible in $\Delta[n]\otimes\Delta[m]$. So if
  $\arrow\mu:[r]->[n+m].$ is $k$-divided then the face $\shuffle{i}
  \cdot\mu$ is thin $\Delta[n]\otimes\Delta[m]$ and it follows that
  $\mu$ itself is thin in $P$ and consequently
  $\Delta^a_k[n+m]\subseteq_e P$ as required.

  \vspace{1ex}
  \noindent{\bf lemma~\ref{minor.fill}(\ref{minor.cond.i})}
  The shuffle $\shuffle{i}$ is a simplex of $H_{i-1}(\Delta[n]\otimes
  \Delta[m])$ and satisfies the mediation condition at $k$, so we may
  apply lemma~\ref{H.analysis}(\ref{H.analysis.2}) to show that the
  faces $\shuffle{i}\cdot\face^{n+m}_{k-1}$ and $\shuffle{i}\cdot
  \face^{n+m}_{k+1}$ are both elements of $H_i(\Delta[n]\otimes
  \Delta[m])$ and thus that the operators $\face^{n+m}_{k-1}$ and
  $\face^{n+m}_{k+1}$ are in the regular subset $X\subseteq_r P$ as
  required.

  \vspace{1ex}
  \noindent{\bf lemma~\ref{minor.fill}(\ref{minor.cond.ii})} Suppose
  that $\arrow\mu:[r]->[n+m].$ is a $k$-divided $r$-simplex in $P$ and
  $l$ is the unique integer with $\mu(l)=k$ then the face
  $\shuffle{i}\cdot\mu$ is a mediator simplex as witnessed by $l$ so
  we may apply lemma~\ref{H.analysis}(\ref{H.analysis.1}) to show that
  $\shuffle{i}\cdot\mu$ is in $H_i(\Delta[n]\otimes\Delta[m])$
  whenever $(\shuffle{i}\cdot\mu)\cdot\face^r_l=\shuffle{i}\cdot(\mu
  \circ\face^r_l)$ is in there and thus we see that $\mu$ is in the
  subset $X\subseteq_r P$ whenever $\mu\circ\face^r_l$ is in there as
  required.
    
  \vspace{1ex}
  \noindent{\bf corollary~\ref{minor.fill}(\ref{mfc.t-ext})} If $\mu$
  and $l$ are as in the last paragraph and $\mu\circ\face^r_l$ is thin
  in $P$ then the face $(\shuffle{i}\cdot\mu)\cdot\face^r_l=
  \shuffle{i}\cdot(\mu\circ\face^r_l)$ is thin in
  $\Delta[n]\otimes\Delta[m]$ and we may apply
  lemma~\ref{mediator.admissible}(\ref{mediator.admissible.2}) to the
  mediator simplex $\shuffle{i}\cdot\mu$ to show that
  $(\shuffle{i}\cdot\mu)\cdot\face^r_{l-1}=\shuffle{i}\cdot(
  \mu\circ\face^r_{l-1})$ and
  $(\shuffle{i}\cdot\mu)\cdot\face^r_{l+1}=\shuffle{i}\cdot(
  \mu\circ\face^r_{l+1})$ are both thin in there and
  thus that the simplices $\mu\circ\face^r_{l-1}$ and
  $\mu\circ\face^r_{l+1}$ are both thin in $P$ as required.
 
  \vspace{1ex} 
  
  So applying corollary~\ref{minor.fill.cor} as foreshadowed, we see
  that the upper horizontal inclusion in display~\eqref{psqu.1} is an
  f-extension and consequently, by applying
  observation~\ref{pasting.sq}(\ref{pasting.sq.a}), we find that the
  map of display~\eqref{major.fill.inc.1}, which appears as the lower
  horizontal in our pasting square~\eqref{psqu.1}, is an f-extension
  as required.
\end{proof}

\begin{lemma}\label{prim.f.tens} For each stratified  set
  $Y$ the horizontal arrows in the commutative square
  \begin{equation}\label{f.tens.dia}
    \let\labelstyle=\textstyle
    \xymatrix@R=5ex@C=10em{
      {\Lambda^a_k[n]\pretens Y}
      \ar@{u(->}[r]^{\subseteq_r\pretens Y} 
      \ar@{^{(}->}[d]_{\subseteq_r} &
      {\Delta^a_k[n]\pretens Y} 
      \ar@{^{(}->}[d]^{\subseteq_r} \\
      {\Lambda^a_k[n]\otimes Y}
      \ar@{u(->}[r]_{\subseteq_r\otimes Y} &
      {\Delta^a_k[n]\otimes Y}}
  \end{equation}
  are both f-extensions. Combining this result with
  lemma~\ref{prim.t.tens} and corollary~\ref{prim.d.tens} we see that
  for each complicial set $A$ the stratified set $A\Leftarrow Y$ is
  also complicial.
\end{lemma}

\begin{proof} First, it is worth commenting that the equivalence of
  the clauses in the statement follows by applying
  observation~\ref{partial.adj} parts~(\ref{partial.adj.a})
  and~(\ref{partial.adj.b}) to the partial right adjoint of
  theorem~\ref{otens.t-ext.pres}. Furthermore, since the vertical maps
  in display~\eqref{f.tens.dia} are both t-extensions (and thus
  f-extensions) by lemma~\ref{pretens.density} we may apply the
  composition and cancellation results of
  observation~\ref{l-inv}(\ref{l-inv.comp}) to show that its upper
  horizontal is a f-extension iff its lower horizontal is.
  
  Recall, from observation~\ref{deltat.dense}, that we may represent
  an arbitrary stratified set $Y$ as a weighted colimit
  $\colim(Y,\Delta)$ of standard (possibly thin) simplices.  We know
  that the functors $\Lambda^a_k[n]\pretens{-}$ and
  $\Delta^a_k[n]\pretens{-}$ preserve this colimit, by
  lemma~\ref{colim.pres}, from which it follows that the inclusion $\inc
  \subseteq_r\pretens Y: \Lambda^a_k[n]\pretens Y ->\Delta^a_k[n]
  \pretens Y.$ is isomorphic to the induced map $\inc
  \colim(Y,\subseteq_r\pretens\Delta(-)):
  \colim(Y,\Lambda^a_k[n]\pretens\Delta(-))->\colim(Y,\Delta^a_k[n]
  \pretens\Delta(-)).$. So applying
  observation~\ref{l-inv}(\ref{l-inv.colim}) we see that this latter
  map is an f-extension if $\subseteq_r\pretens\Delta[m]_?$ is an
  f-extension for each $[m]_?\in\tDelta$.
  
  In other words, it is sufficient for us to establish the result
  given in the statement of this lemma for each standard simplex
  $Y=\Delta[m]_?$, which we do by induction on the dimension $m$ of
  that simplex:
  
  \vspace{1ex}

  \noindent{\bf\boldmath Base case.}
  In this case the stratified set $Y$ in display~(\ref{f.tens.dia})
  is $\Delta[0]$ which is the identity for $\otimes$, thus there is
  nothing to prove.\vspace{1ex}

  \noindent{\bf\boldmath Inductive case.} Assume the induction
  hypothesis that the result given in the statement of the lemma holds
  for any standard (thin) simplex $Y=\Delta[p]_?$ of dimension $p<m$.
  In fact, arguing as above we know that this hypothesis immediately
  implies that the result of the statement holds for any
  $(m-1)$-skeletal stratified set $Y$, since any such may be expressed
  as a canonical colimit of standard simplices of dimension $<m$ as
  discussed in observation~\ref{(co)skeletal.strat}.
  
  \def\cylshuff{\pair<\partproj^{n,m}_1;\partproj^{n,m}_2>}
  \def\ycylshuff{i}

  In order to prove that the inclusion $\inc
  \subseteq_r\otimes\Delta[m]_?:\Lambda^a_k[n]\otimes\Delta[m]_?->
  \Delta^a_k[n]\otimes\Delta[m]_?.$ is also an f-extension, first note
  that we may apply observation~\ref{tens.thin.refl} to the regular
  inclusion $\Lambda^a_k[n]\subseteq_r\Delta^a_k[n]$ to show that
  $\Lambda^a_k[n]\otimes\Delta[m]_?$ is a regular subset of
  $\Delta^a_k[n]\otimes\Delta[m]_?$ for which
  $\subseteq_r\otimes\Delta[m]_?$ is the inclusion map. Furthermore it is
  clear that a simplex $\abpair$ is in this regular subset
  iff $\im(\alpha)\cup\{k\}\neq [n]$. Our approach to showing that
  this regular inclusion is an f-extension will be to decompose it
  into a sequence of regular subsets
  \begin{equation}\label{reg.seq}
    \Lambda^a_k[n]\otimes\Delta[m]_?
    \mathop{\subseteq_r}_{\text{\eqref{prim.f.tens.1}}}
    U\mathop{\subseteq_r}_{\text{\eqref{prim.f.tens.2}}}
    V\mathop{\subseteq_r}_{\text{\eqref{prim.f.tens.3}}}
    W\mathop{\subseteq_r}_{\text{\eqref{prim.f.tens.4}}}
    \Delta^a_k[n]\otimes\Delta[m]_?
  \end{equation}
  given by
  \begin{equation}\label{reg.seq.2}
    \begin{split}
      \abpair \in U \iff {} & 
      \im(\alpha)\cup\{k\}\neq [n]\orelse\im(\beta)\neq [m]\\
      \abpair \in V \iff {} & 
      \im(\alpha)\cup\{k\}\neq [n]\orelse
      \im(\beta)\neq [m]\orelse\\
      & \mkern30mu\left(k\notin \im(\alpha)\also
      \pair<n;0>\notin\im\abpair\right)\\
      \abpair \in W \iff {} &
      \im(\alpha)\cup\{k\}\neq [n]
      \orelse\im(\beta)\neq [m]
      \orelse\pair<n;0>\notin\im\abpair
    \end{split}
  \end{equation}
  and to show that each one of the inclusions in this decomposition is
  an f-extension, from which it then follows, by
  observation~\ref{l-inv}(\ref{l-inv.comp}), that their composite is
  also an f-extension as required. To complete our argument, the
  proofs in the following numbered paragraphs establish that the
  correspondingly numbered inclusions in display~\eqref{reg.seq} are
  f-extensions as stated:
  \begin{enumerate}
  \item\label{prim.f.tens.1}\vspace{1ex} By
    observation~\ref{tens.thin.refl}, the regular subset
    $\boundary\Delta[m]\subseteq_r\Delta[m]_?$ give rise to a regular
    subset $\Delta^a_k[n]\otimes\boundary\Delta[m]$ of
    $\Delta^a_k[n]\otimes\Delta[m]_?$ and it is easily demonstrated
    that
    \begin{equation*}
      \begin{split}
        \Lambda^a_k[n]\otimes\boundary\Delta[m] & {} =
        (\Lambda^a_k[n]\otimes\Delta[m]_?)\cap
        (\Delta^a_k[n]\otimes\boundary \Delta[m])\also\\
        U & {} = (\Lambda^a_k[n]\otimes\Delta[m]_?)\cup
        (\Delta^a_k[n]\otimes\boundary \Delta[m])
      \end{split}
    \end{equation*}
    so we get a pasting square
    \begin{equation}\label{prim.f.tens.1.psq}
    \let\labelstyle=\textstyle
    \xymatrix@R=2em@C=4em{
      {\Lambda^a_k[n]\otimes\boundary \Delta[m]\mskip4mu}
      \ar@{u(->}[r]^{\subseteq_r}\ar@{u(->}[d]_{\subseteq_r} &
      {\Delta^a_k[n]\otimes\boundary \Delta[m]}
      \ar@{u(->}[d]^-{\subseteq_r} \\
      {\Lambda^a_k[n]\otimes\Delta[m]_?\mskip4mu}
      \ar@{u(->}[r]_-{\subseteq_r} & \poexcursion
      {U}}
    \end{equation}
    in $\Strat$. However, $\boundary\Delta[m]$ is clearly
    $(m-1)$-skeletal, since all of its simplices of dimension $>m-1$
    are degenerate and thus, as discussed above, our induction
    hypothesis implies that the upper horizontal map in this diagram
    is an f-extension. It follows, by applying
    observation~\ref{pasting.sq}(\ref{pasting.sq.b}) that the lower
    horizontal inclusion in this diagram is also an f-extension as
    required.
    
  \item\label{prim.f.tens.2}\vspace{1ex} The stratified map $\arrow
    \Delta(\face^n_k):\Delta[n-1] ->\Delta^a_k[n].$ is both regular
    and an inclusion, the former fact following from the observation
    that elements in the image of $\Delta(\face^n_k)$ do not have $k$
    as a vertex and thus, by the definition of $\Delta^a_k[n]$, are
    thin if and only if they are degenerate. It follows, by
    observation~\ref{tens.thin.refl}, that the same properties hold
    for the stratified map:
    \begin{equation*}
      \xymatrix@R=1ex@C=8em{
        {\Delta[n-1]\otimes\Delta[m]_?}
        \ar[r]^{\textstyle\Delta(\face^n_k)\otimes\Delta[m]_?} &
        {\Delta^a_k[n] \otimes\Delta[m]_?}} 
    \end{equation*}
    Notice that the image of the subset
    $H(\Delta[n-1]\otimes\Delta[m]_?)$ under this map is precisely the
    regular subset of $\Delta^a_k[n]\otimes\Delta[m]_?$ of those
    simplices $\abpair$ such that $k\notin\im(\alpha)$ and
    $\pair<n;0>\notin\im\abpair$. It is therefore clear
    that the union of that image and the regular subset $U$ is
    precisely the subset $V$ of
    display~\eqref{reg.seq.2}. Furthermore, it is also easily seen
    that the inverse image of $U\subseteq_r V$ under the restricted
    regular inclusion $\inc\Delta(\face^n_k)\otimes\Delta[m]_?:
    H(\Delta[n-1]\otimes\Delta[m]_?)->V.$ is simply the regular subset
    $\boundary H(\Delta[n-1]\otimes\Delta[m]_?)\subseteq_r
    H(\Delta[n-1]\otimes\Delta[m]_?)$. Consequently we have a pasting
    square 
    \begin{displaymath}
      \let\labelstyle=\textstyle
      \xymatrix@R=5ex@C=5em{
        {\boundary H(\Delta[n-1]\otimes\Delta[m]_?)\mkern4mu}
        \ar@{u(->}[r]^{\subseteq_r}
        \ar@{u(->}[d]_{\Delta(\face^n_k)\otimes\Delta[m]_?} &
        {H(\Delta[n-1]\otimes\Delta[m]_?)}
        \ar@{u(->}[d]^{\Delta(\face^n_k)\otimes\Delta[m]_?} \\
        {U\mkern4mu}\ar@{u(->}[r]_-{\subseteq_r} &
        \poexcursion {V} }
    \end{displaymath}
    in $\Strat$.  Now, by lemma~\ref{major.fill} we know that the
    upper horizontal map in this diagram is an f-extension and,
    consequently, can apply
    observation~\ref{pasting.sq}(\ref{pasting.sq.a}) to infer that the
    lower horizontal map is an f-extension as well.
  \item\label{prim.f.tens.3}\vspace{1ex} Consider the regular subsets
    $H(\Delta^a_k[n]\otimes\Delta[m]_?)$ and $V$ of
    $\Delta^a_k[n]\otimes\Delta[m]_?$ and observe that it is easily
    demonstrated, directly from the definitions given in
    display~\eqref{reg.seq.2}, that we have 
    \begin{equation*}
      \begin{split}
        \boundary H(\Delta^a_k[n]\otimes\Delta[m]_?) & {} = 
        V\cap H(\Delta^a_k[n]\otimes\Delta[m]_?) \also \\
        W & {} = V\cup H(\Delta^a_k[n]\otimes\Delta[m]_?)
      \end{split}
    \end{equation*}
    and thus that we have a pasting square
    \begin{displaymath}
      \let\labelstyle=\textstyle
      \xymatrix@R=5ex@C=5em{
        {\boundary H(\Delta^a_k[n]\otimes\Delta[m]_?)\mkern4mu}
        \ar@{u(->}[r]^{\subseteq_r}
        \ar@{u(->}[d]_{\subseteq_r} &
        {H(\Delta^a_k[n]\otimes\Delta[m]_?)}
        \ar@{u(->}[d]^{\subseteq_r} \\
        {V\mkern4mu}
        \ar@{u(->}[r]_-{\subseteq_r} &
        \poexcursion {W} }
    \end{displaymath}
    in $\Strat$. The upper horizontal map in this diagram is an
    f-extension by lemma~\ref{major.fill}, so we may apply
    observation~\ref{pasting.sq}(\ref{pasting.sq.b}) to infer that the
    lower horizontal inclusion is an f-extension as well. 
  \item\label{prim.f.tens.4}\vspace{1ex} From the description given in
    display~\eqref{reg.seq.2} it is clear that only two non-degenerate
    simplices of $\Delta^a_k[n]\otimes\Delta[m]_?$ are {\bf not} in
    its regular subset $W$, those being the maximal shuffle
    $\cylshuff$ and its $k^{\text{th}}$ face $\cylshuff
    \cdot\face^{n+m}_k$.  The shuffle $\cylshuff\in\Delta^a_k[n]
    \otimes\Delta[m]_?$ gives rise to a stratified inclusion
    $\inc\ycylshuff:\Delta[n+m]->\Delta^a_k[n]\otimes \Delta[m]_?.$,
    by Yoneda's lemma, which we may extend to a regular map on its
    entire coimage $\Delta[n+m]\subseteq_e P\defeq
    \coim_e(\ycylshuff)$ as in notation~\ref{substrat.reg}. This
    provides us with a stratified isomorphism between $P$ and the
    regular stratified subset $\stratc{\{\cylshuff\}}\subseteq_r
    \Delta^a_k[n]\otimes\Delta[m]_?$ generated by our shuffle.
    Furthermore, by observation~\ref{pasting.sq}(\ref{pasting.sq.a})
    we know that if we define $X\subseteq_r P$ to be the inverse image
    of the regular subset $W\subseteq_r \Delta^a_k[n]\otimes
    \Delta[m]_?$ under $\ycylshuff$ then we get a pasting square
  \begin{equation}\label{psqu.2}
    \let\labelstyle=\textstyle
    \xymatrix@R=3em@C=5em{
      {X}\ar@{u(->}[r]^{\subseteq_r}
      \ar@{u(->}[d]_{\ycylshuff} &
      {P}\ar@{u(->}[d]^{\ycylshuff} \\
      {W}\ar@{u(->}[r]_-{\subseteq_r} &
      \poexcursion {\Delta^a_k[n]\otimes\Delta[m]_?}}
  \end{equation}
  in $\Strat$. More explicitly, an $r$-simplex (operator)
  $\arrow\mu:[r]->[n+m].$ is thin in $P$ iff the face
  $\cylshuff\cdot\mu$ is thin in $\Delta^a_k[n]\otimes\Delta[m]_?$ and
  it is in the regular subset $X\subseteq_r P$ iff $\cylshuff\cdot\mu$
  is in $W$.
  
  Since $\cylshuff$ and its face $\cylshuff\cdot\face^{n+m}_k$ are the
  only simplices of $\Delta^a_k[n]\otimes\Delta[m]$ which are not in
  $W$, it follows that the regular subset $X\subseteq_r P$ has the
  horn $\Lambda_k[n+m]$ as its underlying simplicial set. Furthermore,
  arguing as in item~(\ref{prim.t.tens.admprf}) of the proof of
  lemma~\ref{prim.t.tens}, we know that the shuffle
  $\pair<\partproj^{n,m}_1;\partproj^{n,m}_2>$ is $k$-admissible in
  $\Delta^a_k[n]\otimes\Delta[m]_?$ so if $\arrow\mu:[r]->[n+m].$ is 
  a $k$-divided simplex (operator) then $\cylshuff\cdot\mu$ is thin in
  $\Delta^a_k[n]\otimes\Delta[m]_?$ and so $\mu$ is thin in $P$ and it
  follows that $\Delta^a_k[n+m]\subseteq_e P$.
  
  These observations immediately serve to establish the conditions of
  lemma~\ref{minor.fill} and condition~(\ref{mfc.adm}) of
  corollary~\ref{minor.fill.cor} for $X\subseteq_r P$. So to complete
  the verification of the remaining condition~(\ref{mfc.t-ext}) of
  that corollary all we need do is demonstrate that if $\face^{n+m}_k$
  is thin in $P$ then so are $\face^{n+m}_{k-1}$ and
  $\face^{n+m}_{k+1}$, or equivalently that if
  $\cylshuff\cdot\face^{n+m}_k$ is thin in
  $\Delta^a_k[n]\otimes\Delta[m]_?$ then so are
  $\cylshuff\cdot\face^{n+m}_{k-1}$ and
  $\cylshuff\cdot\face^{n+m}_{k+1}$. However the face
  $\cylshuff\cdot\face^{n+m}_k$ is equal to the cylinder
  $\pair<\face^n_k\circ\partproj^{n-1,m}_1;\partproj^{n-1,m}_2>$,
  since $k<n$, and this is thin in $\Delta^a_k[n]\otimes\Delta[m]_?$
  iff $\id_{[m]}$ is thin in $\Delta[m]_?$ (since $\face^n_k$ is not
  thin in $\Delta^a_k[n]$) which must therefore be equal to
  $\Delta[m]_t$. Of course, $\Th_{m-1}(\Delta[m])=\Delta[m]_t$ and
  $\Th_{n-1}(\Delta[n]) \subseteq_e\Delta^a_k[n]$ so we may apply
  lemma~\ref{tensor.Th} to show that
  $\Th_{n+m-2}(\Delta[n]\otimes\Delta[m])\subseteq_e
  \Th_{n-1}(\Delta[n])\otimes\Th_{m-1}(\Delta[m])\subseteq_e
  \Delta^a_k[n]\otimes\Delta[m]_t$ and thus that all
  $(n+m-1)$-dimensional simplices of
  $\Delta^a_k[n]\otimes\Delta[m]_t$, including the faces
  $\cylshuff\cdot\face^{n+m}_{k-1}$ and
  $\cylshuff\cdot\face^{n+m}_{k+1}$, are all thin as required.
  
  Finally applying corollary~\ref{minor.fill.cor} we may infer that
  the horizontal inclusion at the top of display~\eqref{psqu.2} is an
  f-extension so, by observation~\ref{pasting.sq}(\ref{pasting.sq.a}),
  it follows that the horizontal inclusion at the bottom of that
  pasting square is also an f-extension as required.\qedhere
 \end{enumerate} 
\end{proof}

To conclude this section we present an important result which
parallels theorem~\ref{biclosed.precomp}, but applies to the category
of complicial sets $\Comp$ rather than the larger category of
pre-complicial sets $\Precomp$.

\begin{thm}\label{otens.f-ext.pres}\label{biclosed.comp}
  For any complicial set $A$ and any stratified sets $X$
  and $Y$, the stratified sets $A\Leftarrow Y$ and $X\Rightarrow A$
  are both complicial.  It follows that if $\arrow f:Z-> Z'.$ is an
  f-extension then so are the maps
  \begin{displaymath}
    \xymatrix@R=1ex@C=10em{
      Z\otimes Y\ar[r]^{\textstyle f\otimes Y} &
      Z'\otimes Y \\
      X\otimes Z\ar[r]^{\textstyle X\otimes f} &
      X\otimes Z'}
  \end{displaymath}
  Furthermore, by Day's reflection theorem (theorem~\ref{day.refl}), we
  may reflect the monoidal structure on $\Strat$ to a monoidal
  bi-closed structure on $\Comp$ with tensor product $A\otimes_c
  B\defeq \refl_c(A\otimes B)$ (where $\refl_c$ is reflector associated with
  $\Precomp$) and left and right closures $A\Rightarrow{*}$ and
  ${*}\Leftarrow B$ respectively.
\end{thm}

\begin{proof} Almost identical to that of theorem~\ref{otens.t-ext.pres},
  except that we start by observing that the functor $\arrow
  {-}\otimes Y:\Strat->\Strat.$ carries the primitive
  f-extensions $\arrow d:\Delta[1]_t->\Delta[0].$ and $\inc
  \subseteq_r:\Lambda^a_k[n]->\Delta^a_k[n].$ to stratified maps which are
  also f-extensions by applying lemmas~\ref{prim.d.tens}
  and~\ref{prim.f.tens} respectively. The remainder of the argument
  presented there carries through unaltered, aside from the
  substitution of the words ``complicial'' and ``f-extension'' for
  ``pre-complicial'' and ``t-extension''.
\end{proof}

\subsection{Superstructures of Complicial Sets}

Of course, we are motivated to think of the category of complicial
sets $\Comp$ as being analogous to the category of $\inf$-categories
$\InfCat$. On doing so, the functor
$\arrow\Sup_n(\cdot):\Strat->\Strat.$ takes on a particular importance
since it allows us to identify a subcategory of $\Comp$ which
corresponds to the full subcategory $\nCat$ of $n$-categories in
$\InfCat$. In pursuing this line of inquiry we'll need the following
lemma:

\begin{lemma}\label{comp.superst} 
  If $A$ is a complicial set then so is its $n$-dimensional
  superstructure $\Sup_n(A)$. Equivalently, if the stratified map
  $\arrow f:X ->Y.$ is an f-extension then so is the stratified map
  $\arrow \Th_n(f): \Th_n(X)->\Th_n(Y).$.
\end{lemma}

\begin{proof} Extending the argument used to prove
  lemma~\ref{precomp.superst}, all we need do is prove that each map
  obtained by applying $\Th_n$ to a primitive f-extension is again an
  f-extension. We consider each of the classes of primitive
  f-extensions in turn:

  \vspace{2ex}
  \noindent{\boldmath $\arrow d:\Delta[1]_t->\Delta[0].$} All of the
  simplices in $\Delta[0]$ and $\Delta[1]_t$ of dimension $>0$ are
  already thin, therefore $\Th_n(\Delta[0])=\Delta[0]$ and
  $\Th_n(\Delta[1]_t)=\Delta[1]_t$ for any $n\in\mathbb{N}$. It
  follows that the functor $\Th_n$ maps $\arrow
  d:\Delta[1]_t->\Delta[0].$ to itself.

  \vspace{2ex}
  \noindent{\boldmath $\overinc \subseteq_r:
    \Lambda^a_k[m]->\Delta^a_k[m].$} Again, we have two cases:
  \begin{itemize}
  \item {\boldmath $n\geq m-1$} In this case, all of the simplices of
    $\Lambda^a_k[m]$ and $\Delta^a_k[m]$ of dimension $>n$ are already
    thin, therefore $\Th_n(\Lambda^a_k[m])=\Lambda^a_k[m]$,
    $\Th_n(\Delta^a_k[m]) =\Delta^a_k[m]$ and, consequently, $\Th_n$
    maps $\overinc \subseteq_r: \Lambda^a_k[m]->\Delta^a_k[m].$ to itself.
  \item {\boldmath $n<m-1$} This is a simple matter of applying
    corollary~\ref{minor.fill.cor} to the regular subset
    $\Th_n(\Lambda^a_k[m])\subseteq_r\Th_n(\Delta^a_k[m])$. To do so
    observe that the conditions of lemma~\ref{minor.fill} follow
    trivially for any stratified subset whose underlying simplicial set
    is a $k$-horn and that condition~(\ref{mfc.adm}) of
    corollary~\ref{minor.fill.cor} simply asks for the standard
    inclusion $\Delta^a_k[m]\subseteq_e\Th_n(\Delta^a_k[m])$. Finally,
    the only simplex of $\Th_n(\Delta^a_k[m])$ which satisfies the
    postulates of corollary~\ref{minor.fill.cor}
    condition~(\ref{mfc.t-ext}) is the $m$-simplex $\id_{[m]}$ of
    which all $(m-1)$-dimensional faces, including $\face^m_{k-1}$,
    $\face^m_k$ and $\face^m_{k+1}$, are thin as required.
    \end{itemize}
    In either case the resulting map is an f-extension as required.
\end{proof}

\begin{defn}\label{n-comp.def} 
  We define an {\em $n$-complicial\/} set to be a complicial set $A$
  which is $n$-trivial and let $\Comp_n=\Comp\cap\Strat_n$ denote the
  full subcategory of $\Comp$ whose objects are these $n$-complicial
  sets.

  Applying lemma~\ref{precomp.triv} and arguing
  just as in observation~\ref{precomp_n.mon.bic}, we know that
  $X\Rightarrow A$ and $A\Leftarrow Y$ are $n$-complicial whenever $A$
  is and that, just as in theorem~\ref{biclosed.comp}, we may reflect
  the monoidal structure of $\Strat$ to a monoidal biclosed
  structure on the reflective full subcategory $\Comp_n$.
\end{defn}

\begin{lemma}\label{0.comp.eq.set} If $A$ is a complicial set then its
  $n$-skeleton $\Sk_0(A)$ and 0-superstructure $\Sup_0(A)$ coincide.
  Furthermore, any 0-skeletal stratified set is 0-complicial and,
  consequently, the category $\Comp_0$ is actually the full
  subcategory of 0-skeletal stratified sets in $\Strat$ (and is
  therefore equivalent to $\Set$).
\end{lemma}

\begin{proof}
  Any simplex of dimension $>0$ in the 0-skeleton $\Sk_0(A)$ is
  degenerate and thus thin in $A$, so it follows that
  $\Sk_0(A)\subseteq_r\Sup_0(A)$. Conversely, if $a$ is a simplex in
  $\Sup_0(A)$ of dimension $>0$ then it is pre-degenerate at $0$ (see
  definition~\ref{defn.predegen}), since every one of its faces of
  dimension $>0$ must be thin for it to be in $\Sup_0(A)$, and it
  follows, by lemma~\ref{lemma.predegen}, that $a$ is actually
  degenerate at $0$.  In other words, we see (by definition) that
  $\Sup_0(A)$ is 0-skeletal and thus it must be a regular subset of
  $\Sk_0(A)$, so the first of our results follows.
  
  From observation~\ref{conn.cpt} we know that all of the simplices
  and horns which constitute the domains and codomains of the
  primitive f-extensions are connected and therefore that each of
  these is bijective on components. However, by the same observation,
  a stratified set is 0-skeletal iff it is orthogonal to every
  stratified map which is bijective on components and it follows that
  every 0-skeletal stratified set is orthogonal to every primitive
  f-extension. In other words, every 0-skeletal stratified set is
  complicial furthermore, by definition, any simplex of a 0-skeletal
  set of dimension $>0$ is degenerate and thus thin, thereby
  demonstrating that any such stratified set is actually 0-complicial.
  Therefore, we may conclude that $\Comp_0$ coincides with the full
  subcategory of 0-skeletal sets in $\Strat$ and thus (from
  observation~\ref{conn.cpt}) that the functor
  $\arrow\dis:\Set->\Strat.$ provides us with the required equivalence
  between $\Set$ and $\Comp_0$.
\end{proof}



\section{The Path Category Construction}\label{comp.pathcat.sec}

\subsection{The Complicial Category of Prisms}

\begin{notation}[complicial categories] 
  We will often use the terms {\em complicial category\/} and {\em
    complicial functor\/} to refer to the objects and arrows of
  $\Cat(\Comp)$ (respectively). An {\em $n$-arrow\/} (resp. {\em
    $n$-object}) of a complicial category $\mathbb{C}$ is simply
  defined to be an $n$-simplex in its complicial set of arrows $C$
  (resp.\ objects $\scobj{C}$).
  
  When manipulating complicial categories we will, in general, assume
  that all standard categorical notions and constructions have been
  internalised to $\Comp$ in the canonical way. When referring to
  these we may not stress the use of the ``complicial'' qualifier if
  it is clear that it may be inferred from the context.
 
  We will also tend to freely apply qualifiers from our theory of
  stratified and complicial sets to complicial categories where the
  meaning is clear. For instance, if we say that $\mathbb{D}$ is a
  regular subcategory of the complicial category $\mathbb{C}$, it
  should be clear that we mean that the former is a complicial
  subcategory of the latter for which $\scarr{D}$ is, in fact, a
  regular subset of $\scarr{C}$.
\end{notation}

\begin{lemma}[the category of prisms in a complicial set]
  \label{defn.pathcat}
  If $\mathbb{T}_{\Cat}$ is the LE-theory of categories discussed in
  observation~\ref{le-th.cat} then the functor $\arrow
  \Cpathcat=\Th_1\circ\Delta:\Delta->\Strat.$ is a finitely presented
  f-almost $\mathbb{T}_{\Cat}$-coalgebra. Consequently we may apply
  the internal version of Kan's construction,
  observation~\ref{internal.kan}, to build a right
  adjoint, finitely accessible functor
  \begin{equation*}
    \let\labelstyle=\textstyle
    \xymatrix@R=2ex@C=10em{
      {\Comp} \ar[r]^<>(0.5){\pathcat\defeq\ihoml{\Cpathcat}} &
      {\Cat(\Comp)\hbox to 0pt{$\mathsurround=0pt 
        {}=\Alg<\mathbb{T}_{\Cat};\Comp>$\hss}}
    }
  \end{equation*}
  which carries a complicial set $A$ to a complicial category
  $\pathcat(A)$ which we call its {\bf\em category of prisms}.
\end{lemma}

\begin{proof} 
  Following the argument given in observation~\ref{conc.tensor.coalg},
  we may apply Yoneda's lemma to show that the stratified map
  $\arrow\colim(\subseteq_s,\Delta):\colim(\Lambda_k[n],\Delta)->
  \colim(\Delta[n],\Delta).$ is actually isomorphic the inclusion
  $\overinc \subseteq_r:\Lambda_k[n]->\Delta[n].$ under the usual
  identification of $\Simp$ as the category of minimally stratified
  sets in $\Strat$.  Furthermore, we know that $\Th_1$ is a left
  adjoint functor and thus that it preserves colimits, so it follows
  we have a family of isomorphisms $\colim(X,\Th_1\circ\Delta) \cong
  \Th_1(\colim(X,\Delta))\cong \Th_1(X)$ which is natural in
  $X\in\Simp$, and thus that
  $\arrow\colim(\subseteq_s,\Th_1\circ\Delta):\colim(\Lambda_k[n],
  \Th_1\circ\Delta)->\colim(\Delta[n],\Th_1\circ\Delta).$ is
  isomorphic to the inclusion $\overinc\subseteq_r:\Th_1(\Lambda_k[n])
  ->\Th_1(\Delta[n]).$. Applying this observation to the colimit
  characterisation of L-almost coalgebras from
  observation~\ref{func.from.coalg}, it follows that
  $\Th_1\circ\Delta$ is an f-almost $\mathbb{T}_{\Cat}$-coalgebra if
  and only if each of the latter inclusions is an f-extension.
  However, we know that $\Th_1(\Delta[n])=\Th_1(\Delta^a_k[n])$ and
  $\Th_1(\Lambda_k[n])= \Th_1(\Lambda^a_k[n])$ so these inclusions may
  equally well be obtained by applying $\Th_1$ to the primitive
  f-extensions of definition~\ref{prim.f-ext} and thus, applying
  lemma~\ref{comp.superst}, we see that they too are f-extensions as
  required.
\end{proof}

\begin{obs}[handedness conventions]
  \label{handedness.conv}
  Notice that we could equally well have chosen to define $\pathcat$
  to be the functor $\arrow\ihomr{\Cpathcat}:\Comp->\Cat(\Comp).$
  obtained using the right handed, rather than left handed, version of
  Kan's construction (as discussed in observation~\ref{internal.kan}).
  However, we have chosen to use the left handed construction in this
  section because it is compatible with subsequent arguments involving
  the decalage construction on $\Comp$.
  
  Where necessary we might differentiate these two possibilities using
  superscripts $\pathcat\lhv\defeq \ihoml{\Cpathcat}$ and
  $\pathcat\rhv\defeq\ihomr{\Cpathcat}$ and adopt this superscript
  convention to indicate the handedness convention for the other
  category and double category constructions introduced in this
  section. However, for the remainder of this section and the next we
  elide these superscripts and assume that all constructions will be
  made using the left closure. In the final section we shall be a
  little more careful, since there it turns out that we will need to
  deploy the right handed dual.
  
  Of course every result proved here for these left handed
  constructions also applies, dually, for the right handed version.
  Indeed, the dual construction on stratified simplicial sets allows
  us to directly relate these constructions. For instance we may show
  that if $A$ is a complicial set then $\pathcat\rhv(A^\circ)$ is
  actually canonically isomorphic to the complicial category obtained
  by applying the duality functor $\arrow ({-})^\circ:\Comp->\Comp.$
  point-wise to $\pathcat\lhv(A)$ and then applying the usual
  categorical dual.
\end{obs}

\begin{obs}[prism categories under the microscope]\label{pathcat.micro} 
  Observation~\ref{conc.tensor.coalg} demonstrates that the prism
  category functor bears a more explicit description as the external
  Kan functor $\homout{\Cpathcat\otimes \Delta}$ from $\Comp$ to
  $\Al{\mathbb{T}_{\Cat}\otimes\mathbb{T}_{\Comp}} \cong\Cat( \Comp)$.
  Combining this with the explicit description of the relationship
  between $\mathbb{T}_{\Cat}$-algebras and categories given in
  observation~\ref{le-th.cat} we gain the following explicit
  information about the structure of $\pathcat(A)$ for a complicial
  set $A$:
  \begin{enumerate}
  \item\label{pathcat.micro.1} The complicial set of arrows of
    $\pathcat(A)$ is $\ihoml{\Cpathcat}(A)_1 = \Delta[1]
    \Rightarrow A$, alternatively using the explicit presentation
    $\homout{\Cpathcat\otimes\Delta}(A)$ we see that the $n$-arrows
    of $\pathcat(A)$ may be identified with stratified maps of the
    form:
    \begin{displaymath}
      \let\labelstyle=\textstyle
      \xymatrix@R=2ex@C=10em{
        {\Delta[1]\otimes\Delta[n]}\ar[r]^p & {A} }
    \end{displaymath}
    Furthermore the action of a simplicial operator
    $\arrow\alpha:[m]->[n].$ on an $n$-arrow expressed in this form is
    given by $p\cdot\alpha=p\circ(\Delta[1]\otimes\Delta(\alpha))$.
  \item\label{pathcat.micro.2} Such an $n$-arrow in $\pathcat(A)$
    is thin if and only if it extends to a stratified map:
    \begin{displaymath}
      \let\labelstyle=\textstyle
      \xymatrix@R=2ex@C=10em{
        {\Delta[1]\otimes\Delta[n]_t}\ar[r]^p & {A} }
    \end{displaymath} 
    but we know, by observation~\ref{tensor.precomp.test}, that to
    demonstrate this all we need to do is show that $p$ maps the
    mediator simplices and crushed cylinders of
    $\Delta[1]\otimes\Delta[n]_t$ to thin simplices in $A$.
    However, we already know that this is true for all mediator
    simplices, since these are thin in $\Delta[1]\otimes\Delta[n]$,
    and it is an easy matter to demonstrate that the only
    non-degenerate crushed cylinders of
    $\Delta[1]\otimes\Delta[n]_t$ are the two $n$-simplices
    $\pair<\vertex^1_0\circ\eta^n;\id_{[n]}>$ and
    $\pair<\vertex^1_1\circ\eta^n;\id_{[n]}>$ and the $(n+1)$-simplex
    $\pair<\partproj^{\pair<1;n>}_1;\partproj^{\pair<1;n>}_2>$. It
    follows that $p$ is a thin $n$-arrow iff it maps these three
    simplices to thin simplices in $A$.
  \item\label{pathcat.micro.3} The complicial set of objects of
    $\pathcat(A)$ is (isomorphic to) $A$ itself and the following
    diagram displays its various source, target and identity maps:
    \begin{displaymath}
      \xymatrix@R=2ex@C=15em{
        *+[l]!<-1em,0ex>{A}\ar[r]|{\idt\defeq\Delta(\eta^1)\Rightarrow A} &
        *+[r]!<1em,0ex>{\Delta[1]\Rightarrow A}
        \ar@/_2.5ex/[l]_{\s\defeq\Delta(\vertex^1_0)\Rightarrow A}
        \ar@/^2.5ex/[l]^{\t\defeq\Delta(\vertex^1_1)\Rightarrow A}
      }
    \end{displaymath}
    More explicitly, if $\arrow p:\Delta[1]\otimes\Delta[n]->A.$ is an
    $n$-arrow of $\pathcat(A)$, presented as
    in~(\ref{pathcat.micro.1}), then we have $\s(p) \defeq
    p\pair<\vertex^1_0\circ\eta^n;\id_{[n]}>$ and $\t(p) \defeq
    p\pair<\vertex^1_1\circ\eta^n;\id_{[n]}>$ in $A$. In the other
    direction, if $a\in A_n$ then the identity $n$-arrow
    $\arrow\idt(a):\Delta[1]\otimes\Delta[n]->A.$ maps an $r$-simplex
    $\pair<\alpha;\beta>\in\Delta[1]\otimes\Delta[n]$ to
    $\idt(a)\pair<\alpha;\beta> \defeq a\cdot\beta$, 
    
    In order to explicitly describe the corresponding complicial
    subset of identities within $\Delta[1]\Rightarrow A$, notice that
    the stratified map $\arrow\Delta(\eta^1):\Delta[1]->\Delta[0].$
    extends to one with domain $\Delta[1]_t$ and it is this which we
    called $d$ in definition~\ref{prim.f-ext}(\ref{prim.f-ext.one}).
    Consequently, if we think of $\Delta[1]_t\Rightarrow A$ as a
    regular subset of $\Delta[1]\Rightarrow A$ (as in
    lemma~\ref{exp.inc}) then the stratified map $\idt$ from above
    factors as:
    \begin{displaymath}
      \let\labelstyle=\textstyle
      \xymatrix@R=2ex@C=8em{
        {A}
        \ar[r]^-{d\Rightarrow A} &
        {\Delta[1]_t\Rightarrow A}
        \ar@{u(->}[r]^{\subseteq_r} &
        {\Delta[1]\Rightarrow A} }
    \end{displaymath}
    However, $d$ is a primitive f-extension and so, by
    theorem~\ref{day.refl}(\ref{day.refl.2}), it follows that
    $d\Rightarrow A$ is a stratified isomorphism. In other words, the
    image of $\idt$, which we know to be the set of identities in
    $\pathcat(A)$, is precisely the regular subset
    $\Delta[1]_t\Rightarrow A\subseteq_r\Delta[1]\Rightarrow A$ and so an
    $n$-arrow $p\in\pathcat(A)$ is an identity if and only if it extends
    to a stratified map:
    \begin{displaymath}
      \let\labelstyle=\textstyle
      \xymatrix@R=2ex@C=10em{
        {\Delta[1]_t\otimes\Delta[n]}\ar[r]^p & {A} }
    \end{displaymath} 
  \item\label{pathcat.micro.4} A composable pair of $n$-arrows $\arrow
    p,q:\Delta[1]\otimes\Delta[n]->A.$ corresponds to a unique
    stratified map $\overarr:\Lambda_1[2]\otimes\Delta[n]->A.$ that we
    may uniquely extend along the f-extension $\overinc\subseteq_r:
    \Lambda_1[2]\otimes\Delta[n]->\Delta[2]_t\otimes\Delta[n].$ to a
    stratified map $\arrow w\pair<p;q>:\Delta[2]_t\otimes
    \Delta[n]->A.$ and which is thus completely determined by the fact
    that $w\pair<p;q>\circ(\Delta(\face^2_2)\otimes\Delta[n])=p$ and
    $w\pair<p;q>\circ (\Delta(\face^2_0)\otimes\Delta[n])=q$. Now the
    composite of these $n$-arrows is given by $q\comp p \defeq
    w\pair<p;q>\circ (\Delta(\face^2_1)\otimes\Delta[n])$ and we say
    that this composite is {\em witnessed\/} by the stratified map
    $w\pair<p;q>$.
  \end{enumerate}
\end{obs}

\begin{obs}[complicially enriched categories]\label{comp.enr} 
  A {\em complicially enriched\/} category in $\Comp$ is a 
  complicial category $\mathbb{C}$ whose set of objects $\scobj{C}$ is a
  0-complicial set.
  
  Applying lemma~\ref{0.comp.eq.set}, this is equivalent to saying
  that $\scobj{C}$ is 0-skeletal and thus isomorphic to the discrete
  stratified set on its set of 0-simplices. Consequently, it is easily
  seen that a complicially enriched category gives rise to an enriched
  category, in the sense of Kelly's book~\cite{Kelly:1982:ECT}, whose
  set of objects is the set of 0-simplices of $\scobj{C}$. Conversely,
  if $\mathcal{D}$ is a category (with a small set of objects) which
  is enriched in the cartesian category $\triple<\Comp;\times;
  \Delta[0]>$ then the disjoint union $\bigsqcup_{d,d'\in\obj(
    \mathcal{D})}\mathcal{D}(d,d')$ in $\Comp$ inherits an internal
  category structure from $\mathcal{D}$ which makes it into the
  complicial set of arrows of a complicially enriched category in our
  sense.  Indeed, these constructions demonstrate that the
  (2-)categories of complicially enriched categories in our sense, on
  the one hand, and categories enriched in the monoidal category
  $\triple<\Comp;\times;\Delta[0]>$ in the sense of
  Kelly~\cite{Kelly:1982:ECT}, on the other, are equivalent.

  Notice that the inclusion of $\CompCat$ in $\Cat(\Comp)$ has a
  right adjoint
  \begin{displaymath}
    \let\labelstyle=\textstyle
    \xymatrix@R=5ex@C=7.5em{
      {\CompCat}\ar@{u(->}@/^2.5ex/[rr]^{\text{inclusion}} & {\bot} &  
      {\Cat(\Comp)}\ar@/^2.5ex/[ll]^{\efunc}}
  \end{displaymath}
  where the complicially enriched category $\efunc(\mathbb{C})$ can be
  constructed as the complicial full subcategory of $\mathbb{C}$
  whose complicial set of objects is $\Sup_0(\scobj{C})\subseteq_r
  \scobj{C}$. In other words, $\efunc(\mathbb{C})$ is the subcategory
  of $\mathbb{C}$ whose complicial set of arrows is the regular subset
  of $\mathbb{C}_a$ obtained by pulling the regular subset
  $\Sup_0(\scobj{C})\times\Sup_0(\scobj{C})\subseteq_r
  \scobj{C}\times\scobj{C}$ back along the stratified map $\arrow
  \pair<\s;\t>:\scarr{C}->\scobj{C}\times\scobj{C}.$. Furthermore, it
  is clear that if $\arrow f:\mathbb{C}->\mathbb{C}'.$ is a complicial
  functor then $\efunc(f)$ is constructed by restricting $f$ to a
  complicial functor $\arrow f: \efunc(\mathbb{C})-> \efunc(
  \mathbb{C}').$.
\end{obs}

\begin{defn}[the category of paths in a complicial set]
  \label{pathecat.defn} 
  Define $\arrow \pathecat:\Comp->\CompCat.$ to be the composite of
  the prism category functor $\arrow \pathcat:\Comp->\Cat(\Comp).$ and
  the coreflection functor $\arrow \efunc: \Cat(\Comp)->\CompCat.$. In other
  words, $\pathecat(A)$ is the complicial full subcategory of
  $\pathcat(A)$ whose set of objects is the superstructure
  $\Sup_0(A)\subseteq_r A$.
\end{defn}

\begin{obs}[an explicit description of the $n$-arrows of $\pathecat(A)$]
  \label{enrich.path.expl}
  By definition an $n$-arrow $p$ of $\pathcat(A)$ is in the complicial
  subcategory $\pathecat(A)$ if and only if $\s(p)$ and $\t(p)$ are
  both simplices in $\Sup_0(A)$. However we know, by
  observation~\ref{pathcat.micro}(\ref{pathcat.micro.3}), that for an
  explicitly presented $n$-arrow $\arrow p:\Delta[1]\otimes\Delta[n]
  ->A.$ we have $\s(p)=p\pair<\vertex^1_0\circ\eta^n;\id_{[n]}>$ and
  $\t(p)=p\pair<\vertex^1_1 \circ\eta^n;\id_{[n]}>$, so it follows
  that $p$ is in $\pathecat(A)$ iff all faces of these $n$-simplices
  with dimension $>0$ are thin in $A$.
\end{obs}

\begin{obs}[discrete path categories]\label{disc.pathcat}
  A trivial observation, which we shall use a couple of times in
  forthcoming arguments, is that if $A$ is 0-complicial then the path
  category $\pathcat(A)$ is complicially enriched (and consequently
  $\pathcat(A)=\pathecat(A)$) and discrete as a category (cf.\ 
  definition~\ref{2cat.qua.doub}). 
  
  The first of these observations is an immediate consequence of the
  fact that the complicial set of objects of $\pathcat(A)$ is $A$. The
  latter we prove by recalling that the primitive f-extension $\arrow
  d:\Delta[1]_t->\Delta[0].$ was constructed by extending the stratified
  map $\arrow\Delta(\eta^1):\Delta[1]->\Delta[0].$ through the 
  inclusion $\overinc\subseteq_e:\Delta[1]->\Delta[1]_t.$. In other
  words, $\Delta(\eta^1)$ is the composite of an f-extension and an
  entire inclusion which is orthogonal to every $0$-trivial stratified
  set, so it is itself orthogonal to every $0$-complicial
  set. However, since $A$ is $0$-complicial it follows, by applying
  the fact that $\Comp_0$ is closed in $\Precomp$ under right closure
  (cf.\ definition~\ref{n-comp.def}), that the identities map
  $\idt=A\Leftarrow\Delta(\eta^1)$ of $\pathcat(A)$ is an isomorphism
  and so this category is discrete as required.
\end{obs}

\subsection{Path Categories and Superstructures}

\begin{obs}[superstructures of complicially enriched
  categories]\label{super.compl.enrich} Observe that (for
  $m\in\mathbb{N}$) the complicial superstructure functor
  $\Sup_m(\cdot)$ of definition~\ref{superstr.defn} preserves finite
  limits, since it is right adjoint to $\Th_m$, so we know that it may
  be lifted to an endo-functor on $\Cat(\Comp)$ as discussed in
  definition~\ref{le-th.defn}.  Concretely, if $\mathbb{C}$ is a
  complicial category then its structural maps $\idt$, $\s$ and $\t$
  restrict to maps between $\Sup_m(\scobj{C})\subseteq_r \scobj{C}$
  and $\Sup_m(\scarr{C})\subseteq_r \scarr{C}$ additionally
  $\Sup_m(\scarr{C})$ is closed in $\scarr{C}$ under the action of
  $\comp$. In other words, $\Sup_m(\scarr{C})$ is the underlying
  complicial set of arrows of a (regular) complicial subcategory of
  $\mathbb{C}$. Furthermore, the action of this lifted functor on
  complicial functors is simply one of restricting them to these
  regular sub-categories.
  
  Of course, complicial sub-categories of complicially enriched
  categories are also complicially enriched and so our lifted
  superstructure functor restricts to an endo-functor on $\CompCat$. 
    
  At the risk of confusing our notation a little, we actually call the
  subcategory obtained in this way the {\em $(m+1)$-dimensional
    superstructure\/} of the complicially enriched category
  $\mathbb{C}$ and denote it using the notation $\Sup_{m+1}(\sc{C})$
  rather than following the conventions established for lifted
  functors in definition~\ref{le-th.defn}. While this ``dimension
  shift'' might cause a little consternation at first it is important
  since, as we shall see in a page or two, this subcategory really is
  the appropriate analogue of the $(m+1)$-dimensional superstructure
  construction as generalised to the realm of complicially enriched
  categories.
  
  Finally, it is convenient to round off this sequence of
  superstructure functors by defining another endo-functor
  $\Sup_0(\cdot)$ on $\CompCat$, letting $\Sup_0(\mathbb{C})$ be
  the discrete subcategory of all identities in $\mathbb{C}$. In
  other words, $\Sup_0(\mathbb{C})$ is the smallest subcategory of
  $\mathbb{C}$ which contains all of its objects.
  
  Taking successive superstructures of a complicially enriched
  category $\mathbb{C}$ it is immediate that we obtain a canonical
  filtered family $\Sup_0(\sc{C})\subseteq_r\Sup_1(\sc{C})
  \subseteq_r...\subseteq_r\Sup_n(\sc{C})\subseteq_r...$ of regular
  subcategories the union of which is $\sc{C}$ itself.
\end{obs}

\begin{obs}\label{super.subcats.expl}
  If $A$ is a complicial set then the path category
  $\pathecat(\Sup_m(A))$ is also (trivially) a complicial subcategory
  of $\pathecat(A)$. Of course, an $n$-arrow $\arrow
  p:\Delta[1]\otimes\Delta[n]->A.$ of $\pathecat(A)$ is in
  $\pathecat(\Sup_m(A))$ iff it factors through the inclusion
  $\Sup_m(A)\subseteq_r A$ and, in turn, we may apply the adjunction
  $\Th_m\dashv\Sup_m(\cdot)$ to show that this happens precisely when
  $p$ extends to a stratified map:
  \begin{displaymath}
    \xymatrix@C=8em@R=1ex{
      {\Th_m(\Delta[1]\otimes\Delta[n])}\ar[r]^<>(0.5){\textstyle p} &
      {A}
    }
  \end{displaymath}
  
  It is instructive to compare this to the closely related complicial
  subcategory $\Sup_m({\pathecat(A)})$ for $m>0$. By Yoneda's lemma,
  and the definition of superstructures of complicially enriched
  categories, an $n$-arrow in here corresponds to a stratified map
  \begin{displaymath}
    \xymatrix@C=8em@R=1ex{
      {\Delta[n]}\ar[r] & 
      {\Sup_{m-1}({\pathecat(A)_a})\subseteq_r
        \Sup_{m-1}(\Delta[1]\Rightarrow A)}
    }
  \end{displaymath}
  which we may transpose under the adjunction $\Th_m\dashv\Sup_m(\cdot)$ and
  then under the partial adjunction $\Delta[1]\otimes{-}\dashv
  {\Delta[1]\Rightarrow\ast}$ to show that it yields a stratified map:
  \begin{displaymath}
    \xymatrix@C=8em@R=1ex{
      {\Delta[1]\otimes\Th_{m-1}(\Delta[n])}\ar[r]^<>(0.5){\textstyle p} &
      {A}
    }
  \end{displaymath}      
  In other words, our $n$-arrow $\arrow p:\Delta[1]\otimes\Delta[n]->
  A.$ in $\pathecat(A)$ is an arrow in the
  subcategory $\Sup_m({\pathecat(A)})$ iff it extends to a
  stratified map with domain $\Delta[1]\otimes\Th_{m-1}(\Delta[n])$.
  
  Given these representations of the arrows of these two
  subcategories, the following lemma is now easily established:
\end{obs}

\begin{lemma}\label{sup.comp.e.cat} If $A$ is a complicial set then
  the complicial sub-categories $\pathecat(\Sup_m(A))$ and
  $\Sup_m({\pathecat(A)})$ of $\pathecat(A)$ are in fact identical
  for each $m\geq 0$.
\end{lemma}

\begin{proof} First let us dispose with the special case $m=0$, for
  which the definition of $\Sup_0(\cdot)$ given in
  observation~\ref{super.compl.enrich} was selected. By
  observation~\ref{disc.pathcat} we know that $\pathecat(\Sup_0(A))$
  is a discrete subcategory of $\pathecat(A)$ whose complicial set of
  objects is $\Sup_0(A)\subseteq_r A$. However, by definition we also
  know that $\Sup_0({\pathecat(A)})$ is the discrete subcategory of
  $\pathecat(A)$ with the same complicial set of objects and thus
  these two sub-categories are identical.
  
  So from now on we assume that $m>0$ and prove the stated equality as
  two inclusions: \vspace{1ex}
 
  \noindent{\boldmath $\Sup_m({\pathecat(A)})\subseteq_r 
    \pathecat(\Sup_m(A))$:} We know that any $r$-simplices of
  $\Delta[1]$ with $r>1$ are degenerate, and thus thin, so it follows
  that $\Delta[1]=\Th_1(\Delta[1])$ and applying
  lemma~\ref{tensor.Th} we see that:
  \begin{displaymath}
    \Th_m(\Delta[1]\otimes\Delta[n])\subseteq_e
    \Th_1(\Delta[1])\otimes\Th_{m-1}(\Delta[n])=
    \Delta[1]\otimes\Th_{m-1}(\Delta[n])
  \end{displaymath}
  Thus, from the explicit descriptions furnished by
  observation~\ref{super.subcats.expl}, we see that if an $n$-arrow
  $p$ of $\pathecat(A)$ is an element of
  $\Sup_m({\pathecat(A)})$ then it extends to a stratified map
  with domain $\Delta[1]\otimes\Th_{m-1}(\Delta[n])$ which may then be
  restricted to one with domain $\Th_m(\Delta[1]\otimes\Delta[n])$
  thus demonstrating that $p$ is also in
  $\pathecat(\Sup_m(A))$. 
  \vspace{1ex}

  \noindent{\boldmath $\pathecat(\Sup_m(A))\subseteq_r
    \Sup_m({\pathecat(A)})$:} By definition we know that the
  complicial set of arrows of $\Sup_m({\pathecat(A)})$ is simply the
  complicial $(m-1)$-superstructure $\Sup_{m-1}({\pathecat(A)_a})$ so
  it follows that the inclusion we wish to prove here holds if and
  only if the complicial set of arrows of $\pathecat(\Sup_m(A))$ is
  $(m-1)$-complicial.  In other words, we
  need to demonstrate that for each $n>m-1$ any $n$-arrow $p$ of
  $\pathecat(\Sup_m(A))$ is thin.  Consulting
  observation~\ref{pathcat.micro}(\ref{pathcat.micro.2}), we know that
  this will be the case if $\arrow p:\Delta[1]\otimes\Delta[n]->A.$
  maps the $n$-simplices $\pair<\vertex^1_0\circ\eta^n;\id_{[n]}>$ and
  $\pair<\vertex^1_1\circ\eta^n;\id_{[n]};>$ and the $(n+1)$-simplex
  $\pair<\partproj^{\pair<1;n>}_1;\partproj^{\pair<1;n>}_2>$ to thin
  simplices in $A$. However, since $p$ is in $\pathcat(\Sup_m(A))$ we
  know, from observation~\ref{super.subcats.expl}, that it extends to a
  stratified map $\arrow p:\Th_m(\Delta[1]\otimes\Delta[n])-> A.$ and
  since $n+1>m$ it follows that $p$ carries the $(n+1)$-simplex above
  to a thin simplex in $A$. Furthermore, since $p$ is in
  $\pathecat(A)$ we may apply observation~\ref{enrich.path.expl} to
  show that it maps both of the $n$-simplices above to thin simplices
  in $A$ as required.
\end{proof}
 
\subsection{A Complicial Double Category with Connections}

By iterating the category of prisms construction $\pathcat$ we
may construct an important complicial double category with connections.

\begin{obs}[double prism categories]\label{iter.path}
  Of course, for each stratified set $X$ the (partially) right adjoint
  functor $X\Rightarrow{*}$ preserves the (small) limits of $\Comp$.
  It follows, by construction and the point-wise nature of limits in
  $\Cat(\Comp)$, that the functor $\arrow\pathcat:
  \Comp->\Cat(\Comp).$ preserves all (small) limits.  Consequently we
  may lift $\pathcat$ to a LE-functor
  \begin{equation*}
    \let\labelstyle=\textstyle
    \xymatrix@R=2ex@C=10em{
      {\Cat(\Comp)} \ar[r]^<>(0.5){\pathcat} &
      {\Cat(\Cat(\Comp))\hbox to 0pt{$\mathsurround=0pt 
        {}\defeq\Double(\Comp)$\hss}}
    }
  \end{equation*}
  as discussed in definition~\ref{le-th.defn}. This constructs a
  vertically presented double category $\pathcat(\mathbb{C})$ by
  applying $\pathcat$ point-wise to the structural components of the
  complicial category $\mathbb{C}$.  More explicitly, the double
  category $\pathcat(\mathbb{C})$ has:
  \begin{itemize}
  \item complicial set of squares $\Delta[1]\Rightarrow\scarr{C}$, 
  \item complicial sets of vertical arrows $\scarr{C}$ and horizontal
    arrows $\Delta[1]\Rightarrow\scobj{C}$,
  \item complicial set of objects $\scobj{C}$,
  \item vertical categories $\mathbb{C}$ of arrows and
    $\Delta[1]\Rightarrow\mathbb{C}$ of squares,
  \item horizontal categories $\pathcat(\scobj{C})$ of arrows and
    $\pathcat(\scarr{C})$ of squares,
  \end{itemize}
  Composing the functors $\arrow\pathcat:\Comp->\Cat(\Comp).$ and
  $\arrow\pathcat:\Cat(\Comp)->\Double(\Comp).$ we get a functor
  $\arrow\pathsq:\Comp->\Double(\Comp).$ and call $\pathsq(A)$ the
  {\em double prism category\/} of the complicial set $A$.
\end{obs}

\begin{obs}[double prism categories under the microscope]
  \label{dpathcat.micro} 
  From observation~\ref{coalg.tensor} we see that we may present
  $\pathsq$ more symmetrically as $\arrow\ihoml{\Cpathcat\otimes
    \Cpathcat}:\Comp->\Al{\mathbb{T}_{\Cat}\otimes\mathbb{T}_{\Cat}}(
  \Comp).$ (whose codomain is isomorphic to $\Double(\Comp)$).
  Furthermore if we now apply observation~\ref{conc.tensor.coalg}, as
  we did in observation~\ref{pathcat.micro}, we see that this in turn
  bears an explicit description as $\arrow\homout{\Cpathcat\otimes
    \Cpathcat\otimes \Delta}: \Comp->\Al{\mathbb{T}_{\Cat}\otimes
    \mathbb{T}_{\Cat}\otimes \mathbb{T}_{\Comp}}.$
  the external Kan functor.  It follows that the complicial double
  category $\pathsq(A)$ associated with a complicial set $A$ bears the
  following explicit description:
  \begin{enumerate}
  \item\label{dpathcat.micro.1} The complicial set of squares of
    $\pathsq(A)$ is $\ihoml{\Cpathcat\otimes
      \Cpathcat}(A)\pair<[1];[1]> = (\Delta[1]\otimes
    \Delta[1])\Rightarrow A$. Alternatively, using the explicit
    external presentation $\homout{\Cpathcat\otimes
      \Cpathcat\otimes \Delta}(A)$ we see that the $n$-squares
    (that is the $n$-simplices in the complicial set of squares) of
    $\pathsq(A)$ may be described as stratified maps of the form:
    \begin{equation}\label{nsquare.disp}
      \let\labelstyle=\textstyle
      \xymatrix@R=2ex@C=10em{
        {\Delta[1]\otimes\Delta[1]\otimes\Delta[n]}\ar[r]^\lambda & {A} }
    \end{equation}
    Furthermore, the action of a simplicial operator
    $\arrow\alpha:[m]->[n].$ on an $n$-square expressed in this form
    is given by $\lambda\cdot\alpha=\lambda\circ(\Delta[1]\otimes
    \Delta[1] \otimes\Delta(\alpha))$.
  \item\label{dpathcat.micro.2} Consulting the internal presentation
    $\ihoml{\Cpathcat\otimes \Cpathcat}(A)$ of $\pathsq(A)$, it
    is clear that its horizontal and vertical categories of arrows of
    are both equal to $\pathcat(A)$ and that its complicial set of
    objects is $A$ itself. Furthermore, arguing just as in
    observation~\ref{pathcat.micro}(\ref{pathcat.micro.4}) we see that
    its sets of horizontal and vertical identity squares in are simply
    the regular subsets $(\Delta[1]\otimes\Delta[1]_t)\Rightarrow A$
    and $(\Delta[1]_t\otimes \Delta[1])\Rightarrow A$ of
    $(\Delta[1]\otimes\Delta[1]) \Rightarrow A$
    respectively.
    
    Applying the external presentation $\homout{\Cpathcat\otimes
      \Cpathcat\otimes \Delta}(A)$ we see that the source and
    target maps for the horizontal and vertical categories of squares,
    when applied to the $n$-square in display~(\ref{nsquare.disp}),
    may be given concretely by
    \begin{equation*}
      \begin{aligned}
        \s_h(\lambda) & {}=\lambda\circ(
        \Delta[1]\otimes\Delta(\vertex^1_0)\otimes\Delta[n]) & \mkern20mu
        \s_v(\lambda) & {}=\lambda\circ(
        \Delta(\vertex^1_0)\otimes\Delta[1]\otimes\Delta[n]) \\
        \t_h(\lambda) & {}=\lambda\circ(
        \Delta[1]\otimes\Delta(\vertex^1_1)\otimes\Delta[n]) &
        \t_v(\lambda) & {}=\lambda\circ(
        \Delta(\vertex^1_1)\otimes\Delta[1]\otimes\Delta[n]) \\
      \end{aligned}
    \end{equation*}
    and if $\arrow p:\Delta[1]\otimes\Delta[n]->A.$ is an $n$-arrow in
    $\pathcat(A)$ then the horizontal and vertical identity squares on
    this are given by:
    \begin{equation*}
      \idt_h(p)=p\circ(\Delta[1]\otimes\Delta(\eta^1)\otimes\Delta[n])
      \mkern40mu
      \idt_v(p)=p\circ(\Delta(\eta^1)\otimes\Delta[1]\otimes\Delta[n])
    \end{equation*}
    Furthermore the horizontal and vertical identity squares of
    $\pathsq(A)$ are those $n$-squares of the form shown in
    display~(\ref{nsquare.disp}) which extend to a stratified map with
    domain $\Delta[1]\otimes\Delta[1]_t\otimes \Delta[n]$ or
    $\Delta[1]_t\otimes\Delta[1]\otimes\Delta[n]$ respectively.
  \item\label{dpathcat.micro.3} As in
    observation~\ref{pathcat.micro}(\ref{pathcat.micro.4}) if
    $\lambda$ and $\lambda'$ are a horizontally composable pair of
    $n$-squares then we may combine them to give a single stratified
    map $\overarr: \Delta[1]\otimes\Lambda_1[2]\otimes\Delta[n]->A.$.
    and extend this along the f-extension $\overinc\subseteq_r:
    \Delta[1]\otimes\Lambda_1[2]\otimes\Delta[n]->\Delta[1]\otimes
    \Delta[2]_t\otimes\Delta[n].$ to form a {\em horizontal witness\/}
    $\arrow w_h\pair<\lambda;\lambda'>:\Delta[1]\otimes\Delta[2]_t
    \otimes\Delta[n]->A.$. This witness is uniquely determined by the
    equations $w_h\pair<\lambda;\lambda'>\circ(\Delta[1]\otimes\Delta(
    \face^2_2) \otimes\Delta[n])=\lambda$ and $w_h\pair<\lambda;
    \lambda'>\circ(\Delta[1]\otimes\Delta(\face^2_0) \otimes
    \Delta[n])=\lambda'$ and we use it to define the horizontal
    composite by $\lambda'\comp_h\lambda= w_h\pair<\lambda;\lambda'>
    \circ( \Delta[1]\otimes\Delta(\face^2_1) \otimes\Delta[n])$.
    Dually, if $\lambda$ and $\lambda'$ are vertically composable then
    we form a {\em vertical witness\/} $\arrow w_v\pair<\lambda;
    \lambda'>:\Delta[2]_t\otimes\Delta[1]\otimes\Delta[n]->A.$ which is
    uniquely determined by the equations $w_v\pair<\lambda;\lambda'>
    \circ(\Delta(\face^2_2)\otimes\Delta[1]\otimes\Delta[n])=\lambda$ and
    $w_v\pair<\lambda;\lambda'>\circ(\Delta(\face^2_0)\otimes\Delta[1]\otimes
    \Delta[n])=\lambda'$ and we use it to define the vertical
    composite by $\lambda'\comp_v\lambda= w_v\pair<\lambda;\lambda'>
    \circ(\Delta(\face^2_1)\otimes\Delta[1]\otimes\Delta[n])$.
  \end{enumerate}
\end{obs}

The following theorem makes a pivotal contribution to understanding the
structure of the double category $\pathsq(A)$:
 
\begin{thm}[$\pathsq(A)$ as a double category with connections]
  \label{pathsq.conn}
  If $A$ is a complicial set then the double prism category
  $\pathsq(A)$ possesses a canonical thinness structure $t\pathsq(A)$
  which extends $\pathsq$ to a functor
  $\overarr:\Comp->\Double_T(\Comp).$. Furthermore, each
  $\strat{\pathsq(A)}$ is a double category with connections in the
  sense of lemma~\ref{dcat.with.con}.
\end{thm}

\begin{proof}
  In order to construct the claimed thinness structure
  $t\pathsq(A)\subseteq\pathsq(A)$ (cf.\ 
  definition~\ref{strat.double.cat}), start by considering the
  f-almost $\mathbb{T}_{\Cat}\otimes \mathbb{T}_{\Cat}$-coalgebra
  $\Cpathcat\otimes \Cpathcat$ which we used to define
  $\pathsq$.  Notice that the functor $\arrow \Th_1:\Strat->\Strat.$
  preserves f-extensions, by lemma~\ref{comp.superst}, and preserves
  colimits, since it is left adjoint to $\Sup_1(\cdot)$, from which it
  follows that it preserves all f-almost colimits. Consequently, the
  composite $\arrow \Th_1 \circ (\Cpathcat\otimes
  \Cpathcat):\Delta\times\Delta ->\Strat.$ is also an f-almost
  $\mathbb{T}_{\Cat}\otimes \mathbb{T}_{\Cat}$-coalgebra, to which we
  may apply the construction of observation~\ref{func.from.coalg} in
  order to construct a second functor $t\pathsq\defeq
  \ihoml{\Th_1\circ(\Cpathcat\otimes \Cpathcat)}$ from $\Comp$
  to $\Double(\Comp)$.  
  
  Now, every stratified set $X$ is an entire subset of the associated
  set $\Th_1(X)$ so we may construct a canonical map
  $\overarr:\Cpathcat\otimes \Cpathcat->\Th_1 \circ (\Cpathcat\otimes
  \Cpathcat).$ of our f-almost coalgebras (that is to say an arrow of
  $\CoAl{\mathbb{T}_{\Cat}\otimes
    \mathbb{T}_{\Cat}}_{\refl_f}(\Strat)$) whose components are the
  entire inclusions $\overinc\subseteq_e:\Cpathcat([n])\otimes
  \Cpathcat([m]) -> \Th_1(\Cpathcat([n])\otimes \Cpathcat([m])).$.  As
  discussed in observation~\ref{func.from.coalg}, this map of f-almost
  coalgebras gives rise to a derived natural transformation
  $\arrow\ihoml{\subseteq_e}:t\pathsq->\pathsq.$ whose component at
  the complicial set $A$ is obtained by applying ${*}\Rightarrow A$ to
  the entire inclusions of the last sentence and thus, by
  lemma~\ref{exp.inc}, identifies $t\pathsq(A)$ as a regular
  sub-double category of $\pathsq(A)$.
  
  Indeed we have $\Th_1 (\Cpathcat([1])\otimes \Cpathcat([1]))
  =\Th_1(\Delta[1] \otimes\Delta[1])$ so it follows that the set of
  squares of $t\pathcat(A)$ is the regular subset
  $\Th_1(\Delta[1]\otimes\Delta[1])\Rightarrow A$ of the set of
  squares $(\Delta[1]\otimes\Delta[1])\Rightarrow A$ of $\pathcat(A)$.
  Also, by observation~\ref{tensor.Th}, we know that
  $\Th_1(\Delta[1]\otimes\Delta[1])$ is an entire subset of both
  $\Delta[1]_t\otimes\Delta[1]$ and $\Delta[1]\otimes\Delta[1]_t$ so
  it follows, from observation~\ref{dpathcat.micro}, that the sets of
  horizontal and vertical identity squares in $\pathsq(A)$ are
  actually regular subsets of the set of squares of $t\pathsq(A)$. In
  other words, the pair $\strat{\pathsq(A)}$ is indeed a double
  category with thinness as required.

  Finally, by construction the regular inclusions
  $\overinc\subseteq_r:t\pathsq(A)->\pathsq(A).$ are natural in $A$.
  In other words, if $\arrow f:A->A'.$ is a complicial map then
  the double functor $\arrow\pathsq(f):\pathsq(A)->\pathsq(A').$ lifts
  to a thinness preserving double functor
  $\arrow\pathsq(f):\strat{\pathsq(A)} ->\strat{\pathsq(A')}.$, thus
  extending $\pathsq$ to a functor from $\Comp$ to the category of
  double categories with thinness $\Double_T(\Comp)$ as stated.
  
  \vspace{1ex}
  
  To complete this proof, we must demonstrate that each double
  category with thinness $\strat{\pathsq(A)}$ also satisfies the
  conditions given in lemma~\ref{dcat.with.con}:\vspace{1ex}

  \noindent {\em Conditions~\ref{dcat.with.con}(\ref{con.one})
    and~\ref{dcat.with.con}(\ref{con.two}).\/} These two conditions
  have dual proofs, so we simply prove the latter regarding vertical
  arrows. Consider the following commutative square of inclusion maps
  \begin{displaymath}
    \xymatrix@R=4ex@C=12em{
      {\Delta[1] + \Delta[1]}
      \ar@{u(->}[r]^{\langle\Delta(\vertex^1_0)\otimes\Delta[1],
                     \Delta(\vertex^1_1)\otimes\Delta[1]\rangle}
      \ar@{u(->}[d]_{\subseteq_e} &
      {\Th_1(\Delta[1]\otimes\Delta[1])}
      \ar@{u(->}[d]^{\subseteq_e} \\
      {\Delta[1]_t + \Delta[1]_t}
      \ar@{u(->}[r]_{\langle\Delta(\vertex^1_0)\otimes\Delta[1],
                     \Delta(\vertex^1_1)\otimes\Delta[1]\rangle} &
      {\Delta[1]\otimes\Delta[1]_t}
    }
  \end{displaymath}
  in which the right hand vertical is an entire subset inclusion
  obtained by applying lemma~\ref{tensor.Th} (with $n=1$ and $m=0$ we
  get $\Th_1(\Delta[1]\otimes \Delta[1]) \subseteq_e
  \Th_1(\Delta[1])\otimes\Th_0(\Delta[1]) =
  \Delta[1]\otimes\Delta[1]_t$) and the horizontal maps are those
  induced by the indicated stratified maps and the universal
  properties of the sums on the left hand side.  Notice that the only
  simplices which are thin in $\Delta[1]\otimes\Delta[1]_t$ but not in
  $\Th_1(\Delta[1]\otimes\Delta[1])$ are the 1-simplices
  $\pair<\vertex^1_0\circ\eta^1;\id_{[1]}>$ and
  $\pair<\vertex^1_1\circ\eta^1;\id_{[1]}>$, which are precisely the
  images of the 1-simplex $\id_{[1]}\in\Delta[1]$ under the maps
  $\Delta(\vertex^1_0)\otimes\Delta[1]$ and
  $\Delta(\vertex^1_1)\otimes\Delta[1]$ respectively and it follows,
  by observation~\ref{pasting.sq}(\ref{pasting.sq.a}), that our square
  is a pasting square. Of course, if we apply the contravariant
  functor ${*}\Rightarrow A$ to this square it will carry colimits in
  the diagram to limits in $\Comp$, giving us a pullback square
  \begin{displaymath}
    \xymatrix@R=4ex@C=8em{
      *+[l]!<-2em,0ex>{(\Delta[1]\otimes\Delta[1]_t)\Rightarrow A}\pbexcursion
      \ar@{u(->}[r]^{\subseteq_r}\ar[d]_{\pair<\s_v;\t_v>} &
      {\Th_1(\Delta[1]\otimes\Delta[1])}\Rightarrow A
      \ar[d]^{\pair<\s_v;\t_v>} \\
      {(\Delta[1]_t\Rightarrow A)\times(\Delta[1]_t\Rightarrow A)}
      \ar@{u(->}[r]_{\subseteq_r} &
      {(\Delta[1]\Rightarrow A)\times(\Delta[1]\Rightarrow A)}
    }
  \end{displaymath}
  In other words, the regular subset $(\Delta[1]\otimes
  \Delta[1]_t)\Rightarrow A$ of horizontal identity squares in
  $\pathsq(A)$ is precisely the subset of $t\pathsq(A)$ of those
  squares $\lambda$ for which $\s_v(\lambda)$ and $\t_v(\lambda)$ are
  in the complicial subset of identities $\Delta[1]_t\Rightarrow A$ of
  the category $\pathcat(A)$ of horizontal arrows, the desired result
  follows.  \vspace{0.5ex}
 
  \noindent {\em Conditions~\ref{dcat.with.con}(\ref{con.three})
    and~\ref{dcat.with.con}(\ref{con.four}).\/} These two conditions
  have dual proofs, so we simply prove the latter. It is easily
  demonstrated that we have a stratified map
  \begin{displaymath}
    \let\labelstyle=\textstyle
    \xymatrix@R=1ex@C=10em{
      {\Th_1(\Delta[1]\otimes\Delta[1])}
      \ar[r]^-{\max} &
      {\Delta[1]} 
    }
  \end{displaymath}
  where $\max\pair<\alpha;\beta>$ is the point-wise
  maximum of the simplicial operators $\arrow\alpha,\beta:[r]->[1].$,
  and that this satisfies the following identities:
  \begin{displaymath}
    \begin{array}{r@{}lcr@{}l}
      \max\circ(\Delta[1]\otimes\Delta(\vertex^1_0)) & {} =
      \id_{\Delta[1]} & \mkern 30mu &
      \max\circ(\Delta(\vertex^1_0)\otimes\Delta[1]) & {} =
      \id_{\Delta[1]} \\
      \max\circ(\Delta[1]\otimes\Delta(\vertex^1_1)) & {} = 
      \Delta(\vertex^1_1\circ\eta^1) &&
      \max\circ(\Delta(\vertex^1_1)\otimes\Delta[1]) & {} =
      \Delta(\vertex^1_1\circ\eta^1) 
    \end{array}
  \end{displaymath}
  Applying the contravariant functor ${*}\Rightarrow A$ we get a
  stratified map $\max\Rightarrow A$ whose domain $\Delta[1]\Rightarrow A$
  is the set of vertical arrows of $\pathsq(A)$ and whose codomain is
  $t\pathsq(A)$. Furthermore, applying ${*}\Rightarrow A$ to the
  identities above, and using the definitions of the various source
  and target maps of $\pathcat(A)$ and $\pathsq(A)$, we get the
  identities:
  \begin{displaymath}
    \begin{array}{r@{}lcr@{}l}
      \s_h\circ(\max\Rightarrow A) & {} = \id_{\Delta[1]\Rightarrow A} & 
      \mkern 30mu &
      \s_v\circ(\max\Rightarrow A) & {} = \id_{\Delta[1]\Rightarrow A} \\
      \t_h\circ(\max\Rightarrow A) & {} = \idt\circ\t  &&
      \t_v\circ(\max\Rightarrow A) & {} = \idt\circ\t 
    \end{array}
  \end{displaymath}
  In other words if $q$ is a vertical arrow of $\pathsq(A)$ then the
  thin square $\nu_q\defeq(\max\Rightarrow A)(q)$ has horizontal and
  vertical source $q$ and horizontal and vertical targets which are
  identities in $\pathcat(A)$. This makes it a suitable candidate to
  be the connection square of
  condition~\ref{dcat.with.con}(\ref{con.four}) as required.
\end{proof}



\section{Decalage for Complicial Sets and Complicially Enriched
  Categories}\label{pathecat.ff.sec}

\subsection{A Decalage Construction on Complicial Sets}
\label{dec.on.comp.subsect}

We show that we may obtain the appropriate generalisation of the
decalage construction to $\Comp$ as a sub-functor of
$\Delta[1]\Rightarrow{*}$.

\begin{defn}[decalage on $\Comp$]\label{dec.on.comp.defn}
  If $A$ is a complicial set then define $\Dec(A)$, its {\em
    decalage}, to be the complicial set obtained by taking the
  following pullback in $\Comp$: 
  \begin{equation}\label{dec.on.comp.pb}
    \let\labelstyle=\textstyle
    \xymatrix@R=2em@C=3em{
      {\Dec(A)}\pbexcursion\ar@{u(->}[r]^-{\subseteq_r}\ar[d] &
      {\Delta[1]\Rightarrow A}\ar[d]^{\Delta(\vertex^1_0)\Rightarrow A} \\
      {\Sup_0(A)}\ar@{u(->}[r]_-{\subseteq_r} & {A}
    }
  \end{equation}
  As a functor, $\Dec$ is the point-wise pullback, in the functor
  category $\funcat[\Comp,\Comp]$, of the natural inclusion
  $\overinc\subseteq_r: \Sup_0(\cdot)->\id_{\Comp}.$ along the natural
  transformation $\arrow \Delta(\vertex^1_0)\Rightarrow{*}:
  \Delta[1]\Rightarrow{*} ->\id_{\Comp}.$. In other words, for each
  stratified map $\arrow f:A->B.$ the map $\arrow \Delta[1]\Rightarrow
  f: \Delta[1]\Rightarrow A->\Delta[1]\Rightarrow B.$ restricts to a
  map $\arrow\Dec(f):\Dec(A)->\Dec(B).$, thereby making our
  construction into a functor $\arrow\Dec:\Comp->\Comp.$ for which the
  family of regular subset inclusions $\overinc\subseteq_r:\Dec(A)->
  \Delta[1]\Rightarrow A.$ constitute a natural transformation from
  $\Dec$ to $\Delta[1]\Rightarrow{*}$.
  
  Using the explicit description of the $n$-simplices of
  $\Delta[1]\Rightarrow A$ of observation~\ref{expl.power}, we
  may describe $\Dec(A)$ as the regular subset of
  $\Delta[1]\Rightarrow A$ on those simplices $\arrow
  p:\Delta[1]\otimes\Delta[n]->A.$ for which the simplex
  $(\Delta(\vertex^1_0)\Rightarrow A)(p) =
  p\pair<\vertex^1_0\circ\eta^n;\id_{[n]}>$ is in the regular subset
  $\Sup_0(A)\subseteq_r A$. 
  
  Another way of re-stating this characterisation is to construct an
  entire superset $(\Delta[1]\otimes\Delta[n])'$ of
  $\Delta[1]\otimes\Delta[n]$ by making thin every simplex
  $\pair<\alpha;\beta>\in\Delta[1]\otimes\Delta[n]$ which has
  dimension $r>0$ and $\alpha(i)=0$ for all $i\in[r]$. Then an
  $n$-simplex $\arrow p:\Delta[1]\otimes\Delta[n]->A.$ in
  $\Delta[1]\Rightarrow A$ is in the regular subset $\Dec(A)$ if and
  only if it extends to a stratified map $\arrow
  p:(\Delta[1]\otimes\Delta[n])'->A.$. 
  
  We make $\Dec$ into a copointed endo-functor on $\Comp$ by defining
  the counit $\arrow \cface:\Dec->\id_{\Comp}.$ to be the natural
  transformation obtained by composing the inclusion
  $\overinc\subseteq_r:\Dec-> \Delta[1]\Rightarrow{*}.$ with the
  natural transformation $\arrow\Delta(\vertex^1_1)\Rightarrow{*}:
  \Delta[1]\Rightarrow{*}-> \id_{\Comp}.$. In terms of the explicit
  presentation of simplices of $\Dec(A)$, this counit maps an
  $n$-simplex $\arrow p:\Delta[1]\otimes\Delta[n]->A.$ of $\Dec(A)$ to
  the simplex $\cface_A(p) = (\Delta(\vertex^1_1)\Rightarrow A)(p) =
  p\pair<\vertex^1_1\circ\eta^n;\id_{[n]}>$ of $A$.
\end{defn}

\begin{notation}
  In the subsequent arguments we will have occasion to consider, and
  relate, three distinct decalage constructions, each one on a
  different category. In particular, we now have the classical
  decalage comonad on $\Simp$, which we discussed in
  lemma~\ref{nerve.dec}, and the construction on $\Comp$, introduced in
  the previous definition. We will be adding to these later in this
  subsection by describing an analogous copointed endo-functor on
  $\CompCat$.

  In order to reduce notational clutter, we will overload the symbols
  $\Dec$ and $\cface$ in order to use them to denote the decalage
  construction on whichever of these categories is currently under
  consideration. Consequently, we will also be taking care to ensure
  that the reader may disambiguate this notation contextually.
\end{notation}

\begin{thm}\label{dec.on.comp}
  Suppose that $\arrow \forget:\Comp->\Simp.$ denotes the forgetful functor,
  which maps a complicial set to its underlying simplicial set,
  then there exists a natural isomorphism
  \begin{equation}\label{dec.on.comp.diag}
    \let\labelstyle=\textstyle
    \xymatrix@C=3.5em@R=2em{
      {\Comp}\ar[d]_{\forget}\ar[r]^{\Dec}_{}="one" &
      {\Comp}\ar[d]^{\forget} \\
      {\Simp}\ar[r]_{\Dec}^{}="three" & 
      {\Simp}
      \ar@{<=}_{\textstyle\varphi}^{\cong}
        "three"!<0ex,2ex>; "one"!<0ex,-2ex>
    }
  \end{equation}
  in $\Cat$ which makes the pair $\pair<\forget;\varphi>$ into a strong
  transformation of copointed endo-functors from the decalage
  construction on $\Comp$ to that on $\Simp$. 
\end{thm}

Most of the remainder of this subsection is devoted to proving this
theorem, which we do via a sequence of observations and lemmas:

\begin{obs}[strong transformations of copointed endo-functors]
  \label{st.of.cp.endo}
  Transformations of copointed endo-functors are defined analogously
  to the comonad transformations discussed in
  observation~\ref{comonad.trans}. That is, the square of
  display~(\ref{dec.on.comp.diag}) represents an endo-1-cell of the
  2-category $\Cylinder$ on the forgetful functor $\arrow
  \forget:\Comp->\Simp.$ and we require that the pair consisting of
  the counits associated with the functor $\Dec$ on the categories
  $\Comp$ and $\Simp$ (respectively) is a cylinder (2-cell) from this
  endo-1-cell to the identity 1-cell on $\forget$.
  
  Explicitly, this latter condition simply requires that the natural
  transformation $\varphi$ must satisfy the equation
  $\forget\cface = (\cface\forget)\cdot\varphi$. Here it is worth observing that
  the use context of the symbol $\cface$ on either side of this
  equation implies that these instances represent different natural
  transformations. Since $\forget$ is a functor from $\Comp$ to $\Simp$ we
  know that $\cface$ in the context $\forget\cface$ must be the counit for
  $\Dec$ on $\Comp$ and, conversely, that in the context $\cface \forget$ it
  instead represents the counit of $\Dec$ on $\Simp$.
  
  The strong transformation postulated in theorem~\ref{dec.on.comp}
  simply demonstrates that the construction in
  definition~\ref{dec.on.comp.defn} really deserves to be
  characterised as a generalisation of the decalage construction to
  the category of complicial sets.
  
  Assuming for the remainder of this observation that
  theorem~\ref{dec.on.comp} holds and consulting
  observation~\ref{conn.cpt}, we know that the set of connected
  components of a complicial set is dependent only on the structure of
  it underlying simplicial set. In other words, we have a commuting
  triangle of functors
  \begin{displaymath}
    \let\labelstyle=\textstyle
    \xymatrix@R=2em@C=2em{
      {\Comp}\ar[rr]^{\forget}\ar[dr]_{\cpt_0}^{}="one" & &
      {\Simp}\ar[dl]^{\cpt_0}_{}="two" \\
      & {\Set} & 
      \ar@{=} "two"+<-2.3em,0.3em>;"one"+<2.3em,0.3em>
    }
  \end{displaymath}
  and so we may apply the semi-simplicial version of
  lemma~\ref{nerves.and.comonad.tr} to the strong transformation
  postulated in lemma~\ref{dec.on.comp} to give a natural isomorphism:
  \begin{displaymath}
    \let\labelstyle=\textstyle
    \xymatrix@R=2em@C=2em{
      {\Comp}\ar[rr]^{\forget}\ar[dr]_{\Ndec}^{}="one" & &
      {\Simp}\ar[dl]^{\Ndec}_{}="two" \\
      & {\SSimp} & 
      \ar@{}|\cong "two"+<-2.3em,0.3em>;"one"+<2.3em,0.3em>
    }
  \end{displaymath}
  Applying the semi-simplicial version of lemma~\ref{nerve.dec} we
  may show that the nerve functor $\arrow \Ndec: \Simp->\SSimp.$ is
  isomorphic to the forgetful functor from $\Simp$ to $\SSimp$.  It
  follows from the triangle above that the nerve functor $\arrow
  \Ndec:\Comp->\SSimp.$ associated with the construction introduced in
  definition~\ref{dec.on.comp.defn} is isomorphic to the forgetful
  functor from the category of complicial sets $\Comp$ to the category
  of semi-simplicial sets $\SSimp$.
\end{obs}

\begin{obs}[constructing the natural transformation $\varphi$ of 
  theorem~\ref{dec.on.comp}]\label{ladder.of.swords.2}\label{en.defn}
  We define the isomorphism $\varphi$ in two steps, here we construct
  a (candidate) inverse $\arrow\varphi^{*}:\Dec\circ
  \forget->\forget\circ\Dec.$ and in subsequent lemmas we demonstrate
  that this is indeed a natural isomorphism.
  
  Consulting the construction of definition~\ref{defn.dec}, and
  applying Yoneda's lemma for stratified sets, we see that if $A$ is a
  complicial set then an $n$-simplex of $\Dec(\forget(A))$ may be identified
  with a stratified map $\arrow x:\Delta[n+1]->A.$.  Furthermore,
  under this identification the right action of simplicial operator
  $\phi$ on our $n$-simplex $x$ in $\Dec(\forget(A))$ is given by
  $x\cdot\phi=x\circ \Delta(\phi\oplus[0])$.  
  
  Of course we also know, from definition~\ref{dec.on.comp.defn}, that
  an $n$-simplex of $\forget(\Dec(A))$ may be identified with a stratified
  map $\arrow p:(\Delta[1]\otimes\Delta[n])'->A.$ and that under this
  identification $p\cdot\phi=p\circ(\Delta[1]\otimes \Delta(\phi))$ as
  discussed in observation~\ref{expl.power}.
  
  Thus we are drawn to considering the relationship between
  $\Delta[n+1]$ and $(\Delta[1]\otimes\Delta[n])'$, so starting with an
  $r$-simplex $\pair<\alpha;\beta>$ of $\Delta[1]\otimes\Delta[n]$ we
  define the $r$-simplex $e_n\pair<\alpha;\beta>$ of $\Delta[n+1]$ by
  \begin{displaymath}
    (e_n\pair<\alpha;\beta>)(i) = \left\{
      \begin{array}{ll}
        0 & \text{if $\alpha(i)=0$,} \\
        \beta(i)+1 & \text{if $\alpha(i)=1$.}
      \end{array}
    \right.
  \end{displaymath}
  The point-wise nature of this construction ensures that $e_n$
  respects the right actions of $\Delta$ on
  $\Delta[1]\otimes\Delta[n]$ and $\Delta[n+1]$ since these are given
  by pre-composition, in other words $e_n$ is a simplicial map.
  Furthermore, $e_n$ maps any simplex $\pair<\alpha;\beta>$ which has
  $\alpha(i)=0$ for some $i>0$ to a degenerate simplex in
  $\Delta[n+1]$.  It follows, directly from the analysis of
  lemma~\ref{tensor.simp.set} and the definition of
  $(\Delta[1]\otimes\Delta[n])'$ in definition~\ref{defn.dec}, that
  $e_n$ carries thin simplices in $(\Delta[1]\otimes\Delta[n])'$ to
  degenerate simplices in $\Delta[n+1]$ and is thus a stratified map
  $\arrow e_n: (\Delta[1]\otimes\Delta[n])'->\Delta[n+1].$.  We may
  now use this to define a function
  \begin{displaymath}
    \let\labelstyle=\textstyle
    \xymatrix@R=3ex@C=11em{
      {\Dec(\forget(A))_n}\ar[r]^-{(\varphi_A)_n}_{}="one"
      \ar@{{}={}}[d] & {\forget(\Dec(A))_n}\ar@{{}={}}[d]\\
      {\Strat(\Delta[n+1],A)}\ar[r]_-{\Strat(e_n,A)}^{}="two" &
      {\Strat((\Delta[n]\otimes\Delta[1])',A)}
      \ar@{{}={}}"one"-<0ex,2ex>;"two"+<0ex,2ex>_{\text{def}}
    }
  \end{displaymath}
  which is natural in $A$ by construction. Observe also that if
  $\arrow\phi:[m]->[n].$ is a simplicial operator then it is a matter
  of routine calculation to demonstrate that the following square of
  stratified maps commutes:
  \begin{equation}\label{e.simpop.nat}
    \let\labelstyle=\textstyle
    \xymatrix@=3em{
      {(\Delta[1]\otimes\Delta[m])'}
      \ar[r]^<>(0.5){e_m}\ar[d]_{\Delta[1]\otimes\Delta(\phi)} &
      {\Delta[m+1]}
      \ar[d]^{\Delta(\phi\oplus[0])} \\
      {(\Delta[1]\otimes\Delta[n])}\ar[r]_<>(0.5){e_n} & {\Delta[n+1]}}
  \end{equation}
  Returning to the explicit descriptions furnished by the first couple
  of paragraphs of this observation, it is clear that we may use the
  commutativity of these squares to show that the collection of maps
  $(\varphi^{*}_A)_n$ ($n\in\mathbb{N}$) respect the actions of
  $\Delta$ on $\Dec(\forget(A))$ and $\forget(\Dec(A))$ and thus provide us with
  the components of a simplicial map $\arrow
  \varphi^{*}_A:\Dec(\forget(A))->\forget(\Dec(A)).$ as desired. To do so, simply
  consider an arbitrary simplicial operator $\phi$ and observe that we
  have $(\varphi^{*}_A)_n(x)\cdot \phi = (x\circ e_n)\cdot\phi =
  x\circ e_n \circ (\Delta(\phi)\otimes\Delta[1]) =^* x\circ
  \Delta([0]\oplus\phi) \circ e_m = (x\cdot\phi)\circ e_m =
  (\varphi^{*}_A)_m(x\cdot \phi)$ wherein the starred equality follows
  from (\ref{e.simpop.nat}).
\end{obs}

\begin{obs}
  \label{ladder.of.swords}
  By definition the component $(\varphi^{*}_A)_n=\Strat(e_n,A)$ is an
  isomorphism iff $e_n$ is orthogonal to $A$, so it follows that
  $\arrow \varphi^{*}_A:\Dec(\forget(A))->\forget(\Dec(A)).$ is an isomorphism for each
  complicial set $A$ iff the stratified map $\arrow
  e_n:(\Delta[n]\otimes\Delta[1])'->\Delta[n+1].$ is an f-extension
  for each $n\in\mathbb{N}$. Consequently, all that remains in proving
  theorem~\ref{dec.on.comp} is to establish this latter condition.
  
  To this end, define $(\Delta[1]\otimes\Delta[n])''$ to be the entire
  superset of $\Delta[1]\otimes\Delta[n]$ in which we make thin all
  those simplices $\pair<\alpha;\beta>$ for which $\alpha(i)=0$ for
  some $i>0$. As discussed in the last observation, the map $e_n$
  actually maps each of the thin simplices in
  $(\Delta[1]\otimes\Delta[n])''$ to a degenerate (and thus thin)
  simplex in $\Delta[n+1]$. However a little more is actually true, in
  fact $e_n$ extends to a {\em regular\/} stratified map
  $\overarr:(\Delta[1]\otimes\Delta[n])'' ->\Delta[n+1].$.

  We'll actually prove that $\arrow e_n:(\Delta[1]\otimes\Delta[n])'
  ->\Delta[n+1].$ is an f-extension by decomposing our proof into two
  steps 
  \begin{itemize}
  \item first show that the entire inclusion $(\Delta[1]\otimes\Delta[n])'
    \subseteq_e(\Delta[1]\otimes\Delta[n])''$ is a t-extension,
  \item then prove that the regular map $\arrow e_n:
    (\Delta[1]\otimes\Delta[n])''->\Delta[n+1].$ is an f-extension,
  \end{itemize}
  and finally appeal to the fact that the composite of these two, 
  $\arrow e_n:(\Delta[1]\otimes\Delta[n])'->\Delta[n+1].$ itself, is
  therefore an f-extension.
\end{obs}

\begin{lemma}
  The entire inclusion $\overinc\subseteq_e:
  (\Delta[1]\otimes\Delta[n])'->(\Delta[1]\otimes\Delta[n])''.$ is
  a t-extension (and thus an f-extension).
\end{lemma}

\begin{proof} 
  We aim to apply the explicit characterisation of t-extensions given in
  observation~\ref{t-ext.expl.obs}, so we define a sequence
  of subsets $\{tX_k\}_{k\in[n]}$ of $\Delta[1]\otimes\Delta[n]$ by:
  \begin{displaymath}
    tX_k \defeq (t\Delta[1]\otimes t\Delta[n])\cup
    \left\{\,\pair<\alpha;\beta>\in\Delta[1]\otimes\Delta[n]\,\big|\,
    \dim\pair<\alpha;\beta>\geq n-k+1 \wedge \alpha(n-k+1)=0 \,\right\}
  \end{displaymath}
  It is a straightforward matter, using the fact that each $\alpha$ in
  the definition above is order preserving, to show that $tX_0$ is the
  set of thin simplices of $(\Delta[1]\otimes\Delta[n])'$, that $tX_n$
  is the set of thin simplices of $(\Delta[1]\otimes\Delta[n])''$ and
  that for all $0\leq k < n$ we have $tX_k\subseteq tX_{k+1}$.
  
  To verify the remaining condition of
  observation~\ref{t-ext.expl.obs}, suppose that $0\leq k <n$ and that
  $\pair<\alpha;\beta>$ is an $r$-simplex in $tX_{k+1}\setminus
  tX_{k}$ then from the definitions of these sets we know that $r \geq
  n-k > 0$, $\alpha(n-k)=0$ and $\alpha(n-k+1)=1$. Furthermore, we
  also know that $\beta(n-k+1)>\beta(n-k)$ since otherwise $(n-k)$
  would witness that fact that $\pair<\alpha;\beta>$ was a mediator
  simplex and thus an element of $t\Delta[1]\otimes t\Delta[n]$ (cf.\ 
  observation~\ref{med.cyl}).
  
  So define a new $(r+1)$-simplex $\pair<\alpha';\beta'>\defeq
  \pair<\alpha\circ\degen^r_{n-k};\beta\circ\degen^r_{n-k+1}>$ and
  observe that:
  \begin{itemize}
  \item The integer $(n-k+1)$ witnesses the fact that
    $\pair<\alpha';\beta'>$ is a mediator simplex and is thus
    $(n-k+1)$-admissible.
  \item The $(n-k)^{\text{th}}$ face
    $\pair<\alpha';\beta'>\cdot\face^{r+1}_{n-k}$ is equal to
    $\pair<\alpha;\beta\circ\degen^r_{n-k+1}\circ\face^{r+1}_{n-k}>$
    which is a mediator simplex witnessed by $n-k$ and is thus an
    element of $t\Delta[1]\otimes t\Delta[n]\subseteq tX_k$.
  \item The $(n-k+1)^{\text{th}}$ face
    $\pair<\alpha';\beta'>\cdot\face^{r+1}_{n-k+1}$ is equal to our
    original simplex $\pair<\alpha;\beta>$.
  \item The $(n-k+2)^{\text{th}}$ face
    $\pair<\alpha';\beta'>\cdot\face^{r+1}_{n-k+2}$ is equal to
    $\pair<\alpha\circ\degen^r_{n-k}\circ\face^{r+1}_{n-k+2};\beta>$
    which is an element of $tX_k$ since
    $\alpha\circ\degen^r_{n-k}\circ\face^{r+1}_{n-k+2}(n-k+1)
    =\alpha(n-k) = 0$.
  \end{itemize}
  In other words, $\pair<\alpha';\beta'>$ witnesses the extension of
  thinness to $\pair<\alpha;\beta>$ as required by
  observation~\ref{t-ext.expl.obs}, which we may now apply to
  establish the proposition.
\end{proof}

\begin{lemma}
  The stratified map $\arrow e_n:(\Delta[1]\otimes\Delta[n])''
    ->\Delta[n+1].$ is an f-extension.
\end{lemma}

\begin{proof}
  Start by considering the $(n+1)$-simplex
  $\pair<\rho^{n+1}_1;\degen^n_0>$ of $(\Delta[1]\otimes\Delta[n])''$,
  where the simplicial operators $\arrow \rho^r_k:[r]->[1].$ are as
  defined in observation~\ref{prim.d.tens.note}. From the definition
  of $e_n$, given in observation~\ref{en.defn}, we see that it carries
  $\pair<\rho^{n+1}_1;\degen^n_0>$ to the simplex $\id_{[n+1]}$ in
  $\Delta[n+1]$ so, applying Yoneda's lemma, we see that the map
  $\arrow i_n=\yoneda{\pair< \rho^{n+1}_1;\degen^n_0>}:\Delta[n+1]->
  (\Delta[1]\otimes \Delta[n])''.$ is a right inverse to $\arrow
  e_n:(\Delta[1]\otimes\Delta[n])''->\Delta[n+1].$.
  
  Consequently, we may argue as we did in the proof of
  corollary~\ref{prim.d.tens} and show that our desired result holds
  iff whenever $A$ is a complicial set and $\arrow
  g,g':(\Delta[1]\otimes\Delta[n])''->A.$ are stratified maps with
  $g\circ i_n = g'\circ i_n$ which, by Yoneda's lemma, is the same as
  saying that $g\pair<\rho^{n+1}_1;\degen^n_0> =
  g'\pair<\rho^{n+1}_1;\degen^n_0>$ then $g=g'$.
  
  But since the underlying simplicial set of
  $(\Delta[1]\otimes\Delta[n])''$ is $\Delta[1]\times\Delta[n]$ we
  know, by observation~\ref{shuffles}, that every one of its simplices
  is a face of some shuffle. It follows that the stratified maps $g$
  and $g'$ coincide iff they act identically on each of the shuffles
  $\pair<\rho^{n+1}_{k+1};\degen^n_k>$ of
  $\Delta[n]\times\Delta[1]$ ($k=0,...,n$).
  
  Now fix $0<k\leq n$ and consider the simplex
  $\pair<\rho^{n+1}_{k+1};\degen^n_k>$.  Suppose that
  $\arrow\alpha:[r]->[n+1].$ is a face operator whose image contains
  $k-1$ and $k$, let $l\in[r]$ be the unique integer such that
  $\alpha(l)=k$ and consider the face
  $\pair<\rho^{n+1}_{k+1};\degen^n_k>\cdot\alpha=
  \pair<\rho^{n+1}_{k+1}\circ\alpha;\degen^n_k\circ\alpha>$. It is
  clear that $l > 0$ and that $\rho^{n+1}_{k+1}\circ\alpha(l)=0$ so it
  follows that this face is a thin simplex in
  $(\Delta[1]\otimes\Delta[n])''$. Quantifying over these face
  operators $\alpha$ it follows that our
  $\pair<\rho^{n+1}_{k+1};\degen^n_k>$ of
  $(\Delta[1]\otimes\Delta[n])''$ is pre-degenerate at $k-1$ in there.

  Returning to our stratified map $\arrow g:(\Delta[1]\otimes
  \Delta[n])''->A.$, we know that it preserves pre-degeneracy and so
  we may infer that the simplex $g\pair<\rho^{n+1}_{k+1};\degen^n_k>$
  in $A$ is pre-degenerate at $k-1$. However since $A$ is complicial,
  and thus well tempered by lemma~\ref{lemma.predegen}, we may infer
  that this simplex is in fact degenerate at $k-1$. 

  \begin{figure}[h]
    \begin{displaymath}
      \def\vertchar{\scriptscriptstyle\bullet}
      \xymatrix@!0@C=2.5em@R=2.5em{
        & *[o]{\vertchar} & *[o]{\vertchar} \\ 
        & *[o]{\vertchar} & *[o]{\vertchar} \\ 
        & *[o]{\vertchar} & *[o]{\vertchar} \\ 
        & *[o]{\vertchar} & *[o]{\vertchar} \\ 
        [n]\ar[uu] & *[o]{\vertchar} & *[o]{\vertchar} \\
        & [1] \ar[r] & 
        \ar@{-}"5,2";"3,2" \ar@{-}"3,2";"3,3" \ar@{-}"3,3";"1,3" 
        \ar@{..}"3,2";"2,3" \ar@{--}"3,2";"2,2" \ar@{--}"2,2";"2,3"
        \save "3,3"+<5em,1.25em>*+{\scriptstyle
          \pair<\rho^{n+1}_k;\degen^n_{k-1}>}
        \ar@{.>}"3,3"+<0em,1.25em>\restore
        \save "3,2"+<-5em,1.25em>*+{\scriptstyle
          \pair<\rho^{n+1}_{k+1};\degen^n_k>}
        \ar@{.>}"3,2"+<0em,1.25em>\restore
        \save "4,3"+<-0.75em,1em>*{\scriptstyle
          \text{vertex $k$}}="one"
        \ar@{.>}"one";"3,3"\ar@{.>}"one";"2,2" \restore
      }
    \end{displaymath}
    \caption{An adjacent pair of shuffles in
      $\Delta[1]\otimes\Delta[n]$}
    \label{shuffle.pic.4}
  \end{figure}
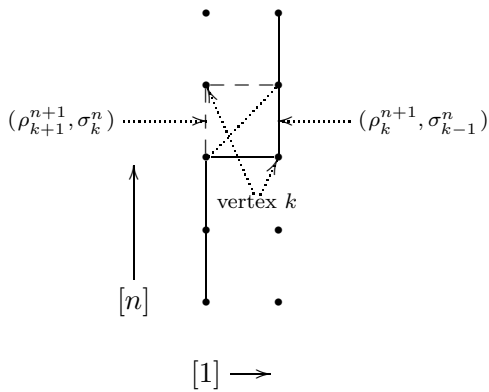
  
  Figure~\ref{shuffle.pic.4} depicts our shuffle
  $\pair<\rho^{n+1}_{k+1};\degen^n_k>$ and the immediately
  preceding one $\pair<\rho^{n+1}_k;\degen^n_{k-1}>$. It clearly
  illustrates that they share the same $k^{\text{th}}$ $n$-face and
  that this is the simplex $\pair<\rho^n_k;\id_{[n]}>$, facts which
  may be verified by straightforward calculations with simplicial
  operators.  From this, and the fact that $g$ is a simplicial map, we
  immediately see that $g\pair<\rho^n_k;\id_{[n]}>$ is the
  $k^{\text{th}}$ $n$-face of $g\pair<\rho^{n+1}_k;\degen^n_{k-1}>$
  and that the degeneracy result of the previous paragraph implies
  that the simplex $g\pair<\rho^{n+1}_{k+1};\degen^n_k>$ is equal to
  its degenerated face $g\pair<\rho^n_k;\id_{[n]}>
  \cdot\degen^n_{k-1}$.  Combining these two, we get the equality
  $g\pair<\rho^{n+1}_{k+1};\degen^n_k> =
  g\pair<\rho^{n+1}_k;\degen^n_{k-1}>\cdot(\face^{n+1}_k\circ
  \degen^n_{k-1})$ which demonstrates that the action of $g$ on the
  shuffle $\pair<\rho^{n+1}_{k+1};\degen^n_k>$ is completely
  determined by its action on the adjacent shuffle
  $\pair<\rho^{n+1}_k;\degen^n_{k-1}>$.
  
  Finally we can apply this result repeatedly for $k=n,n-1,...,1$ to
  show that the action of $g$ on each shuffle is determined by its
  action on the shuffle that immediately precedes it (under the linear
  ordering $\triangleleft$) and that, ultimately, the action of $g$ on
  all shuffles, and thus on all simplices of
  $(\Delta[1]\otimes\Delta[n])''$, is completely determined by its
  action on the $\triangleleft$-minimal shuffle
  $\pair<\rho^{n+1}_1;\degen^n_0>$.  Of course, we can also apply
  exactly the same argument to $g'$ and we know from our initial
  assumption that these two stratified maps act identically on
  $\pair<\rho^{n+1}_1;\degen^n_0>$, it follows therefore that they
  must be identical maps as required.
\end{proof}

\begin{proof}(of theorem~\ref{dec.on.comp}). 
  Finally this follows on observing that the
  composite of the two f-extensions provided by these lemmas is the
  stratified map $\arrow
  e_n:(\Delta[n]\otimes\Delta[1])'->\Delta[n+1].$ which is therefore
  also an f-extension as required by
  observation~\ref{ladder.of.swords}. That the pair $\pair<\forget;\varphi>$
  satisfies the cylinder condition required of a transformation of
  copointed endo-functors is a matter of trivial verification, directly
  from the definition of $\varphi^{*}$, which we leave to the reader.
\end{proof}

\subsection{A Path Construction on Complicially Enriched Categories}

Now we turn to our primary motivation for introducing the double
category with connection $\strat{\pathsq(A)}$ in the previous section,
that is to build an analogue of the functor $\Delta[1]\Rightarrow{*}$
on $\CompCat$. In the next subsection we use this analogue to
emulate the construction of the last subsection on $\CompCat$.

\begin{obs}[what do we mean by an ``analogue'' of 
  {$\Delta[1]\Rightarrow{*}$}]\label{interv.analogue.defn} First we
  should make precise the sense in which we shall be ``building an
  analogue of $\arrow\Delta[1]\Rightarrow{*}:\Comp-> \Comp.$ on
  $\CompCat$''. In fact, our intention is to construct a functor
  $\arrow\interv:\CompCat->\CompCat.$ and a pair of natural
  transformations $\arrow\iend^0,\iend^1:\interv->\id_{\CompCat}.$
  and that this data should come equipped with a natural isomorphism
  \begin{equation}\label{dp.square}
    \let\labelstyle=\textstyle
    \xymatrix@R=3em@C=5em{
      {\Comp}\ar[r]^{\Delta[1]\Rightarrow{*}}^{}="one"
      \ar[d]_{\pathecat} &
      {\Comp}
      \ar[d]^{\pathecat} \\
      {\CompCat}
      \ar[r]_{\interv}^{}="two" &
      {\CompCat}
      \ar@{=>}^{\textstyle\theta}_{\cong}
       "one"-<0em,1.25em>; "two"+<0em,1.25em>
    }
  \end{equation}
  making the triangles
  \begin{equation}\label{dp.triangles}
    \let\labelstyle=\textstyle
    \xymatrix@C=1em@R=2.5em{
      *+[l]!<-1em,0ex>{\pathecat(\Delta[1]\Rightarrow A)}
      \ar[rr]^-{\theta_A}_-{\cong}
      \ar[rd]_{\scriptstyle\pathecat(\Delta(\vertex^1_0)\Rightarrow A)}  & &
      *+[r]!<1em,0ex>{\interv(\pathecat(A))}
      \ar[ld]^{\iend^0_{\pathecat(A)}} \\ 
      & {\pathecat(A)} &
    }\mkern30mu
    \xymatrix@C=1em@R=2.5em{
      *+[l]!<-1em,0ex>{\pathecat(\Delta[1]\Rightarrow A)}
      \ar[rr]^-{\theta_A}_-{\cong}
      \ar[rd]_{\scriptstyle\pathecat(\Delta(\vertex^1_1)\Rightarrow A)}  & &
      *+[r]!<1em,0ex>{\interv(\pathecat(A))}
      \ar[ld]^{\iend^1_{\pathecat(A)}} \\ 
      & {\pathecat(A)} &
    }
  \end{equation}
  commute for each complicial set $A$.
\end{obs}

Our approach to building such a functor will be to start by
identifying $\pathecat(\Delta[1]\Rightarrow A)$ as a substructure of
$\pathsq(A)$. We then show that when we may apply the functor
$\arrow\pathcat:\Comp->\Cat(\Comp).$ ``point-wise'' to a complicially
enriched category, as described in definition~\ref{le-th.defn}, we
obtain a double category the reflection dual of which is in fact
obtain a 2-category and that in the particular case of
$\pathcat(\pathecat(A))$ this is also a substructure of $\pathsq(A)$.
Finally we show that we may apply the double category of squares
construction to $\refld{\pathcat(\pathecat(A))}$ to re-construct a part
of $\pathsq(A)$ which contains $\pathecat(\Delta[1]\Rightarrow A)$.

When handling path categories it will be convenient to restrict our
attention to a certain sub-double category of $\pathsq(A)$:

\begin{defn}\label{pathesq.defn}
  Let $\pathesq(A)$ denote the sub-double category of $\pathsq(A)$
  consisting of those squares $\lambda$ for which each one of the four
  ``corner'' objects $\s(\s_h(\lambda))=\s(\s_v(\lambda))$,
  $\s(\t_h(\lambda))=\t(\s_v(\lambda))$,
  $\t(\s_h(\lambda))=\s(\t_v(\lambda))$ and
  $\t(\t_v(\lambda))=\t(\t_h(\lambda))$ are simplices in the
  0-superstructure $\Sup_0(A)$ of $A$.  Alternatively, a square
  $\lambda$ is in $\pathesq(A)$ iff $\s_h(\lambda)$ and
  $\t_h(\lambda)$ are both in the complicial subcategory
  $\pathecat(A)\subseteq_r\pathcat(A)$ or, equivalently, iff
  $\s_v(\lambda)$ and $\t_v(\lambda)$ are both in $\pathecat(A)$.
  
  That this set of squares is indeed closed in $\pathsq(A)$ under its
  horizontal and vertical category structures is a trivial consequence
  of the fact that the corner objects of composites, sources and
  targets of squares are all corners of one or other of the original
  squares being operated upon (cf.  observation~\ref{dcat.piccy}). So
  if the original squares all have corners in $\Sup_0(A)$ then so do
  their composites etcetera.
  
  Furthermore, let $t\pathesq(A)$ denote the intersection
  $\pathesq(A)\cap t\pathsq(A)$ then, by the same reasoning as the
  last paragraph applied to theorem~\ref{pathsq.conn}, it follows that
  $\strat{\pathesq(A)}$ is also a double category with connection.
\end{defn}

\begin{obs}[{$\pathecat(\Delta[1]\Rightarrow A)$} as a 
  substructure of {$\pathesq(A)$}]\label{study.pecat.d1} Suppose that
  $\mathbb{D}$ is a complicial double category, then we may define an
  associated double category $\efunc_h(\mathbb{D})$ by applying the
  co-reflector $\arrow \efunc:\Cat(\Comp)->\CompCat.$ point-wise to the
  vertical presentation of $\mathbb{D}$ (cf.\ 
  definition~\ref{le-th.defn}). In other words, $\efunc_h(\mathbb{D})$ is
  the regular sub-double category of $\mathbb{D}$ of those squares and
  horizontal arrows whose horizontal source and target are both in the
  superstructure $\Sup_0({\arr_v(\mathbb{D})})\subseteq_r\arr_v(
  \mathbb{D})$. Observe also that, by definition, its horizontal
  categories of arrows and squares are actually complicially enriched,
  since they are $\efunc(\arr_h(\mathbb{D}))$ and $\efunc(\sq_h(\mathbb{D}))$
  respectively.
  
  From observation~\ref{iter.path}, we know that the vertical
  presentation of $\pathsq(A)$ is constructed by applying the left
  exact functor $\arrow\pathcat:\Comp->\Cat(\Comp).$ point-wise to the
  prism category $\pathcat(A)\in\Cat(\Comp)$. It follows, therefore,
  that the vertical presentation of $\efunc_h(\pathsq(A))$ is obtained by
  applying $\pathecat= \efunc\circ\pathcat$ point-wise to
  $\pathcat(A)$, or in other words it is simply the double category
  we would usually refer to as $\pathecat(\pathcat(A))$.
  
  Since the complicial sets of arrows and objects of $\pathcat(A)$ are
  $\Delta[1]\Rightarrow A$ and $A$ respectively it follows that the
  horizontal categories of squares and arrows of
  $\pathecat(\pathcat(A))$ are $\pathecat(\Delta[1]\Rightarrow A)$ and
  $\pathecat(A)$ respectively. Additionally, the source and target
  maps of $\pathcat(A)$ are $\s=\Delta(\vertex^1_0)\Rightarrow A$ and
  $\t=\Delta(\vertex^1_1)\Rightarrow A$ therefore the vertical
  presentation of $\pathecat(\pathcat(A))$ has vertical source and
  target functors given by $\s_v=\pathecat(\Delta(\vertex^1_0)
  \Rightarrow A)$ and $\t_v=\pathecat(\Delta(\vertex^1_1)\Rightarrow
  A)$. These observations explain our interest in the sub-double
  category $\pathecat(\pathcat(A))$, in short it provides us with a
  convenient encapsulation of the structures on the left hand side of
  the triangles in display~(\ref{dp.triangles}) of
  observation~\ref{interv.analogue.defn} in terms of the vertical
  presentation of the double category $\pathsq(A)$. 
  
  Finally, consulting definition~\ref{pathesq.defn} we see that the
  fact that $\pathecat(\pathcat(A))$ has $\pathecat(A)$ as its
  horizontal category of arrows implies that it is actually contained
  within the sub-double category $\pathesq(A)$ of $\pathsq(A)$ and it
  follows immediately that we actually have
  $\pathecat(\pathcat(A))=\efunc_h(\pathesq(A))$.
\end{obs}

\begin{obs}[2-categories from complicially enriched categories]
  \label{2cat.from.compenr}
  Suppose that $\mathbb{C}$ is a complicially enriched category and
  consider the complicial double category $\pathcat(\mathbb{C})$
  discussed in observation~\ref{iter.path}. We know that its vertical
  presentation is constructed by applying the functor
  $\arrow\pathcat:\Comp-> \Cat(\Comp).$ ``point-wise'' to the
  structural components of $\mathbb{C}$ and in particular that its
  complicial category of horizontal arrows
  $\arr_h(\pathcat(\mathbb{C}))$ is equal to $\pathcat(\scobj{C})$. Of
  course since $\mathbb{C}$ is complicially enriched we know, by the
  definition in observation~\ref{comp.enr}, that its set of objects
  $\scobj{C}$ is 0-complicial and so we can infer, by
  observation~\ref{disc.pathcat}, that $\pathcat(\scobj{C})$ is the
  discrete category on $\scobj{C}$ itself. It follows therefore that
  the reflection dual $\refld{\pathcat(\mathbb{C})}$ is a complicial
  2-category.
  
  Consider now the path category $\pathecat(A)$.  We know that it is a
  regular subcategory of $\pathcat(A)$, so it follows that the
  complicial double category $\pathcat(\pathecat(A))$ is a regular
  sub-double category of $\pathsq(A)$.  More precisely, from
  observation~\ref{iter.path} we know that the horizontal category of
  squares of $\pathcat(\pathecat(A))$ is formed by applying the
  functor $\arrow\pathcat:\Comp->\Cat(\Comp).$ to the complicial set
  $\pathecat(A)_a$ which was defined as a pullback in
  observation~\ref{comp.enr}. However $\pathcat$ preserves pullbacks
  and it follows that the complicial category
  $\sq_h(\pathcat(\pathecat(A)))=\pathcat(\pathecat(A)_a)$ is the
  regular subcategory of $\sq_h(\pathsq(A))=\pathcat(\pathcat(A)_a)$
  obtained by pulling the discrete subcategory
  $\pathcat(\Sup_0(A))\times\pathcat(\Sup_0(A))$ of
  $\pathcat(A)\times\pathcat(A)$ back along the map $\arrow
  \pair<\pathcat(\s);\pathcat(\t)>=\pair<\s_v;\t_v>:\sq_h(\pathsq(A))->
  \pathcat(A)\times\pathcat(A).$.
  
  However, observation~\ref{disc.pathcat} tells us that
  $\pathcat(\Sup_0(A))$ is actually the discrete category on
  $\Sup_0(A)$, so it follows that $\pathcat(\pathecat(A))$ is simply
  the regular sub-double category of $\pathsq(A)$ consisting of those
  squares $\lambda$ whose vertical source $\s_v(\lambda)$ and target
  $\s_v(\lambda)$ are each horizontal identities in $\pathcat(A)$ on
  some object in $\Sup_0(A)$. Taking reflection duals, it follows that
  we can summarise this description of $\pathcat(\pathecat(A))$ using
  the sub-double category $\pathesq(A)$ of
  definition~\ref{pathesq.defn} and the 2-category of globs
  construction of observation~\ref{globs.defn} and simply observe that
  the 2-categories $\refld{\pathcat(\pathecat(A))}$ and $\glob(\refld{
    \pathesq(A)})$ are identical (as substructures of
  $\refld{\pathsq(A)}$).
\end{obs}

\begin{lemma}\label{recons.lemma}
  We may construct a functor
  $\arrow\pathcat^{*}:\CompCat->\Twocat(\Comp).$ and a family of
  isomorphisms $\pathcat^{*}(\pathecat(A))\cong\glob(\pathesq(A))$
  that is natural in the complicial set $A$. Furthermore, applying the
  double category of pasting squares construction (cf.
  observation~\ref{dcat.psq}) we obtain the natural isomorphism
  depicted in the following diagram:
  \begin{equation}\label{recons.iso}
    \let\labelstyle=\textstyle
    \xymatrix@=3em{
      {\Comp} \ar[r]^{\pathecat}
      \ar[d]_{\pathesq}^{}="one" & {\CompCat}
      \ar[d]^{\pathcat^{*}}_{}="two" \\
      {\Double(\Comp)} &
      {\Twocat(\Comp)}\ar[l]^{\squares}
      \ar@{=>}^{\cong}_{\phi} "one"+<3em,0em>;"two"-<3em,0em> }
  \end{equation}
\end{lemma}

\begin{proof} Taking our lead from the last observation, we define 
  $\arrow\pathcat^{*}:\CompCat->\Twocat(\Comp).$ to be the functor
  which takes a complicial category $\mathbb{C}$, maps it to the
  double category $\pathcat(\mathbb{C})$, applies the reflection dual
  to obtain a 2-category $\refld{\pathcat(\mathbb{C})}$ {\bf and then}
  applies the 2-categorists 2-cellular dual to obtain
  $\pathcat^{*}(\mathbb{C})\defeq\co{(\refld{\pathcat(\mathbb{C})})}$.
  Our reason for applying this final duality will make itself clear
  momentarily, for now it is enough to observe that this construction,
  as the composite of three functorial constructions, is itself
  functorial.
  
  Everything else follows from the fact that $\strat{\pathesq(A)}$ is
  a complicial double category with connection for each complicial set
  $A$. In particular, observation~\ref{duals.and.psquares} provides us
  with a natural isomorphism
  $\glob(\refld{\pathesq(A)})\cong\co{\glob( \pathesq(A))}$ and
  observation~\ref{2cat.from.compenr} shows that its left hand side is
  actually identical to $\refld{\pathcat(\pathecat(A))}$.  So applying
  the 2-cellular dual $\co{({-})}$, and using the fact that it is
  involutive, we get the first of the natural isomorphisms asked for
  in the statement of this lemma $\pathcat^{*}(\pathecat(A)) \defeq
  \co{(\refld{\pathcat( \pathecat(A))})}\cong\glob( \pathesq(A))$,
  which explains why we defined $\pathcat^{*}$ as we did.
  
  Applying the pasting squares construction
  $\arrow\squares:\Twocat(\Comp)->\Double(\Comp).$ to this natural
  isomorphism we get $\squares(\pathcat^{*}(\pathecat(A)))\cong
  \squares(\glob(\pathesq(A)))$ and then lemma~\ref{dcat.with.con}
  provides another isomorphism
  $\squares(\glob(\pathesq(A)))\cong\pathesq(A)$ which is also natural
  in $A$. Finally, composing the two we get the natural isomorphism
  $\phi$ asked for in display~(\ref{recons.iso}) of the statement.
\end{proof}

\begin{obs}[a more concrete description of
  $\arrow\phi:\pathesq->\squares\circ\pathcat^{*}\circ\pathecat.$]
  \label{recons.concrete}
  Tracing through the construction of lemma~\ref{dcat.with.con} and
  those of the last few observations, it is clear that the double
  categories $\pathesq(A)$ and $\squares(\pathcat^{*}( \pathecat(A)))$
  as we've constructed them share precisely the same categories of
  horizontal and vertical arrows, namely $\pathecat(A)$, and that the
  double functor $\arrow \phi_{A}:\pathesq(A)->\squares(\pathcat^{*}(
  \pathecat(A))).$ acts trivially on these. In other words, it carries
  each square $\lambda$ in $\pathesq(A)$ to a pasting square of the
  form:
  \begin{equation*}
    \let\labelstyle=\textstyle
    \xymatrix@=3em{
      \bullet \ar[r]^{\s_v(\lambda)}_{}="one"\ar[d]_{\s_h(\lambda)} & 
      \bullet\ar[d]^{\t_h(\lambda)} \\
      \bullet \ar[r]_{\t_v(\lambda)}^{}="two" & \bullet 
      \ar@{<=}^{\lambda^{*}} "two"+<0.25em,1em>;"one"+<0.25em,-1em>
    }
  \end{equation*}
  Returning again to lemma~\ref{dcat.with.con}, we see that the 2-cell
  $\lambda^{*}$ is constructed by composing $\lambda$ above and below
  with compatible connection squares in $\pathesq(A)$. By
  observation~\ref{2cat.from.compenr} this vertical composite
  $\nu_{\t_v(\lambda)}\comp_v\lambda\comp_v \mu_{\s_v(\lambda)}$ is
  actually in the sub-double category $\pathcat(\pathecat(A))$ and is
  easily seen to provide a 2-cell with the appropriate orientation in
  the dual 2-category $\pathcat^{*}(\pathecat(A))$.
\end{obs}

\begin{obs}\label{interv.analogue.obs}
  Given the result of lemma~\ref{recons.lemma} and the discussion in
  observation~\ref{study.pecat.d1}, the result alluded to in
  observation~\ref{interv.analogue.defn} is now a matter of mere
  formality. To illuminate that comment, define the functor
  $\arrow\interdc:\CompCat->\Double(\Comp).$ to be the composite
  \begin{equation*}
    \let\labelstyle=\textstyle
    \xymatrix@R=1ex@C=4em{
      {\CompCat}\ar[r]^{\pathcat^{*}} &
      {\Twocat(\Comp)}\ar[r]^{\squares} &
      {\Double(\Comp)}\ar[r]^{\efunc_h} &
      {\Double(\Comp)}}
  \end{equation*}
  and observe that on applying the functor $\efunc_h$ to the natural isomorphism 
  $\phi$ of lemma~\ref{recons.lemma}, and appealing to the equality
  $\pathecat(\pathcat(A))=\efunc_h(\pathesq(A))$ of
  observation~\ref{study.pecat.d1}, we get a natural isomorphism:
  \begin{equation}\label{theta.dc}
    \let\labelstyle=\textstyle
    \xymatrix@=3em{
      {\Comp}\ar[r]^<>(0.5){\pathcat}_{}="1"\ar[d]_{\pathecat} & 
      {\Cat(\Comp)}\ar[d]^{\pathecat} \\
      {\CompCat}\ar[r]_<>(0.5){\interdc}^{}="2" & {\Double(\Comp)}
      \ar@{=>}"1"-<0em,1em>;"2"+<0em,1em>^{\theta}_{\cong}}
  \end{equation}
  Examining the composite at the top right of this diagram we see, by
  observation~\ref{study.pecat.d1}, that if $A$ is a complicial set
  then the double category $\pathecat(\pathcat(A))$ has horizontal
  categories of arrows $\pathecat(A)$ and squares
  $\pathecat(\Delta[1]\Rightarrow A)$ and that its vertical source and
  target maps are $\pathecat(\Delta(\vertex^1_0)\Rightarrow A)$ and
  $\pathecat(\Delta(\vertex^1_1)\Rightarrow A)$ respectively. As also
  observed there, these are simply the various structures which appear
  on the right of the triangles in display~(\ref{dp.triangles}) of
  observation~\ref{interv.analogue.defn}.
  
  Taking this as our lead, we turn to the lower left composite of
  display~(\ref{theta.dc}). Consulting
  observation~\ref{interv.analogue.defn}, it is clear that if
  $\mathbb{C}$ is a complicially enriched category then the natural
  candidate for $\interv(\mathbb{C})\in\CompCat$ is the horizontal
  category of squares of $\interdc(\mathbb{C})$, which we know to be
  complicially enriched (cf.\ observation~\ref{study.pecat.d1}).
  Furthermore, by construction, the horizontal category of arrows of
  $\interdc(\mathbb{C})$ is actually $\mathbb{C}$ itself, and so the
  choice $\interv(\mathbb{C})\defeq\sq_h(\interdc(\mathbb{C}))$ also
  presents us with obvious candidates for the complicial functors
  $\arrow\iend^0_{\mathbb{C}},\iend^1_{\mathbb{C}}:
  \interv(\mathbb{C})->\mathbb{C}.$ which we can take to be the
  vertical source and target maps of $\interdc(\mathbb{C})$
  respectively. Now suppose that $\arrow f:\mathbb{C}->\mathbb{C}'.$
  is a functor in $\CompCat$ then the double functor $\arrow
  \interdc(f): \interdc(\mathbb{C})->\interdc(\mathbb{C}).$ acts like
  $f$ on horizontal categories of arrows and so if we define
  $\interv(f)$ to be its action on horizontal categories of squares
  then the naturality conditions for $\iend^0_{\mathbb{C}}$ and
  $\iend^1_{\mathbb{C}}$ simply reduce to the fact that the double
  functor $\interdc(f)$ preserves vertical sources and targets.
  
  Finally, having made this choice for $\interv$ and the associated
  natural transformations $\iend^0$ and $\iend^1$, the isomorphism
  $\theta$ of display~(\ref{theta.dc}) provides us with a natural
  family of isomorphisms
  \begin{equation*}
    \let\labelstyle=\textstyle
    \xymatrix@R=1ex@C=1em{
      {\pathecat(\Delta[1]\Rightarrow A)}\ar@{=}[r] &
      {\sq_h(\pathecat(\pathcat(A)))}\ar@{=>}[rrr]^{\theta_A}_{\cong} &&&
      {\sq_h(\interdc(\pathecat(A)))}\ar@{=}[r] &
      {\interv(\pathecat(A))}}
  \end{equation*}
  as required by display~(\ref{dp.square}) of
  observation~\ref{interv.analogue.defn}. Furthermore, each $\theta_A$
  is actually a double functor so, in particular, the fact that it
  preserves vertical sources and targets, when combined with our
  identification of the vertical sources and targets of
  $\pathecat(\pathcat(A))$, the definitions of
  $\iend^0_{\pathecat(A)}$ and $\iend^1_{\pathecat(A)}$ and the
  observation that $\theta$ acts as the identity on horizontal
  categories of arrows (cf.\ observation~\ref{recons.concrete}),
  establishes the triangle identities of
  display~(\ref{dp.triangles}).\qed
\end{obs}

\begin{obs}[a more explicit description of $\interv(\mathbb{C})$]
  \label{interv.analogue.expl}
  If $\mathbb{C}$ is a complicially enriched category then
  $\interv(\mathbb{C})$ has $\Sup_0(\scarr{C})$ as its complicial set
  of objects and $n$-arrows which are pasting squares of the form
  \begin{equation}\label{sq.interv.1}
    \xymatrix@=2.5em{
      {\bullet}\ar[r]^{q}_{}="1"\ar[d]_{p} &
      {\bullet}\ar[d]^{p'} \\
      {\bullet}\ar[r]_{q'}^{}="2" & {\bullet}
      \ar@{=>}"1"-<0em,1em>;"2"+<0em,1em>^{\varphi}
    }
  \end{equation}
  where $p$ and $p'$ are $n$-arrows in
  $\Sup_0(\sc{C})\subseteq_r\sc{C}$ (which are it source and target
  object in $\interv(\mathbb{C})_o=\Sup_0(\scarr{C})$), $q$ and $q'$
  are $n$-arrows in $\sc{C}$ and $\varphi$ is an $n$-dimensional
  2-cell in $\pathcat^{*}(\mathbb{C})$ (which has $\sc{C}$ as its
  underlying category of 1-cells).  Such pasting squares are thin in
  $\interv(\sc{C})_a$ precisely when $q$, $q'$ and $\varphi$ are all
  thin in their respective complicial sets, simplicial operators act
  on them point-wise and they compose by horizontal pasting.
  Furthermore the functors $\iend^0_{\sc{C}}$ and $\iend^1_{\sc{C}}$
  simply project a pasting square to its vertical source $q$ and
  target $q'$ respectively. 
\end{obs}

\subsection{A Decalage Construction on Complicially Enriched
  Categories}

Finally we may transfer the decalage construction of
subsection~\ref{dec.on.comp.subsect} to $\CompCat$ using the
endo-functor $\interv$ which we constructed in
observation~\ref{interv.analogue.obs}, as follows:

\begin{defn}[decalage on $\CompCat$]\label{dec.on.compcat.defn}
  If $\mathbb{C}$ is a complicially enriched category then define
  $\Dec(\mathbb{C})$, its decalage, to be the complicially enriched
  category obtained by taking the pullback
  \begin{equation}\label{dec.on.compcat.pb}
    \let\labelstyle=\textstyle
    \xymatrix@=3em{
      {\Dec(\mathbb{C})}\pbexcursion\ar@{u(->}[r]^-{\subseteq_r}\ar[d] &
      {\interv(\mathbb{C})}\ar[d]^{\iend^0_{\mathbb{C}}} \\
      {\Sup_0(\mathbb{C})}\ar@{u(->}[r]_-{\subseteq_r} & {\mathbb{C}}
    }
  \end{equation}
  in $\CompCat$. As a functor, this construction is the point-wise
  pullback, in the functor category $\funcat[\CompCat,\CompCat]$, of
  the natural inclusion $\overinc\subseteq_r:\Sup_0(\cdot)->
  \id_{\CompCat}.$ along the natural transformation $\arrow
  \iend^0:\interv->\id_{\CompCat}.$.  We make $\Dec$ into a copointed
  endo-functor on $\CompCat$ by defining the counit $\arrow
  \cface:\Dec->\id_{\CompCat}.$ to be the natural transformation
  obtained by composing the natural inclusion
  $\overinc\subseteq_r:\Dec->\interv.$ with the natural transformation
  $\arrow \iend^1: \interv -> \id_{\CompCat}.$.
\end{defn}

\begin{obs}[an explicit description of {$\Dec$} on {$\CompCat$}]
  \label{dec.on.compcat.expl}
  Applying the explicit description of $\interv(\sc{C})$ given in
  observation~\ref{interv.analogue.expl} and the definition of
  $\Sup_0(\sc{C})$, we see that $\Dec(\sc{C})$ is the complicial
  subcategory of $\interv(\sc{C})$ consisting of those pasting
  squares of the form depicted in display~(\ref{sq.interv.1}) for
  which the 1-cell $q$ is an identity $n$-arrow in $\mathbb{C}$. In
  other words, contracting the 1-cell $q$ in that diagram to a point,
  we may depict the $n$-arrows of $\Dec(\sc{C})$ as triangles:
  \begin{equation}\label{tri.dec.1}
    \let\labelstyle=\textstyle
    \xymatrix@C=1em@R=2.6em{
      & {\bullet}\ar[dl]_{p}^{}="1"\ar[dr]^{p'}_{}="2" & \\
      {\bullet}\ar[rr]_{q'} && {\bullet}
      \ar@{{}=>}_{\varphi} "2"+<-0.5em,-0.5em>;"1"+<0.5em,-0.5em>\\
    }
  \end{equation}
  With respect to this description, the component
  $\arrow\cface_{\sc{C}}:\Dec(\sc{C})->\sc{C}.$ of the counit simply
  maps the $n$-arrow depicted in~(\ref{tri.dec.1}) to the $n$-arrow
  $q'$ in $\mathbb{C}$.
\end{obs}

We may easily establish a close relationship between our decalage
constructions on $\Comp$ and $\CompCat$:

\begin{thm}\label{dec.on.compcat}
  The natural isomorphism $\theta$ of
  observation~\ref{interv.analogue.defn}, as constructed in
  observation~\ref{interv.analogue.obs}, restricts to give rise to a
  natural isomorphism
  \begin{equation}\label{dec.on.compcat.diag}
    \let\labelstyle=\textstyle
    \xymatrix@R=2.5em@C=3em{
      {\Comp}\ar[d]_{P}\ar[r]^{\Dec}_{}="one" &
      {\Comp}\ar[d]^{P} \\
      {\CompCat}\ar[r]_{\Dec}^{}="three" & 
      {\CompCat}
      \ar@{=>}^{\textstyle\psi}_{\cong}
        "one"-<0em,1em>; "three"+<0em,1em>
    }
  \end{equation}
  which makes the pair $\pair<P;\psi>$ into a strong
  transformation of copointed endo-functors from the decalage
  construction on $\Comp$ to that on $\CompCat$.
\end{thm}

\begin{proof}
  This is simple calculation with pullbacks, which is best summarised
  in the following cubical diagram:
  \begin{equation*}
    \let\labelstyle=\textstyle
    \xymatrix@=1.5em{
      {\pathecat(\Dec(A))}\pbexcursion
      \ar@{^{(}->}[rr]^{\pathecat(\subseteq_r)}
      \ar@{-->}[dr]^(0.6){\psi_A}_{\cong}\ar[dd] & &
      {\pathecat(\Delta[1]\Rightarrow A)}
      \ar[dr]^{\theta_A}_(0.4){\cong}
      \ar@{-}[d] & \\
      & {\Dec(\pathecat(A))}\pbexcursion
      \ar@{^{(}->}[rr]^(0.4){\subseteq_r}\ar[dd] & 
      {}\ar[d]^(0.4){\scriptstyle\pathecat(\Delta(\vertex^1_0) 
        \Rightarrow A)} &
      {\interv(\pathecat(A))}\ar[dd]^{\iend^0_{\pathecat(A)}} \\
      {\pathecat(\Sup_0(A))}
      \ar@{^{(}-}[r]\ar@{=}[dr] & {}\ar[r]_(0.2){\subseteq_r} & 
      {\pathecat(A)}\ar@{=}[dr] & \\
      & {\Sup_0({\pathecat(A)})}\ar@{^{(}->}[rr]_{\subseteq_r} & &
      {\pathecat(A)} 
      \save "1,4"+<6em,0em> *+{\pathecat(A)} 
      \ar@{<-}"1,3"_{\scriptstyle\pathecat(\Delta(\vertex^1_1)
        \Rightarrow A)}
      \ar@{<-}"2,4"^{\iend^1_{\pathecat(A)}}
      \restore
    }
  \end{equation*}
  Here the front face is simply the pullback used to define
  $\Dec(\pathecat(A))$ in display~(\ref{dec.on.compcat.diag}) of
  definition~\ref{dec.on.compcat}, the right-hand face is the first of
  the commuting triangles in display~(\ref{dp.triangles}) of
  observation~\ref{interv.analogue.defn} and the equality in the lower
  left corner is that supplied by lemma~\ref{sup.comp.e.cat}.  The
  back face is also a pullback since it is obtained by applying the
  left exact functor $\arrow\pathecat:\Comp->\CompCat.$ to the
  pullback used to define $\Dec(A)$ in
  display~(\ref{dec.on.comp.diag}) of definition~\ref{dec.on.comp}.
  
  It follows therefore, by the pullback property of the front and back
  faces, that we get the induced isomorphism $\arrow\psi_A:
  \pathecat(\Dec(A))->\Dec(\pathecat(A)).$ which is uniquely
  determined by the property that it makes the upper face commute.
  Furthermore, it is clear that the naturality of $\theta_A$
  immediately implies that $\psi_A$ is also natural in the complicial
  set $A$.
  
  Finally, consulting the definitions of the copointings of these
  decalage constructions given in definitions~\ref{dec.on.comp}
  and~\ref{dec.on.compcat}, we see that we may combine the upper
  face of our cube with the triangle to its right, which is the second
  of the triangles from display~(\ref{dp.triangles}) of
  observation~\ref{interv.analogue.defn}, to demonstrate that $\psi_A$
  satisfies the compatibility condition required of a transformation
  of copointed endo-functors.
\end{proof}

\subsection{Semi-Simplicial Reconstruction}

Using the decalage constructions of the last few subsections we may
now apply the semi-simplicial version of
lemma~\ref{nerves.and.comonad.tr} to prove an exceedingly important
result by which we may re-construct the underlying semi-simplicial
structure of a complicial set $A$ from the complicially enriched
category $\pathecat(A)$. To do this, we must first introduce a
suitable connected components functor on the category of complicially
enriched categories.

\begin{defn}[connected components of complicially enriched categories]
  \label{conn.cpt.compenr}
  If $\mathbb{C}$ is a complicially enriched category then let
  $\cpt_0(\mathbb{C})$, its {\em set of connected components}, 
  be the set obtained by forming the coequaliser
  \begin{displaymath}
    \xymatrix@R=2ex@C=8em{
      {(\scarr{C})_0}\ar@<1ex>[r]^{\textstyle \s}
      \ar@<-1ex>[r]_{\textstyle \t} & {(\scobj{C})_0} 
      \ar@{->>}[r]^{\textstyle q_{\mathbb{C}}} & {\cpt_0(\mathbb{C})}
      }
  \end{displaymath}
  in $\Set$. In other words, we start with the set of objects of
  $\mathbb{C}$ and form $\cpt_0(\mathbb{C})$ from it by identifying
  any pair of objects which occurs as the source and target of the
  same 1-arrow in $\mathbb{C}$.

  If $\arrow f:\mathbb{C}->\mathbb{C}'.$ is a complicial
  functor then we know, by definition~\ref{cat.two.sort}, that it
  preserves sources and targets, in other words its underlying
  maps $\arrow f:\scarr{C}->\scarr{C}'.$ and $\arrow
  f_o:\scobj{C}->\scobj{C}'.$ satisfy the equalities
  $f_o\circ\s=\s\circ f$ and $f_o\circ\t=\t\circ f$ and therefore that
  they induce a unique function $\arrow
  \cpt_0(f):\cpt_0(\mathbb{C})->\cpt_0(\mathbb{C}').$ satisfying the
  equation $\cpt_0(f)(q_{\mathbb{C}}(c)) =q_{\mathbb{C}'}(f_o(c))$ for
  each $c\in (C_o)_0$. This action on complicial functors clearly
  makes $\cpt_0$ into a functor $\arrow \cpt_0:\CompCat->\Set.$.
\end{defn}

\begin{obs}[relating $\cpt_0$ on $\Comp$ and $\CompCat$]
  \label{conn.cpt.rel}
  We would like to use our connected components functors $\cpt_0$ on
  $\Comp$ and $\CompCat$ in applying (the semi-simplicial version of)
  lemma~\ref{nerves.and.comonad.tr} to the strong transformation of
  copointed decalage functors provided by
  theorem~\ref{dec.on.compcat}.  However, in order to do so we need to
  be able to demonstrate that we have a family of isomorphisms
  $\cpt_0(\pathecat(A))\cong\cpt_0(A)$ which is natural in the
  complicial set $A$.  These may be constructed with the aid of the
  following diagram
  \begin{displaymath}
    \xymatrix@R=6ex@C=8em{
      {(\pathecat(A)_a)_0}\ar@<1ex>[r]^{\textstyle \s}\ar[d]_{\cong}
      \ar@<-1ex>[r]_{\textstyle \t} & {(\pathecat(A)_o)_0} 
      \ar@{->>}[r]^{\textstyle q_{\pathecat(A)}}\ar[d]_{\cong} & 
      {\cpt_0(\pathecat(A))}\ar@{-->}[d]_{\cong} \\
      {A_1}\ar@<1ex>[r]^{\textstyle A_{\vertex^1_0}}
      \ar@<-1ex>[r]_{\textstyle A_{\vertex^1_1}} & {A_0} 
      \ar@{->>}[r]^{\textstyle q_{A}} & {\cpt_0(A)}
      }
  \end{displaymath}
  in which the upper and lower horizontal ``forks'' are the defining
  coequalisers for $\cpt_0(\pathecat(A))$
  (definition~\ref{conn.cpt.compenr}) and $\cpt_0(A)$
  (definition~\ref{conn.cpt}) respectively.
  
  The vertical isomorphism in the centre is, in fact, an equality
  since we know that the complicial set of objects $\pathecat(A)_o$ is
  $\Sup_0(A)\subseteq_r A$ and that $\Sup_0(A)$ and $A$ have the same
  sets of 0-simplices so it follows that
  $(\pathecat(A)_o)_0=A_0$. To construct the vertical
  isomorphism to the left, we start by observing that every 0-arrow of
  $\pathcat(A)$ is (trivially) in its subcategory $\pathecat(A)$ and
  so, by observation~\ref{pathcat.micro}, we know that the elements of
  $(\pathecat(A)_a)_0$ bear a canonical representation as
  stratified maps $\arrow p:\Delta[1]\cong\Delta[1]\otimes
  \Delta[0]->A.$. In other words, $(\pathecat(A)_a)_0$ may be
  identified with the homset $\Strat(\Delta[1],A)$ and the vertical
  arrow on the left of our diagram is simply Yoneda's isomorphism.  It
  is now a routine matter to check that the left-hand side of our
  diagram commutes serially (cf.\ \cite{Maclane:1971:CWM}) and it
  follows that the coequaliser property of our forks ensures that we
  get the induced isomorphism indicated as a dotted vertical arrow to
  the right of the diagram. The naturality of this family of
  isomorphisms in $A$ is now a trivial consequence of their
  definition.
\end{obs}

\begin{thm}[semi-simplicial reconstruction]
  \label{semisimp.recons}
  Let $\arrow\recons:\CompCat->\SSimp.$ denote the semi-simplicial nerve
  functor $\Ndec$ associated with the decalage construction on
  $\CompCat$  (of definition~\ref{dec.on.compcat.defn}) under the
  connected components functor $\arrow \cpt_0:\CompCat->Set.$ 
  (of definition~\ref{conn.cpt.compenr}). Then we have a natural
  isomorphism
  \begin{equation}\label{recons.triangle}
    \let\labelstyle=\textstyle
    \xymatrix@C=1em@R=2.6em{
      {\Comp}\ar[rr]^<>(0.5){\pathecat}\ar[dr]_{\forget}^<>(0.4){}="1" & &
      {\CompCat}\ar[dl]^{\recons}_<>(0.4){}="2" \\
      & {\SSimp}\ar@{}"1";"2"|{\cong}}
  \end{equation}
  where $\forget$ on the left-hand side is simply the forgetful functor from
  $\Comp$ to $\SSimp$. Consequently, we often refer to $\recons$ as the {\bf
  semi-simplicial reconstruction functor}.
\end{thm}

\begin{proof} 
  Observation~\ref{conn.cpt.rel} and theorem~\ref{dec.on.compcat}
  provide us with the data which allows us to apply the
  semi-simplicial version of lemma~\ref{nerves.and.comonad.tr} to
  construct a natural isomorphism:
  \begin{equation*}
    \let\labelstyle=\textstyle
    \xymatrix@C=1em@R=2.6em{
      {\Comp}\ar[rr]^<>(0.5){\pathecat}\ar[dr]_{\Ndec}^<>(0.4){}="1" & &
      {\CompCat}\ar[dl]^{\Ndec}_<>(0.4){}="2" \\
      & {\SSimp}\ar@{}"1";"2"|{\cong}}
  \end{equation*}
  However, by definition the functor on the right-hand diagonal is
  simply our semi-simplicial reconstruction functor
  $\arrow\recons:\CompCat->\SSimp.$ and by
  observation~\ref{st.of.cp.endo} we know that to one on the left-hand
  diagonal is (isomorphic to) the forgetful functor
  $\arrow\forget:\Comp->\SSimp.$.  In other words, the natural
  isomorphism thus constructed may be composed with the isomorphism
  $\arrow\forget\cong\Ndec:\Comp->\SSimp.$ to provide the natural
  isomorphism of display~(\ref{recons.triangle}) in the statement of
  the theorem as required.
\end{proof}

We now proceed to establish some important properties of
semi-simplicial reconstruction:

\begin{lemma} 
  \label{recons.pres.limits}
  The decalage functor $\arrow\Dec:\CompCat->\CompCat.$
  preserves all (small) limits. Furthermore, its composite with the
  connected components functor $\arrow\cpt_0\circ\Dec:\CompCat->\Set.$
  simply maps each complicially enriched category $\mathbb{C}$ to a
  set which is naturally isomorphic to its set of objects
  $(\scobj{C})_0$ and so this functor also preserves all (small)
  limits.
  
  It follows, therefore, that the semi-simplicial reconstruction
  functor $\arrow\recons:\CompCat->\SSimp.$ also preserves all limits
  and that, as a consequence, it preserves monomorphisms.
\end{lemma}

\begin{proof}
  The limit of a diagram in $\CompCat$ is formed by forgetting the
  category structures on its vertices, calculating the limit of the
  resulting diagram in $\Comp$ (which is a reflective full subcategory
  of $\Strat$ and thus closed in there under all limits) and then
  defining the category structure on the consequent complicial set
  ``point-wise''. It is therefore a matter of routine verification,
  using the explicit description of $\Dec(\sc{C})$ given in
  observations~\ref{interv.analogue.expl}
  and~\ref{dec.on.compcat.expl}, to show that $\Dec$ preserves limits
  that are constructed in this way.
  
  To prove the second of these results, we fix a complicially enriched
  category $\mathbb{C}$ and exploit the explicit description of
  $\Dec(\mathbb{C})$ given in observation~\ref{dec.on.compcat.expl} to
  construct the following diagram of sets and functions:
  \begin{equation}\label{cpt.dec.compcat}
    \xymatrix@R=2ex@C=8em{
      {(\Dec(\mathbb{C})_a)_0}\ar[r]|{\textstyle \s}
      \ar@<-2ex>[r]_{\textstyle \t} & 
      {(\Dec(\mathbb{C})_o)_0} 
      \ar@<-2ex>[l]_{\textstyle g}
      \ar@<-1ex>[r]_{\textstyle q} & 
      {(\scobj{C})_0}\ar@<-1ex>[l]_{\textstyle f}
    }
  \end{equation}
  Here the maps $\s$ and $\t$ are the source and target
  maps of $\Dec(\mathbb{C})$ which simply take an arrow depicted as a
  triangle in display~(\ref{tri.dec.1}) to the 0-arrows $p$ and $p'$
  respectively. To construct $q$, $f$ and $g$, notice first that, by
  construction, we have $\Dec(\mathbb{C})_o = \Sup_0(\mathbb{C}_a)$
  and we know that $\Sup_0(\mathbb{C}_a)$ and $\mathbb{C}_a$ have the
  same sets of 0-simplices so we have that
  $(\Dec(\mathbb{C})_o)_0=(\mathbb{C}_a)_0$. It follows that
  we may take $q$ and $f$ to be the source and identity functions
  {\em of the category $\mathbb{C}$\/} respectively. Finally, the
  action of $g$ is best described in diagrammatic terms:
  \begin{equation*}
    \xymatrix@C=1em@R=2.6em {
      {\bullet}\save []+<0.5em,1.5em>*+{(\Dec(\mathbb{C})_o)_0}="1"\restore
      \ar[d]_{p}^{} & {\mkern80mu} & & 
      {\bullet}\save []+<0.5em,1.5em>*+{(\Dec(\mathbb{C})_a)_0}="2"\restore
      \ar@{=}[dl]_{\id}^(0.6){}="3"\ar[dr]^{p}_(0.6){}="4" & \\
      {\bullet} & & {\bullet}\ar[rr]_{p} & & {\bullet}
      \ar@{}"3";"4"|{=}\ar@{|->}"1,1"-<-1em,1.5em>;"1,3"-<1em,1.5em>
      \ar"1";"2"^<>(0.5){\textstyle g}
    }
  \end{equation*}
  It is now a routine matter to check that in the following diagram
  \begin{equation*}
    \let\labelstyle=\textstyle
    \xymatrix@R=2.5em@C=3em{
      {(\Dec(\sc{C})_o)_0}
      \ar[r]^{g}\ar[d]_{q} &
      {(\Dec(\sc{C})_a)_0}
      \ar[r]^{\t}\ar[d]^{\s} &
      {(\Dec(\sc{C})_o)_0}\ar[d]^{q} \\
      {(\scobj{C})_0}\ar[r]_<>(0.5){f} &
      {(\Dec(\sc{C})_o)_0}\ar[r]_<>(0.5){q} &
      {(\scobj{C})_0}
    }
  \end{equation*}
  both squares commute and that its upper and lower horizontal
  composites are identities. Consequently
  display~(\ref{cpt.dec.compcat}) is a split coequaliser and, comparing
  it with the coequaliser used to define $\cpt_0$ in
  definition~\ref{conn.cpt.compenr}, we see that it demonstrates that
  $(\scobj{C})_0$ is isomorphic to $\cpt_0(\Dec(\sc{C}))$ as
  postulated. The naturality of this isomorphism follows easily as
  does the stated result regarding preservation of limits (since
  limits of complicially enriched categories are constructed
  point-wise in $\Comp$).
  
  Finally, returning to the construction of
  $\arrow\recons:\CompCat->\SSimp.$, as exposed in
  observation~\ref{comonad.lact} and
  definition~\ref{defn.comonad.nerve}, we see that for each
  $n\in\mathbb{N}$ the functor $\recons({-})_n$ which applies
  $\recons$ and then extracts the set of $n$-simplices of the result
  is, by definition, equal to $(\cpt_0\circ\Dec)\circ\Dec^n$. However,
  by the two preservation results we have already established it
  follows that this latter functor preserves all limits.  Therefore,
  quantifying over $n\in\mathbb{N}$ and appealing to the fact that
  limits are constructed point-wise (dimension by dimension) in
  $\SSimp$, we see that $\recons$ also preserves all limits.
  
  Finally the comment regarding preservation of monomorphisms by
  $\recons$ is a standard categorical result which applies to all
  pullback preserving functors (see~\cite{Maclane:1971:CWM}).
\end{proof} 

\begin{thm}
  The semi-simplicial reconstruction functor
  $\arrow\recons:\CompCat->\SSimp.$ is faithful.
\end{thm}

\begin{proof} 
  This result hinges on the observation that if $\sc{C}$ is a
  complicially enriched category then the path category
  $\pathecat(\scarr{C})$ constructed from its complicial set of arrows
  is actually (isomorphic to) a regular subcategory of the
  complicially enriched decalage category $\Dec(\sc{C})$.  In
  particular it is easily seen, directly from the definition of
  $\pathcat^{*}(\sc{C})$ in lemma~\ref{recons.lemma} and the explicit
  description of $\Dec(\sc{C})$ given in
  observations~\ref{interv.analogue.expl}
  and~\ref{dec.on.compcat.expl}, that $\pathecat(\scarr{C})$ is
  naturally isomorphic to the subcategory of $\Dec(\sc{C})$ of those
  triangles of the form depicted in display~(\ref{tri.dec.1}) for
  which $q'$ is actually an identity arrow in $\sc{C}$.
 
  So if $\sc{C}$ is a complicially enriched category then we may apply
  the semi-simplicial reconstruction functor $\recons$ to the regular
  inclusion $\overinc:\pathecat(\scarr{C})-> \Dec(\mathbb{C}).$ and
  appeal to the monomorphism preservation result of
  theorem~\ref{recons.pres.limits} to demonstrate that we thereby
  obtain a monomorphism $\overinc:\recons(\pathecat(\scarr{C}))->
  \recons(\Dec(\sc{C})).$ in $\SSimp$. Now, by
  theorem~\ref{semisimp.recons} we have an isomorphism
  $\forget(\scarr{C})\cong\recons(\pathecat( \scarr{C}))$ and
  furthermore, on consulting the construction of $\recons$ as
  described in observation~\ref{comonad.lact} and
  definition~\ref{defn.comonad.nerve}, we easily see that we also have
  $\recons(\Dec(\sc{C}))=\Dec(\recons(\sc{C}))$, where the functor
  $\Dec$ on the right-hand side is simply the usual decalage functor
  on $\SSimp$. Of course all of these inclusions, isomorphisms and
  equalities are natural in the complicially enriched category
  $\sc{C}$ and on composing them we get a family of monomorphisms
  $\inc m_{\scarr{C}}:\forget(\scarr{C})->
  \Dec(\recons(\mathbb{C})).$ in $\SSimp$ which is also natural in
  $\sc{C}\in\CompCat$.

  We now come to the crux of our argument, so suppose that $\arrow
  f,g:\sc{C}->\sc{C}'.$ are two functors in $\CompCat$ and observe
  that the naturality of the family of monomorphisms constructed in
  the previous paragraph provides us with a serially commutative
  diagram:
  \begin{equation*}
    \let\labelstyle=\textstyle
    \xymatrix@R=3em@C=5em{
      {\forget(\scarr{C})} \ar[d]<-1ex>_{\forget(f)}
      \ar[d]<1ex>^{\forget(g)}\ar@{^{(}->}[r]^<>(0.5){m_{\sc{C}}} &
      {\Dec(\recons(\sc{C}))}
      \ar[d]<-1ex>_{\Dec(\recons(f))}
      \ar[d]<1ex>^{\Dec(\recons(g))} \\
      {\forget(\scarr{C}')} \ar@{^{(}->}[r]_<>(0.5){m_{\sc{C}'}} &
      {\Dec(\recons(\sc{C}'))}}
  \end{equation*}
  So to prove that $\recons$ is faithful suppose that
  $\recons(f)=\recons(g)$ and observe that then the vertical arrows on
  the right of this square are equal and so it follows, by serial
  commutativity, that the composites $m_{\sc{C}'}\circ\forget(f)$ and
  $m_{\sc{C}'}\circ\forget(g)$ in the lower left of the square must
  also be equal. However $m_{\sc{C}'}$ is a monomorphism so it follows
  that $\forget(f)=\forget(g)$ and we know that
  $\arrow\forget:\Comp->\SSimp.$ is faithful, since it only forgets
  about degeneracy actions and thinness without discarding any
  simplices, so it follows therefore that $f$ and $g$ coincide as maps
  on complicial sets of arrows and thus are identical functors as
  required.
\end{proof}

Finally we may apply this last result to establish the theorem toward
which we have been working throughout this entire section. Later on we
shall see that this result leads directly to a proof that the
\inf-categorical nerve functor constructed by Street in
\cite{Street:1987:Oriental} actually provides us with an equivalence
between the categories of \inf-categories and complicial sets.
However, before stating and proving this pivotal result we pause to
recall a simple categorical result which will be used in its proof and
again in the next section:

\begin{lemma}\label{ff.comp}
  Suppose that $\arrow\ffunc:\mathcal{C}->\mathcal{D}.$ and
  $\arrow\gfunc:\mathcal{D}->\mathcal{E}.$ is a composable pair of
  functors then
  \begin{enumerate}[(i)]
  \item\label{ff.comp.1} if the composite $G\circ F$ is fully faithful
    and $G$ itself is faithful then $F$ is also fully faithful, and
  \item\label{ff.comp.2} if the composite $G\circ F$ is an equivalence
    and $G$ itself is fully faithful then $F$ is also an equivalence.
  \end{enumerate}
\end{lemma}

\begin{proof} This is an entirely trivial categorical result which we
    leave up to the reader to verify.
\end{proof}

\begin{thm}\label{ff.pathecat} 
  The path category functor $\arrow\pathecat:\Comp->\CompCat.$ is
  fully faithful.
\end{thm}

\begin{proof}
  To prove this result we need to work a little harder to reconstruct
  the thinness and degeneracy information built into a complicial set
  $A$ from its path category $\pathecat(A)$. To do so we'll exploit
  the fully faithful representation of well-tempered stratified sets
  as filtered semi-simplicial sets which we studied in
  lemma~\ref{well.tempered.rep}. Our argument can be summed up in the
  following diagram
  \begin{equation*}
    \let\labelstyle=\textstyle
    \xymatrix@C=1.5em@R=3em{
      {\Comp}\ar[rr]^{\pathecat}\ar@{}[drr]|{=}\ar[d]_{\Sup_{\bullet}(\cdot)} &&
      {\CompCat}\ar[d]^{\Sup_{\bullet}(\cdot)} \\
      {\funcat[\mathbb{N},\Comp]}
      \ar[rr]^{\funcat[\mathbb{N},\pathecat]}
      \ar[dr]_{\funcat[\mathbb{N},\forget]}^(0.4){}="1" &&
      {\funcat[\mathbb{N},\CompCat]}
      \ar[dl]^{\funcat[\mathbb{N},\recons]}_(0.4){}="2" \\
      & {\funcat[\mathbb{N},\SSimp]} &
      \ar@{}"1";"2"|{\cong}}
  \end{equation*}
  in which the lower triangle is obtained by applying the 2-functor
  $\funcat[\mathbb{N},{-}]$ to the triangle in
  display~(\ref{recons.triangle}) of theorem~\ref{semisimp.recons}.
  The left-hand vertical functor labelled $\Sup_{\bullet}(\cdot)$ maps
  a complicial set $A$ to its filtered family of complicial
  superstructures $\Sup_0(A)\subseteq_r\Sup_1(A)\subseteq_r...
  \subseteq \Sup_n(A) \subseteq_r...$. It follows that the composite
  of the functors down the left hand side is simply the fully faithful
  representation of lemma~\ref{well.tempered.rep} restricted to the
  category $\Comp$ which is a subcategory of $\Strat_W$ by
  lemma~\ref{lemma.predegen}.
  
  The right hand vertical, also called $\Sup_{\bullet}(\cdot)$, maps a
  complicially enriched category $\sc{C}$ to the filtered family of
  superstructures $\Sup_0(\sc{C}) \subseteq_r\Sup_1(\sc{C})
  \subseteq_r...\subseteq_r\Sup_n(\sc{C})\subseteq_r...$ in $\CompCat$
  as defined in observation~\ref{super.compl.enrich}.  Now we may
  immediately recast the result of lemma~\ref{sup.comp.e.cat}, which
  showed that $\pathecat(\Sup_n(A))=\Sup_n({\pathecat(A)})$ for
  each $n\in\mathbb{N}$ and $A\in\Comp$, to demonstrate the
  commutativity of the upper square in our diagram.
  
  Furthermore, suppose that we are given two functors $\arrow
  f,g:\sc{C}->\sc{C}'.$ in $\CompCat$ then, by definition, we have
  $\Sup_{\bullet}(f)= \Sup_{\bullet}(g)$ if and only if $f$ and $g$
  agree when restricted to the superstructure
  $\Sup_n(\sc{C})\subseteq_r\sc{C}$ for each $n\in\mathbb{N}$. However
  it is clear that any complicially enriched category $\mathbb{C}$ is
  equal to the union of its superstructures and so it follows that $f$
  and $g$ must be equal as functors on $\mathbb{C}$.  In other words,
  the functor $\Sup_{\bullet}(\cdot)$ on the right-hand side of our diagram
  is faithful as is the diagonal functor $\funcat[\mathbb{N},\recons]$
  on that same side, since it acts by applying the faithful
  semi-simplicial reconstruction functor $\recons$ point-wise to
  filtered complicial sets and maps. 
  
  Now observe that our diagram provides us with an isomorphism between
  the fully faithful composite
  $\funcat[\mathbb{N},\forget]\circ\Sup_{\bullet}(\cdot)$ on its
  left-hand side and the composite
  $\funcat[\mathbb{N},\recons]\circ\Sup_{\bullet}(\cdot)\circ\pathecat$
  on its right-hand side, which is thus also fully faithful.
  Furthermore we have shown that $\funcat[\mathbb{N}, \recons]
  \circ\Sup_{\bullet}(\cdot)$ is a composite of faithful functors and
  is thus itself faithful, so finally applying
  lemma~\ref{ff.comp}(\ref{ff.comp.1}) we may infer that $\pathecat$
  is fully faithful as postulated.
\end{proof}



\section{Street's \inf-Categorical Nerve Construction}
\label{street.nerve.sec}

\subsection{Parity Complexes}

Street's nerve construction \cite{Street:1987:Oriental} proceeds by
first providing an explicit construction of the free \inf-category
whose generators are the faces of a simplex and whose relations are
suitably re-interpreted and oriented versions of the face relations.
Later he introduced structures called {\em parity
  complexes\/}~\cite{Street:1991:Parity} in order to generalise this
work to encompass cubes and a range of other polytopes which may be
obtained as products or joins of simplices and globs (oriented
globes).

We will need to understand and calculate with these free,
geometrically derived \inf-categories in the sequel and so we briefly
review the important parts of the theory of parity complexes. For a
more detailed analysis of the combinatorics of these structures we
refer the reader to \cite{Street:1991:Parity} and
\cite{Street:1994:Parity}.

\begin{defn}[pre-parity complexes]\label{pre-parity defn}
  A {\em pre-parity complex\/} is a graded set
  $C=\bigcup_{n=0}^{\infty} C_n$ that comes equipped with a pair of
  operations which map each element $x\in C_n$ (for $n>0$) to disjoint
  non-empty finite subsets $x^{-}, x^{+}\subseteq C_{n-1}$. We say
  that the elements of $C_n$ are $n$-dimensional. In order to simplify
  some calculations and definitions we will generally adopt the
  convention that if $x$ is a 0-dimensional element then we take
  $x^{-}$ and $x^{+}$ to be the empty set.
  
  Following Steet we will generally call the elements of $x^{-}$ and
  $x^{+}$ the {\em negative faces\/} and {\em positive faces\/} of $x$
  respectively. We also reserve the symbols $\parvone$
  and $\parvtwo$ to vary over the set of parity symbols $\{{+},{-}\}$ and
  use the notation $\neg\parvone$ to denote the opposite parity to
  $\parvone$. We also say that an integer {\em $i$ is of parity
    $\parvtwo$\/} if $i$ is even and $\parvtwo$ is the parity ${+}$ or
  $i$ is odd and $\parvtwo$ is the parity ${-}$.
  
  We consider any subset $S$ of $C$ to be graded according to the
  grading of $C$ and define its {\em $n$-superstructure\/} by
  $\Sup_{n}(S)\defeq \bigcup_{m=0}^n S_m$. We will also tend to use
  the notation $S^{\neg n}$ to denote the subset obtained by omitting
  the $n$-dimensional elements from $S$ (that is $S^{\neg
    n}=S\smallsetminus S_n$). 
  
  We say that a subset $D$ is a sub-pre-parity complex of $C$ (denoted
  $D\subseteq_p C$) if it is closed in $C$ under face operations, that
  is if $\parvone$ is an arbitrary parity symbol then for all $s\in D$
  we have $s^\parvone\subseteq D$. The family of sub-pre-parity
  complexes of $C$ is closed under unions and intersections and it
  follows that each subset $S$ of $C$ is contained in a unique
  smallest sub-pre-parity complex $D$, called the pre-parity complex
  {\em generated\/} by $S$.
  
  When manipulating pre-parity complexes we will often have use for
  functions $f$ which map elements one pre-parity complex to {\em
    sets\/} of the elements of another (or possibly the same)
  pre-parity complex. We will depict such functions using the
  crossed arrow notation $\spanarr f:C->D.$ and if $S\subseteq C$ then
  we'll take the notation $f(S)$ to mean the union $\bigcup_{x\in S}
  f(x)$. In particular, if $\parvtwo$ is a parity symbol then under this
  notational convention we have $S^\parvtwo=\bigcup_{x\in S} x^\parvtwo$ the set
  of those elements which are $\parvtwo$-parity faces of some element of
  $S$.
  
  Street introduces the notations $S^{\mp} \defeq S^{-}\smallsetminus
  S^{+}$ and $S^{\pm}\defeq S^{+}\smallsetminus S^{-}$ for the sets of
  negative (resp. positive) faces of elements of $S$ which are not
  positive (resp. negative) faces of any element of $S$. He also
  introduces a binary perpendicularity relation $S\perp T$ which holds
  when $(S^{+}\cap T^{+})\cup (S^{-}\cap T^{-})=\emptyset$, that is to
  say when none of the $\parvtwo$-parity faces of some simplex in $T$ is
  also an $\parvtwo$-parity face of some simplex in $S$. 
  
  If $x,y\in C$ then we say that $x\perp y$ when the corresponding
  singleton sets are perpendicular. We say that a subset $S\subseteq
  C$ is {\em well-formed\/} if it has at most one 0-dimensional
  element and if whenever we have $x,y\in S$ with $x\neq y$ then
  $x\perp y$.
  
  Finally, for $x,y\in C$ we write $x < y$ if $x^{+}\cap y^{-}\neq
  \emptyset$ and $x\prec y$ if $x\in y^{-}$ or $y\in x^{+}$ then for
  $S\subseteq C$ we define $\triangle_S$ and $\filledtri_S$ to be the
  pre-orders obtained as the reflexive transitive closures of $<$ and
  $\prec$ on $S$ (respectively). Notice that $\triangle_S$ is a
  sub-order of $\filledtri_S$ and that if $T\subseteq S$ then
  $\triangle_T$ and $\filledtri_T$ are sub-orders of the restrictions
  of $\triangle_S$ and $\filledtri_S$ to $T$ (respectively). We also
  adopt the assumption that whenever order properties of a subset
  $S\subseteq C$ are referred to that it will be implicitly understood
  that the order $\triangle_S$ is intended.
\end{defn}

\begin{defn}[movement]\label{movement.defn}
  Suppose that $S$, $M$ and $P$ are subsets of a pre-parity complex
  $C$ then we say that $S$ moves $M$ to $P$ when we have
  \begin{equation*}
    P = (M\cup S^{+})\smallsetminus S^{-} \text{ and } M=(P\cup
    S^{-})\smallsetminus S^{+}
  \end{equation*}
  and we denote this relationship by $\arrow S:M->P.$.  The movement
  concept is fundamental to Street's use of parity complexes in
  studying \inf-categories. The reader may find a study of the basic
  properties of this concept in propositions 2.1-2.4 of section 2
  of~\cite{Street:1991:Parity}.
\end{defn}

\begin{defn}[parity complexes]\label{parity.defn}
  In general, we will not be that interested here in the detailed
  combinatorics involved in the definition and theory of parity
  complexes. However, for the record we will recall that a {\em parity
    complex\/} is a pre-parity complex $C$ which satisfies the
  following axioms:

  \begin{tabular}{ll}
    {\bf Axiom 1} & $x^{--}\cup x^{++} = x^{+-}\cup x^{-+}$,\\
    {\bf Axiom 2} & $x^{-}$ and $x^{+}$ are both well formed,\\
    {\bf Axiom 3(a)} & $x\triangle_C y\triangle_C x$ 
    implies $x=y$, and \\
    {\bf Axiom 3(b)} & $x \triangle_C y$, $x\in z^\parvone$ and 
    $y\in z^\parvtwo$ imply that $\parvone=\parvtwo$.\\  
  \end{tabular}

  \noindent In general, much of the theory of parity complexes follows from
  these axioms, however at various stages Street found that he needed
  to introduce auxiliary assumptions to make certain arguments work.
  Notably, in order to correct the proof of his ``excision of
  extremals'' result in \cite{Street:1994:Parity} he needed to assume
  that the ordering $\filledtri_C$ on $C$ was anti-symmetric and
  that each element $x\in C$ satisfied a certain {\em globularity
    condition}.  Furthermore, to extend his results to products of
  parity complexes (see later) he needed to assume that each factor
  complex {\em and its odd dual\/} satisfied these conditions.
  
  However, we will not concern ourselves with the details of these
  various conditions here. Suffice it to say that Street demonstrated
  that every parity complex we shall meet in the sequel satisfies
  every one of these conditions. Consequently, from here on we shall
  take it as understood that when we say that such and such is a
  parity complex we tacitly assume that it satisfies all of the
  technical conditions required to make Street's excision of extremals
  argument work.

  Indeed all of the examples discussed in~\cite{Street:1991:Parity}
  satisfy a stronger totality condition with regard to the
  $\filledtri_C$ relation which is of interest in its own right, but we
  shall not pursue that point here.
  
  Finally, notice that if $D$ is a sub-pre-parity complex of the
  parity complex $C$ then, since the orderings $\triangle_D$ and
  $\filledtri_D$ are sub-orders of (the restrictions of) $\triangle_C$
  and $\filledtri_C$ (respectively), it is clear that $D$ is also a
  parity complex.
\end{defn}

\begin{obs}[\inf-categories from parity complexes]
  \label{norient.defn}
  If $C$ is a any graded set we may define a simple \inf-category
  $\norient(C)$ with:
  \begin{itemize}
  \item cells $\pair<M;P>$ where $M$ and $P$ are finite subsets of
  $C$,
  \item $n$-source and $n$-target operations given by
    \begin{equation*}\begin{split}
      \s_n\pair<M;P> & {} = \pair<\Sup_{n}(M); M_n\cup \Sup_{n-1}(P)> \\
      \t_n\pair<M;P> & {} = \pair<\Sup_{n-1}(M)\cup P_n;\Sup_{n}(P)>
    \end{split}\end{equation*}
    for $n\geq 0$, and
  \item $n$-composition of an $n$-compatible pair given by:
    \begin{equation*}
      \pair<M;P>\comp_n\pair<N;Q> = \pair<M\cup N^{\neg n};
      P^{\neg n}\cup Q>
    \end{equation*}
  \end{itemize}
  The idea here is that $M_n$ and $P_n$ are the sets consisting of the
  $n$-dimensional elements of $C$ which we consider as comprising the
  surface elements of the $n$-dimensional source and target
  hemispheres of the cell $\pair<M;P>$. Given this interpretation, the
  meaning of the source and target maps should be clear. We form
  composites by taking unions of the sets of faces and eliding the
  common $n$-dimensional face along which the cells we are composing
  agree. Notice that a pair $\pair<M;P>$ is an $n$-cell in
  $\norient(C)$ if and only if $M_n=P_n$ and $M_m=P_m=\emptyset$ for
  all $m>n$.

  While this is, of itself, a somewhat routine and uninspiring
  construction, its interest comes from Street's observation that the
  free \inf-category (see definition~\ref{free.inf.cat} below)
  generated by the elements of the parity complex $C$ may be obtained
  as a sub-\inf-category $\oriental(C)$ of $\norient(C)$. To be
  precise, we define $\oriental(C)$ to be the subset of $\norient(C)$
  of those pairs $\pair<M;P>$ for which:
  \begin{itemize}
  \item both of $M$ and $P$ are both well-formed, non-empty subsets of
    $C$, and
  \item each of the subsets $M$ and $P$ moves $M$ to $P$, that is
    symbolically we have $\arrow M:M -> P.$ and $\arrow P:M->P.$.
  \end{itemize}
  It is easily seen that $\oriental(C)$ is closed in $\norient(C)$
  under sources and targets, however much of the detailed argument
  in~\cite{Street:1991:Parity} is devoted to demonstrating that it is
  also closed under compositions in $\norient(C)$ and is thus a
  sub-\inf-category of it.
  
  Before stating Street's result, it is worth noting that each
  $n$-dimensional element $x\in C_n$ gives rise to an inductively
  defined pair of subsets given by
  \begin{equation*}
    \begin{array}{lp{12em}}
      \mu(x)_m = \pi(x)_m = \emptyset & for $m>n$ \\
      \mu(x)_n = \pi(x)_n = \{x\} & \\
      \mu(x)_m = \mu(x)_{m+1}^{\mp} \text{ and }
      \pi(x)_m = \pi(x)_{m+1}^{\pm} & for $n> m\ge 0$
    \end{array}
  \end{equation*}
  which provide for us a cell $\atom{x}\defeq\pair<\mu(x);\pi(x)>$ in
  $\norient(C)$. In some places we adopt the convention that we may
  modify our use of the symbols $\mu$ and $\pi$ by superscripting them
  with parity symbols. In particular, if $\chi$ represents a symbol in
  the set $\{\mu,\pi\}$ then $\chi^{+}$ represents the same symbol
  and $\chi^{-}$ represents the opposite one.
  
  We say that $x$ is {\em relevant\/} if this pair is in fact a cell
  in $\oriental(C)$ and one of the basic postulates of Street's
  theory, which we are tacitly assuming, is that every element of a
  parity complex $C$ should be relevant. Cells of the form $\atom{x}$
  are called {\em atoms}.
\end{obs}

\begin{defn}[freely generated \inf-categories 
  (Street~\cite{Street:1987:Oriental})]\label{free.inf.cat} If
  $\mathbb{C}$ is an \inf-category and $G\subseteq\mathbb{C}$ is a set
  of its cells then we let $\Sup_n(G)\defeq G\cap\Sup_n(\mathbb{C})$
  and we grade $G$ by letting $G_0=\Sup_0(G)$ and
  $G_{n+1}\defeq\Sup_{n+1}(G)\smallsetminus \Sup_n(\mathbb{C})$ for
  each $n\in\mathbb{N}$.
  
  We say that $\mathbb{C}$ is {\em weakly generated\/} by $G$ when for
  each $n\in\mathbb{N}$ the set $\Sup_n(\mathbb{C})\cup G_{n+1}$
  separates \inf-functors with domain $\Sup_{n+1}(\mathbb{C})$, in the
  sense that if $\arrow f,g:\Sup_{n+1}(\mathbb{C})->\mathbb{D}.$ is a
  pair of \inf-functors then we may infer that $f=g$ whenever
  $f(c)=g(c)$ for all $c\in\Sup_n(\mathbb{C})\cup G_{n+1}$.
  
  We say that $\mathbb{C}$ is {\em generated\/} by $G$ when for each
  $n\in\mathbb{N}$ the smallest sub-\inf-category of $\mathbb{C}$
  which contains the set $\Sup_n(\mathbb{C})\cup G_{n+1}$ is the
  superstructure $\Sup_{n+1}(\mathbb{C})$ itself.
  
  We say that $\mathbb{C}$ is {\em freely generated\/} by $G$ when,
  for all \inf-categories $\mathbb{D}$, for all $n\in\mathbb{N}$, for
  all \inf-functors $\arrow f:\Sup_n(\mathbb{C}) ->\mathbb{D}.$ and
  for all functions $\arrow g:G_{n+1} ->\mathbb{D}.$ such that
  $\s_n(g(c)) = f(\s_n(c))$ and $\t_n(g(c))=f(\t_n(c))$ for each $c\in
  G_{n+1}$, there exists a unique $(n+1)$-functor $\arrow
  h:\Sup_{n+1}(\mathbb{C})->\mathbb{D}.$ whose restriction to
  $\Sup_n(\mathbb{C})$ is $f$ and whose restriction to $G_{n+1}$ is
  $g$.  Street depicts this extension property as a serially
  commutative diagram
  \begin{displaymath}
    \let\labelstyle=\scriptstyle
    
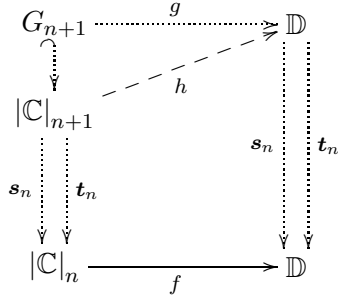
\xymatrix@R=1.5em@C=6em{
      {G_{n+1}}\ar@{.>}[r]^{g}\ar@{_{(}.>}[d] & 
      {\mathbb{D}}\ar@{.>}@<-1ex>[ddd]_{\s_n}\ar@{.>}@<1ex>[ddd]^{\t_n} \\
      {\Sup_{n+1}(\mathbb{C})}\ar@{.>}@<-1ex>[dd]_{\s_n}
      \ar@{.>}@<1ex>[dd]^{\t_n} \ar@{-->}[ur]_{h} & \\
      & \\
      {\Sup_n(\mathbb{C})}\ar[r]_{f} & {\mathbb{D}}}
  \end{displaymath}
  in which $f$ is a \inf-functor and $h$ is its unique \inf-functorial
  ``lift'' extending the action of the function $g$ on the
  $(n+1)$-dimensional generators in $G$. 
  
  Notice that none of these generation properties place any
  restriction on the set $G_0$. However, it will be convenient in what
  follows to assume that each one includes the postulate that
  $G_0=\Sup_0(\mathbb{C})$.
  
  It is clear that the free generation and generation both imply weak
  generation. A little less straightforwardly, we may also prove that
  if $G\subseteq\mathbb{C}$ freely generates $\mathbb{C}$ then it
  generates it. To prove this assume that the subset $G$ freely
  generates $\mathbb{C}$ and let $\mathbb{D}$ be the smallest
  sub-\inf-category of $\mathbb{C}$ which contains $\Sup_n(\mathbb{C})
  \cup G_{n+1}$. Applying the free generation property, we may extend
  the inclusion \inf-functor $\overinc:\Sup_n(\mathbb{C})->
  \mathbb{D}.$ using the subset inclusion function $\overinc:
  G_{n+1}->\mathbb{D}.$ to give an \inf-functor $\arrow
  h:\Sup_{n+1}(\mathbb{C})->\mathbb{D}.$ which, by definition, maps
  elements of $\Sup_n(\mathbb{C})\cup G_{n+1}$ to themselves.
  Composing $h$ with the inclusion $\mathbb{D}\subseteq
  \Sup_{n+1}(\mathbb{C})$ and applying the separation property of $G$
  as a weak generator we see that the resulting \inf-functor is equal
  to the identity on $\Sup_{n+1}(\mathbb{C})$ and it follows that the
  inclusion $\mathbb{D}\subseteq\Sup_{n+1}(\mathbb{C})$ must in fact
  be an equality as required.
\end{defn}
  
\begin{thm}[Street~\cite{Street:1991:Parity}
  and~\cite{Street:1994:Parity} and
  originally~\cite{Street:1987:Oriental}]\label{free.ocat.on.parity}
  The set $\oriental(C)$ is a sub-\inf-category of $\norient(C)$ and,
  furthermore, it is the freely generated \inf-category on the set of
  atoms $\left\langle C\right\rangle\defeq\{\atom{x}\mid x\in
  C\}$.
\end{thm}

\begin{proof} See loc.\ cit.
\end{proof}

 The following corollary is a triviality, but never-the-less it
 provides us with a useful decomposition result for free
 \inf-categories constructed from parity complexes:

 \begin{lemma}\label{widepo.parity}
   If $C$ is a parity complex and $C^{(i)}\subseteq_p C$ ($i\in I$) is
   a family of sub-parity complexes with $C=\bigcup_{i\in I} C^{(i)}$
   then the diagram of sub-\inf-category inclusions
   
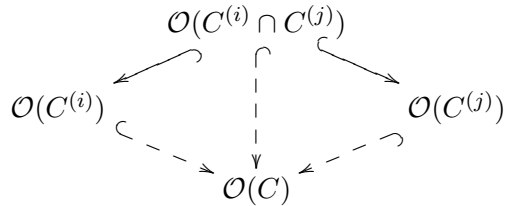
\begin{equation*}
     \xymatrix@R=0.3em@C=0.45em{
       && {\oriental(C^{(i)}\cap C^{(j)})}
       \ar@{u(->}[ddll]\ar@{u(-->}[dddd]\ar@{u(->}[ddrr] && \\
       && && \\
       {\oriental(C^{(i)})}\ar@{u(-->}[ddrr] && && 
       {\oriental(C^{(j)})}\ar@{u(-->}[ddll] \\
       && && \\
       && {\oriental(C)} && }
   \end{equation*}
   displays $\oriental(C)$ as the wide pushout of its
   sub-\inf-categories $\oriental(C^{(i)})$ in $\InfCat$.
 \end{lemma}

 \begin{proof}[Proof (sketch)]
   Suppose that we are given \inf-functors $\arrow
   f_i:\oriental(C^{(i)})->\mathbb{C}.$ forming a cocone under this
   diagram and that we've constructed a functor $\arrow f:\Sup_n(
   {\oriental(C)})->\mathbb{C}.$ which coincides with $f_i$ on
   $\Sup_n({\oriental(C^{(i)})})$ for each $i\in I$. If $x$ is an
   $(n+1)$-dimensional element of $C$ then $x$ is an element of some
   $C^{(i)}$ and we have $f(\s_n(\atom{x}))=f_i(\s_n(\atom{x}))
   =\s_n(f_i(\atom{x}))$ and $f(\t_n(\atom{x}))=
   f_i(\t_n(\atom{x}))=\t_n(f_i(\atom{x}))$. Furthermore, if we have
   another $j$ such that $x\in C^{(j)}$ then by the cocone property of
   the family $f_i$ we know that $f_i(\atom{x})=f_j(\atom{x})$, and it
   follows that we may use the free generation property of
   $\oriental(C)$ to extend $f$ uniquely to a functor on $\Sup_{n+1}(
   {\oriental(C)})$ for which $f(\atom{x})=f_i(\atom{x})$ for all
   $i\in I$ and $x\in C^{(i)}_{n+1}$. But if $f$ and $f_i$ agree on
   these atoms they must agree in the whole of
   $\Sup_{n+1}({\oriental(C^{(i)})})$. Applying this extension result
   on successive skeleta, we may construct a unique \inf-functor whose
   domain is their union $\oriental(C)$ and which restricts to $f_i$
   on each $\oriental(C^{(i)})$ as required.
 \end{proof}
 
The next two, very simple, technical results provide us with a
convenient way to construct isomorphisms with freely generated
\inf-categories. 

\begin{lemma}
  \label{gen.epi}
  Let $\epi e:\mathbb{C}->\mathbb{D}.$ be an epimorphism of
  \inf-categories and let $G\subseteq\mathbb{C}$ be a subset which
  weakly generates $\mathbb{C}$. Consider a subset
  $H\subseteq\mathbb{D}$, which we grade as in
  definition~\ref{free.inf.cat}, and suppose that $e(G_0)\subseteq
  H_0$ and $e(G_{n+1}) \smallsetminus \Sup_n(\mathbb{D})\subseteq
  H_{n+1}$ (for each $n\in\mathbb{N}$) then $H$ weakly generates
  $\mathbb{D}$.
\end{lemma}

\begin{proof}
  First note that $e$ restricts to an epimorphism $\epi
  e:\Sup_n(\mathbb{C})->\Sup_n(\mathbb{D}).$ for each $n\in\mathbb{N}$
  (by observation~\ref{ncells.ncats}). Epimorphisms between
  0-categories (sets) are all surjective and $G_0=\Sup_0(\mathbb{C})$
  so from the condition on $H_0$ in the statement it is immediate that
  we must have $H_0=\Sup_0(\mathbb{D})$. To establish the separation
  property of weak generation, consider a pair of \inf-functors
  $\arrow f,g:\Sup_{n+1}( \mathbb{D})->\mathbb{E}.$ with $f(d)=g(d)$
  for all $d\in\Sup_n(\mathbb{D})\cup H_{n+1}$. Of course, we know
  that $e(\Sup_n(\mathbb{C}))\subseteq\Sup_n(\mathbb{D})$ (since
  \inf-functors carry $n$-cells to $n$-cells) and the condition on
  $H_{n+1}$ in the statement is equivalent to
  $e(G_{n+1})\subseteq\Sup_n(\mathbb{D})\cup H_{n+1}$, therefore
  $e(\Sup_n(\mathbb{C})\cup G_{n+1})=e(\Sup_n(\mathbb{C}))\cup
  e(G_{n+1})\subseteq \Sup_n(\mathbb{D}) \cup H_{n+1}$.  It follows,
  from the fact that $f$ and $g$ coincide on $\Sup_n(\mathbb{D})\cup
  H_{n+1}$ by assumption, that $f(e(c))=g(e(c))$ for each
  $c\in\Sup_{n+1}(G)$.  However $G$ weakly generates $\mathbb{C}$ so
  we may infer that $f\circ e=g\circ e$ and use the epimorphism
  property of $\epi e:\Sup_{n+1}(\mathbb{C})->
  \Sup_{n+1}(\mathbb{D}).$ to show then that $f=g$ as required.
\end{proof}

\begin{obs}
  \label{gen.epi.obs}
  A partial converse to the last lemma is the simple observation that
  if $\arrow f:\mathbb{C}->\mathbb{D}.$ is an \inf-functor and $H$ is
  a weak generator of $\mathbb{D}$ for which $H\subseteq
  f(\mathbb{C})$ then $f$ is an epimorphism. In many cases we are
  presented with a jointly epimorphic family of \inf-functors
  $\{\arrow e_i:\mathbb{C}_i-> \mathbb{D}.\}_{i\in I}$, rather than a
  single epimorphism, and a weak generator $G^{(i)}$ for each
  $\mathbb{C}_i$. In particular, this is the case when $\mathbb{D}$ is
  the colimit of some diagram and the \inf-functors $e_i$ are the
  components of the colimiting cocone. Of course, in such a case the
  coproduct \inf-category $\coprod_{i\in I}\mathbb{C}_i$ is weakly
  generated by the subset $\coprod_{i\in I}G^{(i)}$ and we may apply
  the result of the last lemma to the induced epimorphism $\epi
  \langle e_i\rangle_{i\in I}: \coprod_{i\in I}
  \mathbb{C}_i->\mathbb{D}.$.
\end{obs}

\begin{lemma}\label{fg.isom.lemma}
  Suppose that $\mathbb{C}$ and $\mathbb{D}$ are \inf-categories and
  that we have sets of cells $G\subseteq\mathbb{C}$ and
  $H\subseteq\mathbb{D}$ such that $G$ {\bf\em freely} generates
  $\mathbb{C}$ whereas $H$ {\bf\em weakly} generates $\mathbb{D}$. Assume
  further that we have an \inf-functor $\arrow
  f:\mathbb{D}->\mathbb{C}.$ which restricts to a bijection between
  the sets $H_n$ and $G_n$ for all $n\in\mathbb{N}$ then $f$ is
  actually an {\bf\em isomorphism} of \inf-categories.
\end{lemma}

\begin{proof}
  Before proving this result, its worth making a slight methodological
  point. Its important here that we've only made the assumption that
  $H$ weakly generates $\mathbb{D}$. In all of the applications we
  have in mind $\mathbb{D}$ is obtained by taking some kind of
  quotient of a freely generated \inf-category. In such situations its
  usually relatively easy to construct functors with domain
  $\mathbb{D}$ and even to prove weak generation by a set of cells,
  however it is rarely easy to directly establish free generation.
  
  Our proof proceeds by induction on superstructures, demonstrating
  that $f$ restricts to an isomorphism between each corresponding pair
  of superstructures. For the base case, definition~\ref{free.inf.cat}
  explicitly postulates that $G_0=\Sup_0(\mathbb{C})$ and
  $H_0=\Sup_0(\mathbb{D})$ and we assumed that $f$ restricts to a
  bijection from $H_0$ to $G_0$ which thus provides a isomorphism
  between the discrete \inf-categories $\Sup_0(\mathbb{C})$ and
  $\Sup_0(\mathbb{D})$.
  
  To establish the induction step, we fix an $n\in\mathbb{N}$ and make
  the inductive hypothesis that the restriction $\arrow
  f:\Sup_n(\mathbb{D})->\Sup_n( \mathbb{C}).$ is an isomorphism of
  \inf-categories, suppose that it has inverse $\arrow g:
  \Sup_n(\mathbb{C})->\Sup_n(\mathbb{D}).$ in $\InfCat$ and let
  $\arrow f^{-1}:G_{n+1}->H_{n+1}.$ denote the postulated inverse to
  $f$ on sets of $(n+1)$-dimensional generators.  Now consider an
  arbitrary $c\in G_{n+1}$ and observe that we have
  $\s_n(f^{-1}(c))=g(f(\s_n(f^{-1}(c))))
  =g(\s_n(f(f^{-1}(c))))=g(\s_n(c))$ where the first equality holds
  since $g$ and $f$ are inverses on $n$-superstructures (by the
  inductive hypothesis), the second holds because $f$ is an
  \inf-functor (and thus commutes with $n$-sources) and the last
  follows from the fact that $f$ and $f^{-1}$ are inverses on sets of
  generators. We can argue similarly with $n$-targets and so it
  follows that we may apply the free generation property of $G$ to
  extend $g$ to an unique \inf-functor $\arrow
  h:\Sup_{n+1}(\mathbb{C})-> \Sup_{n+1}(\mathbb{D}).$ with
  $h(c)=f^{-1}(c)$ for all $c\in G_{n+1}$. By definition, we easily
  see now that $f(h(c))=c$ for all $c\in\Sup_{n}(\mathbb{C})\cup
  G_{n+1}$ and $h(f(d))=d$ for all $d\in\Sup_{n}(\mathbb{D})\cup
  H_{n+1}$, so applying the weak generation properties of $G$ and $H$
  we infer that $f\circ h=\id_{\Sup_{n+1}(\mathbb{C})}$ and $h\circ
  f=\id_{\Sup_{n+1}( \mathbb{D})}$. Thus the restriction $\arrow f:
  \Sup_{n+1}(\mathbb{D})->\Sup_{n+1}(\mathbb{C}).$ is also an
  isomorphism of \inf-categories, with inverse $h$.
  
  Finally, our inductive argument has shown that $f$ restricts to an
  isomorphism between each corresponding pair of superstructures and,
  using the fact that an \inf-category is the union of its
  superstructures, it follows that it is itself an isomorphism of
  \inf-categories.
\end{proof}

In the sequel, we'll use a few properties of the \inf-category
$\oriental(C)$ generated by a parity complex $C$. These are an
immediate consequence of Street's work in~\cite{Street:1991:Parity}
and~\cite{Street:1994:Parity} but are not explicitly remarked upon
there.

\begin{obs}\label{simple.cell.results}
  If $S$ is a finite, non-empty subset of parity complex $C$ then
  axiom 3(a) of definition~\ref{parity.defn} ensures that we may find
  a $\triangle_S$ minimal element $s$ in $S$. By minimality it follows
  that $s^{-}\cap c^{+}=\emptyset$ for all $c\in S$ and thus we know
  that $s^{-}\subseteq S^\mp$. Arguing dually, we also see that we
  have an element $t\in S$ with $t^{+}\subseteq S^\pm$. In particular,
  these observations imply that both of the sets $S^\mp$ and $S^\pm$
  are non-empty.
  
  Applying this result to a cell $\pair<M;P>$ of $\oriental(C)$ we see
  that if $M_n$ ($n>0$) is non-empty then it contains elements $s$ and
  $t$ such that $s^{-}\subseteq M_n^\mp$ and $t^{+}\subseteq M_n^\pm$.
  However, the movement property of our cell implies that
  $M_n^\mp\subseteq M_{n-1}$ and $M_n^\pm\subseteq P_{n-1}$, so it
  follows that $s^{-}\subseteq M_{n-1}$ and $t^{+}\subseteq P_{n-1}$
  and in particular that $M_{n-1}$ and $P_{n-1}$ are also non-empty.
  Of course, we may apply exactly the same argument to $P_n$ to obtain
  precisely the same result under the assumption that it is non-empty.
  
  Now, every cell $\pair<M;P>$ has a maximum $n$ for which $M_n$
  or $P_n$ is non-empty (since $M$ and $P$ are finite and non-empty)
  and by the movement property for $P_{n+1}= M_{n+1}=\emptyset$ we
  have $M_n=P_n$, it follows that $\pair<M;P>$ is a non-trivial
  $n$-cell. Applying the result of the last paragraph repeatedly, it
  follows also that $M_r$ and $P_r$ are non-empty for every $r<n$.
\end{obs}

\begin{lemma}\label{elt.in.cell} 
  Suppose that $C$ is a parity complex and that there is some $x\in
  C_n$ ($n>0$) which generates it, in the sense that $C$ is the
  smallest sub-parity complex of itself which contains $x$, then if we
  are given an element $y\in C_m$ ($m>0$) we may construct an
  $(m+1)$-cell $\pair<M;P>$ in $\oriental(C)$ with $y\in M_m$ and
  $\t_m\pair<M;P>=\t_m(\atom{x})$.
\end{lemma}

\begin{proof}
  The sub-parity complex generated by a subset $S\subseteq C$ is
  simply the smallest substructure of $C$ which is closed under the
  face operations and which contains $S$.  It follows, therefore, that
  $x\in C_n$ generates $C$ if and only if $C_n=\{x\}$ and for all
  $y\in C_m$ with $m<n$ there exists a $z\in C_{m+1}$ such that $y\in
  z^{-}\cup z^{+}$. Notice also that this also implies that
  $C_m=\emptyset$ for $m>n$.
  
  Armed with this observation, we prove our result by ``downward''
  induction on $m$. The base case $m=n$ is trivial, since in that case
  $y=x$ and we may take the atom $\atom{x}$ itself as our $(m+1)$-cell
  $\pair<M;P>$.
  
  So assume the inductive hypothesis, that our result is true at
  dimensions greater than $m$, and consider an element $y\in C_m$. We
  know, by our initial observation, that there is some $z\in C_{m+1}$
  for which $y\in z^{-}\cup z^{+}$ and we may apply our inductive
  hypothesis to $z$ to construct an $(m+2)$-cell $\pair<N;Q>$ in
  $\oriental(C)$ which has $z\in N_{m+1}$ and
  $\t_{m+1}\pair<N;Q>=\t_{m+1}(\atom{x})$. However, suppose that $y\in
  z^{+}$ then, by the movement property of $N_{m+1}$, we know that
  $Q_m=(N_m\cup N_{m+1}^{+})\smallsetminus N_{m+1}^{-}$ so either
  $y\in Q_m$ or $y\in N_{m+1}^{-}$. In the first case there is nothing
  more to do since we can take $\pair<M;P>=\t_m\pair<N;Q>$, in the
  second case we can replace $z$ by the element of $N_{m+1}^{-}$ which
  counts $y$ amongst its negative faces.
  
  Thus from now on we may assume that $y\in z^{-}$, define
  $X=\{w\in N_{m+1}\mid z\triangle_N w\}$ and observe that since $N$
  is well formed, by the well formedness property of the cell
  $\pair<N;Q>$, then so is the subset $X$.  Furthermore we have:

  \vspace{1ex}
  
  \noindent{\boldmath $X^\pm\subseteq N_{m+1}^\pm\subseteq Q_m$}, 
  certainly we know that the second of these inclusions holds, by the
  movement property of the cell $\pair<N;Q>$, so to prove the former
  suppose that $u\in X^{+}\smallsetminus N_{m+1}^\pm$ then we have
  \begin{equation*}
    X^{+}\smallsetminus N_{m+1}^\pm = 
    X^{+}\smallsetminus(N_{m+1}^{+}\smallsetminus N_{m+1}^{-}) = 
    (X^{+}\smallsetminus N_{m+1}^{+})\cup (X^{+}\cap N_{m+1}^{-}) =
    X^{+}\cap N_{m+1}^{-}
  \end{equation*}
  where the last of these equalities holds since $X\subseteq N_{m+1}$
  which implies that $X^{+}\subseteq N_{m+1}^{+}$, so it follows that
  we have a $v\in N_{m+1}$ and a $w\in X$ with $u\in w^{+}\cap v^{-}$.
  However, the fact that this latter intersection is non-empty tells
  us that $w\triangle_N v$ and we know that $z\triangle_N w$, by the
  definition of $X$ of which $w$ is an element, so by transitivity we
  get $z\triangle_N v$ and thus we have $v\in X$, from which it follows
  that $u\in v^{-}\subseteq X^{-}$. In summary we have shown that 
  $(X^{+}\smallsetminus N_{m+1}^\pm)\smallsetminus X^{-}
  =X^\pm\smallsetminus N_{m+1}^\pm$ is empty, or in other words that 
  $X^\pm\subseteq N_{m+1}^\pm$ as required.

  \vspace{1ex}
  
  \noindent{\boldmath $z^{-}\subseteq (Q_m\cup X^{-})\smallsetminus
    X^{+}$}, we know that $\triangle_N$ is reflexive so $z\in X$ and
  it follows that $z^{-}\subseteq Q_m\cup X^{-}$, therefore all we
  have to prove is that $z^{-}\cap X^{+}$ is empty. So suppose the
  converse, for a contradiction, and observe that it would imply that
  there is some $w\in X$ with $z^{-}\cap w^{+}\neq\emptyset$, in which
  case we have $w\triangle_N z$, by definition, and $w\neq z$, since
  all elements have disjoint sets of negative and positive faces.
  However, by the definition of $X$, of which $w$ is an element, we
  also know that $z\triangle_N w$ and when combined these three facts
  contradict axiom 3(a) of definition~\ref{parity.defn} as required.

  \vspace{1ex}
  
  Now, the first of these results implies that we may apply lemma~3.2
  of~\cite{Street:1991:Parity} to the $m$-cell $\t_m\pair<N;Q>
  =\pair<\Sup_{m-1}(N)\cup Q_m;\Sup_{m}(Q)>$ and the subset $X$.
  Part~(c) of that lemma allows us to construct an $(m+1)$-cell
  $\pair<M;P>\defeq\pair<\Sup_{m-1}(N)\cup Y\cup X;\Sup_{m}(Q)\cup X>
  \in\oriental(C)$, where $Y$ is the subset $(Q_m\cup
  X^{-})\smallsetminus X^{+}$ discussed in the second of the results
  above. However, that analysis tells us that we have $z^{-}\subseteq
  Y$ and, since we originally selected $z$ so that $y\in z^{-}$, it
  follows that $y\in Y=M_m$. Finally, by construction, it is clear that
  $\t_m\pair<M;P>=\t_m\pair<N;Q>$ and since we originally selected 
  $\pair<N;Q>$ to have $\t_{m+1}\pair<N;Q>=\t_{m+1}(\atom{x})$ it
  follows that $\t_m\pair<M;P>=\t_m(\atom{x})$ as required.
\end{proof}

\begin{cor}\label{top.diml.unique} 
  Suppose that $C$ is the parity complex of lemma~\ref{elt.in.cell}
  then the atom $\atom{x}$ is the {\bf\em unique} non-trivial $n$-cell
  of the $n$-category $\oriental(C)$.
\end{cor}

\begin{proof}
  Suppose that $\pair<M;P>$ is a non-trivial $n$-cell of
  $\oriental(C)$ then since $x$ is the unique
  $n$-dimensional element of $C$ we must have
  $M_n=P_n=\{x\}$. Now, for a contradiction, assume that
  $\pair<M;P>\neq\atom{x}$, and observe that this assumption
  implies that there is some $m<n$ and some cell which is not an
  $m$-cell but which is $m$-composable with $\atom{x}$. 
  
  To establish this fact, observe that we may apply theorem~4.1
  of~\cite{Street:1991:Parity} to construct an $m<n$ and two $n$-cells
  $\pair<N;Q>$ and $\pair<L;R>$ which are not $m$-cells and for which
  $\pair<M;P>=\pair<N;Q>\comp_m \pair<L;R>$. Notice that one of these
  new cells must be a non-trivial $n$-cell, since $\pair<M;P>$ itself
  is, and indeed we may assume without loss of generality that this is
  actually equal to the atom $\atom{x}$.  Were this not the case
  then we could simply iterate our argument and use the rank of cells,
  as in theorem~4.2 of~\cite{Street:1991:Parity}, to show that this
  process must eventually terminate at the atom $\atom{x}$.
  
  Without loss of generality we may assume that the composite we
  constructed is $\atom{x}\comp_m\pair<L;R>$, were it the reverse
  composite we would simply argue dually.  We know that $\pair<L;R>$
  is not an $m$-cell and so, by observation~\ref{simple.cell.results},
  we may infer that there is some $y\in R_{m+1}$ such that
  $y^{+}\subseteq R_m$.  However $R_m=\mu(x)_m$, since $\atom{x}$ and
  $\pair<L;R>$ are $m$-composable, so we may infer that
  $y^{+}\subseteq\mu(x)_m$.
  
  Applying lemma~\ref{elt.in.cell} we obtain an $(m+2)$-cell
  $\pair<M;P>$ which has $y\in M_{m+1}$ and which has
  $\t_{m+1}\pair<M;P>=\t_{m+1}(\atom{x})$. By the movement
  property for this cell we know that $M_{m+1}$ moves $M_m=\mu(x)_m$
  to $P_m=\pi(x)_m$ so in particular $\mu(x)_m = (\pi(x)_m\cup
  M_{m+1}^{-}) \smallsetminus M_{m+1}^{+}$, however we have $y\in
  M_{m+1}$ thus $y^{+}\subseteq M_{m+1}^{+}$ and it follows that
  $\mu(x)_m\cap y^{+} =((\pi(x)_m\cup M_{m+1}^{-})\smallsetminus
  M_{m+1}^{+})\cap y^{+} =\emptyset$. This, provides us with the
  desired contradiction, since in the last paragraph we constructed
  $y$ so that $y^{+}\subseteq \mu(x)_m$ and we know that $y^{+}$
  itself is non-empty.
\end{proof}
  
  We may represent standard simplices, cubes and product polytopes as
  parity complexes and thereby construct appropriate \inf-categorical
  models of such structures.
  
\begin{obs}[\inf-categorical simplex models]\label{simplex.models}
  Following Street in~\cite{Street:1987:Oriental}, we may define the
  standard \inf-simplex $\OSimp$ to be the graded set for which
  $\OSimp_r$ is the set of $(r+1)$-element (indexed) subsets
  $\svec{v}=\{v_0 < v_1 < ... < v_r\}\subseteq\mathbb{N}$ which we
  often simply write as $v_0...v_r$.  In essence each $v_i$ is
  thought of as a distinct vertex of an infinite dimensional simplex
  and $v_0...v_r$ is thought of as its $r$-dimensional face
  spanning the listed vertices.  
  
  Alternatively, we may represent an element $\svec{v}\in\OSimp_r$ as
  the strictly order preserving map $\face_{\svec{v}}$ from $[r]$ to
  $\mathbb{N}$ which maps $i$ to $v_i$. In general, we will identify
  these representations, passing between them as the mood takes us and
  without comment.
  
  Of course, if $\arrow \alpha:[p]->[r].$ is an arbitrary face
  operator then we may pre-compose it with
  $\face_{\svec{v}}\in\OSimp_r$ to give an element
  $\face_{\svec{v}}\circ\alpha\in\OSimp_p$. Using our vectorial
  notation we write $\svec{v}\cdot\alpha$ or $v_{\alpha(0)}
  v_{\alpha(1)}...  v_{\alpha(p)}$ for the subset corresponding to
  the map $\face_{\svec{v}}\circ\alpha$, and this action makes
  $\OSimp$ into a semi-simplicial set. In terms of our vertex-wise
  notation we will sometimes use the traditional algebraic
  topologist's device $v_0...\hat{v}_i...v_r$ where this denotes the
  face of $\svec{v}$ obtained by omitting the vertex $v_i$, in other
  words it is an alternative way of writing $\svec{v}\cdot\face^r_i$.
  
  Now we may use this semi-simplicial structure to define a pre-parity
  complex structure on $\OSimp$ by:
  \begin{equation*}
    \svec{v}^\parvtwo = \{ \svec{v}\cdot\face^r_i \mid i\in[r] 
    \text{ and $i$ is of parity $\parvtwo$}\}
  \end{equation*}
  As demonstrated in \cite{Street:1991:Parity} and
  \cite{Street:1994:Parity}, this actually makes $\OSimp$ into a
  parity complex satisfying all of the conditions required for us to
  be able to apply theorem~\ref{free.ocat.on.parity} and construct a
  free \inf-category $\oriental(\OSimp)$.
  
  We obtain an \inf-categorical model of the $n$-simplex by
  restricting our attention to the sub-parity complex $\OSimp[n]$ of
  $\OSimp$ consisting of those elements $\svec{v}$ for which each
  $v_i$ is actually an element of the ordered set $[n]$. In other
  words, this is the same as saying that the corresponding map
  $\face_{\svec{v}}$ restricts to a face operator with codomain $[n]$.
  When manipulating the elements of $\OSimp[n]$ as face operators well
  often apply a superscripted $n$, as in $\face^n_{\svec{v}}$, to
  remind ourselves of this restricted codomain. The resulting
  sub-$n$-category $\oriental(\OSimp[n])\subseteq\oriental(\OSimp)$
  was dubbed the $n^{\text{th}}$ oriental in
  \cite{Street:1987:Oriental}.
  
  Notice that if $\svec{v}$ is an $r$-dimensional element of
  $\OSimp[n]$ then $\svec{v}^{-}\cup\svec{v}^{+}$ is the set of all
  $(r-1)$-dimensional subsets of $\svec{v}$. Applying this result
  recursively, it follows therefore that the sub-parity complex of
  $\OSimp[n]$ generated by $\svec{v}$ consists of the set of all those
  elements $\svec{w}\in\OSimp[n]$ which have $\svec{w}\subseteq
  \svec{v}$. In particular, it follows that the $n$-dimensional
  element $01...n$ generates $\OSimp[n]$ as a parity complex. So we
  may apply corollary~\ref{top.diml.unique} to show that the atom
  $\atom{01...n}$ is the unique non-trivial $n$-cell of the
  $n$-category $\oriental(\OSimp[n])$.
\end{obs}

\begin{obs}[products of parity complexes and models of cubes]
  \label{parity.prods.and.cubes}
  The product $C\times D$ of two pre-parity complexes $C$ and $D$ is
  the graded set defined by
  \begin{equation*}
    (C\times D)_n = \bigcup_{p+q=n} C_p\times D_q
  \end{equation*}
  in which the $\parvtwo$-parity faces of an element $\pair<x;y>\in
  C_p\times D_q$ are given by
  \begin{equation*}
    \pair<x;y>^{\parvtwo} = x^{\parvtwo}\times\{y\} \cup
    \{x\}\times y^{\parvtwo(p)}
  \end{equation*}
  where $\parvtwo(p)$ is equal to $\parvtwo$ if $p$ is even and is the opposite
  parity $\neg\parvtwo$ if $p$ is odd. 
  
  In~\cite{Street:1991:Parity} and~\cite{Street:1994:Parity} it is
  also shown that if $C$ and $D$ are parity complexes which satisfy
  the fairly ubiquitous conditions alluded to in
  definition~\ref{parity.defn} then $C\times D$ is also a parity
  complex which satisfies those conditions. This implies that we can
  apply theorem~\ref{free.ocat.on.parity} to an arbitrary $n$-fold
  product $C_1\times C_2\times...\times C_n$ of such parity complexes
  to obtain a free \inf-category $\oriental(C_1\times
  C_2\times...\times C_n)$.
  
  A particular example of this is the \inf-categorical model of the
  $n$-cube, which is defined to be the free \inf-category on the
  $n$-fold product of the 1-simplex parity complex $\OSimp[1]$ with
  itself.
  
  Notice that if the element $x\in C_n$ generates $C$ and $y\in D_m$
  generates $D$ (as in lemma~\ref{elt.in.cell}) then we may show that
  $\pair<x;y>$ generates the product $C\times D$, as the reader may
  readily verify using the characterisation of generation given in the
  first paragraph of the proof of lemma~\ref{elt.in.cell}.  In
  particular we know, from observation~\ref{simplex.models}, that the
  element $01...n$ generates $\OSimp[n]$ so it follows that
  $\pair<01...n;01...m>$ generates the product
  $\OSimp[n]\times\OSimp[m]$ and thus that the atom
  $\atomp<01...n;01...m>$ is the unique non-trivial
  $(n+m)$-cell of the $(n+m)$-category
  $\oriental(\OSimp[n]\times\OSimp[m])$ (by
  corollary~\ref{top.diml.unique}).
\end{obs}

\begin{obs}[defining functors between free \inf-categories on parity
  complexes] \label{norient.funct.cons} In the sequel we will need to
  construct a number of \inf-functors between \inf-categories of the
  form $\oriental(C)$. To this end, consider a pair of graded sets $C$
  and $D$ and a function $\spanarr f: C-> D.$ which maps elements of
  $C$ to finite subsets of $D$ and respects dimensions in the sense
  that for all $x\in C_n$ we have $f(x)\subseteq D_n$. We call such a
  function a {\em graded set morphism}.
  
  As usual, if $S\subseteq C$ then we'll use $f(S)$ to denote the
  union $\bigcup_{x\in D} f(x)$ and observe that since $f$ respects
  dimensions we have $f(S_n)=f(S)_n$, $f(\Sup_{n}(S))=\Sup_n({f(S)})$ and
  $f(S^{\neg n})=f(S)^{\neg n}$ for each $n\in\mathbb{N}$.
  Furthermore it is clear that if $S$ is finite then so is $f(S)$ and
  that if $T\subseteq C$ is another such subset then we have $f(S\cup
  T)=f(S)\cup f(T)$.
  
  Notice that these are precisely the operations and conditions which
  are used to define the \inf-category $\norient(C)$. So since they
  are all preserved by $f$ it follows easily that the function which
  maps a pair $\pair<M;P>$ in $\norient(C)$ to the pair
  $\pair<f(M);f(P)>$ in $\norient(D)$ is an \inf-functor, for which we
  shall use the notation $\norient(f)$. 
  
  Indeed we may construct a category $\GrSet$ with objects that are
  graded sets and arrows which are graded set morphisms.
  Composites and identities in this category are defined in the
  obvious way, that is $\id_C(x)=\{x\}$ for a graded set $C$ and
  $g\circ f(x)=g(f(x))=\bigcup_{y\in f(x)} g(y)$ ($\forall x\in C$)
  for a composable pair $\spanarr f:C->D.$ and $\spanarr g:D->E.$ .
  Now it is clear that the construction of $\norient(f)$ given in the
  last paragraph makes $\norient$ into a functor from this category
  $\GrSet$ to the category of \inf-categories $\InfCat$.
  
  In the case where $C$ and $D$ are parity complexes, we would now
  like to determine some easily verifiable conditions under which the
  induced \inf-functor $\arrow\norient(f):\norient(C)-> \norient(D).$
  restricts to give an \inf-functor $\arrow\oriental(f):
  \oriental(C)->\oriental(D).$. A simple, but adequate,
  characterisation of those functions for which this holds is provided
  by lemma~\ref{inf.funct.parity}.
\end{obs}

\begin{obs}\label{ind.cell.char}
  Before stating the next lemma, it is worth commenting on the fact
  that the conditions in its statement are motivated by a simple
  inductive characterisation of the $(n+1)$-cells $\pair<M;P>$ of
  $\oriental(C)$. This states that the $(n+1)$-cell
  $\pair<M;P>\in\norient(C)$ is in $\oriental(C)$ if and only if
  $\s_n\pair<M;P>$ and $\t_n\pair<M;P>$ are both in $\oriental(C)$ and
  $M_{n+1}$ (which is equal to $P_{n+1}$) is well-formed and moves
  $M_n$ to $P_n$.
  
  To prove this characterisation first observe that, since we know
  that $\oriental(C)$ is closed under sources and targets in
  $\norient(C)$, it is enough to fix an $(n+1)$-cell
  $\pair<M;P>\in\norient(C)$ whose $n$-source and $n$-target are both
  in $\oriental(C)$ and show that the residual well-formedness and
  movement condition characterises when our $(n+1)$-cell may actually
  be found to be in $\oriental(C)$.
  
  Notice that since $\pair<M;P>$ is an $(n+1)$-cell in $\norient(C)$
  we know that $M=\Sup_{n}(M)\cup M_{n+1}$ and $P=\Sup_{n}(P)\cup
  M_{n+1}$ and, furthermore, it is clear that well-formedness and
  movement are properties which respect dimension. It follows
  therefore that $M$ and $P$ are well-formed iff $\Sup_{n}(M)$,
  $\Sup_{n}(P)$ and $M_{n+1}$ are all well-formed and that they each
  move $M$ to $P$ iff $\Sup_{n}(M)$ and $\Sup_{n}(P)$ each move
  $\Sup_{n-1}(M)$ to $\Sup_{n-1}(P)$ and $M_{n+1}$ moves $M_n$ to
  $P_n$.
  
  However we have $\s_n\pair<M;P>=\pair<\Sup_{n}(M);M_n\cup
  \Sup_{n-1}(P)>\in\oriental(C)$ so it follows, by the definition of
  the cells that occur in $\oriental(C)$, that $\Sup_{n}(M)$ is
  non-empty and well-formed and that it moves $\Sup_{n-1}(M)$ to
  $\Sup_{n-1}(P)$.  Similarly, considering $\t_n\pair<M;P>=
  \pair<\Sup_{n-1}(M)\cup P_n; \Sup_{n}(P)>\in\oriental(C)$ we see,
  dually, that $\Sup_{n}(P)$ is also non-empty and well-formed and
  furthermore that it too moves $\Sup_{n-1}(M)$ to $\Sup_{n-1}(P)$. So
  finally, excluding these observations from the list of conditions in
  the last paragraph, we see that the remaining conditions required to
  show that $\pair<M;P>$ is a cell of $\oriental(C)$ are simply those
  postulated in the statement above.\qed
\end{obs}

\begin{lemma}\label{inf.funct.parity}
  Suppose that $C$ and $D$ are parity complexes and $\spanarr f:C->D.$
  is a graded set morphism as discussed in
  observation~\ref{norient.funct.cons}. Then the \inf-functor
  $\arrow\norient(f):\norient(C)-> \norient(D).$ restricts to give an
  \inf-functor $\arrow\oriental(f):\oriental(C)-> \oriental(D).$ if
  and only if $f$ satisfies the following conditions:
  \begin{enumerate}[(a)]
  \item\label{func.par.a} for all $x\in C_0$ the subset $f(x)\subseteq
    D_0$ is a singleton, and
  \item\label{func.par.b} for all $n\in\mathbb{N}$ and $x\in C_{n+1}$
    the subset $f(x)\subseteq D_{n+1}$ is well formed and moves $f(x^{-})$
    to $f(x^{+})$ or symbolically $\arrow f(x):f(x^{-})->f(x^{+}).$.
  \end{enumerate}
  Consequently we have a subcategory $\Parity$ of $\GrSet$, whose
  objects are parity complexes and whose arrows are those graded set
  morphisms which satisfy these two conditions, which we shall often
  refer to as {\bf\em parity complex morphisms}.
  
  Clearly then, our free \inf-category construction lifts to a functor
  $\arrow\oriental:\Parity->\InfCat.$ which acts on an arrow
  $f\in\Parity$ by restricting the action of $\norient(f)$ as
  described.
\end{lemma}

\begin{proof} 
  The necessity of this equivalence follows immediately from the
  following observations:
  \begin{enumerate}[(i)]
  \item\label{func.par.obs.a} If $x\in C_0$ then
    $\norient(f)(\atom{x})=\pair<f(x);f(x)>$ which is a 0-cell of
    $\oriental(D)$ if and only if $f(x)$ is a singleton subset of
    $D_0$.
  \item\label{func.par.obs.b} If $x\in C_{n+1}$ then
    $\norient(f)(\atom{x})= \pair<f(\mu(x));f(\pi(x))>$ for which we
    have $f(\mu(x))_{n+1}=f(\pi(x))_{n+1}=f(x)$,
    $f(\mu(x))_n=f(x^{-})$ and $f(\pi(x))_n=f(x^{+})$. It follows, by
    observation~\ref{ind.cell.char}, that $\norient(f)(\atom{x})$ is
    in $\oriental(D)$ if and only if $\s_n(\norient(f)(\atom{x}))=
    \norient(f)(\s_n(\atom{x}))$ and $\t_n(\norient(f)(\atom{x}))=
    \norient(f)(\t_n(\atom{x}))$ are in $\oriental(D)$ and $f(x)$ is a
    well formed set that moves $f(x^{-})$ to $f(x^{-})$.
  \end{enumerate}
  We also prove sufficiency from these observations by inductively
  demonstrating that $\norient(f)$ restricts to an \inf-functor
  $\overarr:\Sup_n({\oriental(C)})->\oriental(D).$ for each dimension
  $n$ in turn and arguing that it therefore restricts in this way to
  the whole of $\oriental(C)$, since this is equal to the union of its
  superstructures.  To that end we assume that $f$ satisfies
  conditions (\ref{func.par.a}) and (\ref{func.par.b}) and proceed by
  induction on $n$.
  
  The base case $n=0$ of our induction is easy, since we know that the
  cells of $\Sup_0({\oriental(C)})$ are simply the atoms $\atom{x}$
  for $x\in C_0$ so by applying observation~(\ref{func.par.obs.a}) and
  condition~(\ref{func.par.a}) we see that $\norient(f)$ restricts to
  an \inf-functor $\overarr:\Sup_0({\oriental(C)})->\oriental(D).$.
  
  So assume the inductive hypothesis, that $\norient(f)$ restricts to
  an \inf-functor from $\Sup_n({\oriental(C)})$ to $\oriental(D)$, and
  consider an arbitrary atom in $\Sup_{n+1}({\oriental(C)})$. Either
  this is in $\Sup_n({\oriental(C)})$ itself, in which case
  $\norient(f)$ maps it into $\oriental(D)$ by hypothesis, or it is of
  the form $\atom{x}$ with $x\in C_{n+1}$. In the latter case the
  $n$-source and $n$-target of $\atom{x}$ are, of course, both cells
  of $\Sup_n({\oriental(C)})$ so by the inductive hypothesis it
  follows that $\norient(f)(\s_n(\atom{x}))$ and
  $\norient(f)(\t_n(\atom{x}))$ are both $n$-cells in $\oriental(D)$.
  Combining this information with condition~(\ref{func.par.b}) gives
  us enough to apply observation~(\ref{func.par.obs.b}) and infer that
  $\norient(f)(\atom{x})$ is in $\oriental(D)$.  It follows therefore
  that $\norient(f)$ maps every atom of dimension $\leq n+1$ into
  $\oriental(D)$.
  
  Now, by theorem~\ref{free.ocat.on.parity} we know that the
  sub-\inf-category $\Sup_{n+1}({\oriental(C)})\subseteq\norient(C)$
  is freely generated by the atoms of dimension $\leq n+1$ and, in
  particular, this implies that every cell of
  $\Sup_{n+1}({\oriental(C)})$ may be obtained as some finite
  \inf-categorical composite of such atoms. So since $\norient(f)$
  preserves these composites and maps all atoms of dimension $\leq
  n+1$ into $\oriental(D)$ it follows that it restricts to an
  \inf-functor $\overarr:\Sup_{n+1}({\oriental(C)})->\oriental(D).$ as
  required.
\end{proof}

\begin{obs}\label{func.movement.prop}
  As an aside, it is worth observing that the movement property of
  condition~(\ref{func.par.b}) of the last lemma may be summarised
  conveniently in the statement that the single equation
  \begin{equation*}
    (f(x^\parvtwo)\cup f(x)^{\neg\parvtwo})\smallsetminus
    f(x)^\parvtwo = f(x^{\neg\parvtwo})
  \end{equation*}
  should hold for an arbitrary parity $\parvtwo\in\{{-},{+}\}$.
\end{obs}

The following lemma and its corollary often allows us to determine the
action of an \inf-functor of the form $\oriental(f)$ on the atoms
$\atom{x}$ in its domain without requiring us to perform detailed
calculations directly on the combinatorially complex sets $\mu(x)$
and $\pi(x)$.

\begin{lemma}\label{elt.in.cell.morp}
  Suppose that $\spanarr f:C->D.$ is a morphism of parity complexes
  and that $x\in C_n$ generates the sub-parity complex $E$ of $C$. If
  we let $m\defeq\max\{i\leq n\mid f(E_i)\neq\emptyset\}$ then the cell
  $\oriental(f)(\atom{x})$ is a non-trivial $m$-cell in
  $\oriental(D)$.
\end{lemma}

\begin{proof}
  Notice that our choice of $m$ ensures that $f$ maps every element of
  $E$ of dimension $>m$ to the empty set, so it follows that
  $\oriental(f)$ maps every cell of $\oriental(E)\subseteq
  \oriental(C)$ to an $m$-cell of $\oriental(D)$. In particular since
  $x$ generates $E$ we know that the associated atom $\atom{x}$ is a
  cell in the sub-\inf-category $\oriental(E)\subseteq\oriental(C)$
  and so it follows that $\oriental(f)(\atom{x})$ is indeed an
  $m$-cell in $\oriental(D)$.  It remains to show that
  $\oriental(f)(\atom{x})$ is actually a {\em non-trivial\/} $m$-cell.
  
  To do this we select a $y\in E_m$ such that $f(y)\neq\emptyset$ and
  apply lemma~\ref{elt.in.cell} to obtain an $(m+1)$-cell $\pair<M;P>$
  in $\oriental(E)$ which has $y\in M_m$ and
  $\t_m\pair<M;P>=\t_m(\atom{x})$. Now, arguing as in the previous
  paragraph we know that $\oriental(f)\pair<M;P>$ is an $m$-cell and
  it follows that $\t_m(\oriental(f)(\atom{x})) = \oriental(f)(
  \t_m(\atom{x}))=\oriental(f)( \t_m\pair<M;P>)=\t_m(\oriental(f)
  \pair<M;P>)=\oriental(f)\pair<M;P>$. In particular this equation
  implies that $f(\pi(x)_m)=f(M_m)$ and, by construction, we know that
  $y\in M_m$ and it follows that $f(M_m)\neq\emptyset$ because
  $y$ was originally chosen so that $f(y)\neq\emptyset$. So
  $f(\pi(x)_m)\neq\emptyset$ and it follows that
  $\oriental(f)(\atom{x})=\pair<f(\mu(x));f(\pi(x))>$ has a non-empty
  set of $m$-dimensional elements and thus is a non-trivial $m$-cell
  in $\oriental(D)$ as postulated.
\end{proof}

\begin{cor}\label{morp.action.1}
  Under the conditions described in lemma~\ref{elt.in.cell.morp}, if
  $F$ is the sub-parity complex of $D$ generated by $f(E)$ and it may
  also be generated by a single element $y\in F$ then this element is
  $m$-dimensional, where $m$ is defined as in the statement of
  lemma~\ref{elt.in.cell.morp}, and $\oriental(f)(\atom{x})=\atom{y}$.
\end{cor}

\begin{proof}
  We may construct the sub-parity complex $F$ generated by a subset
  $S\subseteq D$ by closing $S$ up in $D$ under its face operations,
  which carry each element to subsets of elements of lower dimension.
  It follows therefore that the maxima of the dimensions of elements
  in $S$ and in $D$ must be equal. In particular, in
  lemma~\ref{elt.in.cell.morp} we defined $m$ to be the maximum
  dimension of elements in $f(E)$, so we may infer that the elements
  of $F$ have maximum dimension $F$. Furthermore, since the single
  element $y$ also generates $F$ it follows the dimension of $y$ must
  also be $m$.
  
  Now, applying lemma~\ref{top.diml.unique} we see that the atom
  $\atom{y}$ is the unique non-trivial $m$-cell of the
  sub-\inf-category $\oriental(F)$ of $\oriental(D)$. However,
  applying lemma~\ref{elt.in.cell.morp} we know that
  $\oriental(f)(\atom{x})$ is also a non-trivial $m$-cell in
  $\oriental(F)$ and so it follows that
  $\oriental(f)(\atom{x})=\atom{y}$ as required.
\end{proof}

\begin{obs}[the product of parity complexes as a bifunctor]
  \label{parity.prod.func}
  In fact, we may extend the product of parity complexes to a
  bifunctor from $\Parity\times\Parity$ to $\Parity$ by defining the
  product $\spanarr f\times g:C\times D-> E\times F.$ of a pair of
  arrows $\spanarr f:C->E.$ and $\spanarr g:D->F.$ in $\Parity$ to be
  the function given by $(f\times g)\pair<x;y>=f(x)\times g(y)$.
  Notice that condition~(\ref{func.par.a}) of
  lemma~\ref{inf.funct.parity} holds trivially for this definition
  and, indeed, it is only slightly more difficult to prove that it
  also satisfies condition~(\ref{func.par.b}). To do this we consider
  an arbitrary element $\pair<x;y>\in C_p\times D_q$ and show that:
  
  \vspace{1ex}

  \noindent {\bf\boldmath $(f\times g)\pair<x;y>$ is well formed:}
  notice first that we know that $f(x)\subseteq E_p$ and
  $g(y)\subseteq F_q$ are well formed, since $f$ and $g$ are arrows of
  $\Parity$, and suppose that $\pair<u;v>$ and $\pair<w;z>$ are
  arbitrary but distinct elements of the subset $(f\times
  g)\pair<x;y>=f(x)\times g(y)$ of $E_p\times F_q$. We know that
  $\pair<u;v>^\parvtwo = u^\parvtwo\times\{v\}\cup\{u\}\times v^{\parvtwo(p)}$ and
  $\pair<w;z>^\parvtwo = w^\parvtwo\times\{z\}\cup\{w\}\times z^{\parvtwo(p)}$ and
  that the left hand factors in these unions are subsets of
  $E_{p-1}\times F_q$ whereas their right hand factors are subsets of
  the disjoint set $E_p\times F_{q-1}$. So it follows that if
  $\pair<u;v>^\parvtwo\cup\pair<w;z>^\parvtwo \neq 0$ then either $v=z$, $u\neq
  w$ and $u^\parvtwo\cap w^\parvtwo\neq\emptyset$, which contradicts the
  well-formedness of $f(x)$, or $u=w$, $v\neq z$ and $v^{\parvtwo(p)}\cap
  z^{\parvtwo(p)}\neq\emptyset$, which contradicts the well-formedness of
  $g(y)$.  It follows therefore that
  $\pair<u;v>^\parvtwo\cup\pair<w;z>^\parvtwo= \emptyset$ and thus that
  $f(x)\times g(y)$ is well formed as required.

  \vspace{1ex}

  \noindent {\bf\boldmath $(f\times g)\pair<x;y>$ moves $(f\times
    g)(\pair<x;y>^{-})$ to $(f\times g)(\pair<x;y>^{+})$:} it is
  easily show that $((f\times g)\pair<x;y>)^\parvtwo = f(x)^\parvtwo\times
  g(x)\cup f(x)\times g(x)^{\parvtwo(p)}$ and that $(f\times
  g)(\pair<x;y>^\parvtwo) = f(x^\parvtwo)\times g(y)\cup f(x)\times
  g(y^{\parvtwo(p)})$ and it is again the case that the left factors of
  these unions are in the set $E_{p-1}\times F_q$ whereas their right
  factors are in the disjoint set $E_p\times F_{q-1}$. It follows that
  the movement property of observation~\ref{func.movement.prop} holds
  if and only if it can be verified independently on corresponding
  factors in these unions.  However we have
  \begin{equation*}\begin{split}
    (f(x^\parvtwo)\times g(x) \cup f(x)^{\neg\parvtwo}\times g(x))\smallsetminus
    (f(x)^{\parvtwo}\times g(x)) & {} = \\ 
    ((f(x^\parvtwo)\cup f(x)^{\neg\parvtwo})\smallsetminus
    f(x)^\parvtwo)\times g(x) & {} = f(x^{\neg\parvtwo})\times g(x) \\
    (f(x)\times g(x^{\parvtwo(p)})\cup f(x)\times
    g(x)^{\neg\parvtwo(p)})\smallsetminus (f(x)\times
    g(x)^{\parvtwo(p)}) & {} = \\ 
    f(x)\times ((g(x^{\parvtwo(p)})\cup g(x)^{\neg\parvtwo(p)})\smallsetminus
    g(x)^{\parvtwo(p)}) & {} = f(x)\times g(x^{\neg\parvtwo(p)})
  \end{split}\end{equation*}
  where the final equalities in these calculations are simply obtained
  by applying the operations ${-}\times g(y)$ and $f(x)\times{-}$ to
  the movement properties of the parity complex morphisms $f$ and $g$
  at the elements $x\in C_p$ and $y\in D_q$ respectively.\qed
\end{obs}

\begin{obs}[$\OSimp$ as a functor]
  \label{osimp.functor}
  Suppose that $\arrow\alpha:[m]->[n].$ is a simplicial operator and
  $\svec{v}$ is an element of $\OSimp[m]$ then we define
  $\alpha(\svec{v})$ of $\OSimp[n]$ to be the element determined by
 the equation $\face^n_{\alpha(\svec{v})}=(\alpha\circ\face^m_{
    \svec{v}})^f$. In other words the face operator corresponding to
  $\alpha(\svec{v})$ is obtained by composing the operators $\alpha$
  and $\face^m_{\svec{v}}$ and then applying the face-degeneracy
  factorisation of observation~\ref{face.degen.fac}. In terms of our
  vectorial notation, this operation corresponds to applying $\alpha$
  to the vertices of $\svec{v}$ to construct a, possibly repetitious,
  list of vertices $\alpha(v_0),\alpha(v_1),..., \alpha(v_r)$ from
  which we discard any repeats before re-indexing to obtain
  $\alpha(\svec{v})$. Our intuition is that $\alpha(\svec{v})$ is the
  face of $\OSimp[n]$ which spans the set of vertices obtained by
  applying $\alpha$ to the vertices of the face $\svec{v}$ of
  $\OSimp[m]$.
  
  Notice that while this operation is nicely functorial, in the sense
  that $(\beta\circ\alpha)(\svec{v})=\beta( \alpha( \svec{v}))$, it is
  not, in general, the case that it respects dimensions since the
  dimension of $\alpha(\svec{v})$ will be strictly less than that of
  $\svec{v}$ whenever $\alpha(v_i)=\alpha(v_{i+1})$ for some $i$.
  However, using this operation we may still define an arrow
  $\spanarr\OSimp(\alpha):\OSimp[m]-> \OSimp[n].$ of $\Parity$ by
  letting:
  \begin{equation*}
    \OSimp(\alpha)(\svec{v}) = \left\{
      \begin{array}{lp{20em}}
        \left\{\alpha(\svec{v})\right\} &
        if the dimensions of $\svec{v}$ and $\alpha(\svec{v})$ are the same, \\
        \emptyset & otherwise.
      \end{array}
      \right.
  \end{equation*}
  The verification of conditions~(\ref{func.par.a})
  and~(\ref{func.par.b}) of lemma~\ref{inf.funct.parity} for this
  function is a simple matter, the only interesting part of which is
  the demonstration that $\OSimp(\alpha)(\svec{v})$ moves
  $\OSimp(\alpha)(\svec{v}^{-})$ to $\OSimp(\alpha)(\svec{v}^{+})$ for
  each element $\svec{v}\in\OSimp[m]$. Assuming that $\svec{v}$ has
  dimension $r$ we argue by cases:
  \begin{itemize}
  \item {\bf\boldmath $\alpha(\svec{v})$ has dimension $r$}, in which
    case it is easily seen that $\OSimp(\alpha)(\svec{v}^\parvtwo)=
    \alpha(\svec{v})^\parvtwo$ (for either parity $\parvtwo$), but we know that
    every element moves its set of negative faces to its set of
    positive faces so it is immediate that
    $\OSimp(\alpha)(\svec{v})=\{\alpha(\svec{v})\}$ has the required
    movement property.
  \item {\bf\boldmath $\alpha(\svec{v})$ has dimension $r-1$}, in which
    case we have $\alpha(v_i)=\alpha(v_{i+1})$ for exactly one $0\leq
    i<r$ and it follows that $\alpha(\svec{v}\cdot\face^r_j)$ has
    dimension $r-2$ for all $j\neq i,i+1$ and that
    $\alpha(\svec{v}\cdot\face^r_i) =\alpha(\svec{v}\cdot\face^r_{i+1})=
    \alpha(\svec{v})$. However the faces $\svec{v}\cdot\face^r_i$ and
    $\svec{v}\cdot\face^r_{i+1}$ of $\svec{v}$ are of opposite parity
    and so it follows that $\OSimp(\alpha)(\svec{v}^{-})=
    \OSimp(\alpha)(\svec{v}^{+})= \{\alpha(\svec{v})\}$ and it is clear
    that the set $\OSimp(\alpha)(\svec{v})$ which is empty moves any
    set to itself as required.
  \item {\bf\boldmath $\alpha(\svec{v})$ has dimension $\leq r-2$}, in
    which case each element $\alpha(\svec{v}\cdot \face^r_i)$ must also
    have dimension $\leq r-2$, so $\OSimp(\alpha)(\svec{v})=
    \OSimp(\alpha)(\svec{v}^{-})= \OSimp(\alpha)(\svec{v}^{+})=
    \emptyset$ and it follows that the movement property is trivial in
    this case.
  \end{itemize}
  It is also a matter of trivial verification to show that this action
  on simplicial operators preserves composition, giving a functor
  $\arrow \OSimp:\Delta->\Parity.$.
\end{obs}

\subsection{Collapsers and Stratified Parity Complexes}

\begin{obs}[collapsing and collapsers]\label{collapsing.quotients} 
  Suppose that $\mathbb{C}$ is an \inf-category and that
  $S=\bigcup_{n>0} S_n$ is a graded subset of it's cells. Then we say
  that an \inf-functor $\arrow f:\mathbb{C}->\mathbb{D}.$ {\em
    collapses\/} (the cells in) $S$ iff for each $n\in\mathbb{N}$ it
  maps each cell $c$ in the subset $S_{n+1}$ of $\mathbb{C}$ to an
  $n$-cell $f(c)$ in $\mathbb{D}$.  Since $\InfCat$ possesses all
  colimits, it follows that we may form a quotient \inf-category
  $\mathbb{C}/S$ of $\mathbb{C}$ which represents the \inf-functors
  that collapse $S$
  \begin{equation*}
    \let\labelstyle=\textstyle
    \xymatrix@C=1em@R=2.6em{
      {\mathbb{C}}\ar@{-|>}[rr]^-{q_S}\ar[dr]_{f}
      & & {\mathbb{C}/S}\ar@{-->}[dl]^{\exists!\,\hat{f}} \\
      & {\mathbb{D}} &
    }
  \end{equation*}
  in the sense that we have a quotient \inf-functor $q_S$, as depicted
  in the diagram above, which collapses $S$ and through which any
  other \inf-functor $f$ which collapses $S$ factors uniquely. An
  \inf-functor which possesses this universal property is said to be a
  {\em collapser\/} of (the cells in) $S$.
  
  At a couple of places in the sequel, we will have use for a standard
  and very simple factorisation result with respect to these
  collapsers. Suppose that $S$ and $T$ are two graded subsets of the
  cells of $\mathbb{C}$ with $T\subseteq S$ then in the diagram
  \begin{equation*}
    \let\labelstyle=\textstyle
    \xymatrix@C=1em@R=2.6em{
      {\mathbb{C}} \ar@{-|>}[rr]^{q_T} 
      \ar@{-|>}[dr]_{q_S} && {\mathbb{C}/T}
      \ar@{-->}[dl]^{\exists!\, q} \\
      & {\mathbb{C}/S} & 
    }
  \end{equation*}
  the \inf-functor $q_S$ collapses the cells in $T$, since it
  collapses those in $S$ and $T\subseteq S$, so it factors through
  $q_T$ to give a \inf-functor which is itself a collapser. To be
  precise, $q$ is a collapser of the graded subset
  $q_T(S\smallsetminus T)\subseteq\mathbb{C}/T$ given by
  \begin{math}
    q_T(S\smallsetminus T)_r=\{ q_T(x) \mid x\in S_r\also x\notin T_r\}
  \end{math}
  for each $r>0$.
\end{obs}

\begin{defn}[stratified parity complexes]\label{strat.parity}
  By analogy with stratified sets, we define a {\em stratified parity
    complex\/} to be a pair $\strat{C}$ in which $C$ is a parity
  complex and $tC$ is a subset of its elements, called {\em thin
    elements}, all of which must be of dimension $>0$. A {\em
    stratified morphism\/} $\spanarr f:\strat{C}->\strat{D}.$ between
  such structures is simply a parity complex morphism $\spanarr
  f:C->D.$ satisfying the thinness preservation axiom that
  $f(tC)\subseteq tD$. Equivalently this axiom states that for all
    $x\in tC$ we have $y\in tD$ for all $y\in f(x)$. 
  
  Suppose that $\spanarr g:\strat{D}->\strat{E}.$ is another stratified
  morphism then their
 composite $\spanarr g\circ f:C->E.$ in $\Parity$
  also satisfies the thinness preservation axiom (with respect to $tC$
  and $tE$) and so we may lift the composition of $\Parity$ to
  construct a category $\SParity$ of stratified parity complexes and
  stratified morphisms between them. Notice that we may identify
  $\Parity$ with the full subcategory of those stratified parity
  complexes in $\SParity$ whose set of thin elements is empty, which
  we shall do routinely in the sequel.
  
  In general, when discussing stratified parity complexes we will
  adopt the same notational conventions that we have been applying to
  stratified sets, including dropping explicit mention of their sets
  of thin elements (cf.\ definition~\ref{strat.defn}). Having done
  this, we also adopt $\arrow \forget:\SParity->\Parity.$ to denote
  the functor which forgets stratifications and use this when we wish
  to notationally distinguish between a stratified parity complex $C$
  and its underlying parity complex $\forget(C)$.
  
  It should also be clear that most of the basic concepts,
  constructions and standard structures introduced for stratified sets
  may be naturally transferred to stratified parity complexes. We
  leave the details of this process up to the reader, since we will be
  content to restrict ourselves to the relatively limited collection
  of results and generalisations discussed in the next few
\end{defn}

\begin{obs}[extending $\oriental$ to $\SParity$]
  \label{sparity.ext.oriental}
  If $C$ is a stratified parity complex then we define a
  graded set $\atom{tC}$ of cells of the freely generated
  \inf-category $\oriental(\forget(C))$ by
  \begin{math}
    \atom{tC}_r = \{ \atom{x} \mid x\in tC\cap C_r \}
  \end{math}
  ($r=1,2,...$) and let $\oriental(C)$ to be the collapsed
  \inf-category $\oriental(\forget(C))/\atom{tC}$. We'll generally
  denote the associated collapser by $\cover q_t:\oriental(\forget(C))
  ->\oriental(C).$ and extend our atom notation by letting $\atom{x}
  \in\oriental(C)$ denote the cell obtained by applying this collapser
  to the atom in $\oriental(\forget(C))$ which corresponds to the
  element $x\in C$. Of course this collapser is an epimorphism, so we
  may apply lemma~\ref{gen.epi} to show that the set
  \begin{equation*}
    \atom{C}\defeq \bigcup_{n\in\mathbb{N}}\{ \atom{x} 
    \mid \text{$x\in C_n$ and $\atom{x}$ is not an $(n-1)$-cell 
      in $\oriental(C)$.}\}
  \end{equation*}
  is a weak generator of $\oriental(C)$.
  
  To extend this construction to a functor on $\SParity$, suppose that
  we are given a stratified morphism $\spanarr f:C->D.$ and consider
  the \inf-functor $\arrow\oriental(f):\oriental(\forget(C))->
  \oriental(\forget(D)).$ between the free \inf-categories on
  underlying parity complexes.  If $x$ is an $n$-element in $C$ then
  $\atom{x}$ is of course a non-trivial $n$-cell in
  $\oriental(\forget(C))$ and we know that $\oriental(f)$ maps it to
  an $n$-cell $\pair<N;P>$ in $\oriental(\forget(D))$ which has
  $N_n=P_n=f(x)$.  So applying Street's Excision of Extremals
  algorithm, theorem 4.1 of~\cite{Street:1991:Parity}, as discussed in
  the proof of theorem 4.2 of that paper, we see that the cell
  $\pair<N;P>$ may be expressed as a composite of atoms on elements in
  the sub-parity complex of $D$ generated by $N\cup P$. In particular,
  all atoms in this decomposition are $n$-cells and the only ones that
  are non-trivial are those which are atoms on an element in $f(x)$.
  Now if $x$ is thin then the thinness preservation property of $f$
  tells us that $f(x)$ is a set of thin $n$-elements in $D$, so we
  know that for each $y\in f(x)$ the quotient functor $\arrow
  q_t:\oriental(\forget(D))->\oriental(D).$ collapses the
  atom $\atom{y}$ to an $(n-1)$-cell and therefore that it collapses
  each cell in the atomic decomposition of $\pair<N;P>$, and thus this
  cell itself, to an $(n-1)$-cell in $\oriental(D)$.
  
  In other words, we have shown that the composite \inf-functor
  $\arrow q_t\circ\oriental(f): \oriental(\forget(C))->
  \oriental(D).$ collapses the cells in the graded set
  $\atom{tC}$. From this we can infer, by the universal property of
  the quotient $\oriental(D)$, that this composite induces the
  unique (dashed) \inf-functor which makes the following square
  commute:
  \begin{equation*}
    \let\labelstyle=\textstyle
    \xymatrix@C=5em@R=2em{
      {\oriental(\forget(C))}\ar[r]^{\oriental(f)}
      \ar@{-|>}[d]_{q_t} &
      {\oriental(\forget(D))}\ar@{-|>}[d]^{q_t} \\
      {\oriental(C)}\ar@{-->}[r]_{\exists!\,\oriental(f)} &
      {\oriental(D)}
    }
  \end{equation*}
  So we define $\arrow\oriental(f):\oriental(C)-> \oriental(D).$ to be
  this induced \inf-functor and it follows, from its uniqueness
  property, the functoriality of $\oriental$ on $\Parity$ and the fact
  that stratified morphisms compose as in the underlying category
  $\Parity$, that this action is functorial with respect to the
  category structure on $\SParity$ as required.
  
  Suppose now that $D$ is an entire sub-complex of the stratified
  parity complex $C$, in other words these share the same underlying
  parity complex $\forget(C)$ and have $tD\subseteq tC$. Then
  $\oriental(C)$ and $\oriental(D)$ are constructed by taking
  collapsers of $\oriental(\forget(C))$ with regard to the sets
  $\atom{tC}$ and $\atom{tD}$ of cells respectively. The latter of
  these is a subset of the former, so we may apply the factorisation
  result of observation~\ref{collapsing.quotients} to show that the
  induced \inf-functor $\arrow\oriental(\subseteq_e):\oriental(D)->
  \oriental(C).$ is a collapser of the graded set given by
  $\atom{tC\smallsetminus tD}_r= \{\atom{x}\mid x\in tC\cap C_n\also
  x\notin tD\}$.

  If we consider $\Parity$ to be the full subcategory of trivially
  stratified parity complexes in $\SParity$, as mentioned in the last
  observation, then the functor $\oriental$ on $\SParity$ is indeed an
  extension of our original free \inf-category functor on $\Parity$.
\end{obs}

\begin{obs}
  \label{simplex.to.sparity}
  Suppose that $N$ is a stratified set in the full subcategory
  $\Simplex$ of definition~\ref{thinner.simplices}, which according to
  convention has underlying simplicial set $\Delta[n]$. Then we may
  define a corresponding stratified parity complex $\tilde{N}$ which
  has underlying parity complex $\OSimp[n]$ and in which an element
  $\svec{v}$ is thin iff the corresponding face $\face^n_{\svec{v}}\in
  \Delta[n]$ is thin in $N$.
  
  Furthermore, if $\arrow f:N->M.$ is a stratified map in $\Simplex$
  then we may apply the simplicial Yoneda lemma to show that there is
  a unique simplicial operator $\arrow\alpha:[n]->[m].$ for which
  $\arrow\Delta(\alpha):\Delta[n]->\Delta[m].$ is the simplicial map
  underlying $f$. Consider $\spanarr\OSimp(\alpha):
  \OSimp[n]->\OSimp[m].$ the corresponding parity complex morphism and
  recall that, from the definition given in
  observation~\ref{osimp.functor}, $\OSimp(\alpha)(\svec{v})=
  \{\svec{w}\}$ if and only if $\face^m_{\svec{w}}=\alpha\circ
  \face^n_{\svec{v}}= \Delta(\alpha)(\face^n_{\svec{v}})$ and that
  otherwise $\OSimp(\alpha)(\svec{v})=\emptyset$. In the former case,
  if $\svec{v}$ is thin in $\tilde{N}$ then we know, by definition,
  that $\face^n_{\svec{v}}$ is thin in $N$ and thus that
  $\Delta(\alpha)(\face^n_{\svec{v}})=\face^m_{\svec{w}}$ is thin in
  $M$, since $\Delta(\alpha)$ underlies the stratified map $f$, so
  finally we see that $\svec{w}$ is thin in $\tilde{M}$. It follows
  that $\OSimp(\alpha)$ provides us with a stratified parity complex
  morphism $\arrow \tilde{f}:\tilde{N}->\tilde{M}.$ and that this
  construction enriches our operation on the objects of $\Simplex$ to
  a functor from that category into $\SParity$.
\end{obs}

\begin{obs}[collapsed orientals]\label{coll.oriental}
  Of course, if we consider $\tDelta$ to be a full subcategory of
  $\Strat$ then it is, in fact, a subcategory of $\Simplex$ and it
  follows that we may restrict the functor of the last observation to
  a functor $\arrow\Delta:\tDelta->\SParity.$. This functor extends
  the functor $\arrow \OSimp:\Delta->\Parity.$ of
  observation~\ref{osimp.functor} by mapping $[n]_t$ to the stratified
  parity complex $\OSimp[n]_t$ obtained from $\OSimp[n]$ by making
  thin its unique $n$-dimensional element $01...n$.
    
  Now, we've seen that the atom $\atom{01...n}$ is the only
  non-trivial $n$-cell of the oriental $\oriental(\OSimp[n])$ so it is
  clear that $\oriental(\OSimp[n]_t)$ is an $(n-1)$-category, since it
  is obtained by collapsing this unique $n$-cell. We call
  $\oriental(\OSimp[n]_t)$ the {\em $n^{\text{th}}$ collapsed
    oriental}. Our intuition here is that this is the appropriate
  \inf-categorical model of the standard thin $n$-simplex.
\end{obs}

\begin{obs}[extending a few other constructions to $\SParity$]
  \label{sparity.prod.th}
  If $C$ and $D$ are stratified parity complexes then their product is
  obtained by taking the underlying product of parity complexes, just
  as in observation~\ref{parity.prods.and.cubes}, and defining its set
  of thin elements to be:
  \begin{equation*}
    t(C\times D) \defeq \{ \pair<x;y>\in C\times D \mid
    \text{either $x\in tC$ or $y\in tD$}\}
  \end{equation*}
  Furthermore, if $\spanarr f:C->E.$ and $\spanarr g:D->F.$ are
  stratified morphisms it is a matter of routine verification to check
  that their product as parity complex morphisms satisfies the
  thinness preservation axiom and thus that it provides a stratified
  morphism $\spanarr f\times g: C\times D->E\times F.$. It follows
  that we may canonically extend the product bifunctor on $\Parity$ to
  a bifunctor $\arrow \times:\SParity \times\SParity ->\SParity.$
  simply by lifting its action on underlying parity complex morphisms.
  
  Another construction on stratified parity complexes that we shall
  have use for in the sequel is defined by analogy to the functors
  $\arrow\Th_n:\Strat->\Strat.$ introduced in
  definition~\ref{superstr.defn}. If $C$ is a stratified parity
  complex and $n\in\mathbb{N}$ then $\Th_n(C)$ is defined to be the
  stratified parity complex obtained from $C$ by making thin all
  elements of dimension greater than $n$. Of course, it is clear that
  any stratified parity complex morphism $\arrow f:C->D.$ lifts to a
  stratified morphism from $\Th_n(C)$ to $\Th_n(D)$, which we call
  $\Th_n(f)$, making $\Th_n$ into an endo-functor on $\SParity$. The
  identity morphism on the underlying parity complex of $C$ provides
  us with a canonical stratified morphism $\spanarr i_C^n:
  C->\Th_n(C).$ and the corresponding \inf-functor
  $\cover\oriental(i_C^n):\oriental(C)->\oriental(\Th_n(C)).$ is the
  collapser of the graded set of atoms of $C$ on elements of dimension
  greater than $n$.
\end{obs}

\subsection{\inf-Categorical Nerve Constructions}

\begin{obs}
  \label{osimp(alpha).act.on.atoms}
  Composing the functors $\arrow\OSimp:\Delta->\Parity.$ and
  $\arrow\oriental:\Parity->\InfCat.$ we obtain a functor which
  represents the category of simplicial operators as a subcategory of
  $\InfCat$. Now, if $\arrow\alpha:[m]->[n].$ is a simplicial
  operator, then the corresponding \inf-functor
  $\arrow\oriental(\OSimp( \alpha)):\oriental( \OSimp[m])->
  \oriental(\OSimp[n]).$ is completely determined by its action on
  atoms, since they freely generate $\oriental(\OSimp[m])$. In the
  sequel, it will useful to have a more explicit description of this
  action on atoms.
  
  Indeed, it is the case that $\oriental(\OSimp(\alpha))$ is the
  unique \inf-functor which maps each atom
  $\atom{\svec{v}}\in\oriental(\OSimp[m])$ to the atom
  $\atom{\alpha(\svec{v})}\in \oriental(\OSimp[n])$.  Unfortunately,
  this fact does not follow immediately from our argument so far, but
  it can be established by a somewhat tedious direct calculation.
  However we choose to take a different route and use an argument
  based on corollary~\ref{morp.action.1}. So let
  $E=\{\svec{w}\in\OSimp[m] \mid \svec{w} \subseteq\svec{v}\}$ be the
  sub-parity complex of $\OSimp[m]$ generated by $\svec{v}$ and notice
  that if $\svec{w}\subseteq\svec{v}$ then $\alpha(\svec{w})
  \subseteq\alpha(\svec{v})$ so it is clear that $\alpha$ maps
  elements of $E$ into the sub-parity complex $F=\{\svec{u}\in
  \OSimp[n]\mid \svec{u} \subseteq\alpha(\svec{v})\}$ generated by
  $\alpha(\svec{v})$ in $\OSimp[n]$ and it follows that
  $\OSimp(\alpha)$ restricts to a parity complex morphism from $E$ to
  $F$.  Now suppose that $\svec{u}=u_0u_1...u_r\in F_r$ (for some
  dimension $r$) then from the definition of $F$ we know that the
  vertices of $\svec{u}$ are all vertices of $\alpha(\svec{v})$ and
  are thus of the form $\alpha(v_j)$ for some $j$.  Consequently we
  may construct an element $\svec{w}\in E_r$ by letting
  $w_i=\min\{v_j\mid j\in[m]\also \alpha(v_j)= u_i\}$ for which we
  know, by construction, that $\alpha(\svec{w})=\svec{u}$.  It follows
  that we have $\OSimp(\alpha)(\svec{w})=\{\svec{u}\}$ and so
  $\svec{u}\in\OSimp(\alpha)(E_r)$, or in other words we have
  succeeded in showing that the subset $\OSimp(\alpha)(E)$ is actually
  equal to $F$ itself. Now we may apply corollary~\ref{morp.action.1}
  to conclude that $\OSimp(\alpha)(\atom{\svec{v}})=\atom{\alpha(
    \svec{v})}$ as postulated.
  
  We can also go further and extend this analysis to free
  \inf-categories of the form $\oriental(\OSimp[n]\times\OSimp[m])$,
  which we will study in some detail in the next subsection. In this
  case if $\arrow\phi:[n]->[r].$ and $\arrow\psi:[m]->[s].$ are a pair
  of simplicial operators then we may easily adapt the argument of the
  last paragraph to show that if $\pair<\svec{v};
  \svec{w}>\in\OSimp[n]\times\OSimp[m]$ then
  $\oriental(\OSimp(\phi)\times\OSimp(\psi))(\atomp<\svec{v};\svec{w}>)
  = \atomp<\phi(\svec{v});\psi(\svec{w})>$.
\end{obs}

\begin{obs}[Street's \inf-categorical nerve]
  \label{inf.cat.nerve}
  Our primary reason for extending $\OSimp$ to a functor
  $\overarr:\Delta->\Parity.$ and forming the composite functor
  $\arrow\oriental\circ\OSimp:\Delta->\InfCat.$ is that we may now
  follow Street~\cite{Street:1987:Oriental} and define an adjoint pair
  \begin{equation}\label{simp.nerve}
    \let\labelstyle=\textstyle
    \xymatrix@R=1ex@C=14em{
      {\InfCat}\ar[r]^{\bot}_{\inerve} &
      {\Simp}\ar@/_3ex/[l]_{\inladj}}
  \end{equation}
  by Kan's construction (observation~\ref{func.from.coalg}). 
  
  In other words, the right adjoint $\inerve$ is constructed by
  ``homming out'' of the functor $\oriental\circ\OSimp$, that is to
  say if $\mathbb{C}$ is an \inf-category then the simplicial set
  $\inerve(\mathbb{C})$ has $\inerve(\mathbb{C})_{n}\defeq
  \InfCat(\oriental(\OSimp[n]),\mathbb{C})$ as its set of
  $n$-simplices and a right action given by the pre-composite
  $x\cdot\alpha\defeq x\circ\oriental(\OSimp(\alpha))$ where
  $x\in\inerve(\mathbb{C})_{n}$ is an $n$-simplex and
  $\arrow\alpha:[m]->[n].$ is a compatible simplicial operator. 
  
  The left adjoint $\inladj$ may be constructed as the left Kan
  extension of the composite functor
  $\arrow\oriental\circ\OSimp:\Delta-> \InfCat.$ along the Yoneda
  embedding $\arrow\Delta:\Delta->\Simp.$, that is we can write
  $\inladj(X)$ as the weighted colimit
  $\colim(X,\oriental\circ\OSimp)$ in $\InfCat$.  A cocone under the
  diagram $\oriental\circ\OSimp$ weighted by $X$ is simply a family of
  \inf-functors $\arrow f_x:\oriental(\OSimp[n])->\mathbb{C}.$ for
  $x\in X$ (where $n=\dim(x)$) satisfying the naturality condition
  that for each simplicial operator $\arrow\alpha:[m]->[n].$ which is
  compatible with $x\in X$ we have $f_{x\cdot\alpha}=f_x\circ
  \oriental(\OSimp(\alpha))$.
  
  In future calculations we'll use the notation $\arrow\coliminj_x:
  \oriental(\OSimp[n])->\inladj(X).$ to denote the components of the
  colimiting cocone.  Of course any other such cocone $f_x$ induces a
  \inf-functor $\arrow f:\inladj(X)->\mathbb{C}.$ which is uniquely
  determined by the equality $f\circ\coliminj_x=f_x$ for each $x\in
  X$. In particular, if $\arrow f:X->Y.$ is a simplicial map then the
  \inf-functor $\arrow \inladj(f):\inladj(X)-> \inladj(Y).$ is
  characterised by the fact that it is the unique \inf-functor which
  makes the triangle
  \begin{equation}\label{inladj.act.on.arr}
    \let\labelstyle=\textstyle
    \xymatrix@C=0.8em@R=2.6em{
      & {\oriental(\OSimp[n])}\ar[dl]_{\coliminj_x}
      \ar[dr]^{\coliminj_{f(x)}} & \\
      {\inladj(X)}\ar[rr]_{\inladj(f)} & &
      {\inladj(Y)} }
  \end{equation}
  commute for each simplex $x\in X$ (where $n=\dim(x)$).
\end{obs}

\begin{obs}[$\inladj(X)$ as a (freely) generated \inf-category]
  \label{fg.inladj.obs}
  Suppose that $x$ is an $n$-simplex in $X$ then we adopt the notation
  $\satom{x}$ to denote the cell $\coliminj_x(\atom{01...n})\in
  \inladj(X)$. At the risk of confusing our nomenclature a little
  we'll tend to refer to $\satom{x}$ as the {\em atom of $\inladj(X)$
    associated with the simplex $x$}.
  
  Notice that if $\arrow\alpha:[m]->[n].$ is a simplicial operator
  then we have
  \begin{equation}\label{inladj.inj.beh}
    \coliminj_x(\atom{\alpha(01...m)})=
    (\coliminj_x\circ\oriental(\OSimp(\alpha)))(\atom{01...m})=
    \coliminj_{x\cdot\alpha}(\atom{01...m})\defeq\satom{x\cdot\alpha}
  \end{equation}
  where the first equality follows by
  observation~\ref{osimp(alpha).act.on.atoms} and the second expresses
  the naturality of the cocone $\coliminj_x$.  As a consequence if
  $\arrow \alpha:[m]->[n].$ is a face operator then
  \begin{equation*}
    \satom{x\cdot\alpha}=\coliminj_x(\atom{\alpha(01...m)})=
    \coliminj_x(\atom{01...n})\defeq\satom{x} 
  \end{equation*}
  where the second of these equalities holds because $\alpha$ is
  surjective and thus has $\alpha(01...m)=01...n$.  
  
  Conversely, we may show that whenever $x$ is a non-degenerate
  simplex of $X$ then $\satom{x}$ is actually a non-trivial $n$-cell
  of $\inladj(X)$.  This can be done by employing the explicit
  description of $\inladj(X)$ as a colimit to construct an
  \inf-functor $k_X$ from there into $\norient(\tilde{X})$, the
  \inf-category constructed from the graded set $\tilde{X}$ of
  non-degenerate simplices in $X$ as in
  observation~\ref{norient.defn}, which maps each atom $\satom{x}$ to
  a cell $\pair<M;P>\in\norient(\tilde{X})$ defined by:
  \begin{equation*}
    M = \left\{ x\cdot\face^n_{\svec{v}}\mid \svec{v}\in\mu(01...n)
    \right\}\cap\tilde{X} \mkern40mu
    P = \left\{ x\cdot\face^n_{\svec{v}}\mid \svec{v}\in\pi(01...n)
    \right\}\cap\tilde{X} 
  \end{equation*}
  Then, in the case where $x$ is non-degenerate we see that
  $M_n=P_n=\{x\}$ so it follows that $k_X(\satom{x})$ is a non-trivial
  $n$-cell in $\norient(\tilde{X})$ and thus (cf.\ 
  observation~\ref{ncells.ncats}) that $\satom{x}$ itself is a
  non-trivial $n$-cell in $\inladj(X)$ as required. We leave the
  detailed verification of this result as an exercise for the reader.
  
  Indeed we may actually show that $\inladj(X)$ is freely generated by
  its set of atoms:
  \begin{equation*}
    \satom{X}\defeq\left\{\satom{x}\mid x\in X\right\}
    =\left\{\satom{x}\mid\text{$x$ is a non degenerate
        simplex of $X$}\right\}
  \end{equation*}
  However this result will not be required here and so we leave its
  (relatively routine) proof as an exercise for the reader.  Instead
  we will content ourselves with simply showing that $\satom{X}$
  weakly generates $\inladj(X)$. To do so we apply
  lemma~\ref{gen.epi}, via the associated
  observation~\ref{gen.epi.obs}, to the jointly epimorphic family
  $\arrow \coliminj_x:\oriental(\OSimp[n])->\inladj(X).$ of colimiting
  cocone components.  Now if $\svec{v}$ is an $m$-dimensional element
  of $\OSimp[n]$ we have $\svec{v}=\face^n_{\svec{v}}(01...m)$ so it
  follows, from equation~(\ref{inladj.inj.beh}), that
  $\coliminj_x(\atom{ \svec{v}}) =\satom{x\cdot\face^n_{ \svec{v}}}$.
  In other words, each \inf-functor $\coliminj_x$ maps the generators
  (atoms) of its domain $\oriental(\OSimp[n])$ into $\satom{X}$,
  therefore the conclusion of lemma~\ref{gen.epi} is that this set
  weakly generates $\inladj(X)$ as required.
  
  In particular, it follows that if $\arrow f:X->Y.$ is a simplicial
  map then the characterisation in display~(\ref{inladj.act.on.arr})
  can be re-expressed as stating that $\arrow \inladj(f):\inladj(X)->
  \inladj(Y).$ is the unique \inf-functor which is determined by the
  fact that for each non-degenerate simplex $x\in X$ it maps the cell
  $\satom{x}$ in $\inladj(X)$ to the cell $\satom{f(x)}$ in
  $\inladj(Y)$.
\end{obs}

\begin{obs}\label{inf.cat.nerve.strat}
  It turns out, however, that this nerve construction discards a
  little too much information to make it a useful representation of
  \inf-categories as simplicial structures. To be precise, we do not
  preserve enough about the identities in the \inf-categories we are
  representing and, consequently, there is no hope that we might be
  able to re-construct the compositional structures of an
  \inf-category from its simplicial nerve.
  In~\cite{Street:1987:Oriental} Street follows the lead set by
  Roberts in~\cite{Roberts:1977:Complicial}
  and~\cite{Roberts:1978:Complicial} and suggests that the appropriate
  approach to rectifying this deficiency would be to extend the nerve
  construction to the category of stratified sets
  \begin{equation}\label{strat.nerve}
    \let\labelstyle=\textstyle
    \xymatrix@R=1ex@C=14em{
      {\InfCat}\ar[r]^{\bot}_{\inerve} &
      {\Strat}\ar@/_3ex/[l]_{\inladj}}
  \end{equation}
  in which we can use appropriately selected sets of thin simplices to
  preserve enough information about identities between the composites
  that form their sources and targets.
  
  To make this construction precise, we proceed as in
  observation~\ref{inf.cat.nerve} by forming the composite of the
  extended functors $\arrow\OSimp:\tDelta->\SParity.$ of
  observation~\ref{coll.oriental} and $\arrow\oriental:\SParity->
  \InfCat.$ of observation~\ref{sparity.ext.oriental}. This composite
  extends $\arrow\oriental\circ\OSimp:\Delta->\InfCat.$ by mapping the
  thin simplex $[n]_t$ to the collapsed oriental
  $\oriental(\OSimp[n]_t)$ and the operator
  $\arrow\thop^n:[n]->[n]_t.$ to the canonical collapser $\cover
  q:\oriental(\OSimp[n])->\oriental( \OSimp[n]_t).$. Of course, each
  of these quotient maps is an epimorphism and so it follows that the
  composite $\arrow \oriental\circ\OSimp:\tDelta->\InfCat.$ is
  actually a $\mathbb{T}_{\Strat}$-coalgebra. It follows that we may
  construct the adjunction shown in display~\eqref{strat.nerve} by
  applying Kan's construction, observation~\ref{func.from.coalg}, to
  this coalgebra.
  
  More explicitly, we may form the stratified nerve
  $\inerve(\mathbb{C})\in\Strat$ of an \inf-category $\mathbb{C}$ by
  first forming its nerve in $\Simp$ and then making thin any
  $n$-simplex $x\in\inerve(\mathbb{C})$ for which the corresponding
  \inf-functor $\arrow x:\oriental(\OSimp[n])-> \mathbb{C}.$ carries
  the atom $\atom{01...n}$ to an $(n-1)$-cell in $\mathbb{C}$. Going
  in the other direction, to form $\inladj(X)$ for a stratified set
  $X$ we first apply the left adjoint of display~(\ref{simp.nerve}) to
  the underlying simplicial set $\forget(X)\in\Simp$ to give an
  \inf-category $\inladj(\forget(X))$ then we define a graded subset
  of the cells of $\inladj(\forget(X))$ by
  \begin{equation*}
    \satom{tX}_n\defeq\left\{\satom{x}\mid \text{$x\in X$ is a 
        non-degenerate, thin $n$-simplex}\right\}
  \end{equation*}
  and let $\inladj(X)$ be the quotient $\inladj(\forget(X))/
  \satom{tX}$ of observation~\ref{collapsing.quotients}.  In other
  words, we may construct $\inladj(\forget(X))$ as a colimit of
  orientals and then collapse those (images of) orientals in there
  which correspond to thin simplices in $X$. We'll generally denote
  the associated collapser by $\cover q_t: \inladj(\forget(X)) ->
  \inladj(X).$ and extend our atom notation by letting
  $\satom{x}\in\inladj(X)$ denote the $n$-cell obtained by applying
  this collapser to the corresponding atom in $\inladj(\forget(X))$
  associated with the $n$-simplex $x\in X$.
  
  Of course, the collapser $\cover q_t:\inladj(\forget(X))
  ->\inladj(X).$ is an epimorphism and it carries the atom
  $\satom{x}\in\inladj(\forget(X))$ associated with any thin
  $n$-simplex to an $(n-1)$-cell in $\inladj(X)$.  Therefore we may
  apply lemma~\ref{gen.epi} to the weak generating set
  $\satom{\forget(X)}\subseteq\inladj(\forget(X))$ to show that the
  set
  \begin{equation}\label{satom.wgen}
    \satom{X} \defeq \bigcup_{n\in\mathbb{N}}\left\{\satom{x}\left|\;
        \text{$x\in X_n$ and $\satom{x}$ is not an $(n-1)$-cell in 
          $\inladj(X)$.}
      \right.\right\}
  \end{equation}
  weakly generates the \inf-category $\inladj(X)$. Again, it follows
  that if $\arrow f:X->Y.$ is a stratified map then the \inf-functor
  $\arrow\inladj(f):\inladj(X)->\inladj(Y).$ is completely
  characterised by the fact that for each simplex $x\in X$ it maps the
  cell $\satom{x}\in\inladj(X)$ to the cell $\satom{f(x)}\in
  \inladj(Y)$.
  
  Suppose that $Y$ is an entire subset of the stratified set $X$, in
  other words these stratified sets share the same underlying
  simplicial set and $tY\subseteq tX$. Then $\inladj(X)$ and
  $\inladj(Y)$ are constructed by taking collapsers of
  $\inladj(\forget(X))$ with regard to the sets $\satom{tX}$ and
  $\satom{tY}$ of cells respectively. The latter of these is a subset
  of the former, so we may apply the factorisation result of
  observation~\ref{collapsing.quotients} to show that the induced
  \inf-functor $\arrow\inladj(\subseteq_e):\inladj(Y)->\inladj(X).$ is
  a collapser of the graded set given by $\satom{tX\smallsetminus tY}_r=
  \{\satom{x}\mid x\in tX\cap X_n\also x\notin tY\}$.
  
  It might appear that we have introduced some ambiguity in our
  notation by using $\inladj$ to denote the left adjoints from $\Simp$
  and $\Strat$ into $\InfCat$ and at the same time identifying the
  category $\Simp$ with the full subcategory of $\Strat$ on those
  objects which are minimally stratified. However, we know that each
  thin simplex of a minimally stratified set $X$ is degenerate,
  therefore $\satom{tX}=\emptyset$ and we see that, in this case, the
  quotient $\inladj(X)$ is actually (isomorphic to)
  $\inladj(\forget(X))$ itself.
\end{obs}

\begin{lemma}\label{inerv.rest}
  For each $n\in\mathbb{N}$ the functors $\inerve$ and $\inladj$
  restrict to give an adjunction
  \begin{equation}\label{inerv.rest.diag}
    \let\labelstyle=\textstyle
    \xymatrix@R=1ex@C=14em{
      {\nCat}\ar[r]^{\bot}_{\nnerve} &
      {\Strat_n}\ar@/_3ex/[l]_{\nnladj}}
  \end{equation}
  between the categories of $n$-categories and $n$-trivial
  stratified sets.
\end{lemma}

\begin{proof}
  If $\mathbb{C}$ is an $n$-category and $r>n$ then an $r$-simplex of
  $\inerve(\mathbb{C})$ is simply an \inf-functor $\arrow
  x:\oriental(\OSimp[r])->\mathbb{C}.$. However all cells in
  $\mathbb{C}$ are $n$-cells and so it follows that any such
  \inf-functor must collapse the oriental that is its domain and thus,
  by definition, that our $r$-simplex is thin in
  $\inerve(\mathbb{C})$. In other words, we see easily that
  $\inerve(\mathbb{C})$ is indeed $n$-trivial as stated.
  
  To prove that $\inladj$ also restricts as described, consider an
  $n$-trivial stratified set $X$ and recall that we know,
  from observation~\ref{deltat.dense}, that it admits a canonical
  representation as the colimit of a diagram
  $\arrow\dgm_{X}:\groth(X)->\Strat.$ of standard (thin) simplices.
  Furthermore the left adjoint $\Th_n$ to $\Sup_n(\cdot)$ preserves
  this colimit and we know that $\Th_n(X)=X$, since in any
  $n$-trivial stratified set every simplex of dimension
  greater than $n$ is already thin, so it follows that we may express
  $X$ as a colimit of stratified sets of the form
  $\Th_n(\Delta[r]_?)$. Applying the left adjoint functor $\inladj$ to
  this we see that $\inladj(X)$ may be constructed as a colimit of
  \inf-categories of the form $\inladj(\Th_n(\Delta[r]_{?}))$.
  
  Now $\inladj(\Th_n(\Delta[r]_?))$ is constructed from
  $\inladj(\Delta[r])$, the freely generated \inf-category
  $\oriental(\OSimp[r])$, by collapsing a graded set of cells which
  contains all atoms on elements of $\OSimp[r]$ with dimension greater
  than $n$. However any $r$-cell in $\oriental(\OSimp[r])$ may be
  expressed as a composite of atoms and all atoms of dimension $r$
  greater than $n$ have been collapsed to $(r-1)$-cells in
  $\inladj(\Th_n(\Delta[r]_{?}))$ so we may show, by a simple
  inductive proof on superstructures, that
  $\inladj(\Th_n(\Delta[r]_{?}))$ is in fact an $n$-category. Finally
  we know that $\nCat$ is closed in $\InfCat$ under colimits so it
  follows that $\inladj(X)$, as a colimit of $n$-categories, is
  actually an $n$-category itself as required.
\end{proof}

\begin{obs}\label{inladj.sup.compare}
  We might also expect that a strong relationship should hold between
  $\inladj(\Sup_n(X))$, the \inf-category associated with the
  $n$-superstructure of a stratified set $X$, and the
  $n$-superstructure $\Sup_n({\inladj(X)})$.  To construct a
  comparison between these two \inf-categories, we apply
  $\inladj$ to the inclusion $\overinc\subseteq_r:\Sup_n(X)->X.$ to
  obtain an \inf-functor $\arrow\inladj(\subseteq_r):\inladj(
  \Sup_n(X))->\inladj(X).$ which is characterised by the fact that for
  each simplex $x\in\Sup_n(X)$ it maps the atom
  $\satom{x}\in\inladj(\Sup_n(X))$ to the corresponding atom
  $\satom{x}\in\inladj(X)$. Now applying lemma~\ref{inerv.rest} we
  find that $\inladj(\Sup_n(X))$ is an $n$-category and therefore that
  we may restrict our \inf-functor to a canonical comparison
  $\overarr:\inladj(\Sup_n(X))->\Sup_n( {\inladj(X)}).$.
  
  Indeed, this comparison \inf-functor can be shown to be an
  isomorphism. However, a direct proof of this fact is a little
  involved and for our purposes here all we need is the observation
  that it is an epimorphism. To prove this, start with the fact that
  the set $\satom{X}$ of observation~\ref{inf.cat.nerve.strat} weakly
  generates $\inladj(X)$, from which we may infer immediately that its
  subset $\Sup_n({\inladj(X)})\cap \satom{X}$ of $n$-cells weakly
  generates the superstructure $\Sup_n({\inladj(X)})$. By the
  definition of $\satom{X}$ given in display~(\ref{satom.wgen}) each
  cell in $\Sup_n({\inladj(X)})\cap \satom{X}$ is an atom
  $\satom{x}\in\Sup_n( {\inladj(X)})$ on some non-thin $r$-simplex
  $x\in X$ with $r\leq n$, furthermore each such simplex is in the
  regular subset $\Sup_n(X)\subseteq_r X$ and thus provides an atom
  $\satom{x}\in\inladj(\Sup_n(X))$ which maps to
  $\satom{x}\in\Sup_n({\inladj(X)})$ under our comparison
  \inf-functor.  In other words every one of the cells in our weak
  generator $\Sup_n({\inladj(X)})\cap \satom{X}$ for
  $\Sup_n({\inladj(X)})$ is mapped to by some cell in
  $\inladj(\Sup_n(X))$ under the comparison
  $\overarr:\inladj(\Sup_n(X))->\Sup_n( {\inladj(X)}).$, which is
  therefore an epimorphism by observation~\ref{gen.epi.obs}.
\end{obs}

\begin{thm}[Street~\cite{Street:1988:Fillers}]\label{nerves.are.complicial}
  The \inf-categorical nerve functor $\arrow \inerve:\InfCat->\Strat.$
  takes each \inf-category to a complicial set. It follows that the
  left adjoint $\arrow\inladj:\Strat->\InfCat.$ carries f-extensions
  to isomorphisms and that we may restrict the adjunctions of
  displays~(\ref{strat.nerve}) and~(\ref{inerv.rest.diag}) to:
  \begin{equation*}
    \let\labelstyle=\textstyle
    \xymatrix@R=1ex@C=10em{
      {\InfCat}\ar[r]^{\bot}_{\inerve} &
      {\Comp}\ar@/_3ex/[l]_{\inladj}}\mkern20mu
    \text{and}\mkern20mu
    \xymatrix@R=1ex@C=10em{
      {\nCat}\ar[r]^{\bot}_{\nnerve} &
      {\Comp_n}\ar@/_3ex/[l]_{\nnladj}}
  \end{equation*}
\end{thm}

\begin{proof}
  For a proof of this result see loc.\ cit., which essentially
  proceeds by establishing an optimal decomposition result for the
  cells of an oriental.
  
  It is worth observing, however, that the result Street actually
  proves there looks, on the face of it, to be a little stronger than
  our compliciality condition. In particular his {\em admissible
    horn\/} notion is somewhat more liberal than the one we gave in
  notation~\ref{stan.strat}, since it requires fewer faces to be thin
  for admissibility to hold. In turn, this means that his complicial
  set notion explicitly insists that a strictly greater number of
  horns should have unique fillers. It follows therefore that Street's
  result, which establishes compliciality in his apparently stronger
  sense, immediately implies compliciality in our sense which is
  actually the notion originally proposed by Roberts
  in~\cite{Roberts:1978:Complicial}.
  
  It should be noted that there is no a priori reason to believe that
  Street's compliciality notion isn't strictly stronger than ours.
  However our slightly weaker notion is adequate for all of our work
  in earlier sections and its adoption frees us from having to verify
  Street's more complex admissibility condition in our calculations.
  
  Indeed the equivalence of these notions is actually a simple
  corollary of the ultimate result in this work, which establishes
  Street's primary conjecture in~\cite{Street:1987:Oriental} that the
  nerve functor actually provides us with a canonical equivalence
  between the categories $\InfCat$ of \inf-categories and $\Comp$ of
  complicial sets.
\end{proof}

\begin{obs}[$\inerve$ as a finitely accessible functor]
  \label{inerve.fin.acc}
  It is, of course, the case that $\InfCat$ is an LFP-category. Limits
  and filtered colimits are formed in $\InfCat$ in the familiar way
  for such categories, that is we take the limit or filtered colimit
  of the associated diagram of underlying sets in $\Set$ and then
  enrich that with an \inf-category structure derived point-wise from
  those of the nodes of our original diagram. All finite
  \inf-categories and quotients thereof are finitely presentable so,
  in particular, the orientals $\oriental(\OSimp[n])$ and collapsed
  orientals $\oriental(\OSimp[n])_t$ are all finitely presentable.  It
  follows that the coalgebra $\arrow\oriental\circ\OSimp:\tDelta->
  \InfCat.$, which we used in observation~\ref{inf.cat.nerve.strat} to
  define $\inladj\dashv\inerve$ via Kan's construction, is finitely
  presented and thus that $\arrow \inerve:\InfCat->\Strat.$ preserves
  all filtered colimits (cf.\ observation~\ref{func.from.coalg}).
\end{obs}

\subsection{Relating Parity Complex Products and the Tensor Product of
  Stratified Sets}\label{tensor.rel.sect}

Earlier, we claimed in passing that the tensor product of complicial
sets was the appropriate generalisation of Gray's lax tensor
product~\cite{Gray:1974:FormalCat} to the category of complicial sets.
The main result of this subsection gives this bold claim validity by
demonstrating a strong relationship between our tensor product and 
the free \inf-categories generated by products of parity complexes.

\begin{lemma}
  \label{diag.into.simpprod}
  We may define a canonical parity complex morphism
  $\spanarr\diag:\OSimp->\OSimp\times\OSimp.$ which maps each $r$-simplex
  $\svec{v}$ to the subset: 
  \begin{equation*}
    \diag(\svec{v})\defeq\left\{\pair<\svec{v}\cdot\partinj^{\vec{r}}_1;
      \svec{v}\cdot\partinj^{\vec{r}}_2>\mid \vec{r}\text{ is a
      partition of }r\right\}
  \end{equation*}
  Alternatively, if we express $\svec{v}$ using the vertex-wise
  notation $v_0...v_r$ then we can write this as: 
  \begin{equation*}
    \diag(v_0...v_r)\defeq\left\{
      \pair<v_0...v_s;v_s...v_r>\in\OSimp_s\times\OSimp_{r-s}
      \mid s=0,1,...,r\right\}
  \end{equation*}
\end{lemma}

\begin{proof} 
  Notice first that $\diag$ is indeed a morphism of graded sets, since
  for each $s=0,1,...r$ we know that $\pair<v_0...v_s;v_s...v_r>$ is
  an element of $\OSimp_s\times\OSimp_{r-s}$ which is subset of
  $(\OSimp\times\OSimp)_r\defeq\bigcup_{i=0}^r
  \OSimp_i\times\OSimp_{r-i}$ and so it follows that
  $\diag(v_0...v_r)\subseteq(\OSimp\times\OSimp)_r$. Furthermore,
  we have $\diag(v_0)=\{\pair<v_0;v_0>\}$ and so
  condition~(\ref{func.par.a}) of lemma~\ref{inf.funct.parity} holds
  for our morphism and it remains to prove that
  condition~(\ref{func.par.b}) of that lemma also holds.
  
  To that end, we start by working out what the sets of faces of an
  element $\pair<v_0...v_s;v_s...v_r>$ of $\diag(v_0...v_r)$ look
  like.  Consulting the formulae for the sets of faces of simplices
  and products we see that its faces of parity $\parvtwo$ are
  $\pair<v_0...\hat{v}_i...v_s;v_s...v_r>$ if $i$ is of parity $\parvtwo$
  and $\pair<v_0...v_s;v_s...\hat{v}_i...v_r>$ if $(i-s)$ is of parity
  $\parvtwo(s)$. However in the latter case it is clear that this
  condition also reduces to simply insisting that $i$ is of parity
  $\parvtwo$, therefore we have:
  \begin{equation}\label{faces.of.diagv}
    \begin{split}
      \pair<v_0...v_s;v_s...v_r>^\parvtwo = {} &\left\{\pair<v_0...\hat{v}_i...v_s;
        v_s...v_r>\mid 0<s \also 0\leq i\leq s\also\text{$i$
          is of parity $\parvtwo$}\right\}\cup \\
      & \left\{\pair<v_0...v_s;v_s...\hat{v}_i...v_r>\mid s<r\also s\leq
        i\leq r\also\text{$i$ is of parity $\parvtwo$}\right\}
    \end{split}
  \end{equation}
  Notice also that the first factor in this union is a subset of
  $\OSimp_{s-1}\times\OSimp_{r-s}$ and that its second factor is a
  subset of the disjoint set $\OSimp_s\times\OSimp_{r-s-1}$.
  
  So suppose that we have a distinct pair of indices $s,t\in[r]$,
  assume w.l.o.g.\ that $0\leq s<t\leq r$ and consider the
  intersection $\pair<v_0...v_s;v_s...v_r>^\parvtwo\cap
  \pair<v_0...v_t;v_t...v_r>^\parvone$. We know that the first of the
  factors here is a subset of $(\OSimp_{s-1}\times\OSimp_{r-s})\cup(
  \OSimp_s\times \OSimp_{r-s-1})$ and that its second is a subset of
  $(\OSimp_{t-1}\times\OSimp_{r-t})\cup(\OSimp_t\times
  \OSimp_{r-t-1})$ therefore it follows, by the pairwise disjointness
  of the sets $\OSimp_i\times\OSimp_{r-i-1}$ for $i=0,1,...,r-1$, that
  our intersection is empty if $t>s+1$ and in the case where $t=s+1$
  we have:
  \begin{equation*}
    \begin{split}
      \pair<v_0...v_s;v_s...v_r>^\parvtwo\cap
      \pair<v_0...v_{s+1};v_{s+1}...v_r>^\parvone = {} &
      \left\{\pair<v_0...v_s;v_s...\hat{v}_i...v_r>\mid 
        \text{$i$ is of parity $\parvtwo$}\right\}\cap\\
      & \left\{\pair<v_0...\hat{v}_j...v_{s+1};v_{s+1}...v_r>\mid
        \text{$j$ is of parity $\parvone$}\right\}
    \end{split}
  \end{equation*}
  Notice, however, that those elements  in the first factor which have
  $i>  s$ cannot occur  in the  second factor,  since all  elements in
  there mention the vertex $v_i$  in their second ordinate. Dually any
  element in the  second factor with $j<s+1$ cannot  also occur in the
  first factor.  Thus it follows easily that the only possible element
  in this intersection must be obtained by dropping $v_s$ in the first
  factor and $v_{s+1}$ in the second  factor and thus it is clear that
  we have:
  \begin{equation*}
    \pair<v_0...v_s;v_s...v_r>^\parvtwo\cap
    \pair<v_0...v_{s+1};v_{s+1}...v_r>^\parvone  =
    \left\{
      \begin{array}{lp{7em}}
        \{\pair<v_0...v_s;v_{s+1}...v_r>\} &
        $s$ is of parity $\parvtwo$ and $\parvone=\neg\parvtwo$,\\
        \emptyset & otherwise.
      \end{array}
    \right.
  \end{equation*}
  From this calculation it immediately follows that
  $\diag(v_0...v_r)$ is a well formed subset of
  $\OSimp\times\OSimp$, since the only case which had
  $\pair<v_0...v_s;v_s...v_r>^\parvtwo\cap\pair<v_0...v_t;
  v_t...v_r>^\parvone\neq\emptyset$ also had $\parvone=\neg\parvtwo$, that
  is to say in that case $\parvone$ and $\parvtwo$ were opposing, rather
  than equal, parity symbols.
  
  To prove the movement property for $\diag(v_0....v_r)$ start by
  considering $\diag(v_0...v_r)^\parvtwo$ which is obtained by taking
  the union of the sets in equation~(\ref{faces.of.diagv}) for
  $s=1,2,...,r$ and thus may be written
  \begin{equation}\label{diag(v)xi}
    \begin{split}
      \diag(v_0...v_r)^\parvtwo = {} &
      \left\{\pair<v_0...\hat{v}_i...v_s;v_s...v_r>\mid 
        0\leq i\leq s\leq r\also i\text{ is of parity }\parvtwo\right\}\cup \\
      & \left\{\pair<v_0...v_s;v_s...\hat{v}_i...v_r>\mid 
        0\leq s\leq i\leq r\also i\text{ is of parity }\parvtwo\right\}
    \end{split}
  \end{equation}
  what's more our analysis of $\pair<v_0...v_s;v_s...v_r>^\parvtwo\cap
  \pair<v_0...v_t;v_t...v_r>^\parvone$ in the last paragraph
  demonstrates that:
  \begin{equation}\label{int.diag(v)xi}
    \diag(v_0...v_r)^{-}\cap\diag(v_0...v_r)^{+} =
    \left\{\pair<v_0...v_s;v_{s+1}...v_r>\mid s=0,1,...,r-1\right\}
  \end{equation}
  On the other hand, consider a face $v_0...\hat{v}_i...v_r$ of
  $v_0...v_r$ and observe that we have
  \begin{equation*}
    \diag(v_0...\hat{v}_i...v_r) = 
    \left\{\pair<v_0...\hat{v}_i...v_s;
      v_s...v_r>\mid i< s\leq r\right\}\cup 
    \left\{\pair<v_0...v_s;v_s...\hat{v}_i...v_r>\mid 0\leq s< i
    \right\}
  \end{equation*}
  so it follows that $\diag(v_0...v_r^{\parvtwo})$, which is the union
  of these sets taken over those $i=0,1,...,r$ with parity $\parvtwo$, is
  given by
  \begin{equation}\label{diag(vxi)}
    \begin{split}
      \diag(v_0...v_r^\parvtwo) = {} &
      \left\{\pair<v_0...\hat{v}_i...v_s;v_s...v_r>\mid 
        0\leq i< s\leq r\also i\text{ is of parity }\parvtwo\right\}\cup \\
      & \left\{\pair<v_0...v_s;v_s...\hat{v}_i...v_r>\mid 
        0\leq s< i\leq r\also i\text{ is of parity }\parvtwo\right\}
    \end{split}
  \end{equation}
  which differs from the expression in~(\ref{diag(v)xi}) only to the
  extent that it insists on strict inequalities between $i$ and $s$.
  Notice that each element of $\diag(v_0...v_r^\parvtwo)$ has a vertex
  which is held in common by both of its ordinates, it follows that no
  element of the form $\pair<v_0...v_s;v_{s+1}...v_r>$ can be in there
  so we may infer from
  equations~(\ref{diag(v)xi}),~(\ref{int.diag(v)xi})
  and~(\ref{diag(vxi)}) that we can express $\diag(v_0...v_r)^\parvtwo$
  as the disjoint union:
  \begin{equation*}\label{diagvxi.disj.union}
    \diag(v_0...v_r)^\parvtwo=\diag(v_0...v_r^\parvtwo)\sqcup
    \left\{\pair<v_0...v_s;v_{s+1}...v_r>\mid s=0,1,...,r-1\right\}
  \end{equation*}
  Finally, from this latter expression and
  equation~(\ref{int.diag(v)xi}) it follows immediately that
  \begin{equation*}
    \diag(v_0...v_r^{-}) = \diag(v_0...v_r)^\mp 
    \mkern20mu\text{and}\mkern20mu
    \diag(v_0...v_r^{+}) = \diag(v_0...v_r)^\pm
  \end{equation*}
  and it is trivially the case that any subset $\diag(v_0...v_r)$ moves
  $\diag(v_0...v_r)^\mp$ to $\diag(v_0...v_r)^\pm$.
\end{proof}

\begin{obs}\label{diag.nattrans}
  It is clear that for any $n\in\mathbb{N}$ the parity
  complex morphism of the previous lemma restricts to a morphism
  $\spanarr \diag_n:\OSimp[n]->\OSimp[n]\times \OSimp[n].$ in
  $\Parity$.  Furthermore, it is only slightly less routine to
  demonstrate that these form a natural family, in the sense that for
  each simplicial operator $\arrow\alpha:[m]->[n].$ we have a
  commuting square
  \begin{equation*}
    \xymatrix@R=3em@C=4em{
      {\OSimp[m]}\ar[r]^<>(0.5){\diag_m}
      \ar[d]_{\OSimp(\alpha)} &
      {\OSimp[m]\times\OSimp[m]} 
      \ar[d]^{\OSimp(\alpha)\times\OSimp(\alpha)} \\
      {\OSimp[n]}\ar[r]_<>(0.5){\diag_n} &
      {\OSimp[n]\times\OSimp[n]} 
    }
  \end{equation*}
  in $\Parity$. To demonstrate this we pick an element
  $\svec{v}=v_0...v_r$ of $\OSimp[m]$ and observe that
  $\OSimp(\alpha)(\svec{v})=\emptyset$ iff there is some $i$ such that
  $\alpha(v_i)=\alpha(v_{i+1})$ which in turn holds iff for each
  $s\in[r]$ either $\OSimp(\alpha)(v_0...v_s)$ or
  $\OSimp(\alpha)(v_s...v_r)$ is the empty set. It follows then that
  one leg of this diagram carries $\svec{v}$ to the empty set if and
  only if the other one does. In the common non-empty case, it is now
  a routine matter to show that either way we traverse the
  diagram we get the set $\{\pair<\alpha(v_0...v_s);\alpha(v_s...v_r)>
  \mid s=0,1,...,r\}$.
\end{obs}

\begin{obs}
  The last two results provide us with the tools to construct a
  canonical cocone with vertex $\oriental(\OSimp[n]\times\OSimp[m])$
  under the diagram $\arrow \dgm_{\Delta[n]\times
    \Delta[n]}:\groth(\Delta[n]\times\Delta[m])->\InfCat.$ which we
  used to construct $\inladj(\Delta[n]\times\Delta[m])\in \InfCat$ in
  observation~\ref{inf.cat.nerve}.  
  
  The component of this cocone at the $r$-simplex
  $\pair<\alpha;\beta>\in\Delta[n]\times\Delta[m]$, qua object of
  $\groth(\Delta[n]\times\Delta[m])$, is constructed by applying the
  functor $\oriental$ to the following composite of parity complex
  morphisms:
  \begin{equation*}
    \let\labelstyle=\textstyle
    \xymatrix@R=1ex@C=8em{
      {\OSimp[r]}\ar[r]^<>(0.5){\diag_{r}} &
      {\OSimp[r]\times\OSimp[r]}
      \ar[r]^<>(0.5){\OSimp(\alpha)\times\OSimp(\beta)} &
      {\OSimp[m]\times\OSimp[n]}
    }
  \end{equation*}
  Of course we need to check that these collectively satisfy the
  naturality condition for such a cocone. So suppose that the
  simplicial operator $\arrow \gamma:[s]->[r].$ represents an arrow in
  $\groth(\Delta[n]\times\Delta[m])$ from the $s$-simplex
  $\pair<\alpha';\beta'>$ to the $r$-simplex $\pair<\alpha;\beta>$ and
  observe that we have the commutative diagram
  \begin{equation*}
    \xymatrix@R=1.25em@C=2.5em{
      {\OSimp[s]}\ar[rr]^<>(0.5){\diag_s}\ar[dd]_{\OSimp(\gamma)} &&
      {\OSimp[s]\times\OSimp[s]}
      \ar[dr]^<>(0.75)*+{\labelstyle\OSimp(\alpha')\times\OSimp(\beta')} 
      \ar[dd]|{\OSimp(\gamma)\times\OSimp(\gamma)}& \\
      & & & {\OSimp[n]\times\OSimp[m]} \\
      {\OSimp[r]}\ar[rr]_<>(0.5){\diag_r}&& 
      {\OSimp[r]\times\OSimp[r]}
      \ar[ur]_<>(0.75)*+{\labelstyle\OSimp(\alpha)\times\OSimp(\beta)} & \\
    }
  \end{equation*}
  in $\Parity$. Here the square on the left expresses the naturality
  property discussed in observation~\ref{diag.nattrans} and the
  triangle on the right commutes by the functoriality of $\OSimp({-})
  \times\OSimp({-})$ applied to the defining property of $\gamma$ as
  an arrow of $\groth(\Delta[n]\times\Delta[m])$, that being
  $\pair<\alpha';\beta'>=\pair<\alpha;\beta>\cdot\gamma=\pair<\alpha
  \circ\gamma; \beta\circ\gamma>$. So applying the functor $\oriental$
  to this diagram we get the required cocone naturality property and
  it follows that this information induces a comparison \inf-functor
  $\arrow\cprprod^{n,m}: \inladj(\Delta[n]\times\Delta[m])->
  \oriental(\OSimp[n]\times \OSimp[m]).$ which is uniquely determined
  by the property that the square
  \begin{equation*}
    \xymatrix@=2.5em{
      {\oriental(\OSimp[r])}\ar[r]^<>(0.5){\oriental(\diag_r)}
      \ar[d]_{\coliminj_{\pair<\alpha;\beta>}} &
      {\oriental(\OSimp[r]\times\OSimp[r])}
      \ar[d]^{\oriental(\OSimp(\alpha)\times\OSimp(\beta))} \\
      {\inladj(\Delta[n]\times\Delta[m])}
      \ar[r]_{\cprprod^{n,m}} &
      {\oriental(\OSimp[n]\times\OSimp[m])}
    }
  \end{equation*}
  commutes for each simplex $\pair<\alpha;\beta>$ in
  $\Delta[n]\times\Delta[m]$ (where $r=\dim\pair<\alpha;\beta>$).
  Alternatively, using the notation of observation~\ref{fg.inladj.obs},
  we see that we may characterise $\cprprod^{n,m}$ as the unique
  \inf-functor which has
  \begin{equation}\label{cprprod.char}
    \cprprod^{n,m}(\satomp<\alpha;\beta>)= \oriental(
    (\OSimp(\alpha)\times\OSimp(\beta))\circ\diag_r)(\atom{01...r})
  \end{equation}
  for each simplex $\pair<\alpha;\beta>$ in $\Delta[n]\times\Delta[m]$
  (where $r=\dim\pair<\alpha;\beta>$). 
  
  Finally it is easily seen that the family of \inf-functors
  $\cprprod^{n,m}$ is natural, in the sense that the square
  \begin{equation}\label{cprprod.nat}
    \xymatrix@=3em{
       {\inladj(\Delta[n]\times\Delta[m])}
      \ar[r]^{\cprprod^{n,m}} 
      \ar[d]_{\inladj(\Delta(\phi)\times\Delta(\psi))} &
      {\oriental(\OSimp[n]\times\OSimp[m])} 
      \ar[d]^{\oriental(\OSimp(\phi)\times\OSimp(\psi))}\\
       {\inladj(\Delta[r]\times\Delta[s])}
      \ar[r]_{\cprprod^{r,s}} &
      {\oriental(\OSimp[r]\times\OSimp[s])} \\
    }
  \end{equation}
  commutes for each pair of simplicial operators
  $\arrow\phi:[n]->[r].$ and $\arrow\psi:[m]->[s].$. We leave it up to
  the reader to verify this claim, either by a simple diagram chase
  and an appeal to the uniqueness property of our defining colimit or
  by considering the action of these \inf-functors on the atoms of
  $\inladj(\Delta[n]\times\Delta[m])$ and using the explicit
  descriptions given in equation~(\ref{cprprod.char}) and
  observations~\ref{osimp(alpha).act.on.atoms}
  and~\ref{fg.inladj.obs}.
\end{obs}

\begin{obs}\label{prod.simp.anal}
  In section~\ref{precomp.sec} we singled our two classes of simplices
  in $\Delta[n]\times\Delta[m]$, cylinders and mediators, which we
  found to be of great interest in our analysis of the stratified set
  $\Delta[n]\otimes\Delta[m]$.  On order to extend our work in the
  last observation to the \inf-category $\inladj(\Delta[n]\otimes
  \Delta[m])$, we are lead to consider how the \inf-functor $\arrow
  \cprprod^{n,m}:\inladj(\Delta[n]\times \Delta[m])
  ->\oriental(\OSimp[n]\times\OSimp[m]).$ acts on those atoms which
  correspond to simplices in these classes.
  
  In the sequel we'll be interested in analysing cells
  $\cprprod^{n,m}(\satomp<\alpha;\beta>)$ where $\pair<\alpha;\beta>$
  is an $r$-simplex of $\Delta[n]\times\Delta[m]$ and our primary tool
  in making such calculations will the expression in
  equation~(\ref{cprprod.char}). It is therefore worth observing that
  we may combine the definitions in
  observations~\ref{parity.prod.func}, \ref{osimp.functor}
  and~\ref{diag.into.simpprod} to show that we have:
  \begin{equation}\label{expl.diag.calc}
    ((\OSimp(\alpha)\times\OSimp(\beta))\circ\diag_r)(01...r)
    =\{\pair<\alpha(01...s);\beta(s(s+1)...r)>\mid s=0,1,...,r\}\cap
    (\OSimp[n]\times\OSimp[m])_r
  \end{equation}
  In other words, we form $\diag_r(01...r)$ then apply $\alpha$ and
  $\beta$ to respective components of successive pairs in this set and
  finally we intersect the result with $(\OSimp[n]\times\OSimp[m])_r$
  to exclude any elements which ended up having dimension less than
  $r$.

  \vspace{1ex}

  \noindent{\bf Cylinders:} 
  Firstly let $r=n+m$ and consider the unique non-degenerate
  $r$-dimensional cylinder $\pair<\partproj^{n,m}_1;
  \partproj^{n,m}_2>$ in $\Delta[n]\times \Delta[m]$. Consulting the
  definition of these partition operators it is easy to see that
  \begin{equation*}
    \partproj^{n,m}_1(01...n)=01...n\mkern20mu\also\mkern20mu
    \partproj^{n,m}_2(n(n+1)...r)=01...m
  \end{equation*}
  so we see, by equation~(\ref{expl.diag.calc}), that
  $\pair<01...n;01...m>\in((\OSimp(\partproj^{n,m}_1)
  \times\OSimp(\partproj^{n,m}_2))\circ\diag_r)(01...r)$.  It follows
  therefore that $\cprprod^{n,m}(\satomp<\partproj^{n,m}_1;
    \partproj^{n,m}_2>)=\oriental((\OSimp(\partproj^{n,m}_1)
  \times\OSimp(\partproj^{n,m}_2))\circ\diag_r)(\atom{01...r})$ has
  non-empty sets of $r$-dimensional elements and is thus a non-trivial
  $r$-cell in $\oriental(\OSimp[n]\times\OSimp[m])$. However we know,
  by the comment at the end of
  observation~\ref{parity.prods.and.cubes}, that the atom
  $\atomp<01...n;01...m>$ is the unique non-trivial $r$-cell of
  $\oriental(\OSimp[n]\times\OSimp[m])$ and we must have
  $\cprprod^{n,m}(\satomp<\partproj^{n,m}_1;
  \partproj^{n,m}_2>)=\atomp<01...n;01...m>$.
  
  Now, referring back to observation~\ref{med.cyl}, by definition we
  know that the partition $\pair<p;q>$ of $r$ witnesses the
  $r$-simplex $\pair<\alpha;\beta>\in\Delta[n]\times\Delta[m]$ to be a
  cylinder iff there are (necessarily unique) simplicial operators
  $\arrow\alpha':[p]->[n].$ and $\arrow\beta':[q]->[m].$ with
  $\alpha=\alpha'\circ\partproj^{p,q}_1$ and $\beta=\beta'\circ
  \partproj^{p,q}_2$. Furthermore, consulting the definition our
  partition operators $\partproj^{p,q}_i$, it is clear that
  $\pair<\alpha;\beta>$ is non-degenerate iff $\alpha'$ and $\beta'$
  are both face operators.  In that case, we may define
  $\svec{v}=\alpha(01...r)$ and $\svec{w}=\beta(01...r)$, observe that
  these are elements of $\OSimp[n]_p$ and $\OSimp[m]_q$ respectively
  and easily show that $\alpha'=\face^n_{\svec{v}}$ and
  $\beta'=\face^m_{\svec{w}}$. Finally, the naturality property
  depicted in display~(\ref{cprprod.nat}), when applied to the
  operators $\face^n_{\svec{v}}$ and $\face^m_{\svec{w}}$ and
  evaluated at the atom
  $\satomp<\partproj^{p,q}_1;\partproj^{p,q}_2>$, demonstrates that
  \begin{equation*}
    \cprprod^{n,m}(\satomp<\alpha;\beta>)=
    \cprprod^{n,m}(\satomp<\face^n_{\svec{v}}\circ
      \partproj^{p,q}_1; \face^m_{\svec{w}}\circ\partproj^{p,q}_2>)=
    \oriental(\OSimp(\face^n_{\svec{v}})\times\OSimp(\face^m_{\svec{w}}))
    (\cprprod^{p,q}(\satomp<\partproj^{p,q}_1;\partproj^{p,q}_2>))
  \end{equation*}
  but we've just shown that 
  \begin{equation*}
    \cprprod^{p,q}(\satomp<
      \partproj^{p,q}_1;\partproj^{p,q}_2>)=\atomp<01...p;01...q>
  \end{equation*}
  so it follows that 
  \begin{equation*}
    \begin{split}
    \cprprod^{n,m}(\satomp<\alpha;\beta>)&=
      \oriental(\OSimp(\face^n_{\svec{v}})\times\OSimp(\face^m_{\svec{w}}))(
      \atomp<01...p;01...q>) \\
      &=\atomp<\face^n_{\svec{v}}(01...p);\face^m_{\svec{w}}(01...q)>\\
      &=\atomp<\svec{v};\svec{w}>
    \end{split}
  \end{equation*}
  where the second equality follows from
  observation~\ref{osimp(alpha).act.on.atoms} and the last one follows
  directly from the relationship between $\svec{v}$, $\svec{w}$ and
  the corresponding face operators $\face^n_{\svec{v}}$,
  $\face^m_{\svec{w}}$. 
  
  In summary, we have demonstrated that $\cprprod^{n,m}$ restricts to
  a bijection between the set of $r$-dimensional atoms in
  $\oriental(\OSimp[n]\times\OSimp[m])$ and the set of atoms
  $\satomp<\alpha;\beta>$ in $\inladj(\Delta[n]\times\Delta[m])$ for
  which $\pair<\alpha;\beta>$ is a non-degenerate cylinder in
  $\Delta[n]\times\Delta[m]$.

  \vspace{1ex}

  \noindent{\bf Mediators:} Consider an arbitrary $r$-simplex
  $\pair<\alpha;\beta>\in\Delta[n]\times\Delta[m]$ and observe that
  equations~(\ref{cprprod.char}) and~(\ref{expl.diag.calc}) imply that
  $\cprprod^{n,m}(\satomp<\alpha;\beta>)$ is an $(r-1)$-cell
  iff $\pair<\alpha(01...s);\beta(s(s+1)...r)>\notin
  (\OSimp[n]\times\OSimp[m])_r$ (for all $s=0,1,...,r$). Of course,
  this condition holds if and only if either $\alpha(01...s)$ has
  dimension $<s$ (that is $\exists\,0\leq i<s$ s.t.\ 
  $\alpha(i)=\alpha(i+1)$) or $\beta(s(s+1)...r)$ has dimension
  $<(r-s)$ (in other words $\exists\,s\leq j<r$ s.t.\ 
  $\beta(j)=\beta(j+1)$).  
  
  In fact, if we let $s\leq r$ be the largest integer such that
  $\alpha(0)<\alpha(1)<...<\alpha(s)$ then this latter condition
  implies that there is some $s\leq j<r$ such that
  $\beta(j)=\beta(j+1)$ and then, by the maximality of $s$, we also
  have $\alpha(s)=\alpha(s+1)$. Conversely, if we had some $i\leq j$
  such that $\alpha(i)=\alpha(i+1)$ and $\beta(j)=\beta(j+1)$ then it
  is clear that the triviality condition of the last paragraph
  immediately holds.
  
  Notice that we've simply demonstrated that $\cprprod^{n,m}(
    \satomp<\alpha;\beta>)$ is an $(r-1)$-cell if and only if the
  condition of lemma~\ref{deg.ord}(\ref{deg.ord.both}) holds and we
  know, by corollary~\ref{tensor.simp.set}, that this characterises
  the thin simplices in the tensor product
  $\Delta[n]\otimes\Delta[m]$.
  
  Of course, given that this result builds a strong and compelling
  link to our work on tensor products of stratified sets in
  sections~\ref{precomp.sec} and~\ref{comp.sec}, it would be tempting
  to claim that it was this observation which motivated our original
  definition of that tensor product. However this is by no means the
  case, since we arrived at our original definition of
  $\Delta[n]\otimes\Delta[m]$ through purely simplicial means and only
  later derived its close relationship to the free \inf-category
  $\oriental(\OSimp[n]\times\OSimp[m])$.
\end{obs}

These observations lead directly to, and are summarised by, the
primary result of this subsection:

\begin{thm}\label{tensor.simp}
  The \inf-functor $\arrow\cprprod^{n,m}:\inladj(\Delta[n]\times
  \Delta[m])->\oriental(\OSimp[n]\times\OSimp[m]).$ factors through
  the canonical collapser $\arrow q_t:\inladj(\Delta[n]\times
  \Delta[m])->\inladj(\Delta[n]\otimes \Delta[m]).$ and provides a
  comparison \inf-functor
  \begin{equation*}
    \let\labelstyle=\textstyle
    \xymatrix@R=1ex@C=8em{
      {\inladj(\Delta[n]\otimes\Delta[m])}
      \ar[r]^{\cprprods^{n,m}} &
      {\oriental(\OSimp[n]\times\Delta[m])}
    }
  \end{equation*}
  which inherits the naturality properties of $\cprprod^{n,m}$.  If
  $\pair<p;q>$ is the partition which witnesses the non-degenerate
  $r$-simplex $\pair<\alpha;\beta>$ as a cylinder then we have
  elements $\svec{v}=\alpha(01...r)$ in $\OSimp[n]_p$ and
  $\svec{w}=\beta(01...r)$ in $\OSimp[m]_q$ (for which
  $\alpha=\face^n_{\svec{v}}\circ\partproj^{p,q}_1$ and
  $\beta=\face^m_{\svec{w}}\circ\partproj^{p,q}_2$) and
  $\cprprods^{n,m}$ maps the cell $\satomp<\alpha;\beta>$ in
  $\inladj(\Delta[n]\otimes\Delta[m])$ to the atom
  $\atomp<\svec{w};\svec{v}>$ in $\oriental(\OSimp[n]\times
  \OSimp[m])$. It follows that $\cprprods^{n,m}$ restricts to a
  bijection between each set
  \begin{equation*}
    C^{n,m}_r\defeq\left\{
      \satomp<\alpha;\beta>\mid \pair<\alpha;\beta>\text{ is a
      non-degenerate $r$-cylinder in }\Delta[n]\otimes\Delta[m]
    \right\}
  \end{equation*}
  and the set of $r$-dimensional atoms of
  $\oriental(\OSimp[n]\times\OSimp[m])$.  Finally,
  $\inladj(\Delta[n]\otimes\Delta[m])$ is weakly generated by the
  cells in the set $C^{n,m}$ and consequently we may apply
  lemma~\ref{fg.isom.lemma} to show that $\cprprods^{n,m}$ is an
  {\bf\em isomorphism\/} of \inf-categories.
\end{thm}

\begin{proof}
  The latter part of the last observation demonstrated that the
  \inf-functor $\cprprod^{n,m}$ maps atoms $\satomp<\alpha;\beta>$
  associated with those $r$-simplices $\pair<\alpha;\beta>\in\Delta[n]
  \times \Delta[m]$ which are thin in $\Delta[n]\otimes\Delta[m]$ to
  $(r-1)$-cells in $\oriental(\OSimp[n]\otimes\OSimp[m])$. In other
  words, $\cprprod^{n,m}$ collapses the cells in the graded set
  $\satom{t(\Delta[n]\otimes\Delta[m])}\subseteq
  \inladj(\Delta[n]\times\Delta[m])$ and it therefore factors through
  the quotient $\inladj(\Delta[n]\otimes\Delta[m])$ as postulated.
  Checking the naturality properties of the resulting \inf-functors
  $\cprprods^{n,m}$ is a routine matter which we leave up to the
  reader.  Notice also that the first half of the last observation
  demonstrated that $\cprprod^{n,m}$, and thus $\cprprods^{n,m}$, acts
  on cells $\satomp<\alpha;\beta>$ associated with non-degenerate
  cylinders as stated and thus that it restricts to a bijection
  between each $C^{n,m}_r$ and the set of $r$-dimensional atoms of
  $\oriental(\OSimp[n]\times\OSimp[m])$.
  
  So all that remains is to apply induction on $(n+m)$ to prove that
  $C^{n,m}$ weakly generates the \inf-category
  $\inladj(\Delta[n]\otimes\Delta[m])$. The base case is trivial, so
  we adopt the inductive hypothesis, that $C^{p,q}$ weakly generates
  $\inladj(\Delta[p]\otimes\Delta[q])$ whenever $p+q < n+m$, and
  consider family of stratified maps
  \begin{equation}\label{bdary.family.1}
    \begin{aligned}
      & \arrow\Delta(\face^n_i)\otimes\Delta[m]:\Delta[n-1]\otimes\Delta[m]
      ->\Delta[n]\otimes\Delta[m]. && \text{for $i=0,1,...,n$ and}\\
      & \arrow\Delta[n]\otimes\Delta(\face^m_j):\Delta[n]\otimes\Delta[m-1]
      ->\Delta[n]\otimes\Delta[m]. && \text{for $j=0,1,...,m$.}
    \end{aligned}
  \end{equation}
  each of which factors through the regular boundary subset
  $\boundary(\Delta[n]\otimes\Delta[m])\subseteq_r
  \Delta[n]\otimes\Delta[m]$ of
  notation~\ref{Hdelta}(\ref{Hdelta.bdary}) to provide a jointly
  epimorphic family of stratified maps with codomain
  $\boundary(\Delta[n]\otimes \Delta[m])$ The functor $\inladj$
  preserves coproducts and epimorphisms, since it is left adjoint, so
  applying it to this family we obtain a jointly epimorphic family of
  \inf-functors with codomain $\inladj(\boundary(\Delta[n]\otimes
  \Delta[m]))$.
  
  Clearly $\boundary(\Delta[n]\otimes\Delta[m])$ is actually a regular
  subset of the $(n+m-1)$-superstructure $\Sup_{n+m-1}(\Delta[n]
  \otimes\Delta[m])$ and indeed it is easily seen that if
  $H(\Delta[n],\Delta[m])\subseteq_r\Delta[n]\otimes\Delta[m]$ is the
  regular subset defined in notation~\ref{Hdelta}(\ref{Hdelta.H(N,M)})
  then $\Sup_{n+m-1}(\Delta[n]\otimes\Delta[m])=\boundary(\Delta[n]
  \otimes\Delta[m])\cup H(\Delta[m],\Delta[m])$ so we have a pasting
  square:
  \begin{equation*}
    \xymatrix@=3em{
      {\boundary H(\Delta[n],\Delta[m])}\ar@{_{(}->}[d]_{\subseteq_r}
      \ar@{^{(}->}[r]^{\subseteq_r} &
      {H(\Delta[n],\Delta[m])}\ar@{^{(}->}[d]^{\subseteq_r} \\
      {\boundary(\Delta[n]\otimes\Delta[m])}
      \ar@{^{(}->}[r]_{\subseteq_r} &
      {\Sup_{\rlapm{n+m-1}}(\Delta[n]\otimes\Delta[m])}\poexcursion
    }
  \end{equation*}
  By lemma~\ref{major.fill}, we know that the upper horizontal in this
  square is an f-extension and so the lower horizontal inclusion is
  also an f-extension (by
  observation~\ref{pasting.sq}(\ref{pasting.sq.b})). It follows, from
  theorem~\ref{nerves.are.complicial}, that this inclusion is
  perpendicular to the nerve $\inerve(\mathbb{C})$ of each
  \inf-category or equivalently that the \inf-functor
  $\arrow\inladj(\subseteq_r):\inladj(\boundary(\Delta[n]\otimes
  \Delta[m]))->\inladj(\Sup_{n+m-1}(\Delta[n]\otimes\Delta[m])).$ is
  an isomorphism.
  
  We also know, from observation~\ref{inladj.sup.compare}, that the
  regular inclusion $\Sup_{n+m-1}(\Delta[n]\otimes \Delta[m])
  \subseteq_r \Delta[n]\otimes\Delta[m]$ provides an \inf-functor
  $\arrow\inladj(\subseteq_r):\inladj(\Sup_{n+m-1}(
  \Delta[n]\otimes\Delta[m]))->\inladj(\Delta[n]\otimes\Delta[m]).$
  which restricts to an epimorphism with codomain
  $\Sup_{n+m-1}({\inladj(\Delta[n]\otimes \Delta[m])})$. Observe now
  that the functoriality of $\inladj$ implies that the family obtained
  by composing this latter epimorphic \inf-functor with the composite
  of the isomorphism of the previous paragraph and the jointly
  epimorphic family of the paragraph before that gives a family of
  \inf-functors which could otherwise be constructed by applying
  $\inladj$ to the stratified maps in display~\ref{bdary.family.1} and
  factoring each one through $\Sup_{n+m-1}({\inladj(\Delta[n]\otimes
    \Delta[m])}) \subseteq \inladj(\Delta[n]\otimes\Delta[m])$. As the
  composite of an epimorphism, an isomorphism and a jointly epimorphic
  family this too is a jointly epimorphic family.
  
  If $\pair<\alpha;\beta>\in\Delta[n-1]\otimes\Delta[m]$ is a
  non-degenerate cylinder then $\inladj(\Delta(\face^n_i)
  \otimes\Delta[m])(\satomp<\alpha;\beta>) =\satomp<\face^n_i\circ
  \alpha;\beta>$ and $\pair<\face^n_i\circ\alpha;\beta>$ is a
  non-degenerate cylinder in $\Delta[n]\otimes\Delta[m]$, in other
  words the \inf-functor $\inladj(\Delta(\face^n_i)\otimes\Delta[m])$
  maps cells in $C^{n-1,m}$ into $\Sup_{n+m-1}(C^{n,m})$. Similarly,
  each \inf-functor $\inladj(\Delta[n]\otimes\Delta(\face^m_j))$ maps
  cells in $C^{m,n-1}$ into $\Sup_{n+m-1}(C^{n,m})$. Consequently, we
  may apply lemma~\ref{gen.epi}, and the associated
  observation~\ref{gen.epi.obs}, to the jointly epimorphic family of
  the last paragraph and the weakly generating sets
  $C^{n-1,m}\subseteq\inladj(\Delta[n-1]\otimes\Delta[m])$ and
  $C^{n,m-1}\subseteq\inladj(\Delta[n]\otimes\Delta[m-1])$ of the
  inductive hypothesis to show that $\Sup_{n+m-1}(C^{n,m})$ weakly
  generates $\Sup_{n+m-1}( {\inladj(\Delta[n]\otimes\Delta[m])})$.
  
  Finally, the only cell in the standard weak generator
  $\satom{\Delta[n] \otimes\Delta[m]}$ of
  $\inladj(\Delta[n]\otimes\Delta[m])$ (cf.\ 
  observation~\ref{inf.cat.nerve.strat}) which is not in
  $\Sup_{n+m-1}({\inladj( \Delta[n]\otimes\Delta[m])})$ is the atom
  $\satomp<\partproj^{n,m}_1; \partproj^{n,m}_2>$ because all
  simplices of $\Delta[n]\otimes\Delta[m]$ above dimension $(n+m)$ are
  degenerate and its only non-thin $(n+m)$-simplex is the minimal
  shuffle $\pair<\partproj^{n,m}_1; \partproj^{n,m}_2>$. Therefore
  since $\Sup_{n+m-1}(C^{n,m})$ weakly generates $\Sup_{n+m-1}(
  {\inladj(\Delta[n]\otimes\Delta[m])})$ we may infer that
  $C^{n,m}=\Sup_{n+m-1}(C^{n,m})\cup\{\satomp<\partproj^{n,m}_1;
  \partproj^{n,m}_2>\}$ weakly generates $\inladj(\Delta[n]\otimes
  \Delta[m])$ as required.
\end{proof}

\begin{cor}\label{tensor.simp.th}
  For each pair of stratified sets $N$ and $M$ in $\Simplex$ the
  isomorphism of theorem~\ref{tensor.simp} induces an isomorphism 
  $\arrow\cprprods^{N,M}:\inladj(N\otimes M)->\oriental(\tilde{N}
  \times\tilde{M}).$ which is uniquely determined by the fact that it
  makes the square
  \begin{equation}\label{cprprodsNM.square}
    \let\labelstyle=\textstyle
    \xymatrix@R=2.5em@C=5em{
      {\inladj(\Delta[n]\otimes\Delta[m])}
      \ar[r]^{\cprprods^{n,m}}_{\cong}
      \ar@{-|>}[d]_{\inladj(\subseteq_e)} &
      {\oriental(\OSimp[n]\times\OSimp[m])}
      \ar@{-|>}[d]^{\oriental(\subseteq_e)} \\
      {\inladj(N\otimes M)}
      \ar@{-->}[r]_{\exists!\,\cprprods^{N,M}}^{\cong} &
      {\oriental(\tilde{N}\times\tilde{M})}
    }
  \end{equation}
  commute. In other words, $\cprprods^{N,M}$ is characterised by the
  fact that it maps the cell $\satomp<\alpha;\beta>\in\inladj(N\otimes
  M)$ associated with a non-thin $r$-cylinder $\pair<\alpha;\beta>\in
  N\otimes M$ to the cell $\atomp<\svec{v};\svec{w}>\in\oriental(
  \tilde{N}\times\tilde{M})$ with $\svec{v}=\alpha(01...r)$ and
  $\svec{w}=\beta(01...r)$.  It follows, therefore, that this family
  of isomorphisms is natural in $N,M\in\Simplex$.
\end{cor}

\begin{proof}
  First observe that both of the vertical \inf-functors in
  display~\eqref{cprprodsNM.square} are collapsers.  The right hand
  one is simply the collapser we used to construct
  $\oriental(\tilde{N}\times\tilde{M})$ from the free \inf-category
  $\oriental(\OSimp[n]\times\OSimp[m])$ in
  observation~\ref{sparity.ext.oriental}, in other words this is the
  collapser of the graded subset of this latter \inf-category of those
  atoms associated with elements which are thin in
  $\tilde{N}\times\tilde{M}$.  The left hand one features in the
  commutative square
  \begin{equation*}
    \let\labelstyle=\textstyle
    \xymatrix@R=2.5em@C=5em{
      {\inladj(\Delta[n]\pretens\Delta[m])}
      \ar[r]^{\inladj(\subseteq_e)}_{\cong}
      \ar@{-|>}[d]_{\inladj(\subseteq_e)} &
      {\inladj(\Delta[n]\otimes\Delta[m])}
      \ar@{-|>}[d]^{\inladj(\subseteq_e)} \\
      {\inladj(N\pretens M)}
      \ar[r]_{\inladj(\subseteq_e)}^{\cong} &
      {\inladj(N\otimes M)}
    }
  \end{equation*}
  where the horizontal \inf-functors are obtained by applying
  $\inladj$ to the t-extensions of lemma~\ref{pretens.density}, and
  are therefore isomorphisms by Street's
  theorem~\ref{nerves.are.complicial}. As discussed in
  observation~\ref{inf.cat.nerve.strat} the left hand vertical in this
  square, and thus its isomorphic right hand vertical, is a collapser
  of the set of cells associated with those simplices which are thin
  in $N\pretens M$ but not thin in $\Delta[n]\pretens\Delta[m]$.
  Returning to definition~\ref{pre.tens}, we see that the simplices
  which satisfy this condition are precisely the non-degenerate
  cylinders that are thin in $N\otimes M$.
  
  Now we know, from theorem~\ref{tensor.simp}, that the isomorphism
  $\cprprods^{n,m}$ at the top of display~\eqref{cprprodsNM.square}
  maps the cell $\satomp<\alpha;\beta>\in\inladj(
  \Delta[n]\otimes\Delta[m])$ associated with a non-degenerate
  cylinder $\pair<\alpha;\beta>$ to a corresponding atom
  $\atomp<\svec{v};\svec{w}>\in\oriental(\OSimp[n]\times\OSimp[m])$
  and that these are related by $\alpha=\face^n_{\svec{w}}\circ
  \partproj^{p,q}_1$ and $\beta=\face^m_{\svec{v}}\circ
  \partproj^{p,q}_2$.  Consulting observation~\ref{med.cyl}, we see
  that our cylinder is thin in $N\otimes M$ if either
  $\face^n_{\svec{v}}$ is thin in $N$ or $\face^m_{\svec{w}}$ is thin
  in $M$ and consequently we may infer, from the definitions of the
  corresponding stratified parity complexes $\tilde{N}$, $\tilde{M}$
  and $\tilde{N}\times\tilde{M}$, that our non-degenerate cylinder
  $\pair<\alpha;\beta>$ is thin in $N\otimes M$ iff the corresponding
  element $\pair<\svec{v};\svec{w}>$ is thin in
  $\tilde{N}\times\tilde{M}$.  It follows that $\cprprods^{n,m}$ sets
  up a graded bijection between the sets of cells with respect to
  which the vertical \inf-functors of
  display~\eqref{cprprodsNM.square} are collapsers and thus that it
  induces the isomorphism at the bottom of that square as postulated.
\end{proof}

\subsection{An Inductive Proof of the Street-Roberts Conjecture}

\begin{obs}[an outline of our proof of the Street-Roberts conjecture] 
  \label{st.rob.prf.outline}
  We now have almost all of the technical tools required to complete
  our proof of the Street-Roberts conjecture, which is that the nerve
  functor $\arrow\inerve:\InfCat->\Comp.$ is an equivalence of
  categories. Our approach will be to inductively prove that it
  restricts to an equivalence $\arrow\nerv_n:\nCat->\Comp_n.$ for each
  $n\in\mathbb{N}$ and then use a simple colimiting argument to extend
  those equivalences to $\arrow\inerve:\InfCat->\Comp.$

  To complete this proof the only real work left is to show that we
  may construct an essentially commutative square
  \begin{equation}\label{path.cons.comp}
    \let\labelstyle=\textstyle
    \xymatrix@R=2.5em@C=3em{
      {\InfCat}\ar[r]^<>(0.5){\pathecat}
      \ar[d]_{\inerve}\ar@{}[dr]|{\cong} &
      {\InfCatCat} \ar[d]^{\ECat\inerve} \\
      {\Comp}\ar[r]_<>(0.5){\pathecat\rhv} & {\CompCat}
    }
  \end{equation}
  in which the horizontal functor along the top is the equivalence
  discussed in observation~\ref{infcat.enr} and the functor along the
  bottom is the right handed version of the path category construction
  discussed in sections~\ref{comp.pathcat.sec}
  and~\ref{pathecat.ff.sec}. Furthermore, the right vertical here is
  the functor obtained by applying the (product preserving) nerve
  functor $\inerve$ point-wise to the hom-\inf-categories of each
  $\InfCat$-enriched category.
  
  Once we have constructed this square, the remainder of the proof is
  straightforward. First we combine lemma~\ref{inerv.rest},
  lemma~\ref{sup.comp.e.cat} and the discussion in
  observation~\ref{infcat.enr} to show that for each $n\in\mathbb{N}$
  the square in display~(\ref{path.cons.comp}) restricts to:
  \begin{equation}\label{path.cons.comp.2}
    \let\labelstyle=\textstyle
    \xymatrix@R=2.5em@C=3em{
      {\ECat{(n+1)}}\ar[r]^<>(0.5){\pathecat}
      \ar[d]_{\nerv_{n+1}}\ar@{}[dr]|{\cong} &
      {\ECat{(\nCat)}} \ar[d]^{\ECat{\nnerve}} \\
      {\Comp_{n+1}}\ar[r]_<>(0.5){\pathecat\rhv} & {\ECat{\Comp_n}}
    }
  \end{equation}
  Applying induction we use this to show that $\inerve$ restricts
  to an equivalence $\equiv\nnerve:\nCat-> \Comp_n.$ for each
  $n\in\mathbb{N}$. The base case $n=0$ follows immediately from
  lemma~\ref{0.comp.eq.set}. For the inductive step we adopt the
  hypothesis that $\arrow\nnerve:\nCat->\Comp_n.$ is an equivalence,
  from which it is immediate that right hand vertical functor
  $\ECat\nnerve$ of display~(\ref{path.cons.comp.2}) is an
  equivalence.  Furthermore, as remarked in
  observation~\ref{infcat.enr}, the upper horizontal arrow $\pathecat$
  in this square is also an equivalence, so it follows that its upper
  right composite, and thus its isomorphic lower left composite, is an
  equivalence. However we know, by theorem~\ref{ff.pathecat}, that the
  lower horizontal functor $\pathecat\rhv$ is fully faithful, which
  allows us to apply lemma~\ref{ff.comp}(\ref{ff.comp.2}) and
  demonstrate that the left hand vertical $\nerv_{n+1}$ is an
  equivalence as required.
  
  Of course, this result immediately implies that the adjunction
  $\arrow \nnladj \dashv\nnerve:\nCat->\Comp_n.$ is an adjoint
  equivalence (cf.~\cite{Maclane:1971:CWM}) for each $n\in\mathbb{N}$.
  The unit and counit of each such adjunction is obtained by
  restricting the unit $\nattrans \eta:\id_{\Comp}->
  \inerve\circ\inladj.$ and counit $\nattrans
  \epsilon:\inladj\circ\inerve-> \id_{\InfCat}.$ of the adjunction
  $\arrow \inladj\dashv\inerve:\InfCat-> \Comp.$. So it follows that
  each $\nnladj\dashv\nnerve$ is an adjoint equivalence iff the unit
  component $\eta_A$ is an isomorphism whenever $A$ is $n$-complicial
  and the counit component $\epsilon_{\mathbb{C}}$ is an isomorphism
  whenever $\mathbb{C}$ is an $n$-category.
  
  To extend this result to all components of $\eta$ and $\epsilon$, we
  exploit the fact that each \inf-category or complicial set is equal
  to the union of its ascending sequence of superstructures. To
  elaborate, for each \inf-category $\mathbb{C}$ the cocone of
  inclusions $\inc i_n:\Sup_n(\mathbb{C})->\mathbb{C}.$ displays
  $\mathbb{C}$ as the colimit $\colim_{n\in\mathbb{N}}\Sup_n(
  \mathbb{C})$ of the ascending chain
  \begin{math}
    \xymatrix@R=0ex@C=2em@M=1pt{
      {\Sup_0(\mathbb{C})}\ar@{^{(}->}[r] &
      {\Sup_1(\mathbb{C})}\ar@{^{(}->}[r] & {} \ar@{.}[r] &
      {\Sup_n(\mathbb{C})}\ar@{^{(}->}[r] & {} \ar@{.}[r] & {}
    }
  \end{math}
  of superstructures in $\InfCat$. Of course all colimits are preserved
  by the left adjoint $\inladj$, furthermore ascending chains are
  filtered diagrams and thus, by observation~\ref{inerve.fin.acc}, the
  nerve functor $\inerve$ also preserves the colimits of such chains.
  It follows therefore that the cocone $\arrow \inladj(\inerve(i_n)):
  \inladj(\inerve(\Sup_n(\mathbb{C})))->\inladj(\inerve(\mathbb{C})).$
  displays $\inladj(\inerve(\mathbb{C}))$ as the colimit of the chain:
  \begin{math}
    \xymatrix@R=0ex@C=2em@M=1pt{
      {\inladj(\inerve(\Sup_0(\mathbb{C})))}\ar[r] &
      {\inladj(\inerve(\Sup_1(\mathbb{C})))}\ar[r] & {} \ar@{.}[r] &
      {\inladj(\inerve(\Sup_n(\mathbb{C})))}\ar[r] & 
      {} \ar@{.}[r] & {}
    }
  \end{math}
  Applying the naturality of $\epsilon$ to the inclusions
  $\overinc:\Sup_n(\mathbb{C})-> \mathbb{C}.$ we get a commutative
  square
  \begin{equation*}
    \let\labelstyle=\textstyle
    \xymatrix@C=4em@R=2em{
      {\inladj(\inerve(\Sup_n(\mathbb{C})))}
      \ar[r]^-{\epsilon_{\Sup_n(\mathbb{C})}}
      \ar[d]_{\inladj(\inerve(i_n))} &
      {\Sup_n(\mathbb{C})}
      \ar@{^{(}->}[d]^{i_n} \\
      {\inladj(\inerve(\mathbb{C}))}
      \ar[r]_-{\epsilon_{\mathbb{C}}} & {\mathbb{C}}
    }
  \end{equation*}
  for each $n\in\mathbb{N}$, and this family displays
  $\epsilon_{\mathbb{C}}$ as the \inf-functor induced between these
  colimits by the family of \inf-functors
  $\{\epsilon_{\Sup_n(\mathbb{C})}\}_{n\in\mathbb{N}}$.  However,
  we've already shown that the component of $\epsilon$ on each
  $n$-category $\Sup_n(\mathbb{C})$ is an isomorphism and so the
  induced \inf-functor $\epsilon_{\mathbb{C}}$ is also an isomorphism.
  Applying a dual argument to the ascending chain of superstructures
  of a complicial set $A$ we see that each $\eta_A$ is also an
  isomorphism. It follows that $\inladj\dashv\inerve$ is an adjoint
  equivalence as required.
  
  We now turn to completing our proof by establishing the isomorphism
  depicted in display~(\ref{path.cons.comp}), the construction of
  which relies primarily on the parity complex calculations laid out
  in the next few observations.
\end{obs}

\begin{notation}[intervals]
  In the sequel if $s\leq t$ are integers then we will use the
  following traditional interval notation for sequences of integers:
  \begin{equation*}
    \begin{aligned}\relax
      [s,t] & {} \defeq \{i\in\mathbb{N}\mid s\leq i\leq t\} & \mkern20mu
      (s,t) &  {} \defeq \{i\in\mathbb{N}\mid s< i< t\} \\ 
      [s,t) & {} \defeq \{i\in\mathbb{N}\mid s\leq i< t\} &
      (s,t] &  {} \defeq \{i\in\mathbb{N}\mid s< i\leq t\} 
    \end{aligned}
  \end{equation*}
\end{notation}

\begin{obs}[two point suspensions of parity complexes]
  \label{two.pt.susp}
  If $C$ is a parity complex and $n\in\mathbb{N}$ then we may define a
  parity complex $\Sigma_n C$ with underlying graded set
  \begin{equation*}
    (\Sigma_n C)_0 = \{0,1,...,n\}\mkern40mu
    (\Sigma_n C)_{r+1} = C_r\times\{1,...,n\}
  \end{equation*}
  and face operations:
  \begin{equation*}
    \begin{aligned}
      \pair<x;i>^{-} & {}= \{ i-1 \} & \mkern40mu & \text{and}\\
      \pair<x;i>^{+} & {}= \{ i \} && \text{when $\pair<x;i>\in(\Sigma_n
        C)_1$,} \\
      \pair<x;i>^\parvone & {}= x^\parvone\times\{ i \} && \text{when
      $\pair<x;i>\in(\Sigma_n C)_{r+1}$ with $r>0$.}
    \end{aligned}
  \end{equation*}
  It is a routine matter to check that this is indeed a parity complex
  and that it satisfies all of the extra technical conditions
  discussed in definition~\ref{parity.defn}. In outline, all of the
  basic parity complex axioms and the globularity condition involve
  relationships between 2 or 3 adjacent dimensions so they clearly
  hold for $\Sigma_n C$ above dimension 0, where it is no more than a
  disjoint union of copies of $C$ which have all been shifted up a
  dimension. Instances of these conditions which involve faces at
  dimension 0 all admit simple direct verifications. Finally the
  remaining technical condition, which involves the dimension
  traversing ordering $\filledtri$ on $\Sigma_n C$ and its odd dual,
  follows directly from the simple observation that if
  $\pair<c;i>\filledtri \pair<x';i'>$ in $\Sigma_n C$ then either
  $i<i'$ or we have $i=i'$ and $x\filledtri x'$ in $C$.

  We can extend this construction to a functor $\arrow \Sigma_{*} C:
  \Delta->\Parity.$ which maps each ordinal $[n]\in\Delta$ to $\Sigma_n C$ and
  each simplicial operator $\arrow\alpha:[n]->[m].$ to a parity
  complex morphism $\spanarr\Sigma_\alpha C:\Sigma_n C->\Sigma_m C.$
  given by:
  \begin{equation*}
    \begin{aligned}
      (\Sigma_\alpha C)(i) & {}= \{\alpha(i)\} & \mkern30mu &
      \text{for $i\in (\Sigma_n C)_0$, } \\
      (\Sigma_\alpha C)\pair<x;i> & {}= \{ x \} \times (\alpha(i-1),
      \alpha(i)] && \text{for $\pair<x;i>\in(\Sigma_n C)_{r+1}$ with
        $r\in\mathbb{N}$.}
    \end{aligned}
  \end{equation*}
  Showing that this morphism of graded sets does satisfy the
  conditions given in lemma~\ref{inf.funct.parity} and indeed that
  this construction is functorial in $\alpha\in\Delta$ are matters of
  routine verification which we leave to the reader.  Notice also that
  if $\spanarr f:C->D.$ is a parity complex morphism then we may
  define a morphism of graded sets $\spanarr \Sigma_n f: \Sigma_n
  C->\Sigma_n D.$ by
  \begin{equation*}
    \begin{aligned}
      (\Sigma_n f)(i) & {}= \{ i \} & \mkern30mu &
      \text{for $i\in (\Sigma_n C)_0$, } \\
      (\Sigma_n f)\pair<x;i> & {}= f(x)\times \{ i \} && \text{for
        $\pair<x;i>\in(\Sigma_n C)_{r+1}$ with $r\in\mathbb{N}$.}
    \end{aligned}
  \end{equation*}
  which is again easily shown to be a parity complex morphism and
  functorial in $f\in\Parity$. Combining these observations it is
  clear that we've succeeded in naturally extending our construction
  to a functor $\arrow\Sigma:\Delta\times \Parity->\Parity.$.
  
  To aid the intuition it is appropriate to think of $\Sigma_1 C$ as
  being the two point suspension of $C$ and to view $\Sigma_n C$ as an
  oriented path of $n$ such suspensions abutting at dimension
  $0$, for instance we might picture $\Sigma_3 \OSimp[1]$ as:
  \begin{equation*}
    \xymatrix@R=2ex@C=10em{
      {0} \ar@/^2ex/[r]^{\pair<0;1>}_{}="01"
      \ar@/_2ex/[r]_{\pair<1;1>}^{}="11"
      \ar@{=>}"01";"11"^{\pair<01;1>} &
      {1} \ar@/^2ex/[r]^{\pair<0;2>}_{}="02"
      \ar@/_2ex/[r]_{\pair<1;2>}^{}="12"
      \ar@{=>}"02";"12"^{\pair<01;2>} &
      {2} \ar@/^2ex/[r]^{\pair<0;3>}_{}="03"
      \ar@/_2ex/[r]_{\pair<1;3>}^{}="13"
      \ar@{=>}"03";"13"^{\pair<01;3>} & {3}
    }
  \end{equation*}
  Indeed, in line with this intuition it is possible to show that the
  functor $\arrow\oriental(\Sigma_{*} C):\Delta->\InfCat.$ is actually
  a cocategory ($\mathbb{T}_{\Cat}$-coalgebra).  To do so we apply the
  alternate, finite sketch based, presentation of
  $\mathbb{T}_{\Cat}$-(co)algebras given in
  observation~\ref{le-th.cat} and consider the commutative square
  \begin{equation}\label{orient.sigma.cocat.po}
    \xymatrix@R=2.5em@C=4em{
      {\oriental(\Sigma_0 C)}\ar[r]^{\oriental(\Sigma_{\vertex^p_p} C)}
      \ar[d]_{\oriental(\Sigma_{\vertex^q_0} C)}\ar@{-->}[dr] &
      {\oriental(\Sigma_p C)} 
      \ar@{-->}[d]^{\oriental(\Sigma_{\partinj^{p,q}_1} C)} \\
      {\oriental(\Sigma_q C)}
      \ar@{-->}[r]_-{\oriental(\Sigma_{\partinj^{p,q}_2} C)} & 
      {\oriental(\Sigma_{p+q} C)}
    }
  \end{equation}
  obtained by applying the functor $\oriental(\Sigma_* C)$ to the
  square in display~\eqref{cat.th.sk}. It is clear, from the
  definition of $\Sigma_\alpha C$ for a simplicial operator $\alpha$,
  that the lower horizontal in this square restricts to an isomorphism
  between $\oriental(\Sigma_q C)$ and the sub-\inf-category of
  $\oriental(\Sigma_{p+q} C)$ freely generated by the sub-parity
  complex consisting of those 0-elements $i$ with $i\leq p$ and
  $(r+1)$-elements $\pair<x;i>$ with $i\leq p$. Dually, its right hand
  vertical restricts to an isomorphism between $\oriental(\Sigma_p C)$
  and the sub-\inf-category of $\oriental(\Sigma_{p+q} C)$ freely
  generated by the sub-parity complex of 0-elements $i$ with $i\geq p$
  and $(r+1)$-elements $\pair<x;i>$ with $i> p$. Furthermore, the
  union of these sub-parity complexes is $\Sigma_{p+q} C$ itself and
  the diagonal map in our square restricts to an isomorphism with the
  sub-\inf-category freely generated by the intersection of these
  sub-parity complexes, which consists of the solitary 0-cell
  $\atom{p}$.  Finally, it follows that we may apply
  lemma~\ref{widepo.parity} to show that the square in
  display~\eqref{orient.sigma.cocat.po} is a pushout in $\InfCat$ as
  required.
  
  To identify the functor $\arrow\pathecat_C:\InfCat->\Cat.$ obtained
  by applying Kan's construction (observation~\ref{func.from.coalg})
  to this cocategory we start by studying $\oriental(\Sigma_1 C)$. It
  is immediate, from the definition of $\Sigma_1 C$, that each cell
  $\pair<N;P>$ in the freely generated \inf-category $\oriental(C)$
  gives rise to a corresponding cell $\pair<N;P>^{*}\defeq
  \pair<\{0\}\cup (N\times \{1\}); \{1\}\cup (P\times\{1\})>$ in
  $\oriental(\Sigma_1 C)$ and, indeed, that the only cells of
  $\oriental(\Sigma_1 C)$ that are not of this form are the $0$-cells
  $\atom{0}$ and $\atom{1}$.  Furthermore it is also clear that cells
  $\pair<N;P>$ and $\pair<M;Q>$ are $r$-composable in $\oriental(C)$
  iff $\pair<N;P>^{*}$ and $\pair<M;Q>^{*}$ are $(r+1)$-composable in
  $\oriental(\Sigma_1 C)$ and that in this case we have
  $\pair<M;Q>^{*}\comp_{r+1} \pair<N;P>^{*}=(\pair<M;Q>\comp_r
  \pair<N;P>)^{*}$.  However, if $\mathbb{C}$ is any \inf-category
  then we may combine this analysis with the definition of the
  enriched category $\pathecat(\mathbb{C})\in\InfCatCat$ given in
  observation~\ref{infcat.enr} to show that \inf-functors $\arrow
  f:\oriental(\Sigma_1 C)->\mathbb{C}.$ with $u=f(\atom{0})$ and
  $v=f(\atom{1})$ correspond bijectively to \inf-functors $\arrow
  \hat{f}:\oriental(C)->\pathecat(\mathbb{C} )(u,v).$ where these are
  related by the equality $\hat{f}\pair<N;P>=f\pair<N;P>^{*}$ for each
  cell $\pair<N;P>\in\oriental(C)$.
  
  In order to describe the composition of $\pathecat_C(\mathbb{C})$ in
  these terms we need to study the \inf-functor $\arrow
  w:\oriental(\Sigma_2 C)-> \mathbb{C}.$ that witnesses the composite
  of a pair of arrows $\arrow f,g:\oriental(\Sigma_1 C)->\mathbb{C}.$
  and which, by the definition given in observation~\ref{le-th.cat}, is
  the unique such \inf-functor with $w\circ( \Sigma_{\face^2_2} C)=f$,
  $w\circ(\Sigma_{\face^2_0} C)=g$ and $w\circ(\Sigma_{\face^2_1}
  C)=g\comp f$. Now, a simple calculation demonstrates that if
  $\pair<N;P>$ is a cell in $\oriental(C)$ then we have the equality
  \begin{equation*}
    (\Sigma_{\face^2_0} C)\pair<N;P>^{*} \comp_0 
    (\Sigma_{\face^2_2} C)\pair<N;P>^{*}=
    (\Sigma_{\face^2_1} C)\pair<N;P>^{*}
  \end{equation*}
  in $\oriental(\Sigma_2 C)$, to which we may apply the \inf-functor
  $w$, the witnessing equalities of the last sentence and the defining
  relationship of the last paragraph to show that:
  \begin{equation*}
    \begin{aligned}
      \widehat{g\comp f}\pair<N;P> & {} = (g\comp f)\pair<N;P>^{*} \\
      & {} = w((\Sigma_{\face^2_1} C)\pair<N;P>^{*}) 
      = w((\Sigma_{\face^2_2} C)\pair<N;P>^{*}) \comp_0
      w((\Sigma_{\face^2_0} C)\pair<N;P>^{*}) \\
      & {} = g\pair<N;P>^{*}\comp_0 f\pair<N;P>^{*}
      = \hat{g}\pair<N;P>\comp_0 \hat{f}\pair<N;P>
    \end{aligned}
  \end{equation*}
  In summary, our analysis shows that the objects of
  $\pathecat_C(\mathbb{C})$ correspond to 0-cells of $\mathbb{C}$, its
  homset $\pathecat_C(\mathbb{C})(u,v)$ is isomorphic to the set
  $\InfCat(\oriental(C),\pathecat(\mathbb{C})(u,v))$ and under these
  isomorphisms its composition corresponds to the point-wise
  composition of functors.
  
  In other words, we have established a canonical isomorphism, natural
  in the \inf-category $\mathbb{C}\in\InfCat$, between the category
  $\pathecat_C(\mathbb{C})$ and the one obtained by applying the
  finite limit preserving functor $\InfCat(\oriental(C),{-})$
  point-wise to the $\InfCat$-enriched category
  $\pathecat(\mathbb{C})$.
\end{obs}

\begin{obs}[generalising $\Sigma$ to stratified parity complexes]
  We may extend $\Sigma$ once more to make it into a functor
  $\arrow\Sigma:\Delta\times\SParity->\SParity.$. If $C\in\SParity$
  then we define $\Sigma_n C$ to be the stratified parity complex
  whose underlying parity complex is $\Sigma_n \forget(C)\in\Parity$,
  and in which an element $\pair<x;i>$ is defined to be thin iff $x$
  is thin in $C$. Under this choice of stratification, it is clear
  from the definitions given in the last observation that
  $\Sigma_\alpha C$ preserves thinness for each simplicial operator
  $\alpha\in\Delta$ and that $\Sigma_m f$ preserves thinness whenever
  $f$ is a stratified morphism.
  
  Notice that if $x$ is an element of $C$ then the atoms $\atom{x}$ of
  $\oriental(\forget(C))$ and $\atomp<x;i>$ of $\oriental(\Sigma_1
  \forget(C))$ are related by $\atomp<x;i>= \atom{x}^{*}$, where
  $({-})^{*}$ is the operation on cells which we defined in the last
  observation. It follows, therefore, that an \inf-functor $\arrow
  g:\oriental(\Sigma_1\forget(C))-> \mathbb{C}.$ collapses the atoms
  on thin elements in $\Sigma_1 C$ iff the corresponding \inf-functor
  $\arrow\hat{g}:\oriental(\forget(C)) ->\pathecat(\mathbb{C})(u,v).$
  collapses the atoms on thin elements in $C$. However, we know that
  $\oriental(\Sigma_1 C)$ and $\oriental(C)$ are obtained by taking
  collapsers of these sets of cells, and thus that the bijection of
  the last observation extends to the stratified context.  
  
  Furthermore, we may extend the work of the last observation to show
  that the functor $\arrow\oriental(\Sigma_{*} C):\Delta->\InfCat.$ is a
  cocategory for each stratified parity complex $C$. To be precise, we
  argue that in this generalised context the square depicted in
  display~\eqref{orient.sigma.cocat.po} is constructed by taking
  collapsers of the nodes in the corresponding square for the
  underlying parity complex $\forget(C)$.  However we already know,
  from the last observation, that this latter square is a pushout and
  so it easily follows, from the fact that as colimits the pushout and
  collapser constructions commute, that the square for $C$ is also a
  pushout as required. As before, we adopt the notation $\pathecat_C$
  to denote the Kan functor associated with this cocategory. 
  
  By definition $\Sigma_n \forget(C)$ is the underlying complex of
  $\Sigma_n C$ and so we have a natural family of collapsers $\cover
  q_t:\oriental(\Sigma_n \forget(C))->\oriental(\Sigma_n C).$ and we
  may apply Kan's construction (observation~\ref{func.from.coalg}) to
  this epimorphic coalgebra map to construct a monomorphic natural
  transformation $\overinc:\pathecat_C->\pathecat_{\forget(C)}.$.  For
  each \inf-category $\mathbb{C}$, this identifies
  $\pathecat_C(\mathbb{C})$ with a subcategory of
  $\pathecat_{\forget(C)}(\mathbb{C})$ consisting of those arrows
  $\arrow p:\oriental(\Sigma_1 \forget(C))->\mathbb{C}.$ which
  collapse the atoms on thin elements in $\Sigma_1 C$.
\end{obs}

\begin{obs}[a coalgebraic description of the composite $\ECat{\inerve}\circ\pathecat$]
  \label{coalg.leg.1}
  Consider the functor $\arrow\oriental(\Sigma_{*}\OSimp(-)):
  \tDelta\times\Delta->\InfCat.$ and recall that the last two
  observations demonstrate that
  $\arrow\oriental(\Sigma_{*}\OSimp[n]_?):\Delta->\InfCat.$ is a
  cocategory for each $[n]_?\in\tDelta$. They also demonstrate that
  for each $\mathbb{C}\in\InfCat$ the presheaf
  $\arrow\InfCat(\oriental(\Sigma_1 \OSimp(-)),\mathbb{C}):
  \op{\tDelta}->\Set.$ is naturally isomorphic to the complicial set
  $\coprod_{u,v\in\Sup_0(\mathbb{C})}\inerve(\pathecat(
  \mathbb{C})(u,v))$.  To be precise, this latter fact follows from
  our observation that \inf-functors $\arrow
  g:\oriental(\Sigma_1\OSimp[n]_?)->\mathbb(C).$ with $g(\atom{0})=u$
  and $g(\atom{1})=v$ naturally correspond to \inf-functors
  $\arrow\hat{g}:\oriental(\OSimp[n]_?)-> \pathecat(\mathbb{C})(u,v).$
  which themselves correspond to (thin) $n$-simplices in
  $\inerve(\pathecat(\mathbb{C})(u,v))$ since $\inerve$ is the Kan
  functor on the $\mathbb{T}_{\Comp}$-coalgebra $\arrow\oriental\circ
  \OSimp:\tDelta->\InfCat.$.
  
  It follows therefore that the functor $\arrow\oriental(\Sigma_{*}
  \OSimp(-)):\tDelta\times\Delta->\InfCat.$ is a $\mathbb{T}_{\Comp}
  \otimes\mathbb{T}_{\Cat}$-coalgebra and, as usual, we use the
  notation $\arrow\kan{\oriental(\Sigma_{*}\OSimp(-))}:\InfCat
  ->\Cat(\Comp).$ for the corresponding Kan functor. Applying the
  analysis of observation~\ref{two.pt.susp} again, we know that for
  each $\mathbb{C}\in\InfCat$ the category of (thin) $n$-arrows in the
  complicial category $\kan{\oriental(\Sigma_{*}\OSimp(-))}(
  \mathbb{C})$ is isomorphic to the category obtained by applying the
  representable $\InfCat(\oriental(\OSimp[n]_?),{-})$ point-wise to
  the \inf-category enriched category $\pathecat(\mathbb{C})$.
  Abstracting over $[n]_?\in\tDelta$, we therefore see that
  $\kan{\oriental(\Sigma_{*}\OSimp(-))}( \mathbb{C})$ is isomorphic to
  the complicially enriched category constructed by applying $\inerve$
  point-wise to $\pathecat(\mathbb{C})$. In other words, the
  composite functor $\arrow\ECat{\inerve}\circ\pathecat:\InfCat->
  \CompCat.$ may be represented as the Kan functor on our coalgebra
  $\oriental(\Sigma_{*}\OSimp(-))$.
\end{obs}

\begin{obs}[a coalgebraic description of the composite $\pathcat\rhv\circ\inerve$]\label{coalg.leg.2}
  As discussed in observation~\ref{pathcat.micro}, the prism functor
  $\arrow \pathcat\rhv:\Comp->\Cat(\Comp).$ bears an external
  description as the functor obtained by applying Kan's construction
  (observation~\ref{func.from.coalg}) to the f-almost
  $\mathbb{T}_{\Comp} \otimes\mathbb{T}_{\Cat}$-coalgebra $\arrow
  \Delta\otimes\Cpathcat: \tDelta\times\Delta->\Strat.$.  Of course
  the left adjoint functor $\arrow\inladj:\Strat->\InfCat.$ preserves
  colimits and, by Street's theorem~\ref{nerves.are.complicial}, it
  also carries f-extensions to isomorphisms. So it follows that it
  carries f-almost colimits to colimits and consequently that we may
  apply it point-wise to the f-almost coalgebra
  $\Delta\otimes\Cpathcat$ to give a genuine $\mathbb{T}_{\Comp}
  \otimes\mathbb{T}_{\Cat}$-coalgebra $\arrow
  \inladj\circ(\Delta\otimes\Cpathcat):\tDelta\times\Delta->
  \InfCat.$.  Furthermore, the adjunction $\inladj\dashv\inerve$
  provides us with isomorphisms
  \begin{equation*}
    \begin{aligned}
      \homout{\inladj\circ(\Delta\otimes\Cpathcat)}(\mathbb{C}) & {} \defeq
      \InfCat(\inladj(\Delta({-})\otimes\Cpathcat({*})), \mathbb{C}) \\
      & {} \cong \Strat(\Delta({-})\otimes\Cpathcat({*}),
      \inerve(\mathbb{C})) \defeq \pathcat\rhv(\inerve(\mathbb{C}))
    \end{aligned}
  \end{equation*}
  in $\Al{\mathbb{T}_{\Comp}\otimes\mathbb{T}_{\Cat}}\cong\Cat(\Comp)$
  which are natural in $\mathbb{C}\in\InfCat$.  In other words, we've
  shown that we may represent the composite
  $\arrow\pathcat\rhv\circ\inerve:\InfCat-> \Cat(\Comp).$ as the Kan
  functor on the coalgebra $\inladj\circ(\Delta\otimes\Cpathcat)$.
  
  Now applying corollary~\ref{tensor.simp.th} it is easy to see that
  each \inf-category $\inladj(\Delta[n]_?\otimes\Cpathcat[m])$ is
  canonically isomorphic to $\oriental(\OSimp[n]_?\times
  \Th_1(\OSimp[m]))$ (see observation~\ref{sparity.prod.th} for a
  discussion of the functor $\Th_n$ on $\SParity$) and that these
  isomorphisms are natural in $[n]_?\in\tDelta$ and $[m]\in\Delta$. In
  other words, we could equally well describe $\pathcat\rhv\circ\inerve$
  as the Kan functor on the coalgebra $\arrow\oriental(\OSimp(-)\times
  \Th_1(\OSimp(*))): \tDelta\times\Delta->\InfCat.$.
\end{obs}

\begin{obs}[and its sub-functor $\arrow \pathecat\rhv\circ\inerve:
  \InfCat->\CompCat.$] 
  \label{coalg.leg.2a}
  Let $(\Delta[n]_?\otimes\Delta[1])'$ denote the
  stratified set constructed from $\Delta[n]_?\otimes\Delta[1]$ by
  making thin all non-degenerate $r$-simplices ($r=1,2,...$) of the
  form $\pair<\alpha;\vertex^1_i\circ\eta^r>$. Then we know, from
  observation~\ref{enrich.path.expl}, that if $A$ is a complicial set
  then a (thin) $r$-arrow of $\arrow p:\Delta[n]_?\otimes\Delta[1]->
  A.$ of $\pathcat\rhv(A)$ is in the subcategory $\pathecat\rhv(A)$
  iff it may be lifted to a stratified map with domain
  $(\Delta[n]_?\otimes \Delta[1])'$.
  
  We may define a corresponding stratified parity complex
  $(\OSimp[n]_?\times\OSimp[1])'$ by making thin all elements of
  $\OSimp[n]_?\times\OSimp[1]$ of the form $\pair<v_0...v_p;w_0>$ with
  $p>0$. Arguing as in corollary~\ref{tensor.simp.th} we have a square
  \begin{equation}\label{cprprods.pathecat.sq}
    \let\labelstyle=\textstyle
    \xymatrix@R=2.5em@C=5em{
      {\inladj(\Delta[n]_?\otimes\Delta[1])}
      \ar[r]^{\cong}
      \ar@{-|>}[d]_{\inladj(\subseteq_e)} &
      {\oriental(\OSimp[n]_?\times\OSimp[1])}
      \ar@{-|>}[d]^{\oriental(\subseteq_e)} \\
      {\inladj(\Delta[n]_?\otimes \Delta[1])'}
      \ar@{-->}[r]_{\exists!}^{\cong} &
      {\oriental(\OSimp[n]_?\times\OSimp[1])'}
    }
  \end{equation}
  in which the upper horizontal is the isomorphism of that corollary
  and the verticals are collapsers of the sets of cells associated
  with the extra thin elements in the dashed versions of these
  stratified structures (as discussed in
  observations~\ref{sparity.ext.oriental}
  and~\ref{inf.cat.nerve.strat}). From the description of the upper
  horizontal given in corollary~\ref{tensor.simp.th} it is clear that
  it restricts to a bijection between the sets of cells collapsed by
  these vertical \inf-functors and thus that it induces the (dashed)
  isomorphism at the bottom of this square.
  
  Returning to the last observation, we know that (thin) $n$-arrows of
  $\pathcat\rhv(\inerve(\mathbb{C}))$ correspond naturally to
  \inf-functors $\arrow p:\oriental(\OSimp[n]_?\times\OSimp[1])->
  \mathbb{C}.$. It follows therefore, from the comment in the first
  paragraph of this observation and the square in
  display~\eqref{cprprods.pathecat.sq}, that this \inf-functor
  represents a (thin) $n$-arrow in the subcategory
  $\pathecat\rhv(\inerve(\mathbb{C}))$ iff it factors through the
  collapser $\cover \oriental(\subseteq_e):
  \oriental(\OSimp[n]_?\times\OSimp[1])->\oriental(\OSimp[n]_?\times
  \OSimp[1])'.$.
\end{obs}
  
\begin{obs}[comparing {$C\times\Th_1(\OSimp[m])$} and {$\Sigma_m
    C$}]\label{coalg.comparison} If $C$ is a parity complex, we may
  define a parity complex morphism $\spanarr k^{C,m}:
  C\times\OSimp[m]->\Sigma_m C.$ by:
  \begin{equation*}
      k^{C,m}\pair<x;w_0...w_q> =
      \begin{cases}
        \{ w_0 \} & \text{when $q=0$ and $x\in C_0$,} \\
        \{ x \}\times(w_0,w_1] & \text{when $q=1$,} \\
        \emptyset & \text{otherwise.}
      \end{cases}
  \end{equation*}
  This clearly respects dimensions and it is a matter of simple
  case-wise verification to check that it satisfies the movement
  condition given in lemma~\ref{inf.funct.parity}.  In a similar
  fashion it is also easily demonstrated that this family of parity
  complex morphisms is natural in $C\in\Parity$ and $[m]\in\Delta$. We
  leave the completely routine verification of these facts for the
  delectation of the reader.
  
  Having defined this parity complex morphism, it is worth spending a
  few moments studying the action of the associated \inf-functor
  $\arrow\oriental(k^{C,m}):\oriental(C\times \OSimp[m])->
  \oriental(\Sigma_m C).$ on the atoms of its domain.  As we may
  easily establish directly, or verify by arguing as
  in observation~\ref{osimp(alpha).act.on.atoms} using the result of
  corollary~\ref{morp.action.1}, we have
  \begin{itemize}
  \item $\oriental(k^{C,m})(\atomp<x;w_0>)=\atom{w_0}$, and 
    
  \item $\oriental(k^{C,m})(\atomp<x;w_0...w_q>)= \atomp<x;w_q>
    \comp_0 \atomp<x;w_q-1>\comp_0 \cdots\comp_0\atomp<x;w_0+1>$
    whenever $q>0$ and in particular if $q=1$ and $w_1=w_0+1$ 
    this is equal to the atom $\atomp<x;w_1>$.
  \end{itemize}
  
  Now suppose that $C$ is a stratified parity complex and consider the
  product $C\times\Th_1(\OSimp[m])$ in $\SParity$.  By definition, the
  element $\pair<x;w_0...w_q>$ is thin in $C\times\Th_1(\OSimp[m])$
  iff either $w_0...w_q$ is thin in $\Th_1(\OSimp[n])$, which happens
  when $q>1$ in which case $k^{\forget(C),m}\pair<x;w_0...w_q>$ is
  empty, or $x$ is thin in $C$, in which case its dimension is $>0$
  and so either $k^{\forget(C),m}\pair<x;w_0...w_q>$ is empty or $q=1$
  and it is equal to the set $\{x\}\times(w_0,w_1]$ of thin elements
  in $\Sigma_m C$.  It follows that $k^{\forget(C),m}$ satisfies the
  thinness preservation condition making it into a stratified morphism
  $\spanarr k^{C,m}:C\times\Th_1(\OSimp[m])->\Sigma_m C.$. Forgetting
  stratifications, we see that the naturality of this stratified
  morphism in $C\in\SParity$ and $[n]\in\Delta$ follows directly from
  its established naturality on underlying parity complexes.
\end{obs}

\begin{lemma}
  \label{k.iso.stuff}
  The stratified morphism $\spanarr k:\OSimp[n]_?\times\OSimp[1]->
  \Sigma_1 \OSimp[n]_?.$ of the last observation lifts to one with
  domain $(\OSimp[n]_?\times\OSimp[1])'$ and the \inf-functor obtained
  by applying $\oriental$ to this lifted morphism provides an
  isomorphism between $\oriental(\OSimp[n]_?\times\OSimp[1])'$ and
  $\oriental(\Sigma_1 \OSimp[n]_?)$.
\end{lemma}

\begin{proof}
  The first part simply follows from the definition of $k$ which tells
  us that $k\pair<v_0...v_p;w_0>=\emptyset$ whenever $p>0$. 
  
  To prove the remainder, first consider the case $[n]_?=[n]$ and
  recall, from the last observation, that
  $\arrow\oriental(k):\oriental( \OSimp[n]\times \OSimp[1])'->
  \oriental(\Sigma_1\OSimp[n]).$ maps
  $\atomp<0;w_0>\in\oriental(\OSimp[n]\times\OSimp[1])'$ to
  $\atom{w_0}\in\oriental(\Sigma_1\OSimp[n])$ and $\atomp<\svec{v};
  01>\in\oriental(\OSimp[n]\times\OSimp[1])'$ to
  $\atomp<\svec{v};1>\in \oriental(\Sigma_1\OSimp[n])$. It follows
  that this \inf-functor restricts to a bijection between the subset
  \begin{math}
    S=\{ \atomp<0;w_0>\mid w_0=0,1 \} \cup \{\atomp<\svec{v};01> \mid 
    \svec{v}\in\OSimp[n]\}
  \end{math}
  of cells in $\oriental(\OSimp[n]\times\OSimp[1])'$ and the set of
  atoms of the freely generated category
  $\oriental(\Sigma_1\OSimp[n])$. Notice however that any element
  $\pair<v_0v_1;w_0>$ in $(\OSimp[n]\times\OSimp[1])'$ is thin and so
  it follows that $\atomp<v_0v_1;w_0>$ is a 0-cell in
  $\oriental(\OSimp[n]\times \OSimp[1])'$ and therefore that its
  0-source $\atomp<v_0;w_0>$ and 0-target $\atomp<v_1;w_0>$ are equal.
  From this it follows that any cell $\atomp<v_0;w_0>$ is in $S$ and
  thus that this subset is a weak generator for
  $\oriental(\OSimp[n]\times\OSimp[1])'$ since it contains the weak
  generator $\atom{(\OSimp[n]\times\OSimp[1])'}$ discussed in
  observation~\ref{sparity.ext.oriental}. Now we may apply
  lemma~\ref{fg.isom.lemma} to show that
  $\arrow\oriental(k):\oriental( \OSimp[n]\times \OSimp[1])'->
  \oriental(\Sigma_1\OSimp[n]).$ is an isomorphism as required.

  To extend this to the case $[n]_?=[n]_t$ consider the commutative
  square
  \begin{equation*}
    \let\labelstyle=\textstyle
    \xymatrix@R=2.5em@C=5em{
      {\oriental(\OSimp[n]\times\OSimp[1])'}
      \ar[r]^-{\oriental(k)}_-{\cong}
      \ar@{-|>}[d]_{\oriental(\subseteq_e)} &
      {\oriental(\Sigma_1\OSimp[n])}
      \ar@{-|>}[d]^{\oriental(\subseteq_e)} \\
      {\oriental(\OSimp[n]_t\times\OSimp[1])'}
      \ar[r]_-{\oriental(k)} &
      {\oriental(\Sigma_1\OSimp[n]_t)}
    }
  \end{equation*}
  and again, as discussed in observation~\ref{sparity.ext.oriental},
  the verticals here are collapsers of the cells $\atomp<01...n;01>\in
  \oriental(\OSimp[n]\times\OSimp[1])'$ and $\atomp<01...n;1>\in
  \oriental(\Sigma_1\OSimp[n])$ respectively. The horizontal
  isomorphism at the top of this square maps these cells to each other
  and it follows that we may apply the collapser property of the
  right hand vertical to induce an inverse to the lower horizontal 
  as required.
\end{proof}

Finally we arrive at our ultimate destination.

\begin{thm}[the Street-Roberts conjecture]\label{street.roberts.thm}
  We may construct the essentially commutative square illustrated in
  display~(\ref{path.cons.comp}) and thus apply the argument described
  in observation~\ref{st.rob.prf.outline} to show that Street's nerve
  functor $\arrow\inerve:\InfCat->\Comp.$ provides an equivalence
  between the categories of \inf-categories and complicial sets.
\end{thm}

\begin{proof}
  Observation~\ref{coalg.leg.1} showed that we may represent $\arrow
  \ECat{\inerve}\circ\pathecat:\InfCat->\CompCat.$ as the Kan functor
  associated with a coalgebra $\arrow \Sigma_{*} \OSimp(-):
  \tDelta\times\Delta->\InfCat.$. Similarly
  observation~\ref{coalg.leg.2} demonstrated that the composite
  $\arrow \pathcat\rhv\circ\inerve:\InfCat->\Cat(\Comp).$ may be
  represented by a coalgebra $\arrow\OSimp(-)\times\Th_1(\OSimp(*)):
  \tDelta\times\Delta->\InfCat.$. Furthermore,
  observation~\ref{coalg.comparison} provides us with a coalgebra map
  from the latter to the former which we can use to construct a
  comparison 2-cell
  \begin{equation}\label{path.cons.comp2}
    \let\labelstyle=\textstyle
    \xymatrix@R=5em@C=6em{
      {\InfCat}\ar[r]^-{\pathecat}
      \ar[d]_{\inerve}^{}="4" &
      {\InfCatCat}\ar@{}[d]_{}="3" \\ {\Comp} & {}
      \save []+<0em,+2.5em>*+{\CompCat}="1"\restore
      \save []+<-3em,0em>*+{\Cat(\Comp)}="2"\restore
      \ar"2,1";"2"_-{\pathcat\rhv}
      \ar"1,2";"1"^{\ECat\inerve} 
      \ar@{^{(}->}"1";"2"
      \ar@{=>}"3"-<4em,0em>;"4"+<4em,0em>
    }
  \end{equation}
  by pre-composition, as discussed in
  observation~\ref{func.from.coalg}. Now applying
  observation~\ref{coalg.leg.2a} and lemma~\ref{k.iso.stuff}, we see
  that for each \inf-category $\mathbb{C}$ the component $\overarr:
  \ECat{\inerve}(\pathecat(\mathbb{C}))->\pathcat\rhv(\inerve(
  \mathbb{C})).$ of this acts as an isomorphism between the complicial
  set of arrows of $\ECat{\inerve}(\pathecat(\mathbb{C}))$ and that of
  the subcategory $\pathecat\rhv(\inerve(\mathbb{C}))$ of
  $\pathcat\rhv(\inerve(\mathbb{C}))$, and thus restricts to an
  isomorphism of those categories. In other words, the natural
  transformation in display~\eqref{path.cons.comp2} restricts to
  provide the isomorphism of display~\eqref{path.cons.comp} as
  required.
\end{proof}


\newpage
\bibliographystyle{abbrv}
\bibliography{cattheory}

\end{document}